\documentclass[notitlepage]{book}

\usepackage{etex}
\usepackage{makeidx}
\makeindex

\usepackage[usenames,dvipsnames]{xcolor}
\usepackage[bookmarksopen,bookmarksdepth=2]{hyperref}
\hypersetup{colorlinks=true,citecolor=NavyBlue,linkcolor=BrickRed,urlcolor=Green} 
\usepackage{tikz-cd}

\RequirePackage{amsmath}
\RequirePackage{amssymb}
\RequirePackage{amsxtra}
\RequirePackage{amsfonts}
\RequirePackage{latexsym}
\RequirePackage{euscript}
\RequirePackage{amscd}
\RequirePackage{amsthm}
\RequirePackage{xypic}
\usepackage{stmaryrd}
\usepackage{mathtools}

\usepackage{mathrsfs}
\xyoption{all}

\usepackage{enumitem}
\DeclareMathOperator{\Max}{Max}
\def\acong{\buildrel \text{almost }  \over \cong}
\renewcommand{\mathbb}{\mathbf}

\newcommand{\BT}{\mathrm{BT}}
\newcommand{\Lotimes}{{\otimes}^{\mathbb L} \, }
\newcommand{\Acirc}{A^\circ}
\newcommand{\Bcirc}{B^\circ}
\newcommand{\Ccirc}{C^\circ}

\newcommand{\Gmhat}{\widehat{\mathbb G}_m}
\newcommand{\cYha}{\cY_{h}^a}

\newcommand{\AKnrhat}{\widehat{\A}_K^{\nr}}
\newcommand{\AKnrhatA}{\widehat{\A}_{K,A}^{\nr}}
\newcommand{\AKnrhatcK}{\widehat{\A}_{K,\cK}^{\nr}}

\newcommand{\lambdau}{\underline{\lambda}}
\newcommand{\bE}{\mathbf{E}}

\newcommand{\Kbasic}{K^{\mathrm{basic}}}
\newcommand{\piflat}{\pi^\flat}
\newcommand{\piflatprime}{(\pi')^\flat}

\newcommand{\Sets}{\underline{{Sets}}}
\newcommand{\cyc}{\operatorname{cyc}}
\newcommand{\BK}{\operatorname{BK}}
\newcommand{\Kcyc}{K_{\cyc}}

\newcommand{\Aplus}{\A^+}
\newcommand{\AAA}{\A}

\newcommand{\Cone}{\operatorname{Cone}}
\newcommand{\Fib}{\operatorname{Fib}}
\newcommand{\Lift}{\operatorname{Lift}}
\newcommand{\pro}{{\operatorname{pro\, -}}}
\def\im{\mathop{\mathrm{im}}\nolimits}
\newcommand{\negligible}{\mathrm{small}}
\newcommand{\OEA}{\cO_{\mathcal{E},A}}
\newcommand{\OEApiflat}{\cO_{\mathcal{E},\pi^\flat,A}}
\newcommand{\hOEA}{\widehat{\cO}_{\mathcal{E},A}}

\newcommand{\tOEA}{\widetilde{\cO}_{\mathcal{E},A}}

\newcommand{\OEB}{\cO_{\mathcal{E},B}}

\newcommand{\cXbar}{\overline{\cX}}

\newcommand{\Frac}{\operatorname{Frac}}

\newcommand{\gMt}{\gM^{\inf}}
\newcommand{\gNt}{\gN^{\inf}}
\newcommand{\gFt}{\gF^{\inf}}

\newcommand{\ybar}{\overline{y}}
\newcommand{\Ghat}{\widehat{G}}

\usepackage{dsfont}

\def\A{\mathbb A}

\def\C{\mathbb C}
\def\F{\mathbb F}
\def\Khat{\widehat{K}}

\def\Lhat{\widehat{L}}

\def\Q{\mathbb{Q}}

\def\Z{\mathbb{Z}}
\def\Fbar{\overline{\F}}

\def\Qbar{\overline{\Q}}

\def\Zbar{\overline{\Z}}
\def\Zhat{\widehat{\Z}}

\def\m{\mathfrak m}
\newcommand{\st}{\mathrm{st}}

\def\Gr{\mathrm{Gr}}

\def\chibar{\overline{\chi}}

\def\id{\mathrm{id}}

\def\alg{\mathrm{alg}}

\def\red{\mathrm{red}}

\def\ab{\mathrm{ab}}
\def\nr{\mathrm{ur}}
\def\fl{\mathrm{fl}}

\def\ss{\mathrm{ss}}

\def\ur{\mathrm{ur}}

\def\nil{\mathop{\mathrm{nil}}\nolimits}

\def\GL{\operatorname{GL}}

\def\Gal{\mathrm{Gal}}

\def\Aut{\mathrm{Aut}}

\def\Ext{\mathrm{Ext}}

\def\End{\mathrm{End}}

\def\Hom{\mathop{\mathrm{Hom}}\nolimits}
\def\RHom{\mathop{\mathrm{RHom}}\nolimits}

\def\Spec{\mathop{\mathrm{Spec}}\nolimits}
\def\Spf{\mathop{\mathrm{Spf}}\nolimits}
\def\Frob{\mathop{\mathrm{Frob}}\nolimits}
\def\Supp{\mathop{\mathrm{Supp}}\nolimits}

\def\Ind{\mathop{\mathrm{Ind}}\nolimits}

\def\Fil{\mathop{\mathrm{Fil}}\nolimits}

\def\rhobar{\overline{\rho}}

\def\cotimes{\operatorname{\widehat{\otimes}}}

\def\crys{\mathrm{crys}}

\def\dR{\mathrm{dR}}
\def\pst{\mathrm{pst}}

\def\WD{\mathrm{WD}}

\def\m{\mathfrak{m}}

\def\iso{\buildrel \sim \over \longrightarrow}

\swapnumbers

\newcommand{\onto}{\twoheadrightarrow}

\newcommand{\into}{\hookrightarrow}

\newcommand{\To}{\longrightarrow}
\newcommand{\isoto}{\stackrel{\sim}{\To}}
\newcommand{\isofrom}{\stackrel{\sim}{\longleftarrow}}

\newcommand{\textD}{\mathrm{D}}

\newcommand{\textB}{\mathrm{B}}
\newcommand{\Bcris}{\textB_{\cris}}

 \newcommand{\BdR}{\textB_{\dR}}

\newcommand{\Bst}{\textB_{\st}}

\newcommand{\DdR}{\textD_{\dR}}

\newlength{\ownl}

\newcommand{\ad}{{\operatorname{ad}\,}}

\newcommand{\Fitt}{{\operatorname{Fitt}\,}}

\newcommand{\Id}{{\operatorname{Id}}}
\renewcommand{\Im}{{\operatorname{Im}\,}}

\newcommand{\rec}{{\operatorname{rec}}}

\newcommand{\rk}{{\operatorname{rk}\,}}

\newcommand{\tr}{{\operatorname{tr}\,}}

\newcommand{\Gm}{{\mathbb{G}_m}}

\newcommand{\cris}{{\operatorname{cris}}}

\newcommand{\disc}{{\operatorname{disc}}}

\newcommand{\semis}{{\operatorname{ss}}}

\newcommand{\univ}{{\operatorname{univ}}}
\newcommand{\B}{{\mathbb{B}}}

\newcommand{\E}{{\mathbb{E}}}
\newcommand{\G}{{\mathbb{G}}}

\newcommand{\cC}{\mathcal{C}}
\renewcommand{\cD}{\mathcal{D}}
\newcommand{\cE}{\mathcal{E}}
\newcommand{\cF}{\mathcal{F}}
\newcommand{\cG}{\mathcal{G}}
\renewcommand{\cH}{\mathcal{H}}

\newcommand{\cI}{\mathcal{I}}

\newcommand{\cK}{\mathcal{K}}

\newcommand{\cO}{\mathcal{O}}
\renewcommand{\O}{\cO}

\renewcommand{\cR}{\mathcal{R}}

\newcommand{\cT}{\mathcal{T}}
\newcommand{\cU}{\mathcal{U}}

\newcommand{\cW}{\mathcal{W}}
\newcommand{\cWW}{\mathcal{W}}
\newcommand{\cX}{\mathcal{X}}
\newcommand{\cY}{\mathcal{Y}}
\newcommand{\cZ}{\mathcal{Z}}

\newcommand{\gF}{{\mathfrak{F}}}

\newcommand{\gL}{{\mathfrak{L}}}
\newcommand{\gM}{{\mathfrak{M}}}
\newcommand{\gN}{{\mathfrak{N}}}

\newcommand{\gP}{{\mathfrak{P}}}
\newcommand{\gQ}{{\mathfrak{Q}}}

\newcommand{\gS}{{\mathfrak{S}}}

\newcommand{\gu}{{\mathfrak{u}}}

\newcommand{\barK}{\overline{{K}}}

\newcommand{\bark}{\overline{{k}}}

\newcommand{\tA}{\widetilde{{\A}}}

\newcommand{\tB}{\widetilde{{B}}}
\newcommand{\tC}{\widetilde{{C}}}

\newcommand{\tM}{\widetilde{{M}}}

\newcommand{\tT}{\widetilde{{T}}}

\newcommand{\tX}{\widetilde{{X}}}

\newcommand{\tZ}{\widetilde{{Z}}}

\newcommand{\tb}{\widetilde{{b}}}
\newcommand{\tc}{\widetilde{{c}}}

\newcommand{\te}{\widetilde{{e}}}
\newcommand{\tif}{\widetilde{{f}}}

\newcommand{\tx}{\widetilde{{x}}}

\newcommand{\tcF}{\widetilde{{\cF}}}
\newcommand{\tcG}{\widetilde{{\cG}}}
\newcommand{\tcH}{\widetilde{{\cH}}}

\newcommand{\tcK}{\widetilde{{\cK}}}

\newcommand{\alphabar   }{\overline{\alpha  }}
\newcommand{\betabar         }{\overline{\beta}}

\newcommand{\epsilonbar    }{\overline{\epsilon}}

 \newcommand{\thetabar    }{\overline{\theta}}

 \newcommand{\sigmabar   }{\overline{\sigma}}

 \newcommand{\psibar   }{\overline{\psi}}

 \newcommand{\tGamma     }{\widetilde{\Gamma}}
 \newcommand{\Gammat     }{\widetilde{\Gamma}}

\newcommand{\gammat   }{\widetilde{\gamma}}

 \newcommand{\varphit   }{\widetilde{\varphi}}

\def\RCS$#1: #2 ${\expandafter\def\csname RCS#1\endcsname{#2}}
\RCS$Revision: 85 $
\RCS$Date: 2011-12-16 18:54:17 -0600 (Fri, 16 Dec 2011) $

\newcommand{\s}{\mathcal{S}} 
 
\newcommand{\mf}{\mathfrak}

\newcommand{\rbar}{{\overline{r}}}

\newcommand{\HT}{\operatorname{HT}}

 \newcommand{\Qp}{{\Q_p}}

\newcommand{\Zp}{{\Z_p}}
 
\newcommand{\Qpbar}{{\overline{\Q}_p}}

\newcommand{\Zpbar}{{\overline{\Z}_p}}

\newcommand{\Fpbar}{{\overline{\F}_p}}
\newcommand{\Fpbartimes}{{\overline{\F}_p^\times}}
\newcommand{\Fnew}{\F'}

\newcommand{\Fp}{{\F_p}}

\newif\iffinalrun
\iffinalrun
\else
\fi

\iffinalrun
  \newcommand{\need}[1]{}
  \newcommand{\mar}[1]{}
\else
  \newcommand{\need}[1]{{\tiny *** #1}}
  \newcommand{\mar}[1]{\marginpar{\raggedright\tiny fixme #1}}
\fi

\newcommand{\Ainf}{\mathbf{A}_{\inf}}
\newcommand{\AAinf}[1]{\mathbf{A}_{\inf,#1}}

\newtheorem{theorem}[subsection]{Theorem}
\newtheorem{thm}[subsection]{Theorem}
\newtheorem{lemma}[subsection]{Lemma}
\newtheorem{lem}[subsection]{Lemma}

\newtheorem{cor}[subsection]{Corollary}
\newtheorem{conj}[subsection]{Conjecture}

\newtheorem{prop}[subsection]{Proposition}

\newtheorem{atheorem}[section]{Theorem}
\newtheorem{athm}[section]{Theorem}
\newtheorem{alemma}[section]{Lemma}
\newtheorem{alem}[section]{Lemma}

\newtheorem{acor}[section]{Corollary}
\newtheorem{aprop}[section]{Proposition}

\theoremstyle{definition}
\newtheorem{df}[subsection]{Definition}
\newtheorem{defn}[subsection]{Definition}
\newtheorem{situation}[subsection]{Situation}

\newtheorem{adefn}[section]{Definition}
\theoremstyle{remark}
\newtheorem{remark}[subsection]{Remark}
\newtheorem{rem}[subsection]{Remark}
\newtheorem{aremark}[section]{Remark}
\newtheorem{arem}[section]{Remark}

\newtheorem{example}[subsection]{Example}
\newtheorem{aexample}[section]{Example}

\newtheorem{hyp}[subsection]{Hypothesis}

\def\numequation{\addtocounter{subsection}{1}\begin{equation}}
\def\nummultline{\addtocounter{subsection}{1}\begin{multline}}
\def\anumequation{\addtocounter{section}{1}\begin{equation}}
\renewcommand{\theequation}{\arabic{chapter}.\arabic{section}.\arabic{subsection}}

\newenvironment{dedication}
  {
   \thispagestyle{empty}
   \vspace*{\stretch{1}}
   \itshape             
   \raggedleft          
  }
  {\par 
   \vspace{\stretch{3}} 
   \clearpage           
  }

\title
{Moduli
  stacks of \'etale $(\varphi,\Gamma)$-modules and the existence of
  crystalline lifts}
\author{
Matthew Emerton
  \and
Toby Gee
}
\date{}
  
\begin{document}
  


\frontmatter
Moduli stacks of \'etale $(\varphi,\Gamma)$-modules and the existence
of crystalline lifts

\cleardoublepage
\thispagestyle{empty}

\maketitle
 \chapter{Dedication}
  \begin{dedication}
In memory of Jean-Marc Fontaine.
\end{dedication}

\thispagestyle{empty}
\cleardoublepage
\thispagestyle{empty}


\setcounter{tocdepth}{1}
\tableofcontents
  
\mainmatter             


\chapter{Introduction}
In this book we construct moduli stacks of \'etale $(\varphi,\Gamma)$-modules
(projective, of some fixed rank, 
and with coefficients in $p$-adically complete rings),
and establish some of their basic properties.  We also present some first
applications of this construction to the theory of Galois representations.

\section{Motivation}\label{subsec: motivation}Mazur's theory of
deformations of Galois representations~\cite{MR1012172}
is modelled on the geometric study of
infinitesimal neighbourhoods of points in 
moduli spaces
via 
formal deformation theory.
In the mid-2000s,
Kisin suggested 
that some kind of
moduli spaces of local Galois representations should
exist; that is, there should be formal algebraic stacks
over~$\Zp$ whose
closed points correspond to representations
$\rhobar:G_K\to\GL_d(\Fpbar)$, and
whose versal rings at such points
should recover appropriate Galois deformation rings. This
expectation is borne out by the results of this book. (In fact, Kisin
was motivated by calculations of crystalline deformation rings
for~$\GL_2(\Qp)$ that had been carried out by Berger--Breuil using the
$p$-adic local Langlands correspondence, and suggested that the versal
rings should give crystalline deformation rings.  Thus his suggestion
is realized by the stacks~$\cX_d^{\crys,\underline{\lambda}}$ of
Theorem~\ref{thm:intro statement on crystalline moduli} below.)

A natural way to construct such a stack would be to  consider a
literal moduli stack of continuous representations
$\rho:G_K\to\GL_d(A)$, for $K$ a $p$-adic field
and $A$ a $p$-adically complete $\Z_p$-algebra; indeed such stacks
were constructed by Carl Wang-Erickson~\cite{MR3831282}. However,
the stacks constructed in this way are less ``global'' than one would
wish, and in particular the corresponding families of mod~$p$
representations $\rhobar:G_K\to\GL_d(\Fpbar)$ have constant semisimplification.

In this book, we instead consider moduli stacks of \'etale
$(\varphi,\Gamma)$-modules. These
contain Wang-Erickson's stacks as substacks, and coincide with them on
the level of $\Fpbar$-points, but their geometry is quite different; in particular,
we see much larger families,  
exhibiting some unexpected features (for example, irreducible
representations arising as limits of reducible representations). The
relationship between the theory and constructions that we develop
here, and the usual formal deformation theory of Galois
representations, is the same as that between the theory of moduli
spaces of algebraic varieties, and the formal deformation theory of
algebraic varieties: the latter gives valuable local information about
the former, but moduli spaces, when they can be constructed, capture
global aspects of the situation 
inaccessible to the purely
infinitesimal tools of formal deformation theory.



\section{Our main theorems}
Our goal in this book is to construct, and establish the basic
properties of, moduli stacks of \'etale $(\varphi,\Gamma)$-modules.
More precisely, if we fix a finite extension $K$ of $\Q_p$, and a
non-negative integer $d$ (the rank), then we let $\cX_d$ denote the
category fibred in groupoids over $\Spf \Z_p$ whose groupoid of
$A$-valued points, for any $p$-adically complete $\Z_p$-algebra $A$,
is equal to the groupoid of rank~ $d$ projective
\'etale $(\varphi,\Gamma)$-modules with $A$-coefficients.
(See Section~~\ref{subsec: phi Gamma
  coefficients}
below for a definition of these.) 
Our first main theorem is the following. (See
Corollary~\ref{cor:Xd is formal algebraic} and
Theorem~\ref{thm:reduced dimension}.)

\begin{thm}\label{thm:intro statement of basic properties of X}
	The category fibred in groupoids $\cX_d$ is a 
	Noetherian formal algebraic
	stack.  Its underlying reduced substack~$\cX_{d,\red}$
	{\em (}which is an algebraic stack{\em )}
	is of finite type
	over $\F_p$, and is equidimensional of dimension
	$[K:\Q_p] d(d-1)/2$. The irreducible components of 
	$\cX_{d,\red}$
	admit a natural labelling by Serre
        weights. 
\end{thm}

%

We will elaborate on the labelling of components by Serre weights
further below.  For now, we mention that,
under the usual
correspondence between \'etale $(\varphi,\Gamma)$-modules
and Galois representations,
the groupoid of $\Fpbar$-points of~$\cX_d$,
which coincides with the groupoid of $\Fpbar$-points
of the underlying reduced substack $\cX_{d,\red}$,
is naturally equivalent to the groupoid of continuous representations
$\rhobar:G_K\to\GL_d(\Fpbar)$.  (More generally, if $A$ is any finite $\Z_p$-algebra,
then the groupoid $\cX_d(A)$ is canonically equivalent
to the groupoid of continuous representations $G_K \to \GL_d(A)$.) It is expected 
 that our labelling of the irreducible components can be refined (by
 adding further labels to some of the components) to
 give a geometric description of the weight part of Serre's
 conjecture, so that $\rhobar$ corresponds to a point in a component
of $\cX_{d,\red}$ which is labeled by the Serre weight
$\underline{k}$ if and only if $\rhobar$ admits $\underline{k}$ as a
Serre weight; we discuss this expectation, and what is known
about it, in Section~\ref{subsec:
  geometric BM} below (and in more detail in Chapter~\ref{sec: BM}).

Again using the
correspondence between \'etale $(\varphi,\Gamma)$-modules
and Galois representations,
we see that the universal lifting ring of a 
representation $\rhobar$ as above will provide a versal ring
to~$\cX_d$ at the corresponding $\Fbar_p$-valued
point. Accordingly we expect that the stacks~$\cX_d$ will have
applications to the study of Galois representations and
their deformations. As a first
example of this, we prove the following result on
the existence of crystalline lifts; although the statement of this
theorem involves a fixed~$\rhobar$, we do not know how to prove it
 without using the stacks~$\cX_d$, over which~$\rhobar$ varies.
\begin{thm}[Theorem~\ref{thm: strong existence of crystalline lifts}]
	\label{thm:intro crystalline lifts}
	If $\rhobar: G_K \to \GL_d(\Fbar_p)$ is a continuous
	representation, then $\rhobar$
	has a lift $\rho^{\circ}: G_K \to \GL_d(\Zbar_p)$ for which
	the associated $p$-adic representation $\rho:G_K \to \GL_d(\Qbar_p)$
	is crystalline of regular Hodge--Tate weights. 

We can furthermore 
        ensure that~$\rho^{\circ}$ is potentially diagonalizable.
\end{thm}
(The notion of a potentially diagonalizable representation was
introduced in~\cite{BLGGT}, and is recalled as Definition~\ref{defn:
  pot diag} below.) In combination with
potential automorphy 
 theorems, this has the following application to the globalization of
local Galois representations. 

\begin{thm}[Corollary~\ref{cor: existence of
  global lift}]
  \label{thm: intro existence of global lifts}Suppose that~$p\nmid 2d$, and
  fix $\rhobar:G_K\to\GL_d(\Fpbar)$. Then there is an imaginary CM
  field~$F$ and an irreducible conjugate self dual automorphic Galois representation
  $\rbar:G_F\to\GL_d(\Fpbar)$ such that  for every $v|p$, we have
  $F_v\cong K$ and either~$\rbar|_{G_{F_v}}\cong \rhobar$ or~$\rbar|_{G_{F_{v^c}}}\cong \rhobar$.
\end{thm}

Another key result of the book is the following theorem,
describing moduli stacks of \'etale $(\varphi,\Gamma)$-modules corresponding to
{\em crystalline} and {\em semistable} Galois representations.

\begin{thm}[Theorem~\ref{thm: existence of ss stack}]\label{thm:intro statement on crystalline moduli}
	If $\underline{\lambda}$ is a collection
	of labeled Hodge--Tate weights,
	and if $\cO$ denotes the ring of integers
	in a finite extension $E$ of $\Q_p$ containing the Galois
	closure of $K$ {\em(}which will serve as the ring of coefficients{\em )},
	then there is a closed substack $\cX_d^{\crys,\underline{\lambda}}$
	of $(\cX_{d})_{\cO}$ which is a $p$-adic formal algebraic stack and is flat 
	over~$\cO$, and which 
	is characterized as being the unique closed substack
	of $(\cX_{d})_{\cO}$ which is flat over $\cO$ and whose groupoid 
	of $A$-valued points,
	for any finite flat $\cO$-algebra~$A$,
	is equivalent
	{\em (}under the equivalence between \'etale
	$(\varphi,\Gamma)$-modules and continuous
	$G_K$-representations{\em )} to the groupoid of continuous
	representations
	$G_K\to \GL_d(A)$ which become {\em crystalline}
	after extension of scalars to $A\otimes_{\cO} E$,
	and whose labeled Hodge--Tate weights are
	equal to~$\underline{\lambda}$.

	Similarly, 
	there is a closed substack $\cX_d^{\semis,\underline{\lambda}}$
	of $(\cX_{d})_{\cO}$ which is a $p$-adic formal algebraic stack and is flat 
	over~$\cO$, and which 
	is characterized as being the unique closed substack
	of $(\cX_{d})_{\cO}$ which is flat over $\cO$ and whose groupoid 
	of $A$-valued points,
	for any finite flat $\cO$-algebra~$A$,
	is equivalent
	to the groupoid of continuous representations
	$G_K\to \GL_d(A)$ which become {\em semistable}
	after extension of scalars to $A\otimes_{\cO} E$,
	and whose labeled Hodge--Tate weights are
	equal to $\underline{\lambda}$.
\end{thm}

\begin{remark}
  In fact, Theorem~\ref{thm: existence of ss stack} also proves
  the analogous result for potentially crystalline and potentially
  semistable representations of arbitrary inertial type, but for
  simplicity of exposition we restrict ourselves to the crystalline
  and semistable cases in this introduction.
\end{remark}

A crucial distinction between the stacks $\cX_d$ and their
closed substacks
$\cX_d^{\crys,\underline{\lambda}}$
and
$\cX_d^{\semis,\underline{\lambda}}$
is that while $\cX_d$ is a formal algebraic stack lying over $\Spf \Z_p$,
it is not actually a $p$-adic formal algebraic stack (in the sense of Definition~\ref{adefn: p adic formal algebraic stack}); 
see Proposition~\ref{prop: X is not padic formal algebraic}.
On the other hand,
the stacks
$\cX_d^{\crys,\underline{\lambda}}$
and
$\cX_d^{\semis,\underline{\lambda}}$
{\em are} $p$-adic formal algebraic stacks, which implies 
that their mod $p^a$ reductions are in fact algebraic stacks.
This gives in particular a strong interplay between the structure of the
mod~ $p$ fibres of crystalline and semistable
lifting rings and the geometry of the underlying reduced substack~$\cX_{d,\red}$.  This plays an important role in determining the structure of this reduced
substack, and also in the proof of
Theorem~\ref{thm:intro crystalline lifts}.
As we explain in more detail in Section~\ref{subsec:
  geometric BM} below, 
it also allows us to reinterpret the Breuil--M\'ezard conjecture
in terms of the interaction between the
structure of the mod $p$ fibres of the stacks
$\cX_d^{\crys,\underline{\lambda}}$
and
$\cX_d^{\semis,\underline{\lambda}}$,
and the geometry of $\cX_{d,\red}$.

\section{\texorpdfstring{$(\varphi,\Gamma)$}{(phi,Gamma)}-modules
  with coefficients}\label{subsec: phi Gamma coefficients}
There is quite a lot of evidence, for example from Colmez's work on
the $p$-adic local Langlands correspondence~\cite{MR2642409}, and work
of Kedlaya--Liu
\cite{MR3379653}, 
that rather than considering families of representations of~$G_K$, it
is more natural to consider families of
\'etale~$(\varphi,\Gamma)$-modules.

The theory
of \'etale~$(\varphi,\Gamma)$-modules for~$\Zp$-representations was introduced
by Fontaine in~\cite{MR1106901}. There are various possible definitions that can
be made, with perfect, imperfect or overconvergent coefficient rings,
and different choices of~$\Gamma$; we discuss the various variants
that we use, and the relationships between them, at some length in the
body of the book.  For the purpose of this introduction we simply
let~$\A_K=W(k)((T))^\wedge$, where~$k$ is a finite extension of~$\Fp$
(depending on~$K$), and the hat denotes the $p$-adic completion. This
ring is endowed with a Frobenius~$\varphi$ and an action of a
profinite group~$\Gamma$ (an open subgroup of~$\Z_p^\times$) that
commutes with~$\varphi$; the formulae for $\varphi$ and for this action can be rather
complicated for general~$K$, although they admit a simple description
if~$K/\Qp$ is abelian.  (See Definition~\ref{defn: basic} and the
surrounding material.) 

An \'etale $(\varphi,\Gamma)$-module is then, by
definition, a finite $\A_K$-module endowed with commuting semilinear
actions of~$\varphi$ and~$\Gamma$, with the property that the
linearized action of~$\varphi$ is an isomorphism. There is a natural
equivalence of categories between the category of \'etale
$(\varphi,\Gamma)$-modules and the category of continuous
representations of~$G_K$ on finite $\Zp$-modules.

Let~$A$ be a $p$-adically complete $\Zp$-algebra. We
let~$\A_{K,A}:=A\cotimes_{\Zp}\A_K$ (where the completed tensor
product is taken with respect to the $p$-adic topology on~$A$
and the so-called weak topology on~$\A_K$),
and define an \'etale
$(\varphi,\Gamma)$-module with $A$-coefficients 
just 
as in the case~$A=\Zp$ described above, but now using the ring~$\A_{K,A}$.
 In the case that~$A$ is
finite as a $\Zp$-module there is again an equivalence of categories
with the category of continuous representations of~$G_K$ on finite
$A$-modules, but for more general~$A$ no such equivalence
exists. Our moduli stack~$\cX_d$ is defined to be the stack
over~$\Spf\Zp$ with the property that~$\cX_d(A)$ is the groupoid
of projective \'etale $(\varphi,\Gamma)$-modules of rank~$d$ with $A$-coefficients. (That
this is indeed an \'etale stack, indeed even an \emph{fpqc} stack,
follows from results of
Drinfeld.) Using the machinery of our 
paper~\cite{EGstacktheoreticimages} we are able to show that~$\cX_d$
is an Ind-algebraic stack, but to prove Theorem~\ref{thm:intro
  statement of basic properties of X} we need to go further and make a
detailed study of its special fibre and of the underlying reduced
substack. This study is guided by ideas coming from Galois
deformation theory and the weight part of Serre's conjecture,
in a manner that we now describe.

\section{Families of extensions}\label{subsec: families of extensions}
As we have already explained, over a general base~$A$ there is no
longer an equivalence between $(\varphi,\Gamma)$-modules and
representations of~$G_K$. Perhaps surprisingly, from the point of view
of applications of our stacks to the study of $p$-adic Galois
representations, this is a feature rather than a bug. For example, an
examination of the known results on the reductions modulo~$p$ of 2-dimensional
crystalline representations of~$G_{\Qp}$ (see for
example~\cite[Thm.\ 5.2.1]{MR2906353}) suggests that any
moduli space of mod~$p$ representations of~$G_K$ should have the
feature that the representations are generically reducible, but can
specialise to irreducible representations. A literal moduli space of
representations of a group cannot behave in this way (essentially
because Grassmannians are proper), but it turns out that the
underlying reduced substack~$\cX_{d,\red}$ of~$\cX_d$ does have this
property. (See also Section~\ref{sec:
  comparison to CWE} and Remark~\ref{rem: rank 1 versus rank d} for 
further discussions of the relationship between our stacks of
$(\varphi,\Gamma)$-modules and stacks of representations of a Galois
or Weil--Deligne group.) 

More precisely, the results of~\cite[Thm.\ 5.2.1]{MR2906353}, together
with the weight part of Serre's conjecture, suggest that each
irreducible component of~$\cX_{d,\red}$ should contain a dense set of
$\Fpbar$-points which are successive extensions of characters
$G_K\to\Fpbar$, with the extensions being as non-split as possible. This
turns out to be the case. The restrictions of these characters to the
inertia subgroup~$I_K$ are constant on the irreducible components, and
the discrete data of these characters, together with some further
information about peu- and tr\`es ramifi\'ee extensions, determines
the components. This discrete data can be conveniently and naturally
organised in terms of ``Serre weights'' $\underline{k}$, which are
tuples of integers which biject with the isomorphism classes of the
irreducible $\Fpbar$-representations of~$\GL_d(\cO_K)$. The
relationship between Serre weights and Galois representations is
important in the $p$-adic Langlands program, and in proving automorphy
lifting theorems, and we discuss it further in Section~\ref{subsec:
  geometric BM}. 

Having guessed that the $\Fpbar$-points of~$\cX_d$ should be arranged
in irreducible components in this way, an inductive strategy to prove
this suggests itself. It is easy to see that irreducible
representations of~$G_K$ are ``rigid'', in that there are up to twist
by unramified characters only finitely many in each dimension;
furthermore, it is at least intuitively clear that each such family of
unramified twists of a $d$-dimensional irreducible representation
should give rise to a zero-dimensional substack of~$\cX_d$ (there is a
$\Gm$ of twists, but also a $\Gm$ of automorphisms). On the other
hand, given
characters~$\chibar_1,\dots,\chibar_d:G_K\to\Fbar_p^\times$, a Galois
cohomology calculation suggests that there should be a substack
of~$\cX_{d,\red}$ 
of dimension~$[K:\Qp]d(d-1)/2$ given by successive
extensions of unramified twists of the~$\chibar_i$. Accordingly, one
could hope to construct the stacks 
corresponding to the Serre weights~$\underline{k}$ by inductively
constructing families of extensions of representations.

To confirm this expectation, we use the machinery originally developed
by Herr~\cite{MR1693457,MR1839766}, who gave an explicit complex
which is defined in terms of~$(\varphi,\Gamma)$-modules and computes Galois
cohomology. This definition goes over unchanged to the case with
coefficients, and with some effort we are able to adapt Herr's
arguments to our setting, and to prove finiteness and base change
properties (following Pottharst~\cite{MR3117501}, we in fact find it
helpful to think of the Herr complex of a $(\varphi,\Gamma)$-module
with $A$-coefficients as a perfect complex of $A$-modules). Using the Herr
complex, we can inductively construct irreducible closed
substacks~$\cX_{d,\red}^{\underline{k}}$ of~$\cX_{d,\red}$ of
dimension~$[K:\Qp]d(d-1)/2$ whose generic $\Fpbar$-points correspond
to successive extensions of characters as described above (the
restrictions of these characters to~$I_K$ being
determined by~$\underline{k}$). 
Furthermore, by a rather involved induction, we can show that the
union of the~$\cX_{d,\red}^{\underline{k}}$, together possibly with a closed
substack of~$\cX_{d,\red}$ of dimension strictly less than~$[K:\Qp]d(d-1)/2$,
exhausts~$\cX_{d,\red}$. In particular, each ~$\cX_{d,\red}^{\underline{k}}$ is an irreducible component
of~$\cX_{d,\red}$, and any 
irreducible component that is not one of the $\cX_{d,\red}^{\underline{k}}$
is of strictly
smaller dimension than these components.

One way to show that the~$\cX_{d,\red}^{\underline{k}}$ exhaust the
irreducible components of~$\cX_{d,\red}$ would be to show that every
representation $G_K\to\GL_d(\Fpbar)$ occurs as an $\Fbar_p$-valued point of
some~$\cX_{d,\red}^{\underline{k}}$.  We expect this to be 
difficult to show directly;
indeed, already for $d=2$ the
paper~\cite{CEGSKisinwithdd} shows that the closed points
of~$\cX_{d,\red}^{\underline{k}}$ are governed by the weight part of
Serre's conjecture, and the explicit description of this conjecture is
complicated (see e.g.\ \cite{BuzzardDiamondJarvis,MR3589336}). 
Furthermore it seems hard to explicitly understand the way
in which families of reducible $(\varphi,\Gamma)$-modules degenerate
to irreducible ones, or to reducible representations with different
restrictions to~$I_K$
(phenomena which are implied by the weight part of Serre's conjecture).

Instead, our approach is to show by a consideration of versal rings that~$\cX_{d,\red}$ is equidimensional of
dimension~$[K:\Qp]d(d-1)/2$; this suffices, since our inductive construction showed
that any other irreducible component would necessarily have dimension
strictly less than~$[K:\Qp]d(d-1)/2$. Our proof of this
equidimensionality relies on Theorems~\ref{thm:intro crystalline
  lifts} and~\ref{thm:intro statement on crystalline moduli}, as we
explain in Remark~\ref{rem: dimension and codimension} below.






\section{Crystalline lifts}\label{subsec: crystalline lifts}Theorem~\ref{thm:intro crystalline
  lifts} solves a problem that has been considered by various
authors, in particular~\cite{MullerThesis,2015arXiv150601050G}.
It admits a well-known inductive approach (which is taken in~\cite{MullerThesis,2015arXiv150601050G}): one writes $\rhobar$ as 
a successive extension of irreducible representations, lifts each of these
irreducible representations 
to a crystalline representation, and then attempts to lift the various
extension classes.  The difficulty that arises in this approach (which
has proved an obstacle to obtaining general statements along the lines
of Theorem~\ref{thm:intro crystalline lifts} until now) is showing that the 
mod $p$ extension classes that appear in this description of $\rhobar$ 
can actually be lifted to crystalline extension classes in characteristic
zero.  The basic source of the difficulty is that the local Galois $H^2$
can be non-zero, and non-zero classes in $H^2$ obstruct the lifting 
of extension classes (which can be interpreted as classes lying in
$H^1$). In fact, the difficulty is not so much in obtaining
crystalline extension classes, as in lifting to any classes in
characteristic zero; indeed, 
it was not previously known that an arbitrary~$\rhobar$ had \emph{any} lift to
characteristic zero at all. (Subsequently a different proof of the
existence of such a lift has been found by B\"ockle--Iyengar--Pa{\v{s}}k{\=u}nas~\cite{BIP}.)

Our proof of Theorem~\ref{thm:intro crystalline lifts} relies on the
inductive strategy described in the preceding paragraph, but we are able to
prove the following key result, which controls the obstructions that
can be presented by $H^2$, and is a consequence of
Theorems~\ref{thm:Xdred is algebraic} and~\ref{thm:reduced dimension}
(see also Remark~\ref{rem: dimension and codimension}).

\begin{prop}
	\label{prop:H2 dimension}
	The locus of points $\rhobar \in \cX_{d,\red}(\Fbar_p)$ 
	at which \[\dim H^2(G_K,\rhobar) \geq r\] is Zariski closed
	in $\cX_{d,\red}(\Fbar_p)$, and is of codimension $\geq r$.
      \end{prop}

Let~$R^{\square}_{\rhobar}$ denote the universal lifting ring
of~$\rhobar$, with universal lifting~$\rho^{\univ}$. For each regular
tuple of labeled Hodge--Tate weights $\underline{\lambda}$, we
let~$R_{\rhobar}^{\crys,\underline{\lambda}}$ denote the quotient
of~$R^{\square}_{\rhobar}$ corresponding to crystalline lifts
of~$\rhobar$ with Hodge--Tate weights~$\underline{\lambda}$ (of
course, this quotient is zero unless~$\rhobar$ admits such a
crystalline lift). Then $H^2(G_K, \rho^{\univ})$ is 
an $R^{\square}_{\rhobar}$-module, and Proposition~\ref{prop:H2 dimension}
implies the following corollary.

\begin{cor}
	\label{cor:H2 dimension}
	For any regular tuple of labeled Hodge--Tate weights $\underline{\lambda}$
	the locus of points 
	$x\in \Spec R^{\crys,\underline{\lambda}}_{\rhobar}/p$ for which
	$$\dim_{\kappa(x)}
 H^2(G_K,\rho^{\univ})	\otimes_{R^{\square}_{\rhobar}}\kappa(x)
	\geq r$$
	has codimension $\geq r$.
\end{cor}

\begin{remark}
	\label{rem:H2 dimension}
	Tate local duality, 
together with the compatibility of $H^2$ with base-change,
shows that
\begin{multline*}
	\dim_{\kappa(x)}
	 \bigl( H^2(G_K,\rho^{\univ})\otimes_{R^{\square}_{\rhobar}}\kappa(x)
\bigr)
	=
	\dim_{\kappa(x)}
	 H^2\bigr(G_K,\rho^{\univ}\otimes_{R^{\square}_{\rhobar}}\kappa(x)\bigl)
\\
	=
	\dim_{\kappa(x)}
	\Hom_{G_K}\bigl((\rho^{\univ})^{\vee}\otimes_{R^{\square}_{\rhobar}}\kappa(x),\epsilonbar\bigr)
\end{multline*}
	(where $\epsilonbar$ denotes the mod $p$ cyclotomic character, thought
	of as taking values in $\kappa(x)^{\times}$).
	Thus the statement of Corollary~\ref{cor:H2 dimension} is
	related to the way in which 
        $\Spec R_{\rhobar}^{\crys,\underline{\lambda}}/p$ intersects
	the reducibility locus in $\Spec R^{\square}_{\rhobar}$.
\end{remark}

Given Corollary~\ref{cor:H2 dimension},
we prove Theorem~\ref{thm:intro crystalline
	lifts}
by working purely within the context of formal lifting rings.  
However we don't know how to
prove the corollary while staying within that context.
Indeed, as Remark~\ref{rem:H2 dimension} indicates, 
this corollary is related to the way in which the special fibre
of a potentially crystalline deformation ring intersects another natural
locus in $\Spec R_{\rhobar}^{\square}$ (namely, the reducibility locus).
Since the special fibre of a potentially crystalline lifting ring
is not directly defined in deformation-theoretic terms, such questions
are notoriously difficult to study directly.
Our proof of the corollary proceeds
differently, by 
replacing a computation on the special fibre of the potentially
crystalline deformation ring by a computation on $\cX_{d,\red}$;
this latter space has a concrete description in terms of families of varying
~$\rhobar$, whose $H^2$ we are  able to compute, as a result of the
inductive construction of families of extensions described in Section~\ref{subsec: families of extensions}.

In order to deduce Corollary~\ref{cor:H2 dimension} from
Proposition~\ref{prop:H2 dimension}, it is crucial that we know that
the natural morphism $ \Spf R^{\crys,\underline{\lambda}}/p\to\cX_d$
is effective, in the sense that it arises from a morphism
$ \Spec R^{\crys,\underline{\lambda}}/p\to\cX_d$. More concretely, the
universal representation~$\rho^\univ$ gives an \'etale
$(\varphi,\Gamma)$-module over each Artinian quotient
of~$R_{\rhobar}^{\square}$. By passing to the limit over these
quotients, we obtain a
``universal formal \'etale $(\varphi,\Gamma)$-module''
over the completion of $(k\otimes_{\Zp}R_{\rhobar}^{\square}/p)((T))$
with respect to the maximal ideal~$\m$
of~$R_{\rhobar}^{\square}$. Since the special fibre of~$\cX_d$ is
formal algebraic but not algebraic (see Section~\ref{subsec: further
  questions} below), there is no corresponding $(\varphi,\Gamma)$-module
with $R_{\rhobar}^{\square}/p$-coefficients;
the  $\varphi$ and $\Gamma$ actions 
on the universal formal \'etale $(\varphi,\Gamma)$-module
 involve Laurent
tails of unbounded degree (with the coefficients of~$T^{-n}$ tending
to zero~$\m$-adically as $n\to\infty$).

The assertion that $ \Spf R^{\crys,\underline{\lambda}}/p\to\cX_d$ is
effective is equivalent to showing that the base change of the universal formal
\'etale $(\varphi,\Gamma)$-module to~$R^{\crys,\underline{\lambda}}/p$
arises from a genuine $(\varphi,\Gamma)$-module, i.e.\ from one that involves
only Laurent
tails of bounded degree. We deduce this from Theorem~\ref{thm:intro
  statement on crystalline moduli}. Indeed, the ring~$R^{\crys,\underline{\lambda}}/p$ is a versal ring for the special
fibre of the $p$-adic formal algebraic
stack~$\cX_{d}^{\crys,\lambdau}$, and (by the very definition of a
$p$-adic formal algebraic stack) this special fibre is an algebraic
stack; and the versal rings for algebraic stacks are always effective.

\begin{rem}
  \label{rem: dimension and codimension}As our citation of both
  Theorems~\ref{thm:Xdred is algebraic} and~\ref{thm:reduced
    dimension} for the proof of Proposition~\ref{prop:H2 dimension}
  may indicate, our proof of Proposition~\ref{prop:H2 dimension} is
  somewhat intricate. Indeed, in Theorem~\ref{thm:Xdred is algebraic},
  we show that~$\cX_{d,\red}$ has dimension at most~$[K:\Qp]d(d-1)/2$,
  and that the locus considered in Proposition~\ref{prop:H2 dimension}
  has dimension at most $[K:\Qp]d(d-1)/2-r$. This is in fact enough to
  deduce Corollary~\ref{cor:H2 dimension}, as ~$\Spec
  R^{\crys,\lambdau}/p$ is known to be equidimensional.

  Given Corollary~\ref{cor:H2 dimension}, we prove
  Theorem~\ref{thm:intro crystalline lifts}. In combination with the
  effective versality of the crystalline deformation rings discussed
  above we are then able to deduce the equidimensionality
  of~$\cX_{d,\red}$, and then also prove Proposition~\ref{prop:H2
    dimension} as stated.
\end{rem}

\section{Crystalline and semistable moduli stacks}\label{subsec: intro crystalline
  stacks}


We now explain the proof of Theorem~\ref{thm:intro statement on
  crystalline moduli}; the proof is essentially identical in the
crystalline and semistable cases, so we concentrate on the crystalline
case. To prove the theorem, it is necessary to have a criterion for a
$(\varphi,\Gamma)$-module to come from a crystalline Galois
representation. In the case that~$K/\Qp$ is unramified, it is possible
to give an explicit criterion in terms of Wach
modules~\cite{MR1415732}, but no such direct description is known for
general~$K$. Instead, following Kisin's construction of the
crystalline deformation rings~$R_{\rhobar}^{\crys,\lambda}$
in~\cite{MR2373358}, we use the theory of Breuil--Kisin modules. More
precisely, Kisin shows that crystalline representations of~$G_K$ have
finite height over the (non-Galois) Kummer extension~$K_\infty/K$
obtained by adjoining a compatible system of $p$-power roots of a
uniformizer of~$K$; here being of finite height means that the
corresponding \'etale $\varphi$-modules admit certain $\varphi$-stable
lattices, called Breuil--Kisin modules.

While not every representation of finite height over~$K_\infty$ comes
from a crystalline representation, we are able to show in
Appendix~\ref{app: BKF pst} (jointly written by T.G.\ and Tong Liu)
that a representation $G_K\to\GL_d(\Zpbar)$ is crystalline if and only
if it is of finite height for \emph{every} choice of~$K_\infty$, and
if the corresponding Breuil--Kisin modules satisfy certain natural
compatibilities. (These compatibilities are best expressed in terms of
Bhatt--Scholze's prismatic site, as in~\cite{bhatt2021prismatic}, but
we do not make use of that perspective in this book. Instead, we write
down explicit conditions on the corresponding Breuil--Kisin--Fargues
modules; recall that Breuil--Kisin--Fargues modules are a variant of
Breuil--Kisin modules introduced by Fargues, see e.g.\ \cite[\S 4]{2016arXiv160203148B}.)

We use this description of the crystalline representations to
prove the existence of the stacks~$\cX_d^{\crys,\lambdau}$. 
The proof
that~$\cX_d^{\crys,\lambdau}$ is a $p$-adic formal algebraic stack
relies on an
analogue of results of Caruso--Liu~\cite{MR2745530} on extensions of
the Galois action on Breuil--Kisin modules, which roughly speaking
says that the action of~$G_{K_\infty}$ determines the action of~$G_K$
up to a finite amount of ambiguity. More precisely, given a
Breuil--Kisin module over a $\Z/p^a$-algebra for some~$a\ge 1$, there
is a finite subextension $K_s/K$ of~$K_\infty/K$ depending only
on~$a$, $K$ and the height of the Breuil--Kisin module, such that there is a canonical action
of~$G_{K_s}$ on the corresponding Breuil--Kisin--Fargues
module. This canonical action is constructed by Frobenius
amplification, and in the case that the Breuil--Kisin module arises from the
reduction modulo~$p^a$ of a crystalline representation of~$G_K$, the canonical
action coincides with the restriction to~$G_{K_s}$ of the $G_K$-action
on the Breuil--Kisin--Fargues module. 
(In \cite{MR2745530} a version
of this canonical action is used to prove ramification bounds on the
reductions modulo~$p^a$ of crystalline representations; in
Chapter~\ref{sec: the rank one case}, we use analogous arguments in
the setting of~$(\varphi,\Gamma)$-modules to relate our stacks to
stacks of Weil group representations in the rank one case.)

There is one significant technical difficulty, which is that we need
to define morphisms of stacks that correspond to the restriction of
Galois representations from~$K$ to~$K_\infty$. In order to do this we
have to compare $(\varphi,\Gamma)$-modules with $A$-coefficients
(which are defined via the cyclotomic extension
$K(\zeta_{p^\infty})/K$) to $\varphi$-modules with $A$-coefficients
defined via the extension~$K_\infty/K$. We do not know of a direct way to do this; we
proceed by proving a correspondence between $\varphi$-modules over
Laurent series rings with $\varphi$-modules over the perfections of
these Laurent series rings,
and proving the following descent result which may be of independent
interest; in the statement,
$\C$ denotes the completion of an algebraic closure of~$\Qp$.

\begin{theorem}[Theorem~\ref{thm:descending projective
    modules}]\label{thm:descending projective modules intro version}
	  Let $A$ be a finite type $\Z/p^a$-algebra, for some $a \geq 1$. Let
  $F$ be a closed perfectoid subfield of $\C$, with tilt~$F^\flat$, a
  closed perfectoid subfield of $\C^\flat$.  Write
  $W(F^\flat)_A := W(F^\flat)\otimes_{\Z_p}A .$ 

Then the inclusion
  $W(F^\flat)_A \to W(\C^\flat)_A$ is a faithfully flat morphism of
  Noetherian rings, and the functor
  $M \mapsto W(\C^\flat)_A\otimes_{W(F^\flat)_A} M$ induces an
  equivalence between the category of finitely generated projective
  $W(F^\flat)_A$-modules and the category of finitely generated
  projective $W(\C^\flat)_A$-modules endowed with a continuous
  semi-linear $G_F$-action.
\end{theorem} The existence of the required morphism of stacks follows
easily from two applications of Theorem~\ref{thm:descending projective modules intro
  version}, applied with~$F$ equal to respectively the completion of~$K_\infty$ and
the completion of~$K(\zeta_{p^\infty})$. Furthermore, this
construction gives an alternative description of our stacks, as moduli
spaces of $W(\C^\flat)_A$-modules endowed with commuting
  semi-linear actions of~$G_{K}$ and~$\varphi$.  It seems plausible
  that this description will be useful in future work, as it connects
  naturally to the theory of Breuil--Kisin--Fargues modules (and
  indeed we use this connection in our construction of the potentially
  semistable moduli stacks). Note though that the description in terms
  of~$(\varphi,\Gamma)$-modules is important (at least in our
  approach) for establishing the basic finiteness properties of our stacks.

\section{The geometric Breuil--M\'ezard conjecture and the weight
  part of Serre's conjecture}\label{subsec:
  geometric BM}\sectionmark{The geometric Breuil--M\'ezard conjecture}
We will now briefly explain our results and conjectures relating our
stacks to the Breuil--M\'ezard conjecture and the weight part of
Serre's conjecture. Further explanation and motivation can be found
throughout Chapter~\ref{sec: BM}. Some of these results were previewed
in~\cite[\S 6]{2015arXiv150902527G}, and the earlier sections of that
paper (in particular the introduction) provide an overview of the
weight part of Serre's conjecture and its connections to the
Breuil--M\'ezard conjecture that may be helpful to the reader who is
not already familiar with them. As in the rest of this introduction,
we ignore the possibility of inertial types, and we also restrict to
crystalline representations for the purpose of exposition. Everything
in this section extends to the more general setting of potentially
semistable representations, and indeed as we explain in
Section~\ref{subsec: CEGS} when discussing the
papers~\cite{CEGSKisinwithdd} and~\cite{geekisin}, the additional
information provided by non-trivial inertial types is very important.

Let $\rhobar: G_K \to \GL_d(\F)$ be a continuous representation (for some finite
extension $\F$ of $\F_p$), and let $R_{\rhobar}^\square$ be the
corresponding universal lifting ring. The corresponding formal scheme $\Spf R_{\rhobar}^\square$ doesn't carry
a lot of evident geometry in and of itself; for example, its underlying reduced
subscheme is simply the closed point $\Spec \F$, corresponding 
to $\rhobar$ itself. 
On the other hand, $\cX_d$ has a quite non-trivial underlying reduced
substack~$\cX_{d,\red}$, which parameterizes all the $d$-dimensional
residual representations of $G_K$.  It is natural to ask whether this
underlying reduced substack has any significance in formal deformation
theory.  More precisely, we could ask for the meaning of the fibre
product $ \Spf R_{\rhobar}^\square \times_{\cX_d} \cX_{d,\red}.$

This fibre product is a reduced closed formal subscheme of $\Spf R_{\rhobar}^\square$
of dimension $d^2+[K:\Q_p]d(d-1)/2$.   It arises (via completion at
the closed point) from a closed subscheme of $\Spec R_{\rhobar}^\square$ (as 
does any closed formal subscheme of the $\Spf$ of a complete Noetherian local
ring), whose irreducible components,
when thought of as cycles on $\Spec R_{\rhobar}^\square$, are precisely the cycles
that (conjecturally) appear in the geometric Breuil--M\'ezard
conjecture of~\cite{emertongeerefinedBM}. 
More precisely, we obtain the following qualitative version of the
geometric Breuil--M\'ezard
conjecture~\cite[Conj.~4.2.1]{emertongeerefinedBM}.

\begin{theorem}[Theorem~\ref{thm: qualitative BM}]\label{thm: intro geom BM without multiplicities}
	If $\rhobar: G_K \to \GL_d(\F)$ is a continuous
        representation, 
	then there are finitely many cycles of dimension
        $d^2+[K:\Q_p]d(d-1)/2$ in $\Spec R^{\square}_{\rhobar}/p$,
        such that for any regular tuple of labeled Hodge--Tate
        weights $\underline{\lambda}$, 
        the special fibre
        $\Spec R^{\crys,\underline{\lambda}}_{\rhobar}/p$ is
        set-theoretically supported on the union of some number of
        these
        cycles.
\end{theorem}

The cycles in the statement of the theorem are precisely
those arising from the fibre products $\Spf R_{\rhobar}^\square\times_{\cX_d}
\cX_{d,\red}^{\underline{k}}$, where~$\underline{k}$ runs over the
Serre weights. While Theorem~\ref{thm: intro geom BM without
  multiplicities} is a purely local statement, we do not know how to
prove it without using the stacks~$\cX_d$.

The full geometric Breuil--M\'ezard conjecture
of~\cite{emertongeerefinedBM} makes precise predictions about the
multiplicities of the cycles of the special fibres of
$\Spec R^{\crys,\underline{\lambda}}_{\rhobar}/p$; passing from cycles
to Hilbert--Samuel multiplicities then recovers the original
Breuil--M\'ezard conjecture~\cite{BreuilMezard} (or rather a natural
generalisation of it to~$\GL_d$), which we refer to as the ``numerical
Breuil--M\'ezard conjecture''. In particular, the multiplicities are
expected to be computed in terms of
quantities~$n_{\underline{k}}^\crys(\lambda)$ that are defined as
follows: one associates an irreducible algebraic
representation~$\sigma^{\crys}(\lambda)$ of~$\GL_d/K$ to~$\lambdau$,
defined to have highest weight (a certain shift of)~$\lambdau$. The
semisimplification of the reduction mod~$p$
of~$\sigma^{\crys}(\lambda)$ can be written as a direct sum of
irreducible representations of~$\GL_d(k)$,
and~$n_{\underline{k}}^\crys(\lambda)$ is defined to be the
multiplicity with which the Serre weight~$\underline{k}$ occurs.

In Chapter~\ref{sec: BM} we explain
that as we run over all~$\rhobar$, the geometric Breuil--M\'ezard
conjecture is equivalent to the following analogous conjecture for the
special fibres of our crystalline and semistable stacks. Here by a
``cycle'' in~$\cX_{d,\red}$ we mean a formal $\Z$-linear combination
of its irreducible components~$\cX_d^{\underline{k}}$.

\begin{conj}[Conjecture~\ref{conj: geometric BM}]
  \label{conj: intro version of geometric BM}There are
  cycles~$Z_{\underline{k}}$ in~$\cX_{d,\red}$ with the property that
  for each regular tuple of labeled Hodge--Tate weights~$\lambdau$,
  the underlying cycle of the special fibre
  of~$\cX_d^{\crys,\lambdau}$ is
  $\sum_{\underline{k}}n_{\underline{k}}^\crys(\lambda)\cdot
  Z_{\underline{k}}$.
\end{conj}
In fact we expect that the cycles~$Z_{\underline{k}}$ are effective,
i.e.\ that they are a linear combination of the irreducible
components~$\cX_d^{\underline{k}'}$ with non-negative integer
coefficients. Since there are infinitely many possible~$\lambdau$, the
cycles~$Z_{\underline{k}}$, if they exist, are hugely overdetermined
by Conjecture~\ref{conj: intro version of geometric BM}.

As first explained in~\cite{KisinFM}, the (numerical) Breuil--M\'ezard
conjecture has important consequences for automorphy lifting theorems;
indeed proving the conjecture is closely related to
proving automorphy lifting theorems in situations with arbitrarily
high weight or ramification at the places dividing~$p$. Conversely,
following~\cite{MR2827797}, one can use automorphy lifting theorems to
deduce cases of the Breuil--M\'ezard conjecture. Automorphy lifting
theorems involve a fixed~$\rhobar$, and in fact we can deduce
Conjecture~\ref{conj: intro version of geometric BM} from the
Breuil--M\'ezard conjecture for a finite set of suitably
generic~$\rhobar$.

In particular, we are able to combine results in the literature to
show that for~$\GL_2$ the cycles~$Z_{\underline{k}}$ in
Conjecture~\ref{conj: intro version of geometric BM} must have a
particularly simple form: we necessarily
have~$Z_{\underline{k}}=\cX_{d,\red}^{\underline{k}}$
unless~$\underline{k}$ is a so-called ``Steinberg'' weight, in which
case~ $Z_{\underline{k}}$ is the sum of~$\cX_{d,\red}^{\underline{k}}$
and one other irreducible component. (More precisely, what we show,
following~\cite{CEGSKisinwithdd,geekisin}, is that with these cycles~
$Z_{\underline{k}}$, Conjecture~\ref{conj: intro version of geometric
  BM} holds for all ``potentially Barsotti--Tate'' representations.)

  The weight part of Serre's conjecture predicts the weights in which
  particular Galois representations contribute to the mod~$p$
  cohomology of locally symmetric spaces. Following~\cite{geekisin},
  this conjecture is closely related to the Breuil--M\'ezard
  conjecture; indeed, if Conjecture~\ref{conj: intro version of
    geometric BM} holds, then the set of Serre weights associated to a
  representation $\rhobar:G_K\to\GL_d(\Fpbar)$ should be precisely the
  weights~$\underline{k}$ for which~$Z_{\underline{k}}$ is supported
  at~$\rhobar$. In other words, if we refine our labelling of the
  irreducible components of~$\cX_{d,\red}$ by labelling each component
  by the union of the weights~$\underline{k}$ for which that component
  contributes to~$Z_{\underline{k}}$, then we expect the set of Serre
  weights for~$\rhobar$ to be the union of the sets of weights
  labelling the irreducible components containing~$\rhobar$. This
  expectation holds for~$\GL_2$ by the main results
  of~\cite{CEGSKisinwithdd}.


  \section{Further questions}\label{subsec: further questions}
  There are many other questions one could ask about the
  stacks~$\cX_d$, which we hope to return to in future
  papers. For example, we show in Proposition~\ref{prop: X is not
    padic formal algebraic} that~$\cX_d$ is not a $p$-adic formal
  algebraic stack. Indeed, if it were $p$-adic formal algebraic, then
  its special fibre would be an algebraic stack, whose dimension would
  be equal to the dimension of its underlying reduced substack. In
  turn, this would imply that the versal rings~$R_{\rhobar}^\square$
  would have dimension equal to the dimensions of the crystalline
  deformation rings~$R_{\rhobar}^{\crys,\lambdau}$, and this is known
  not to be true. In fact, it is a folklore conjecture, recently
  proved by B\"ockle--Iyengar--Pa{\v{s}}k{\=u}nas~\cite{BIP},  that the
  lifting rings~$R_{\rhobar}^\square$ are $\Zp$-flat local complete intersections
  of dimension~$1+d^2+[K:\Qp]d^2$, which should imply that the
  stacks~$\cX_d$ are $\Zp$-flat local complete intersections of
  dimension~$1+[K:\Qp]d^2$ (a notion that we do not attempt to make
  precise for formal algebraic stacks).

It is natural to ask about the rigid analytic generic fibre of $\cX_d$; this 
should exist as a rigid analytic stack in an appropriate sense.  The
 generic fibres of the substacks $\cX_d^{\underline{k}}$ should admit 
morphisms to the stacks of Hartl and Hellmann~\cite{2013arXiv1312.6371H}
(although these morphisms won't be isomorphisms, since for any finite 
extension $E$ of~$\Q_p$, the $\cO_E$-points of 
$\cX_d^{\underline{k}}$,  which would coincide with the $E$-points
of its generic fibre, correspond to lattices in crystalline representations,
whereas the stacks of~\cite{2013arXiv1312.6371H} parameterize crystalline
or semistable representations themselves). 

We expect that the~$\cX_d$ will have a role to play in generalisations
of the $p$-adic local Langlands correspondence. For example, we expect
that when~$K=\Qp$ the $p$-adic local Langlands correspondence
for~$\GL_2(\Qp)$ can be extended to give rise to sheaves
of~$\GL_2(\Qp)$-representations on~$\cX_2$. More generally, we expect
 that there will be a $p$-adic analogue of the work of Fargues--Scholze 
on the local Langlands correspondence \cite{fargues2021geometrization} 
involving the stacks~$\cX_d$.


\section{Previous work}
The description of local Galois representations in terms of \'etale $(\varphi,\Gamma)$-modules
is due to Fontaine~\cite{MR1106901}.  The importance of ``height''
as an aspect of the theory was already emphasized in~\cite{MR1106901}, 
and was further developed by Wach~\cite{MR1415732}, who explored 
the relationship between the finite height condition and crystallinity 
of Galois representations in the absolutely unramified context.

The use of what are now called Breuil--Kisin modules as a tool
for studying crystalline and semistable representations 
for general (i.e.\ not necessarily absolutely unramified) $p$-adic fields
(a study which, apart from its intrinsic importance, 
is crucial for treating potentially crystalline or
semistable representations, even in the absolutely unramified context)
was due originally to Breuil~\cite{Breuilunpublished}  
and was extensively developed by Kisin~\cite{KisinModularity,
  MR2373358}, who used them to study Galois deformation rings.

The algebro-geometric and  moduli-theoretic perspectives that already played
key roles in Kisin's work were further developed by
Pappas and Rapoport \cite{MR2562795}, who introduced moduli
stacks of Breuil--Kisin modules and of \'etale $\varphi$-modules;
it this work of Pappas and Rapoport, which can be very roughly thought
of as constructing moduli stacks of representations of the absolute Galois
groups of certain perfectoid fields, which is the immediate launching
point for our work in this book, as well as for our 
paper~\cite{EGstacktheoreticimages}.
(We should also mention Drinfeld's work~\cite{MR2181808}, which underpins the
verification of the stack property for the constructions of~\cite{MR2562795},
as well as for those of the present book.) Our use of moduli stacks
of Breuil--Kisin--Fargues modules (in the construction of the
potentially crystalline and semistable substacks) was in part inspired
by the work of Fargues and Bhatt--Morrow--Scholze (see in particular
\cite[\S 4]{2016arXiv160203148B}), which taught us not to be afraid of~$\Ainf$.

Moduli stacks parameterizing crystalline and semistable representations 
have already been constructed by Hartl and Hellmann~\cite{2013arXiv1312.6371H};
as remarked upon above, these stacks should have a relationship
to the stacks $\cX_d^{\underline{k}}$ that we construct. See also the
related papers of Hellmann~\cite{MR3566522,MR3103130}.

As far as we are aware, the first construction of moduli stacks of
representations of~$G_K$ in which the residual
representation~$\rhobar$ can vary is the work of Carl
Wang-Erickson~\cite{MR3831282} mentioned above, which constructs and studies such
stacks in the case that~$\rhobar$ has fixed semisimplification. These
are literally
moduli stacks of representations of~$G_K$; 
they are isomorphic to certain 
substacks of our stacks
$\cX_d$, as we explain in Section~\ref{sec:
  comparison to CWE}.



\section{An outline of the book}\label{subsec: outline}We finish this
introduction with a brief overview of the contents of this
book. The reader may also wish to refer to the introductions to each
chapter, as well as to the overview of this book provided by the notes~\cite{emerton2020moduli}.  

In Chapter~\ref{section: coefficient rings} we recall various of the
coefficient rings used in the theories of $(\varphi,\Gamma)$-modules
and Breuil--Kisin modules, and introduce versions of these rings with
coefficients in a $p$-adically complete $\Zp$-algebra. We also prove
almost Galois descent results for projective modules, and deduce
Theorem~\ref{thm:descending projective modules intro version}.

In Chapter~\ref{sec: phi modules and phi gamma modules} we recall the
results of~\cite{EGstacktheoreticimages} on moduli stacks of
$\varphi$-modules, and use them to define our stacks~$\cX_d$ of
\'etale $(\varphi,\Gamma)$-modules. With some effort, we prove
that~$\cX_d$ is an Ind-algebraic stack. Chapter~\ref{sec: crystalline
  and semistable} defines various moduli stacks of Breuil--Kisin and
Breuil--Kisin--Fargues modules, and uses them 
to construct our moduli
stacks of potentially semistable and potentially crystalline
representations, and in particular to prove Theorem~\ref{thm:intro
  statement on crystalline moduli}.

Chapter~\ref{sec: families of extensions} develops the theory of the
Herr complex, proving in particular that it is a perfect complex and
is compatible with base change. We show how to use the Herr complex to
construct families of extensions of $(\varphi,\Gamma)$-modules, and we
use these families to define the irreducible
substack~$\cX_{d,\red}^{\underline{k}}$ corresponding to a Serre
weight~$\underline{k}$. By induction on~$d$ we prove that~$\cX_d$ is a
Noetherian formal algebraic stack, and establish a version of
Proposition~\ref{prop:H2 dimension} (although as discussed in
Remark~\ref{rem: dimension and codimension}, we do not prove
Proposition~\ref{prop:H2 dimension} as stated at this point in the
argument). 

It may help the reader for us to point out
 that Chapters~\ref{sec: crystalline and semistable}
and~\ref{sec: families of extensions} are essentially independent of one another,
and are of rather different flavours.  
Chapter~\ref{sec: crystalline and semistable} involves an interleaving of stack-theoretic
arguments with ideas from $p$-adic Hodge theory and the theory of Breuil--Kisin
modules, while in 
Chapter~\ref{sec: families of extensions},
once we complete our analysis of the Herr complex,
our perspective begins to shift: although at a technical level we of course continue
to work with $(\varphi,\Gamma)$-modules, we begin to think in terms 
of Galois representations and Galois cohomology, and the more foundational
 arguments of the preceding chapters recede somewhat into the background.
 
In Chapter~\ref{sec:properties} we combine the results of
Chapters~\ref{sec: crystalline and semistable} and ~\ref{sec: families
  of extensions} with a geometric argument on the local
deformation ring to prove Theorem~\ref{thm:intro crystalline
  lifts}. Having done this, we are then able to improve on the results
on~$\cX_d$ established in the earlier chapters by proving
Theorem~\ref{thm:intro statement of basic properties of X}. We also
deduce Theorem~\ref{thm: intro existence of global lifts}, as well as
determining the closed points of~$\cX_d$, and describing the
relationship of our stacks with Wang--Erickson's stacks of Galois representations.

Chapter~\ref{sec: the rank one case} gives explicit descriptions of
various of our moduli stacks in the case~$d~=~1$, relating them to
moduli stacks of Weil group representations. Chapter~\ref{sec: BM}
explains our geometric version of the Breuil--M\'ezard conjecture, and
proves some results towards it, particularly in the case~$d=2$.

Finally the appendices for the most part establish various technical
results used in the body of the book. We highlight in particular
Appendix~\ref{app: formal algebraic stacks}, which summarises the
theory of formal algebraic stacks developed
in~\cite{Emertonformalstacks}, and Appendix~\ref{app: BKF pst}, which
combines the theory of Breuil--Kisin--Fargues modules with Tong Liu's
theory of $(\varphi,\Ghat)$-modules to give a new characterisation of
integral lattices in potentially semistable representations, of which we
make crucial use 
in Chapter~\ref{sec: crystalline and semistable}.

\section{Acknowledgements}\label{subsec: acknowledgements}We would like to thank Robin Bartlett,
Laurent Berger, Bhargav Bhatt, Xavier Caruso, Pierre Colmez, Yiwen
Ding, Florian Herzig, Ashwin
Iyengar, Mark Kisin, Tong Liu, George Pappas, Dat Pham, L\'eo Poyeton, Michael Rapoport, and
Peter Scholze for helpful correspondence and conversations.
We would particularly like to
thank Tong Liu for his contributions to Appendix~\ref{app: BKF
  pst}. We would also like to thank the organisers (Johannes Ansch\"utz,
Arthur-C\'esar Le Bras, and Andreas Mihatsch) 
of the 2019 Hausdorff School on ``the Emerton--Gee stack and related
topics'', as well as all the participants,
both for the encouragement to finish this book and for the
many helpful questions and corrections resulting from our lectures. 

Our mathematical debt to the late Jean-Marc Fontaine will be obvious
to the reader. This book benefited from several conversations with
him over the years~2011-2018, 
and from the interest he showed in our
results; in particular, his explanations to us of the relationship
between framing Galois representations and Fontaine--Laffaille modules
during his visit to Northwestern University in the spring of 2011 provided
an important clue as to the correct definitions of our stacks.

We owe special thanks to Colette and Therese for all of their patience and
support during the writing and revision of this book.


The first author was supported in part by the
  NSF grants DMS-1303450, DMS-1601871, and DMS-1902307. The second author was 
  supported in part by a Leverhulme Prize, EPSRC grant EP/L025485/1, Marie Curie Career
  Integration Grant 303605,
  ERC Starting Grant 306326, ERC Advanced Grant 884596, and a Royal Society Wolfson Research
  Merit Award. This project has received funding from the European Research Council (ERC) under the European Union’s Horizon 2020 research and innovation programme (grant agreement No. 884596).

\section{Notation and conventions}\label{subsec: notation and conventions}
\subsection*{$p$-adic Hodge theory}Let $K/\Qp$ be a finite
extension. If $\rho$ is a de Rham representation of $G_K$ on a
$\Qpbar$-vector space~$W$, then we will write $\WD(\rho)$ for the
corresponding Weil--Deligne representation of $W_K$ (see e.g.\
\cite[App. B]{CDT}), and if $\sigma:K \into \Qpbar$ is a continuous
embedding of fields then we will write $\HT_\sigma(\rho)$ for the
multiset of Hodge--Tate numbers of $\rho$ with respect to $\sigma$,
which by definition contains $i$ with multiplicity
$\dim_{\Qpbar} (W \otimes_{\sigma,K} \widehat{\barK}(i))^{G_K} $. Thus
for example if~$\epsilon$ denotes the $p$-adic cyclotomic character,
then $\HT_\sigma(\epsilon)=\{ -1\}$.

By a \emph{$d$-tuple of labeled Hodge--Tate
  weights~$\underline{\lambda}$}, \index{labeled Hodge--Tate
  weights} \index{$\underline{\lambda}$}
we mean a tuple of integers~$\{\lambda_{\sigma,i}\}_{\sigma:K\into\Qpbar,1\le i\le d}$
with $\lambda_{\sigma,i}\ge \lambda_{\sigma,i+1}$ for all~$\sigma$ and all $1\le
i\le d-1$. We will also refer to~$\underline{\lambda}$ as a
\index{Hodge type} \index{inertial type}
\emph{Hodge type}. By an \emph{inertial type~$\tau$} we mean a representation
$\tau:I_K\to\GL_d(\Qpbar)$ which extends to a  representation
of~$W_K$ with open kernel (so in particular, $\tau$ has finite
image). 

Then we say that~$\rho$ has Hodge type~$\underline{\lambda}$ (or
labeled Hodge--Tate weights~$\underline{\lambda}$) if for
each~$\sigma:K\into\Qpbar$ we
have~$\HT_\sigma(\rho)=\{\lambda_{\sigma,i}\}_{1\le i\le
  d}$, and we say that~$\rho$ has inertial
type~$\tau$ if~$\WD(\rho)|_{I_K}\cong\tau$.

We often somewhat abusively write that a representation
$\rho:G_K\to\GL_d(\Zp)$ is crystalline (or potentially crystalline, or
semistable, or\dots) if the corresponding
representation $\rho:G_K\to\GL_d(\Qp)$ is crystalline (or potentially
crystalline, or\dots).

\subsection*{Serre weights and Hodge--Tate weights}By a \emph{Serre
weight} $\underline{k}$ we mean a tuple of \index{Serre weight} \index{\underline{k}}
integers~$\{k_{\sigmabar,i}\}_{\sigmabar:k\into\Fpbar,1\le i\le d}$
with the properties that \begin{itemize}
\item $p-1\ge k_{\sigmabar,i}-k_{\sigmabar,i+1}\ge 0$ for each $1\le i\le
  d-1$, and
\item $p-1\ge k_{\sigmabar,d}\ge 0$, and not every~$k_{\sigmabar,d}$
  is equal to~$p-1$.
\end{itemize}
The set of Serre weights is in bijection with the set of irreducible
$\Fpbar$-representations of $\GL_d(k)$, via passage to highest
weight vectors (see for example the appendix to~\cite{MR2541127}). 


Each embedding~$\sigma:K\into\Qpbar$ induces an
embedding~$\sigmabar:k\into\Fpbar$; if $K/\Qp$ is ramified, then
each~$\sigmabar$ corresponds to multiple embeddings~$\sigma$.  We say
that~$\underline{\lambda}$ is a lift of  ~$\underline{k}$ 
if for each
embedding~$\sigmabar:k\into\Fpbar$, we can choose an
embedding~$\sigma:K\into\Qpbar$ lifting~$\sigmabar$, with the
properties that:
\begin{itemize}
\item $\lambda_{\sigma,i}=k_{\sigma,i}+d-i$, and
\item if $\sigma':K\into\Qpbar$ is any other lift of~$\sigmabar$,
  then $k_{\sigma',i}=d-i$.
\end{itemize}

\subsection*{Lifting rings}Let $K/\Qp$ be a finite extension, and
let $\rhobar:G_K\to\GL_d(\Fpbar)$ be a continuous representation. Then
the image of~$\rhobar$ is contained in~$\GL_d(\F)$ for any
sufficiently large finite
extension~$\F/\F_p$. Let~$\cO$ be the ring of integers in some
finite extension~$E/\Qp$, and suppose that the residue field of~$E$
is~$\F$. 
Let~$R^{\square,\cO}_{\rhobar}$ be the universal \index{$R^{\square,\cO}_{\rhobar}$}
lifting~$\cO$-algebra of~$\rhobar$; by definition, this (pro-)represents the
functor given by lifts of~$\rhobar$ to representations $\rho: G_K \to
\GL_d(A)$, for~$A$ an Artin local $\cO$-algebra with residue field~$\F$. 
The precise choice of~$E$ is unimportant, in the sense that
if~$\cO'$ is the ring of integers in a finite extension~$E'/E$, then
by~\cite[Lem.\ 1.2.1]{BLGGT} we have
$R^{\square,\cO'}_{\rhobar}=R^{\square,\cO}_{\rhobar}\otimes_{\cO}\cO'$.


Fix some Hodge type~$\underline{\lambda}$ and inertial
type~$\tau$. If~$\cO$ is chosen large enough that the inertial
type~$\tau$ is defined over~$E=\cO[1/p]$, and large enough that~$E$
contains the images of all embeddings~$\sigma:K\into\Qpbar$, then we
have the usual lifting $\cO$-algebras
~$R^{\crys,\underline{\lambda},\tau,\cO}_{\rhobar}$ \index{$R^{\crys,\underline{\lambda},\tau,\cO}_{\rhobar}$}
and~$R^{\semis,\underline{\lambda},\tau,\cO}_{\rhobar}$. By \index{$R^{\semis,\underline{\lambda},\tau,\cO}_{\rhobar}$}
definition, these are the unique $\cO$-flat quotients
of~$R_{\rhobar}^{\square,\cO}$ with the property that if~$B$ is a
finite flat~$E$-algebra, then an $\cO$-algebra homomorphism
$R_{\rhobar}^{\square,\cO}\to B$ factors through ~$R^{\crys,\underline{\lambda},\tau,\cO}_{\rhobar}$
(resp.\ through~$R^{\semis,\underline{\lambda},\tau,\cO}_{\rhobar}$)
if and only if the corresponding representation of~$G_K$ is
potentially crystalline (resp.\ potentially semistable) of Hodge
type~$\underline{\lambda}$ and inertial type~$\tau$. If~$\tau$ is
trivial, we will sometimes omit it from the notation.
By
the main theorems of~\cite{MR2373358}, 
these rings are (when they are nonzero) equidimensional of
dimension \[1+d^2+\sum_{\sigma}\#\{1\le i<j\le
  d|\lambda_{\sigma,i}>\lambda_{\sigma,j}\}.\]
Note that this quantity is at most $1+d^2+[K:\Qp]d(d-1)/2$, with
equality if and only if~$\underline{\lambda}$ is \emph{regular}, in
the sense that $\lambda_{\sigma,i}> \lambda_{\sigma,i+1}$ for all~$\sigma$ and all $1\le
i\le d-1$. As above, we have
$R^{\crys,\underline{\lambda},\tau,\cO'}_{\rhobar}=R^{\crys,\underline{\lambda},\tau,\cO}_{\rhobar}\otimes_{\cO}\cO'$,
and similarly for~$R^{\semis,\underline{\lambda},\tau,\cO'}_{\rhobar}$. 
By \cite[Thm.\ 3.3.8]{MR2373358} 
the localized rings 
$R^{\crys,\underline{\lambda},\tau,\cO}_{\rhobar}[1/p]$
are regular, and thus the rings
$R^{\crys,\underline{\lambda},\tau,\cO}_{\rhobar}$
(which embed into their localizations away from~$p$,
since they are $\cO$-flat) are reduced.

\subsection*{Algebra}Our conventions typically
follow~\cite{stacks-project}. In particular, if~$M$ is an abelian
topological group with a linear topology, then as
in~\cite[\href{https://stacks.math.columbia.edu/tag/07E7}{Tag
  07E7}]{stacks-project} we say that~$M$ is {\em complete} if the
natural morphism $M\to \varprojlim_i M/U_i$ is an isomorphism,
where~$\{U_i\}_{i \in I}$ is some (equivalently any) fundamental
system of neighbourhoods of~$0$ consisting of subgroups. Note that in
some other references this would be referred to as being~{\em complete
  and separated}.

If~$R$ is a ring, we write 
$D(R)$ for the (unbounded) derived category of $R$-modules. 
We say that a complex~$P^\bullet$ is \emph{good} if it is a bounded
complex of finite projective $R$-modules; then an object $C^\bullet$ of $D(R)$ is called a \emph{perfect
  complex} if there is a quasi-isomorphism
$P^\bullet \rightarrow C^\bullet$ where $P^\bullet$ is good. In fact, $C^\bullet$ is
perfect if and only if it is isomorphic in $D(R)$ to a good complex
$P^\bullet$: if we have another complex
$D^\bullet$ and quasi-isomorphisms $P^\bullet \rightarrow D^\bullet$,
$C^\bullet \rightarrow D^\bullet$, then there is a quasi-isomorphism
$P^\bullet \rightarrow C^\bullet$
(\cite[\href{http://stacks.math.columbia.edu/tag/064E}{Tag
  064E}]{stacks-project}). \index{good complex} \index{perfect complex}

\subsection*{Stacks}Our conventions on algebraic stacks and formal
algebraic stacks are those of~\cite{stacks-project}
and~\cite{Emertonformalstacks}. We recall some terminology and results
in Appendix~\ref{app: formal algebraic
  stacks}. Throughout the book, if~$A$ is a topological ring
and~$\cC$ is a stack we write~$\cC(A)$ for~$\cC(\Spf A)$; if~$A$ has
the discrete topology, this is equal to~$\cC(\Spec A)$. 
\chapter{Rings and coefficients}\label{section: coefficient rings}In this
chapter we study various rings which will be the coefficients of the
$\varphi$-modules and $(\varphi,\Gamma)$-modules
that we study in the subsequent chapters.
Throughout the chapter,
we fix a finite extension $K$ of~$\Q_p$,
which we regard as a subfield of some fixed algebraic closure $\Qbar_p$
of~$\Q_p$.

We begin by recalling some rings used in integral $p$-adic Hodge
theory (in particular, in defining $(\varphi,\Gamma)$-modules and
Breuil--Kisin modules), before introducing versions of them with
coefficients in certain $p$-adically complete $\Zp$-algebras. We then
prove (almost) descent results that allow us to relate $\varphi$-modules over
various different rings,
before introducing the
$(\varphi,\Gamma)$-modules which we work with throughout the book.

\section{Rings}\label{subsec:rings}
\subsection{Perfectoid fields and their tilts}
As usual, we let $\C$ denote the completion of the algebraic closure \index{$\C$}
$\Qbar_p$ of $\Q_p$.
It is a perfectoid field,
whose tilt
$\C^{\flat}$ is a complete  \index{$\C^{\flat}$}
non-archimedean valued perfect field of characteristic~$p$.
If $F$ is a perfectoid closed subfield of $\C$,
then its tilt $F^{\flat}$ is a closed, and perfect,
subfield of $\C^{\flat}$.

We let $\cO_{\C}$ denote the ring of integers in $\C$. \index{$\cO_{\C}$}
Its tilt $\cO_{\C}^{\flat}$ \index{$\cO_{\C}^{\flat}$} is then the ring of integers
in $\C^{\flat}$.
Similarly, if $\cO_F$ denotes the ring of integers
in $F$, then $\cO_F^\flat$ is the ring
of integers in~$F^\flat$.

We may form the rings of Witt vectors
$W(\C^\flat)$ and $W(F^\flat)$, \index{$W(\C^\flat)$}
and the rings of Witt vectors
$W(\cO_{\C^\flat})$ and $W(\cO_F^\flat)$;
	following the standard convention,
	we typically denote $W(\cO_{\C^\flat})$  by $\Ainf$.\index{$\Ainf$}

Each of these rings of Witt vectors
is a $p$-adically complete ring,
but we always consider them as topological rings \index{weak topology}
by endowing them with a finer topology, the so-called {\em weak
	topology}, which admits the following description: 
If $R$ is any of $\C^\flat$, $F^\flat$, $\cO_\C^\flat$,
or $\cO_F^\flat$, endowed with its natural (valuation)
topology,
then there is a canonical identification (of sets):
$$W_a(R) \iso R\times \cdots \times R \text{ ($a$ factors)},$$
and we endow $W_a(R)$ with the product topology.
We then endow $W(R) := \varprojlim_a W_a(R)$ with the inverse 
limit topology.

This topology admits the following more concrete description 
(in the general case of a perfectoid $F$; setting $F = \C$ 
recovers that particular case):
If $x$ is any element of $\cO_F^\flat$ of positive
valuation, and if $[x]$ denotes the Teichm\"uller 
lift of $x$, then we endow
$W_a(\cO_F^\flat)$ with the $[x]$-adic topology,
so that $W(\cO_F^\flat)$ is then endowed with
the $(p,[x]$)-adic topology.
The topology on $W_a(F^\flat)$ is then characterized 
by the fact that $W_a(\cO_F^\flat)$ is an open subring
--- concretely $W_a(F^\flat) = W_a(\cO_F^\flat)[\frac{1}{[x]}] 
= \varinjlim_i W_a(\cO_F^\flat)$ (the transition maps being
given by multiplication by~$[x]$, and each transition map
being an open and closed embedding),
while the topology on $W(F^\flat)$ is the
inverse limit topology --- concretely, $W(F^\flat)$
is the $p$-adic completion~$\widehat{W(\cO_F^\flat)[\frac{1}{[x]}]}$.

Apart from the case of $\C$ itself, 
there are two main examples of perfectoid $F$ that will be of importance
to us; see Example~\ref{ex:Sen} for a justification of the claim that
these fields are indeed perfectoid.

\begin{example}[The cyclotomic case]
	\label{ex:cyclo}
	We write $K(\zeta_{p^{\infty}})$ to denote the 
	extension of~ $K$ obtained by adjoining all $p$-power
	roots of unity.  It is an infinite degree Galois extension
	of $K$, whose Galois group is naturally
	identified with an open subgroup of $\Z_p^{\times}$.
	We let $K_{\cyc}$ denote the unique subextension of \index{$K_{\cyc}$}
	$K(\zeta_{p^{\infty}})$ whose Galois group
	over $K$ is isomorphic to $\Z_p$ (so $K_{\cyc}$ is
	the ``cyclotomic $\Z_p$-extension'' of $K$).
	If we let $\Khat_{\cyc}$ denote the closure \index{$\Khat_{\cyc}$}
	of $K_{\cyc}$ in~$\C$,
	then $\Khat_{\cyc}$ is a perfectoid subfield of~$\C$.
\end{example}

\begin{example}[The Kummer case]
	\label{ex:kummer}
	If we choose a uniformizer $\pi$ of $K$,
	as well as a compatible system of $p$-power roots
	$\pi^{1/p^n}$ of $\pi$ (here, ``compatible'' has
	the obvious meaning, namely that $(\pi^{1/p^{n+1}})^p = \pi^{1/p^n}$),
	then we define $K_{\infty} = K(\pi^{1/p^{\infty}}) :=
	\bigcup_n K(\pi^{1/p^n}).$ \index{$K_{\infty}$}
	If we let $\Khat_{\infty}$ denote the closure
	of $K_{\infty}$ in~$\C$,
	then $\Khat_{\infty}$ is again a perfectoid subfield of~$\C$.
\end{example}

\begin{rem}
  \label{rem: Krasner Ax Tate Sen}Let~$L$ be an algebraic extension
  of~$\Qp$, and let~$\Lhat$ be the closure of~$L$ in~$\C$. By
  Krasner's lemma (see~\cite[Prop.\ 9.1.16]{Gabber--Ramero} for full
  details), the field~$\Lhat\otimes_L\Qpbar$ is an algebraic closure
  of~$\Lhat$, so that the absolute Galois groups of~$L$ and~$\Lhat$
  are canonically identified. The action of~$G_L$ on~$\Qpbar$ extends
  to an action on~$\C$, and by a theorem of
  Ax--Tate--Sen~\cite{MR0263786}, we have~$\C^{G_L}=\Lhat$. We will in
  particular make use of these facts in the cases ~$L=K_{\cyc}$
  and~$L=K_\infty$.
\end{rem}


\subsection{Fields of norms}
If $L$ is
an infinite strictly arithmetically profinite (strictly APF) extension \index{strictly APF}
of $K$ in $\Qbar_p$, 
then we may form the field of norms $X_K(L)$, \index{field of norms}
as in~\cite{MR526137} and~\cite[\S 2]{MR719763}.
This is a complete discretely valued field of characteristic $p$,
whose residue field is canonically identified with that of
$L$. 
(We don't define the notion of strictly APF here, but simply refer to
\cite[Def.~1.2.1]{MR719763}
for the definition.  We will only
apply these notions to Examples~\ref{ex:cyclo} and~\ref{ex:kummer},
in which case the end result of the field of norms
construction can be spelt out quite explicitly.)
The following theorem is well known.
\begin{thm}
	\label{thm:field of norms}
	Let $L$ be an infinite strictly APF extension of $K$ in $\Qbar_p$.
	\begin{enumerate}
		\item The closure $\widehat{L}$ of $L$ in $\C$
			is perfectoid.
		\item 
			There is a canonical
	embedding $X_K(L) \hookrightarrow (\widehat{L})^\flat$,
	which
	identifies the target with the completion of the perfect
	closure of the source.
\end{enumerate}
\end{thm}
\begin{proof}
  To show that $\widehat{L}$ is perfectoid, it suffices to show that
  it is not discretely valued, and that the absolute Frobenius on
  $\cO_L/p\cO_L$ is surjective.  
  This surjectivity
  follows from \cite[Cor.~4.3.4]{MR719763}. It follows immediately
  from the description in~\cite[\S 1.4]{MR719763} of APF extensions as
  towers of elementary extensions that~$L$, and thus~$\widehat{L}$, is
  not discretely valued.  The canonical embedding of~(2) is
  constructed in \cite[\S 4.2]{MR719763}, and its claimed property is
   proved in \cite[Cor.~4.3.4]{MR719763}.
\end{proof}

\begin{example}
	\label{ex:Sen}
	A theorem of Sen \cite{MR0319949} 
	shows that
	if the Galois group of the Galois closure of $L$ 
	over $K$ is a $p$-adic Lie group, and the induced extension of
        residue fields is finite, then
	$L$ is a strictly APF extension of $K$.
	(Sen's theorem shows that the Galois
       closure of $L$ is strictly APF over $K$; see
       \cite[Ex. 1.2.2]{MR719763}. It then follows
       from 
       \cite[Prop.~1.2.3~(iii)]{MR719763}
       that $L$ itself is strictly APF over $K$.)
       As a consequence, we see that the theory of the field of norms,
       and in particular Theorem~\ref{thm:field of norms},
       applies in the
       cases of Examples~\ref{ex:cyclo} and~\ref{ex:kummer}.
       \end{example}

       \subsection{Thickening fields of norms}
       In the context of Theorem~\ref{thm:field of norms},
       given that one may thicken $(\widehat{L})^\flat$ to
       the flat $\Z/p^a$-algebra $W_a\bigl((\widehat{L})^\flat\bigr)$,
       it is natural to ask whether $X_K(L)$ admits a similar
       such thickening, such that the embedding $X_K(L) \hookrightarrow
       (\widehat{L})^\flat$
       lifts to an embedding of the
       corresponding thickenings.  One would furthermore like
       such a lifted embedding to be compatible with various auxiliary
       structures, such as the action of Frobenius on
       $W_a\bigl((\widehat{L})^\flat\bigr)$,
       or the action of the Galois group
       $\Gal(L/K)$, in the case when $L$ is Galois over $K$
       (i.e.\ one would like the image
       of this embedding to be stable under these actions).

       Since $X_K(L)$ is imperfect, there is no canonical
       thickening of $X_K(L)$ over $\Z/p^a$, and
       as far as we know, there is no simple or general answer to the 
       question of whether such thickenings and embeddings exist
       with desirable 
       extra properties, such as being Frobenius or Galois stable.
       However, in the cases of Examples~\ref{ex:cyclo} and~\ref{ex:kummer},
       such thickenings and thickened embeddings can be constructed
       directly, as we now recall.

       \subsection{The cyclotomic case}\label{subsubsec: cyclotomic}
       The extension $K(\zeta_{p^{\infty}})$ of $K$ is infinite and
       strictly APF, as well as being Galois over $K$.  If we write
       $\tGamma_K := \Gal(K(\zeta_{p^{\infty}})/K)$, then the
       cyclotomic character induces an embedding
       $\chi:\tGamma_K \hookrightarrow \Z_p^{\times}$.  Consequently,
       there is an isomorphism
       $\tGamma_K \cong \Gamma_K \times \Delta$, where
       $\Gamma_K \cong \Z_p$ and $\Delta$ is finite.  We have
       $K_{\cyc} = (K(\zeta_{p^{\infty}}))^{\Delta}$. Later in
       the book, $K$ will typically be fixed, and we will
       write~$\Gamma$ for~$\Gamma_K$.

       Suppose for a moment that $K = \Q_p$.
       If we choose a compatible system of $p^n$th roots of $1$,
       then these give rise in the usual way to an element
       $\varepsilon \in (\widehat{\Qp(\zeta_{p^{\infty}})})^{\flat}.$
              If we identify the field of norms $X_{\Q_p}(\Qp(\zeta_{p^{\infty}}))$ 
       with a subfield of $(\widehat{\Qp(\zeta_{p^{\infty}})})^\flat$ via the embedding
       $X_{\Q_p}(\Qp(\zeta_{p^{\infty}})) \hookrightarrow (\widehat{\Qp(\zeta_{p^{\infty}})})^\flat,$
       then $\varepsilon - 1 $ is a uniformizer of $X_{\Q_p}(\Qp(\zeta_{p^{\infty}}))$.  
       As usual, let $[\varepsilon]$ denote the Teichm\"uller 
       lift of $\varepsilon$ to an element of
       $W\bigl( \cO_{\widehat{\Qp(\zeta_{p^{\infty}})}}^\flat).$
       There is then a continuous embedding
       $$\Z_p[[T]] \hookrightarrow 
       W( \cO_{\widehat{\Qp(\zeta_{p^{\infty}})}}^\flat)$$
       (the source being endowed with its $(p,T)$-adic topology,
       and the target with its weak topology),
       defined via $T \mapsto [\varepsilon]- 1$.
       We denote the image of this embedding by $(\A'_{\Q_p})^+$.
       This embedding extends to an embedding
       $$\widehat{\Z_p((T))} \hookrightarrow W\bigl( (\widehat{\Qp(\zeta_{p^{\infty}})})^\flat
       \bigr)$$
       (here the source is the $p$-adic completion of the Laurent
       series ring $\Z_p((T))$),
       whose image we denote by $\A'_{\Q_p}$.

       We now return to the case of general $K$. We will compare 
       this case with the case of $\Q_p$. Since $K\Qp(\zeta_{p^\infty}) = K(\zeta_{p^\infty})$,
       the theory of the field of norms gives an identification
       $X_K(K(\zeta_{p^\infty})) = X_{\Q_p}(K(\zeta_{p^\infty}))$
       \cite[Rem.~2.1.4]{MR719763},
       and shows that this field 
       is a separable extension of $X_{\Q_p}(\Qp(\zeta_{p^\infty}))$
       \cite[Thm.~3.1.2]{MR719763},
       if we regard both as embedded in $\C^\flat$ via
       the embeddings of Theorem~\ref{thm:field of norms}. To ease
       notation, from now on we
       write~$\E'_K=X_K(K(\zeta_{p^\infty}))$,
       and~$\E_K:=(\E'_K)^\Delta$, so that ~$\E'_K/\E_K$ is a
       separable extension. We write $(\E'_K)^+$, $\E_K^+$ for the
       respective rings of integers. We write~$\varphi$ for the
       ($p$-power) Frobenius on~$E_K'$, $E_K$, $(\E'_K)^+$ and~ $\E_K^+$.

       We note now that $B'_{\Q_p} := \A'_{\Q_p}[1/p]$
       is a discretely valued field
       admitting $p$ as a uniformizer,
       with residue field $\E'_{\Q_p}$. 
       There is then a unique finite unramified extension of
       $\B'_{\Q_p}$ contained in the field $W(\C^\flat)[1/p]$
       with residue field $\E'_K$.
       We denote this extension by $\B'_K$; it is again
       a discretely valued field, admitting $p$ as a uniformizer,
       and we let $\A'_K$ denote its ring of integers;
       equivalently, we have $\A'_K = \B'_K \cap W(\C^\flat).$
       We see that $\A'_K$ is a discrete valuation ring,
       admitting $p$ as a uniformizer, and that
       $\A'_K/p \A'_K = \E'_K.$ There is a natural lift of the
       Frobenius~$\varphi$ from~$\E'_K$ to~$\A'_K$.

       The action of $\tGamma_K$ on $\E'_K=X_K(K(\zeta_{p^\infty}))$ induces
       an action of $\tGamma_K$ on $\A'_K$, and we write
       $\A_K := (\A'_K)^{\Delta}$. \index{$\A_K$}  This is the ring of integers in
       the discretely valued field $\B_K := \A_K[1/p] = (\B'_K)^{\Delta}$,
       and has residue field equal to $\E_K.$ The actions of
       $\tGamma_K$ on~$\A'_K$ and of~$\Gamma_K$ on~$\A_K$ commute with~$\varphi$.

       If we let $T'_K$ denote a lift of a uniformizer of
       $\E'_K$, and let $k'_{\infty}$ denote the residue field
       of $K(\zeta_{p^\infty})$, then there is an isomorphism
       $\widehat{W(k'_{\infty})((T'_K))} \iso \A'_K.$  (Here the source
       denotes the $p$-adic completion.)  Similarly if~$k_\infty$
       denotes the residue field of~$\Kcyc$, and~$T_K$ is a lift of a uniformizer of
       $\E_K$, then there is an isomorphism
       $\widehat{W(k_{\infty})((T_K))} \iso \A_K.$

       For a general extension~$K/\Qp$ it is hard to give an explicit
      formula for the actions of~$\varphi$ and~$\Gammat_K$ on
       $W(k_\infty')((T'_K))^{\wedge}$, but in some of our arguments
       it is useful to reduce to a special case where we can use
       explicit formulae. If~$K=K_0$ (that is, if~$K/\Qp$ is
       unramified) then we have such a description as follows: 
       we have $\E'_{K_0}=k((\varepsilon-1))$,
       and ~$\A'_{K_0}=W(k((\varepsilon-1)))$, and
       we set~$T'_{K_0}=[\varepsilon]-1$, with the square brackets
       denoting the Teichm\"uller lift.

The actions of~$\varphi$ and~$\gamma\in\Gammat_{K_0}$
on~$T'_{K_0}\in\A_{K_0}'$ are given by the explicit formulae
\numequation\label{eqn: phi on T basic case with
  epsilon}\varphi(T'_{K_0})=(1+T'_{K_0})^p-1, \end{equation}
\numequation\label{eqn: gamma on T basic case with
  epsilon}\gamma(1+T'_{K_0})=(1+T'_{K_0})^{\chi(\gamma)},
\end{equation}where~$\chi:\Gammat_{K_0}\to\Z_p^\times$ denotes the
cyclotomic character. We set $(\A_{K_0}')^+=W(k)[[T'_{K_0}]]$, which is
visibly~$(\varphi,\Gammat_{K_0})$-stable. We set
$T_{K_0}=\tr_{\A_{K_0}'/\A_{K_0}}(T'_{K_0})$ and
$\A_{K_0}^+=W(k)[[T_{K_0}]]$; then we
have~$\A_{K_0}=W(k)((T_{K_0}))^{\wedge}$, and~$\A_{K_0}^+$ is
$(\varphi,\Gamma_{K_0})$-stable (if $p>2$ this is~\cite[Prop.\
A.3.2.3]{MR1106901}, and if~$p=2$ the same statements hold, as
explained in~\cite[\S1.1.2.1]{MR1693457}). After possibly
replacing~$T_{K_0}$ by ~$(T_{K_0}-\lambda)$ for some~$\lambda\in
pW(k)$, by~\cite[Prop.\ 4.2, 4.3]{MR3322781} we can and do assume that
$\varphi(T_{K_0})\in T_{K_0}\A_{K_0}^+$, and that $g(T_{K_0})\in
T_{K_0}\A_{K_0}^+$ for all~$g\in\Gamma_{K_0}$.

\begin{defn}\index{basic (field~$K$)}
  \label{defn: basic}We say that~$K$ is \emph{basic} if it is
  contained in~$K_0(\zeta_{p^\infty})$. 
  (Although it won't play a role in what follows,
 we remark that this implies in particular that $K$ is abelian over $\Q_p$.)
\end{defn}

If~$K$ is basic then we have
$K(\zeta_{p^\infty})=K_0(\zeta_{p^\infty})$, so
that~$\bE'_K=\bE'_{K_0}$, and we take $\A'_K=\A'_{K_0}$,
$T'_K=T'_{K_0}$, $(\A_K')^+=(\A_{K_0}')^+$, and similarly for~$\A_K$,
with the actions of~$\Gammat_K$ and~$\Gamma_K$ being the restrictions
of the actions of~$\Gammat_{K_0}$ and~$\Gamma_{K_0}$. 

If~$K$ is not basic then 
as explained above, we may still choose some~$T_K'$ so that
$\A_K'=W(k_\infty')((T'_K))^{\wedge}$, and some~$T_K$ so that $\A_K=W(k_\infty)((T_K))^{\wedge}$. 
Having done so, we set $(\A_K')^+=W(k'_\infty)[[T'_K]]$,  $\A_K^+=W(k_\infty)[[T_K]]$, where the topology on~$W(k_\infty)[[T]]$
is as usual the~$(p,T)$-adic topology.\index{$\A_K^+$} 

\begin{rem}\label{rem: not phi stable T}
For an arbitrary~$K$, it is not necessarily possible to choose~$T_K$ so
  that~$\A_K^+$ is $\varphi$-stable and $\Gamma_K$-stable (and similarly
  for~$(\A_K')^+$). 
  In fact, it is not even possible to choose~$T_K$ so
  that~$\A_K^+$ is $\varphi$-stable. Indeed, by~\cite[Prop.\
  4.2, 4.3]{MR3322781}, 
  if~$\A_K^+$ is $\varphi$-stable then it is also~$\Gamma_K$-stable. As
  explained in~\cite[\S1.1.2.2]{MR1693457}, it is possible to choose $T_K$
  so that $\A_K^+$ is $\varphi$ and $\Gamma_K$-stable precisely
  when 
  $K$ is contained in an unramified extension of~$K_0(\zeta_{p^\infty})$,
  i.e.\ when $K$ is 
  abelian over~$\Qp$.
\end{rem}


       \subsection{The Kummer case}\label{subsubsec: Kummer}
       
       Set~$\gS=W(k)[[u]]$, \index{$\gS$} with a~$\varphi$-semi-linear
       endomorphism~$\varphi$ determined by~$\varphi(u)=u^p$, and
       let~$\cO_{\cE}$ \index{$\cO_{\cE}$} be the $p$-adic completion of~$\gS[1/u]$. The
       extension $K_\infty= K(\pi^{1/p^{\infty}})$ of~$K$ is infinite
       and strictly APF (but not Galois). The choice of compatible
       system of $p$-power roots of~$\pi$ gives an element
       $\pi^{1/p^\infty}\in \cO_{\widehat{K_\infty}}^\flat$, and there
       is a continuous $\varphi$-equivariant
       embedding \[\gS\into W(\cO_{\widehat{K_\infty}}^\flat)\]
       which lifts the morphism $X_K(K_{\infty}) \hookrightarrow
(\widehat{K}_{\infty})^{\flat},$
       defined by
       sending $u\mapsto [\pi^{1/p^\infty}]$.
This embedding extends
       to a continuous $\varphi$-equivariant
       embedding \[\cO_{\cE}\into W((\widehat{K_\infty})^\flat).\]

\subsection{Graded and rigid analytic techniques}In
Appendix~\ref{sec: rigid analytic perspective} we prove a number of
results using rigid analysis and graded techniques, which we will
apply in the next section in order to prove some basic facts about our
coefficient rings. In order to apply these results to the various
rings introduced above, we need to verify the hypotheses introduced in
\ref{subsec:setting}, which are as follows. There is an Artinian 
 local ring $R$, which in our present context we take to be the ring~$\Z/p^a$.
Then we work with a $\Z/p^a$-algebra $C^+$, and an element $u \in C^+$, satisfying the
following properties:
\begin{enumerate}[label=(\Alph*)]
	\item\label{item: Cplus A} $u$ is a regular element (i.e.\ a non-zero divisor) of $C^+$.
	\item\label{item: Cplus B}  $C^+/u$ is a flat (equivalently, free) $\Z/p^a$-algebra.
	\item\label{item: Cplus C}  $C^+/p$ is a rank one complete valuation ring, and
		the image of $u$ in $C^+/p$ (which
		is necessarily non-zero, by \eqref{item: Cplus A}) is of
		positive valuation (i.e.\ lies in the maximal
		ideal of $C^+/p$).
\end{enumerate}
The following lemma gives the examples of this construction that we
will use in the next section.

\begin{lem}
  \label{lem: actual examples of Cplus and C that we use}The following
  pairs~$(C^+,u)$ satisfy axioms \emph{\ref{item: Cplus A}--\ref{item:
    Cplus C}} above.
  \begin{enumerate}
  \item $C^+=W_a(\cO_F^\flat)$, where~$F$ is any perfectoid
    field, and~$u$ is any element of~$W_a(\cO_F^\flat)$
    whose image in $\cO_F^\flat$ is of positive
    valuation. 
  \item $C^+=\A_K^+/p^a$, $u=T_K$.
  \item $C^+=W_a(k)[[u]]$.
  \end{enumerate}

\end{lem}
\begin{proof}
	Axioms~\ref{item: Cplus A} and~\ref{item:
    Cplus C} are clear from the definitions in each case, so we need
  only verify that $C^+/u$ is flat over~$\Z/p^a$ in each case.
  This is immediate in cases~(2) and (3), so we focus on case~(1).
  Rather than checking~\ref{item: Cplus B} directly in this case,
  we instead check that $W(\cO_F^\flat)/ u$ is flat over $\Z_p$;
  base-changing to $\Z/p^a$ then gives us~\ref{item: Cplus B}.
  For this, it suffices to observe that multiplication by $u$ is injective
  on each term of the short exact sequence
  $$0 \to W(\cO_F^\flat) \buildrel p \cdot \over \longrightarrow
  W(\cO_F^\flat) \longrightarrow \cO_F^\flat \to 0 \, ;$$
  the snake lemma then shows that multiplication by $p$ is injective
  on $W(\cO_F^{\flat})/u$, as required.
\end{proof}

We also note the following lemma.

\begin{lem}
  \label{lem: boundedness of Galois actions in our concrete
    settings}The actions of~$G_K$ on~$W_a(\C^\flat)$, 
  and of~$\Gamma_K$
  on~$\A_K/p^a$ and $\tA_K/p^a$, are continuous and bounded in the sense
  of Definition~{\em \ref{adefn: bounded group action}}.
\end{lem}
\begin{proof}It suffices to show that the action of $G_K$
  on~$W_a(\C^\flat)$ is continuous and bounded, as the other cases follow by
  restricting this action to the corresponding subgroup and subring.
  By Remark~\ref{arem: bounded doesn't depend on u}, we may assume
  that ~$u=[v]$ for some~$v\in\cO_{\C}^\flat$ with positive valuation;
  since the action of~$G_K$ on $\cO_{\C}^\flat$
  preserves the valuation, we then
  have~$G_K\cdot u^{M}W_a(\C^\flat) =  u^{M}W_a(\C^\flat)$ for
  all~$M\in\Z$, so the action is continuous and bounded.
\end{proof}

\section{Coefficients}\label{subsec: coefficients}
Our main concern throughout this book will be families of
\'etale $(\varphi,\Gamma)$-modules; an auxiliary role will also be played by
families of (various flavours of) \'etale $\varphi$-modules. Such a family
will be parameterized by an algebra of coefficients, typically 
denoted by $A$.   Since we will work throughout
in the context of formal algebraic stacks over $\Spf \Z_p$ (or closely 
related contexts), we will always assume that our coefficient ring $A$ is
a $p$-adically complete $\Zp$-algebra.
Frequently, we will work modulo some a fixed power $p^a$ of $p$,
and thus assume that $A$ is actually a $\Z/p^a$-algebra
(and sometimes we impose further conditions on $A$,
such as that of being Noetherian, or even of finite type over $\Z/p^a$). 

It is often convenient to introduce an auxiliary base ring for our
coefficients,
which we will take to be the ring of integers $\cO$ in a finite extension
$E$ of $\Q_p$; in this case, we will let $\varpi$ denote a uniformizer
of $\cO$, and $A$ will be taken to be a $p$-adically complete
(or, equivalently, a $\varpi$-adically complete) $\cO$-algebra,
or, quite frequently, an $\cO/\varpi^a$-algebra
(perhaps Noetherian, or even of finite type), for some power
$\varpi^a$ of $\varpi$.

For the moment, we put ourselves in the most general case;
that is, we assume that $A$ is a $p$-adically complete $\Z_p$-algebra,
and define versions of the various rings considered in
Section~\ref{subsec:rings} ``relative to $A$''.

Let $F$ be a perfectoid closed subfield of $\C$,
and if $a \geq 1$,
let $v$ denote an element of the maximal ideal
of $W_a(\cO_F^\flat)$
whose image
in $\cO_F^\flat$ is non-zero.
We then set
$$
W_a(\cO_F^\flat)_A =
W_a(\cO_F^\flat) \cotimes_{\Z_p}  A
:= \varprojlim_i \bigl( W_a(\cO_F^\flat)\otimes_{\Z_p}A\bigr)/v^i$$
(so that the indicated completion is the $v$-adic completion).
Note that any two choices of $v$ induce the same topology
on $W_a(\cO_F^\flat)\otimes_{\Z_p} A$, so
that 
$W_a(\cO_F^\flat) \cotimes_{\Z_p}  A$
is well-defined independent of the choice of $v$.
We then define
$$W_a(F^\flat)_A =  W_a(F^\flat) \cotimes_{\Z_p} A
:= W_a(\cO_F^\flat)_A[1/v];$$
this ring is again well-defined independently of the choice of $v$.

There are natural reduction maps
$$W_a(\cO_F^\flat) \cotimes_{\Z_p}  A \to
W_b(\cO_F^\flat) \cotimes_{\Z_p}  A$$
and
$$W_a(F^\flat) \cotimes_{\Z_p}  A \to
W_b(F^\flat) \cotimes_{\Z_p}  A,$$
if $a \geq b$,
so that we may define
$$
W(\cO_F^\flat)_A =
W(\cO_F^\flat)\cotimes_{\Z_p} A := \varprojlim_a W_a(\cO_F^\flat)_A,
$$
and similarly
$$
W(F^\flat)_A =
W(F^\flat)\cotimes_{\Z_p} A := \varprojlim_a W_a(F^\flat)_A.
$$

For certain choices of $F$, we introduce alternative notation
for the preceding constructions. 
In the case when $F$ is equal to $\C$ itself,
we write
$$\AAinf{A} :=  W(\cO_\C^\flat)_A.$$
In the case when $F$ is equal to $\widehat{K}_{\cyc}$,
we write 
$$\tA_{K,A} := W\bigl((\widehat{K}_{\cyc})^\flat\bigr)_A.$$
In the case when $F$ is equal to $\widehat{K}_{\infty}$,
we write
$$\tOEA := W\bigl((\widehat{K}_{\infty})^\flat\bigr)_A.$$

We next introduce various imperfect coefficient rings with $A$-coefficients:
If $T_K$ denotes a lift to $\A_K$ of a uniformizer of $\E_K$, 
defining a subring~$\A_K^+$ of $\A_K$ as in Section~\ref{subsec:rings}, 
then we write
\[\A^+_{K,A}=\A^+_K\cotimes_{\Zp}A := \varprojlim_{n} \A^+_K/(p,T_K)^n
	\otimes_{\Z_p} A = \varprojlim_{m} (\varprojlim_n \A^+_K/(p^m,T_K^n)
	\otimes_{\Z_p} A), \]
and
\[\A_{K,A} =\A_K\cotimes_{\Zp}A :=
	\varprojlim_{m} \bigl(\varprojlim_n (\A^+_K/(p^m,T_K^n)
	\otimes_{\Z_p} A) [1/T_K]\bigr). \]
Similarly, if $u$ denotes the usual element of $\gS$, whose
reduction modulo $p$
is a uniformizer of $X_K(K_{\infty})$, then we write 
\[\gS_{A}=\gS\cotimes_{\Zp}A := \varprojlim_{n} \gS/(p,u)^n
	\otimes_{\Z_p} A = \varprojlim_{m} (\varprojlim_n \gS/(p^m,u^n)
	\otimes_{\Z_p} A), \]
and
\[\OEA =\cO_\cE \cotimes_{\Zp}A :=
	\varprojlim_{m} \bigl(\varprojlim_n (\gS/(p^m,u^n)
	\otimes_{\Z_p} A) [1/u]\bigr). \]
Note that by definition the ring~$\gS_{A}$ is a power
series ring over~$W(k)\otimes_AA$, and $\OEA$ is the $p$-adic
completion of the corresponding Laurent series ring;  and similarly
for the rings~ $\A^+_{K,A}$ and
~ $\A_{K,A}$.

\subsection{A digression on flatness and completion}
We record some results giving sufficient conditions 
for flatness to be preserved after passage to certain inverse limits,
beginning with the
      following results from~\cite{MR2774689}.

      \begin{prop}
        \label{prop: FGK flat}Let~$R$ be a ring which is $x$-adically
        complete for some $x\in R$, and let ~$S$ be an~$R$-algebra
        which is also $x$-adically complete, and for which 
$R[1/x]$ and $S[1/x]$ are both Noetherian.
If for each~$k~\ge~1$, the induced morphism $R/x^k\to S/x^kS$ is {\em (}faithfully{\em )}
flat, then the morphism $R\to S$ is {\em (}faithfully{\em )} flat.
      \end{prop}
      \begin{proof}
The flatness claim is a special case of the discussion at the beginning
        of~\cite[\S 5.2]{MR2774689}.
The faithful flatness claim 
        is a special case of~\cite[Prop.\ 5.2.1~(2)]{MR2774689}.
It also follows from the flatness claim together with
Lemma~\ref{lem:faithful flatness} below.
      \end{proof}

\begin{remark}
Under the faithful flatness hypothesis of the preceding proposition,
it is proved in~\cite[Prop.\ 5.2.1~(1)]{MR2774689}
that Noetherianness of $R[1/x]$ is a consequence of
Noetherianness of $S[1/x]$.   This is a kind of {\em fpqc} descent
result for this property, which however we won't need in the present book.
\end{remark}

We next present
the following variation on
a result of Bhatt--Morrow--Scholze.

\begin{prop}
\label{prop:BMS flatness}
Let $R$ be a Noetherian ring which is $x$-adically complete, for some element
$x \in R$, let $S$ be an $x$-adically complete $R$-algebra,
and suppose that $S/x^n$ is {\em (}faithfully{\em )} flat over $R/x^n$ for every~$n~\geq~1$.
Then $S$ is {\em (}faithfully{\em )} flat over~$R$.
\end{prop}

\begin{remark}
In~\cite[Rem.\ 4.31]{2016arXiv160203148B}, the authors prove
the above proposition in the particular
case when $R$ and $S$ are flat $\Z_p$-algebras, and $x$ equals~$p$.
In this case they only need to assume that $S/p$ is flat over~$R/p$.
(Since $R$ and $S$ are both flat over~$\Z_p$,
it follows automatically that each $S/p^n$ is flat over~$R/p^n$.)  
Their proof makes use of the following key ingredients: that for a $p$-adically complete 
ring, any pseudo-coherent complex is derived complete; that for a complex 
over a $p$-torsion
free ring, the derived $p$-adic completion may be computed ``naively'';
and the Artin--Rees lemma for the $p$-adically complete and Noetherian ring~$R$.

Our argument is identical to theirs, once we confirm that these ingredients remain
available (with $p$-adic completions being replaced by $x$-adic completions). 
The first statement holds quite generally~\cite[\href{https://stacks.math.columbia.edu/tag/0A05}{Tag
    0A05}]{stacks-project}, and of course Artin--Rees holds
for any Noetherian ring~$R$.  Thus the main point is 
to verify that we may compute derived $x$-adic completions for the (generally
non-Noetherian ring) $S$ ``naively''.   
By~\cite[\href{https://stacks.math.columbia.edu/tag/0923}{Tag 0923}]{stacks-project}, 
this is possible provided that the $x$-power torsion in $S$ is bounded.  So
our task is to verify this boundedness (under our given hypotheses).
\end{remark}

We begin with the following general criterion for a ring to have bounded torsion.

\begin{lemma}
\label{lem:bounded torsion}
If $x$ is an element of the ring  $R$, then the following are equivalent:
\begin{enumerate}
\item
$R[x^m] = R[x^{\infty}]$
\item 
The morphism $R[x^i] \to R/x^m$ is injective for some~$i\geq 1$.
\item 
The morphism $R[x^i] \to R/x^m$ is injective for every~$i\geq 1$.
\end{enumerate}
\end{lemma}
\begin{proof}
 This is straightforward. Indeed, if the morphism $R[x^i] \to R/x^m$ is
 injective for some~$i\geq 1$, then it is in particular injective
 if~$i=1$. Suppose then that injectivity holds for~$i=1$.
 If $y\in R$ is such that~$x^{m+n}y=0$ for some~$n\ge 1$,
 then $x^{m+n-1}y$ is in the kernel of the morphism $R[x] \to R/x^m$,
 so $x^{m+n-1}y=0$, and by an easy induction we have $x^my=0$, as
 required.
 
 Conversely, if $t$ is in the kernel of the morphism $R[x^i] \to R/x^m$
 for some~$i$, then we can write $t=x^my$, and we have $x^it=0$; so
 $x^{m+i}y=0$. If $R[x^m] = R[x^{\infty}]$, then we have $x^my=0$, so
 $t=0$, as required.
\end{proof}

\begin{lemma}
\label{lem:torsion exact one}
If $x$ is any element of the ring~$R$, then for any $i,n\ge 1$ we have the exact sequence
$$R[x^i] \to R/x^n \buildrel x^i \over \longrightarrow R/x^{n+i} \to R/x^i \to 0.$$
\end{lemma}
\begin{proof}
This is immediately verified.
\end{proof}

\begin{lemma}
\label{lem:torsion exact two}
If $x$ is an element of the ring~$R$, and if $R[x^m] = R[x^{\infty}]$,
then 
$$0 \to R[x^i] \to R/x^n \buildrel x^i \over \longrightarrow R/x^{n+i} \to R/x^i \to 0$$
is an exact sequence, for any~$n~\geq~m$ and any $i \geq 1$.
\end{lemma}
\begin{proof}
This follows from Lemmas~\ref{lem:bounded torsion} and~\ref{lem:torsion exact
one}.
\end{proof}

\begin{proof}[Proof of Proposition~\ref{prop:BMS flatness}]
Since $R$ is Noetherian, we see that $R[x^m] = R[x^{\infty}]$ for some $m \geq 0$.
Thus, by Lemma~\ref{lem:torsion exact two},
for each $n~\geq~m$ and each $i \geq 1$,
we obtain an exact sequence
$$0 \to R[x^i] \to R/x^n \buildrel x^i \over \longrightarrow R/x^{n+i} \to R/x^i \to 0$$
Tensoring with the flat $R/x^{n+i}$-algebra $S/x^{n+i},$
we obtain an exact sequence
$$0 \to S\otimes_R R[x^i] \to S/x^n \buildrel x^i \over \longrightarrow S/x^{n+i} \to S/x^i \to 0$$
Passing to the inverse limit over $n$ (and taking into account that all the 
transition morphisms are surjective),
we obtain an exact sequence of $R/x^{n+i}$-modules
$$0 \to S\otimes_R R[x^i] \to S \buildrel x^i \over \longrightarrow S \to S/x^i \to 0$$
In particular, we find that $S\otimes_R R[x^i] \iso S[x^i]$,
for every~$i$, and so in particular
$S[x^m] = S[x^{\infty}]$, so that $S$ has bounded $x$-power torsion.

We now follow the proof of~\cite[Rem.\ 4.31]{2016arXiv160203148B}.
By (for example) \cite[\href{https://stacks.math.columbia.edu/tag/00M5}{Tag
    00M5}]{stacks-project}, it is enough to show that if~$M$ is a
  finitely generated $R$-module, then $M\otimes^{\mathbf{L}}_RS$ has
  cohomology concentrated in degree~$0$. Since~$R$ is Noetherian,
  \cite[\href{https://stacks.math.columbia.edu/tag/066E}{Tag
    066E}]{stacks-project} shows that $M$ is pseudo-coherent when
regarded as a complex
  of~$R$-modules in a single degree,
and so $M\otimes^{\mathbf{L}}_RS$ is a pseudo-coherent complex
  of~$S$-modules,
  by~\cite[\href{https://stacks.math.columbia.edu/tag/0650}{Tag
    0650}]{stacks-project}. Also,
since~$R$ and~$S$ are $x$-adically
  complete, it then follows
  from~\cite[\href{https://stacks.math.columbia.edu/tag/0A05}{Tag
    0A05}]{stacks-project} that both~$M$ and
  $M\otimes^{\mathbf{L}}_RS$ are derived $x$-adically complete.

Since $S$ has bounded $x$-power torsion,
it follows
from~\cite[\href{https://stacks.math.columbia.edu/tag/0923}{Tag 0923}]{stacks-project} 
that we may compute derived $x$-adic completions naively, i.e.\ via
$R\lim(\text{--}\otimes_S S/x^n).$  Thus
$$M\otimes^{\mathbf{L}}_R S
\iso R\lim( M\otimes^{\mathbf{L}}_R S/x^n) 
\iso R\lim\bigl( (M\otimes^{\mathbf{L}}_R R/x^n) \otimes^{\mathbf{L}}_{R/x^n} S/x^n\bigr).$$
Artin--Rees allows us to replace the pro-system
 $\{M\otimes^{\mathbf{L}}_R R/x^n\}$ with the pro-system
 $\{M/x^n\}$, 
and so we find that in fact
\begin{multline*}
M\otimes^{\mathbf{L}}_R S
\iso R\lim(M/x^n \otimes^{\mathbf{L}}_{R/x^n} S/x^n)
\\
\iso R\lim( M\otimes_R S/x^n) \iso \lim M\otimes_R S/x^n,
\end{multline*}
the penultimate isomorphism holding
since $S/x^n$ is flat over $R/x^n$, 
and the final isomorphism following from the fact that the transition morphisms
in the pro-system $\{ M\otimes_R S/x^n\}$ are surjective.
Thus indeed 
$M\otimes^{\mathbf{L}}_R S$ has cohomology supported in a single degree,
as required.

The claim about faithful flatness follows from the immediately
following Lemma~\ref{lem:faithful flatness}.
\end{proof}

\begin{lemma}
\label{lem:faithful flatness}
If $R\to S$ is a flat morphism of rings, if $I$ is an ideal in $R$ for which
$R$ is $I$-adically complete, and if the morphism
$R/I \to S/I$ is faithfully flat, then the morphism $R\to S$ is faithfully flat.
\end{lemma}
\begin{proof}
Since $R$ is $I$-adically complete, we see that $I$ is contained in the
Jacobson radical of~$R$.  Thus $\Max\Spec R/I = \Max\Spec R$,
and so our assumption that $R/I \to S/I$ is faithfully flat shows
that $\Max\Spec R$
 is contained in the image of $\Spec S/I \subseteq \Spec S$ in~$\Spec R$.
Since flat morphisms satisfy going down, we find that in fact
$\Spec S \to \Spec R$ is surjective, as claimed.
\end{proof}

      We also note the following lemma.

      \begin{lemma}
	      \label{lem:Witt flatness}
      If $F \hookrightarrow F'$ is an inclusion of
      perfectoid fields in characteristic $p$,
	      then the induced morphism $W_a(\cO_F)
	      \to W_a(\cO_{F'})$ 
	      is faithfully flat, for any~$a~\geq~1$.
      \end{lemma}
      \begin{proof}
Since $\cO_F$ is a B\'ezout ring,
being flat over $\cO_F$ is the same as being
torsion free, and so the inclusion
$\cO_{F} \hookrightarrow \cO_{F'}$
is a flat morphism,
and thus even faithfully flat, being a local morphism of local rings. 
A standard grading argument, applying Lemma~\ref{lem:graded flat} to 
the $p$-adic filtrations on source and target,
then shows that the inclusion
$W_a(\cO_{F}) \hookrightarrow W_a(\cO_{F'})$,
is faithfully flat, for each~$a~\geq~1$.
\end{proof}

\subsection{Back to coefficient rings}
      We now return to our discussion 
of coefficient rings, and record that, at least if ~$A$ is a finite
type 
      $\Z/p^a\Z$-algebra, the various natural maps between these rings
      are in fact (faithfully) flat injections.
We also show that the various maps of coefficient rings induced by
a (faithfully) flat morphism of finite type $\Z/p^a$-algebras
are again (faithfully) flat.

\begin{prop}
  \label{prop: maps of coefficient rings are faithfully flat injections}
Suppose that~$A\to B$ is a flat homomorphism
  of finite type $\Z/p^a\Z$-algebras for
  some $a\ge 1$. Then all the maps in the following diagram are flat.
  Furthermore the vertical arrows are all injections,
while the  horizontal arrows are all faithfully flat {\em (}and
so in particular also injections{\em )}.
If $A\to B$ is furthermore faithfully flat, then the same is true of the
diagonal arrows.

\[\xymatrix{&\A_{K,B}^+\ar[rr]\ar[dd]&&\AAinf{B}\ar[dd]&&\gS_B
    \ar[dd]\ar[ll]\\
    \A_{K,A}^+\ar[rr]\ar[dd]\ar[ur]&&\AAinf{A}\ar[dd]\ar[ur]&&\gS_A
    \ar[dd]\ar[ll]\ar[ur]\\
&\A_{K,B}\ar[r]&\tA_{K,B}\ar[r]& W(\C^\flat)_B &
\widetilde{\cO}_{\cE,B}\ar[l]& \OEB\ar[l]\\
\A_{K,A}\ar[r]\ar[ur]&\tA_{K,A}\ar[r]\ar[ur]& W(\C^\flat)_A\ar[ur] &
\widetilde{\cO}_{\cE,A}\ar[l]\ar[ur]& \OEA\ar[l]\ar[ur]}\]
  
\end{prop}
\begin{proof}From left to right, the vertical maps are given by
  inverting the elements $T$, $v$, and~$u$ respectively (where as
  above, $v$ is any element of the maximal ideal of~$W_a(\cO^\flat_\C)$, whose
  image in~$\cO^\flat$ is nonzero; in particular,
  we can take $v$ equal to either $T$ or $u$).
  Since localizations are flat, 
  to prove the proposition for these maps it
  is enough to note that the sources of the vertical maps are
  respectively $T$-, $v$-, and $u$-torsion free, by
  Lemma~\ref{lem:non-zero divisor}. 

  We now turn to the horizontal maps, where we will make repeated use
  of Proposition~\ref{prop: FGK flat}. It evidently suffices to treat the maps
  with $A$-coefficients. Since
  each of $\A_{K,A}$, $\tA_{K,A}$, 
  $W(\C^\flat))_A$,
  $\widetilde{\cO}_{\cE,A}$, and $\OEA$
  is Noetherian by Proposition~\ref{prop:
    rigid analysis FGK}~(1), it follows from Proposition~\ref{prop: FGK flat} (and
  the flatness of the vertical maps) that we need only show that the maps
  $\gS_A/u^i\to
  (W_a(\cO_{\widehat{K_\infty}}^\flat)\otimes_{\Zp}A)/u^i\to
  \AAinf{A}/u^i$ and
  $\A^+_{K,A}/T^i\to
  (W_a(\cO_{\widehat{\Kcyc}}^\flat)\otimes_{\Zp}A)/T^i\to
  \AAinf{A}/T^i$ are faithfully flat for each~$i~\ge~1$. 

  For the faithful flatness of
  $\gS_A/u^i
  \to(W_a(\cO_{\widehat{K_\infty}}^\flat)\otimes_{\Zp}A)/u^i$,
  it is
  enough to show that $\gS/p^a\to W_a(\cO_{\widehat{K_\infty}}^\flat)$
  is faithfully flat. By a standard grading argument (see
  Lemma~\ref{lem:graded flat}), it suffices in turn to prove that
  $k[[u]]\to \cO_{\widehat{K_\infty}}^\flat$ is faithfully flat, so in
  turn it is enough to check that~$\cO_{\widehat{K_\infty}}^\flat$ is
  $u$-torsion free, which is clear. To see that
  $(W_a(\cO_{\widehat{K_\infty}}^\flat)\otimes_{\Zp}A)/u^i\to
  \AAinf{A}/u^i$ is faithfully flat, it suffices to note that
  $W_a(\cO_{\widehat{K_\infty}}^\flat)\to W_a(\cO_{\C}^\flat)$ is
  faithfully flat,
  by Lemma~\ref{lem:Witt flatness}.
  The faithful flatness of the remaining horizontal maps
  is proved in exactly the same way.

Finally, the (faithful) flatness of the diagonal maps is now immediate
from another application of Proposition~\ref{prop: FGK flat}, together
with the (faithful) flatness of $A\to B$.
\end{proof}


 Recall that if~$\Acirc$ is a $p$-adically
    complete $\cO$-algebra, then~$\Acirc$ is said to be \emph{topologically of finite
    type over~$\cO$} if it can be written as a quotient of a
  restricted formal power series ring in finitely many variables~$\cO\langle\langle
  X_1,\dots,X_n\rangle\rangle$; equivalently, if and only if~$\Acirc\otimes_{\cO}k$ is a finite
    type $k$-algebra (\cite[\S0, Prop.\
    8.4.2]{MR3752648}). In particular, if~$\Acirc$ is a
    $\cO/\varpi^a$-algebra for some~$a\ge 1$, then~$\Acirc$ is
    topologically of finite type over~$\cO$ if and only if it is of
    finite type over~$\cO/\varpi^a$.
    Since $\cO$ is finite over $\Z_p$, it is equivalent for $A$ to be
    topologically of finite type over $\cO$ or over $\Z_p$, and so for the
    moment we consider the maximally general case of a topologically finite type
    $\Z_p$-algebra. 

\begin{remark} We don't know if the analogue of 
  Proposition~\ref{prop: maps of coefficient rings are faithfully flat injections}
holds when $A$ and $B$ are taken to be merely $p$-adically complete 
and topologically of finite type over~$\Z_p$,
rather than of finite type over $\Z/p^a$ for some~$a \geq 1$.   Each of the various
coefficient
rings over $A$ and $B$ is (by definition) formed by first forming the
corresponding coefficient ring over each $A/p^a$ or~$B/p^a$, and then taking
an inverse limit.  Since the formation of inverse limits is left exact,
we see that the horizontal and vertical arrows in the diagram are injective,
but we don't know in general that the various arrows are flat (although
we have no reason to doubt it).  One can however establish some partial results,
using~Proposition~\ref{prop:BMS flatness}.  
We record here one such result,
which we will need later on.
\end{remark}

\begin{prop}
\label{prop:top f.t. flatness}
If $A$ is a $p$-adically complete $\Z_p$-algebra which is topologically of
finite type, 
then the natural morphism $\gS_A \to \AAinf{A}$ is faithfully flat.
\end{prop}
\begin{proof}
This follows from Propositions~\ref{prop:BMS flatness}
and~\ref{prop: maps of coefficient rings are faithfully flat injections},
once we note that $\gS_A$ is Noetherian when $A$ is $p$-adically
complete and topologically of finite type; indeed since $A$ is in
particular Noetherian, so is the power series ring $W(k)\otimes_{\Z_p} A [[u]].$
\end{proof}

We conclude this initial discussion of coefficient rings by 
explaining how
the action of~$\varphi$ on the various rings~$\gS$, $\Ainf$ and so on
extends to a continuous action on the corresponding rings~$\gS_A$,
$\AAinf{A},$ etc.,
and similarly for the various Galois actions.
For this,
it is convenient for us to briefly digress,
and to introduce the
following situation, which will also be useful for us in Chapter~\ref{sec:
phi modules and phi gamma modules}. 

\begin{situation}
  \label{subsubsec:general framework}
  Fix a finite extension~$k/\Fp$ and
  write~$\Aplus:=W(k)[[T]]$. Write~$\A$ for the $p$-adic completion
  of~$\Aplus[1/T]$.


  If $A$ is a $p$-adically complete~$\Zp$-algebra, we write
  $\Aplus_A:=(W(k)\otimes_{\Zp}A)[[T]]$; we equip~$\Aplus_A$ with its
  $(p,T)$-adic topology, so that it is a topological $A$-algebra
  (where~$A$ has the $p$-adic topology). Let $\AAA_A$ be the $p$-adic
  completion of~ $\Aplus_A[1/T]$, which we regard as a topological
  $A$-algebra by declaring $\Aplus_A$ to be an open subalgebra. Note
  that the formation of~$\gS_A$, $\OEA$, $\A_{K,A}^+$ and~$\A_{K,A}$
  above are particular instances of this construction.

  Let~$\varphi$ be a ring endomorphism of~$\A$ which is congruent to
  the ($p$-power) Frobenius endomorphism modulo~$p$. We say
  that~$\A^+$ is $\varphi$-stable if~$\varphi(\A^+)\subseteq \A^+$.

By~\cite[Lem.\ 5.2.2 and
  5.2.5]{EGstacktheoreticimages}, if~$\A^+$ is $\varphi$-stable, then
  $\varphi$ is faithfully flat, and induces the usual 
  Frobenius on~$W(k)$; the same arguments show that this is true
  for~$\varphi$ on~$\A$, even
  if $\A^+$ is not $\varphi$-stable.
\end{situation}

 We have the following variant
  of~\cite[Lem.\ 5.2.3]{EGstacktheoreticimages}.

  \begin{lem}
    \label{lem: phi on powers of T without assuming A plus
      stable}Suppose that we are in Situation~{\em \ref{subsubsec:general
      framework}}. Then for
    each integer $a\ge 1$ there is an integer~$C\ge 0$ such that for
    all~$n\ge 1$, we have $\varphi(T^n)\in T^{pn-C}\A^++p^a\A$. In
    particular, the action of~$\varphi$ on~$\A$ is continuous.
  \end{lem}
  \begin{proof}
    For some~$h\ge 0$ we have $\varphi(T)\in T^{-h}\A^++p^a\A$. Write
    $\varphi(T)=T^p+pY$, so that~$\varphi(T^n)=(T^p+pY)^n$. If we
    expand using the binomial theorem, then every term on the right
    hand side is either divisible by~$p^a$, or is a multiple of
    $T^{p(n-r)}Y^r$ for some $0\le r \le a-1$. It follows that we can
    take~$C=(a-1)(p+h)$. 
  \end{proof}

  \begin{lem}\label{lem: varphi extends to A}
	  If~$A$ is a $p$-adically complete $\Zp$-algebra,
    then in Situation {\em \ref{subsubsec:general framework}}, the
    endomorphism~$\varphi$ of~$\A$ extends uniquely to an $A$-linear continuous
  endomorphism of~$\AAA_A$ \emph{(}which we continue to denote
  by~$\varphi$\emph{)}. If~$\A^+$ is $\varphi$-stable, then $\A^+_A$
  is $\varphi$-stable for all~$A$.  
  \end{lem}
  \begin{proof}Since~$\varphi$ is $\Z_p$-linear
    by definition, and since the topologies on~$\A^+_A$ and~$\AAA_A$
 are defined by passage to the limit modulo~$p^a$
    as $a\to\infty$, we are immediately reduced to the case that~$A$
    is a $\Z/p^a\Z$-algebra. 
In this case the result is immediate from Lemma~\ref{lem: phi on powers of T without assuming A plus
      stable} and Lemma~\ref{lem:extending endomorphisms}. 
  \end{proof}

  \begin{lem}\label{lem: varphi extends to A and Galois
      etc}Let~$A$ be a $p$-adically complete $\Zp$-algebra.
    \begin{enumerate}
\item The endomorphism
	\label{item: continuity of phi 2} $\varphi$
  of~$\gS$ {\em (}resp.\ $\cO_{\cE}$, $\A_{K}$,
  $W(\cO_F^\flat)$, $W(F^\flat)${\em )}
	extends uniquely to an $A$-linear continuous
  endomorphism
  of~$\gS_A$ {\em (}resp.\ $\OEA$, $\A_{K,A}$,
  $W(\cO_F^\flat)$, $W(F^\flat)_A${\em )}, which we again denote by~$\varphi$.
  If~$K$ is furthermore basic,
  then the endomorphism $\varphi$ of $\A_{K,A}$ preserves~$\A^+_{K,A}$.
\item\label{item: continuity of Galois action 3}
	The continuous $G_K$-action on 
  $W(\cO_{\C}^\flat)$ and $W(\C^\flat)$ {\em  (}
  resp.\ the continuous~$\Gamma_K$-action on~$\A_{K}$ and~$\tA_{K}${\em )}
  extends to a continuous $A$-linear action
	of~$G_K$ on
  $W(\cO_{\C}^\flat)_A$ and $W(\C^\flat)_A$ {\em (}
  resp.\ of~$\Gamma_K$ on~$\A_{K,A}$ and~$\tA_{K,A}${\em )}.
    \end{enumerate}
  \end{lem}
  \begin{proof}As in the proof of Lemma~\ref{lem: varphi extends to
      A}, we can immediately reduce to the case that~$A$ is a
    $\Z/p^a$-algebra. Part~\eqref{item: continuity of phi 2} then
    follows from Lemma~\ref{lem:extending endomorphisms} as in the
    proof of Lemma~\ref{lem: varphi extends to
      A},
while part~\eqref{item: continuity of Galois action 3}
    follows from Lemma~\ref{lem:extending group actions},     bearing in mind
    Lemmas~\ref{lem: actual examples of Cplus and C that we use}
and~\ref{lem: boundedness of Galois actions in our concrete
  settings}.\end{proof}

For our final result of this section,
we compute
the $\varphi$-invariants in $W(\C^\flat)_A$.

\begin{lem}
\label{lem:phi invariants}
If $A$ is any $p$-adically complete $\Z_p$-algebra,
then for each $a \geq 1$,
the natural morphism
$A/p^a \to W_a(\C^\flat)_A$
induces an identification
$A/p^a = \bigl(W_a(\C^\flat)_A\bigr)^{\varphi =  1} $,
and {\em (}consequently{\em )}
the natural morphism
$A \to W(\C^\flat)_A$
induces an identification
$A = \bigl(W(\C^\flat)_A\bigr)^{\varphi = 1}.$ 
\end{lem}
\begin{proof}
Since 
$W(\C^\flat)_A = \varprojlim_a W_a(\C^\flat)_A$,
we see that
$$\bigl(W(\C^\flat)_A\bigr)^{\varphi = 1} =
\varprojlim \bigl(W_a(\C^\flat)_A\bigr)^{\varphi = 1},$$
and thus (as the phrasing of the lemma indicates),
the claim for $W(\C^\flat)_A$ follows from the 
claims for each $W_a(\C^\flat)_A.$
Thus for the remainder of the proof,
we may and do assume that $A$ is a $\Z_p/p^a$-module for some $a \geq 1$,
and prove that
$A = \bigl(W_a(\C^\flat)_A\bigr)^{\varphi = 1}.$
In fact,
for any $\Z/p^a$-module $M$, 
we may define
$$W_a(\C^\flat)_M := 
\Bigl( \varprojlim_n  \bigl( W_a(\cO_\C^\flat)/T^n\bigr)\otimes_{\Z/p^a} M \Bigr)
[1/T]$$
(this extends our preceding definition for $\Z/p^a$-algebras),
and
we will prove that
$M = \bigl(W_a(\C^\flat)_M\bigr)^{\varphi = 1}.$
(In fact this statement for modules, which obviously
implies the statement for algebras, is actually equivalent to it,
as one sees by applying the statement for algebras to 
$\Z_p/p^a$-algebras of the form $\Z_p /p^a \oplus M,$
where $M$ is an arbitrary $\Z_p/p^a$-module equipped with the
square zero multiplication.)

Suppose first that $a = 1$, so that $M$  is simply an $\F_p$-vector space..  
If we choose a basis for~$M$, i.e.\ an isomorphism $\F_p^{\oplus I} \iso M$,
then $(\cO_\cC^{\flat})_M = \prod'_{i \in I} \cO_\cC^\flat,$
where the prime indicates we consider the set of tuples $(x_i)_{i\in I}$
for which $\lim_i x_i = 0,$ where the limit is taken with respect to
the filter of cofinite subsets of~$I$.
Thus we have, correspondingly, that
$\cC^{\flat}_M  = \prod'_{i \in I} \cC^\flat$
(where the prime has the same meaning).
Now $(\cC^\flat)^{\varphi = 1} = \F_p$, and so
$$(\cC^{\flat}_M)^{\varphi= 1}  = \prod_{i \in I}{}' (\cC^\flat)^{\varphi = 1}
= \oplus_{i \in I} \F_p  = M.$$
This establishes the case~$a = 1$.

Now consider the case of general~$a$.  An application of
Lemma~\ref{lem:completion properties}~(1) 
(taking $R=\Z/p^a$ and $C^+=W_a(\cO_{\C^\flat})$, and considering  the
exact sequence  of $\Z/p^a$-modules $0 \to pM\to M \to M/pM \to 0$;
see also Lemma~\ref{lem: actual examples of Cplus and C that we use})
shows that the sequence
$$0 \to  W_{a-1}(\cC^{\flat})_{p M} = W_a(\cC^\flat)_{p M} \to 
W_a(\cC^\flat)_M \to W_a(\cC^\flat)_{M/pM} = \cC^{\flat}_{M/pM} \to 0$$
is short exact.
Suppose that $x \in \bigl(W_a(\C^\flat)_M\bigr)^{\varphi = 1}.$
The case $a = 1$ already proved then shows that 
we may find $m \in M$ so that $y := x - m$  lies in
$W_{a-1}(\cC^{\flat})_{p M}.$    Of course $\varphi(y) = y$,
and so, by induction on $a$, we find that $y \in pM$.  Thus $x = m + y \in M$,
as claimed.
\end{proof}

We will also use the following variant on Lemma~\ref{lem:phi invariants}. Suppose  that ~$A$ is a complete local Noetherian $\cO$-algebra
  with finite residue field, and write $\widehat{W(\C^\flat)}_A$ for
the $\m_A$-adic completion of~$W(\C^\flat)_A$.
\begin{lem}
  \label{lem: phi invariants in W C flat with coefficients}Let~$A$ be
  a complete local Noetherian $\cO$-algebra with finite residue
  field. Then $(\widehat{W(\C^\flat)}_A)^{\varphi=1}=A$. 
\end{lem}
\begin{proof}Since 
$\widehat{W(\C^\flat)}_A = \varprojlim_n W(\C^\flat)_{A/\m_A^n}$,
we have
\[\bigl(\widehat{W(\C^\flat)}_A\bigr)^{\varphi=1} = \varprojlim_n \bigl(W(\C^\flat)_{A/\m_A^n}\bigr)^{\varphi=1},\]
so the result follows from Lemma~\ref{lem:phi invariants} (applied to each~$A/\m_A^n$).
\end{proof}
\section[Almost Galois descent for profinite group actions]{Almost Galois descent for profinite group actions
  \sectionmark{Almost Galois descent}}\sectionmark{Almost Galois descent}
\label{subsec:profinite}
We will be interested in descent results for profinite group actions,
and in this section we establish the key result that we will need.
Our setup is slightly elaborate, but accords with the situations that
will arise in practice. We begin with the case of finite group actions.

Suppose that $R$ is a ring, and that $I$ is an ideal in $R$ such
that $I^2 = I$. Assume further that~$I\otimes_R I$ is a flat $R$-module;
by~\cite[Prop.\ 2.1.7]{MR2004652}, this holds if~$I$ is a filtered
union of principal ideals. In particular, these assumptions hold
if~$R$ is a valuation ring with a non-discrete rank one valuation,
and~$I$ is the maximal ideal of~$R$.

The full subcategory of the category of
$R$-modules whose objects are the modules annihilated by $I$ is a Serre
subcategory, and so we can form the quotient category of {\em almost}
$R$-modules. \index{almost module}

Suppose that $S$ is an $R$-algebra
equipped with an action of a finite
group $G$ by $R$-algebra automorphisms
(i.e.\ the structural morphism
$R \to S$ is equivariant with
respect to the given $G$-action on $S$, and the trivial $G$-action
on $R$).

\begin{df}\index{almost Galois extension}
	We say that the morphism $R \to S$ makes $S$ an
	{\em almost Galois extension of $R$, with Galois group $G$},
	if the natural $G$-equivariant and $S$-linear morphism
$$S\otimes_R S \to \prod_{g \in G} S $$
(here $G$ acts on the source through its action on the second factor,
and on the target by permuting the factors, while $S$ acts on the source
through its action on the first factor and on the target through its
action on each factor) 
defined by $s_1\otimes s_2 \mapsto  \bigl(s_1 g(s_2) \bigr)_{g \in G}$
is an {\em almost isomorphism} (i.e.\ induces an isomorphism in the category
of almost $R$-modules).
\end{df}

\begin{remark}
	\label{rem:base-change}
	Suppose that $R\to R'$ is a morphism, and define $I' = I R'$,
	so that $I'^2 = I'$, allowing us to also define the category
	of almost $R'$-modules.  If $R \to S$ is almost Galois
	with Galois group~$G$, then evidently
	$R' \to S' := R'\otimes_R S$ is almost Galois with Galois
        group~$G$.

        Conversely, if $R\to R'$ is faithfully flat, and $R'\to S'$ is
        almost Galois with Galois group~$G$, then so is $R\to S$;
        indeed we may write $R'\otimes_R (S\otimes_RS)=(R'\otimes_R
        S)\otimes_{R'}(R'\otimes_R S)=S'\otimes_{R'}S'$.
\end{remark}


\begin{lemma}
	\label{lem:almost invariants}Suppose that $R \to S$  makes
$S$ an almost Galois extension of $S$ with Galois group $G$,
and furthermore that $R \to S$ is faithfully flat.
	Then the morphism $R \to S^G$ is an almost isomorphism.
\end{lemma}
\begin{proof}
	Since $R \to S$ is faithfully flat, it suffices to verify that
	the induced morphism 
	\numequation
        \label{eqn:hoped for a.i.}
        S \to (S\otimes_R S)^G
        \end{equation}
        is an almost isomorphism.
        By assumption the morphism $S\otimes_R S \to \prod_{g \in G} S$
        is an almost isomorphism.  One easily verifies that the induced 
        morphism on invariant subrings
        $$(S\otimes_RS)^G \to (\prod_{g \in G} S)^G = S$$
        is then also an almost isomorphism.
        The composite of this map with the morphism~\eqref{eqn:hoped for a.i.}
        is just the identity, and so we find that~\eqref{eqn:hoped for a.i.}
        is indeed an almost isomorphism.
\end{proof}

\begin{lemma}
	\label{lem:almost descent}Suppose that $R \to S$  makes
$S$ an almost Galois extension of $S$ with Galois group $G$,
and furthermore that $R \to S$ is faithfully flat.
Then if $M$ is an $S$-module equipped with a semi-linear $G$-action,
 the induced morphism $S\otimes_R M^G \to M$ is an almost
isomorphism of $S$-modules.
\end{lemma} 
\begin{proof}
	The semi-linear $G$-action on $M$ can be reinterpreted as
	an isomorphism 
	$$(\prod_{g \in G} S) \otimes_S M \cong M \otimes_S (\prod_{g \in G}
	S)$$ (where the tensor product on the left hand side is twisted
        by the $G$-action)
	of $S\otimes_R S$-modules, and hence as an almost isomorphism
	$$ S\otimes_R M \acong M\otimes_R S$$
	of $S\otimes_R S$-modules.  The claim of the lemma now follows
	from faithfully flat descent in the almost category (for which
        see~\cite[\S 3.4.1]{MR2004652}; it is in order to make this citation
	that we have assumed that $I\otimes_R I$ is $R$-flat).
\end{proof}

The following lemma explains our interest in almost Galois extensions.

\begin{lemma}
	\label{lem:almost etale to almost Galois}
	If $F \subseteq F'$ is a finite Galois extension
	of perfectoid fields, 
with Galois group $G$, then the corresponding inclusion of rings
of integers $\cO_F \subseteq \cO_{F'}$ 
realizes $\cO_{F'}$ as an almost 
Galois extension of $\cO_F$, with Galois group $G$
{\em (}the ``almost'' structure being understood with respect
to the maximal ideal of $\cO_F${\em )}.
\end{lemma}
\begin{proof}
	There are various more-or-less concrete ways to see this.
	For example, 
\cite[Prop.~5.23]{ScholzePerfectoid} shows that
the morphism $\cO_F \to \cO_{F'}$ is almost \'etale,
which by definition
\cite[Def.~4.12]{ScholzePerfectoid} means that the surjection
of $\cO_{F'}$-algebras $\cO_{F'}\otimes_{\cO_F} \cO_{F'} \to \cO_{F'}$
may be ``almost'' split, so that $\cO_{F'}$ is almost a direct factor
of the source.  Permuting such a splitting under the action
of the Galois group $G$ then allows one to show that the morphism
$\cO_{F'}\otimes_{\cO_F}\cO_{F'} \to \prod_{g \in G} \cO_{F'}$ is an almost
isomorphism.

However, a more direct (but less explicit) way to prove the lemma is to
use the first equivalence of categories in \cite[Thm.~5.2]{ScholzePerfectoid},
namely the equivalence between the category of perfectoid 
$F$-algebras and the category of perfectoid almost algebras
over~$\cO_F$.  
Under this equivalence, the morphism $\cO_{F'}\otimes_{\cO_F} \cO_{F'}
\to \prod_{g \in G} \cO_{F'}$ is taken to the morphism
${F'}\otimes_F {F'} \to \prod_{g \in G} {F'}$.  This latter morphism
{\em is} an isomorphism, since ${F'}$ is Galois over $F$ with Galois
group~$G$, and so we conclude 
that the former morphism is an almost isomorphism,
as required.
\end{proof}

In the case of perfectoid fields of characteristic~$p$,
we may extend the statement of the preceding lemma to the context
of truncated rings of Witt vectors.

\begin{lem}
	\label{lem:almost Galois Witt version}
	If $F \subseteq F'$ is a finite Galois extension
	of perfectoid fields in characteristic~$p$, 
	with Galois group~$G$,
	then for each $a \geq 1$, the inclusion
	$W_a(\cO_{F}) \hookrightarrow W_a(\cO_{F'})$
	is almost Galois,
	with Galois group $G$
{\em (}the ``almost'' structure being understood with respect
to the maximal ideal of $W_a(\cO_F)${\em )}.
\end{lem}
\begin{proof}
	For simplicity of notation, write $R := W_a(\cO_F)$
	and $S := W_a(\cO_{F'})$.  Each of $R$ and $S$ is flat
	over $\Z/p^a$, and $S$ is also flat over $R$,
	by Lemma~\ref{lem:Witt flatness}.
	Thus $S\otimes_R S$ is flat over $S$, hence over $R$,
	and hence over $\Z/p^a$ as well.

	To see that the natural morphism
	$S\otimes_R S \to \prod_{g \in G} S$
	is an almost isomorphism, it suffices to check the analogous
	condition after passing to associated graded rings for the
	$p$-adic filtration.   Lemma~\ref{lem:graded flat} allows
	us to rewrite this induced morphism on associated graded rings
	in the form
	$$\F_p[T]/(T^a) \otimes_{\F_p} \bigl((S/pS)\otimes_{(R/pR)}(S/pS)\bigr)
	\to \F_p[T]/(T^a) \otimes_{\F_p} \prod_{g \in G} (S/pS)$$
	(here $\F_p[T]/(T^a)$ appears as the associated graded ring
	to $\Z/p^a$ with its $p$-adic filtration),
	which may be identified with the base-change over $\F_p[T]/(T^a)$
	of the natural morphism
	$$(S/pS)\otimes_{(R/pR)}(S/pS) \to \prod_{g \in G} (S/pS).$$
	This latter morphism is an indeed an almost isomorphism,
	by Lemma~\ref{lem:almost etale to almost Galois}.
\end{proof}


We now pass to the profinite setting. We continue to suppose that we are given the ring $R$ endowed with
an idempotent ideal $I$.  We suppose additionally that $v \in I$
is a regular element of $R$ (i.e.\ a non-zero divisor) and 
that $R$ is $v$-adically complete.

We also suppose given a $v$-adically complete $R$-algebra $S$, 
equipped with an action of a profinite group $G$ as $R$-algebra
automorphisms.  We assume that this action is continuous, 
in the sense that the action map $G \times S \to S$ is continuous
when $S$ is endowed with its $v$-adic topology.

We further suppose that
we may write $S = \widehat{\bigcup_n S_n}$ as the $v$-adic completion
of an increasing union of $G$-invariant $v$-adically complete
$R$-subalgebras $S_n$, with $R=S_0$; that the $G$-action
on $S_n$ factors through a finite quotient $G_n :=G/H_n$ of~$G$,
where~$H_n$ is a normal open subgroup of~$G$; and that for each~$i\ge
0$ and each  $m \leq n$,  the morphism $S_m/v^i \to S_n/v^i$ is faithfully flat,
and realizes $S_n/v^i$ as an almost Galois
extension of $S_m/v^i$, having the subgroup $H_m/H_n$ of $G_n$ as Galois group.


\begin{lemma}
	\label{lem:almost invariants two}
	For any value of $m$, the morphism $S_m \to S^{H_m}$
	is an almost isomorphism, as is the induced morphism
	$S_m/v^i \to (S/v^i)^{H_m}$ for any $i \geq 0$.
\end{lemma}
\begin{proof}
	Our hypotheses, together with Lemma~\ref{lem:almost invariants},
        imply that the morphism
	$$S_m/v^i \to (S_n/v^i)^{H_m/H_n}$$
	is an almost isomorphism for each $n \geq m$ and each
	$i \geq 0.$  Passing to the direct limit over $n$
	gives the second claim,
	and then additionally passing to the inverse limit over $i$ 
	gives the first claim.
\end{proof}

\begin{cor}
	\label{cor:invariants}
	The injection $R[1/v] \hookrightarrow S[1/v]^G$ 
	is an isomorphism.
\end{cor}
\begin{proof}
	Note that $S[1/v]^G = (S^G)[1/v],$ since
	$v$ is $G$-invariant (being an element of $R$) 
	and localization is exact.
	The corollary thus follows from the $m = 0$ case of 
	Lemma~\ref{lem:almost invariants two}, as inverting $v$
	converts the almost isomorphism into a genuine isomorphism.
\end{proof}

\begin{df}
	We say that an 
       	$S$-module $M$ is {\em iso-projective of finite rank} \index{iso-projective}
	if it is finitely generated and
	$v$-torsion free, and if $S[1/v]\otimes_S M$ 
	is a projective $S[1/v]$-module.
\end{df}

\begin{lem}
	\label{lem:iso-proj implies complete}
	Any iso-projective $S$-module of finite rank is $v$-adically complete.
\end{lem}
\begin{proof}
 If~$M$ is iso-projective, then we may write $M[1/v]=eF$ where~$F$ is a finite
 free $S[1/v]$-module and $e\in\End_{S[1/v]}F$ is an
 idempotent.
 Write~$F=F_0\otimes_SS[1/v]$, where
 $F_0$ is a finite free $S$-module.
 Then, since $F_0$ is finitely generated,
 we find that $eF_0$ is contained in $v^{-a}F_0$ for some sufficiently large
 value of $a$; the latter $S$-module is $v$-adically separated, and
 thus so is the former.  Since $eF_0$ is furthermore
 the image of the $v$-adically complete $S$-module $F_0$, it is in fact
 $v$-adically complete.  Since $eF_0[1/v] = M[1/v]$, and both $eF_0$
 and $M$ are finitely generated $S$-modules,
 we find that $v^b eF_0 \subseteq M \subseteq v^{-b} e F_0$ 
 for some sufficiently large value of $b$,
 and thus $M$ is also $v$-adically
 complete, as claimed.
%
%
%
\end{proof}

\begin{thm}
	\label{thm:almost descent three}
	If $M$ is 
        an iso-projective $S$-module of finite rank,
	equipped with a semi-linear $G$-action that is continuous
	with respect to the $v$-adic topology on~$M$,
	then the kernel and cokernel 
	of the induced morphism $S \cotimes_R M^G \to M$
	{\em (}the source being the $v$-adically completed tensor
	product{\em )}
	are each annihilated by a 
	power of $v$.
\end{thm}
\begin{proof}
        Write $P := S[1/v]\otimes_S M,$ so that (by assumption) $P$ is a
	finitely presented projective module over $S[1/v]$. 
	Since $M$ is $v$-adically complete
	(by Lemma~\ref{lem:iso-proj implies complete}),
	we have an isomorphism $M \iso \varprojlim_n M/v^n M$.
	This induces a corresponding isomorphism
	$M^G \iso \varprojlim_{n} (M/v^nM)^G.$ 
	(We should also note that the topology on $M^G$ induced
	by the $v$-adic topology on $M$ coincides with the $v$-adic
	topology on $M^G$, as follows immediately from
	that fact that multiplication by $v$ on $M$ is injective
	(as $M$ is iso-projective)
	and commutes with the $G$-action.)
	Thus, to prove the theorem,
	it suffices to show that the kernel and cokernel
	of each of the morphisms
	$$S\otimes_R (M/v^n M)^G \to M/v^n M$$
	is annihilated by a power of $v$ that is bounded independently
	of $n$.

	Multiplication by $v^{-n}$ induces an isomorphism
	$M/v^n M \iso v^{-n} M/M,$ and we will actually prove the 
	equivalent statement that each of the morphisms
	\numequation
	\label{eqn:M morphisms at level n}
	S\otimes_R (v^{-n} M/ M)^G \to v^{-n} M/M
\end{equation}
	has kernel and cokernel
	annihilated by a power of $v$ that is bounded independently
	of~$n$.

	The reason for formulating our argument in terms of
	the quotients $v^{-n}M/M$ is that these may be conveniently
	be regarded as submodules of the quotient $P/M$.
        Our argument will proceed by constructing,
	for each~$n$, a $G$-invariant $S$-submodule $Z_n$ of $P/M$,
	such that $v^{-n}M/M \subseteq Z_n \subseteq v^{-(n+c)} M/M$
	(where $c$ is independent of~$n$),
	and such that each of the morphisms
	\numequation
	\label{eqn:Z morphisms at level n}
	S\otimes_R Z_n^G \to Z_n
\end{equation}
	has kernel and cokernel
	annihilated by a power of $v$ that is bounded independently
	of~$n$.  This implies the corresponding statement for the
	morphisms~(\ref{eqn:M morphisms at level n}), and thus establishes
	the theorem.
	The remainder of the argument is devoted to constructing the
	submodules $Z_n$, and proving the requisite properties
	of the morphisms~(\ref{eqn:Z morphisms at level n}).

	Since $P$ is a projective $S[1/v]$-module, we may 
	choose a finite rank free module $F$ over $S[1/v]$, and an
	idempotent $e \in \End_{S[1/v]}(F)$, such that $P = eF.$
	We choose a free $R$-submodule $F_0$ of $F$ such that
	$S[1/v]\otimes_R F_0 \iso F.$  (More concretely, $F_0$ is simply
	the $R$-span of some chosen $S[1/v]$-basis of $F$.)
	The endomorphism $e$ may not preserve the $S$-submodule
	$S\otimes_R F_0$ of $F$, but if we choose $a$ sufficiently large,
	then $e' := v^a e$ will preserve $S\otimes_R F_0.$

	Since $M$ is finitely generated over $S$, we may and do assume
	that we have chosen $F_0$ in such a manner 
	that $M \subseteq v^a(S\otimes_R F_0)$.  (Simply replace $F_0$ by $v^{-b} F_0$
	for some sufficiently large value of $b$.)
	Then in fact $M \subseteq e'(S\otimes_R F_0).$   Furthermore,
	since $M$ and $e'(S\otimes_R F_0)$ both span $P$ as an $S[1/v]$-module,
	we find that
	$v^c e'(S\otimes_R F_0) \subseteq M$ for some sufficiently large 
	value of $c$.

	For any $n \geq 0$ we have $S/v^n = \bigcup_{i\geq 0} S_i/v^n.$ 
	Thus $e' \bmod v^n$, which is an endomorphism of $(S/v^n)\otimes_R F_0$,
	descends to an endomorphism of $(S_i/v^n)\otimes_R F_0$ for
        all sufficiently large~$i$.  Replacing
	$(S_n)$ by an appropriately chosen subsequence, we may and do assume that
	in fact for each~$n$, $e'$ descends to an endomorphism $e'_n$ of
	$(S_n/v^n)\otimes_R F_0$.

	Since $G$ acts continuously on $P$ and preserves $M$,
	it acts continuously on $P/M$.  The topology on $P/M$ is
	discrete, and thus any finite subset of $P/M$ is fixed
	by some open subgroup of $G$.  In particular,
	we find that
	\[v^{-n} e'F_0 /
	\bigl( v^{-n}e' F_0 \cap M\bigr)\]
	(which we regard in the natural way as a submodule of $P/M$,
	and which we note is finitely generated over $R$)
	is pointwise fixed by some open subgroup of $G$.  Since the $H_n$
	form a cofinal sequence of open subgroups of $G$, if we again
	replace $S_n$ and $H_n$ by appropriately chosen subsequences,
	we may and do assume that in fact
	$v^{-n} e'F_0/ \bigl( v^{-n} e'F_0 \cap M\bigr)$
	is pointwise fixed by $H_n$.

	Since $M \supseteq v^c e'(S\otimes_R F_0),$ we see that
	$v^{-n} e'F_0/ \bigl( v^{-n} e'F_0 \cap M\bigr)$ is annihilated
	by $v^{n+c}$, and thus that 
	$$G \Bigl(
	v^{-n} e'F_0/ \bigl( v^{-n} e'F_0 \cap M\bigr)\Bigr) 
	\subseteq v^{-(n+c)}M/M \subseteq 
	v^{- (n+c)} e'(S\otimes_R F_0)/ M.$$
	Since $G/H_n$ is finite, and since $H_n$ fixes $v^{-n} e'F_0/ \bigl( v^{-n} e'F_0 \cap M\bigr)$
 pointwise, we see that
	in fact
	$$G \Bigl(
	v^{-n} e'F_0/ \bigl( v^{-n} e'F_0 \cap M\bigr)\Bigr) 
	\subseteq
	v^{- (n+c)} e'(T\otimes_R F_0)/ M,$$
	for some finitely generated $R$-subalgebra $T$ of $S/v^{n+c} S.$
	Passing to a subsequence one more time, we may assume
	that $T$ is contained in $S_n$.


	We let $Y_n$, respectively $Z_n$,
       	denote the $S_n$-submodule, respectively the $S$-submodule, of $P/M$ 
	generated by the $G$-translates of 
	$v^{-n} e'F_0/ \bigl( v^{-n} e'F_0 \cap M\bigr)$.
	It remains to prove the requisite properties of $Z_n$.
	We begin by noting the inclusions
	$$v^{-n} e'(S\otimes_R F_0) / M \subseteq Z_n \subseteq 
	v^{-(n+c)}M/M \subseteq v^{-(n+c)} e' (S \otimes_R F_0)/M,$$
	which imply that
	$$v^c Z_n \subseteq v^{-n} e'(S\otimes_R F_0)/M \subseteq Z_n,$$
	and thus that 
	\numequation
	\label{eqn:chain of inclusions}
	v^c Z_n^{H_n} \subseteq
	(v^{-n} e'(S\otimes_R F_0)/M)^{H_n} \subseteq Z_n^{H_n}.
        \end{equation}

	Lemma~\ref{lem:almost invariants two} shows that
	$S_n/v^n  \to (S/v^n)^{H_n}$ is an almost isomorphism.
	Thus $(v^{-n} S_n/S_n) \otimes_R F_0 \to
	(v^{-n} S/S \otimes_R F_0)^{H_n}$ is an almost
        isomorphism, if we declare that $H_n$ acts trivially on $F_0$,
	and thus via its action on the first factor in the 
	target tensor product.
	Since $(e')^2 = v^a e',$ we find that the cokernel
	of the inclusion
	$$
	e' \bigl( (v^{-n} S\otimes_R F_0)/ (S\otimes_R F_0) \bigr)^{H_n}
	\hookrightarrow
	\bigl(e' ( v^{-n} S\otimes_R F_0)/ (S\otimes_R F_0) \bigr)^{H_n}
	$$
	is annihilated by $v^a$, and thus that the natural morphism
	$$
	e' (v^{-n} S_n\otimes_R F_0)/ (S_n\otimes_R F_0) 
	\to 
	\bigl(e' ( v^{-n} S\otimes_R F_0)/ (S\otimes_R F_0) \bigr)^{H_n}
	$$
	has its kernel annihilated by $I$, and its cokernel 
	annihilated by $v^a I.$

	If we let $X_n$ denote the image of $e'(v^{-n} S_n\otimes_R F_0)$
	in 
	$\bigl(e' ( v^{-n} S\otimes_R F_0)/ M \bigr)^{H_n},$
	then certainly
	$$X_n \subseteq Y_n \subseteq Z_n^{H_n}.$$
	It follows from the conclusion of the preceding
	paragraph, along with the fact that $v^c e'(S_0\otimes_R F_0)
	\subseteq M,$
        that the cokernel of the inclusion
	$$X_n \hookrightarrow
	\bigl(e' ( v^{-n} S\otimes_R F_0)/ M \bigr)^{H_n}$$
	is annihilated by $v^{a+c} I,$
	while the chain of inclusions~(\ref{eqn:chain of inclusions})
	shows that the cokernel of the inclusion
	$$
	\bigl(e' ( v^{-n} S\otimes_R F_0)/ M \bigr)^{H_n}
	\subseteq Z_n^{H_n} 
	$$
	is annihilated by $v^c$.  
	Thus the cokernel of the inclusion $X_n \subseteq Z_n^{H_n}$
	is annihilated by $v^{a+2c} I,$ and hence so is the cokernel
	of the inclusion $Y_n \subseteq Z_n^{H_n}$.
	Passing to $G_n$-invariants in the inclusion just mentioned,
	we find that there is an inclusion
	$Y_n^{G_n}\subseteq Z_n^G$,
	whose cokernel is annihilated by $v^{a+2c}I$.
	Extending scalars to $S_n$, we obtain a morphism
	$$S_n\otimes_R Y_n^{G_n} \to S_n\otimes_R Z_n^{G}$$
	(which is in fact an embedding,
	since $S_n/v^{n+c}$ is flat over $R/v^{n+c}$,
	although we don't need this here),
	whose cokernel is annihilated by $v^{a+2c}I$.
	On the other hand,
        Lemma~\ref{lem:almost descent} implies that the natural morphism
	$S_n\otimes_R Y_n^{G_n} \to Y_n$ is an almost isomorphism.

	Putting all these results together,
	we find that the kernel and cokernel of the natural morphism
	$$S_n\otimes_R Z_n^{G} \to Z_n^{H_n}$$
	are each annihilated by $v^{a+2c} I$.
	Indeed, the composite 
        $S_n\otimes_R Y_n^{G_n} \to S_n\otimes_R Z_n^{G} \to
        Z_n^{H_n}$ factors through the almost isomorphism
        $S_n\otimes_R Y_n^{G_n} \to Y_n$, 
	and we have shown that the natural morphism $Y_n \to Z_n^{H_n}$
	is an injection whose cokernel is killed by $v^{a+2c}I$,
	while the cokernel of
        $S_n\otimes_R Y_n^{G_n} \to S_n\otimes_R Z_n^{G}$ is also
        killed by~$v^{a+2c}I$.


	Since $Y_n \subseteq Z_n^{H_n},$ and since $Y_n$ generates $Z_n$
	over $S$ by their very definitions, we see that the natural map
	\numequation
	\label{eqn:natural surjection}
	S\otimes_{S_n}  Z_n^{H_n} \to Z_n
        \end{equation}
	is surjective.
	We bound the exponent of its kernel as follows:
	The inclusion $X_n \subseteq Z_n^{H_n}$,
        whose cokernel we have shown above to be annihilated by $v^{a+2c}I$,
	induces an inclusion 
	$$S\otimes_{S_n} X_n \subseteq S\otimes_{S_n} Z_n^{H_n},$$
	whose cokernel is again annihilated by $v^{a+2c}I,$
	while the defining surjection
	$$e'\bigl((v^{-n} S_n/v^cS_n)\otimes_R F_0\bigr) \to
	X_n$$
	induces a surjection
	$$e'\bigl((v^{-n} S\otimes_R F_0)/(v^c S\otimes_R F_0)\bigr)
	\to S\otimes_{S_n} X_n .$$
	Now the natural morphism
	\numequation\label{eqn:finalsurjection}e'\bigl((v^{-n}S\otimes_R F_0)/(v^c S\otimes_R F_0)\bigr) \to
	M[1/v]/M\end{equation}
	has kernel annihilated by $v^c$,
	and so we find that the kernel of~(\ref{eqn:natural surjection})
	is annihilated by $v^{a+3c} I$. (Indeed, if $x$ is an element
        of the kernel of~(\ref{eqn:natural surjection}), then for any
        $i\in I$ we can lift $v^{a+2c}ix$ to an element of the kernel
        of~(\ref{eqn:finalsurjection}).)

	Putting together the results of the preceding two paragraphs,
	we find that the cokernel of the natural morphism
	$$S\otimes_R Z_n^G \to Z_n$$ is annihilated by $v^{a+2c} I$,
	while its kernel is annihilated by $v^{2a+5 c} I$.
	Recalling that $v\in I$, and noting that 
	the powers of $v$ just mentioned are independent 
	of $n$, and also that we have the inclusions
	$v^{-n} M /M \subseteq Z_n \subseteq v^{-(n+c)}M/M,$
	we see that the proof of the theorem is completed.
\end{proof}

We will now deduce a descent result for projective modules 
over $S[1/v]$.  For this, we need to make some additional hypotheses,
which we now describe.

\begin{hyp}
	\label{hyp:descent}\leavevmode
	\begin{enumerate} 
	\item $S/vS$ is countable.
	\item $S[1/v]$ is Noetherian.
	\item The morphism $R[1/v] \to S[1/v]$ is faithfully flat.
	\end{enumerate}
\end{hyp}

We write $S[1/v] = \varinjlim_n {v^{-n}} S.$
If we identify $v^{-n}S$ with $S$ via multiplication by~
$v^n$, then each of the transition maps becomes identified with
the closed embedding $vS \hookrightarrow S$.
Thus if we equip $S[1/v]$ with the inductive limit topology,
then $S[1/v]$ becomes a topological ring, which is completely
metrizable (since $S$ is $v$-adically complete) and in fact Polish
(since $S/vS$ is countable, so that $S$ and thus $S[1/v]$ is separable).
Note also that $vS$ is an open additive subgroup of $S[1/v]$
which is closed under multiplication, and consists of topologically
nilpotent elements.
Proposition~\ref{prop:module topologies} then
shows that finitely generated $S[1/v]$-modules
have a canonical topology, with respect to which all $S[1/v]$-homomorphisms
are continuous, with closed image. 

We now establish the following descent result in this situation.

\begin{theorem}
	\label{thm:descent with coefficients}
	If Hypothesis~{\em \ref{hyp:descent}} holds,
	and if
	$M$ is a finitely generated projective $S[1/v]$-module,
	equipped with a continuous
        {\em (}when endowed with its canonical topology{\em )}
        semi-linear $G$-action,
	then $M^G$ is a finitely generated and projective $R[1/v]$-module,
	and the natural morphism
	\numequation
	\label{eqn:natural morphism}
	S[1/v]\otimes_{R[1/v]} M^G \to M
\end{equation}
       	is an isomorphism.
\end{theorem}


\begin{proof}


	Let $M'$
		denote any finitely generated $S$-submodule of $M$ that generates
		$M$ over $S[1/v]$; then the $S$-span $M''$ of $G M'$ 
		is finitely generated (note that since~$M'$ is finitely
                generated, there is an open subgroup~$H$ of~$G$ such
                that $HM'\subseteq M'$, and~$H$ has finite index in the
                profinite group~$G$), so it is an iso-projective $S$-submodule of $M$ which
		is $G$-invariant. 
		Applying Theorem~\ref{thm:almost descent three} to $M''$,
		and noting that $$M^G = (M''[1/v])^G = (M'')^G[1/v]$$
		(because localisation is exact),
		we find
		that the natural morphism~(\ref{eqn:natural morphism})
		has dense image.  Since its target is finitely generated
		over $S[1/v]$, its image is also finitely generated
		(because $S[1/v]$ is Noetherian, by assumption),
		and thus closed in its target;
		combined with the density, we find 
		that~(\ref{eqn:natural morphism}) is surjective.

		Since $S[1/v]\otimes_{R[1/v]} M^G \to M$ is surjective,
		while $M$ is finitely generated over $S[1/v]$,
		we see that if $N$ is any sufficiently large finitely generated 
		$R[1/v]$-submodule of $M^G$, then 
		$S[1/v]\otimes_{R[1/v]}N \to M$ is surjective.
		Choose a finitely generated free $R[1/v]$-module $F$
		that surjects onto $N$, and let $E$ denote the kernel
		of the induced surjection
		$$ S[1/v]\otimes_{R[1/v]} F \to S[1/v]\otimes_{R[1/v]} N
		\to M,$$
		so that we have a short exact sequence
		$$0 \to E \to S[1/v]\otimes_{R[1/v]} F \to M \to 0.$$
		Passing to $G$-invariants,
	        and taking into account Corollary~\ref{cor:invariants},
		we obtain a left exact sequence
		$$0 \to E^G \to F \to M^G,$$
		which (by the choice of $F$)
		induces a short exact sequence
		$$0 \to E^G \to F \to N \to 0.$$
		Tensoring back up with $S[1/v]$ (which is
		flat over $R[1/v]$ by assumption),
		we obtain a morphism of short exact sequences
		$$\xymatrix{ 0 \ar[r] &
			S[1/v]\otimes_{R[1/v]} E^G \ar[r]\ar[d] &
			S[1/v]\otimes_{R[1/v]} F \ar[r]\ar[d] &
			S[1/v]\otimes_{R[1/v]} N \ar[r] \ar[d] & 0 \\
			0 \ar[r] & E \ar[r] & S[1/v]\otimes_{R[1/v]} F \ar[r] &
			M \ar[r] & 0 }
		$$
		Evidently the middle vertical arrow is the identity,
		and thus the natural morphism
		$S[1/v]\otimes_{R[1/v]} E^G \to E$ 
		is injective. Since~$E$ is finitely generated and projective
		(being the kernel of a surjection from a finitely generated
		free module to a projective module), it follows from
                what we have already proved that this morphism is also
                surjective.

                Since the left hand two vertical arrows are
                isomorphisms, so is the third, so that $S[1/v]\otimes_{R[1/v]} N \to M$ is an
                isomorphism.  This is true for any sufficiently large
                choice of $N$, and thus we find (using the faithful
                flatness of $S[1/v]$ over $R[1/v]$, which holds by
                assumption) that all these sufficiently large choices
                of~$N$ coincide, implying that $M^G$ is finitely
                generated and that~(\ref{eqn:natural morphism}) is an
                isomorphism.  The faithful flatness of $S[1/v]$ over
                $R[1/v]$ then implies that $M^G$ is projective over
                $R[1/v]$
                \cite[\href{http://stacks.math.columbia.edu/tag/058S}{Tag
                  058S}]{stacks-project}, as
                required.
%
\end{proof}

The following corollary provides a convenient reformulation
of the preceding theorem.

\begin{cor}
	\label{cor:descent with coefficients}
	If Hypothesis~{\em \ref{hyp:descent}} holds,
	then the functor \[M \mapsto S[1/v]\otimes_{R[1/v]} M\]
	induces an equivalence between the category
	of finitely generated projective $R[1/v]$-modules
	and the category of finitely generated projective
	$S[1/v]$-modules endowed with a continuous semi-linear
	action of $G$.  A quasi-inverse functor is given by
	$N \mapsto N^G$.
\end{cor}
\begin{proof}
	Theorem~\ref{thm:descent with coefficients}
	shows that if $N$ is a finitely generated and projective
	$S[1/v]$-module, endowed with a continuous 
	semi-linear $G$-action, then the natural map
	$S[1/v]\otimes_{R[1/v]} N^G \to N$ is an isomorphism.
	To complete the proof of the corollary, then, it suffices
	to show that if $M$ is a finitely generated and projective 
	$R[1/v]$-module, then the natural map
	$M\to (S[1/v]\otimes_{R[1/v]}M)^G$ is an isomorphism.
	Writing $M$ as the direct summand of a finite
	rank free module, we reduce to the case when $M$ is free,
	which (as was already observed in the proof of
	Theorem~\ref{thm:descent with coefficients})
	is established by Corollary~\ref{cor:invariants}.
\end{proof}

\section{An application of almost Galois descent}\label{subsec:coefficients}
Let $F$ be a closed perfectoid subfield of $\C$, 
with tilt~$F^\flat$, a closed perfectoid subfield of $\C^\flat$.  
Recall that, for any $p$-adically complete $\Z_p$-algebra~$A$,
we defined $W(\cO_F^\flat)_A$ and $W(F^\flat)_A$
in Section~\ref{subsec: coefficients}.

\begin{theorem}
	\label{thm:descending projective modules}
			Let $A$ be a finite type $\Z/p^a$-algebra,
			for some $a \geq 1$.
			The inclusion
			$W(F^\flat)_A  \to W(\C^\flat)_A$
			is a faithfully flat morphism
			of Noetherian rings,
			and
	the functor $M \mapsto W(\C^\flat)_A\otimes_{W(F^\flat)_A} M$
	induces an equivalence between the category of finitely generated 
	projective $W(F^\flat)_A$-modules and the category of
	finitely generated projective $W(\C^\flat)_A$-modules endowed
	with a continuous semi-linear $G_F$-action. A quasi-inverse
        functor is given by $N\mapsto N^{G_F}$.
\end{theorem}

The proof will be an application of the almost Galois descent results
of Section~\ref{subsec:profinite}.  Thus we have to place ourselves
in the framework of that section.  To this end,
we note that $\C$ may be regarded as a completion $\widehat{\overline{F}}$ 
of the algebraic closure $\overline{F}$ of $F$.  We then write
$\overline{F} = \bigcup_{n\geq 1} F_n$ as the increasing union of a sequence
of finite Galois extensions of $F$; to match the notation
of Section~\ref{subsec:profinite}, we also write $F_0~:=~F$.
We write $G_n := \Gal(F_n/F),$ $H_n := \Gal(\overline{F}/F_n)$,
and $G := H_0 = \Gal(\overline{F}/F),$ so that $G \iso \varprojlim_n G_n$
and $G/H_n \iso G_n.$
Note that each of the finite Galois extensions $F_n$ of $F$
is again perfectoid. 

%

Assume now that $A$ is a finitely generated $\Z/p^a$-algebra.
Then
$W(F^\flat)_A  =  W_a(F^\flat)_A,$
and similarly
$ W(\C^\flat)_A =  W_a(\C^\flat)_A$,
so that from here on we may work with rings of $a$-truncated Witt vectors,
rather than with full rings of Witt vectors themselves.
To accord with the notation of Section~\ref{subsec:profinite},
we also denote these rings by $R$ and $S$ respectively.
Recall from Section~\ref{subsec: coefficients}
that $v$ denotes a non-zero element of the maximal ideal
of $W_a(\cO_F^\flat)$,
and that by definition
$$R := W_a(\cO_F^\flat)_A 
= \varprojlim_i \bigl( W_a(\cO_F^\flat)\otimes_{\Z/p^a}A\bigr)/v^i,$$
while
$$S := W_a(\cO_\C^\flat)_A
= \varprojlim_i \bigl(W_a(\cO_\C^\flat)\otimes_{\Z/p^a}A\bigr)/v^i;$$
so that both $R$ and $S$ are $v$-adically complete.


We let $I$ denote the ideal in $R$ generated by the maximal ideal
of $W_a(\cO_F^\flat),$ and let $v$ denote the image in $I$ of the
chosen element $v$ of that maximal ideal (no confusion should result
from this duplication of notation).
We further set
$$S_n := W_a(\cO_{F_n}^{\flat})_A
= \varprojlim_i \bigl(W_a(\cO_{F_n}^\flat)\otimes_{\Z/p^a}A \bigr)/v^i$$
for each $n~\geq~0$
(so in particular $S_0 = R$). 
By construction we have that
$S$ coincides with the $v$-adic completion of $\varinjlim_n S_n$.
As we will see below, the transition morphisms $S_m \to S_n$
are in fact faithfully flat, and thus injective, and so
in fact this direct limit is simply a union;
thus $S = \widehat{\bigcup_n S_n},$
as required for the setup of Section~\ref{subsec:profinite}.
Furthermore,
Lemma~\ref{lem: varphi extends to A and Galois etc}
ensures that the $G$-action on $S$ is continuous.



\begin{prop}
	\label{prop:confirming set-up}
	For each $i\geq 0,$ and each $m \leq n$,
	the morphism $S_m/v^i \to S_n/v^i$ is faithfully flat, and realizes
	$S_n/v^i$ as an almost Galois extension of $S_m/v^i$, with Galois group
	$H_m/H_n$.
\end{prop}
\begin{proof}
Lemma~\ref{lem:Witt flatness}
shows that each of the inclusions
$W_a(\cO_{F_m}^{\flat}) \hookrightarrow W_a(\cO_{F_n}^{\flat})$ 
is faithfully flat,
	thus so is the morphism
	$W_a(\cO_{F_m}^\flat)\otimes_{\Z/p^a}A
	\to W_a(\cO_{F_n}^\flat)\otimes_{\Z/p^a} A$,\
	and hence so are each
	of the morphisms
	\begin{multline*}
	S_m/v^i S_m =
	\bigl(W_a(\cO_{F_m}^\flat)\otimes_{\Z/p^a}A)\bigr)
	/ v^i
	\bigl(W_a(\cO_{F_m}^\flat)\otimes_{\Z/p^a}A)\bigr)
	\\
	\to
	\bigl(W_a(\cO_{F_n}^\flat)\otimes_{\Z/p^a}A)\bigr)
	/v^i
	\bigl(W_a(\cO_{F_n}^\flat)\otimes_{\Z/p^a}A)\bigr)
	=S_n/v^i S_n .
\end{multline*}
An identical argument, taking into account
Lemma~\ref{lem:almost Galois Witt version}
and Remark~\ref{rem:base-change},
shows that $S_m/v^i S_m \to S_n/v^i S_n$ is almost Galois
(with respect to the ideal in $S_m$ generated by
the maximal ideal of $W_a(\cO_{F_m}^{\flat})$,
and hence also with respect to $I S_m$, since the latter
ideal is contained
in the former),
with Galois group $H_m/H_n$.
\end{proof}

\begin{prop}
	\label{prop:rigid facts}\leavevmode
\begin{enumerate}
\item Each $S_m[1/v]$, as well as $S[1/v]$, is Noetherian.
	{\em (}Setting $m = 0$ gives in particular that $R[1/v]$ is
	Noetherian.{\em )}
\item Each of the morphisms $S_m \to S_n$ {\em (}for $m \leq n${\em )},
       as well as each morphism $S_m \to S$, is faithfully flat.
       In particular {\em (}setting $m = 0$ and then inverting~$v${\em )}
       the morphism $R[1/v]\to S[1/v]$ is faithfully flat.
\end{enumerate}
\end{prop}
\begin{proof}
	The Noetherian claims of~(1) follow 
	from~Proposition~\ref{prop: rigid analysis FGK}~(1).
	
	We already noted
	in Proposition~\ref{prop:confirming set-up}
	that the morphisms $S_m/v^i \to S_n/v^i$ are faithfully flat,
	and an identical argument shows that
	each morphism $S_m/v^i \to S/v^i$ is faithfully flat.
	Taking into account that statement of~(1),
	we find that the claims of~(2) follow from Proposition~\ref{prop:
		FGK flat}.
	%
\end{proof}

\begin{proof}[Proof of
	Theorem~{\em \ref{thm:descending projective modules}}]
	Proposition~\ref{prop:confirming set-up} verifies
	that the running assumptions imposed at the beginning
	of Section~\ref{subsec:profinite} are satisfied.
	Proposition~\ref{prop:rigid facts} verifies
	that Hypothesis~\ref{hyp:descent} is satisfied,
	and also establishes the Noetherian and faithful
        flatness claims	of
	Theorem~\ref{thm:descending projective modules}.
	The remainder of the theorem then follows from
	Corollary~\ref{cor:descent with coefficients}.
\end{proof}

\section{\'Etale \texorpdfstring{$\varphi$}{phi}-modules}
\label{subsec: EG stuff}
In this section
we briefly recall and generalize some definitions and results
from~
~\cite{EGstacktheoreticimages}. Let~$R$ be a $\Zp$-algebra, equipped
  with a  ring endomorphism~$\varphi$, which is congruent to the
  ($p$-power) Frobenius modulo~$p$. If~$M$ is an $R$-module, we write \[\varphi^*M:=R\otimes_{R,\varphi}M.\]  
\begin{defn}\index{\'etale $\varphi$-module}
  \label{defn: general phi module} An \emph{\'etale
    $\varphi$-module over~$R$} is a finite $R$-module~$M$, equipped
  with a $\varphi$-semi-linear endomorphism $\varphi_M:M\to M$, which
  has the property that the induced $R$-linear morphism
  \[\Phi_M:\varphi^*M\stackrel{1\otimes\varphi_M}{\longrightarrow}M \]
  is an isomorphism.  A morphism of
   \'etale $\varphi$-modules is a morphism of the underlying
  $R$-modules which commutes with the
  morphisms~$\Phi_{M}$. We say that~$M$ is \emph{projective} (resp.\
  \emph{free}) if it is projective of constant rank (resp.\ free of
  constant rank) as an $R$-module. 
\end{defn}
We will typically apply this definition with~$R$ taken to be one of
the coefficient rings defined in Section~\ref{subsec:
  coefficients}. Of particular interest to us will be the cases
corresponding to imperfect fields of norms, as it is these cases which
fit into the framework of~\cite[\S5]{EGstacktheoreticimages}, and we
will use the results of that paper to 
prove the basic algebraicity properties of our moduli
stacks.

\begin{defn}\label{defn: base change of phi modules}
  Let~$S$ be an $R$-algebra, and let~$\varphi_S$ be a ring
  endomorphism of~$S$, which is congruent to Frobenius modulo~$p$, and
  is compatible with~$\varphi$ on~$R$. Then if~$M$ is an \'etale
  $\varphi$-module over~$R$, the extension of scalars $S\otimes_RM$ is
  naturally an \'etale $\varphi$-module over~$S$, with
  $\varphi_{S\otimes_R M}:=\varphi_S\otimes\varphi_M$.
\end{defn}

\subsection{Multilinear algebra}
We briefly recall the multilinear algebra
of projective \'etale $\varphi$-modules.
Firstly, if
$P$ is a projective \'etale $\varphi$-module over~$R$, then we give
its $R$-dual
$P^\vee := \Hom_{R}(P,R)$ the structure of an
\'etale $\varphi$-module by defining the isomorphism $\varphi^*P^\vee\to
P^\vee$ to be the inverse of the transpose of~$\Phi_P$.
Secondly, if~$M$ and~$N$
are projective \'etale $\varphi$-modules, then we endow $M\otimes_RN$ 
with the structure of an \'etale $\varphi$-module by
defining $\varphi_{M\otimes N}:=\varphi_M\otimes\varphi_N$. 

\begin{lem}
  \label{lem: projective phi module duality}If $M,$
  $N$ are projective \'etale $\varphi$-modules over~$R$ then we have
a natural identification 
$\Hom_{R,\varphi}(M,N)=(M^\vee\otimes_{R}
N)^{\varphi=1}$.
\end{lem}
\begin{proof}We have
  $\Hom_{R}(M,N)=M^\vee\otimes_{R} N$. Given~$f\in \Hom_{R}(M,N)$,
  regarded as an element of~$P:=M^\vee\otimes_{R} N$, we have $f\in
  P^{\varphi=1}$ if and only if $\Phi_P(1\otimes f)=f$, and by the
  definition of the $\varphi$-structure on~$P$, this is equivalent
  to~$f$ intertwining~$\Phi_M$ and~$\Phi_N$, as required.
\end{proof}

\section{Frobenius descent}\label{subsec: perfect
  etale phi modules}
Suppose that~$\A$ is as in Situation~\ref{subsubsec:general framework},
and that we have a continuous $\varphi$-equivariant embedding 
$\A\into W(\C^\flat)$. 
Let~$A$ be a finite type $\Z/p^a$-algebra for some $a\ge
1$. Assume that the induced map $\A_A\to
  W(\C^\flat)_A$ is a faithfully flat injection. 
  Since~$\varphi$ is
bijective on~$W(\C^\flat)$, we have an increasing
union 
\[\A_A\subset\varphi^{-1}(\A_A)\subset\varphi^{-2}(\A_A)\subset\dots\subset
  W(\C^\flat)_A.\]
We let~$\tA_A$ be the closure of 
$\cup_{n\ge 0}\varphi^{-n}(\A_A)$ in~$W(\C^\flat)_A$,
and we
set~$\tA^+_A:=\tA_A\cap\AAinf{A}$. 
  Note that~$\varphi$ extends to a bijection
  on~$\tA_A$, which induces a bijection on~$\tA^+_A$.

  \begin{rem}
    \label{rem: abstract perfection section linked to our concrete
      settings} We will apply the results of this section in the
    setting introduced in Section~\ref{subsec:rings}, taking~$\A=\A_K$
    or~$\A=\gS$. In either case it follows from Proposition~\ref{prop: maps of coefficient rings are faithfully flat
    injections} that $\A_A\to
  W(\C^\flat)_A$ is a faithfully flat injection, and it follows easily from Theorem~\ref{thm:field of
      norms} that in the former case we have~$\tA_A=\tA_{K,A}$, and in
    the latter case we have~$\tA_A=\tOEA$.
  \end{rem}

  \begin{rem}
    \label{rem: warning about A plus not phi stable implies not in A
      tilde plus}If~$\A^+_A$ is not $\varphi$-stable, then $\A^+_A$ is
    not a subring of ~$\tA^+_A$, and in particular it is not equal to
    $\A_A\cap\tA^+_A$. 
  \end{rem}

We will make use of the following variant of Lemma~\ref{lem: phi on powers of T without assuming A plus
      stable}; in view of Remark~\ref{rem: warning about A plus not phi stable implies not in A
      tilde plus} we cannot apply Lemma~\ref{lem: phi on powers of T without assuming A plus
      stable} in the setting of $\tA^+_A$-modules.

  \begin{lem}
    \label{lem: phi on powers of T for A tilde plus}There is an
    integer~$r\ge 0$ such that for all~$s\ge r$ we have
    $T^s\in\tA^+_A$, and an integer~$C\ge 0$ such that for
    all~$n\ge 1$, we have $\varphi(T^n)\in T^{pn-C}\tA^+_A$.
  \end{lem}
  \begin{proof}Since~$T$ is topologically nilpotent in~$\A_A^+$ it is
    topologically nilpotent in~$\tA_A$, so in particular
    $T^s\in\tA^+_A$ for all~$s$ sufficiently large. For the second
    statement, we follow the proof of Lemma~\ref{lem: phi on powers of
      T without assuming A plus stable}, and write
    $\varphi(T)=T^p+pY$ for some~$Y\in\A_A$, so that~$\varphi(T^n)=(T^p+pY)^n$. If we
    expand using the binomial theorem, and recall that~$p^a=0$, then
    we see that $\varphi(T^n)-T^{np}$ is a sum of terms of the form 
    $T^{p(n-r)}Y^r$ for some $0\le r \le a-1$. It therefore suffices
    to  choose~$C$ such that $T^{C}T^{-pr}Y^r\in \tA^+_A$ for $0\le r
    \le a-1$, which we can do by the topological nilpotence of~$T$.
  \end{proof}

\begin{rem}
\label{rem:A-tilde as a localization of  A-tilde-plus}
As already noted in Remark~\ref{rem: warning about A plus not phi stable implies not in A tilde plus}, $\tA^+_A$  need not contain~$\A^+_A$, and so, although $\A_A
= \A_A^+[1/T]$, it doesn't make sense to write $\tA_A =  \tA^+_A[1/T]$.
On the other hand, Lemma~\ref{lem: phi on powers of T for A tilde plus} does
ensure that $T^r \in \tA^+_A$ for some $r > 0$, and then it does
make sense to write, and is true (since $T$ is
topologically nilpotent in $\tA_A$, while $\tA^+_A$ is open 
in $\tA_A$), that $\tA_A = \tA^+_A[1/T^r].$
\end{rem}

Given an \'etale $\varphi$-module $M$ over~$\A_A$, we may form its
\emph{perfection} $\widetilde{M}:=\tA_A\otimes_{\A_A} M$,
\index{perfection (of \'etale $\varphi$-module)}
which is an \'etale $\varphi$-module over ~$\tA_A$.

\subsection{Descending morphisms}


\begin{prop}
  \label{prop: descent of morphisms for phi modules}Let~$M, N$ be
  projective \'etale $\varphi$-modules. Then the natural map
  $\Hom_{\A_A,\varphi}(M,N)\to\Hom_{\tA_A,\varphi}(\widetilde{M},\widetilde{N})$
  is a bijection.
\end{prop}
\begin{proof}
  Write $P:=M^\vee\otimes N$, so that~$P$ is a projective
  \'etale $\varphi$-module over~$\A_A$.  Then by Lemma~\ref{lem: projective phi module duality}, 
  we are reduced to checking that the natural morphism of $A$-modules
  $P^{\varphi=1}\to\widetilde{P}^{\varphi=1}$ is an isomorphism. It
  is certainly injective, so we need only show that it is
  surjective.

Since the formation of $\varphi$-invariants is compatible with direct
sums, it follows from~ \cite[Lem.\ 5.2.14]{EGstacktheoreticimages} that it is enough to consider the case that~$P$ is in fact
free. By Lemma~\ref{lem: phi on powers of T for A tilde plus}
(\emph{cf.}\ Remark~\ref{rem:A-tilde as a localization of  A-tilde-plus})
we can choose a $\varphi$-stable free
$\tA^+_A$-submodule~$\widetilde{\gP}$ of~$\widetilde{P}$, which generates $\widetilde{P}$ over $\A_A$
(indeed, choose any basis of~$\widetilde{P}$, multiply by a sufficiently large
power of~$T$, and let~$\widetilde{\gP}$ be the submodule
generated by this scaled basis).

Consider an element $x\in \widetilde{P}^{\varphi=1}$. Choose an
integer $r\ge 1$ such that if~$s\ge r$ then $\varphi(T^s\tA^+_A)\subseteq T^{s+1}\tA^+_A$ (such
an~$r$ exists by Lemma~\ref{lem: phi on powers of T for A tilde plus}). We may write $x=x_1+x_2$
where $x_2\in T^r\widetilde{\gP}$ and $x_1\in
\varphi^{-n}(\A_A)\otimes_{\A_A}P$ for some $n\ge 0$. Choose~$n$ to be
minimal with this property. If $n>0$, then since
$x=\varphi(x)=\varphi(x_1)+\varphi(x_2)$, and $\varphi(x_2)\in T^r\widetilde{\gP}$,
we see that we may replace~$x_1$ with $\varphi(x_1)\in\varphi^{1-n}(\A_A)\otimes_{\A_A}P$, a contradiction.

Thus $x_1\in P$, and we have $x_1-\varphi(x_1)=\varphi(x_2)-x_2\in
T^r\widetilde{\gP}\cap P$. 
The sum
$x':=x_1+\sum_{i=0}^\infty\varphi^i(\varphi(x_1)-x_1)$ converges to an
element of~$T^r\widetilde{\gP}$, and the sum therefore converges
in~$P$. By definition we have $\varphi(x')=x'$. We have $x-x_1$,
$x'-x_1\in T^r\widetilde{\gP}$, so $x-x'\in T^r\widetilde{\gP}$; then
$x-x'=\varphi(x-x')\in T^{r+1}\widetilde{\gP}$, and iterating gives
$x=x'\in P$, as required.
\end{proof}

\subsection{Descending objects}


We will show that every projective  \'etale $\varphi$-module over~$\tA_A$
arises as the perfection of an \'etale $\varphi$-module over~$\A_A$.
We first note the following analogue of~\cite[Lem.\ 5.2.14]{EGstacktheoreticimages}
for \'etale $\varphi$-modules over~$\tA_A$.

\begin{lemma}\label{lem: projective perfect phi module is summand of free}
If~$M$ is a projective  \'etale $\varphi$-module
  over~$\tA_A$, then~$M$ is a direct summand of a free 
  \'etale $\varphi$-module
  over~$\tA_A$. 
  \end{lemma}
  \begin{proof}
	  This can be proved in an identical fashion to~\cite[Lem.\ 5.2.14]{EGstacktheoreticimages}.
  \end{proof}

\begin{prop}
  \label{prop: descending projective modules from infinite level}Let
  $M$ be a projective  \'etale $\varphi$-module over~$\tA_A$. Then there is a projective \'etale
  $\varphi$-module $M_0$ over~$\A_A$ and an isomorphism
  $M\isoto \widetilde{M}_0$.
\end{prop}
\begin{proof}Suppose firstly that $M$ is  free of some rank~$d$ as an
  $\tA_A$-module. Let~$X$ denote the matrix
  of~$\varphi$ with respect to a choice of basis
  $e_1,\dots,e_d$. As in the proof of Proposition~\ref{prop: descent of morphisms for phi modules}, after possibly
  scaling the~$e_i$ by powers of~$T$, we may assume
  that~$X$ has entries in~$\tA^+_A$. Since $\Phi_{M}$
  is an isomorphism, it follows from Remark~\ref{rem:A-tilde as a localization of  A-tilde-plus})
that there is an integer $h\ge 0$ such that
$X^{-1} \in T^{-h} M_d(\tA^+_A).$

  If we change basis via a matrix~$Y\in M_d(\tA^+_A)\cap\GL_d(\tA_A)$, the new matrix for~$\varphi$ is
  $\varphi(Y)XY^{-1}$. It suffices to show that we can choose~$Y$ so
  that this matrix has entries in~$\A_A$ (as we can then let~$M_0$ be
  the $\A_A$-span of this new basis). 

  Fix some~$H>\lceil (C+h+1)/(p-1)\rceil$, where~$C$ is as in
  Lemma~\ref{lem: phi on powers of T for A tilde plus}, and write $X=X'+X''$ where
$X'\in T^{H+h}M_d(\tA^+_A)$ and
  $X''\in M_d(\varphi^{-n}(\A_A))$ 
  for some sufficiently large~$n$. (The reason for this choice of~$H$ is to be able
  to apply Lemma~\ref{lem: Frobenius conjugacy without phi stable} below.) Taking~$Y=X$, we may replace~$X,X',X''$ 
  by~$\varphi(X),\varphi(X'),\varphi(X'')$ respectively
(note that by Lemma~\ref{lem: phi on powers of T for A tilde plus}, we have
  $\varphi(X')\in T^{H+h}M_d(\tA^+_A)$, because 
  $p(H+h)-C\ge H+h$ by our choice of~$C$),
 which has the
  effect of replacing~$n$ by~$(n-1)$.
Iterating this procedure, we
  can assume that~$n=0$, so~$X''\in M_d(\A_A)$. By Lemma~\ref{lem:
    Frobenius conjugacy without phi stable} below, we can find~$Y$
  such that $\varphi(Y)XY^{-1}=X''$, completing the proof in the case
  that~$M$ is free.




  We now return to the general case in which $M$ is only assumed finitely
  generated projective, rather than free. By Lemma~\ref{lem: projective perfect
	  phi module is summand of free}, we may write $M$ as a direct summand of a free 
  \'etale $\varphi$-module~$F$ over~$\tA_A$. By the case already proved, we may
  write $F\isoto\widetilde{F}_0$ for some free \'etale $\varphi$-module
  $F_0$. By Proposition~\ref{prop: descent of morphisms for phi
    modules}, the idempotent in~$\End(F)$ corresponding to ~$M$ comes
  from an idempotent in~$\End(F_0)$, and we may take~$M_0$ to be the
  \'etale $\varphi$-module corresponding to this idempotent.
  \end{proof}The following lemma and its proof are based
  on~\cite[Prop.\ 2.2]{MR2562795}.
  \begin{lem}
    \label{lem: Frobenius conjugacy without phi stable}Suppose
    that~$X\in M_d(\tA^+_A)\cap\GL_d(\tA_A)$, and that
    $X^{-1}\in T^{-h} M_d(\tA^+_A)$ for some~$h\ge
    0$. Suppose that ~$X''\in M_d(\tA^+_A)\cap\GL_d(\tA_A)$ is such
    that $X^{-1}X''\in 1+T^{\lceil (C+h+1)/(p-1)\rceil}M_d(\tA^+_A)$, where~$C$
    is as in Lemma~{\em \ref{lem: phi on powers of T for A tilde plus}}. 
Then there exists $Y\in M_d(\tA^+_A)\cap\GL_d(\tA_A)$ with $X''=\varphi(Y)XY^{-1}$.
  \end{lem}
  \begin{proof}Define sequences $X_i$, $h_i$ by $X_0=X$,
    $h_0=(X'')^{-1}X$, and for each $i\ge 1$,
    $X_i=\varphi(h_{i-1})X_{i-1}h_{i-1}^{-1}$,
    $h_i=(X'')^{-1}X_i$. Then if we set $y_i=h_ih_{i-1}\cdots h_0$, we
    have $X''=X_ih_i^{-1}=\varphi(y_{i-1})Xy_{i}^{-1}$ for
    each~$i$. We claim that the~$y_i$ tend to a limit~$Y$ as
    $i\to\infty$; then we have $X''=\varphi(Y)XY^{-1}$, as required.

    To see that the ~$y_i$ tend to a limit, it is enough to show that
    $h_i\to 1$ as~$i\to\infty$. To see this, suppose
    that $h_i\in 1+T^{s}M_d(\tA^+_A)$ for some $s\ge 
    (C+h+1)/(p-1)$. We
    have \[X''h_{i+1}=X_{i+1}=\varphi(h_{i})X_{i}h_{i}^{-1}=\varphi(h_{i})X'',\]so that
    $h_{i+1}=(X'')^{-1}\varphi(h_{i})X''$. Using Lemma~\ref{lem: phi on powers of T for A tilde plus} and the assumption that $(X'')^{-1}\in T^{-h} M_d(\tA^+_A)$,
    we see that $h_{i+1}\in 1+T^{ps-C-h}M_d(\tA^+_A)$. Since
    $ps-C-h\ge s+1$, we are done.
      \end{proof}

  \subsection{An equivalence of categories}We summarise the results
  of this section in the following proposition.
  \begin{prop}
	  \label{prop:perfect equivalence}
  Let~$A$ be a finite type $\Z/p^a$-algebra for some ~$a\ge 1$.  Then the functor $M\mapsto\widetilde{M}$ is an equivalence of
    categories from the category of projective \'etale
    $\varphi$-modules over~$\A_A$ to the category of
    projective \'etale $\varphi$-modules over~$\tA_A$.
  \end{prop}
  \begin{proof}
    The functor is essentially surjective by Proposition~\ref{prop:
      descending projective modules from infinite level}, and fully
    faithful by Proposition~\ref{prop: descent of morphisms for phi modules}.
  \end{proof}
\section{
  \texorpdfstring{$(\varphi,\Gamma)$}{(phi,Gamma)}--modules}
\label{subsec: phi gamma over Zp}

By definition, an {\em \'etale $(\varphi,\Gamma_K)$-module} 
is a finitely generated 
$\A_K$-module $M$, equipped with
\begin{itemize}
\item a $\varphi$-linear morphism $\varphi:M\to M$ with the property
  that the corresponding morphism $\Phi_M:\varphi^*M\to M$ is an
  isomorphism (i.e.\ $M$ is given the structure of an \'etale
  $\varphi$-module over $\A_K$), and 
\item a continuous semi-linear action
  of~$\Gamma_K$ that commutes with~$\varphi$.
\end{itemize}


\subsection{The relationship with Galois representations}
\label{subsubsec:Galois reps}
 There is an equivalence of categories between the category of
 continuous representations of~$G_K$ on finite $\Zp$-modules, and the
 category of \'etale $(\varphi,\Gamma_K)$-modules, which is given by 
functors $\mathbb{D}$ and $T$ that are defined as
follows. Let~$\AKnrhat$ denote the $p$-adic
completion of the ring of integers
 of the
maximal unramified extension of~$\A_K[1/p]$ in~$W(\C^\flat)[1/p]$;
this is preserved by the natural actions of~$\varphi$ and~$G_K$
on~$W(\C^\flat)[1/p]$.
Then for a $G_K$-representation $T$, we define 
\[\mathbb{D}(T) :=(\AKnrhat
  \otimes_{\Zp}T)^{G_{\Kcyc}},\]
while for an \'etale $(\varphi,\Gamma_K)$-module $M$ we define
 \[T(M) :=(\AKnrhat \otimes_{\A_K}M)^{\varphi=1}.\]
The action of~$G_{\Kcyc}$ (resp.\ $\varphi$) in the definition
of $\mathbb D$ (resp.\ of $T$)
is the diagonal one (with the $G_K$-action on~$M$
being that inflated from the action of~$\Gamma_K$). There is a
completely analogous theory of~$(\varphi,\Gammat_K)$-modules, and
taking~$\Delta_K$-invariants gives an
equivalence of categories between ~$(\varphi,\Gammat_K)$-modules
and~$(\varphi,\Gamma_K)$-modules. 

The functor $T(M)$ is also defined for finitely generated \'etale $\varphi$-modules
over~$\A_K$ (i.e.\ in the absence of a  $\Gamma$-action on~$M$), 
and in this context yields an equivalence of categories between
finitely generated  \'etale $\varphi$-modules over~$\A_K$
and continuous $G_{K_{\cyc}}$-representations on finite $\Z_p$-modules.


Suppose now that $L/K$ is a finite extension, and
suppose further that $K = K_{\cyc} \cap L$ (or, equivalently,
that the natural embedding  $\Gamma_L \hookrightarrow \Gamma_K$ induced
by the inclusion $K_{\cyc}\subseteq L_{\cyc}$ is an isomorphism).  
We will recall the description of the functor $\Ind_{G_L}^{G_K}$
in terms of $(\varphi,\Gamma)$-modules.
(We will use this description in Section~\ref{subsec:ram bound}.)

By construction, there are inclusions $\A_K \subseteq \A_L \subseteq \AKnrhat$,
which are unramified embeddings of DVRs.  The composite of these embeddings
is $G_K$-equivariant (where $G_K$-acts on $\A_K$ through
its quotient~$\Gamma_K$),
while the second is $G_L$-equivariant (where $G_L$-acts on
$\A_L$ through its quotient~$\Gamma_L$).   Regarding the second of these embeddings
as an $\A_K$-linear morphism of $\A_K$-modules,
it induces an $\AKnrhat$-linear surjection
$\AKnrhat\otimes_{\A_K} \A_L \to \AKnrhat.$
The source of this surjection has a diagonal action of $G_K$
(the $G_K$-action on $A_L$ being via its quotient $\Gamma_K \iso \Gamma_L$),
while the surjection itself is $G_L$-equivariant.
Thus it induces  a $G_K$-equivariant morphism
$$\AKnrhat\otimes_{\A_K} \A_L \to
\Ind_{G_L}^{G_K} \AKnrhat,$$
which is easily seen to be an isomorphism.
Using the description of the induction as a tensor product,
we can express this as an isomorphism
\numequation
\label{eqn:induction formula}
\AKnrhat\otimes_{\A_K} \A_L \iso \Z_p[G_K]\otimes_{\Z_p[G_L]} \AKnrhat.
\end{equation}

If $M$ is an \'etale $(\varphi,\Gamma)$-module over~$\A_L$,
then we let $M'$  denote $M$ regarded as an \'etale $(\varphi,\Gamma)$-module
over the subring $\A_K$ of~$\A_L$. (Recall again that we are assuming
$\Gamma_L \iso \Gamma_K$.)
The isomorphism~\eqref{eqn:induction formula}
then yields an isomorphism 
\nummultline
\label{eqn:induction via phi Gamma mods}
T(M') := (\AKnrhat \otimes_{\A_K} M')^{\varphi = 1}  =
\bigl((\AKnrhat\otimes_{\A_K} \A_L) \otimes_{\A_L} M\bigr)^{\varphi =1}
\\
= 
(\Z_p[G_K]\otimes_{\Z_p[G_L]} \AKnrhat \otimes_{\A_L} M\bigr)^{\varphi = 1}
=\Z_p[G_K]\otimes_{\Z_p[G_L]} (\AKnrhat \otimes_{\A_L} M\bigr)^{\varphi = 1}
\\
= \Ind_{G_L}^{G_K} T(M). 
\end{multline}
In short, induction of Galois representations corresponds to
the restriction of coefficients on the $(\varphi,\Gamma)$-module side.

The preceding discussion also applies in the context of \'etale
$\varphi$-modules (in the sense of Definition~\ref{defn: general phi module}).
Our assumption that $\Gamma_L \iso \Gamma_K$ implies that
$G_{K_{\cyc}}/G_{L_{\cyc}}  \to G_K/G_L$  is a bijection,
so that \eqref{eqn:induction formula} may also be interpreted as an isomorphism
$$
\AKnrhat\otimes_{\A_K} \A_L \iso \Z_p[G_{K_{\cyc}}]\otimes_{\Z_p[G_{L_{\cyc}}]}
\AKnrhat.
$$
Then,
if $M$ is a finitely generated \'etale $\varphi$-module over~$\A_L$,
and if $M'$ denotes its restriction to an \'etale $\varphi$-module over~$\A_K$,
the preceding isomorphism induces an isomorphism
\numequation 
\label{eqn:induction via phi mods}
T(M') = \Ind_{G_{L_{\cyc}}}^{G_{K_{\cyc}}} T(M). 
\end{equation}

We also need a variant of the preceding theory
which again follows from the results
of~\cite{MR1106901}, using the Kummer extension~$K_\infty/K$
introduced in
Example~\ref{ex:kummer} and Section~\ref{subsubsec: Kummer}; the integral version of this theory was first
studied by Breuil and Kisin, see~\cite[\S1]{KisinModularity}. 
Namely,
there is an equivalence of categories between the
category of continuous representations of~$G_{K_\infty}$ on finite
$\Zp$-algebras, and the category of \'etale
$\varphi$-modules 
over~$\cO_{\cE}$, which is given by functors
$\mathbb{D}_\infty$, $T_\infty$ that are defined as
follows: Let~$\cO_{\widehat{\cE^{\nr}}}$ denote the $p$-adic
completion of the ring of integers
in the
maximal unramified extension of~$\Frac(\cO_{\cE})$
in~$W(\C^\flat)[1/p]$;
this is preserved by the natural actions of~$\varphi$ and~$G_{K_{\infty}}$
on~$W(\C^\flat)[1/p]$.
We define
\[\mathbb{D}_\infty(T) :=( \cO_{\widehat{\cE^{\nr}}} \otimes_{\Zp}T)^{G_{K_\infty}},\]
for a $G_{K_\infty}$-representation $T$,
and define
\[T_\infty(M):=(  \cO_{\widehat{\cE^{\nr}}}\otimes_{\A_K}M)^{\varphi=1},\]
for an \'etale $\varphi$-module $M$.

\subsection{Various types of $\varphi$-modules with coefficients}
We now introduce coefficients.

\begin{df}
	\label{def:etale phi Gamma modules}
       
Let~$A$ be a $p$-adically complete $\Zp$-algebra.	A {\em
  projective \'etale $(\varphi,\Gamma_K)$-module of rank~$d$} \index{\'etale $(\varphi,\Gamma_K)$-module}
	with $A$-coefficients is a projective \'etale $\varphi$-module \index{\'etale $(\varphi,\Gamma_K)$-module}
        $M$ of rank~$d$
	over $\A_{K,A}$ equipped
	with a semi-linear action of $\Gamma_K$, which commutes
        with~$\varphi$, and which is furthermore continuous
	when $M$ is endowed with its canonical topology (i.e.\ the
        topology of Remark~\ref{rem:topologies on M}).
\end{df}

If $A$ is a finite
type $\Z/p^a$-algebra, we can use Theorem~\ref{thm:descending projective modules} to give
a useful alternative description of \'etale $(\varphi,\Gamma_K)$-modules with $A$-coefficients in terms of
\'etale $\varphi$-modules with
coefficients in the ring $W(\C^\flat)_A$, as we now explain.



\begin{defn}
  \label{defn: GK module over Ainf}Let~$A$ be a $p$-adically complete $\Zp$-algebra. 
An \'etale \emph{$(\varphi,G_K)$-module \index{\'etale $(\varphi,G_K)$-module}
with $A$-coefficients} (resp.\ an \'etale
\emph{$(\varphi,G_{\Kcyc})$-module \index{\'etale $(\varphi,G_{\Kcyc})$-module}
with $A$-coefficients}, resp.\ an \'etale
\emph{$(\varphi,G_{K_\infty})$-module \index{\'etale $(\varphi,G_{K_\infty})$-module}
with $A$-coefficients}) is by definition a finitely
generated~$W(\C^\flat)_A$-module $M$ equipped with an isomorphism
\[\varphi_M:\varphi^*M\isoto M\]
 of $W(\C^\flat)_A$-modules, 
 and
a~$W(\C^\flat)_A$-semi-linear action of~$G_K$ (resp.\ $G_{\Kcyc}$,
resp.\ $G_{K_\infty}$), which is continuous 
and commutes with~$\varphi_M$.
We say that~$M$ is projective if it is projective of constant rank as
a~$W(\C^\flat)_A$-module.
\end{defn}

If~$A$ is a finite type $\Z/p^a$-algebra, then
there is a functor from the category of finite projective \'etale
$(\varphi,\Gamma_K)$-modules 
with~$A$-coefficients to the category of finite projective
\'etale $(\varphi,G_K)$-modules with~$A$-coefficients, which takes an \'etale
$(\varphi,\Gamma_K)$-module~$M$ to $W(\C^\flat)_A\otimes_{\A_{K,A}}M$, endowed
with the extension of scalars of~$\varphi$, and the diagonal action
of~$G_K$, with the action of~$G_K$ on~$M$ being the action inflated
from~$\Gamma_K$. 

Similarly, there is a functor from the category of finite projective
\'etale $\varphi$-modules over~$\OEA$ 
to the category of finite projective
\'etale $(\varphi,G_{K_\infty})$-modules with~$A$-coefficients, which takes an \'etale
$\varphi$-module~$M$ to $W(\C^\flat)_A\otimes_{\OEA}M$, endowed
with the extension of scalars of~$\varphi$, and the diagonal action
of~$G_{K_\infty}$, with the action of~$G_{K_\infty}$ on~$M$ being the
trivial action.




\begin{prop}
  \label{prop: equivalences of categories to Ainf}Let~$A$ be a
  finite type $\Z/p^a$-algebra for some~$a\ge 1$.
  \begin{enumerate}
  \item The functor $M\mapsto W(\C^\flat)_A\otimes_{\A_{K,A}}M$ is an
    equivalence between the category of finite projective \'etale
$\varphi$-modules over~$\A_{K,A}$ and the category of
finite projective $(\varphi,G_{\Kcyc})$-modules
with $A$-coefficients.

It induces  an
    equivalence of categories between the  category of finite projective \'etale
$(\varphi,\Gamma_K)$-modules with $A$-coefficients and the category of
finite projective \'etale $(\varphi,G_K)$-modules
with $A$-coefficients. 

A quasi-inverse functor is given by the
composite of $N\mapsto N^{G_{\Kcyc}}$ and a quasi-inverse to the functor
of Proposition~{\em \ref{prop:perfect equivalence}}.
\item The functor $M\mapsto W(\C^\flat)_A\otimes_{\OEA}M$ is an
    equivalence of categories between the  category of finite projective \'etale
$\varphi$-modules over~$\OEA$ and the category of
finite projective \'etale $(\varphi,G_{K_\infty})$-modules
with $A$-coefficients. A quasi-inverse functor is given by $N\mapsto
N^{G_{K_\infty}}$ and a quasi-inverse to the functor of  Proposition~{\em \ref{prop:perfect equivalence}}.
  \end{enumerate}

\end{prop}
\begin{proof}
  This is a formal consequence of Theorem~\ref{thm:descending
    projective modules} and Proposition~\ref{prop:perfect
    equivalence}.  
  We begin with~(2). Firstly, by
  Proposition~\ref{prop:perfect equivalence} (and Remark~\ref{rem: abstract perfection section linked to our concrete settings}) the functor
  $M\mapsto \tM=\tOEA\otimes_{\OEA}M$ is an equivalence of categories
  between the category of finite projective \'etale $\varphi$-modules
  over~$\OEA$ and the category of finite projective \'etale
  $\varphi$-modules over~$\tOEA$, so it is enough to show that
  $\tM\mapsto W(\C^\flat)_A\otimes_{\tOEA}\tM$ is an equivalence of
  categories. By Example~\ref{ex:kummer} and Remark~\ref{rem: Krasner
    Ax Tate Sen}, $\widehat{K}_\infty$ is a perfectoid field, and we
  have a canonical identification of Galois groups
  $G_{\widehat{K}_\infty}=G_{K_\infty}$. The result then follows from
  the equivalence of categories given by Theorem~\ref{thm:descending
    projective modules}, as we can think of $\varphi$ as being an
  isomorphism of $\tOEA$-modules $\varphi^*\tM\isoto\tM$.

The same argument shows in~(1) that the functor
$M\mapsto W(\C^\flat)_A\otimes_{\A_{K,A}}M$ is an equivalence of
categories between the category of projective \'etale
$\varphi$-modules over~$\A_{K,A}$, and the category of
projective~\'etale $(\varphi,G_{\Kcyc})$-modules with
$A$-coefficients. Since~$\Gamma_K=G_K/G_{\Kcyc}$, this extends to the
claimed equivalence of categories, noting that by construction the
continuity of the $\Gamma_K$-action on~$M$ is equivalent to the
continuity of the $G_K$-action on~$W(\C^\flat)_A\otimes_{\A_{K,A}}M$
(since
Lemma~\ref{lem: varphi extends to A and Galois etc}
shows  that  the $G_K$-action on~$W(\C^\flat)_A$ is continuous).
\end{proof}




\chapter{Moduli stacks of \texorpdfstring{$\varphi$}{phi}-modules and 
  \texorpdfstring{$(\varphi,\Gamma)$}{(phi,Gamma)}-modules}\chaptermark{Stacks of \texorpdfstring{$\varphi$}{phi}-modules and 
  \texorpdfstring{$(\varphi,\Gamma)$}{(phi,Gamma)}-modules}\label{sec:
phi modules and phi gamma modules}

In this chapter we build on the results
of~\cite{EGstacktheoreticimages} (which in turn built
on~\cite{MR2562795}) to construct our stacks of
$(\varphi,\Gamma)$-modules, and various related stacks of
$\varphi$-modules. We show in particular that our stack~$\cX_{K,d}$ of
$(\varphi,\Gamma)$-modules is Ind-algebraic.

\section{Moduli stacks of \texorpdfstring{$\varphi$}{phi}-modules}\label{subsec: moduli stacks of phi
  modules}

In this section we put ourselves in the context of Situation~\ref{subsubsec:general
	framework};
we also remind the reader that the notion of \'etale $\varphi$-module is
defined in Section~\ref{subsec: EG stuff}.
We furthermore fix
a finite extension~$E/\Qp$ with ring of
integers~$\cO$ and residue field~$\F$; all of our coefficient rings
from now on will be $\cO$-algebras. All of our constructions are
compatible with replacing~$E$ by a finite extension, and we will
typically not comment on this, although see Remark~\ref{rem: changing O} below.

If we fix integers $a,d \geq 1,$
then we may follow~\cite[\S
5]{EGstacktheoreticimages} and define
an \emph{fpqc} stack in groupoids $\cR_d^a$ over
$\Spec \cO/\varpi^a$ as follows: For any $\cO/\varpi^a$-algebra~$A$,
we define $\cR^a_d(A)$ to be the groupoid of
 \'etale $\varphi$-modules over~$\A_A$ which are projective 
 of rank~$d$.
If $A \to B$ is a morphism of $\cO$-algebras, and $M$
is an object of $\cR_d^a(A)$,
then the pull-back of $M$ to $\cR_d^a(B)$
is defined to be the tensor product~$\AAA_B\otimes_{\AAA_A} M$.

A key point is that this definition does not require $\A_A^+$ to be $\varphi$-stable,
although (as far as we know) this hypothesis {\em is} required to make any deductions about
$\cR^a_d$ beyond the fact that it is an {\em fpqc} stack (which relies just
on Drinfeld's general descent results, as described
in \cite{MR2181808} and \cite[\S 5.1]{EGstacktheoreticimages}).

Since we are ultimately interested in questions of algebraicity or of Ind-algebraicity,
from now on we regard $\cR^a_d$ as an {\em fppf} stack over~$\cO/\varpi^a$. 
By~\cite[\href{http://stacks.math.columbia.edu/tag/04WV}{Tag
  04WV}]{stacks-project}, 
we may also regard the stack $\cR_d^a$ as an {\em fppf} stack over $\cO$,
and as $a$ varies, we may form the $2$-colimit $\cR := \varinjlim_a \cR^a_d$,
which is again an {\em fppf} stack over $\cO$.
In fact $\cR_d$ lies over $\Spf \cO := \varinjlim_a \Spec \cO/\varpi^a$,
the formal spectrum of $\cO$ with respect to the $\varpi$-adic,
or equivalently $p$-adic, topology.

We now fix a polynomial $F\in W(k)[T]$ which is congruent to a positive power
of~$T$ modulo~$p$ (for example, an Eisenstein polynomial). 

\begin{defn}
  \label{defn: phi module of finite E-height}
Suppose that~$\A_A^+$ is $\varphi$-stable. Let $h$ be a non-negative
integer, and let~$A$ be a $p$-adically complete $\cO$-algebra. A  \emph{$\varphi$-module of
   $F$-height at most~$h$ over~$\A^+_A$} \index{finite height $\varphi$-module}
is  a pair $(\gM,\varphi_M)$
consisting of a finitely generated $T$-torsion free
$\Aplus_A$-module~$\gM$, and a 
$\varphi$-semi-linear map $\varphi_\gM:\gM\to\gM$, with the further
properties that if we
write \[\Phi_{\gM}:=1\otimes \varphi_{\gM}:\varphi^*\gM\to\gM,\] then 
$\Phi_{\gM}$ is injective, and the cokernel of $\Phi_{\gM}$ is killed by $F^h$.

A \emph{$\varphi$-module of finite $F$-height over~$\A^+_A$} is a $\varphi$-module of
  $F$-height at most~$h$  for some $h\ge 0$. A
  morphism of $\varphi$-modules is a morphism of the underlying
  $\Aplus_A$-modules which commutes with the morphisms~$\Phi_\gM$.

  We say that a $\varphi$-module of finite $F$-height is projective of
  rank~$d$ if it is a finitely generated projective~$\Aplus_A$-module of
  constant rank~$d$.
\end{defn}

If we maintain the assumption that $\A^+$ is~$\varphi$-stable,
and if we fix integers $a,d\ge 1$ and an integer $h\ge 0$,
then we may again follow~\cite[\S
5]{EGstacktheoreticimages} to define an 
\emph{fpqc} stack in groupoids $\cC_{d,h}^a$
over $\Spec \cO/\varpi^a$ as follows: For any $\cO/\varpi^a$-algebra~$A$,
we define $\cC^a_{d,h}(A)$ to be the groupoid of $\varphi$-modules of
$F$-height at most $h$ over~$\A^+_A$ which are projective
 of rank~$d$.
If $A \to B$ is a morphism of $\Z_p$-algebras, and $\mathfrak M$
is an object of $\cC_{d,h}^a(A)$ 
then the pull-back of $\gM$ to $\cC_{d,h}^a(B)$
is defined to be the tensor product
$\Aplus_{B}\otimes_{\Aplus_A} \gM$.

Just as for the stack $\cR_d^a$,
we may and do also regard the stack $\cC_{d,h}^{a}$ as an {\em fppf}
stack over $\cO$,
and we then, allowing $a$ to vary, define $\cC_{d,h}:=\varinjlim_{a}\cC_{d,h}^a$,
obtaining an {\em fppf} stack over $\cO$
which in fact lies over $\Spf \cO$.
There are  canonical morphisms $\cC_{d,h}^a\to\cR_d^a$ and $\cC_{d,h}\to\cR_d$ given 
  by tensoring with $\A_A$ over~$\A_A^+$. 

  \begin{rem}
    \label{rem: changing O}If~$E'/E$ is a finite extension with ring
    of integers~$\cO'$, then by definition we have (with obvious
    notation) $\cC_{d,h,\cO'}=\cC_{d,h}\times_{\cO}\cO'$ (in the case
that $\A^+$ is $\varphi$-stable, so that these stacks are defined) and
    $\cR_{d,\cO'}=\cR_{d,\cO}\times_\cO\cO'$ (in general). 
  \end{rem}

The following lemma provides a concrete interpretation of the
$A$-valued points of the stacks we have defined, when $A$ is a $\varpi$-adically
complete $\cO$-algebra (rather than just an algebra over some $\cO/\varpi^a$).

\begin{lem}
  If $A$ is a $\varpi$-adically complete $\cO$-algebra, then there
is a canonical equivalence between the groupoid of morphisms
$\Spf A\to\cR_{d}$ and the groupoid of rank~$d$
 \'etale $\varphi$-modules over~$\A_A$. 
If $\A_A^+$ is furthermore $\varphi$-stable, 
then there is a canonical equivalence between the groupoid of morphisms
  $\Spf A\to\cC_{d,h}$
and the groupoid of $\varphi$-modules of rank~$d$ and $F$-height at most~$h$
 over~$\A^+_A$.
\end{lem}
\begin{proof}
  This is immediate from Lemma~\ref{lem: can check projectivity of Kisin and etale phi modules modulo
    p^n}.
\end{proof}

We now apply the results of \cite[\S 5]{EGstacktheoreticimages} to deduce
various results about the stacks we have introduced.  
This requires the assumption that $\A^+$ is $\varphi$-stable.

\begin{thm}
  \label{thm: C is a $p$-adic formal algebraic stack and related
    properties}Suppose that~$\A^+$ is $\varphi$-stable, and let $a\ge 1$ be arbitrary.
  \begin{enumerate}
  \item The stack $\cC_{d,h}^a$ is an algebraic stack of finite
    presentation over $\Spec\cO/\varpi^a$, with affine diagonal.
  \item The morphism $\cC_{d,h}^a\to\cR_d^a$ is representable by
    algebraic spaces, proper, and of finite presentation.
  \item The diagonal morphism
    $\Delta:\cR_d^a\to\cR_d^a\times_{\cO/\varpi^a}\cR_d^a$ is
    representable by algebraic spaces, affine, and of finite
    presentation.
  \item  $\cR^a_d$ is a limit preserving Ind-algebraic stack, whose diagonal is
    representable by algebraic spaces, affine, and of finite
    presentation. 
  \end{enumerate}
\end{thm}
\begin{proof}
  Part (1) is~\cite[Thm.\ 5.4.9~(1)]{EGstacktheoreticimages}, and
  parts~(2), (3) and~(4) are proved in~\cite[Thm.\
  5.4.11]{EGstacktheoreticimages}, except for the claim that~$\cR^a_d$
  is Ind-algebraic, which is~\cite[Thm.\
  5.4.20]{EGstacktheoreticimages}. 
  \end{proof}

\begin{cor}
  \label{cor: basic properties of C and R p adic stacks}Suppose
  that~$\A^+$ is $\varphi$-stable. 
  \begin{enumerate}
  \item $\cC_{d,h}$ is a $p$-adic formal algebraic stack of finite
    presentation over $\Spf\cO$, with affine diagonal.
 \item  $\cR_d$ is a limit preserving Ind-algebraic stack, whose diagonal is
    representable by algebraic spaces, affine, and of finite
    presentation.
  \item The morphism $\cC_{d,h}\to\cR_d$ is representable by
    algebraic spaces, proper, and of finite presentation.
  \item The diagonal morphism
    $\Delta:\cR_d\to\cR_d\times_{\Spf\cO}\cR_d$ is
    representable by algebraic spaces, affine, and of finite
    presentation.
   \end{enumerate}
\end{cor}
\begin{proof}
  The first part is immediate from~ Theorem~\ref{thm: C is a $p$-adic
    formal algebraic stack and related properties}~(1) and
  Proposition~\ref{prop: criterion for p-adic formal algebraic
    stack}. Everything else is immediate from Theorem~\ref{thm: C is a
    $p$-adic formal algebraic stack and related properties}.
\end{proof}

\section{Moduli stacks of
  \texorpdfstring{$(\varphi,\Gamma)$}{(phi,Gamma)}-modules}\label{subsec:
defn of Xd}
In this section we begin the study of our main objects of interest,
namely the moduli stacks of \'etale $(\varphi,\Gamma)$-modules.
 As in Chapter~\ref{section: coefficient rings} 
 we fix a finite extension $K/\Q_p$.  As in Section~\ref{subsec: moduli
 stacks of phi modules}, we also fix a finite extension $E$ of $\Q_p$
 with ring of integers $\cO$, which will serve as our ring of coefficients.
As always, $k$ denotes the residue field of the ring of integers of $K$,
and $\F$ denotes the residue field of~$\cO$.



\begin{df} \label{defn: Xd}
	We let $\cX_{K,d}$ denote the moduli stack of projective \index{$\cX_{K,d}$}
	\'etale $(\varphi,\Gamma_K)$-modules of rank~$d$. More
        precisely, if~$A$ is a $p$-adically complete $\cO$-algebra,
        then we define $\cX_{K,d}(A)$ (i.e., the groupoid of morphisms
        $\Spf A\to\cX_{K,d}$) to be the groupoid of 
        projective \'etale~$(\varphi,\Gamma_K)$-modules  of rank~$d$ with~$A$-coefficients,
in the sense of Definition~\ref{def:etale phi Gamma modules},
 with morphisms given by isomorphisms. If $A\to B$ is a morphism of complete
        $\cO$-algebras, and $M$ is an object of $\cX_{K,d}(A)$, then the
        pull-back of~$M$ to~$\cX_{K,d}(B)$ is defined to be the tensor
        product $\A_{K,B}\otimes_{\A_{K,A}}M$. 

        It follows from the results of~\cite{MR2181808},
        and more specifically from~\cite[Thm.\
        5.1.16]{EGstacktheoreticimages}, that~$\cX_{K,d}$ is
        an~\emph{fpqc} stack over~$\cO$. As in Remark~\ref{rem:
          changing O}, the definition of~$\cX_{K,d}$ behaves
        naturally with respect to change of the coefficient ring~$\cO$.\end{df}

One of the main results of this book is that
$\cX_{K,d}$ is a Noetherian formal algebraic stack.  However,
the proof of this is quite involved, and will only be fully
achieved at the conclusion of Chapter~\ref{sec: families of extensions}.
In this section and the two that follow it, we establish
the preliminary result that $\cX_{K,d}$ is an Ind-algebraic stack.

We begin by discussing the moduli stacks of \'etale $\varphi$-modules
over~$\A_{K,A}$, which  will play an auxiliary role in our study 
of~$\cX_{K,d}$. 

\begin{df} \index{$\cR_{K,d}$}
We let~$\cR_{K,d}$ denote the moduli stack of rank $d$ projective \'etale
$\varphi$-modules, defined as in Section~\ref{subsec: moduli stacks of
  phi modules}, taking $\A$ to be $\A_K$.
\end{df}

If~$\A_K^+$ is
$\varphi$-stable, then Corollary~\ref{cor: basic properties of C and R
  p adic stacks}  applies to~$\cR_{K,d}$; this
is in particular the case if~$K/\Qp$ is basic in the sense of
Definition~\ref{defn: basic}. Our first task is to establish the same
results for general~$K$, which we will do by reducing to the basic
case.

\begin{defn}\label{defn: Kbasic} \index{$\Kbasic$}
  If~$K/\Qp$ is any finite extension, we
  set~$\Kbasic:=K\cap K_0(\zeta_{p^\infty})$.
\end{defn}
By definition, $\Kbasic$ is basic, and~$(\Kbasic)_0=K_0$. Note that the natural
  restriction map~$\Gammat_K\to\Gammat_{\Kbasic}$ is an
  isomorphism, and induces an isomorphism
  $\Gamma_K\to\Gamma_{\Kbasic}$. By~\cite[Thm.~3.1.2]{MR719763}, $\A_K$ is a
  free~$\A_{\Kbasic}$-module of
  rank~$[K:\Kbasic]=[K(\zeta_{p^{\infty}}):K_0(\zeta_{p^\infty})]$, 
and the inclusion $\A_{\Kbasic}\subset\A_K$ is $\varphi$-equivariant,
so there is a natural morphism $\cR_{K,d}\to\cR_{\Kbasic,d[K:\Kbasic]}$ given
by forgetting the $\A_{K,A}$-algebra structure on an \'etale
$\varphi$-module.

\begin{remark}
\label{rem:explanation of reduction to basic case}
As noted in Remark~\ref{rem: not phi stable T},
in order to ensure that $\A_K^+$ is $\varphi$-stable,
it suffices for $K$ to be abelian over $\Q_p$.  Thus,
rather than relating the theory for $K$ to that for the field $\Kbasic$ introduced above,
we could just as well relate it to any other subfield $K'$ of $K$ which is abelian
over $\Q_p$, and for which the natural map $\Gamma_K \to \Gamma_{K'}$ is an isomorphism; 
e.g.\ the field $K' := K \cap \Q_p^{\ab}$ (where $\Q_p^{ab}$ denotes
the maximal abelian extension of $K$).  
It doesn't matter (for our purposes) which particular $K'$ we choose; 
$\Kbasic$ is simply a convenient choice.
\end{remark}

\begin{lem}\label{lem: relative representability for R over Kbasic}
  The morphism $\cR_{K,d}\to\cR_{\Kbasic,d[K:\Kbasic]}$ is representable by
  algebraic spaces, affine, and of finite presentation.
\end{lem}
\begin{proof}We can prove the statement after pulling back via a
  morphism $\Spec A\to \cR_{\Kbasic,d[K:\Kbasic]}$, where~$A$ is an
  $\cO/\varpi^a$-algebra for some~$a\ge 1$. This morphism corresponds to
  a projective \'etale $\varphi$-module over~$\A_{\Kbasic,A}$ of
  rank~$d[K:\Kbasic]$, and we need to show that the functor on
  $A$-algebras taking~$B$ to the set of projective \'etale $\varphi$-modules over~$\A_{K,B}$ of
  rank~$d$, whose underlying \'etale $\varphi$-module
  over~$\A_{\Kbasic,B}$ coincides with~$M_B$, is representable by an affine scheme of finite presentation over~$\Spec A$.

  Note that since~$M_B$ is projective and in particular flat
  over~$\A_{\Kbasic,B}$, we have a natural inclusion
  $i:M_B\into \A_{K,B}\otimes_{\A_{\Kbasic,B}}M_B$.  The additional
  structure needed to make~$M_B$ into a projective \'etale
  $\varphi$-module over~$\A_{K,B}$ is the data of a morphism of
  \'etale $\varphi$-modules over~$\A_{\Kbasic,B}$
  \[f: \A_{K,B}\otimes_{\A_{\Kbasic,B}}M_B\to M_B \] satisfying the
  conditions that
\begin{enumerate}
\item the composite $M_B\stackrel{i}{\to}
  \A_{K,B}\otimes_{\A_{\Kbasic,B}}M_B\stackrel{f}{\to} M_B $ is the
  identity morphism, and
\item the kernel of~$f$ is $\A_{K,B}$-stable.
\end{enumerate}
(We can then define the ~$\A_{K,B}$-module structure
on~$M_B$ via the formula $\lambda\cdot m:=f(\lambda\otimes
m)$. The first condition guarantees that this action is compatible
with the existing ~$\A_{\Kbasic,B}$-module structure
on~$M_B$, and the second condition that $(\lambda_1\lambda_2)\cdot
m=\lambda_1\cdot(\lambda_2\cdot m)$.)

By~\cite[Prop.\ 5.4.8]{EGstacktheoreticimages}, the data of a morphism
of \'etale $\varphi$-modules
$f: \A_{K,B}\otimes_{\A_{\Kbasic,B}}M_B\to M_B$ is representable by an
affine scheme of finite presentation over~$\Spec A$, so it is enough
to show that conditions~(1) and~(2) are closed conditions, given by
finitely many equations. To see this, we follow the proof
of~\cite[Prop.\ 5.4.8]{EGstacktheoreticimages}. Exactly as in that
argument, we can reduce to the case that~$M_B$ is free, and after
choosing bases, any~$f$ is determined by the coefficients of finitely
many powers of~$T$ in the Laurent series expansions of the entries of
the matrix given by~$f$. Condition~(1) is then evidently given by
finitely many equations in these coefficients.

To see that the same is true of condition~(2), note that since~$M_B$ is
projective and condition~(1) implies in particular that~$f$ is
surjective, we have a splitting
$\A_{K,B}\otimes_{\A_{\Kbasic,B}}M_B=M_B\oplus\ker(f)$. The projection
onto~$\ker(f)$ is given by $(1-i\circ f)$, so the condition
that~$\ker(f)$ is $\A_{K,B}$-stable is the condition that for
any~$\lambda$ in~$\A_{K,B}$, and any~$m\in M_B$, we
have \[f(\lambda(m-i(f(m))))=0. \]This is evidently a closed
condition, and since $\A_{K,B}$ and~$M_B$ are both finitely
generated~$\A_{\Kbasic,B}$-modules, it is determined by finitely many
equations, as required.
\end{proof}

\begin{cor}\label{cor: R good properties all K}
  The stack $\cR_{K,d}$ is a limit preserving
  Ind-algebraic stack, whose diagonal is representable by algebraic
  spaces, affine, and of finite presentation.
\end{cor}
\begin{proof}
	This follows from
	Lemma~\ref{lem: relative representability for R over Kbasic},
  Corollary~\ref{cor: basic properties of C and R p adic stacks}
(which establishes the claimed properties for $\cR_{\Kbasic, d[K:\Kbasic]}$),
and Corollary~\ref{cor:deducing properties} below.
\end{proof}
 
The following series of results concerning morphisms of stacks
culminates in Corollary~\ref{cor:deducing properties},
which was used in the proof of Corollary~\ref{cor: R good properties all K}.

\begin{lem}
  \label{lem:finiteness of diagonals}
  Let $\cX\to\cY$ be a morphism of stacks
  over a base scheme~$S$, which is representable by algebraic
  spaces. 
   As usual, let $\Delta_f: \cX \to  \cX\times_{\cY} \cX$
	  denote the diagonal of~$f$.
%
If $f$ is of finite type and quasi-separated,
 		  then  $\Delta_f$ is of finite presentation.
\end{lem}
   \begin{proof}
	   This can be checked after
	   pulling back along an arbitrary morphism $T \to \cY$,
	   where $T$ is a scheme, and hence reduced to the case
	   of a morphism from an algebraic space to a scheme.  In this case,
	   the claim of the lemma is proved
  in~\cite[\href{http://stacks.math.columbia.edu/tag/084P}{Tag
  084P}]{stacks-project}.  
   \end{proof}

  \begin{lem}
  \label{lem: abstract pullback diagonals lemma}
  Let $f:\cX\to\cY$ be a
  morphism of stacks over a base scheme~$S$ which
  is representable by
  algebraic spaces, has affine diagonal, and is of finite type.
  Suppose that the
  diagonal of~$\cY$ is representable by algebraic
  spaces, affine, and of finite presentation. Then the
  diagonal of~$\cX$ is also representable by algebraic
  spaces, affine, and of finite presentation.
\end{lem}
\begin{proof}
  We may factor the diagonal of~$\cX$
  as \[\cX\to\cX\times_{\cY}\cX\to\cX\times_{S}\cX.\]Since the
  morphism $\cX\times_{\cY}\cX\to\cX\times_{S}\cX$ is pulled back from
  the diagonal $\cY\to\cY\times_S\cY$, it is enough to show that
  the relative diagonal
  $\Delta_f: \cX\to\cX\times_{\cY}\cX$ is representable by algebraic spaces,
  affine, and of finite presentation.

  Now $\Delta_f$ is representable by algebraic spaces (since $f$ is),
  and affine (by assumption).  Since affine morphisms
  are quasi-compact, we see that $f$ is quasi-separated, as well as
  being of finite type (by assumption); thus $\Delta_f$ is 
  of finite presentation (by Lemma~\ref{lem:finiteness
	  of diagonals}).  \end{proof} 

\begin{cor}\label{cor:deducing properties}
	Let $\cX \to \cY$ be a morphism of stacks which
	is representable by algebraic spaces,
	affine, and of finite presentation.
	If $\cY$ is a limit preserving Ind-algebraic stack,
	whose diagonal is representable by algebraic spaces,
	affine, and of finite presentation,
	then the same is true of $\cX$.
\end{cor}
\begin{proof}
The claimed limit preserving property of $\cX$ follows,
by \cite[Lem.~2.3.20~(3)]{EGstacktheoreticimages},
from that of $\cY$
and the fact that 
$\cX\to \cY$
is of finite presentation,
while the claimed Ind-algebraic property of $\cX$ 
follows from that of $\cY$
and the fact that 
$\cX\to \cY$
is representable by algebraic spaces. 
Finally,
the claimed properties of the diagonal of $\cX$
follow from those of $\cY$
by Lemma~\ref{lem: abstract pullback diagonals lemma}.
(Note that by assumption the hypotheses of
Lemma~\ref{lem: abstract pullback diagonals lemma}
hold, since morphisms of finite presentation are in particular of finite type,
while affine morphisms are separated, and thus have their diagonals
being closed immersions, which in particular are again affine.)
\end{proof}

We now give a concrete description of~$\cR_{K,d}^a$ as an Ind-algebraic stack
which will be useful in Chapter~\ref{sec: the rank one case}.  In order to state it,
we introduce the notation
$\cC_{\Kbasic,d[K:\Kbasic],h}^a$ for the stack
over~$\Spec\cO/\varpi^a$ classifying rank $d$ projective
$\varphi$-modules over~$\A_{\Kbasic,A}^+$ of
$T$-height at most~$h$; by Theorem~\ref{thm: C is a $p$-adic formal algebraic stack and related
    properties}, $\cC_{\Kbasic,d[K:\Kbasic],h}^a$ is an algebraic
stack of finite presentation over~$\Spec\cO/\varpi^a$.

We begin with a lemma which is a variant of one of the steps appearing 
in the proof of~\cite[Thm.~5.4.20]{EGstacktheoreticimages}.

\begin{lemma}
\label{lem:lifting from R to C}
If $T \to \cR_{K,d}^a$ is a morphism
whose source is a Noetherian scheme,
then there is a Noetherian scheme $Z$ and a scheme-theoretically dominant and
surjective morphism
$Z\to T$ such that the composite $Z \to T \to \cR_{K,d}^a 
\to \cR^a_{\Kbasic,d[K:\Kbasic]}$ can be lifted
to a morphism $Z \to \cC_{\Kbasic,d[K:\Kbasic],h}^a$
for some sufficiently large value of~$h$.
\end{lemma}
\begin{proof}
Since $T$ is Noetherian, is is quasi-compact, and so admits a scheme-theoretically
dominant surjection from a Noetherian affine scheme (e.g.\ the disjoint union
of the members of a finite cover of $T$ by affine open subsets).  Thus we are
reduced to the affine case.

Since $T$ is affine, by~\cite[Prop.~5.4.7]{EGstacktheoreticimages},
we may find a scheme-theoretically dominant surjection $Z  = \Spec A\to T$
such that the \'etale $\varphi$-module $M$ over $A$ corresponding
to the composite $Z \to T \to \cR_{K,d}$ is free (of rank $d$)
over $\A_{K,A}$ (rather
than merely projective).  Thus $M$ is also free (of rank $d[K:\Kbasic]$) 
over $\A_{\Kbasic,A}$.  We may then choose a $\varphi$-stable
$\A_{\Kbasic,A}$-basis
of~$M$.  If we let $\gM$ denote the $\A_{\Kbasic,A}^+$-module
spanned by this basis, then $\gM$ is $\varphi$-invariant, and has some height~$h$.
Thus $\gM$ induces the desired morphism $Z \to \cC_{\Kbasic,d[K:\Kbasic],h}^a.$
\end{proof}

We now give the promised description of the Ind-algebraic stack structure on~$\cR_{K,d}$.

\begin{lem}
  \label{lem: explicit Ind description for R pulled back from K basic}
  Fix~$a \geq 1$, and for each $h \geq 0$,
  let $\cR_{K,d,\Kbasic,h}^a$ denote the scheme-theoretic image of the
  base-changed morphism
  \[\cC_{d[K:\Kbasic],h}^a\times_{\cR_{\Kbasic,d[K:\Kbasic]}}\cR_{K,d}\to\cR_{K,d},\]
  so that $\cR_{K,d,\Kbasic,h}^a$ is a closed algebraic substack of $\cR^a_{K,d}$.
  Then the canonical morphism
  $\varinjlim_h\cR_{K,d,\Kbasic,h}^a \to \cR_{K,d}^a$
  is an isomorphism.
\end{lem}
\begin{proof}
  We have to show that any morphism $T \to \cR_{K,d}^a$ whose source
is a scheme factors through
the inductive limit.  Since $\cR_{K,d}^a$ is limit preserving, we may assume
that $T$ is of finite type over $\cO/\varpi^a$.  Let $Z\to T$
be as in the statement of Lemma~\ref{lem:lifting from R to C}.
Then the composite
\numequation
\label{eqn:dominant composite}
Z\to T \to \cR^a_{K,d}
\end{equation}
 lifts to a morphism
$Z \to \cC_{d[K:\Kbasic],h}^a\times_{\cR_{\Kbasic,d[K:\Kbasic]}}\cR_{K,d},$
and hence the morphism~\eqref{eqn:dominant composite} 
factors through 
  $\cR_{K,d,\Kbasic,h}^a$.
Since $Z \to T$ is scheme-theoretically dominant, the original morphism
$T \to \cR_{K,d}^a$ also factors through $\cR_{K,d,\Kbasic,h}^a$,
and hence through the inductive limit.
Thus the lemma is proved.
(We remark that this argument is essentially identical to that used to 
prove~\cite[Thm.~5.4.20]{EGstacktheoreticimages}.)
\end{proof}

%
%
%
%

We now turn to studying $\Gamma_K$-actions on our $\varphi$-modules.
To ease notation, write~$\Gamma=\Gamma_K$ from now on. 
We choose
a topological generator $\gamma$ of $\Gamma$, and let $\Gamma_{\disc} :=
\langle \gamma \rangle;$ so $\Gamma_{\disc} \cong \Z$.
\begin{remark}
	\label{rem:discrete}
Since $\Gamma_{\disc}$ is dense in $\Gamma$, while a projective \'etale $\varphi$-module $M$ is complete with respect
to its canonical topology, in order to endow $M$ with the structure
of an \'etale $(\varphi,\Gamma)$-module, it suffices to equip $M$ with
a continuous action of $\Gamma_{\disc}$ (where we equip  $\Gamma_{\disc}$ 
with the topology induced on it by $\Gamma$).
\end{remark}

In order to study the properties of $\cX_{K,d}$ we will take advantage
of Remark~\ref{rem:discrete}. Accordingly, we now consider the moduli stack of 
projective \'etale $\varphi$-modules of rank~$d$ equipped with
a semi-linear action of $\Gamma_{\disc}$. 
We don't introduce particular
notation for this stack, since the following proposition
identifies it with a fixed point stack $\cR_{K,d}^{\Gamma_{\disc}}$, which
we now define.

There is a canonical action of~$\Gamma_{\disc}$ on~$\cR_d$ (that is, a canonical morphism $\gamma:\cR_d\to\cR_d$):
 if~$M$ is an object of~$\cR_d(A)$, then $\gamma(M)$ is given
by $\gamma^*M:=\A_{K,A}\otimes_{\gamma,\A_{K,A}}M$. (Note that this is naturally a
$\varphi$-module, because the action of~$\gamma$ on~$\A_{K,A}$ commutes
with~$\varphi$.)  Then we set 
\[\cR_d^{\Gamma_{\disc}}:=
\cR_d\underset{\Delta,\cR_d\times\cR_d,\Gamma_\gamma}{\times}\cR_d, \]where
 $\Delta$ is the diagonal of~$\cR_d$ and $\Gamma_\gamma$ is the
graph of~$\gamma$, so that $\Gamma_\gamma(x)=(x,\gamma(x))$. 


\begin{prop}
	\label{prop:Gamma-disc modules as fixed points; etale case}
        The moduli stack of projective
	\'etale $\varphi$-modules of rank~$d$ equipped with a semi-linear action
	of $\Gamma_{\disc}$ is isomorphic to the fixed point stack
	$\cR_d^{\Gamma_{\disc}}.$
\end{prop}
\begin{proof}Using the usual construction of the 2-fibre product, we
  see that $\cR_d^{\Gamma_\disc}$ consists of tuples
  $(x,y,\alpha,\beta)$, with $x,y$ being objects of~$\cR_d$, and
  $\alpha:x\isoto y$ and~$\beta:\gamma(x)\isoto y$ being
  isomorphisms. This is equivalent to the category fibred in groupoids
  given by pairs~$(x,\iota)$ consisting of an object~$x$ of~$\cR_d$
  and an isomorphism~$\iota:\gamma(x)\isoto x$. Thus an object
  of~$\cR_d^{\Gamma_\disc}(A)$ is a projective
  \'etale $\varphi$-module of  rank~$d$ with $A$-coefficients~$M$,
  together with an isomorphism of $\varphi$-modules
  $\iota:\gamma^*M\isoto M$; but this isomorphism is precisely the
  data of a semi-linear action of~$\Gamma_\disc=\langle\gamma\rangle$
  on~$M$, as required.
\end{proof}
Since~$\cR_d$ is an Ind-algebraic stack, so
is~$\cR_d^{\Gamma_\disc}$. More precisely, we have the following
lemma.

\begin{lem}\label{lem: good properties R Gamma disc}$\cR^{\Gamma_{\disc}}_d$ is a limit preserving
  Ind-algebraic stack, whose diagonal is representable by algebraic
  spaces, affine, and of finite presentation.
  \end{lem}
\begin{proof}
The description of $\cR^{\Gamma_{\disc}}_d$ as a fibre product
shows that the forgetful morphism $\cR^{\Gamma_{\disc}}_d \to \cR_d$
(given by forgetting the $\Gamma_{\disc}$-action)
is representable by algebraic spaces, affine, and of finite presentation,
since these properties hold for the diagonal of $\cR_d$,
by Corollary~\ref{cor: R good properties all K}.
The claim of the lemma is now seen to follow via another 
application 
of Corollary~\ref{cor: R good properties all K},
together with Corollary~\ref{cor:deducing properties}.
\end{proof}

Restricting the $\Gamma$-action on an \'etale $(\varphi,\Gamma)$-module
to $\Gamma_{\disc}$ (and taking into account Proposition~\ref{prop:Gamma-disc modules as fixed points; etale case}),
we obtain a morphism $\cX_{K,d} \to \cR_d^{\Gamma_{\disc}},$
which by Remark~\ref{rem:discrete} is fully faithful. 
Thus $\cX_{K,d}$ may be regarded as a substack of
$\cR_d^{\Gamma_{\disc}}$.

As already noted, 
the first step in proving that
$\cX_{K,d}$ is a Noetherian formal algebraic stack is to
show that it is an Ind-algebraic stack.
Although $\cX_{K,d}$ is a substack of the Ind-algebraic stack
$\cR_d^{\Gamma_{\disc}}$,
it is not a closed substack, but should rather be thought of as a
certain formal neighbourhood of $\cX_{d,\red}$ in $\cR_d^{\Gamma_{\disc}}$ (see Remark~\ref{rem: rank one X not closed substack}),
and so even this statement will require additional work to prove. We
begin with the following lemma, which allows us to reduce to the case
that~$K$ is basic.

\begin{lem}
  \label{lem: cartesian diagram for X basic} We have a 2-Cartesian
  diagram
  \[ \xymatrix{\cX_{K,d}\ar[r]\ar[d] &
      \cX_{\Kbasic,d[K:\Kbasic]}\ar[d]\\
      \cR_{K,d}^{\Gamma_{\disc}}\ar[r]&\cR_{\Kbasic,d[K:\Kbasic]}^{\Gamma_{\disc}}} \]where
  the horizontal arrows are the natural maps \emph{(}forgetting the
  $\A_K$-module structure\emph{)}, and the vertical arrows are the
  monomorphisms given by restricting the $\Gamma$-action
  to~$\Gamma_{\disc}$. 
\end{lem}
\begin{proof}
Unwinding the definitions, we need to show that if~$A$ is an
$\cO/\varpi^a$-algebra for some~$a\ge 1$, and~$M$ is an \'etale
$\varphi$-module over $\A_{K,A}$ with a semi-linear action
of~$\Gamma_{\disc}$, then the action of~$\Gamma_{\disc}$ extends to a
continuous action of~$\Gamma$ if and only if the same is true of~$M$
regarded as a module $\A_{\Kbasic,A}$. Since~$\A_{K,A}$ is free of finite
rank over~$\A_{\Kbasic,A}$, this is clear (for example, because it
follows from Lemma~\ref{lem:lattice properties} that the set of
lattices in~$M$ regarded as an $\A_{K,A}$-module is cofinal in the set of
lattices in~$M$ regarded as an $\A_{\Kbasic,A}$-module). 
\end{proof}

As a consequence of Lemmas~\ref{lem: relative representability for R
  over Kbasic} and~\ref{lem: cartesian diagram for X basic}, we can
deduce some properties of~$\cX_d$ from the corresponding properties
of~$\cR_d$.

\begin{prop}
  \label{prop: properties of X pulled back from R}
  \leavevmode
  \begin{enumerate}
  \item The morphism $\cX_{K,d}\to \cX_{\Kbasic,d[K:\Kbasic]}$ is
    representable by algebraic spaces, affine, and of finite presentation.
  \item The diagonal of $\cX_{K,d}$ is representable by algebraic
  spaces, affine, and of finite presentation.
  \end{enumerate}
\end{prop}
\begin{proof}To make the proof easier to read, we write~$\cR_K$
  for~$\cR_{K,d}$ and~$\cR_{\Kbasic}$ for~$\cR_{\Kbasic,d[K:\Kbasic]}$.
We factor the morphism
$\cR_{K}^{\Gamma_{\disc}}\to
  \cR_{\Kbasic}^{\Gamma_{\disc}}$
  as
  \begin{multline*}
\cR_{K}^{\Gamma_{\disc}} :=
\cR_{K}\underset{\Delta,\cR_{K}\times\cR_{K},\Gamma_\gamma}{\times}\cR_{K}
\to
\cR_{K}\underset{\Delta,\cR_{\Kbasic}\times\cR_{\Kbasic},\Gamma_\gamma}{\times}\cR_{K}
\\
\to
\cR_{\Kbasic}\underset{\Delta,\cR_{\Kbasic}\times\cR_{\Kbasic},\Gamma_\gamma}{\times}\cR_{\Kbasic}
=:
  \cR_{\Kbasic}^{\Gamma_{\disc}}.
\end{multline*}
	By Lemma~\ref{lem: relative representability for R over
    Kbasic}, the morphism
$\cR_{K}\to \cR_{\Kbasic}$
  is representable by algebraic spaces, affine, and of finite presentation.
Thus so is the second arrow in the preceding displayed expression.
The first arrow is a base-change of a product of two copies
of the diagonal 
$$\cR_{K} \to \cR_{K} \times_{\cR_{\Kbasic}} 
\cR_{K}.$$
As in the proof of Lemma~\ref{lem: abstract pullback diagonals lemma},
we deduce from 
Lemma~\ref{lem: relative representability for R over Kbasic}
that this diagonal 
  is representable by algebraic spaces, affine, and of finite presentation.
  Thus so is the first arrow in the expression,
  and hence so is the morphism
$\cR_{K}^{\Gamma_{\disc}}\to
  \cR_{\Kbasic}^{\Gamma_{\disc}}$.

It follows from Lemma~\ref{lem:
  cartesian diagram for X basic} that
 $\cX_{K,d}\to \cX_{\Kbasic,d[K:\Kbasic]}$ is then
 also representable by algebraic spaces, affine,
 and of finite presentation.
The claimed properties of the diagonal of~$\cX_{K,d}$ follow from
the corresponding properties of the diagonal
of~$\cR^{\Gamma_{\disc}}_{K,d}$ proved in Lemma~\ref{lem: good properties R Gamma disc},
together with the fact that $\cX_{K,d} \to \cR_{K,d}^{\Gamma_{\disc}}$
is a monomorphism.
\end{proof}
From now on we will typically drop~$K$ from the notation, simply writing
$\cX_d$, $\cR_d$ and so on. We conclude this section by showing that $\cX_d$ is limit preserving,
using the material
on lattices and continuity developed in Appendix~\ref{app: Tate modules and continuity}.
\begin{lem}
  \label{lem: gamma on T}We have ~$\gamma(T)-T\in
  p\A_{K,A}+T^2\A^+_{K,A}$. If~$K$ is basic then $\gamma(T)-T\in
  (p,T)T\A^+_{K,A}$.
\end{lem}
\begin{proof} Given the definitions of $\A_{K,A}$ and $\A^+_{K,A}$
in terms of $\A_K$ and $\A^+_K$, it suffices to prove the lemma
for these latter rings (i.e.\ in the case when $A =  \Z_p$).
We then reduce modulo~$p$, and
  write~$\gamma(T)=\sum_{i=n}^{\infty}a_iT^i$, with
  ~$a_i\in k_\infty$ and~$a_{n}\ne 0$. Since the action of~$\gamma$ is continuous for
  the ~$T$-adic topology, we see in particular that for~$m>0$
  sufficiently large, we have $\gamma(T)^m=\gamma(T^m)\in T
  k_\infty[[T]]$. Since~$\gamma(T)^m=a_{n}^mT^{nm}+\dots$, we see that~$n>0$. 
Since~$\gamma^{p^s}\to 1$ as
  $s\to \infty$, we see that~$n=1$ and that~$a_1^{p^s}=1$ for all sufficiently
  large~$s$, which implies that~$a_1=1$.

  If~$K$ is basic then 
  we have chosen~$T$ so that $\gamma(T)\in T\A_{K,A}^+$, so
  that $\gamma(T)-T\in (p,T)T\A^+_{K,A}$ by the above.
\end{proof}
If~$K$ is basic, it follows from Lemma~\ref{lem: gamma on T}
that~\eqref{eqn: condition on gamma semi-linear} holds, so that by
Lemma~\ref{lem:gamma minus 1 T quasi linear} we can use the material
on $T$-quasi-linear morphisms developed in Appendix~\ref{app: Tate
  modules and continuity} to study the action of~$\gamma$.

\begin{lem}
  \label{lem:X is limit preserving}$\cX_d$ is limit preserving.
\end{lem}
\begin{proof}
By Proposition~\ref{prop: properties of X pulled back from R} it suffices to prove
this in the case that~$K$ is basic (since a morphism which is 
of finite presentation is in particular limit preserving). Since~$\cX_d\hookrightarrow\cR_d^{\Gamma_\disc}$ is fully faithful, and~$\cR_d^{\Gamma_\disc}$ is limit
  preserving by Corollary~\ref{cor: basic properties of C and R p adic stacks}, it suffices to prove that~$\cX_d$ is limit preserving
  on objects. Since~$\cR_d^{\Gamma_\disc}$ is limit preserving on
  objects, we are reduced to showing that if~$T=\varprojlim_iT_i$ is a
  limit of affine schemes, and  $T_{i_0}\to\cR_d^{\Gamma_\disc}$ is a morphism with
  the property that the composite \[T\to T_{i_0}\to\cR_d^{\Gamma_\disc}\]
    factors through~$\cX_d$, then for some~$i\ge i_0$, the composite \[T_i\to T_{i_0}\to\cR_d^{\Gamma_\disc}\]
    factors through~$\cX_d$.
    This follows from the equivalence of conditions~(1)
    and~(7) 
    of Lemma~\ref{lem:testing continuity on M mod T}, together with Lemma
    ~\ref{lem:quasi-linear limit preserving}.
%
\end{proof}

\section{Weak Wach modules}\label{subsec: Wach modules}
In this section we introduce the notion of a weak Wach module of
height at most~$h$ and level at most~$s$. These will play a purely technical
auxiliary role for us, and will be used only in order to show
that~$\cX_d$ is an Ind-algebraic stack.  
We suppose throughout this
section that~$K$ is basic.
\begin{remark}
	\label{rem:testing continuity on M mod u}
        By Lemma~\ref{lem:testing continuity on M mod T}, if $A$ is an $\cO/\varpi^a$-algebra for some~$a\ge 1$, and
        $\gM$ is a rank~$d$ projective $\varphi$-module of
        $T$-height $\leq h$ over $A$, such that~$\gM[1/T]$ is equipped
        with a semi-linear action of $\Gamma_{\disc}$, then this
        action extends to a 
        continuous action of~$\Gamma$ if and only if for some~$s\ge 0$ we have $(\gamma^{p^s}-1)(\gM)\subseteq
        T\gM$.
\end{remark}

\begin{df}
	\label{def:weak Wach modules}Suppose that~$K$ is basic.
	A {\em rank $d$ projective weak Wach module of $T$-height
    $\leq h$} \index{weak Wach module}
	with $A$-coefficients is a rank~$d$ projective $\varphi$-module $\gM$ 
	over $\A^+_{K,A}$, 
        which is of $T$-height $\leq h$, 
        such that~$\gM[1/T]$ is equipped
	with a semi-linear action of $\Gamma$ which is furthermore
        continuous 
	when $\gM[1/T]$ is endowed with its canonical topology  (see
        Remark~\ref{rem:topologies on M}).

	If~$s\ge 0$, then we say that~$\gM$ has
        \emph{level $\le s$} \index{level (of weak Wach module)} if $(\gamma^{p^s}-1)(\gM)~\subseteq~T\gM$.
	Remark~\ref{rem:testing continuity on M mod u} shows
	that any projective weak Wach module is of level $\leq s$
	for some $s \geq 0$.
\end{df}
A special role in the classical theory of $(\varphi,\Gamma)$-modules is played by the weak Wach modules of
  level~$0$.  However, this will not be important for us.



\begin{df}
	We let $\cWW_{d,h}$ denote the moduli stack of rank $d$ \index{$\cWW_{d,h}$}
       	projective weak Wach modules of $T$-height $\leq h$. 
        (That this
        is an \emph{fpqc} stack follows as in Definition~\ref{defn:
          Xd}.) We let~$\cWW_{d,h,s}$ denote the substack of rank $d$ \index{$\cWW_{d,h,s}$}
       	projective weak Wach modules of $T$-height $\leq h$ and
        level~$\le s$.  
\end{df}

We recall from Section~\ref{subsec: moduli stacks of phi modules}
that there is a $p$-adic formal
algebraic stack
$\cC_{d,h}$ classifying rank $d$ projective $\varphi$-modules of
$T$-height at most~$h$.
We consider the fibre product
$\cR_d^{\Gamma_{\disc}}\times_{\cR_d}\cC_{d,h}$, where the
map~$\cR_d^{\Gamma_{\disc}}\to\cR_d$ is the canonical morphism given by
forgetting
the~$\Gamma_{\disc}$ action. By Proposition~\ref{prop:Gamma-disc
  modules as fixed points; etale case}, this is the moduli stack of rank $d$ projective $\varphi$-modules~$\gM$ of
$T$-height at most~$h$, equipped with a semi-linear action
of~$\Gamma_{\disc}$ on~$\gM[1/T]$. 
\begin{lem}
  \label{lem: R Gamma disc times C is p adic formal
    algebraic}$\cR_d^{\Gamma_{\disc}}\times_{\cR_d}\cC_{d,h}$ is a
  $p$-adic formal algebraic stack of finite presentation over~$\Spf\cO$.
\end{lem}
\begin{proof}The map $\cR_d^{\Gamma_{\disc}}\to\cR_d$ is  representable by algebraic spaces
	and
	of finite presentation, being a base-change of the diagonal
        of~$\cR_d$, which is representable by algebraic spaces
	and
	of finite presentation by Corollary~\ref{cor: basic properties of C and R p adic stacks}. 
  Again by  Corollary~\ref{cor: basic properties of C and R p adic stacks}, $\cC_{d,h}$ is a $p$-adic formal
  algebraic stack of finite presentation, and thus
  $\cR_d^{\Gamma_{\disc}}\times_{\cR_d} \cC_{d,h}$ is a $p$-adic
  formal algebraic stack of finite presentation, as claimed.
  \end{proof}
Restricting the $\Gamma$-action on a weak Wach module
to $\Gamma_{\disc}$, 
we obtain a morphism $\cWW_{d,h} \to \cR_d^{\Gamma_{\disc}}\times_{\cR_d}\cC_{d,h}$,
which, by the evident analogue for weak Wach modules
of Remark~\ref{rem:discrete}, is fully faithful. 
Thus $\cWW_{d,h}$ may be regarded as a substack of
$\cR_d^{\Gamma_{\disc}}\times_{\cR_d}\cC_{d,h}$. The following
proposition records the basic properties of the stacks~$\cWW_{d,h}$
and~$\cWW_{d,h,s}$.

\begin{prop}
	\label{prop:inductive description of W}\leavevmode
        \begin{enumerate}
        \item For $s \geq 0,$ the morphism
	$$\xymatrix{\cWW_{d,h,s}
          \ar
          [r]
          & \cR_d^{\Gamma_{\disc}}\times_{\cR_d}\cC_{d,h}\\ }$$ is a
        closed immersion of finite presentation.
        In particular, each of the stacks $\cWW_{d,h,s}$ is a $p$-adic
          formal algebraic stack of finite presentation over
          $\Spf \cO$.
        \item If $s'\ge s$, then the canonical monomorphism $
          \cWW_{d,h,s} \hookrightarrow \cWW_{d,h,s'}$ is a closed
          immersion of finite presentation.
        

       \item 
         The canonical morphism $\varinjlim_s \cWW_{d,h,s} \to \cWW_{d,h}$
         is an isomorphism. In particular, $\cWW_{d,h}$
         is an Ind-algebraic stack.  
      \end{enumerate}

\end{prop}
\begin{proof} 	Since  $\cR_d^{\Gamma_{\disc}}\times_{\cR_d}
	\cC_{d,h}$ classifies
	$\varphi$-modules $\gM$ that are
	projective of rank $d$ and
	of $T$-height $\leq h$, which are endowed with a $\Gamma_{\disc}$-action
	on the underlying \'etale $\varphi$-module, in order to prove~(1),
	we must show that
        the condition that~$(\gamma^{p^s}-1)(\gM)\subseteq T\gM$ is a
        closed condition, and is determined by finitely many
        equations. It suffices to check this after pulling back via an
        arbitrary morphism 
          $\Spec A\to \cR_d^{\Gamma_{\disc}}\times_{\cR_d}
	\cC_{d,h}$,
where~$A$ is an $\cO/\varpi^a$-algebra for some~$a\ge 1$. 


This morphism gives rise to a projective $\A^+_{K,A}$ $\varphi$-module $\gM$,
such that $M:=  \gM[1/T]$ is \'etale, and is furthermore endowed with
a semi-linear $\Gamma_{\disc}$-action.
As already stated above, we must now check that the condition
$(\gamma^{p^s}-1)(\gM) \subseteq T\gM$  is a finitely presented closed condition
(in the sense that it holds after replacing $\gM$  by  $\gM_B := \A^+_{K,B}
\otimes_{\A^+_{K,A}} \gM$ for an $A$-algebra $B$ if and only if
$\Spec B \to \Spec A$ factors through a certain finitely presented closed
subscheme of~$\Spec A$).

Note firstly that for~$n$ sufficiently large, we have
$(\gamma^{p^s}-1)(T^n\gM)\subseteq T\gM$. Indeed,
writing \[\gamma^{p^s}-1=((\gamma-1)+1)^{p^s}-1\]and expanding via the
binomial theorem, it suffices to show that for ~$n$ sufficiently large, we have
$(\gamma-1)^m(T^n\gM)\subseteq T\gM$ for all $0\le m\le p^s$; and this
is immediate from Lemma~\ref{lem: how f behaves on lattices}.

We next 
choose a finitely generated projective $\A_{K,A}^+$-module ~$\gN$ such
that ~$\gF:=\gM\oplus\gN$ is free, and write~$N:=\gN[1/T]$,
$F:=\gF[1/T]$. Extend the morphism~$\gamma^{p^s}-1:M\to M$ to the
morphism~$f:F\to F$ given by
$(\gamma^{p^s}-1,0):M\oplus N\to M\oplus N$, so that the condition
that $(\gamma^{p^s}-1)(\gM)\subseteq T\gM$ is equivalent to asking that
$f(\gF)\subseteq T\gF$.
By our choice of~$n$ above,  the morphism $\gF\to F/T\gF$
induced by~$f$ factors through a morphism $\gF/T^n\gF\to F/T\gF$;
since~$\gF$ is finitely generated, after enlarging~$n$ if necessary,
it factors through
$\gF/T^n\gF\to T^{-n}\gF/T\gF$. 
The vanishing of this morphism is obviously a closed condition; indeed
it is given by the vanishing of finitely many matrix entries, so is
closed and of finite presentation, as required.
      Since~$\cR_d^{\Gamma_{\disc}}\times_{\cR_d}\cC_{d,h}$ is a
      $p$-adic formal algebraic stack of finite presentation over
      $\Spf \cO$, it follows that so is each~$
      \cWW_{d,h,s}$. That~$\cWW_{d,h,s} \hookrightarrow \cWW_{d,h,s'}$
      is a closed immersion of finite presentation follows immediately
      from a consideration of the composite
      \[\cWW_{d,h,s} \hookrightarrow \cWW_{d,h,s'}\hookrightarrow
        \cR_d^{\Gamma_{\disc}}\times_{\cR_d}\cC_{d,h},\]as we have
      just shown that both the composite and the second morphism are
      closed immersions of finite presentation. Thus we have
      proved~(1) and~(2).

For~(3), we need to show that  every morphism $\Spec A\to\cWW_{d,h}$,
where~$A$ is an $\cO/\varpi^a$-algebra for some~$a\ge 1$, factors
through~$\cWW_{d,h,s}$  for~$s$ sufficiently large. As
we already noted in Definition~\ref{def:weak Wach modules}, this follows from
Remark~\ref{rem:testing continuity on M mod u}. Since
each~$\cWW_{d,h,s}$ is an Ind-algebraic stack, so is~$\cWW_{d,h}$. 
\end{proof}

\section{$\cX_{\lowercase{d}}$ is an Ind-algebraic stack}
\label{subsec: X is Ind algebraic}

Continue to assume that~$K$ is basic. By definition, we have a $2$-Cartesian diagram
\numequation
\label{eqn:C to R square}
\xymatrix{\cWW_{d,h} \ar[r]\ar[d] &
	\cR_d^{\Gamma_{\disc}}\times_{\cR_d}\cC_{d,h} \ar[d] \\
\cX_d \ar[r] & \cR_d^{\Gamma_{\disc}} }
\end{equation}
If $h' \geq h$ then
the closed immersion $\cC_{d,h} \hookrightarrow \cC_{d,h'}$ is
compatible with the morphisms from each of its source and target
to $\cR_d$, 
and so 
we obtain a closed immersion 
\numequation
\label{eqn:closed immersion of Wachs}
\cWW_{d,h} \hookrightarrow \cWW_{d,h'}.
\end{equation}

By construction,
the morphisms $\cWW_{d,h}\to\cX_d$  are compatible,
as $h$ varies,
with the closed immersions~(\ref{eqn:closed immersion of Wachs}). 
Thus we also obtain a morphism
\numequation
\label{eqn:Ind W to X}
\varinjlim_h \cWW_{d,h} \to \cX_d.
\end{equation}
Roughly speaking, we will prove that~$\cX_d$ is an Ind-algebraic stack
by showing that it is the ``scheme-theoretic image'' of the morphism
$\varinjlim_h \cWW_{d,h} \to \cR_d^{\Gamma_{\disc}}$ induced
by~\eqref{eqn:Ind W to X}. More precisely, choose $s \geq 0,$ and consider the composite
        \numequation
        \label{eqn:W to R composite}
        \cWW_{d,h,s} \to \cWW_{d,h} \to \cX_d \to \cR_d^{\Gamma_{\disc}}.
        \end{equation}
        This admits the alternative factorization
        $$\cWW_{d,h,s} \to \cWW_{d,h} \to 	\cR_d^{\Gamma_{\disc}}\times_{\cR_d}\cC_{d,h} 
        \to \cR_d^{\Gamma_{\disc}}.$$
        Proposition~\ref{prop:inductive description of W} shows
        that the composite of the first two arrows is a closed 
        embedding of finite presentation,
        while  Corollary~\ref{cor: basic properties of C and R p adic stacks} 
        shows that the third arrow is representable by algebraic spaces, proper, 
        and of finite presentation. 
        Thus~(\ref{eqn:W to R composite}) is representable by
        algebraic spaces, proper, and of
        finite presentation.   

	Fix an integer $a \geq 1,$ and write $\cWW^a_{d,h,s} := \cWW_{d,h,s}
	\times_{\Spf \cO} \Spec \cO/\varpi^a$.   
	Proposition~\ref{prop:inductive description of W}
	shows that $\cWW_{d,h,s}$ is a $p$-adic formal
	algebraic stack of finite presentation over $\Spf \cO$,
	and so $\cWW^a_{d,h,s}$ is an algebraic stack, 
	and a closed substack of $\cWW_{d,h,s}$.

\begin{df}  We let $\cX^a_{d,h,s}$ denote the \index{$\cX^a_{d,h,s}$}
scheme-theoretic image of the composite 
	\numequation
	\label{eqn:another composite}
	\cWW^a_{d,h,s} \hookrightarrow \cWW_{d,h,s}
	\buildrel \text{(\ref{eqn:W to R composite})} \over
	\longrightarrow \cR_d^{\Gamma_{\disc}},
        \end{equation}
defined via the formalism of scheme-theoretic images for morphisms of Ind-algebraic
stacks developed in Appendix~\ref{app: formal algebraic stacks}.
\end{df}

More concretely,
	since $\cR_d^{\Gamma_{\disc}}$ is an Ind-algebraic stack, 
        constructed as the $2$-colimit of a directed system of algebraic
        stacks whose transition morphisms are closed immersions,
	the morphism~\eqref{eqn:another composite},
	which is representable by algebraic spaces, proper, and of finite presentation,
	factors through a closed algebraic substack $\cZ$ of
	$\cR_d^{\Gamma_{\disc}}.$   We then define $\cX^a_{d,h,s}$
	to be the 
        scheme-theoretic image of
	$\cWW^a_{d,h,s}$ in $\cZ$.   Note that
	$\cX^a_{d,h,s}$ 
        is a closed algebraic substack of $\cR_d^{\Gamma_{\disc}}$,
        and is independent of the choice of~$\cZ$.

\begin{remark}
As already observed above, the morphism~\eqref{eqn:another composite} 
factors through $\cX_d$. However, since at this point we don't know that $\cX_d$
is Ind-algebraic, we can't directly define a scheme-theoretic image of 
$\cW^a_{d,h,s}$ in $\cX_d$.  It might be possible to do this using
the formalism of \cite{EGstacktheoreticimages}; since~\eqref{eqn:another composite} 
is proper, this scheme-theoretic image would then coincide with $\cX^a_{d,h,s}$.
We don't do this here; but we do prove somewhat more directly,
in Lemma~\ref{lem: scheme theoretic images factor through X} below,
that $\cX^a_{d,h,s}$ is a substack of~$\cX_d$.
\end{remark}

        As in Definition~\ref{defn: lattice}, a \emph{lattice} $\gM$ \index{lattice}
        in a projective \'etale $\varphi$-module~$M$ is a finitely
        generated ~$\A_{K,A}^+$-submodule of~$M$ whose $\A_{K,A}$-span
        is~$M$. Note that~$\gM$ is not assumed to be projective.
	\begin{lem}
		\label{lem:Artinian points of scheme theoretic images}
		Suppose that $M$ is a projective \'etale $\varphi$-module of rank $d$
		over a finite type Artinian $\cO/\varpi^a$-algebra $A$,
		and that~$M$ is endowed with an action of $\Gamma_{\disc}$,
		such that the corresponding morphism $\Spec A \to
		\cR_d^{\Gamma_{\disc}}$ 
		factors through $\cX^a_{d,h,s}$.
		Then $M$ contains a $\varphi$-invariant
		lattice 
                $\gM$ of $T$-height $\leq h$,
		such that $(\gamma^{p^s}-1) (\gM) \subseteq
		T\gM.$
	\end{lem}
	\begin{proof}Since an Artinian ring is a direct product of
          Artinian local rings, it suffices to treat the case that~$A$
          is local.
          Let the residue field of~$A$ be~$\F'$, a finite
          extension of~$\F$, and write~$\cO'$ for the ring of integers
          in the compositum of~$E$ and~$W(\F')[1/p]$, so that~$\cO'$
          has residue field~$\F'$; note that $A$ is naturally on $\cO'$-algebra. 

          The projective
          \'etale $\varphi$-module $M$ is in fact free of rank~$d$,
          since $\A_{K,A}$ is again a local ring, and we fix an (ordered)
          basis for $M$ as an $\A_{K,A}$-module.
          Write~$M_{\F'} := \A_{K,\F'} \otimes_{\A_{K,A}} M$; the $\A_{K,A}$-basis
          of $M$ gives rise to a corresponding $\A_{K,\F'}$-basis of~$M_{\F'}$.
          Let~$\cC_{\cO'}$ be the category of Artinian local
          $\cO'/\varpi^a$-algebras for which the structure map induces
          an isomorphism on residue fields.  
By a \emph{lifting} of~$M_{\F'}$ to an object ~$\Lambda$
         of~$\cC_{\cO'}$, we mean a triple consisting of an \'etale
 $\varphi$-module~$M_\Lambda$ which is free of rank~$d$, a choice of (ordered)
 $\A_{K,\Lambda}$-basis of~$M_\Lambda$, and an isomorphism $M_\Lambda\otimes_\Lambda
 \F'\cong M_{\F'}$ of \'etale $\varphi$-modules which takes the chosen basis
 of~$M_\Lambda$ to the fixed basis of~$M_{\F'}$. Let~$D$ be the functor
          $\cC_{\cO'}\to\Sets$ taking~$\Lambda$ to the set of isomorphism
          classes of liftings of~$M_{\F'}$ to an \'etale $\varphi$-module~$M_\Lambda$
          with $\Lambda$-coefficients,
	  endowed with an action of~$\Gamma_{\disc}$ lifting
	  that on~$M_{\F'}$.
          Note that $A$  is an object of $\cC_{\cO'}$, and that
our originally chosen basis for $M$ is classified by a 
          continuous morphism $R  \to A.$

          The functor~$D$ is pro-representable by an object~$R$
          of~$\pro\cC_{\cO'}$, by the same argument as in the
          proof 
          of~\cite[Prop.\ 5.3.6]{EGstacktheoreticimages}: namely, 
          by Grothendieck's representability theorem, it suffices to
          prove the compatibility of~$D$ with fibre products
          in~$\cC_{\cO'}$, which is obvious. The universal lifting
          $M_R$ gives a morphism $\Spf R\to\cR_d^{\Gamma_{\disc}}$. 
          The composite  $\Spec A \to \Spf R \to \cR_d^{\Gamma_{\disc}}$ 
          is of course just the morphism that classifies~$M$.

We let~$D'$ denote the
          subfunctor of~$D$ consisting of those lifts for which there
          is a lattice  ~$\gM_\Lambda$ 
          of $F$-height~$\le h$ and level~$\le s$,
          with~$\gM_\Lambda[1/T]=M$. (More precisely, we require that
          $T^h\gM_\Lambda\subseteq (1\otimes\varphi_{M_\Lambda})(\varphi^*\gM_\Lambda)$  and
          $(\gamma^{p^s}-1)(\gM_\Lambda)\subseteq T\gM_\Lambda$.) We claim that the functor~$D'$ is
          pro-representable by a quotient~$S$ of~$R$. To see this, it
          is again enough to show that~$D'$ preserves fibre products,
          and in turn it is enough to show that the property of an
          \'etale $\varphi$-module~$M$ with an action
          of~$\Gamma_{\disc}$ admitting a ~$\varphi$-stable
          lattice~$\gM$ of $T$-height~$\le h$, and such that
          $(\gamma^{p^s}-1)(\gM)\subseteq T\gM$, is stable under
          taking direct sums and subquotients. This is obvious for
          direct sums, and the case of subquotients follows as
          in~\cite[Lem.\ 5.3.10]{EGstacktheoreticimages} (which is the
          same result without the conditions
          on~$\Gamma_{\disc}$). More precisely, once checks that if we
          have a short exact sequence $0\to M'\to M\to M''\to 0$ of
          \'etale $\varphi$-modules with $\Lambda$-coefficients and
          actions of~$\Gamma_{\disc}$, and $\gM$ is a lattice of the
          appropriate kind in~$M$, then the kernel and image of the
          map $\gM\to M''$ give the appropriate lattices~$\gM'$,
          $\gM''$ in~$M'$, $M''$ respectively. (The properties that
          $(\gamma^{p^s}-1)(\gM')\subseteq T\gM'$ and
          $(\gamma^{p^s}-1)(\gM'')\subseteq T\gM''$ follow from a
          short diagram chase, using that the composite
          $\gM\stackrel{(\gamma^{p^s}-1)}{\to}M\to M/T\gM$ vanishes by
          hypothesis.)

 By the definition of~$D'$, the statement of the lemma is equivalent to
          the statement that the 
map $R\to A$ factors
          through~$S$. Let $X=\cW_{d,h,s}\times_{\cR_d^{\Gamma_{\disc}}}\Spf
          R$, a formal algebraic space, and let~$\Spf T$ be the
          scheme-theoretic image of the morphism $X\to\Spf
          R$. Since the morphism~$\Spec A\to\cR_d^{\Gamma_{\disc}}$
          factors through $\cX^a_{d,h,s}$ by hypothesis, the morphism
          $\Spec A\to\Spf R$ factors through~$\Spf T$, so it is enough
to show that~$\Spf T$ is a closed formal subscheme of~$\Spf S$.
By Lemma~\ref{lem: criterion for Artin to map to scheme theoretic
    image}, 
it in turn suffices to          show that if~$A'$ is a finite
          type Artinian local~$R$-algebra for which $R\to A'$ factors
          through a discrete quotient of~$R$, and for which the
          morphism $\cW_{d,h,s}\times_{\cR_d^{\Gamma_{\disc}}}\Spec
          A'\to\Spec A'$ admits a section, then~$R\to A'$ factors
          through~$S$.

          By replacing~$R$
          with~$R\otimes_{W(\F')}\kappa(A')$ we can reduce to the case
          that~$A'$ has residue field~$\F'$ (\emph{cf.}\ \cite[Cor.\ 5.3.18]{EGstacktheoreticimages}), so that $R\to A'$ factors
          through~$S$ if and only if the \'etale $\varphi$-module
          corresponding to $\Spec A'\to \cR_{d}^{\Gamma_{\disc}}$
          admits a $\varphi$-stable lattice~$\gM_{A'}$ of
          $T$-height~$\le h$ for
          which~$(\gamma^{p^s}-1)(\gM_{A'})\subseteq T\gM_{A'}$. But
          a section to the morphism $\cW_{d,h,s}\times_{\cR_d^{\Gamma_{\disc}}}\Spec
          A'\to\Spec A'$ gives us a morphism $\Spec A'\to\cW_{d,h,s}$, and the
          corresponding weak Wach module provides the required lattice.
          \end{proof}
          

        \begin{lem}
  \label{lem: scheme theoretic images factor through X}
  Each~$\cX^a_{d,h,s}$ is a
  closed substack of~$\cX_d$.
\end{lem}
\begin{proof}
  Since $\cX^a_{d,h,s}$ is a closed substack of $\cR_d^{\Gamma_{\disc}}$,
  if it is a substack of $\cX_d$, it will in fact
  be a closed substack.
  Thus it suffices to show that $\cX^a_{d,h,s}$ is indeed a substack
  of~$\cX_d$.
  Since~$\cX^a_{d,h,s}$ is limit preserving, it is enough to check
  that if~$A$ is a finite type $\cO/\varpi^a$-algebra, then for any
  morphism $\Spec A\to\cX^a_{d,h,s}$, the composite morphism
  $\Spec A\to\cX^a_{d,h,s}\to\cR_d^{\Gamma_{\disc}}$ factors
  through~$\cX_d$.  Equivalently,
  if $M$ denotes the \'etale $\varphi$-module over~$A$, 
  endowed with a $\Gamma_{\disc}$-action,
  associated to the given point $\Spec A \to \cR_d^{\Gamma_{\disc}},$
  then we must show that the $\Gamma_{\disc}$-action
  on $M$ is continuous. 

  By Lemma~\ref{lem:testing continuity on M mod T}, to see that the
  $\Gamma_{\disc}$-action on $M$ is continuous, it suffices to produce
  a lattice $\gM \subseteq M$ such that
  $(\gamma^{p^s}-1)(\gM)\subseteq T\gM$. We will produce such a
  lattice by reduction to the Artinian case, as follows. Let
  $\{A_i\}_{i \in I}$ be the directed system of Artinian quotients of
  $A$.  Since $A$ is of finite type over $\cO/\varpi^a$, and so
  Noetherian, the natural map $A \to B := \prod_i A_i$ is injective.
  Lemma~\ref{lem:Artinian points of scheme theoretic images} shows
  that each base-changed module $M_{A_i}$ admits a $\varphi$-invariant
  lattice $\gM_i$ of $T$-height~$\le h$ satisfying the condition
  $(\gamma^{p^s}- 1)(\gM_i) \subseteq T \gM_i.$

  Write~$\gM_B:=\prod_i\gM_i\subset M_B$, and set
  $\gM:=M\cap\gM_B\subset M_B$. Since  $(\gamma^{p^s}- 1)(\gM_B)
  \subseteq T \gM_B$, we have  $(\gamma^{p^s}- 1)(\gM) \subseteq T
  \gM$, so it only remains to check that~$\gM$ is a lattice. To see
  this, let~$\gM'$ be a finite height lattice in~$M$ (i.e.\ a lattice
  in $M$ which is~$\varphi$-stable, and for which the cokernel
  of~$\Phi$ is killed by a power of~$T$; such a lattice exists
  by~\cite[Lem.\ 5.2.15]{EGstacktheoreticimages}). By
  Lemma~\ref{lem:lattice properties}~(3), it suffices to prove that
  there is an integer~$l\ge 0$ such that $T^l\gM'\subseteq
  \gM\subseteq T^{-l}\gM$. Accordingly, if for each~$i$ we
  let~$\gM'_i$ denoted the image of~$\gM_i\otimes_AA_i$ in~$M_{A_i}$,
  we need to show that we have $T^l\gM'_i\subseteq
  \gM_i\subseteq T^{-l}\gM_i'$.

  The existence of such an~$l$ is established in the course of the
  proof of~\cite[Lem.\ 5.3.14]{EGstacktheoreticimages}, although this
  is not explicitly recorded there. For the convenience of the reader,
  we recall the argument in our present setting. Increasing~$h$ if
  necessary, we may assume that~$\gM'$ has $T$-height at most~$h$. It
  follows that~$\gM_i'$ and~$\gM_i$ are each $\varphi$-stable lattices
  of height at most~$h$ in~$M_{A_i}$. We claim that we may
  take~$l=\lfloor (h+ap)/(p-1)\rfloor$.

  To see this, let~$j$ be minimal such that $T^j\gM_i\subseteq\gM_i'$;
  then we have
  \[\Phi_{\gM_i}(\varphi^*\gM_i)\subseteq \gM_i\subseteq
    T^{-j}\gM'_i\subseteq T^{-h-j}\Phi_{\gM'_i}(\varphi^*\gM'_i),\] so
  that $T^{h+j}\varphi^*\gM_i\subseteq\varphi^*\gM'_i$. It follows
  from~\cite[Lem.\ 5.3.13]{EGstacktheoreticimages} that $h+j>(j-a)p$,
  so that $j<(h+ap)/(p-1)$, and~$j\le l$ by definition.  Thus
  $T^l\gM_i\subseteq\gM_i'$, as claimed. Reversing the roles
  of~$\gM_i$ and~$\gM'_i$, we have $T^l\gM'_i\subseteq\gM_i$, and we
  are done.
\end{proof}    

We now prove our first key structural result for $\cX_d$, in the case when $K$ is basic.

\begin{prop}
	\label{prop:X is an Ind-stack basic case}
If~$K$ is basic, then the canonical morphism
	$\varinjlim \cX^a_{d,h,s} \to \cX_d$ is an isomorphism.
	Thus $\cX_d$ is an Ind-algebraic stack, and may in fact
	be written as the inductive limit of
        algebraic stacks of finite presentation, with the transition maps
        being closed immersions.
\end{prop}
\begin{proof}Note firstly that if $a'\ge a$, $h'\ge h$ and $s'\ge s$
  then the canonical morphism~$\cX^a_{d,h,s}\to\cX^{a'}_{d,h',s'}$ is a
  closed immersion by construction. By Lemma~\ref{lem: scheme
    theoretic images factor through X}, each~$\cX^a_{d,h,s}$ is a
  closed substack of~$\cX_d$, so it remains to show that any morphism
  $T\to\cX_d$ whose source is an affine scheme
  factors through some~$\cX^a_{d,h,s}$,
  or equivalently, that the closed immersion
  \numequation
  \label{eqn:closed immersion into test scheme}
  \cX^a_{d,h,s}\times_{\cX_d} T \to T
  \end{equation}
  is an isomorphism,
  for some choice of $h$ and $s$.

  Since~$\cX_d$ is
  limit preserving, by Lemma~\ref{lem:X is limit preserving}, we can
  reduce to the case where $T = \Spec A$ for a Noetherian
  $\cO/\varpi^a$-algebra~$A$.
  If $M$ denotes the \'etale $(\varphi,\Gamma)$-module corresponding
  to the morphism $\Spec A \to \cX_d$,
  then an application of \cite[Prop.~5.4.7]{EGstacktheoreticimages}
  shows that
  we may find a scheme-theoretically dominant morphism 
  $\Spec B \to \Spec A$ such that $M_B$ is free of rank $d$.
  If we show that the composite $\Spec B \to \Spec A \to \cX_d$
  factors through $\cX^a_{d,h,s}$ for some $h$ and $s$,
  then we see that the morphism $\Spec B \to \Spec A$
  factors through the closed subscheme
  $\cX^a_{d,h,s} \times_{\cX_d} \Spec A$ of $\Spec A.$ 
  Since $\Spec B \to \Spec A$ is scheme-theoretically dominant,
  this implies that~(\ref{eqn:closed immersion into test scheme})
  is indeed an isomorphism, as required.

  Since $M_B$ is free,
  we may choose a $\varphi$-invariant free lattice $\gM \subseteq M_B$,
  of height $\leq h$ for some sufficiently large value of $h$.
  Since the $\Gamma_{\disc}$-action on~$M$, and hence
  on~$M_B$, is continuous by assumption,
  Lemma~\ref{lem:testing continuity on M mod T}
  shows that $(\gamma^{p^s}-1)(\gM)\subseteq T\gM$ 
  for some sufficiently large value of $s$.
  Then $\gM$ gives rise to a $B$-valued point
  of $\cW^a_{d,h,s}$, whose image in $\cR^{\Gamma_{\disc}}_d$
  is equal to the \'etale $\varphi$-module $M_B$.
  Thus the morphism $B \to \cX_d$
  corresponding to $M_B$
  does indeed factor through $\cX^a_{d,h,s}$.
%
%
%
%
%
%
\end{proof}  

Finally, we drop our assumption that $K$ is basic.

\begin{prop}
	\label{prop:X is an Ind-stack}
Let~$K$ be an arbitrary finite extension of~$\Qp$. 
	Then $\cX_d$ is an Ind-algebraic stack, and may in fact
	be written as the inductive limit of
        algebraic stacks of finite presentation over~$\Spec\cO$, with the transition maps
        being closed immersions. Furthermore the diagonal of~$\cX_d$
        is representable by algebraic spaces, affine, and of finite presentation.
\end{prop}
\begin{proof}
  This is immediate from Propositions~\ref{prop: properties of X
    pulled back from R} and~\ref{prop:X is an Ind-stack basic case},
  together with Corollary~\ref{cor:deducing properties}.
\end{proof}
\section{Canonical actions and weak Wach modules}\label{subsec:
  canonical weak Wach}
We now explain an alternative perspective on some of the above
results, which gives more information about the moduli stacks of weak
Wach modules. 
The results of this section are only used
in Chapter~\ref{sec: the rank one case}, where we use them to
establish a concrete description of~$\cX_d$ in the case~$d=1$.
For each~$s\ge 0$, we write~$K_{\cyc,s}$
for the unique subfield of~$K_{\cyc}$ which is cyclic over~$K$ of
degree~$p^s$. The following lemma can be proved in exactly the same
way as Lemma~\ref{lem: Caruso Liu Galois action on Kisin} below.
\begin{lem}\label{lem: Caruso Liu Galois action on weak Wach}
Assume that~$K$ is basic.  For any fixed $a,h$,  there is a constant~$N(a,h)$ such that if $N\ge N(a,h)$, there is a
  positive integer~$s(a,h,N)$ with the property that for any finite type
  $\cO/\varpi^a$-algebra~$A$, any finite
  projective $\varphi$-module~$\gM$ over~$\A_{K,A}^+$  of $T$-height at most~$h$, and 
  any~$s\ge s(a,h,N)$, there is a unique continuous action of~ $G_{K_{\cyc,s}}$ on
  $\gMt:=\AAinf{A}\otimes_{\A_{K,A}^+}\gM$ which commutes with~$\varphi$
  and is semi-linear with respect to the natural action of~$G_{K_{\cyc,s}}$
  on~$\AAinf{A},$ with the additional property that for
  all~$g\in G_{K_{\cyc,s}}$ we have $(g-1)(\gM)\subset
  T^N\gMt$. 
\end{lem}
It is possible to
prove the following result purely in the world of $(\varphi,\Gamma)$-modules
(by following the proof of Lemma~\ref{lem: Caruso Liu Galois action on
  Kisin}, interpreting the existence of a semilinear action of~$\Gamma_{K_{\cyc,s}}$ in
terms of linear maps between twists of $\varphi$-modules, satisfying
certain compatibilities),
but we have found it more straightforward to argue with Proposition~\ref{prop: equivalences of
  categories to Ainf}.
 \begin{cor}\label{cor: Caruso Liu Galois action phi Gamma discrete version}
Assume that~$K$ is basic.  For any fixed $a,h$,  there is a constant~$N(a,h)$ such that if $N\ge N(a,h)$, there is a
  positive integer~$s(a,h,N)$ with the property that for any finite type
  $\cO/\varpi^a$-algebra~$A$, any finite
  projective $\varphi$-module~$\gM$ over~$\A_{K,A}^+$  of $T$-height at most~$h$, and 
  any~$s\ge s(a,h,N)$, there is a  unique semi-linear  action of~
  $\langle \gamma^{p^s}\rangle$ on
  $\gM$ which commutes with~$\varphi$,
  and  with the additional property that $(\gamma^{p^s}-1)(\gM)\subseteq
  T^N\gM$.

  In particular, this action gives $\gM[1/T]$ the structure of a projective \'etale $(\varphi,\Gamma_{K_{\cyc,s}})$-module.
\end{cor}
\begin{proof}By Lemma~\ref{lem: Caruso Liu Galois action on weak
    Wach}, there is a a unique continuous semi-linear action of~ $G_{K_{\cyc,s}}$ on
  $\gMt:=\AAinf{A}\otimes_{\A_{K,A}^+}\gM$ which commutes with~$\varphi$,
   with the additional property that for
  all~$g\in G_{K_{\cyc,s}}$ we have $(g-1)(\gM)\subset
  T^N\gMt$. Write~$M:=\A_{K,A}\otimes_{\A_{K,A}^+}\gM$, $\widetilde{M}:=W(\C^\flat)_A\otimes_{\AAinf{A}}\gMt$. By Proposition~\ref{prop: equivalences
  of categories to Ainf}, the~$G_{K_{\cyc,s}}$-action on $\widetilde{M}$
endows~$M$ with the structure of a 
projective $(\varphi,\Gamma_{K_{\cyc,s}})$-module with $A$-coefficients,
and $\widetilde{M}$ with its~$G_{K_{\cyc,s}}$-action is recovered as
$$\widetilde{M}=W(\C^\flat)_A\otimes_{\A_{K,A}}M.$$ 


Since~$\Gamma_{K_{\cyc,s}}$ is topologically generated by~$\gamma^{p^s}$, and  since
$(g-1)(\gM)\subset
  T^N\gMt$ for $g \in G_{K_{\cyc,s}}$, we have $(\gamma^{p^s}-1)(\gM)\subset
  T^N\gMt$. We also have~$(\gamma^{p^s}-1)(\gM)\subseteq M$, and
  since~$M\cap\gMt=\gM$ (as is easily checked, by reducing to the case
  that~$\gM$ is free, and then to the case that it is free of rank
  one), 
  we have $(\gamma^{p^s}-1)(\gM)\subset
  T^N\gM$, as required. For the uniqueness, note that if we
  have two such actions, then taking their difference gives a nonzero
  $\varphi$-linear morphism $(\gamma^{p^s})^*\gM\to T^N\gM$, which is
  impossible for~$N$ sufficiently large (see e.g.\ the first part of
  the proof of Lemma~\ref{lem: Frobenius amplification lemma over
    Ainf} below).

Finally, the claim that this gives $\gM[1/T]$ the structure of a
$(\varphi,\Gamma_{K_{\cyc,s}})$-module is immediate from Remark~\ref{rem:testing continuity on M mod u}.
\end{proof}

We continue to assume that $K$ is basic.
For any~$s~\geq~1$,
we write $\Gamma_{s,\disc}$ to denote the subgroup
of $\Gamma_{\disc}$ generated by $\gamma^{p^s}$.
We write $\cR_d^{\Gamma_{s,\disc}}$ in obvious analogy to the notation
$\cR_d^{\Gamma_{\disc}}$ introduced above.
The result of 
Corollary~\ref{cor: Caruso Liu Galois action phi Gamma discrete version}
may be interpreted as constructing a morphism
\numequation
\label{eqn:CL morphism cyclo case}
\cC^a_{d,h} \to \cR_d^{\Gamma_{s,\disc}}
\end{equation}
for sufficiently large values of $s$ (depending on $a$ and $h$),
lying over the morphism $\cC^a_{d,h} \to \cR_d.$ Of course there is also a morphism
$\cR_d^{\Gamma_{\disc}} \to \cR_d^{\Gamma_{s,\disc}}$
given by restricting the action of $\Gamma_{\disc}$ to $\Gamma_{s,\disc}$.

Since the morphism $\cR_d^{\Gamma_{s,\disc}} \to \cR_d$ given by
forgetting the $\Gamma_{s,\disc}$-action is representable by algebraic 
spaces and separated (indeed, even affine --- it is a base-change
of the morphism $\cR_d \to \cR_d \times_{\cO} \cR_d$ giving
the graph of the action of~$\gamma^{p^s}$, and this latter
morphism is representable by algebraic spaces and affine, since
the diagonal morphism of $\cR_d$ is so),
the diagonal morphism 
$$\cR_d^{\Gamma_{s,\disc}} \to
\cR_d^{\Gamma_{s,\disc}} \times_{\cR_d} \cR_d^{\Gamma_{s,\disc}}
$$
is a closed immersion.
We may thus define a closed substack $\cZ_{d,h,s}^a$ of 
$\cR_d^{\Gamma_\disc}\times_{\cR_d} \cC^a_{d,h}$
via the following $2$-Cartesian diagram:
$$\xymatrix{
\cZ^a_{d,h,s} \ar[r] \ar[d] &  \cR_d^{\Gamma_\disc}\times_{\cR_d} \cC^a_{d,h} \ar[d] \\
\cR_d^{\Gamma_{s,\disc}} \ar[r] & \cR_d^{\Gamma_{s,\disc}} \times_{\cR_d}
\cR_d^{\Gamma_{s,\disc}} }
$$
in which the lower horizontal arrow is the diagonal, and the right-hand vertical
arrow is the product of the restriction morphism and the morphism~\eqref{eqn:CL 
morphism cyclo case}.
In less formal language, the stack $\cZ^a_{d,h,s}$ parameterizes
projective rank~$d$ $\varphi$-modules $\gM$  over $\A_{K,A}^+$ of $T$-height at most~$h$,
for which $\gM[1/T]$ is endowed with a $\Gamma_{\disc}$-action extending
the canonical action of $\Gamma_{s,\disc}$ given by 
Corollary~\ref{cor: Caruso Liu Galois action phi Gamma discrete version}.

\begin{prop}\label{prop:ind structure on X via canonical
    actions}Assume that~$K$ is basic.
Then each $\cZ^a_{d,h,s}$ is contained {\em (}as a closed algebraic substack{\em )} in
$\cX_d\times_{\cR_d} \cC_{d,h},$
and the natural morphism
$\varinjlim_{a,s} \cZ_{d,h,s}^a \to \cX_d\times_{\cR_d} \cC_{d,h}$
is an isomorphism.
\end{prop}
\begin{proof}That $\cZ^a_{d,h,s}$ is a substack of
$\cX_d\times_{\cR_d} \cC_{d,h}$ is immediate from
Remark~\ref{rem:testing continuity on M mod u}; indeed, if we are given $\gM$ over $\A_{K,A}^+$ of $T$-height at most~$h$,
for which $\gM[1/T]$ is endowed with a $\Gamma_{\disc}$-action extending
the canonical action of $\Gamma_{s,\disc}$ given by 
Corollary~\ref{cor: Caruso Liu Galois action phi Gamma discrete
  version}, then we have~$(\gamma^{p^s}-1)(\gM)\subseteq
T^N\gM\subseteq T\gM$. Since $\cZ_{d,h,s}^a$ is a closed substack of 
$\cR_d^{\Gamma_\disc}\times_{\cR_d} \cC^a_{d,h}$, it is a closed
substack of~$\cX_d\times_{\cR_d} \cC_{d,h}$.

To see that the natural morphism
$\varinjlim_{a,s} \cZ_{d,h,s}^a \to \cX_d\times_{\cR_d} \cC_{d,h}$ is an isomorphism,
we need to show that it is surjective, and so we need to show that for any
finite type $\cO/\varpi^a$-algebra~$A$, and any projective $\varphi$-module $\gM$ over $\A_{K,A}^+$ of $T$-height at most~$h$,
for which $\gM[1/T]$ is endowed with a continuous
$\Gamma_{\disc}$-action, there is some~$s\ge 1$ such that 
the restriction of this action to~$\gamma_{s,\disc}$ agrees with the canonical action of $\Gamma_{s,\disc}$ given by 
Corollary~\ref{cor: Caruso Liu Galois action phi Gamma discrete
  version}.

By Remark~\ref{rem:testing continuity on M mod u}, there is
some~$s'\ge 1$ such that $(\gamma^{p^{s'}}-1)(\gM)\subseteq
T\gM$. Let~$N=N(a,h)$ be as in the statement of Lemma~\ref{lem: Caruso
  Liu Galois action on weak Wach}. By the equivalence of conditions
(3) and~(4) of Lemma~\ref{lem:testing continuity on M mod T}, there is
some~$s\ge s'$ such that $(\gamma^{p^s}-1)(\gM)\subseteq
T^N\gM$, as required.
\end{proof}

For arbitrary~$K$ (not necessarily basic) we can deduce the following description
of~$\cX_d^a$ as an Ind-algebraic stack, in the style of Lemma~\ref{lem:
  explicit Ind description for R pulled back from K basic}.

\begin{lem}
  \label{lem: explicit Ind description for X pulled back from K basic}
Fix $a \geq 1$.
For each $h$, and each sufficiently large {\em (}depending on $a$ and $h${\em )}
value of~$s$,
let $\cX_{K,d,\Kbasic,h,s}^a$ denote the scheme-theoretic image of the
  base-changed
  morphism
  \[\cZ_{d[K:\Kbasic],h,s}^a\times_{\cX_{\Kbasic,d[K:\Kbasic]}}\cX^a_{K,d}\to\cX^a_{K,d},\]
so that
$\cX_{K,d,\Kbasic,h,s}^a$ is a closed algebraic substack of $\cX_{K,d}$.
Then the natural morphism
  $\varinjlim_{h,s}\cX_{K,d,\Kbasic,h,s}^a \to \cX_{K,d}^a$
is an isomorphism.
\end{lem}
\begin{proof}
This is proved in the
  same way as Lemma~\ref{lem: explicit Ind description for R pulled
    back from K basic}, bearing in mind Proposition~\ref{prop:ind structure on X
    via canonical actions}. 
\end{proof}

\section{The connection with Galois representations}
\label{subsec:Galois reps}
In Section~\ref{subsubsec:Galois reps}  we  recalled
the  relationship between  \'etale $(\varphi,\Gamma)$-modules
without coefficients (which is to say, with coefficients in~$\Z_p$)
and $p$-adic representations of~$G_K$,
and the similar relationship between \'etale $\varphi$-modules
over $\cO_{\mathcal E}$ and $p$-adic representations of~$G_{K_{\infty}}$.
In this section, we revisit those topics in the context of
\'etale $(\varphi,\Gamma)$-modules (resp.\ \'etale $\varphi$-modules)
and $G_K$-representations (resp.\ $G_{K_{\infty}}$-modules) with coefficients.

\subsection{Galois representations with coefficients}
\label{subsec:Galois reps with coeffs}
 For a general $p$-adically complete
$\Zp$-algebra~$A$, a projective \'etale $(\varphi,\Gamma)$-module with
$A$-coefficients need not correspond to a family
of~$G_K$-representations. It will be useful, though, to have a version of
such a correspondence in the case that~$A$ is complete local Noetherian
with finite residue field. In this case, following ~\cite{MR1805474}
we let~$\widehat{\A}_{K,A}$ denote the $\m_A$-adic completion
of~${\A}_{K,A}$, and we define a \emph{formal \'etale
  $(\varphi,\Gamma)$-module} with $A$-coefficients \index{formal \'etale
  $(\varphi,\Gamma)$-module}
to be an \'etale $(\varphi,\Gamma)$-module over
$\widehat{\A}_{K,A}$ in the obvious sense.

\begin{remark}
	If $A$ is a complete local Noetherian $\cO$-algebra,
	with finite residue field,
	then the groupoid of formal \'etale $(\varphi,\Gamma)$-modules
	is equivalent to the groupoid $\cX_d(A)$ (which we remind the reader
	refers to the groupoid of morphisms $\Spf A \to \cX_d$; here
	$\Spf A$ is taken with respect to the $\mathfrak m_A$-adic
	topology on $A$).  Indeed, by definition, the 
	latter groupoid may be identified with the $2$-limit
	$\varprojlim_i \cX_d(A/\mathfrak m_A^i)$, which is easily
	seen to be equivalent 
	to the groupoid of formal \'etale $(\varphi,\Gamma)$-modules,
	via an application of~\cite[Prop.\
        0.7.2.10(ii)]{MR3075000}. 
\end{remark}

We let $\AKnrhatA$ denote 
  ~$\AKnrhat\cotimes_{\Zp} A$, where the completed
  tensor product is with respect to the usual topology
  on~$\AKnrhat$, and the  $\m_A$-adic topology
  on~$A$.
%
%
The functors~$T_A$, $\mathbf{D}_A$ defined by  \[\mathbb{D}_A(T)=(\AKnrhatA 
  \otimes_{A}T)^{G_{\Kcyc}},\]  \[T_A(M)=(  \AKnrhatA \otimes_{\widehat{\A}_{K,A}}M)^{\varphi=1}\]
then give equivalences of categories between the category of finite projective formal \'etale
$(\varphi,\Gamma)$-modules with $A$-coefficients, and the category of
finite free $A$-modules with a continuous action of~$G_K$.
(In fact the results of~\cite{MR1805474} do
not require the $(\varphi,\Gamma)$-modules to be projective, but we
will for the most part only consider projective modules in this book.)

Continuing to assume that~$A$ is complete local Noetherian with finite
residue field, we let~$\hOEA$ denote the $\m_A$-adic completion
of~$\OEA$,
  and we define
  $\widehat{\cO}_{\widehat{\cE^{\nr}},A}:=\cO_{\widehat{\cE^{\nr}}}\cotimes_{\Zp}A$,
  where the completed
  tensor product is with respect to the usual topology
  on~$\cO_{\widehat{\cE^{\nr}}}$, and the  $\m_A$-adic topology
  on~$A$.
Then the analogous statements to those of the previous
paragraph, relating representations
of~$G_{K_\infty}$ on finite free $A$-modules and \'etale
$\varphi$-modules over~$\hOEA$\footnote{One might refer to
an \'etale $\varphi$-module over~$\hOEA$ as a
``formal \'etale $\varphi$-module over $A$'', but since there
are several different species of \'etale $\varphi$-modules
under consideration throughout the book, we avoid using
this potentially ambiguous terminology.} via the functors
\[\mathbb{D}_{\infty,A}(T)=(\widehat{\cO}_{\widehat{\cE^{\nr}},A} 
  \otimes_{A}T)^{G_{K_\infty}},\]  \[T_{\infty,A}(M)=(  \widehat{\cO}_{\widehat{\cE^{\nr}},A} \otimes_{\hOEA}M)^{\varphi=1},\] can be proved in exactly the same way as in~\cite{MR1805474}
(by passage to the limit over~$A/\m_A^n$).

We will also occasionally
apply these statements in the case~$A=\Fpbar$.
Their validity in this case 
follows from their validity in the case when~$A$ is
a finite extension of~$\Fp$, 
and the fact that both Galois representations
and  projective \'etale $(\varphi,\Gamma)$-modules over~$\Fpbar$ 
arise as base changes from such finite contexts; for
Galois representations, this  
follows from the compactness of~$G_K$ and~$G_{K_{\infty}}$,
and for
\'etale $(\varphi,\Gamma)$-modules, it is an immediate consequence of
Lemma~\ref{lem:X is limit preserving}.

We can use the equivalence between Galois representations
and $(\varphi,\Gamma)$-modules as a tool to deduce facts about
the finite type points of $\cX_{d}$.  
      More precisely,
if~$\F'/\F$ is a finite extension, then the groupoid of points 
$x\in \cX_{d}(\F')$ is canonically equivalent to the groupoid of
 Galois representations $\rhobar: G_K \to \GL_d(\F')$. 
For 
this reason, we will often denote such a point $x$ simply
by the corresponding Galois representation~$\rhobar$.

Suppose now that~$\F'/\F$ is a finite extension, and
that~$x:\Spec\F'\to\cX_d$ is a finite type point, 
with  corresponding Galois representation 
$\rhobar:G_K\to\GL_d(\F')$. Let $\widehat{(\cX_d)}_x$ be the category
of Definition~\ref{adefn: deformation category} (with~$\cF$ there
being $\cX_d$). Let ~$\cO'=\cO\otimes_{W(\F)}W(\F')$ be the ring of integers in
the finite extension~$E'=W(\F')E$. 
\begin{prop}
	\label{prop:versal rings}
There is a morphism
$\Spf R_{\rhobar}^{\square,\cO'}
\to\cX_d$
which is versal at the point $x$ corresponding to~$\rhobar$,
an isomorphism 
$\Spf R_{\rhobar}^{\square,\cO'}
\times_{\cX_d}
\Spf R_{\rhobar}^{\square,\cO'}
\iso
\widehat{\GL}_{d,R_{\rhobar}^{\square,\cO'},Z(\rhobar)},$
where
$\widehat{\GL}_{d,R_{\rhobar}^{\square,\cO'},Z(\rhobar)}$
denotes the completion of $(\GL_{d})_{R_{\rhobar}^{\square,\cO'}}$
along the closed subgroup of $(\GL_{d})_{\F'}$
given by the centraliser of~$\rhobar$,
and an isomorphism 
$\Spf R_{\rhobar}^{\square,\cO'}
\times_{\widehat{(\cX_d)}_x}
\Spf R_{\rhobar}^{\square,\cO'}
\iso
\widehat{\GL}_{d,R_{\rhobar}^{\square,\cO',1}},$
where
$\widehat{\GL}_{d,R_{\rhobar}^{\square,\cO'},1}$
denotes the completion of $(\GL_{d})_{R_{\rhobar}^{\square,\cO'}}$
along the identity of $(\GL_{d})_{\F'}$.
\end{prop}
\begin{proof}
The existence of the
morphism follows from the theorem of Dee
recalled above, and the descriptions of the two fibre products are
clear from its very definition.
To see that this morphism is versal, it suffices to show that 
if $\rho: G_K \to \GL_d(A)$ is a representation with $A$ a
finite Artinian $\cO$-algebra, and if $\rho_B: G_K \to \GL_d(B)$
is a second representation, with $B$ a finite Artinian $\cO$-algebra
admitting a surjection onto $A$, such that the base change
$\rho_A$ of $\rho_B$ to $A$ is isomorphic to $\rho$
(more concretely, so that there exists $M \in \GL_d(A)$ with
$\rho = M \rho_A M^{-1}$),
then we may find $\rho': G_K \to \GL_d(B)$ which lifts $\rho$,
and is isomorphic to $\rho_B$.  But this is clear: the natural
morphism $\GL_d(B) \to \GL_d(A)$ is surjective, 
and so if $M'$ is any lift of $M$ to an element of $\GL_d(B)$,
then we may set $\rho' = M' \rho_B (M')^{-1}$.
\end{proof}

Because of the equivalence between $(\varphi,\Gamma)$-modules and
Galois representations with~$\Zp$-coefficients, 
and because of the traditional notation $\rho$ for Galois
representations, we will often denote a family of rank $d$
projective \'etale
$(\varphi,\Gamma)$-modules over~$T$, i.e.\ a morphism
$T \to \cX_d$, by $\rho_T$, or some similar notation.    We caution the
reader that this notation is chosen purely for its psychological 
suggestiveness; a family of $(\varphi,\Gamma)$-modules over a general
base $T$ does not admit a literal interpretation in terms of Galois
representations.
\subsection{Galois representations associated to
  \texorpdfstring{$(\varphi,G_K)$}{(phi,Galois)}-modules}\label{subsubsec:
Galois reps for GK phi}
At times
it will be convenient to consider the Galois representations
associated to \'etale $(\varphi,G_K)$-modules. Let~$A$ be a finite
$\cO/\varpi^a$-algebra for some~$a\ge 1$, and let ~$M$ be a finite
projective \'etale $(\varphi,G_K)$-module (resp.\
\'etale $(\varphi,G_{K_\infty})$-module) with $A$-coefficients. Then
we may apply the equivalence of 
Proposition~\ref{prop: equivalences of categories to Ainf}
to obtain an \'etale $(\varphi,\Gamma)$-module (resp.\
an \'etale $\varphi$-module) from $M$;
applying
the functor~$T_A$ (resp.\
$T_{\infty,A}$) of Section~\ref{subsec:Galois reps with coeffs} to this 
latter object yields a 
$G_K$-representation (resp.\
$G_{K_\infty}$-representation)
on a finite free~$A$-module,
which we denote simply by~$T_A(M)$
(resp.\ ~$T_{\infty,A}(M)$).

By passage to the limit over~$a$, we can extend these
functors to the case where~$A$ is a finite $\cO$-algebra, and in
particular to the case that~$A=\cO'$ is the ring of integers in a
finite extension~$E'/E$. As in  Section~\ref{subsec:Galois reps with coeffs},
we can and do
then further
extend these functors to the case~$A=\Fpbar$.

While we have defined the functors~$T_A(M)$ (resp.\ $T_{\infty,A}(M)$) 
via our equivalences of categories,
we note that both also admit a more direct description: 
namely $T_A(M)=M^{\varphi=1}$
(resp.\ $T_{\infty,A}(M)=M^{\varphi=1}$). Similarly, the composites of
the functors $\mathbb{D}_A$ and $\mathbb{D}_{\infty,A}$ with the equivalences of 
Proposition~\ref{prop: equivalences of categories to Ainf} admit the
following simple description: if $T_A$ is a finite free $A$-module
with a continuous action of $G_K$ (resp.\ $G_{K_\infty}$), then the
corresponding \'etale $(\varphi,G_K)$-module (resp.\
$(\varphi,G_{K_\infty})$-module) is given by
$W(\C^\flat)_A\otimes_AT_A$ (with~$\varphi$ acting on the first
factor in the tensor product, and $G_K$ (resp.\ $G_{K_\infty}$) acting diagonally).


\subsection{Galois representations over certain
  fields}\label{subsubsec: Galois repns into certain fields}
There is one more context in which we will need to consider the correspondence 
between Galois representations and \'etale $(\varphi,\Gamma)$-modules.
Namely, suppose that $A$ is a complete Noetherian local $\Z_p$-algebra
which is a domain,
and which has finite residue field.
Let $\cK$ denote the fraction field of~$A$, and as usual, let $\mathfrak m$
denote the maximal ideal of~$A$.
We say that a representation $\rho: G_K \to \GL_d(\cK)$ is {\em continuous}
if there exists a finitely generated $A$-submodule $L$ of $\cK^d$ which spans
$\cK^d$ over $\cK$ (a ``lattice''), such that $G_K$ preserves $L$ and
acts continuously (when $L$ is given its $\mathfrak m$-adic topology).

We write $\widehat{\A}_{K,\cK}:= \widehat{\A}_{K,A}\otimes_A \cK$,
and also write
$\AKnrhatcK := \AKnrhatA\otimes_{A}\cK.$
We have the obvious notion of a projective $(\varphi,\Gamma)$-module with
coefficients in $\widehat{\A}_{K,\cK}$.  We say that such a $(\varphi,\Gamma)$-module
$D$ is {\em \'etale} if there exists a (not necessarily projective) \'etale
$(\varphi,\Gamma)$-module $D_A$ over $\widehat{\A}_{K,A}$ contained in $D$ such
that the evident morphism
$\widehat{\A}_{K,\cK} \otimes_{\widehat{\A}_{K,A}} D_A \to D$ is an isomorphism.

Then, analogously to the case of $A$ itself, we have
functors~$T_{\cK}$, $\mathbf{D}_{\cK}$ defined by  \[\mathbb{D}_{\cK}(T)=(\AKnrhatcK
  \otimes_{\cK}T)^{G_{\Kcyc}},\]  \[T_{\cK}(M)=(  \AKnrhatcK \otimes_{\widehat{\A}_{K,\cK}}M)^{\varphi=1}\]
which give equivalences of categories between the category of finite projective formal \'etale
$(\varphi,\Gamma)$-modules with $\cK$-coefficients, 
and the category of finite dimensional $\cK$-vector spaces with a continuous action of~$G_K$.
(This is proved by passing to lattices on each side, and using 
the results of~\cite{MR1805474}; note that the lattices 
involved need not be projective in general, and so we apply those results 
in their full generality.)

This formalism is most often applied in the literature in the case when $A = \Z_p$,
so that $\cK = \Q_p$.  However, we will not consider 
\'etale $(\varphi,\Gamma)$-modules with $\Q_p$-coefficients in this book. Rather
we will apply the preceding formalism only once, namely in our analysis
of the closed $\Fbar_p$-points of $\cX_d$,
which we make in Section~\ref{subsec:closed points};
and in this application,
we will take $A$ to be a complete local domain of characteristic~$p$.

\section{\texorpdfstring{$(\varphi,G_K)$}{(phi,Galois}-modules
and restriction}
\label{subsec:restriction}
Our goal in this section is to define certain morphisms of stacks
which are the analogues,
for families of $(\varphi,\Gamma)$-modules,
of the restriction functors on Galois representations
$\rho \mapsto \rho_{| G_{K_{\infty}}}$
and $\rho \mapsto \rho_{|G_L}$ (for any finite extension $L$ of~$K$).


\begin{defn}\index{$\cR_{\BK,d}$}
We let~$\cR_{\BK,d}$ denote the moduli stack
of rank~$d$ projective \'etale $\varphi$-modules over~$\OEA$; that is,
the stack~$\cR_d$ constructed in Section~\ref{subsec: moduli stacks of phi
  modules} in the case when~$\A_A=\OEA$.
\end{defn}

Recall that we also have the stack~$\cX_d$ of \'etale $(\varphi,\Gamma)$-modules
of Definition~\ref{defn: Xd}. We now construct a 
morphism~$\cX_d\to\cR_{\BK,d}$ 
which corresponds to the restriction of Galois representations from~$G_K$
to~$G_{K_\infty}$. 
\begin{prop}
  \label{prop: natural morphism X to R}There is a canonical
  morphism~$\cX_{d}\to\cR_{\BK,d}$. If~$A$ is a complete local Noetherian
  $\cO$-algebra with finite residue field, or equals  $\Fpbar$, then the
  morphism~$\cX_d(A)\to\cR_{\BK,d}(A)$ is given by
  restriction of the corresponding representation of~$G_K$
  to~$G_{K_\infty}$.
\end{prop}
\begin{proof}
  If~$A$ is a finite type $\cO/\varpi^a$-algebra, for some~$a\ge
  1$, then we define a morphism of groupoids $\cX_{d}(A) \to\cR_{\BK,d}(A)$
  via the equivalences of categories of Proposition~\ref{prop:
    equivalences of categories to Ainf}; that is, given a finite
  projective \'etale $(\varphi,\Gamma)$-module with $A$-coefficients,
  we form the corresponding finite projective \'etale $(\varphi,G_K)$-module,
  which yields a finite projective \'etale $(\varphi,G_{K_\infty})$-module by
  restricting the~$G_K$-action to~$G_{K_\infty}$, and thus gives an \'etale
  $\varphi$-module over~$\OEA$. Since both~$\cX_d$ and~$\cR_{\BK,d}$
  are limit preserving (by  Corollary~\ref{cor: basic properties of C and R p adic stacks} and Lemma~\ref{lem:X is limit preserving}), it follows from~\cite[Lem.\ 2.5.4, Lem.\
  2.5.5~(1)]{EGstacktheoreticimages} that this construction determines
  a morphism ~$\cX_{d}\to\cR_{\BK,d}$.

  Suppose now that ~$A$ is a complete local Noetherian $\cO$-algebra
  with finite residue field. 
  Let $M$ be a formal
  \'etale $(\varphi,\Gamma)$-module corresponding to a morphism
  $\Spf A\to\cX_{d}$, and~$M_\infty$ be the \'etale $\varphi$-module
  over~$\hOEA$ corresponding to the composite
  $\Spf A\to \cX_d\to\cR_{\BK,d}$.
  The relationship between $M$ and $M_{\infty}$ is
  expressed as an isomorphism
  \numequation
  \label{eqn:module isomorphism}
  \widehat{W(\C^\flat)}_A \otimes_{\widehat{\A}_{K,A}} M
  \iso 
    \widehat{W(\C^\flat)}_A\otimes_{\hOEA}M_\infty.
  \end{equation}
  
  Recall that the Galois representation associated to $M$ is 
  defined via $T_A(M) := (\AKnrhatA\otimes_{\widehat{\A}_{K,A}}M)^{\varphi = 1},$
  and that the evident
  $(\varphi,G_{K})$-equivariant $\AKnrhatA$-linear morphism
  \[\AKnrhatA\otimes_{A}T_A(M)\to
    \AKnrhatA\otimes_{\widehat{\A}_{K,A}}M \]
  is then an isomorphism (see for
  example~\cite[Prop.\ 2.1.26]{MR1805474}).   
  Thus there is an induced natural
  $(\varphi,G_{K})$-equivariant
  isomorphism
  \[\widehat{W(\C^\flat)}_A\otimes_{A}T_A(M)\isoto
    \widehat{W(\C^\flat)}_A\otimes_{\widehat{\A}_{K,A}}M, \] where we recall
  that we write $\widehat{W(\C^\flat)}_A$ for
the $\m_A$-adic completion of~$W(\C^\flat)_A$.  Similarly, we obtain a natural
  $(\varphi,G_{K_\infty})$-equivariant isomorphism
  \[ \widehat{W(\C^\flat)}_A\otimes_{A}T_{\infty,A}(M_\infty)\isoto
    \widehat{W(\C^\flat)}_A\otimes_{\hOEA}M_\infty.\]
Combining these two isomorphisms with~\eqref{eqn:module isomorphism},
we obtain a natural $(\varphi,G_{K_\infty})$-equivariant isomorphism
\[\widehat{W(\C^\flat)}_A\otimes_{A}T_A(M)\isoto
  \widehat{W(\C^\flat)}_A\otimes_{A}T_{\infty,A}(M_\infty). \]
Each of $T_A(M)$ and $T_{\infty,A}(M_{\infty})$ 
is a free $A$-module with trivial~$\varphi$-action,
and so if we pass to $\varphi$-invariants in this isomorphism,
and take into account
Lemma~\ref{lem: phi invariants in W C flat with coefficients},
we obtain an isomorphism
$T_A(M)|_{G_{K_\infty}}\isoto T_{\infty,A}(M_\infty)$, 
which is what we wanted to show.
Finally, if $A=\Fpbar$, 
the result follows by taking direct limits.
\end{proof}
\begin{prop}
  \label{prop: diagonal of X to R}The diagonal morphism~$\Delta:
  \cX_d\to\cX_d\times_{\cR_{\BK,d}}\cX_d$
  induced by the morphism of
  Proposition~{\em \ref{prop: natural morphism X to R}}
  is a closed immersion.
\end{prop}
\begin{proof}	
	The product $\cX_d\times_{\cR_{\BK}} \cX_d$ can be described explicitly
	as follows: its $A$-valued points are pairs $(M_1,M_2)$ of \'etale
	$(\varphi,G_K)$-modules with a $G_{K_{\infty}}$-equivariant isomorphism
	between them.  The diagonal is defined by $M \mapsto (M,M),$ with
	the isomorphism being the identity.
%
To see that this defines a closed
	immersion,  we have to check that if we 
	are given a pair~$(M_1,M_2)$ of $(\varphi,G_K)$-modules with $A$-coefficients,
	together with a $G_{K_{\infty}}$-equivariant isomorphism $f:M_1\isoto M_2$,
	then the locus where~$f$ becomes $G_K$-equivariant
	is closed.   That this is so follows from Corollary~\ref{cor:equivariance locus}.  
	(That corollary shows that, for each $g \in G_K$,
       	there is an ideal $J_g \subseteq A$
	which cuts out the condition for the given
       isomorphism to commute with~$g$.   The ideal $J := \sum_g J_g$
       then cuts out the condition for the given isomorphism to be
       $G_K$-equivariant.)
%
\end{proof}

We now construct the restriction maps corresponding to finite extensions
of~$K$.

\begin{lem}
  \label{lem: restriction from X to Xs}Let~$L/K$ be a finite
  extension. There is a canonical
  morphism~$\cX_{K,d}\to\cX_{L,d}$, which is representable by
  algebraic spaces, affine, 
  and of finite
  presentation. If~$A$ is a complete local Noetherian
  $\cO$-algebra with finite residue field,
  or equals~$\Fpbar$, then the
  morphism~$\cX_{K,d}(A)\to\cX_{L,d}(A)$ is given by
  restriction of the corresponding representation of~$G_K$
  to~$G_{L}$.
\end{lem}
\begin{proof}
As in the case of Proposition~\ref{prop: natural morphism X to R},
the morphism is defined first in the case of $\cO$-algebras that
are of finite type over $\cO/\varpi^a$ for some $a$ by
  applying the equivalences of categories of Proposition~\ref{prop:
    equivalences of categories to Ainf}:
we send an \'etale $(\varphi,G_K)$-module~$M$ 
to an \'etale $(\varphi,G_L)$-module~$M_L$ via restricting the continuous
$G_K$-action to a continuous $G_L$-action.
Since $\cX_{K,d}$ and $\cX_{L,d}$ 
  are limit preserving (by Lemma~\ref{lem:X is limit preserving}),
  it follows from~\cite[Lem.\ 2.5.4, Lem.\ 2.5.5~(1)]{EGstacktheoreticimages}
  that this construction determines a morphism $\cX_{K,d} \to \cX_{L,d}$.

Before verifying the various claimed properties of this morphism,
we confirm that it has the claimed effect on associated Galois representations.
To this end, we note that
if~$A$ is a complete local Noetherian
  $\cO$-algebra, then, as in the proof of Proposition~\ref{prop:
    natural morphism X to R}, it follows from the definition that we
  have a 
  natural $(\varphi,G_{L})$-equivariant isomorphism
  \[\widehat{W(\C^\flat)}_A\otimes_{A}T_A(M)\isoto
    \widehat{W(\C^\flat)}_A\otimes_{A}T_{A}(M_L), \] and by Lemma~\ref{lem: phi
    invariants in W C flat with coefficients},
  taking~$\varphi$-invariants induces an isomorphism
  $T_A(M)|_{G_L}\isoto T_{A}(M_L)$. We deduce the same statement if
    ~$A=\Fpbar$ 
    by taking direct limits.  
%

We now verify the claimed properties of the restriction morphism $\cX_{K,d}
\to \cX_{L,d}$.
Suppose firstly that~$L/K$ is Galois, and let $\{g_i\}$
be a set of (finitely many) coset representatives for $G_L$ in
$G_K$. For each~$i$, we may give~$g_i^*M$ the structure of a
$(\varphi,G_L)$-module by letting each $h\in G_L$ act as~$g_i^{-1}h g_i$
acts on~$M$; then to extend the action of~$G_L$ to an action of~$G_K$
is to give an isomorphism of~$(\varphi,G_L)$-modules $g_i^*M \to M$
for each index~$i$, satisfying a slew of compatibilities. (Since $G_L$
is open in $G_K$, the continuity of the $G_K$-action is automatic,
given that the $G_L$-action that it is extending is continuous.)

We let $\cY_i$ denote the stack classifying objects $M_A$ 
of~$\cX_{L,d}(A)$,
endowed with an isomorphism 
$g_i^*M_A\isoto M_A$.  If we regard $g_i^*$  as an automorphism
of~$\cX_{L,d}(A)$,
then we  may form  its  graph~$\Gamma_i$,
and we then have an isomorphism of stacks
$$\cY_i \iso
\cX_{L,d}(A)\times_{\Delta,
\cX_{L,d}(A) \times
\cX_{L,d}(A),\Gamma_i}
\cX_{L,d}(A).$$
(Here $\Delta$ denotes the diagonal of~$\cX_{L,d}$;
\emph{cf.}\ the definition of $\cR_d^{\Gamma_{\disc}}$ in 
Section~\ref{subsec: defn of Xd}
above, together 
with Proposition~\ref{prop:Gamma-disc modules as fixed points; etale case}.)
The projection onto the second factor $\cY_i \to 
\cX_{L,d}(A)$,
which corresponds to
forgetting the isomorphism, is a base-change of the
diagonal~$\Delta,$
and so is representable by algebraic spaces, affine and of finite presentation by Proposition~\ref{prop: properties of X pulled back from R}.

We can rephrase our interpretation of objects 
of~$\cX_{K,d}(A)$ 
as objects of $\cX_{L,d}(A)$
endowed with isomorphisms
$g_i^*M_A\isoto M_A$ satisfying certain compatibilities
as the existence of  a closed immersion
$$
\cX_{K,d}\into
\cY_1 \times_{\cX_{L,d}}\times \cdots \times_{\cX_{L,d}}
\cY_n
$$  
(the point  being that the compatibilities
arise as
base-changes of the double diagonal, and so impose closed conditions).
Since the morphisms ~$\cY_i\to\cX_{L,d}$ are
representable by algebraic spaces, affine and of finite presentation,
it follows that the same is true of the morphism
$\cX_{K,d}\to\cX_{L,d}$, as claimed.


Finally, let~$L/K$ be a general finite extension. Let~$M/K$ be the
normal closure of~$L/K$, so that we have
morphisms \[\cX_{K,d}\to\cX_{L,d}\to\cX_{M,d}.\]
By what we have already proved, both the second arrow and the composite
of the two arrows are representable by algebraic
spaces and affine, and of finite presentation.
An affine morphism is separated, and thus has affine and finite type diagonal.
One sees (e.g.\ by
Lemma~\ref{lem:finiteness of diagonals})
that the diagonal of an affine morphism of finite type (and so, in particular,
the diagonal of an affine morphism of finite presentation)
is furthermore of finite presentation.
A standard graph argument,
in which we factor the first morphism as
$$\cX_{K,d} \hookrightarrow \cX_{K,d} \times_{\cX_{M,d}} \cX_{L,d} 
\to \cX_{L,d}$$
(the first arrow in this factorization
being the graph of the first morphism,
which is a base change of the diagonal
$\cX_{L,d} \to \cX_{L,d} \times_{\cX_{M,d}} \cX_{L,d} $,
and the second arrow being the projection, which is a base change
of the composite $\cX_{K,d} \to \cX_{L,d} \to \cX_{M,d}$),
then shows that the first morphism is also
representable by algebraic spaces and affine, and of finite presentation.
%
\end{proof}


\section{Tensor products and duality} \label{subsec: tensor product
  of phi gamma and duality}
If~$M$ is a projective \'etale $(\varphi,\Gamma)$-module with $A$-coefficients, we define the dual
$(\varphi,\Gamma)$-module by \[M^\vee:=\Hom_{\A_{K,A}}(M,\A_{K,A}); \] if $A$
is a finite $\cO$-module, then there is a natural isomorphism
$T(M^\vee)\cong T(M)^\vee$. Similarly, for
  each~$M$, we let~$M(1):=M\otimes_{\A_{K,A}}\A_{K,A}(1)$ denote the Tate
  twist, where $\A_{K,A}(1)$ denotes the free $(\varphi,\Gamma)$-module of
  rank~$1$ with a generator~$v$ on which $\Gamma$ acts via the
  cyclotomic character and~$\varphi(v)=v$.  
  (Note that if~$A=\cO$ then $T(\A_{K,A}(1))=\cO(1)$, the usual
  Tate twist on the Galois side.)  We also define the Cartier dual
$M^*:=M^\vee(1)$, which again in the case that $A$ is a finite
$\cO$-algebra is compatible with the usual notion for Galois
modules. If $T\to\cX_d$ corresponds to a family $\rho_T$,
then we adopt the natural convention of writing ~$\rho_T(1)$, $\rho_T^\vee$ and~$\rho_T^*$ for the
corresponding families.

Given two projective \'etale $(\varphi,\Gamma)$-modules
$M_1$, $M_2$ with $A$-coefficients of respective ranks~$d_1,d_2$, we may form the tensor
product $(\varphi,\Gamma)$-module $M_1\otimes M_2$ (given by the
tensor product on  underlying
$\A_{K,A}$-modules, and the tensor products of the actions of each of~$\varphi$
and $\Gamma$). If~$A$ is a finite $\cO$-module, then we have a natural
isomorphism \[T(M_1\otimes M_2)\cong T(M_1)\otimes T(M_2).\] The
tensor product induces a natural morphism \[\cX_{d_1}\times_{\cO}
  \cX_{d_2}\to\cX_{d_1d_2};\] in particular, twisting by (the
$(\varphi,\Gamma)$-module corresponding to) any
character~$G_K\to\cO^\times$ gives an automorphism of each~$\cX_d$. Given morphisms $S\to\cX_{d_1}$,
$T\to\cX_{d_2}$, denoted 
by $\rho_S$ and $\rho'_T$, we write $\rho_{S}\boxtimes
\rho'_T$ for the family corresponding to the composite 
\[S\times_{\cO}
  T\to\cX_{d_1}\times_{\cO}\cX_{d_2}\to\cX_{d_1d_2}.\] If $S=T$ then we
write $\rho_T\otimes\rho'_T$ for the composite  \[T\to T\times_{\cO}
  T\to\cX_{d_1}\times_{\cO}\cX_{d_2}\to\cX_{d_1d_2},\]where the first map is
the diagonal embedding.

\chapter{Crystalline and semistable moduli stacks}\label{sec:
  crystalline and semistable}
In this chapter we build on the results of Appendix~\ref{app: BKF pst} to
define the potentially semistable and potentially crystalline
substacks of~$\cX_{K,d}$. Using the results of the earlier chapters,
and Kisin's results on potentially semistable lifting
rings~\cite{MR2373358}, we show that these stacks are $p$-adic formal
algebraic stacks, and compute the dimensions of their special fibres.
\section{Notation} Recall from Section~\ref{subsubsec: Kummer}
that the embedding~$\gS_A\into\AAinf{A}$ depends on the choice of a
uniformizer~$\pi$ of~$K$, and of a compatible
family~$\pi^{1/p^\infty}$ of $p$-power roots of unity
in~$\overline{K}$. In this chapter we will want to consider all
possible choices of~$\pi$ and~$\pi^{1/p^\infty}$, and consequently we
introduce some notation to do so. We write~$\pi^\flat\in\cO_\C^\flat$ for
the element determined by~$\pi^{1/p^\infty}$. For each~$s\ge 0$ we
write~$K_{\pi^\flat,s}$ for $K(\pi^{1/p^s})$,
and~$K_{\pi^\flat,\infty}$ for~$\cup_sK_{\pi^\flat,s}$. We
write~$\gS_{\pi^\flat,A}$ for the image of~$\gS_A$ in ~$\AAinf{A}$ via
the homomorphism determined by~$u\mapsto[\pi^\flat]$,
and~$\OEApiflat$ for the image of~$\OEA$ in~$W(\C^\flat)_A$
under the corresponding homomorphism. When~$\piflat$ is fixed, we will often write~$u$
for~$[\pi^\flat]$ when discussing~$\gS_{\piflat,A}$.  We write~$E_\pi$ for the
Eisenstein polynomial corresponding to~$\pi$. We also have a natural
map~$\theta:\AAinf{A}\to\cO_{\C,A}$ (where the target is the tensor
product  $\cO_{\C}\cotimes A$, which is 
the quotient of the source by the principal ideal 
generated by the kernel of the usual map $\theta:\Ainf\to\cO_{\C}$).

\begin{rem}\label{rem: we don't twist our embeddings by phi}In
  contrast to many papers in the literature (e.g.\
  \cite{2016arXiv160203148B,liulattice2}) we do not twist the
  embedding $\gS\to\Ainf$ by~$\varphi$. While such a twist is
  important in applications involving comparisons to Fontaine's period
  rings, it is more convenient for us to use the embedding we have
  given here. 
  Since~$\varphi$ is an automorphism of~$\Ainf$, it is
  in any case straightforward to pass back and forth between these two
  conventions. 
\end{rem}


\section[Breuil--Kisin modules and Breuil--Kisin--Fargues
  modules]{Breuil--Kisin modules and Breuil--Kisin--Fargues
  modules \sectionmark{Breuil--Kisin--Fargues
  modules}}\sectionmark{Breuil--Kisin--Fargues
  modules}\label{subsec: BK and BKF modules}
\begin{defn}
Fix a choice of~$\pi^\flat$.	Let $A$ be a $p$-adically complete
$\cO$-algebra which is topologically of finite type. 
We define a {\em projective Breuil--Kisin module} \index{Breuil--Kisin module} 
(resp.\ a {\em projective Breuil--Kisin--Fargues module}) \index{Breuil--Kisin--Fargues module}  {\em of height at
  most~$h$ with $A$-coefficients} to be a finitely generated projective $\gS_{\pi^\flat,A}$-module
  (resp.\ $\AAinf{A}$-module)
  $\gM$, equipped with a~$\varphi$-semi-linear morphism
  $\varphi:\gM\to\gM$, with the property that the corresponding
  morphism 
  $\Phi_{\gM}:\varphi^*\gM\to\gM$ is
  injective, with cokernel killed by~$E_\pi^h$. 
If~$\gM$ is a Breuil--Kisin module   then
$\AAinf{A}\otimes_{\gS_{\pi^\flat,A}}\gM$ is a Breuil--Kisin--Fargues module.
\end{defn}
\begin{rem}
  \label{rem: all BK BKF are effective}Note that our definition of a
  Breuil--Kisin--Fargues module is less general than the definition
  made in~\cite{2016arXiv160203148B}, in that we require~$\varphi$ to
  take~$\gM$ to itself; this corresponds to only considering Galois
  representations with non-negative Hodge--Tate weights. This
  definition is convenient for us, as it allows us to make direct
  reference to the literature on Breuil--Kisin modules. The
  restriction to non-negative Hodge--Tate weights is harmless in our
  main results, as we can reduce to this case by twisting by a large
  enough power of the cyclotomic character (the interpretation of
  which on Breuil--Kisin--Fargues modules is explained in~\cite[Ex.\
  4.24]{2016arXiv160203148B}).
\end{rem}

\begin{defn}
  \label{defn: BKF module with GK action coefficients}Let $A$ be a $p$-adically complete
$\cO$-algebra which is topologically of finite type.
A Breuil--Kisin--Fargues $G_K$-module with $A$-coefficients is \index{Breuil--Kisin--Fargues module!$G_K$-module}
  a Breuil--Kisin--Fargues module with $A$-coefficients, equipped with a continuous
  semilinear action of~$G_K$ which commutes with~$\varphi$.
\end{defn}
Note that if~$\gMt$ is a Breuil--Kisin--Fargues $G_K$-module with
$A$-coefficients, then $W(\C^\flat)_A\otimes_{\AAinf{A}}\gMt$ is
naturally an \'etale $(\varphi,G_K)$-module in the sense of
Definition~\ref{defn: GK module over Ainf}.  The following definition
is motivated by Theorem~\ref{athm: admits all descents if and only if
  semistable} and Corollary~\ref{acor: admits all descents if and only
  if potentially semistable}.
\begin{defn}
  \label{defn: descending BKF to BK coefficients}
Let $A$ be a $p$-adically complete
$\cO$-algebra which is topologically of finite type over~$\cO$ (and recall
then that $\AAinf{A}$ is faithfully flat over $\gS_{\pi^\flat,A}$,
by Proposition~\ref{prop:top f.t. flatness}).
We say that a
  Breuil--Kisin--Fargues $G_K$-module with $A$-coefficients
$\gMt$
\emph{descends for~$\pi^\flat$} or \emph{descends
    to~$\gS_{\piflat,A}$} if there is a
  Breuil--Kisin module~$\gM_{\pi^\flat}$ with
  $\gM_{\pi^\flat}\subseteq(\gMt)^{G_{K_{\pi^\flat,\infty}}}$ and
for which the natural map
  $\AAinf{A}\otimes_{\gS_{\pi^\flat,A}}\gM_{\pi^\flat} \to \gMt$ is an isomorphism.

We say
  that~$\gMt$ \emph{admits all descents} if it descends for every
  \index{Breuil--Kisin--Fargues module!admitting all descents}
  choice of~$\pi^\flat$ (for every choice of~$\pi$), and if
  furthermore 

  \begin{enumerate}
  \item\label{item: M mod u descends} The $W(k)\otimes_{\Zp}A$-submodule $\gM_{\piflat}/[\piflat]\gM_{\piflat}$ of
  $W(\overline{k})_A\otimes_{\AAinf{A}}\gMt$ is independent of the
    choice of~$\pi$ and~$\piflat$.
  \item\label{item: M mod E descends} The $\cO_K\otimes_{\Zp}A$-submodule $\varphi^*\gM_{\piflat}/E_{\piflat}\varphi^*\gM_{\piflat}$
    of  $\cO_{\C,A}\otimes_{\AAinf{A},\theta}\varphi^*\gMt$  is independent of the
    choice of~$\pi$ and~$\piflat$.
  \end{enumerate}
  If~$\gMt$ admits all descents, then we say that it is
  \emph{crystalline} if for each~$\piflat$, and each~$g\in G_K$, we
  have \numequation\label{eqn: condition on BK GK for
    crystalline}(g-1)(\gM_{\piflat})\subseteq \varphi
  ^{-1}(\mu)[\piflat]\gMt \end{equation} (where ~$\mu=[\varepsilon]-1$
for some compatible choice of roots of unity
$\varepsilon=(1,\zeta_p,\zeta_{p^{2}},\dots)\in\cO_\C^\flat$).

 If~$L/K$ is a finite Galois extension, then we say that~$\gMt$
 \emph{admits all descents over~$L$} if the corresponding
 Breuil--Kisin--Fargues $G_L$-module (obtained by restricting
 the~$G_K$-action to~$G_L$) admits all descents.
\end{defn}

\begin{rem}
  \label{rem: descents are automatically projective}There is
  considerable redundancy and rigidity in Definition~\ref{defn:
    descending BKF to BK coefficients}. Note in particular that if 
  for some~$\piflat$ there exists a $\gS_{\piflat,A}$-module
  $\gM_{\piflat}$ with
  $\gMt=\AAinf{A}\otimes_{\gS_{\pi^\flat,A}}\gM_{\pi^\flat}$, then
  since~$\gS_{\piflat,A}\to\AAinf{A}$ is faithfully flat,
 the $\gS_{\piflat,A}$-module $\gM_{\piflat}$ is automatically finite projective
  (by~\cite[\href{https://stacks.math.columbia.edu/tag/058S}{Tag
    058S}]{stacks-project}).
\end{rem}
We will find the following basic lemmas useful.
\begin{lem}
  \label{lem: descent for GK action depends only on pi}Fix a choice
  of~$\pi$. Suppose that $\gMt$ is a
Breuil--Kisin--Fargues $G_K$-module with~$A$ coefficients, where~$A$
is a $p$-adically complete
$\cO$-algebra which is topologically of finite type over~$\cO$. Then~$\gMt$ descends for some
choice of~$\pi^\flat$ \emph{(}for our fixed~$\pi$\emph{)} if and only if it descends
for all such choices.
\end{lem}
\begin{proof}
  We can choose an element~$g\in G_K$ with
  $g(\pi^\flat)=(\pi^\flat)'$, so
  that~$\gS_{(\pi^\flat)',A}=g(\gS_{\pi^\flat,A})$.  Then if~$\gM_\pi$
  is a descent to~$\gS_{\pi,A}$, $g(\gM_\pi)$ is a descent
  to~$\gS_{(\pi^\flat)',A}$.
\end{proof}

\begin{lem}
  \label{lem: uniqueness of descent for G K infty Ainf}Let $A$ be a $p$-adically complete
$\cO$-algebra which is topologically of finite type over~$\cO$. 
  Suppose that $\gMt$ is a
  Breuil--Kisin--Fargues  $G_K$-module. Then if $\gMt$ descends for~$\pi^\flat$, the Breuil--Kisin
  module~$\gM_{\pi^\flat}$ is uniquely determined.
\end{lem}
\begin{proof}
Assume to begin with that $A$ is an $\cO/\varpi^a$-algebra of finite type,
for some $a \geq 1$.
  By Proposition~\ref{prop: equivalences of categories to Ainf},
the $\varphi$-module 
  $(W(\C^\flat)_{A}\otimes_{\AAinf{A}}\gMt)^{G_{K_{\pi^\flat,\infty}}}$ uniquely descends to an
  \'etale $\varphi$-module~$M$ over~$\cO_{\mathcal{E},\pi^\flat,A}$.
Thus
if~$\gM_1$, $\gM_2$ are two descents of~$\gMt$
  to~$\gS_{\piflat,A}$, then $\cO_{\mathcal{E},\pi^\flat,A}\otimes_{\gS_{\piflat,A}}\gM_1$
and  $\cO_{\mathcal{E},\pi^\flat,A}\otimes_{\gS_{\piflat,A}}\gM_2$ coincide; they are both
equal to~$M$.

This allows us to show that 
  $\gM_1+\gM_2$ is also a descent of~$\gMt$.  For this, we have to show that
the natural map $ \AAinf{A}\otimes_{\gS_{\piflat,A}}(\gM_1  + \gM_2) \to \gMt$
is an isomorphism. 
(Note that it will then follow automatically that $\gM_1+\gM_2$
is projective over $\gS_{\piflat,A}$,
  by Remark~\ref{rem: descents are automatically projective}.)
That it is a surjection is clear, since it is already a surjection
when restricted to each summand.  To see that it is an injection, we may check
that its composite with
the injection $\gMt \hookrightarrow  W(\C^\flat)_A\otimes_{\AAinf{A}}\gMt$ is
injective.  But this composite may be factored as the composite of the embedding
$$\AAinf{A}\otimes_{\gS_{\piflat,A}}(\gM_1 + \gM_2) \hookrightarrow
 \AAinf{A} \otimes_{\gS_{\piflat,A}}M $$
(obtained by tensoring the embedding $\gM_1 + \gM_2 \hookrightarrow M$ with
the extension $\AAinf{A}$ of $\gS_{\piflat, A}$, which is flat by
Proposition~\ref{prop:top f.t. flatness}),
and the isomorphism 
\begin{multline*}
 \AAinf{A} \otimes_{\gS_{\piflat,A}}M\cong
 \bigl( \AAinf{A}\otimes_{\gS_{\pi^{\flat},A}} \cO_{\mathcal{E},\pi^{\flat},A} 
  \bigr)
\otimes_{\cO_{\mathcal{E},\pi^\flat,A}} M
\\
\cong  W(\C^{\flat})_A  \otimes_{\cO_{\mathcal{E},\piflat,A}}M
\cong W(\C^\flat)_A \otimes_{\AAinf{A}}\gMt.
\end{multline*}
Since $\gM_1 + \gM_2$ is also a descent,
  we may replace~$\gM_2$ by~$\gM_1+\gM_2$, and 
  therefore assume that
  $\gM_1\subseteq\gM_2$.
  Since~$\AAinf{A}\otimes_{\gS_{\piflat,A}}(\gM_2/\gM_1)=0$, and
  as ~$\gS_{\piflat,A}\to\AAinf{A}$ is faithfully flat
  (again by Proposition~\ref{prop:top f.t. flatness}),
  we see that~$\gM_2/\gM_1=0$, as required.

Suppose now that $A$ is $p$-adically complete and topologically of finite type.
Let $\gM_1$ and $\gM_2$ be two descents of $\gMt$.  Since $\gS_A \to \AAinf{A}$
is faithfully flat
(again by Proposition~\ref{prop:top f.t. flatness}),
we find that $\gM_1/\varpi^a$ and $\gM_1/\varpi^a$ both embed into
$\gMt/\varpi^a$, for any~$a~\geq~1$, and both provide descents for $\piflat$ of $\gMt/\varpi^a$.
Applying the uniqueness result we've already proved (our coefficients now being
$A/\varpi^a$), we find that $\gM_1/\varpi^a$ and $\gM_2/\varpi^a$ coincide
as submodules of $\gMt/\varpi^a$.  Passing to the inverse limit over~$a$,
we find that $\gM_1$ and $\gM_2$ coincide as submodules of~$\gMt$,
proving the desired uniqueness in general.
\end{proof}

\section{Canonical extensions of $G_{K_{\infty}}$-actions}
\label{subsec:canonical actions}

In this section we will consider a fixed choice of~$\pi^\flat$, which
we will accordingly drop from the notation. For each~$s\ge 1$,
write~$K_s=K(\pi^{1/p^s})$. We now use a variant of the arguments of
~\cite{MR2745530}, which show that the trivial action
of~$G_{K_\infty}$ on a Breuil--Kisin module with
$\cO/\varpi^a$-coefficients can be extended to some~$G_{K_s}$. We
begin with some preliminary lemmas, the first of which is an analogue
for Breuil--Kisin modules of~\cite[Lem.\
5.2.14]{EGstacktheoreticimages}, and is proved in a similar way.

\begin{lem}
  \label{lem: adding projective Kisin to get free}Let $A$ be a $p$-adically complete
$\cO$-algebra which is topologically of finite type, and let~$\gM$ be a finite projective Breuil--Kisin module with $A$-coefficients of height at
  most~$h$. Then~$\gM$ is a direct summand of a finite free Breuil--Kisin
  module with $A$-coefficients of height at
  most~$h$.
\end{lem}
\begin{proof}
  Since the map $\Phi_{\gM}:\varphi^*\gM\to\gM$ is an injection with
  cokernel killed by~$E^h$, there is a
  map~$\Psi_{\gM}:\gM\to\varphi^*\gM$ with $\Psi_{\gM}\circ\Phi_{\gM}$
  and $\Phi_{\gM}\circ\Psi_{\gM}$ both being given by multiplication
  by~$E^h$. Let~$\gP$ be a finite projective $\gS_A$-module
  such that~$\gF:=\gM\oplus\gP$ is a finite free~$\gS_A$
  module. Set~$\gQ:=\gP\oplus\gM\oplus\gP$, a finite
  projective~$\gS_A$-module.

  Note that~$\gM\oplus\gQ\cong\gF\oplus\gF$ is a finite
  free~$\gS_A$-module, so it suffices to show that~$\gQ$ can
  be endowed with the structure of a Breuil--Kisin module of height at
  most~$h$. Since~$\gF$ is free, we can choose an
  isomorphism~$\Phi_{\gF}:\varphi^*\gF\isoto\gF$, and we then
  define~$\Phi_{\gQ}:\varphi^*\gQ\to\gQ$ as the composite
  \begin{align*}
    \varphi^*\gQ
    &=\varphi^*\gP\oplus\varphi^*\gM\oplus\varphi^*\gP\\&\isoto
    \varphi^*\gM\oplus\varphi^*\gP\oplus\varphi^*\gP\\&=\varphi^*\gF\oplus\varphi^*\gP\\&\stackrel{\Phi_{\gF}}{\to}\gF\oplus\varphi^*\gP\\&\isoto
    \gP\oplus\gM\oplus\varphi^*\gP\\&\stackrel{\Psi_{\gM}}{\to} \gP\oplus\varphi^*\gM\oplus\varphi^*\gP\\&= \gP\oplus\varphi^*\gF\\&\stackrel{\Phi_{\gF}}{\to}\gP\oplus\gF=\gQ.
  \end{align*}Every morphism in this composite other than~$\Psi_{\gM}$
  is an isomorphism. Since
  $\Psi_{\gM}\circ\Phi_{\gM}=\Phi_{\gM}\circ\Psi_{\gM}=E^h$, we see
  that~$\Psi_{\gM}$ is an injection with cokernel killed
  by~$E^h$, so the same is true of~$\Phi_{\gQ}$, as required.
\end{proof}

The following lemma is proved by a standard Frobenius amplification
argument, which is in particular almost identical to the proof
of~\cite[Lem.\ 4.1.9]{CEGSKisinwithdd}.
\begin{lem}
  \label{lem: Frobenius amplification lemma over Ainf}Let~$\gMt$, $\gNt$
  be projective Breuil--Kisin--Fargues modules  of height at
  most~$h$, where $A$ is a finite type
  $\cO/\varpi^a$-algebra. Suppose that~$N\ge e(a+h)/(p-1)$. Let \[f:\gMt\to
    \gNt\] be an $\AAinf{A}$-linear map,
  and suppose that for all~$m\in\gMt$,
  $(f\circ\Phi_{\gMt}-\Phi_{\gNt}\circ\varphi^*f)(m)\in
  u^{pN}\gNt$.

  Then there is a unique $\AAinf{A}$-linear, $\varphi$-linear morphism
  $f':\gMt\to
    \gNt$ with the property that for all~$m\in\gMt$,
  $(f'-f)(m)\in u^N\gNt$; in fact, we have $(f'-f)(m)\in u^{N+1}\gNt$.
\end{lem}
\begin{proof}We firstly prove uniqueness. Suppose that~$f''$ also
  satisfies the properties that~$f'$ does, and write~$g=f'-f''$, so that by
  assumption we have
  that~$\im g\subseteq u^N \gNt$ and
  $g\circ\Phi_{\gMt}=\Phi_{\gNt}\circ\varphi^*g$. We claim that this implies
  that~$\im g\subseteq u^{N+1} \gNt$; if
  this is the case, then by induction on~$N$ we have
  $\im g\subseteq u^N \gNt$ for all~$N$,
  and so~$g=0$ (note that as~$\gNt$ is a finite
  projective~$\AAinf{A}$-module,
  it is $u$-adically complete and separated).

  To prove the claim, note that
  since~$\im g\subseteq u^N \gNt$ we
  have~$\im \varphi^*g\subseteq u^{pN}
  \gNt$. By~\cite[Lem.\
  5.2.6]{EGstacktheoreticimages} (and its proof), and the assumption
  that~$\gMt$ has height at most~$h$, we
  have~$\im \Phi_{\gMt}\supseteq u^{e(a+h-1)}\gMt$.  Since
  $g\circ\Phi_{\gMt}=\Phi_{\gNt}\circ\varphi^*g$, it follows that
  \[u^{e(a+h-1)}\im g\subseteq\im g\circ\Phi_{\gMt}= \im \Phi_{\gNt}\circ\varphi^*g
    \subseteq u^{pN} \gNt\]and since
  $\gNt$ is $u$-torsion free, we have
  $\im g\subseteq u^{pN-e(a+h-1)}
  \gNt$. Our assumption on~$N$ implies
  that $pN-e(a+h-1)\ge N+1$, as required.

  We now prove the existence of~$f'$. For any~$\AAinf{A}$-linear
  map $h: \gMt\to
  \gNt$ we
  set~$\delta(h):=\Phi_{\gNt}\circ\varphi^*h-h\circ\Phi_{\gMt}:\varphi^*\gMt\to\gNt$. We
  claim that for
  any~$s:\varphi^*\gMt\to
  u^{pN}\gNt$, we can find~$t:\gMt\to
  u^{N+1}\gNt$ with~$\delta(t)=s$. Given this claim, the existence of~$f'$ is immediate,
taking~$s=\delta(f)$ and~$f':=f-t$, so that~$\delta(f')=0$.

To prove the claim, we first prove the weaker claim that for
any~$s$ as above we can
find~$h:\gMt\to
u^{N+1}\gNt$
with~$\im(\delta(h)-s)\subseteq u^{p(N+1)}\gNt$. 
Admitting this second claim, we prove
the first claim by successive approximation: we set~$t_0=h$, so
that~$\im(\delta(t_0)-s)\subseteq u^{p(N+1)}\gNt$. Applying the second
claim again with~$N$ replaced by~$N+1$, and $s$ replaced
by~$s-\delta(t_0)$, we
find~$h:\gMt\to
u^{N+2}\gNt$
with~$\im(\delta(h)-(s-\delta(t_0)))\subseteq u^{p(N+2)}\gNt$. 
Setting~$t_1=t_0+h$, and
proceeding inductively, we obtain a Cauchy sequence converging (in
the~$u$-adically complete finite $\AAinf{A}$-module
$\Hom_{\AAinf{A}}(\gMt,\gNt)$)
to the required~$t$.

Finally, we prove this second
claim.
If~$h:\gMt\to
u^{N+1}\gNt$ then 
$\im\varphi^*h\subseteq u^{p(N+1)}\gNt$, so it suffices to show that we can
find~$h$ such that $h\circ\Phi_{\gMt}=-s$; but this is immediate,
because the cokernel of~$\Phi_{\gMt}$ is killed by~$u^{e(a+h-1)}$, and
$pN-e(a+h-1)\ge N+1$.
\end{proof}


The following lemma is proved by a reinterpretation of some of the arguments
of~\cite[\S2]{MR2745530}.  

\begin{lem}\label{lem: Caruso Liu Galois action on Kisin}
  For any fixed $a,h$, and any $N\ge e(a+h)/(p-1)$, there is a
  positive integer~$s(a,h,N)$ with the property that for any finite type
  $\cO/\varpi^a$-algebra~$A$, any
  projective Breuil--Kisin
  module~$\gM$  of height at most~$h$, and 
  any~$s\ge s(a,h,N)$, there is a unique continuous action of~ $G_{K_s}$ on
  $\gMt:=\AAinf{A}\otimes_{\gS_{A}}\gM$ which commutes with~$\varphi$
  and is semi-linear with respect to the natural action of~$G_{K_s}$
  on~$\AAinf{A},$ with the additional property that for
  all~$g\in G_{K_s}$ we have $(g-1)(\gM)\subset
  u^N\gMt$. 
\end{lem}
\begin{proof}By
  definition, to give a semi-linear action of~$G_{K_s}$  on~$\gMt$ is
  to give for each~$g\in G_{K_s}$ a morphism of Breuil--Kisin--Fargues
  modules \[\beta_g:g^*\gMt\to\gMt\] with the property that for
  all~$g,h\in G_{K_s}$, we have $\beta_{gh}=\beta_h\circ
  h^*\beta_g$. 
  Now, the requirement that
  $(g-1)(\gM)\subseteq u^N\gMt$ implies that if we have two such
  morphisms~$\beta_g,\beta'_g$ then
  $(\beta_g-\beta'_g)(g^*\gMt)\subseteq u^N\gMt$, so
  that~$\beta'_g=\beta_g$ by the uniqueness assertion of
  Lemma~\ref{lem: Frobenius amplification lemma over Ainf}.

  We now use the existence part of Lemma~\ref{lem: Frobenius
    amplification lemma over Ainf} to construct the required
  action. Since the homomorphism $\gS_A\into\AAinf{A}$ takes~$u$ to the
  Teichm\"uller lift of~$(\pi^{1/p^n})_n$, 
we can
  and do choose~$s(a,h,N)$ sufficiently large that for all~$s\ge
  s(a,h,N)$ and~$g\in
  G_{K_s}$, we have~$g(u)-u\in u^{pN}\AAinf{A}$; thus for
  all~$\lambda\in \gS_A$, we have $g(\lambda)-\lambda\in
  u^{pN}\AAinf{A}$. 

We claim that for each~$g\in G_{K_s}$, we may define
  an~$\AAinf{A}$-linear map  \[\alpha_g:g^*\gMt\to
    \gMt\] with the following two properties:
\begin{itemize}
	\item
$(\alpha_g\circ\Phi_{g^*\gMt}-\Phi_{\gMt}\circ\varphi^*\alpha_g)(\varphi^*g^*\gMt)\subseteq u^{pN}\gMt.$
\item
For each $m \in \gM,$ we have $\alpha_g(1\otimes m) - m \in u^N \gMt.$
\end{itemize}
It suffices to construct such maps in the case when~$\gM$ is free;
indeed, by Lemma~\ref{lem: adding projective Kisin to get free},
  we may write~$\gF=\gM\oplus\gP$ where~$\gF,\gP$ are respectively free and
  projective Breuil--Kisin modules of height at most~$h$, and given a
  morphism~$\alpha_g:g^*\gF^{\inf}\to\gF^{\inf}$ satisfying our
  requirements (where we write $\gF^{\inf}=\AAinf{A}\otimes_{\gS_{A}}\gF$),
  we may obtain the desired
  morphism~$\alpha_g:g^*\gMt\to\gMt$ by projection onto the
  corresponding direct summands.

If then~$\gF$ is a free Breuil--Kisin module of height at most~$h$, we choose a
basis~$f_1,\dots,f_d$ for~$\gF$, and then, for each~$g\in G_{K_s}$, we define
  an~$\AAinf{A}$-linear map  \[\alpha_g:g^*\gFt\to
    \gFt\] (which depends on our choice of basis) by
  %
  %
  %
  %
  \[\alpha_g(\sum_i\lambda_i\otimes f_i)=\sum_i\lambda_if_i.\] We now check
  that~$\alpha_g$ has the claimed properties.
  If~$\Phi_{\gF}(1\otimes f_i)=\sum_j\theta_{i,j}f_j$ (so
  that~$\theta_{i,j}\in\gS_A$), then 
  for any $\sum_i\lambda_i\otimes
    f_i\in\varphi^*g^*\gFt, $ we compute that 
    \[(\alpha_g\circ\Phi_{g^*\gFt}-\Phi_{\gFt}\circ\varphi^*\alpha_g)(\sum_i\lambda_i\otimes
    f_i)=\sum_{i,j}\lambda_i(g(\theta_{i,j})-\theta_{i,j})f_j\in
    u^{pN}\gFt,\] 
  while for any $\sum_i \lambda_i f_i \in \gF$
  (note then that each $\lambda_i \in \gS_A$),
  we compute that 
\begin{multline*}
    \alpha_g\bigl(1\otimes (\sum_i \lambda_i f_i)\bigr)
    - \sum_i \lambda_i f_i
\\
    = 
				 \alpha_g(\sum_ig(\lambda_i)\otimes
                                 g_i)-\sum_i\lambda_ig_i = 
				 \sum_i(g(\lambda_i)-\lambda_i)g_i
				 \in u^N \gFt;
				 \end{multline*}
  in both computations we have used the fact that, by our choice of~$s$,
  if~$\lambda\in\gS_A$, 
  then~$g(\lambda)-\lambda\in u^{pN}\gFt \subseteq u^N\gFt$.

 Returning to the general case of a projective Breuil--Kisin
 module~$\gM$ of height at most~$h$, it follows from 
  Lemma~\ref{lem: Frobenius amplification lemma
    over Ainf} that there is a unique 
morphism of Breuil--Kisin--Fargues modules
  \[\beta_g:g^*\gMt\to
 \gMt\]
  with~$\im(\beta_g-\alpha_g)\subseteq
  u^N\AAinf{A}\otimes_{\gS_A}\gM$. Since we have
  $\alpha_{gh}=\alpha_h\circ h^*\alpha_g$, it follows from this
  uniqueness that~$\beta_{gh}=\beta_h\circ h^*\beta_g$, 
  so we have
  constructed an action of~$G_{K_s}$ on~$\gMt$. It remains to check
  that~$(g-1)(\gM)\subseteq u^N\gMt$; for this, note that,
  for any $m \in \gM$, we have
    $$(g-1)m = \beta_g(1\otimes m)  - m = (\beta_g - \alpha_g)(1\otimes m)
    + \alpha_g(1\otimes m ) - m \in u^N \gMt;
    $$
here we have used the fact $\im(\beta_g-\alpha_g)\subseteq
  u^N\AAinf{A}\otimes_{\gS_A}\gM$, together with the second of the
  properties satisfied by the maps $\alpha_g$.
%
\end{proof}

\section{Breuil--Kisin--Fargues $G_K$-modules and canonical
  actions}\label{subsec: BKF action canonical}

In an earlier draft of this book, we claimed to show that for
Breuil--Kisin--Fargues $G_K$-modules admitting descents for
all~$\piflat$, the actions of the~$G_{K_{\piflat,s}}\subseteq G_K$
(for sufficiently large values of~$s$)
are necessarily the canonical actions considered in
Section~\ref{subsec:canonical actions}. We then used this claim to
establish finiteness properties of the stacks of
Breuil--Kisin--Fargues $G_K$-modules that we consider later in this
chapter. We are grateful to Dat Pham for pointing out a mistake in our
argument, and we do not know whether this is in fact
automatic. However, the following result of Caruso--Liu (which
motivated our consideration of canonical actions in the first place)
shows that the Breuil--Kisin--Fargues $G_K$-modules coming from
potentially semistable $G_K$-representations do necessarily have this
property; so we can (and do) impose this as an extra condition in the
definition of our stacks.

The following lemma is essentially the content
of~\cite[\S3]{MR2745530} (although they are concerned with getting
precise bounds on the size of~$s$, which is unimportant for
us). Recall that, by Theorem~\ref{athm: admits all descents if and only
  if semistable},
to any $\cO$-lattice in a
semistable~$E$-representation of~$G_K$ with Hodge--Tate weights
in~$[0,h]$
(for some integer~$h\ge 0$),
there corresponds a
Breuil--Kisin--Fargues $G_K$-module~$\gMt$ of height at most~$h$ which
admits all descents.

\begin{prop}
 \label{prop: Caruso Liu canonical action semistable}For any
 fixed~$K,a,h$ and~$N\ge e(a+h)/(p-1)$ there is a
 constant~$s'(K,a,h,N)\ge s(a,h,N)$ {\em (}where~$s(a,h,N)$ is as in
 Lemma~{\em \ref{lem: Caruso Liu Galois action on Kisin})} with the
 following property: for any Breuil--Kisin--Fargues $G_K$-module~$\gMt$
 corresponding as above to a semistable
$G_K$ representation with Hodge--Tate weights in~$[0,h]$, and for any
choice of~$\piflat$ with corresponding descent~$\gM_{\piflat}$, 
 if~$s\ge s'(K,a,h,N)$, then 
 the restriction to~$G_{K_{\piflat,s}}$ of the action of~$G_K$ 
 on $\gMt\otimes_{\cO}\cO/\varpi^a$ 
 agrees with the
 action obtained from~$\gM_{\piflat}$ by Lemma~{\em \ref{lem: Caruso Liu Galois action on Kisin}}.
\end{prop}
\begin{proof}
 This is essentially the content of~\cite[Prop.\
 3.3.1]{MR2745530}. 
 By the
 statement of Lemma~\ref{lem: Caruso Liu Galois action on Kisin}, it is enough to check that if~$s$ is sufficiently
 large, then for all~$g\in G_{K_{\piflat,s}}$ we have
 $(g-1)(\gM_{\piflat}\otimes_{\cO}\cO/\varpi^a)\subset u^N\Ainf\otimes_\gS\gM_{\piflat}\otimes_{\cO}\cO/\varpi^a$, where the
 action of~$g$ is via
the given action of~$G_K$ on~$\gMt$. 
 This is immediate
 from~\cite[Lem.\ 3.2.1]{MR2745530} (taking~$\mathfrak{L}=\gM_{\piflat}$, and
 noting that in their notation, $\widehat{\mathfrak{L}}$ is an
 extension of scalars of~$\gL$ to a subring
 of~$\AAinf{\cO/\varpi^a}$; this extension of scalars involves a
 twist by~$\varphi$, but since~$\varphi$ is a continuous automorphism
 of~$\Ainf$, this is harmless).
\end{proof}

\section[Stacks of semistable and crystalline Breuil--Kisin---Fargues
modules]{Stacks of semistable and crystalline Breuil--Kisin---Fargues
  modules \sectionmark{Stacks of semistable and crystalline BKF modules}}\sectionmark{Stacks of semistable and crystalline BKF modules}\label{sec:
  semistable stack}
We now define moduli stacks parameterizing the Breuil--Kisin--Fargues $G_K$-modules
admitting all descents that were introduced in Definition~\ref{defn:
  descending BKF to BK coefficients}. As discussed in
Section~\ref{subsec: BKF action canonical}, we also impose a condition
that the~$G_K$-actions agree with the canonical actions introduced in
Section~\ref{subsec:canonical actions}. To this end, we make the
following definition.
\begin{defn}
  \label{defn: canonical actions modulo anything}For each ~$K,a,h$
  and~$N\ge e(a+h)/(p-1)$, we fix a constant $s'(K,a,h,N)$ as in the
  statement of Proposition~\ref{prop: Caruso Liu canonical action
    semistable}. Let~$A$ be a $p$-adically complete $\cO$-algebra
  which is topologically of finite type, and let $\gMt$ be a
  projective Breuil--Kisin--Fargues $G_K$-module with $A$-coefficients
  which admits all descents. Then we say that the $G_K$-action on
  \index{canonical $G_K$-action}
  $\gMt$ is {\em canonical} if for any~$a\ge 1$, any~$\piflat$, and
  any~$s\ge s'(K,a,h,N)$, the restriction to~$G_{K_{\piflat,s}}$ of
  the action of~$G_K$ on $\gMt\otimes_{\cO}\cO/\varpi^a$ agrees with
  the action obtained from~$\gM_{\piflat}$ by Lemma~{\ref{lem:
      Caruso Liu Galois action on Kisin}}.
\end{defn}
\begin{rem}
  \label{rem: why canonical actions}The point of this definition is
  that it is immediate from
  Proposition~\ref{prop: Caruso Liu canonical action semistable} that
  the Breuil--Kisin--Fargues modules coming from lattices in
  semistable $G_K$-representations have canonical $G_K$-actions.
\end{rem}

We now define our moduli stacks of Breuil--Kisin--Fargues $G_K$-modules.
\begin{defn}\label{defn:stacks of semistable cris BKF modules}For any~$h\ge 0$ we let~$\cC_{d,\semis,h}^a$ 
  denote the
 limit preserving category of groupoids over
  $\Spec\cO/\varpi^a$ determined by decreeing, for any finite
  type $\cO/\varpi^a$-algebra~$A$, that
  $\cC_{d,\semis,h}^a(A)$ is the
  groupoid of rank~$d$ projective Breuil--Kisin--Fargues $G_K$-modules with
  $A$-coefficients, which are of height at most~$h$, which admit
  all
  descents, and whose $G_K$-action is canonical.

  We let $\cC_{d,\crys,h}^a$ denote the limit preserving subcategory of
  groupoids of~$\cC_{d,\semis,h}^a$ 
  consisting of those
  Breuil--Kisin--Fargues $G_K$-modules~$\gMt$ which are furthermore crystalline. 
  \index{$\cC_{d,\semis,h}$} \index{$\cC_{d,\crys,h}$}
  We let~$\cC_{d,\semis,h}:=\varinjlim_a\cC_{d,\semis,h}^a$,
  and let~$\cC_{d,\crys,h}:=\varinjlim_a\cC_{d,\crys,h}^a$.
\end{defn}
\begin{rem}
  This definition uniquely determines limit preserving categories
  fibred in
  groupoids over $\cO$ by~\cite[Lem.\
  2.5.4]{EGstacktheoreticimages}.
  We will see shortly that the inductive limits~$\cC_{d,\semis,h}$ and~$\cC_{d,\crys,h}$ are in fact $p$-adic formal algebraic stacks
  (Theorem~\ref{thm: semi stable stack is formal algebraic}).
\end{rem}

If $A$ is a $p$-adically complete $\cO$-algebra,
then as usual we write $\cC_{d,\semis,h}(A)$ 
(resp. $\cC_{d,\crys,h}(A)$) to denote the groupoid
of morphisms from $\Spf A$ (defined using the $p$-adic topology on $A$)
to $\cC_{d,\semis,h}$ (resp.\ $\cC_{d,\crys,h}$).
More concretely, we have
$$\cC_{d,\semis,h}(A) = \varprojlim_a \cC_{d,\semis,h}^a(A/\varpi^a),$$
and similarly for $\cC_{d,\crys,h}(A)$. 
We will typically apply this convention in the case when $A$ is furthermore
topologically of finite type over $\cO$.  In this case each
of the quotients $A/\varpi^a$ is of finite type over $\cO/\varpi^a$
(by the definition of topologically of finite type),
and so the objects of $\cC_{d,\semis,h}^a(A/\varpi^a)$
admit a concrete interpretation as Breuil--Kisin--Fargues $G_K$-modules over~$A$.
The objects of $\cC_{d,\semis,h}(A)$ and of $\cC_{d,\crys,h}(A)$
then similarly admit a concrete interpretation.

\begin{lem}
  \label{lem: finite O algebra points of semistable Kisin stack}Let
  $A$ be a $p$-adically complete
  $\cO$-algebra which is topologically of finite type.
Then $\cC_{d,\semis,h}(A)$ 
  is the groupoid of rank $d$ projective
  Breuil--Kisin--Fargues $G_K$-modules with $A$-coefficients, which
  are of height at most~$h$, which admit
  all
  descents, and whose $G_K$-action is canonical. In addition,
  $\cC_{d,\crys,h}(A)$ is the subgroupoid consisting of those
  Breuil--Kisin--Fargues $G_K$-modules~$\gMt$ which are furthermore
  crystalline.
\end{lem}
\begin{proof}
  This follows from Lemma~\ref{lem: can check projectivity of Kisin and etale phi modules modulo p^n}.
\end{proof}

The following lemma is crucial: it shows that Breuil--Kisin--Fargues modules which
admit all descents give rise to~$(\varphi,\Gamma)$-modules.

\begin{lemma}
\label{lem:C to X map}
There is a natural morphism $\cC_{d,\semis,h}\to\cX_d$, which
on $\Spf A$-points {\em (}for $p$-adically complete $\cO$-algebras $A$ that are topologically
of finite type{\em )} is given by extending scalars to~$W(\C^\flat)_A$.
\end{lemma}
\begin{proof}
As indicated in the statement of the lemma, this morphism is defined, for finite type 
$\cO/\varpi^a$-algebras, via
$\gMt \mapsto W(\C^\flat)_A \otimes_{\AAinf{A}} \gMt,$
with the target object being regarded as an $A$-valued point of $\cX_d$
via the equivalence of Proposition~\ref{prop: equivalences of categories to Ainf}.
Since both $\cC_{d,\semis,h}$ and $\cX_d$ are limit preserving
(the former by its very construction, and the latter by 
Lemma~\ref{lem:X is limit preserving}), it follows from~\cite[Lem.\ 2.5.4, Lem.\
  2.5.5~(1)]{EGstacktheoreticimages} that this construction determines
a morphism $\cC_{d,\semis,h}\to\cX_d$.
\end{proof}

\begin{remark}
\label{rem:restricting to finite type algebras}
The proof of Lemma~\ref{lem:C to X map} illustrates
a general principle (which was also applied in the proofs
of the results in Section~\ref{subsec:restriction}),
namely that to construct a morphism between limit preserving
stacks over~$\Spec \cO/\varpi^a$,
or over~$\Spf \cO$,
it suffices to define the intended morphism on finite
type $\cO/\varpi^a$-algebras (with $a$ either fixed or varying,
depending on which case we are in).

Similarly, if $P$ is any property of such a morphism that is preserved under
base-change, and that is tested by pulling back over $A$-valued points,
then to test for the property~$P$, it suffices to consider such pull-backs
in the case when $A$ is of finite type over~$\cO/\varpi^a$.

We will apply these principles consistently throughout the remainder of this section.
\end{remark}

For any~$h\ge 0$ and any choice of~$\pi^\flat$, we write~$\cC_{\piflat,d,h}$ for the moduli stack of
rank~$d$ projective Breuil--Kisin modules over~$\gS_{\pi^\flat,A}$ of height at
most~$h$, and~$\cR_{\piflat,d}$ for the corresponding stack of rank $d$
projective \'etale
$\varphi$-modules.
We also write $\cC^a_{\piflat,d,h}$ and $\cR^a_{\piflat,d}$ for their
base-change over $\cO/\varpi^a$, for any $a \geq 1$.

Recall that,
by Proposition~\ref{prop: natural morphism X to R}, we have a natural morphism
$\cX_{K,d}\to\cR_{\piflat,d}$, which for each~$s\ge 1$ can be factored,
via Lemma~\ref{lem: restriction from X to Xs},
as\[\cX_{K,d}\to\cX_{K_{\piflat,s},d}\to\cR_{\piflat,d}. \]
\begin{prop}
  \label{prop: morphism C to Xs}For any fixed $a,h$, and any~$N$ and~$s(a,h,N)$ as in
  Lemma~{\em \ref{lem: Caruso Liu Galois action on Kisin}}, then for
  any~$s\ge s(a,h,N)$ there is a canonical morphism
  $\cC^a_{\piflat,d,h}\to\cX^a_{K_{\piflat,s},d}$ obtained from the canonical
  action of Lemma~{\em \ref{lem: Caruso Liu Galois action on Kisin}}. This
  morphism fits into a commutative
  triangle
  \[\xymatrix{& \cC^a_{\piflat,d,h}\ar[dl]\ar[d]\\ \cX^a_{K_{\piflat,s},d}\ar[r]&\cR^a_{\piflat,d}}\]
\end{prop}
\begin{proof}
As discussed in Remark~\ref{rem:restricting to finite type algebras},
it is enough to show that if~$A$ is a finite type
  $\cO/\varpi^a$-algebra, and~$\gM$ is a finite projective \'etale Breuil--Kisin
  module with $A$-coefficients, then we can canonically 
  extend the action of~$G_{K_{\pi^\flat,\infty}}$ on
  $W(\C^\flat)_A\otimes_{\gS_{\piflat,A}}\gM$ to an action
  of~$G_{K_{\pi^\flat,s}}$, if~$s\ge s(a,h,N)$. Such an extension
is provided by Lemma~\ref{lem: Caruso Liu Galois action on Kisin}.
\end{proof}

\begin{lem}\label{lem: C to Xs is proper}
  The morphism $\cC^a_{\piflat,d,h}\to\cX^a_{K_{\piflat,s},d}$ of
  Proposition~{\em \ref{prop: morphism C to Xs}} is representable by
  algebraic spaces, proper, and of finite presentation.
\end{lem}
\begin{proof} This follows by a standard graph argument.  Namely,
	we write this morphism as the
  composite \[\cC^a_{\piflat,d,h}\to\cC^a_{\piflat,d,h}\times_{\cR^a_{\piflat,d}}\cX^a_{K_{\piflat,s},d}
    \to\cX^a_{K_{\piflat,s},d}. \] The first arrow is a closed immersion, 
  being
  the base change of the diagonal morphism $\cX_d\to\cX_d\times_{\cR_{\piflat,d}}\cX_d$, which is a closed immersion
  by Proposition~\ref{prop: diagonal of X to R}. The second arrow
  is representable by  algebraic spaces  and proper, 
because it is the base change of the morphism
  $\cC^a_{\piflat,d,h}\to{\cR^a_{\piflat,d}}$ of Theorem~\ref{thm: C is a $p$-adic formal algebraic stack and related
    properties}.
The  composite morphism  is thus representable by algebraic spaces
and  proper.  Such a morphism is in particular of finite type
and separated (and so also quasi-separated), and hence Lemma~\ref{lem:finite presentation from closed immersion}
below allows us to conclude that it is also of finite presentation
(once we take into account Lemma~\ref{lem:X is limit preserving},
 which shows that its target is limit preserving).
\end{proof}

\begin{lem}
\label{lem:finite presentation from closed immersion}
Let $\cZ \hookrightarrow \cX$ be a morphism of 
stacks over a locally
Noetherian base scheme~$S$ 
which is representable by algebraic spaces, of finite type,
and quasi-separated.
If $\cX$ is limit preserving, 
then this morphism is furthermore of finite presentation.
\end{lem}
\begin{proof}
It is enough to show, for any morphism $T \to \cX$ whose source
is an affine $S$-scheme $T$, that the pulled back morphism
(of algebraic spaces)
$T\times_{\cX} \cZ \to T$ is of finite presentation.  Since $\cX$ is
limit preserving,
 the morphism $T \to \cX$ factors through a morphism
of $S$-schemes
$T \to T'$, where $T'$ is affine and of finite presentation
over~$S$.  Thus, replacing $T$ by~$T'$, we may assume
that $T$ is of finite presentation over~$S$, and consequently Noetherian.

The base-changed morphism $T\times_{\cX} \cZ \to T$ is thus a
finite type and quasi-separated morphism  from
an algebraic space to a Noetherian scheme, and hence is of finite presentation
by~\cite[\href{https://stacks.math.columbia.edu/tag/06G4}{Tag 06G4}]{stacks-project}.
This proves the lemma.
\end{proof}

\begin{df}
	\label{def:2-Cartesian}
The morphisms of Lemma~\ref{lem: restriction from X to Xs} and
Proposition~\ref{prop: morphism C to Xs} allow us to define a
stack~$\cC^a_{\piflat,s,d,h}$ by the requirement that it fits into a  \index{$\cC^a_{\piflat,s,d,h}$}
2-Cartesian diagram \numequation\label{eqn: 2 cartesian canonical C}\xymatrix{\cC^a_{\piflat,s,d,h} \ar[r]\ar[d] &  \cC^a_{\piflat,d,h}\ar[d]\\
  \cX^a_{K,d}\ar[r] & \cX^a_{K_{\piflat,s},d}}\end{equation}
The stack~$\cC^a_{\piflat,s,d,h}$ classifies rank $d$
projective Breuil--Kisin modules $\gM$ of height at most $h$
over $\gS_{\piflat,A}$
equipped with an extension of the canonical $G_{K_{\piflat,s}}$-action
on $W(\C^\flat)_A\otimes_{\gS_{\piflat,A}} \gM$ to an action of $G_K$
(making it an \'etale $(\varphi,G_K)$-module).

Since the lower horizontal arrow in~\eqref{eqn: 2 cartesian canonical C}
is representable by algebraic spaces and of 
finite presentation, by Lemma~\ref{lem: restriction from X to Xs},
and since~$\cC^a_{\piflat,d,h}$ is an
algebraic stack of finite presentation over $\cO/\varpi^a$,
by Theorem~\ref{thm: C is a $p$-adic formal algebraic stack and related
    properties}, we see that~$\cC^a_{\piflat,s,d,h}$   is also an algebraic
stack of finite presentation over $\cO/\varpi^a$.
Since the right hand vertical arrow in this diagram is representable by
algebraic spaces, proper, and of
finite presentation,
by Lemma~\ref{lem: C to Xs is proper}, so is the left
hand vertical arrow.
\end{df}

\begin{lem}
\label{lem:auxiliary object has affine diagonal}
The diagonal of~$\cC^a_{\piflat,s,d,h}$ is affine and of finite
presentation.
\end{lem}
\begin{proof}
As we already noted, 
the morphism $\cC^a_{\piflat,s,d,h} \to \cX_{K,d}$ is representable by algebraic
spaces and proper (and so in particular is both separated and of finite type);
the claim of the lemma thus follows from Proposition~\ref{prop: properties
of X pulled back from R}, together with Lemma~\ref{lem:controlling diagonal} below.
%
\end{proof}

\begin{lem}
\label{lem:controlling diagonal}
Suppose that $f:\cX \to \cY$ is a morphism of stacks over a base scheme $S$
which is representable
by algebraic spaces, separated, and of finite type.  Suppose also that
the diagonal $\Delta_{\cY}: \cY \to \cY\times_S \cY$ is affine and of finite
presentation.  Then the diagonal $\Delta_{\cX}: \cX \to \cX \times_S \cX$
is affine and of finite presentation.
\end{lem}
\begin{proof}
The diagonal morphism 
$$\cX \to \cX\times_S \cX$$
may be factored as the composite of the relative diagonal
\numequation
\label{eqn:first factor}
\cX \to \cX \times_{\cY} \cX
\end{equation}
and the morphism
\numequation
\label{eqn:second factor}
\cX\times_{\cY} \cX \to \cX\times_S \cX,
\end{equation}
which is a base change of the diagonal morphism~$\Delta_{\cY}$. 
Our assumption on $\Delta_{\cY}$, then, implies that
the morphism~\eqref{eqn:second factor} 
is affine and of finite presentation.
Since $\cX \to \cY$ is representable by algebraic spaces and separated,
its diagonal~\eqref{eqn:first factor} is a closed immersion, and thus affine. 
Since it is furthermore
of finite type, the morphism~\eqref{eqn:first factor} is also 
of finite presentation, by Lemma~\ref{lem:finiteness of diagonals}.
Thus $\Delta_{\cX}$ is the composite of morphisms that are affine and
of finite presentation, and the lemma follows.
\end{proof}

\begin{prop}
  \label{prop: BKF closed in BK}For each~$\piflat$ and each~$s\ge s'(K,a,h,N)$, there are  natural
  closed immersions $\cC^a_{d,\crys,h}\to\cC^a_{d,\semis,h}\to
  \cC^a_{\piflat,s,d,h}$. In particular, $\cC^a_{d,\semis,h}$
  and~$\cC^a_{d,\crys,h}$ are
 algebraic stacks of finite presentation over~$\cO/\varpi^a$, and have
 affine diagonals. 
\end{prop}
\begin{proof} 
We have already observed that $\cC^a_{\piflat,s,d,h}$ is of
finite presentation over~$\cO/\varpi^a$, and has affine diagonal by
Lemma~\ref{lem:auxiliary object has affine diagonal}. Once we
construct the claimed closed immersions, 
the claims of finite presentation will follow from
Lemma~\ref{lem:finite presentation from closed immersion}; and since closed
immersions are in particular monomorphisms,
the diagonals of  $\cC_{d,\semis,h}$ and  $\cC_{d,\crys,h}$  are then obtained via base change
from the diagonal of~ $\cC^a_{\piflat,s,d,h}$, and are in particular affine.
We thus focus on constructing these closed immersions.
Throughout the proof, we take into
Remark~\ref{rem:restricting to finite type algebras},
which allows us to restrict our attention to the points
of the various stacks in question that are defined
over finite type $\cO/\varpi^a$-algebras.

Lemma~\ref{lem:C to X map}
  constructs a morphism $\cC^a_{d,\semis,h}\to\cX^a_{K,d}$, given
  by extending scalars to~$W(\C^\flat)$, and it follows from
  Lemma~\ref{lem: uniqueness of descent for G K infty Ainf} that there
  is a natural morphism $\cC^a_{d,\semis,h}\to\cC^a_{\piflat,d,h}$, defined via
  $\gMt \mapsto \gM_{\piflat}$. The composite morphisms
  $\cC^a_{d,\semis,h} \to \cC^a_{\piflat,d,h} \to
  \cX^a_{K_{\piflat,s},d}$ and
  $\cC^a_{d,\semis,h} \to \cX^a_{K,d}\to \cX^a_{K_{\piflat,s},d}$
  coincide by definition (see Definitions~\ref{defn: canonical actions
    modulo anything} and~\ref{defn:stacks of semistable cris BKF modules}). 
  Thus these morphisms induce a morphism
\numequation
\label{eqn:induced map}
\cC^a_{d,\semis,h}\to \cC^a_{\piflat,s,d,h}.
\end{equation}

  To see that~\eqref{eqn:induced map}
is a monomorphism, it is enough to note
  that if~$A$ is any finite type $\cO/\varpi^a$-algebra~$A$,
  and~$\gMt$ is any Breuil--Kisin--Fargues module over $A$ which
  admits all descents and whose $G_K$-action is canonical,
 corresponding to an $A$-valued point of~$\cC^a_{d,\semis,h},$
  then~$\gMt=\AAinf{A}\otimes_{\gS_{\piflat,A}}\gM_{\piflat}$ is
  determined by~$\gM_{\piflat}$, and the $G_K$-action on~$\gMt$ is
  determined by the $G_K$-action on
  $W(\C^\flat)_A\otimes_{\AAinf{A}}\gMt$.

We now show that~\eqref{eqn:induced map}
is a closed immersion. It is enough to do
  this after pulling back to some finite type
  $\cO/\varpi^a$-algebra~$A$, where we need to show that the conditions that
$\gMt=\AAinf{A}\otimes_{\gS_{\piflat,A}}\gM_{\piflat}$ is $G_K$-stable
and admits all descents, and that the $G_K$-action is canonical, are closed conditions. We begin with the first
of these. It is enough to show that for each~$g\in G_K$, the condition
that $g(\gMt)\subset\gMt$ is closed. This condition is equivalent to
the vanishing of the composite morphism \[g^*\gMt\to
  \gMt[1/u]\to\gMt[1/u]/\gMt,\]and this is a closed condition by
Lemma~\ref{lem: inclusion lattices closed condition}. Let~$\cC^a_{d,G_K}$ denote the closed substack of
$\cC^a_{\piflat,s,d,h}$ for which~$\gMt$ is $G_K$-stable. 

 We next show that for each choice of uniformizer~$\pi'$, and
each~$(\pi')^\flat$, the condition that~$\gMt$ admits a descent
to~$\gS_{A,(\pi')^\flat}$ is a closed condition. To this end, note
that by Proposition~\ref{prop: equivalences of categories to Ainf},
$\gMt[1/u]$ descends uniquely to~$\gS_{(\pi')^\flat,A}[1/u]$, so we
have a morphism $\cC^a_{d,G_K}\to\cR^a_{(\pi')^\flat,d}$. Then the
fibre product
$\cC^a_{d,G_K}\times_{\cR^a_{(\pi')^\flat,d}}\cC^a_{(\pi')^\flat,d,h}$
is the moduli space of pairs~$(\gM,\gM')$ where~$\gM := \gM_{\piflat}$ is as above
(an object classified by~$\cC^a_{d,G,K}$),
and~$\gM'$ is a Breuil--Kisin module for~$\gS_{(\pi')^\flat,A}$ of
height at most~$h$ which
satisfies \[W(\C^\flat)_A\otimes_{\gS_{(\pi')^\flat,A}}\gM'=W(\C^\flat)_A\otimes_{\gS_{\pi^\flat,A}}\gM.\]
Consider the substack~$\cC^a_{d,G_K,(\pi')^\flat}$ of this fibre
product satisfying the condition that
\numequation\label{eqn: pi pi prime
  agree}\AAinf{A}\otimes_{\gS_{(\pi')^\flat,A}}\gM'=\gMt.\end{equation} Note that this is exactly the condition that~$\gMt$ admits
a descent for~$(\pi')^\flat$, and that~$\gM'$ is this descent (which
is uniquely determined by Lemma~\ref{lem: uniqueness of descent for G K infty Ainf}),
so the
projection $\cC^a_{d,G_K,(\pi')^\flat}\to \cC^a_{d,G_K}$ is a
monomorphism, and we need to show that it is a closed immersion. Since the projection
$\cC^a_{d,G_K}\times_{\cR^a_{(\pi')^\flat,d}}\cC^a_{(\pi')^\flat,d,h}\to
\cC^a_{d,G_K}$ is proper (being a base change of the proper morphism
$\cC^a_{(\pi')^\flat,d,h}\to\cR^a_{(\pi')^\flat,d}$), and since proper
monomorphisms are closed immersions, it is enough to show that
$\cC^a_{d,G_K,(\pi')^\flat}$ is a closed substack of
$\cC^a_{d,G_K}\times_{\cR^a_{(\pi')^\flat,d}}\cC^a_{(\pi')^\flat,d,h}$;
that is, we must show that~\eqref{eqn: pi pi prime
  agree} is a closed condition. This again follows from
Lemma~\ref{lem: inclusion lattices closed condition}.

We now consider the closed substack  of $\cC^a_{\piflat,s,d,h}$ for
which~$\gMt$ is $G_K$-stable and admits a descent for
each~$(\pi')^\flat$. We need to show that for each~$(\pi')^\flat$, the
further conditions that  \[\gM_{\piflat}/[\piflat]\gM_{\piflat}=\gM_{\piflatprime}/[\piflatprime]\gM_{\piflatprime}\]and  \[\varphi^*\gM_{\piflat}/E_{\piflat}\varphi^*\gM_{\piflat}=\varphi^*\gM_{\piflatprime}/E_{\piflatprime}\varphi^*\gM_{\piflatprime}\]
are closed conditions.

The arguments in both cases are very similar (and in turn are similar
to the proof of Lemma~\ref{lem: inclusion lattices closed condition}),
so we only give the argument in the second case, leaving the first to
the reader. Each side is a projective $\cO_K\otimes_{\Zp}A$-submodule
of the projective $\cO_{\C}\otimes_{\Zp}A$-module
$\cO_{\C}\otimes_{\AAinf{A},\theta}\varphi^*\gMt$, and indeed spans this module
after extension of scalars to~$\cO_{\C}$. By symmetry, it is enough to
show that for each
element~$m\in \varphi^*\gM_{\piflatprime}/E_{\piflatprime}\varphi^*\gM_{\piflatprime}$,
the condition that $m\in \varphi^*\gM_{\piflat}/E_{\piflat}\varphi^*\gM_{\piflat}$ is
closed.

Write $P:=\varphi^*\gM_{\piflat}/E_{\piflat}\varphi^*\gM_{\piflat}$, and choose
a finite projective $\cO_K\otimes_{\Zp}A$-module~$Q$ such
that~$F:=P\oplus Q$ is free. We can think of~$m$ as an element of
$\cO_{\C}\otimes_{\cO_K}P$ and thus as an element of
$\cO_{\C}\otimes_{\cO_K}F$, and it is enough to check that the
condition that~$m\in F$ is closed. Choosing a basis for~$F$, we reduce
to the case that~$F$ is one-dimensional, and thus to showing that the
condition that an element of~$\cO_{\C}\otimes_{\Zp}A$ lies
in~$\cO_{K}\otimes_{\Zp}A$ is closed. Since~$\cO_{\C}$ is a torsion-free
and thus flat~$\cO_K$-module, $\cO_{\C}/\varpi^a$ is flat and thus
free as a module for the Artinian ring $\cO_K/\varpi^a$. Thus
$\cO_{\C}\otimes_{\Zp}A$ is a free~$\cO_{K}\otimes_{\Zp}A$-module, and
the result follows by choosing a basis.

We now need to show that the condition that the~$G_K$-action is
canonical is a closed condition. Explicitly, by Lemma~\ref{lem: Caruso
  Liu Galois action on Kisin}, we need to show that for each~$b\le a$,
each $N\ge e(b+h)/(p-1)$, each~$s\ge s'(K,b,h,N)$, 
each~$\piflat$,
and each~$g\in G_{K_{\piflat,s}}$,
the condition that 
\[(g-1)(\gM_{\piflat}\otimes_{\cO}\cO/\varpi^b)\subset
  u^N(\gMt\otimes_{\cO}\cO/\varpi^b)\] is closed. This
follows from Lemma~\ref{lem: element of lattice closed condition}.

Finally, it remains to show that the monomorphism
~$\cC^a_{d,\crys,h}\to \cC^a_{d,\semis,h}$ is a closed
immersion; that is, for each~$g\in G_K$ and each~$\piflat$ we need to
show that the condition that \[(g-1)(\gM_{\piflat})\subseteq \varphi
  ^{-1}(\mu)[\piflat]\gMt \] is a closed condition. It is enough to show that
for each~$m\in\gM_{\piflat}$, the condition that
$(g-1)(m)\in \varphi^{-1}(\mu)[\piflat]\gMt$ is closed. Since $ \varphi^{-1}(\mu)[\piflat]\gMt$ is a
finite projective $\AAinf{A}$-module, this follows from another
application of Lemma~\ref{lem: element of lattice closed condition}.
\end{proof}




\begin{thm}\label{thm: semi stable stack is formal algebraic}
  $\cC_{d,\semis,h}$ is a $p$-adic formal algebraic stack
of finite presentation,
  as is its closed substack~$\cC_{d,\crys,h}$. In addition,
both of these stacks have affine diagonal,
  and each of the 
  morphisms
$\cC_{d,\semis,h}\to\cX_{K,d}$
 and
$\cC_{d,\crys,h}\to\cX_{K,d}$
is representable by algebraic
  spaces, proper, and of finite presentation. 
\end{thm}
\begin{proof}Since ~$\cC_{d,\crys,h}$ is a closed substack
  of~$\cC_{d,\semis,h}$, it suffices to prove the statements of the lemma
  for~$\cC_{d,\semis,h}$.
(Here we use the fact that a closed immersion is a monomorphism
to see that the diagonal of $\cC_{d,\crys,h}$ is obtained via base change
from the diagonal of $\cC_{d,\semis,h}$, while we use Lemma~\ref{lem:finite
presentation from closed immersion} to transfer finite presentation
properties from $\cC_{d,\semis,h}$ to $\cC_{d,\crys,h}$).

 By Proposition~\ref{prop: criterion for
    p-adic formal algebraic stack}, to see that~$\cC_{d,\semis,h}$ is
  a $p$-adic formal algebraic stack of finite presentation, it is enough to show
  that each~$\cC_{d,\semis,h}^a$ is an algebraic stack of finite
presentation over $\cO/\varpi^a$,
  which was proved in Proposition~\ref{prop: BKF closed in BK}.
  Since each $\cC_{d,\semis,h}^a$ has affine diagonal, by the same
proposition, we see as well that $\cC_{d,\semis,h}$ has affine diagonal.

  That $\cC_{d,\semis,h}\to\cX_{K,d}$ is representable by algebraic
  spaces, proper, and of finite presentation also follows from
  Proposition~\ref{prop: BKF closed in BK}. Indeed, for each~$a$ the
  morphism $\cC_{d,\semis,h}^a\to\cX_{K,d}$ factors as
  $\cC_{d,\semis,h}^a\to\cC^a_{\piflat,s,d,h}\to \cX_{K,d}^a$, where the
  first morphism is a closed immersion of finite type algebraic
  stacks, so it is enough to prove the same properties for the
  morphism~$\cC^a_{\piflat,s,d,h}\to \cX_{K,d}^a$; as explained in
  Definition~\ref{def:2-Cartesian}, this follows from Lemma~\ref{lem:
    C to Xs is proper}.
\end{proof}

\subsection{Potentially crystalline and potentially semistable
	stacks}
We now consider the corresponding potentially semistable and
potentially crystalline versions of these moduli stacks of
Breuil--Kisin--Fargues modules. 

  \begin{df}\label{defn: semistable crys stacks L over K}
For a Galois extension~$L/K$,
define stacks~$\cC_{d,\semis,h}^{L/K}$ and \index{$\cC_{d,\semis,h}^{L/K}$}
\index{$\cC_{d,\crys,h}^{L/K}$}
$\cC_{d,\crys,h}^{L/K}$ as follows: for each~$a\ge 1$ we
let~$\cC_{d,\semis,h}^{L/K,a}$ denote the limit preserving category of
groupoids over $\Spec\cO/\varpi^a$ determined by decreeing,
for any finite type $\cO/\varpi^a$-algebra~$A$,
that $\cC_{d,\semis,h}^{L/K,a}(A)$ is the groupoid of
Breuil--Kisin--Fargues $G_K$-modules with $A$-coefficients, which are
of height at most~$h$, which admit all descents over~$L$, and whose
$G_L$-actions are canonical. 

We let
$\cC_{d,\crys,h}^{L/K,a}$ denote the limit preserving subcategory of
groupoids of~$\cC_{d,\semis,h}^a$ consisting of those
Breuil--Kisin--Fargues $G_K$-modules~$\gMt$ for which the action
of~$G_L$ is crystalline.

We let~$\cC_{d,\semis,h}^{L/K}:=\varinjlim_a\cC_{d,\semis,h}^{L/K,a}$,
and let~$\cC_{d,\crys,h}^{L/K}:=\varinjlim_a\cC_{d,\crys,h}^{L/K,a}$.
\end{df}

The following lemma is proved in exactly the same way as Lemma~\ref{lem: finite O algebra points of  semistable Kisin stack}.

\begin{lem}
  \label{lem: finite O algebra points of potentially semistable Kisin
    stack}Let $A$ be a $p$-adically complete
  $\cO$-algebra which is topologically of finite type over~$\cO$. Then $\cC_{d,\semis,h}^{L/K}(A)$ is the groupoid of
  Breuil--Kisin--Fargues $G_K$-modules with $A$-coefficients, which
  are of height at most~$h$, and which admit all descents over~$L$, and whose
$G_L$-actions are canonical; and
  $\cC_{d,\crys,h}^{L/K}(A)$ is the subgroupoid of those
  Breuil--Kisin--Fargues $G_K$-modules~$\gMt$ which are furthermore
  crystalline.
\end{lem}

\begin{prop}
  \label{prop: pst and pcrys stacks are padic}For any Galois
  extension~$L/K$, any $d$, and any~$h$, 
  ~$\cC_{d,\semis,h}^{L/K}$ is a 
  $p$-adic formal algebraic stack of finite presentation, and
  ~$\cC_{d,\crys,h}^{L/K}$ is a closed substack, which is thus
again a $p$-adic formal algebraic stack of finite presentation.
Both $\cC_{d,\semis,h}^{L/K}$ and $\cC_{d,\crys,h}^{L/K}$ have affine diagonal,
and the natural morphisms
  $\cC_{d,\semis,h}^{L/K}\to\cX_{K,d}$
  and $\cC_{d,\crys,h}^{L/K}\to\cX_{K,d}$
are representable by algebraic
  spaces, proper, and of finite presentation.
\end{prop}
\begin{proof}
  We have a natural
  morphism~$\cC_{d,\semis,h}^{L/K}\to\cC_{d,\semis,h}^{L/L}$ (given by
  restricting the action of~$G_K$ to~$G_L$; the target is of course
  just the stack denoted~$\cC_{d,\semis,h}$ above, but for~$L$ instead
  of~$K$), and a natural morphism
  $\cC_{d,\semis,h}^{L/K}\to\cX_{K,d}$, and thus a natural morphism
  \numequation\label{eqn: mapping pst BKF to fibre
    product}\cC_{d,\semis,h}^{L/K}\to \cC_{d,\semis,h}^{L/L}\times_{\cX_{L,d}}\cX_{K,d}   \end{equation}

 We claim that~\eqref{eqn: mapping pst BKF to fibre
    product} is a closed immersion. Given this, the proposition
  follows. Indeed, by Theorem~\ref{thm: semi stable stack is formal algebraic}, $\cC_{d,\semis,h}^{L/L}$ is a
  $p$-adic formal algebraic stack of finite presentation, and the morphism
  $\cC_{d,\semis,h}^{L/L}\to\cX_{L,d}$ is representable by algebraic
  spaces, proper, and of finite presentation. Furthermore,
  $\cX_{K,d}\to\cX_{L,d}$ is representable by algebraic spaces 
  and of
  finite presentation by Lemma~\ref{lem: restriction from X to Xs}. That~$\cC_{d,\crys,h}^{L/K}$ is a closed substack of
  $\cC_{d,\semis,h}^{L/K}$ follows from Proposition~\ref{prop: BKF
    closed in BK}.
  The claims on the diagonals are proved by appeal to Lemma~\ref{lem:controlling
diagonal} (\emph{cf.}\ the proof of Lemma~\ref{lem:auxiliary object has affine diagonal}).

  It remains to prove the claim. It suffices to prove this after
  pulling back via a morphism $\Spec A\to\cX_{L,d}$ for some finite
  type $\cO/\varpi^a$-algebra~$A$. Unwinding the definitions, we have
  to show that if~$\gMt$ is a Breuil--Kisin--Fargues $G_L$-module with
  $A$-coefficients, with a compatible action of~$G_K$ on
  $W(\C^{\flat})_A\otimes_{\AAinf{A}}\gMt$, then the condition
  that~$\gMt$ is $G_K$-stable is a closed condition. It suffices to
  show that for each~$g\in G_K$, the condition that
  $g(\gMt)\subseteq\gMt$ is a closed condition; as in the proof of
  Proposition~\ref{prop: BKF closed in BK}, this follows from
  Lemma~\ref{lem: inclusion lattices closed condition}, applied to the composite morphism \[g^*\gMt\to
  \gMt[1/u]\to\gMt[1/u]/\gMt.\qedhere\]
\end{proof}

We end this section with the following lemma, which will be used in
the proof of Proposition~\ref{prop: versal rings for pst stacks}.

\begin{lem}
  \label{lem: algebraising C}Let~$R$ be a complete local
  Noetherian~$\cO$-algebra with residue field~$\F$, 
  together with a morphism $\Spf
  R\to\cX_{K,d}$, and
  let~$\widehat{C}:=\cC_{d,\semis,h}^{L/K}\times_{\cX_{K,d}}\Spf R$. Then
  there is a projective morphism of schemes $C\to\Spec R$ whose
  $\m_R$-adic completion is isomorphic to~$\widehat{C}$.
\end{lem}
\begin{proof}
Note that since $\cC_{d,\semis,h}^{L/K}\to{\cX_{K,d}}$ 
  is proper and representable by algebraic spaces, 
  we see that $\widehat{C}$ is a formal algebraic space,
  and  the morphism $\widehat{C}\to\Spf R$ is proper and representable
  by algebraic spaces. In particular, $\widehat{C}$ is a proper
  $\m_R$-adic formal algebraic space over~$\Spf R$.

  To show that we can algebraize~$\widehat{C}$, we will use an
  argument of Kisin, see \cite[Prop.\ 2.1.10]{KisinModularity}, which
  proves a similar statement for moduli of Breuil--Kisin modules. In
  order to use this, we show that we can
  realise~$\cC_{d,\semis,h}^{L/K}$ as a closed substack of a product
  of moduli stacks of Breuil--Kisin modules.

Let~$\pi,\pi'$ be two choices of uniformizers of~$L$, chosen such
that~$\pi/\pi'$ is not a $p$th power, and choose compatible
systems of $p$-power roots~$\piflat,(\pi')^\flat$. Note that the
closure of the subgroup of~$G_L$ generated by~$G_{L_{\piflat,\infty}}$
and~$G_{L_{(\pi')^\flat,\infty}}$ is just~$G_L$, because its fixed
field is contained in $L_{\piflat,\infty}\cap
L_{(\pi')^\flat,\infty}=L$, so that a continuous action of~$G_L$ is
determined by its restrictions to~$G_{L_{\piflat,\infty}}$
and~$G_{L_{(\pi')^\flat,\infty}}$.

Suppose firstly that~$L=K$. We have a natural morphism \numequation\label{eqn: L equals K mapping to two
  uniformizers}\cC_{d,\semis,h}\to\cC_{\piflat,d,h}\times_{\Spf \cO}
\cC_{(\pi')^\flat,d,h},\end{equation} which we claim is a  monomorphism. To see this, it is enough to note
that~$\gMt=\AAinf{A}\otimes_{\gS_{\piflat,A}}\gM_{\piflat}$ is
determined by~$\gM_{\piflat}$, while the $G_K$-action on~$\gMt$ is
determined by the actions of ~$G_{K_{\piflat,\infty}}$
and~$G_{K_{(\pi')^\flat,\infty}}$, which are in turn determined by the
conditions that they act trivially on~$\gM_{\piflat}$
and~$\gM_{(\pi')^\flat}$ respectively.

We may factor \eqref{eqn: L equals K mapping to two uniformizers} as
the composite
\[\cC_{d,\semis,h}\to(\cC_{\piflat,d,h}\times_{\Spf \cO}
  \cC_{(\pi')^\flat,d,h})\times_{(\cR_{\piflat,d}\times_{\Spf \cO}\cR_{(\pi')^\flat,d})}\cX_{{K,d}}\to\cC_{\piflat,d,h}\times_{\Spf \cO}
  \cC_{(\pi')^\flat,d,h}.\] Since the composite is a monomorphism, the
first morphism is a monomorphism.
This first morphism is also proper
(since the morphism $\cC_{d,\semis,h} \to \cX_{K,d}$ obtained by
composing it with the projection to $\cX_{K,d}$ is proper,
by Theorem~\ref{thm: semi stable stack is formal algebraic}, while 
this projection is separated, being a base-change of the product of
the morphisms $\cC_{(\pi)^\flat,d,h}\to\cR_{\piflat,d}$ and
  $\cC_{(\pi')^\flat,d,h}\to\cR_{(\pi')^\flat,d}$, each
of which is proper, and so in
  particular separated, by Theorem~\ref{thm: C is a $p$-adic formal algebraic stack and related
    properties}) 
and so it is in fact a closed immersion. 
Thus~$\widehat{C}$ is a
closed algebraic subspace of
\[(\cC_{\piflat,d,h}\times_{\cR_{\piflat,d}}\Spf
  R)\times_{\Spf  R}(\cC_{(\pi')^\flat,d,h}\times_{\cR_{(\pi')^\flat,d}}\Spf
  R).\] By the Grothendieck Existence theorem for algebraic
spaces~\cite[Thm.\ V.6.3]{MR0302647} (and the symmetry
between~$\piflat$ and~$(\pi')^\flat$), we are therefore reduced to
showing that the morphism
$\cC_{\piflat,d,h}\times_{\cR_{\piflat,d}}\Spf R\to \Spf R$ can be
algebraized to a projective morphism. In the case~$h=1$, this
is~\cite[Prop.\ 2.1.10]{KisinModularity} (bearing in mind the main
result of~\cite{MR1320381}, as in the proof of~\cite[Prop.\
2.1.7]{KisinModularity}),
and the case of general~$h$ can be proved in
exactly the same way, as explained in the proof of~\cite[Prop.\
1.3]{MR2827797}; the key point is that
$\cC_{\piflat,d,h}\times_{\cR_{\piflat,d}}\Spf R$ inherits a natural very
ample (formal) line bundle from the affine Grassmannian.

We now consider the case of a general finite Galois
extension~$L/K$, where we argue as in the proof of Proposition~\ref{lem: restriction from X to Xs}. Let~$\{g_i\}_{i=1,\ldots,n}$ be a set of coset representatives
for~$G_L$ in~$G_K$. 
By definition, to give an object of~$\cC_{d,\semis,h}^{L/K}(A)$ (for some~$A$) is
the same as giving an object of~$\cC_{d,\semis,h}^{L/L}(A)$, together
with an extension of the action of~$G_L$ on the underlying Breuil--Kisin--Fargues
module~$\gMt_A$ to an action of~$G_K$. Now, for each~$i$ we can
give~$g_i^*\gMt_A$ the structure of a Breuil--Kisin--Fargues module by
letting~$h\in G_L$ act as~$g_i^{-1}hg_i$ acts on~$\gMt_A$; it follows
from the definitions that $g_i^*\gMt_A$ is also an object
of~$\cC_{d,\semis,h}^{L/L}(A)$. Then to give an extension of the
action of~$G_L$ on~$\gMt_A$ to an action of~$G_K$ is to give for
each~$i$ an isomorphism $g_i^*\gMt_A\iso\gMt_A$ of objects
of~$\cC_{d,\semis,h}^{L/L}(A)$, satisfying a slew of compatibilities.

We let $\cY_i$ denote the stack classifying objects $\gMt_A$ 
of~$\cC_{d,\semis,h}^{L/L}(A)$, 
endowed with an isomorphism 
$g_i^*\gMt_A\to\gMt_A$.  If we regard $g_i^*$  as an automorphism
of~$\cC_{d,\semis,h}^{L/L}(A)$,
then we  may form  its  graph~$\Gamma_i$,
and we then have an isomorphism of stacks
$$\cY_i \iso
\cC_{d,\semis,h}^{L/L}(A)\times_{\Delta,
\cC_{d,\semis,h}^{L/L}(A) \times
\cC_{d,\semis,h}^{L/L}(A),\Gamma_i}
\cC_{d,\semis,h}^{L/L}(A).$$
(Here $\Delta$ denotes the diagonal;
\emph{cf.}\ the definition of $\cR_d^{\Gamma_{\disc}}$ in 
Section~\ref{subsec: defn of Xd}
above, together 
with Proposition~\ref{prop:Gamma-disc modules as fixed points; etale case}.)
The projection onto the second factor $\cY_i \to 
\cC_{d,\semis,h}^{L/L}(A)$,
which corresponds to
forgetting the isomorphism, is a base-change of the
diagonal~$\Delta,$
and so is affine,
since~$\Delta$ is affine by Theorem~\ref{thm: semi stable stack is formal algebraic}.
In particular,
$\cY_i\times_{\cX_{L,d}}\Spf R$ admits a natural ample formal line
bundle, pulled back from the natural very ample formal line bundle on
$\cC_{d,\semis,h}^{L/L}\times_{\cX_{L,d}}\Spf R$, whose existence we
established above.

We can rephrase our interpretation of objects 
of~$\cC_{d,\semis,h}^{L/K}(A)$ 
as objects of $\cC_{d,\semis,h}^{L/L}(A)$
endowed with isomorphisms
$g_i^*\gMt_A\to\gMt_A$ satisfying compatibilities
as the existence of  a closed immersion
$$
\cC_{d,\semis,h}^{L/K}\into
\cY_1 \times_{\cC_{d,\semis,h}^{L/L}}\times \cdots \times_{\cC_{d,\semis,h}^{L/L}}
\cY_n
$$  
(the compatibilities that cut out 
$\cC_{d,\semis,h}^{L/K}$ arise as
base-changes of the double diagonal, and so impose closed conditions).
%
It follows that $\cC_{d,\semis,h}^{L/K}\times_{\cX_{L,d}}\Spf R$ inherits a natural ample formal line
bundle. As noted in the proof of Lemma~\ref{lem: restriction from X to
  Xs}, since the morphism $\cX_{K,d}\to\cX_{L,d}$ is affine, it  has affine
diagonal, so that the natural morphism
\[\cC_{d,\semis,h}^{L/K}\times_{\cX_{K,d}}\Spf
  R\to\cC_{d,\semis,h}^{L/K}\times_{\cX_{L,d}}\Spf R \]is affine. It
follows that that $\cC_{d,\semis,h}^{L/K}\times_{\cX_{K,d}}\Spf R$ in
turn inherits a natural ample formal line bundle (see e.g.\
\cite[\href{https://stacks.math.columbia.edu/tag/0892}{Tag
  0892}]{stacks-project}), and the  result then follows from another application of \cite[Thm.\
V.6.3]{MR0302647}.
\end{proof}

\section{Inertial types}
\label{subsec:inertial types}
  In this section we examine how to extract the inertial type of the
  Galois representation associated to a Breuil--Kisin--Fargues
  $G_K$-module; in Section~\ref{subsec: HT
  weights} we study the same problem for Hodge--Tate weights, and in
Section~\ref{subsec: moduli stacks of pst GK repns} we use these
results to define our moduli stacks of potentially semistable and
potentially crystalline representations of fixed inertial and Hodge
type.

In contrast to the rest of the book, in this and the following
sections we will need to consider Kisin modules with coefficients
with~$p$ inverted. To this end, we will write ~$\Acirc$ and~$\Bcirc$
for $p$-adically complete flat $\cO$-algebras which are topologically
of finite type over~$\cO$, and
write~$A=\Acirc[1/p]$, $B=\Bcirc[1/p]$. We hope that this change of
notation will not cause any confusion. We will freely use that
if~$\Acirc$ is topologically of finite type over~$\cO$, then
$A:=\Acirc[1/p]$ is Noetherian and Jacobson (\cite[\S0, Prop.\ 9.3.2,
9.3.10]{MR3752648}), and the residue fields of the maximal ideals
of~$A$ are finite extensions of~$K$ (\cite[\S0, Cor.\
9.3.7]{MR3752648}).

Let~$L/K$ be a finite Galois
extension with inertia group~$I_{L/K}$, and 
%
suppose now that~$E$ is large enough that it contains the images of
all embeddings $L\into\Qpbar$, that all irreducible
$E$-representations of~$I_{L/K}$ are absolutely irreducible, and that
every irreducible~$\Qpbar$-representation of~$I_{L/K}$ is defined
over~$E$. Write~$l$ for the residue field of~$L$, and
write~$L_0=W(l)[1/p]$. 

Let~$\gM_{\Acirc}$ be a Breuil--Kisin--Fargues $G_K$-module with~$\Acirc$-coefficients which
  admits all descents over~$L$, and write
  $\overline{\gM}_{\Acirc}:=\gM_{\Acirc,\piflat}/[\piflat]\gM_{\Acirc,\piflat}$ for the
  module considered in Definition~\ref{defn: descending BKF to BK
    coefficients}~\eqref{item: M mod u
    descends} (for some choice of~$\piflat$, with~$\pi$ a uniformizer
  of~$L$; note that by the definition of $\gM_{\Acirc}$ admitting
all descents over~$L$, the quotient $\overline{\gM}_{\Acirc}$
is actually well-defined as a $W(l)\otimes_{\Z_p}A^{\circ}$-submodule
of $W(\overline{k})_A\otimes_{\AAinf{A}} \gM_{\Acirc}$,
independently of the choice $\pi$ or~$\piflat$).
Then~$\overline{\gM}_{\Acirc}$ has a natural $W(l)\otimes_{\Zp}A$-semilinear
  action of~$\Gal(L/K)$, which is defined as follows: if~$g\in \Gal(L/K)$,
  then~$g(\gM_{\Acirc,\piflat})=\gM_{\Acirc,g(\piflat)}$ (see the proof of
  Lemma~\ref{lem: descent for GK action depends only on pi}), so the
  morphism $g:\gM_{\Acirc,\piflat}\to g(\gM_{\Acirc,\piflat})=\gM_{\Acirc,g(\piflat)}$
  induces a morphism \[g: \gM_{\Acirc,\piflat}/[\piflat]\gM_{\Acirc,\piflat} \to
    \gM_{\Acirc,g(\piflat)}/[g(\piflat)]\gM_{\Acirc,g(\piflat)},\]and the source and
  target are both
  canonically identified with~$\overline{\gM}_{\Acirc}$.

  This action of~$\Gal(L/K)$ on~$\overline{\gM}_{\Acirc}$ induces an
  $L_0\otimes_{\Qp}A$-linear action of $I_{L/K}$ on the projective
  $L_0\otimes_{\Qp}A$-module
  $\overline{\gM}_{\Acirc}\otimes_{\Acirc}A$. Fix a choice of
  embedding $\sigma:L_0\into E$, and
  let~$e_\sigma\in L_0\otimes_{\Qp} E$ be the corresponding
  idempotent. Then
  $e_\sigma(\overline{\gM}_{\Acirc}\otimes_{\Acirc}A$) is a projective
  $A$-module of rank~$d$, with an $A$-linear action
  of~$I_{L/K}$. (Indeed, since~$\gM_{\Acirc}$ is a finite
  projective~$\gS_A$-module of rank~$d$, $\overline{\gM}_{\Acirc}$ is
  a finite projective $W(k)\otimes_{\Zp}\Acirc$-module of rank~$d$,
  and $\overline{\gM}_{\Acirc}\otimes_{\Acirc}A$ is a finite
  projective $L_0\otimes_{\Qp}A$-module of rank~$d$. It follows that
  $e_\sigma(\overline{\gM}_{\Acirc}\otimes_{\Acirc}A)$ is a finite
  projective $A$-module of rank~$d$.)

  Note that up to canonical isomorphism, this module does
  not depend on the choice of~$e_\sigma$: indeed the induced action
  of~$\varphi$ commutes with~$I_{L/K}$ and induces isomorphisms
  between the 
  $e_\sigma(\overline{\gM}_{\Acirc}\otimes_{\Acirc}A)$ (with~$\sigma$
  varying), because the cokernel of~$\varphi$ is killed
  by~$E_{\piflat}$, which is a unit in our setting (in
  which~$[\piflat]=0$ and~$p$ is a unit).

  \begin{defn}
    \label{defn: WD representation associated to BKF
      module}Let~$\gMt_{\Acirc}$ be as above. Then we write \index{$\WD(\gMt_{\Acirc})$}
    \[\WD(\gMt_{\Acirc}):=e_\sigma(\overline{\gM}_{\Acirc}\otimes_{\Acirc}A),\]
    a projective~$A$-module of rank~$d$ with an $A$-linear action of $I_{L/K}$.
  \end{defn}

\begin{remark}
\label{rem:inertial types point}
  The point of this definition is that if~$\Acirc$ is a finite flat
  $\cO$-algebra, then it computes inertial types in the following
  sense:  writing~$A:=\Acirc[1/p]$, and~$M_{\Acirc}:=W(\C^\flat)_{\Acirc}\otimes_{\AAinf{\Acirc}}\gMt_{\Acirc}$, we have a potentially semistable representation of~$G_K$ on
  a free $A$-module given by
  $V_A(M_{\Acirc}):=T_{\Acirc}(M_{\Acirc})\otimes_{\Acirc}A$.
As explained in Section~\ref{subsec: Hodge and inertial types}, the
  inertial type $D_{\pst}(V_A(M_{\Acirc}))|_{I_K}$ is then given
  by~$\WD(\gMt_{\Acirc})$ with its action of~$I_{L/K}$.
\end{remark}

  \begin{prop}
  \label{prop: good behaviour of inertial type}Let~$\Acirc$ be a $p$-adically complete flat
$\cO$-algebra which is topologically of finite type over~$\cO$, and write~$A=\Acirc[1/p]$. Let~$\gMt_{\Acirc}$
and~$\gM_{\Acirc}$ be as
above, and fix~$\sigma:L_0\into E$. Then~$\WD(\gMt_{\Acirc})$ 
is a finite
  projective $A$-module of rank~$d$ with an action of~$I_K$, whose formation is
  compatible with base changes $\Acirc\to \Bcirc$ of $p$-adically complete flat
$\cO$-algebras which are topologically of finite type over~$\cO$.
\end{prop}
\begin{proof}Writing~$u$ for~$[\piflat]$, the compatibility of
  formation with base change reduces to the observation
  that the natural map $\gS_A\otimes_AB\to\gS_B$ induces an isomorphism \[\gS_B/u\gS_B \cong (\gS_A/u\gS_A)\otimes_AB. \qedhere\]
\end{proof}
Let~$\tau$ be a $d$-dimensional $E$-representation of~$I_{L/K}$.
\begin{defn}
  \label{defn: having inertial type tau}In the setting of
  Proposition~\ref{prop: good behaviour of inertial type}, we say that
  $\gMt_{\Acirc}$ 
  has inertial type~$\tau$ \index{inertial type}
  if Zariski locally on~$\Spec A$, $\WD(\gMt_{\Acirc})$ is isomorphic to the base change of~$\tau$
  to~$A$.
\end{defn}

\begin{cor}
  \label{cor:inertial type decomposition}In the setting of
  Proposition~{\em \ref{prop: good behaviour of inertial type}}, we  can
  decompose~$\Spec A$ as the disjoint union of open and closed
  subschemes~$\Spec A^\tau$, where~$\Spec A^\tau$ is the locus over
  which $\gMt_{\Acirc}$ has inertial type~$\tau$. Furthermore, the
  formation of this decomposition is compatible with base changes $\Acirc\to \Bcirc$ of $p$-adically complete flat
$\cO$-algebras which are topologically of finite type over~$\cO$.
\end{cor}
\begin{proof}
  By Proposition~\ref{prop: good behaviour of inertial type} (more
  precisely, by the compatibility with base change), we can define
  $\Spec A^\tau$ to be the locus over which $\gMt_{\Acirc}$
  has inertial type~$\tau$. That~$\Spec A$ is actually the
  disjoint union of the $\Spec A^\tau$ follows easily from our
  assumptions on~$E$, and standard facts about the representation
  theory of finite groups in characteristic zero. For lack of a
  convenient reference, we sketch a proof as follows.

  Since~$E$ has characteristic zero, the
  representation~$P:= \oplus_rr$ is a projective generator
  of the category of $E[I_{L/K}]$-modules, where~$r$ runs over a set of
  representatives for the isomorphism classes of irreducible
  $E$-representations of~$I_{L/K}$.  Our assumption that~$E$ is large
  enough that each $r$ is
  absolutely irreducible furthermore shows that
  $\End_{I_{L/K}}(r) = E$ for each~ $r$, so that
  $\End_{I_{L/K}}(P) = \prod_{r} E$.

Standard Morita theory then shows that the functor
$M \mapsto \Hom_{I_{L/K}}(P,M)$ induces an equivalence between the category
of $E[I_{L/K}]$-modules and the category of $\prod_{r} E$-modules.
Of course, a $\prod_{r}E$-module is just given by a tuple
$(N_{r})_{r}$ of $E$-vector spaces,
and in this optic, the functor $\Hom_{I_{L/K}}(P,\text{--})$ 
can be written as $M \mapsto \bigl(\Hom_{I_{L/K}}(r,M)\bigr)_{r}$,
with a quasi-inverse functor being given by $(N_{r}) \mapsto
\bigoplus_{r} r\otimes_{E} N_{r}.$
It is easily seen (just using the fact that $\Hom_{I_{L/K}}(P,\text{--})$
induces an equivalence of categories)
that $M$ is a finitely generated projective $A$-module,
for some $E$-algebra $A$,
if and only if each $\Hom_{I_{L/K}}(r,M)$
is a finitely generated projective $A$-module. Writing the various
~$\tau$ in the form~$\oplus_{r} r^{n_{r}}$, we are done.
\end{proof}
\section{Hodge--Tate weights}\label{subsec: HT
  weights}Let~$\Acirc$ be a $p$-adically complete flat
$\cO$-algebra which is topologically of finite type over~$\cO$, and let 
~$\gMt_{\Acirc}$ be a Breuil--Kisin--Fargues $G_K$-module of
height at most~$h$ with $\Acirc$-coefficients, which admits all
descents. 
We now explain how to interpret the condition that~$\gMt_{\Acirc}$
has a fixed Hodge type, following~\cite[Lem.\ 2.6.1, Cor.\
2.6.2]{MR2373358} (but bearing in mind the corrections to these
results explained in~\cite[A.4]{KisinFM}). We would like to thank Mark
Kisin for a helpful conversation about these results, and for some
suggestions regarding the proof of Proposition~\ref{prop: good behaviour of Kisin
    filtration}.

Fix some choice of~$\piflat$, and write~$\gM_{\Acirc}$
for~$\gM_{\piflat,\Acirc}$, and~$u$ for~$[\piflat]$. For
each~$0\le i\le h$ we define
$\Fil^i\varphi^*\gM_{\Acirc}=\Phi_{\gM_{\Acirc}}^{-1}(E(u)^i\gM_{\Acirc})$,
and we set $\Fil^i\varphi^*\gM_{\Acirc}=\varphi^*\gM_{\Acirc}$
for~$i<0$.

\begin{remark}
\label{rem:Hodge filtration point}
The point of this definition is that it captures the Hodge
filtration. Indeed, if~$\Acirc$ is a finite flat $\cO$-algebra,
then writing~$A=\Acirc[1/p]$, and~$M_{\Acirc}:=W(\C^\flat)_{\Acirc}\otimes_{\AAinf{\Acirc}}\gMt_{\Acirc}$, we have a representation of~$G_K$ on
  a free $A$-module given by
  $V_A(M_{\Acirc}):=T_{\Acirc}(M_{\Acirc})\otimes_{\Acirc}A$. As
explained in Section~\ref{subsec: Hodge and inertial types} (or see
the proof of~\cite[Cor.\ 2.6.2]{MR2373358}), there is a natural
identification of~$\DdR(V_A(M_{\Acirc}))$
with~$(\varphi^*\gM_{\Acirc}/E(u)\varphi^*\gM)\otimes_{\Acirc}A$, under
which $\Fil^i\DdR(V_A(M_{\Acirc}))$ is identified with
\[(\Fil^i\varphi^*\gM_{\Acirc}/E(u)\Fil^{i-1}\varphi^*\gM_{\Acirc})\otimes_{\Acirc}A;\]
in particular, this latter module is a finite projective
$K\otimes_{\Qp}A$-module.
\end{remark}

The following proposition shows that the
final conclusion of the preceding remark
holds for more general choices of~$\Acirc$; the proof
uses that particular case (i.e.\ the case that $\Acirc$ is a finite flat $\cO$-algebra)
as an input.

\begin{prop}
  \label{prop: good behaviour of Kisin
    filtration}Let~$\Acirc$ be a $p$-adically complete flat
$\cO$-algebra which is topologically of finite type over~$\cO$, and write~$A=\Acirc[1/p]$. Let~$\gMt_{\Acirc}$
and~$\gM_{\Acirc}$ be as above.  For each~$0\le i\le h$,
  \[(\Fil^i\varphi^*\gM_{\Acirc}/E(u)\Fil^{i-1}\varphi^*\gM_{\Acirc})\otimes_{\Acirc}A\] is a finite
  projective $K\otimes_{\Qp}A$-module, whose formation is
  compatible with base changes $\Acirc\to \Bcirc$ of $p$-adically complete flat
$\cO$-algebras which are topologically of finite type over~$\cO$.
\end{prop}
\begin{proof}For any morphism of $\cO$-algebras $\Acirc\to \Bcirc$, we write~$\gM_{\Bcirc}:=\gS_{\Bcirc}\otimes_{\gS_{\Acirc}}\gM_{\Acirc}$. We begin by showing (following the proof of \cite[Lem.\
  2.6.1]{MR2373358}) that for any $p$-adically complete $\cO$-algebra
  $\Acirc$ which is topologically of finite type over~$\cO$, both $\gM_{\Acirc}/\im\Phi_{\gM_{\Acirc}}$ and $\varphi^*\gM_{\Acirc}/\Fil^h\varphi^*\gM_{\Acirc}$
  are finite projective $\cO_K\otimes_{\Zp}\Acirc$-modules, whose formation
  is compatible with arbitrary base changes $\Acirc\to \Bcirc$
  (with~$\Bcirc$ also topologically of finite type). 

  Indeed, that~$\gM_{\Acirc}/\im\Phi_{\gM_{\Acirc}}$ is a finite projective
 $\cO_K\otimes_{\Zp}\Acirc$-module whose formation is compatible with base change
  follows as in the proof of~\cite[Lem.\ 2.6.1]{MR2373358} from the
  fact that~$\Phi_{\gM}$ remains injective after any base change
  (because~$E(u)$ remains a non-zero-divisor after any base change),
  and the result for $\varphi^*\gM_{\Acirc}/\Fil^h\varphi^*\gM_{\Acirc}$ follows
  from this and the short exact
  sequence \[0\to\varphi^*\gM_{\Acirc}/\Fil^h\varphi^*\gM_{\Acirc}\to \gM_{\Acirc}/E(u)^h\gM_{\Acirc}\to    \gM_{\Acirc}/\im\Phi_{\gM_{\Acirc}}\to 0.\]

 We claim that for each~$0\le i\le h$, 
  $(\Fil^h\varphi^*\gM_{\Acirc}/E(u)^i\Fil^{h-i}\varphi^*\gM_{\Acirc})\otimes_{\Acirc}A$ is a finite
  projective $K\otimes_{\Qp}A$-module whose formation is compatible with base
  change. Admitting the claim, the proposition follows from the short
  exact sequence
  \nummultline\label{eqn: more Kisin
    fun}
0\longrightarrow\Fil^i\varphi^*\gM_{\Acirc}/E(u)\Fil^{i-1}\varphi^*\gM_{\Acirc}
\\
\stackrel{E(u)^{h-i}}{\longrightarrow}
    \Fil^h\varphi^*\gM_{\Acirc}/E(u)^{h-i+1}\Fil^{i-1}\varphi^*\gM_{\Acirc}\\ \longrightarrow
    \Fil^h\varphi^*\gM_{\Acirc}/E(u)^{h-i}\Fil^{i}\varphi^*\gM_{\Acirc}\longrightarrow 0
\end{multline}(because after tensoring with~$A$, the second and third
terms are projective and compatible with base change, so that the
sequence splits, and the first term is also projective and compatible
with base change).

To prove the claim, we firstly consider for each~$0\le i\le h$ the finite $\cO_K\otimes_{\Zp}\Acirc$-module
  \[\varphi^*\gM_{\Acirc}/(E(u)^i\varphi^*\gM_{\Acirc}+\Fil^h\varphi^*\gM_{\Acirc}).\] Since this is
  the cokernel of the morphism
  \[E(u)^i\varphi^*\gM_{\Acirc}
    \to\varphi^*\gM_{\Acirc}/\Fil^h\varphi^*\gM_{\Acirc}, \]we see that its formation is
  compatible with base change. We have a short exact sequence
of finite type $\Acirc$-modules
  \numequation\label{eqn: first Kisin nightmare}
  \begin{split}
    0\to \Fil^h\varphi^*\gM_{\Acirc}/E(u)^i\Fil^{h-i}\varphi^*\gM_{\Acirc}\to
    \varphi^*\gM_{\Acirc}/E(u)^i\varphi^*\gM_{\Acirc} \\ \to
    \varphi^*\gM_{\Acirc}/(E(u)^i\varphi^*\gM_{\Acirc}+\Fil^h\varphi^*\gM_{\Acirc})\to 0,
  \end{split}
\end{equation}in which the second and third terms are compatible with
base change, and we need to prove that after inverting~$p$, the first
term is projective and is of formation compatible with base change. To
this end, note firstly that as discussed above, if~$\Acirc$ is a
finite flat $\cO$-algebra, then
each~\[(\Fil^i\varphi^*\gM_{\Acirc}/E(u)\Fil^{i-1}\varphi^*\gM_{\Acirc})\otimes_{\Acirc}A\]
is a finite projective $K\otimes_{\Qp}A$-module, and it follows
from \eqref{eqn: more Kisin fun} and an easy induction that the same
is true
of~$(\Fil^h\varphi^*\gM_{\Acirc}/E(u)^i\Fil^{h-i}\varphi^*\gM_{\Acirc})\otimes_{\Acirc}A$. 
  

To simplify notation, we write~\eqref{eqn: first Kisin nightmare} as
\[0\to K(\Acirc)\to M_1\to M_2\to 0,\] and for any $\Acirc$-algebra $C$,
we write $K(C)$ for the kernel of the surjection
$M_1\otimes_{\Acirc}C\to M_2\otimes_{\Acirc}C$. In particular
for~$\Bcirc$ as in the statement of the proposition, if as usual we write 
$B := \Bcirc[1/p]$, then we have
\[K(B)=(\Fil^h\varphi^*\gM_{\Bcirc}/E(u)^i\Fil^{h-i}\varphi^*\gM_{\Bcirc})\otimes_{\Bcirc}B;\]
so we need to show that~$K(A)$ is projective, and
that~$K(B)=K(A)\otimes_{A}
B$. 
  
  Let~$\m$ be a maximal ideal of~$A$, so that $A/\m^i$ is a finite
  $E$-algebra for each~$i$. 
  Then since~$\widehat{A}_\m$ is a flat~$A$-algebra,
since the Artin--Rees Lemma shows that tensoring with
$\widehat{A}_\m$ coincides with  $\m$-adic completion for finite type modules,
and since the formation of kernels commutes with limits,  we see that
  \[K(A)\otimes_{\Acirc}\widehat{A}_\m=K(\widehat{A}_\m)=\varprojlim_iK(A/\m^i).\]
  Since ~$A/\m^i$ is a finite~$E$-algebra, 
  we see that~$K(A/\m^i)$ is a finite projective
  $(K\otimes_{\Qp}A)/\m^i$-module of rank bounded independently
  of~$i$. It follows that $K(A)\otimes_{\Acirc}\widehat{A}_\m$ is a finite
  projective~$\cO_K\otimes_\Zp{\widehat{A}_\m}$-module. Since~$A_\m$
  is Noetherian, 
  $\widehat{A}_\m$ is a
  faithfully flat~$A_\m$-algebra, so that $K(A)\otimes_{\Acirc}{A}_\m$ is a finite
  projective~ $\cO_K\otimes_\Zp A_\m$-module. Since this holds for
  all~$\m$, we see that~$K(A)$ is a finite
  projective~$\cO_K\otimes_\Zp A$-module, as claimed.

  It remains to show that we have $K(B)=K(A)\otimes_{A} B$. We have a
  natural surjective morphism of finite projective
  $K\otimes_{\Qp}B$-modules \numequation\label{eqn: base change of
    K map}K(A)\otimes_A B\to K(B),\end{equation}which we need to show
is an isomorphism. Note firstly that if~$B=A/\m$ for some maximal
ideal~$\m$ of~$A$, then this follows from the previous paragraph. Now
suppose that~$B$ is general. Since the kernel of~\eqref{eqn: base
  change of K map} is in particular a finite projective $B$-module, is
enough to prove that for any maximal ideal~$\m_B$ of~$B$, \eqref{eqn:
  base change of K map} becomes an isomorphism after tensoring
with~$B/\m_B$. 
Suppose that~$\m_B$ lies over a
maximal ideal~$\m$ of~$A$, so that $B/\m_B$ is a finite field
extension of~$A/\m$; 
we need to show that the induced surjection
\[K(A/\m)\otimes_{A/\m} B/\m_B\to K(B/\m_B)\]is an
isomorphism. 
Each side is determined by the de Rham filtration on
the corresponding $G_K$-representation, which is compatible with the
extension of scalars, so we are done.  
\end{proof}
In what follows, we will work in the relative setting of an
extension~$L/K$.  
To this end,
we let~$\Acirc$ be a $p$-adically complete flat $\cO$-algebra, with $A := \Acirc[1/p]$,
let~$L/K$
be a finite Galois extension, and let~$\gMt_{\Acirc}$ be a
Breuil--Kisin--Fargues $G_K$-module of rank $d$ and height at most~$h$ with
$\Acirc$-coefficients, which admits all descents
over~$L$. 
Fix some choice of~$\piflat$ a uniformizer of~$L$, write~$\gM_{\Acirc}$
for~$\gM_{\piflat,\Acirc}$, and~$u$ for~$[\piflat]$.
Applying Proposition~\ref{prop: good behaviour of Kisin filtration},
with $L$ in place of $K$,
we obtain a projective $L\otimes_{\Q_p} A$-module 
	$(\varphi^*\gM_{\Acirc}/E(u)\varphi^*\gM_{\Acirc})
	\otimes_{A^{\circ}} A,$  
	which is filtered by projective submodules.
	This filtered module has a natural action of $\Gal(L/K)$,
	which is semi-linear with respect to the action
	of $\Gal(L/K)$ on $L\otimes_{\Q_p} A$ induced by its action
	on the first factor.
	Since $L/K$ is a Galois extension, the tensor product
	$L\otimes_{\Q_p} A$ is an \'etale $\Gal(L/K)$-extension 
	of $K\otimes_{\Q_p} A$, and so \'etale descent allows
	us to descend 
	$(\varphi^*\gM_{\Acirc}/E(u)\varphi^*\gM_{\Acirc})
	\otimes_{A^{\circ}} A$  
	to a filtered module over $K\otimes_{\Q_p} A$;
	concretely, this descent is achieved by taking $\Gal(L/K)$-invariants.
This leads to the following definition. 

\begin{df}
	\label{def:DdR L/K case}
	In the preceding situation, we write \index{$D_{\dR}(\gMt_{\Acirc})$}
	$$D_{\dR}(\gMt_{\Acirc}) := 
	\bigl((\varphi^*\gM_{\Acirc}/E(u)\varphi^*\gM_{\Acirc})
	\otimes_{A^{\circ}} A\bigr)^{\Gal(L/K)},$$  
		and more generally, for each $i \geq 0,$
		we write
	$$\Fil^i
	D_{\dR}(\gMt_{A^{\circ}}) :=
	\bigl((\Fil^i\varphi^*\gM_{\Acirc}/E(u)\Fil^{i-1}\varphi^*\gM_{\Acirc})\otimes_{\Acirc}A\bigr)^{\Gal(L/K)}$$
(and for $i < 0$, we write 
	$\Fil^i D_{\dR}(\gMt_{A^{\circ}}) :=
	D_{\dR}(\gMt_{A^{\circ}})$).
The property of being a finite rank projective module is preserved
under \'etale descent, and so we find that $D_{\dR}(\gMt_{\Acirc})$
is a rank $d$ projective $K\otimes_{\Q_p} A$-module, filtered by projective
submodules. 

	Since $A$ is an $E$-algebra,
	we have the product decomposition  $K\otimes_{\Q_p} A
	\iso \prod_{\sigma: K \hookrightarrow E} A,$
	and so, if we write $e_{\sigma}$ for the idempotent
	corresponding to the factor labeled by $\sigma$
	in this decomposition,
	we find that 
	$$D_{\dR}(\gMt_{A^{\circ}}) = 
	\prod_{\sigma: K \hookrightarrow E} e_{\sigma} D_{\dR}(\gMt_{A^{\circ}}),$$
	where each
	$e_{\sigma}D_{\dR}(\gMt_{A^{\circ}})$
	is a projective $A$-module of rank~$d$.
	For each $i$, we write
	$$\Fil^i 
	e_{\sigma}D_{\dR}(\gMt_{A^{\circ}})
	=
	e_{\sigma}\Fil^i D_{\dR}(\gMt_{A^{\circ}}).$$
	Each $\Fil^i
	e_{\sigma}D_{\dR}(\gMt_{A^{\circ}})$ is again a projective $A$-module.
\end{df}

The base-change property proved in
	Proposition~\ref{prop: good behaviour of Kisin filtration}
shows
that the various quotients
$\Fil^i
e_{\sigma} D_{\dR}(\gMt_{A^{\circ}})/
\Fil^{i+1}
e_{\sigma} D_{\dR}(\gMt_{A^{\circ}})$
are again projective $A$-modules.  Phrased more geometrically,
then, we see that each $e_{\sigma} D_{\dR}(\gMt_{A^{\circ}})$
gives rise to a vector bundle over $\Spec A$ which is endowed 
with a filtration by subbundles.
We may thus decompose $\Spec A$ into a disjoint union of
open and closed subschemes over which the ranks of the
various subbundles $\Fil^i e_{\sigma} D_{\dR}(\gMt_{\Acirc})$
(or equivalently, the ranks of the various constituents
\[\Fil^i
e_{\sigma} D_{\dR}(\gMt_{A^{\circ}})/
\Fil^{i+1}
e_{\sigma} D_{\dR}(\gMt_{A^{\circ}})\] of the associated graded bundle)
are constant.  

To encode this rank data, and the corresponding decomposition
of~$\Spec A$, it is traditional to use the terminology of
Hodge types, which we now recall.

\begin{df}
	\label{def:Hodge type}\index{Hodge type}
A Hodge type~$\lambdau$ of rank $d$ is by definition a set of tuples of
integers $\{\lambda_{\sigma,j}\}_{\sigma:K\into\Qpbar,1\le j\le d}$
with $\lambda_{\sigma,j}\ge \lambda_{\sigma,j+1}$ for all~$\sigma$ and
all $1\le j\le d-1$. 

If $\underline{D} := (D_{\sigma})_{\sigma:K \hookrightarrow E}$
is a collection
of rank $d$ vector bundles over $\Spec A$, labeled (as indicated) by the embeddings
$\sigma:K \hookrightarrow E$,
then we say that $\underline{D}$ has Hodge type $\lambdau$ if
  $\Fil^i D_{\sigma}$
  has constant rank equal 
  to $\#\{j\mid\lambda_{\sigma,j}\ge i\}.$
\end{df}



\begin{cor}\label{cor:Hodge decomposition L/K case}
	In the preceding context,
	if $\lambdau$ is a Hodge type of rank~$d$,
	then we let 
	$\Spec A^{\lambdau}$ 
	denote the open and closed subscheme of $\Spec A$
	over which the tuple $\bigl(e_{\sigma} D_{\dR}(\gMt_{A^{\circ}})\bigr)$
	of filtered vector bundles is of Hodge type $\lambdau$.
	We have a corresponding decomposition
	\numequation
	\label{eqn:Hodge type decomposition L/K case}
	\Spec A = \coprod_{\lambdau}
	\Spec A^{\lambdau}, 
\end{equation}
labeled by the set of Hodge types $\lambdau$ of rank $d$.  This
decomposition is compatible with base changes $\Acirc\to \Bcirc$ of
$p$-adically complete flat $\cO$-algebras which are topologically of
finite type over~$\cO$.
\end{cor}
\begin{proof}
	This follows from the preceding discussion.
\end{proof}

\begin{remark}
  \label{rem: Hodge filtration not depending on choice of piflat}
   {\em A priori},
the tuple of filtered vector bundles
$\bigl(e_{\sigma} D_{\dR}(\gMt_{\Acirc})\bigr)_{\sigma: K \hookrightarrow E}$
on~$\Spec A$,
and hence the decomposition~\eqref{eqn:Hodge type decomposition L/K case}
of~$\Spec A$,
depends on the descent $\gM_{\piflat, \Acirc}$, 
and hence on the choice of~$\piflat$.
However, 
  Theorem~\ref{thm: admits all descents if and only if semistable}
below
  implies in fact that the decomposition~\eqref{eqn:Hodge type
	  decomposition L/K case} is independent of the choice
  of~$\piflat$; see
  Remark~\ref{rem: Hodge filtration not depending on choice of piflat explained}
  for a more detailed explanation of this.
\end{remark}

As in Section~\ref{subsec:inertial types}, suppose now that~$E$ is large enough that it contains the images of
all embeddings $L\into\Qpbar$, that all irreducible
$E$-representations of~$I_{L/K}$ are absolutely irreducible, and that
every irreducible~$\Qpbar$-representation of~$I_{L/K}$ is defined
over~$E$. We can then immediately combine Corollaries~\ref{cor:inertial type decomposition} and~\ref{cor:Hodge
	  decomposition L/K case}, obtaining the following result.

        \begin{cor}
          \label{cor: Hodge and inertial decomposed together}In the
          preceding situation, 	we have a  decomposition
	\numequation
	\label{eqn:Hodge inertia type decomposition L/K case}
	\Spec A = \coprod_{\lambdau,\tau}
	\Spec A^{\lambdau,\tau}, 
\end{equation}
labeled by the set of Hodge types $\lambdau$ of rank $d$, and the set
of inertial types~$\tau$ of~$I_{L/K}$.  This
decomposition is compatible with base changes $\Acirc\to \Bcirc$ of
$p$-adically complete flat $\cO$-algebras which are topologically of
finite type over~$\cO$. The Breuil--Kisin--Fargues
$G_K$-module~$\gMt_{\Acirc}$
is of Hodge type~$\lambdau$ and inertial type~$\tau$ if and only if $A^{\lambdau,\tau}=A$.
        \end{cor}

The following key theorem relates the constructions of this section, 
and the previous one, to the $p$-adic Hodge theory of
$G_K$-representations. If~$B$ is a finite $E$-algebra, then we say
that a sub-$\cO$-algebra $\Bcirc\subseteq B$ is an \emph{order of~$B$}
if \index{order}
$\Bcirc[1/p]=B$, and~ $\Bcirc$ is a finite $\cO$-algebra.

\begin{thm}
  \label{thm: admits all descents if and only if semistable}
Suppose
  that~$\Acirc$ is a  finite flat $\cO$-algebra,
 let~$M$ be a projective \'etale $(\varphi,G_K)$-module with
  $\Acirc$-coefficients, and write~$V_A(M)=T_{\Acirc}(M)[1/p]$. Let~$L/K$ be a finite Galois extension. 
Then~$V_A(M)|_{G_L}$ 
is
semistable with Hodge--Tate weights in~$[0,h]$ if and only if there
is an order $(\Acirc)'$ of $A := \Acirc[1/p]$ that contains $\Acirc,$
and a Breuil--Kisin--Fargues $G_K$-module $\gMt_{(\Acirc)'}$ with~$(\Acirc)'$-coefficients,
which is of height at most~$h$, which
admits all descents over~$L$, whose $G_L$-action is canonical, and which satisfies
$M_{(\Acirc)'}=W(\C^\flat)_{(\Acirc)'}\otimes_{\AAinf{(\Acirc)'}}\gMt_{(\Acirc)'}$.

Furthermore~$V_A(M)|_{G_L}$ is crystalline if and only if
~$\gMt_{(\Acirc)'}$ is
crystalline as a Breuil--Kisin--Fargues $G_L$-module with~$(\Acirc)'$
coefficients.

In either case, 
the Hodge type is determined by
applying Definition~{\em \ref{def:Hodge type}} to the tuple
$\bigl(e_{\sigma} D_{\dR}(\gMt_{(\Acirc)'})\bigr)_{\sigma:K\hookrightarrow E}$
arising from Definition~{\em \ref{def:DdR L/K case}}, and the inertial type of~$V_A(M)$ is given by
Definition~{\em \ref{defn: having inertial type tau}}.
\end{thm}
\begin{proof}
  Suppose firstly that~$(\Acirc)'$ and~$\gMt_{(\Acirc)'}$ exist. Then
  if we simply forget the $(\Acirc)'$-coefficients and
  consider~$\gMt_{(\Acirc)'}$ as a Breuil--Kisin--Fargues $G_K$-module
  with~$\Zp$-coefficients, the theorem follows from
  Corollaries~\ref{acor: admits all descents if and only if
    potentially semistable} and~\ref{acor: reading off inertial type
    and HT weights from BKF}.  Conversely, suppose
  that~$V_A(M)|_{G_L}$ is semistable with Hodge--Tate weights
  in~$[0,h]$. If we show that~$(\Acirc)'$ and~$\gMt_{(\Acirc)'}$
  exist, then (as discussed
  in Remarks~\ref{rem:inertial types point} and~\ref{rem:Hodge filtration point})
it again follows from Corollary~\ref{acor: reading off inertial
    type and HT weights from BKF} that the inertial and Hodge types
  are given by the claimed recipes; so it suffices to prove this
  existence. We begin with a preliminary reduction.
  
The ring $A$ is a product of Artinian local $E$-algebras,
say $A  = \prod_i A_i$; then if $\Acirc_i$ denotes
the image of $\Acirc$ in~$A_i$, 
the product $\prod_i \Acirc_i$ is an order of $A$ containing~$\Acirc$.
Thus it is no loss of generality to replace $\Acirc$ by this product,
and hence, by working one factor at a time, to assume that $A$ is local.
  The residue field~$E'$ of~$A$ is then a finite extension of~$E$, and $A$ is naturally
  an~$E'$-algebra. The compositum $\cO_{E'}\Acirc$ is an order
in $A$ containing~$\Acirc$ which is furthermore an $\cO_{E'}$-algebra.
Replacing $\Acirc$ by this compositum, we may thus assume
that $\Acirc$ is an $\cO_{E'}$-algebra.
Then, relabelling~$E'$ as~$E$ if necessary, we can and do
  assume that~$E'=E$ . 
Thus, reformulating the preceding discussion slightly,
we have $A_{\red} = E$ and $\Acirc_{\red} = \cO_E$.


  By Corollary~\ref{acor: admits all
    descents if and only if potentially semistable} there is a unique
  Breuil--Kisin--Fargues $G_K$-module~$\gMt$ with $\Zp$-coefficients
  which is of height at most~$h$ and admits all descents over~$L$, and
  satisfies \[M=W(\C^\flat)\otimes_{\Ainf}\gMt.\] For each~$\piflat$
  we have by definition a descent~$\gM_{\piflat}$ of~$\gMt$. Fix some
  choice of~$\piflat$, write~$u$ for~$[\piflat]$, and write~$\gS$
  for~$\gS_{\piflat}$ and~$\gM$ for~$\gM_{\piflat}$.

  We now follow the proof of~\cite[Prop.\ 1.6.4]{MR2373358}. In
  particular, we work for the most part with the Breuil--Kisin
  module~$\gM$, rather than the Breuil--Kisin--Fargues
  module~$\gMt$. It is possible that by using~\cite[Prop.\
  4.13]{2016arXiv160203148B}, we could make our arguments with~$\gMt$
  itself, but since this does not seem likely to significantly
  simplify the proof, and would make it harder for the reader to
  compare to the arguments of~\cite{MR2373358}, we have not attempted
  to do this.

 Note
  firstly that~$\gM\subset\gMt\subset M$ are stable under the action
  of~$\Acirc$ on~$M$. Indeed, since~$\Acirc$ is local, it is enough to
  check stability under the action of~$(\Acirc)^\times$. In the case
  of~$\gMt$ it is immediate from the unicity of~$\gMt$ that we
  have~$\gMt=a\gMt$ for any~$a\in(\Acirc)^\times$, and similarly
  for~$\gM$ it follows from Lemma~\ref{lem: uniqueness of descent for
    G K infty Ainf}.

  In particular~$\gM$ is naturally a
  finite~$\gS_{\Acirc}$-module, which is projective
  as an $\gS$-module.  While~$\gM$ need not be a projective
  $\gS_{\Acirc}$-module, $\cO_{\cE,\Acirc}\otimes_{\gS_{\Acirc}}\gM$ is a projective~$\cO_{\cE,\Acirc}$-module,
  because after the faithfully flat base extension
  $\cO_{\cE,\Acirc}\into W(\C^\flat)_{\Acirc}$ it is identified with~$M$.
It follows
from~\cite[Lem.\ 1.6.1]{MR2373358}
that $\gM[1/p]$ is a projective
  $\gS_{\Acirc}[1/p]$-module, necessarily of rank~$d$, 
and that $\gM[1/u]$ is a projective $\gS_{\Acirc}[1/u]$-module,
again of rank~$d$.

  In fact, it will be useful to note that $\gM[1/p]$ is actually
  free of rank~$d$.  To see this, it suffices to prove it
  after base-changing to $$(\gS_{\Acirc})[1/p]_{\red}
  = (\gS_{\Acirc_{\red}})[1/p] = \gS_{\cO}[1/p] 
  \cong \prod_{\sigma:W(k) \hookrightarrow \cO} \cO[[u]][1/p].$$
  Since $\cO[[u]][1/p]$ is a PID (by the Weierstrass preparation
  theorem), 
  any finitely generated projective module over
  $(\gS_{\Acirc_{\red}})[1/p]$ of constant rank is necessarily free.
For later use,  we note that since
$(\cO_{\cE,\Acirc})_{\red}=
\cO\otimes_{\Z_p} \cO_{\cE} = \prod_{\sigma: W(k) \hookrightarrow \cO}
\cO_{\cE}$,
and $\cO_{\cE}$ is a local ring (indeed a discrete valuation ring,
with uniformizer~$p$), 
any projective module of constant rank over
$\cO_{\cE,\Acirc}$ is also necessarily free.

  Now let~$\gM'_{\cO}$ denote the image of~$\gM$ under the
  projection~$\gM[1/u]\to\gM[1/u]\otimes_AE$ (the map~$A\to E$ being
  the projection from~$A$ to its residue field), and
  let~$\gM'_{\cO}\subset\gM_{\cO}$ be the canonical inclusion (with
  finite cokernel) of~$\gM'_{\cO}$ into the corresponding finite
  projective~$\gS_{\cO}$-module. (The existence of~$\gM_{\cO}$
  follows from the structure theory of $\gS$-modules, see~\cite[Prop.\
  4.3]{2016arXiv160203148B}. Concretely, we have
  $\gM_{\cO}=\gM'_{\cO}[1/u]\cap\gM'_{\cO}[1/p]$.) The module~$\gM_{\cO}$ is again a
  Breuil--Kisin module of height at most~$h$.

  Choose an $\gS_{\cO}$-basis for~$\gM_{\cO}$, and lift it to an
  $\gS_{\Acirc}[1/p]$-basis of $\gM[1/p]$,
  and let~$(\Acirc)''$ be the
  $\Acirc$-subalgebra of~$A$ generated by these matrix coefficients.
  These coefficients have bounded powers of $p$  in their
  denominators (when we think of $A$ as equalling $\Acirc[1/p]$)
  and lie in $\cO$ after passing to the quotient of $A$ by its
  maximal ideal (equivalently its nilradical); in other words
  these entries lie in $\Acirc + p^{-N} \nil(\Acirc)$ for some $N \geq 0$,
  and this  $\Acirc$-submodule of $A$  generates a finite $\Acirc$-subalgebra
  of~$A$.
Thus $(\Acirc)''$  is an  order in~$A$. 
Let~$\gM_{(\Acirc)''}$ be the $\gS_{(\Acirc)''}$-submodule
  of~$\gM[1/p]$ spanned by this basis. Then~$\gM_{(\Acirc)''}$ is
  $\varphi$-stable by construction.
  Again by construction we have that
 $$\gM_{(\Acirc)''}[1/p] 
  = (\Acirc)''\otimes_{\Acirc} \gM[1/p],$$  
  and that
$$\cO\otimes_{(\Acirc)''} \cO_{\cE,(\Acirc)''} \otimes_{\gS_{(\Acirc)''}}
\gM_{(\Acirc)''}
=  \cO_{\cE}\otimes_{\gS} \gM_{\cO}
=  \cO_{\cE}\otimes_{\gS} \gM'_{\cO}$$
(where the first tensor product is with respect to the surjection
$(\Acirc)'' \to  (\Acirc)''_{\red} = \cO$).
Thus $\cO_{\cE,(\Acirc)''} \otimes_{\gS_{(\Acirc)''}} \gM_{(\Acirc)''}$ and 
  $\cO_{\cE,(\Acirc)''}\otimes_{\gS_{\Acirc}} \gM$
  are two rank $d$ projective $\cO_{\cE,(\Acirc)''}$-modules which
  coincide after inverting $p$, and also after reducing modulo
 $\nil\bigl( (\Acirc)'').$ In fact, they are both free of rank~$d$, by
our above observation that a projective $\cO_{\cE,\Acirc}$-module of
constant rank is necessarily free
(applied now  with $(\Acirc)''$  in  place of $\Acirc$).
If we choose bases for each of them which coincide 
modulo $\nil\bigl((\Acirc)''\bigr)$, then we may regard
each of these bases as a basis of $\cO_{\cE,(\Acirc)''}\otimes_{\gS_{\Acirc}}
\gM[1/p]$
over $\cO_{\cE,(\Acirc)''}[1/p]$, and so they differ by a
change-of-basis matrix lying in $1 + \nil\bigl((\Acirc)'') 
M_d\bigl(\cO_{\cE,(\Acirc)''}[1/p]\bigr).$
Just as we argued above in the construction of~$(\Acirc)''$,
we see that we may enlarge $(\Acirc)''$ to an
order $(\Acirc)'$ such that the entries of this change-of-basis matrix
lie in~$(\Acirc)'$.  Consequently, if we set
$\gM_{(\Acirc)'}   := \gS_{(\Acirc)')} \otimes_{\gS_{(\Acirc)''}} 
\gM_{(\Acirc)''}$, then we find that
$$\cO_{\cE,(\Acirc)'} \otimes_{\gS_{\Acirc}} \gM
= \cO_{\cE,(\Acirc)'} \otimes_{\gS_{(\Acirc)''}} \gM_{(\Acirc)'}.$$
Extending scalars further, we find that
\begin{multline*}
M_{(\Acirc)'} :=  W(\C^\flat)_{(\Acirc)'}\otimes_{W(\C^\flat)_{\Acirc}}
M \\
 = W(\C^\flat)_{(\Acirc)'} \otimes_{\gS_{\Acirc}} \gM
= W(\C^\flat)_{(\Acirc)''}\otimes_{\gS_{(\Acirc)'}} \gM_{(\Acirc)'}.
\end{multline*}


  We claim that $\gM_{(\Acirc)'}$ is a Breuil--Kisin module of height at most~$h$,
  i.e.\ that the cokernel of~$\Phi_{\gM_{(\Acirc)'}}$ is killed
  by~$E(u)^h$. This cokernel is, as an $\gS$-module, a successive
  extension of copies of the cokernel of~$\Phi_{\gM_{\cO}}$ (the
  injectivity of~$\Phi$ implies that the formation of the cokernel
  of~$\Phi$ is exact), 
  and is in
  particular $p$-torsion free (by \cite[Lem.\ 1.2.2~(2)]{KisinModularity}). 
  After inverting~$p$ the cokernel is a
  base change of the cokernel of~$\Phi_{\gM}$, and is therefore killed
  by~$E(u)^h$, as required.
  
  We now
  set~$\gMt_{(\Acirc)'}:=\AAinf{(\Acirc)'}\otimes_{\gS_{(\Acirc)'}}\gM_{(\Acirc)'}$,
  a Breuil--Kisin--Fargues $G_K$-module of height at most~$h$ with
  $(\Acirc)'$-coefficients which by construction satisfies
  $M_{(\Acirc)'}=W(\C^\flat)_{(\Acirc)'}\otimes_{\AAinf{(\Acirc)'}}\gMt_{(\Acirc)'}$. It
  remains to prove that $\gMt_{(\Acirc)'}$ admits all descents
  over~$L$ and that the $G_L$-action is canonical. Applying Corollary~\ref{acor: admits all descents if and
    only if potentially semistable} again to~$M_{(\Acirc)'}$ regarded
  as an \'etale $(\varphi,G_K)$-module with $\Zp$-coefficients, we see that
  there is a Breuil--Kisin--Fargues $G_K$-module~$(\gMt)'$ of height
  at most~$h$ with $\Zp$-coefficients which admits all descents
  over~$L$ and satisfies
  $M_{(\Acirc)'}=W(\C^\flat)\otimes_{\Ainf}(\gMt)'$. Write~$(\gM)'$
  for the descent of~$(\gMt)'$, for our particular choice of~$\piflat$;
  then by the uniqueness of Breuil--Kisin modules
  with~$\Zp$-coefficients (\cite[Prop.\ 2.1.12]{KisinCrys}) we have
  $(\gM)'=\gM_{(\Acirc)'}$, so that in
  fact~$(\gMt)'=\gMt_{(\Acirc)'}$. Thus~$\gMt_{(\Acirc)'}$ admits all
  descents over~$L$ to Breuil--Kisin modules with~$\Zp$-coefficients,
  and by another application of Lemma~\ref{lem: uniqueness of descent
    for G K infty Ainf}, these Breuil--Kisin modules
  are~$(\Acirc)'$-stable. Finally, the $G_L$-action is canonical by Proposition~\ref{prop: Caruso Liu canonical action semistable}, and we are done.
\end{proof}

\begin{rem}
  \label{rem: Hodge filtration not depending on choice of piflat explained}It
  follows from
  Theorem~\ref{thm: admits all descents if and only if semistable}
  that the decomposition 
  ~\eqref{eqn:Hodge type decomposition L/K case} is
  independent of the choice of~$\piflat$. Indeed,
  the condition that
  $\Fil^i e_{\sigma} D_{\dR}(\gMt_{\Acirc})$
  has given constant rank can be checked after base changing
  to 
  all~$A/\m$, for~$\m$ a maximal ideal of~$A$, so it follows from
  Proposition~\ref{prop: good behaviour of Kisin filtration} and
  Theorem~\ref{thm: admits all descents if and only if semistable}
  that the direct factor $A^{\lambdau}$ of $A$ is characterized as follows:
  for any $\Acirc$-algebra $\Bcirc$ which
  is finite and flat as an~$\cO$-algebra,
  the canonical morphism $A \to \Bcirc[1/p]$ factors through $A^{\lambdau}$
  if and only if the $G_K$-representation
  $V_B(\gMt_{\Bcirc})$ has Hodge type~$\lambdau$.
\end{rem}

\section[{Moduli stacks of potentially semistable
  representations}]{Moduli stacks of potentially semistable
  representations \sectionmark{Stacks of potentially semistable
  representations}}\sectionmark{Stacks of potentially semistable
  representations}\label{subsec: moduli stacks of pst GK repns}
We are now in
a position to define the main objects of interest in this chapter,
which are the closed substacks of~$\cX_d$ classifying representations
which are potentially crystalline or potentially semistable of
fixed Hodge and inertial types. 

In this section ~$L/K$ will denote a finite Galois
extension with inertia group~$I_{L/K}$; we will always assume
(without specifically recalling this) that~$E$ is large enough that all irreducible
$E$-representations of~$I_{L/K}$ are absolutely irreducible, and that
every irreducible~$\Qpbar$-representation of~$I_{L/K}$ is defined
over~$E$. 
For any such~$L/K$ we have the
stacks   $\cC_{d,\ss,h}^{L/K}$ and $\cC_{d,\crys,h}^{L/K}$
that we defined in Definition~\ref{defn: semistable crys stacks L over K}, 
and we write $\cC_{d,\ss,h}^{L/K,\fl}$ and
$\cC_{d,\crys,h}^{L/K,\fl}$ respectively for their flat parts (in
the sense recalled in Appendix~\ref{app: formal algebraic
  stacks}).

\begin{df} \index{Hodge type!effective}
We say that a Hodge
type~$\underline{\lambda}$ is \emph{effective} 
if $\lambda_{\sigma,i}
\geq 0$ for each $\sigma$ and $i$,
and that $\lambdau$ is \emph{bounded by~$h$} \index{Hodge type!bounded
by~$h$} if
$\lambda_{\sigma,i}\in [0,h]$ for each~$\sigma$ and $i$.
\end{df}

\begin{prop}
  \label{prop: HT weights are a closed condition}
  Let~$L/K$ be a finite
  Galois extension. Then the stacks
  $\cC_{d,\ss,h}^{L/K,\fl}$ and $\cC_{d,\crys,h}^{L/K,\fl}$
  are scheme-theoretic unions of closed substacks
  $\cC_{d,\ss,h}^{L/K,\fl,\lambdau,\tau}$ and
  $\cC_{d,\crys,h}^{L/K,\fl,\lambdau,\tau}$, where~$\lambdau$ runs
  over all effective Hodge types that are bounded by~$h$, and~$\tau$ runs over
  all $d$-dimensional $E$-representations of~$I_{L/K}$. These latter closed substacks are
  uniquely characterised by the following property: if~$\Acirc$ is a
  finite flat~$\cO$-algebra, then an $\Acirc$-point of
  $\cC_{d,\ss,h}^{L/K,\fl}$ \emph{(}resp.\
  $\cC_{d,\crys,h}^{L/K,\fl}$\emph{)} is a point of
  $\cC_{d,\ss,h}^{L/K,\fl,\lambdau,\tau}$ \emph{(}resp.\
  $\cC_{d,\crys,h}^{L/K,\fl,\lambdau,\tau}$\emph{)} if and only if
  the corresponding Breuil--Kisin--Fargues module~$\gMt_{\Acirc}$ has
  Hodge type~$\lambdau$ in the sense given by
applying Definition~{\em \ref{def:Hodge type}} to the tuple
$\bigl(e_{\sigma} D_{\dR}(\gMt_{\Acirc})\bigr)_{\sigma:K\hookrightarrow E}$
arising from Definition~{\em \ref{def:DdR L/K case}},
  and inertial type~$\tau$ in the
sense of Definition~{\em \ref{defn: having inertial type tau}};
or, more succinctly {\em (}and writing $A := \Acirc[1/p]$, as usual{\em )},
if and only if $A^{\lambdau,\tau} = A$. 

Finally, each stack $\cC_{d,\ss,h}^{L/K,\fl,\lambdau,\tau}$ or
$\cC_{d,\crys,h}^{L/K,\fl,\lambdau,\tau}$ is a
$p$-adic
formal algebraic stack of finite presentation which is flat over~$\Spf\cO$
and whose diagonal is affine, and the natural morphisms to $\cX_{d}$ are
representable by algebraic spaces, proper, and of finite
presentation. 
\end{prop}
\begin{proof}
  We give the proof
  for~$\cC_{d,\ss,h}^{L/K,\fl}$, the crystalline case being
  formally identical.
  Since~$\cC_{d,\ss,h}^{L/K,\fl}$ is a flat
  $p$-adic formal algebraic stack of finite presentation
  over~$\Spf\cO$, we can choose a smooth surjection
  \numequation
  \label{eqn:B to C}
  \Spf\Bcirc\to
  \cC_{d,\ss,h}^{L/K,\fl}
  \end{equation}
  where~$\Bcirc$ is topologically of
  finite type over~$\cO$.

As usual, we write $B := \Bcirc[1/p]$.
By Corollary~\ref{cor: Hodge and inertial decomposed together} 
  we may write $\Spec B$ as a disjoint union $\Spec B = \coprod_{\lambdau,
	  \tau} B^{\lambdau,\tau},$
  and so correspondingly factor $B$ as a product
  $B = \prod_{\lambdau,\tau} B^{\lambdau,\tau}.$
  If we let $B^{\lambdau,\tau,\circ}$ denote the image of $\Bcirc$
  in $B^{\lambdau,\tau}$, then we obtain an induced injection
  $B^{\circ} \hookrightarrow \prod_{\lambdau,\tau} B^{\lambdau,\tau,\circ}$,
  which induces a scheme-theoretically dominant morphism
  \numequation
  \label{eqn:dominant B}
  \coprod_{\lambdau,\tau} \Spf B^{\lambdau,\tau,\circ} \to \Spf \Bcirc.
  \end{equation}

  Write $R = \Spf \Bcirc \times_{\cC_{d,\ss,h}^{L/K,\fl}} \Spf \Bcirc,$
	  so that 
	  $[\Spf \Bcirc/R] \iso \cC_{d,\ss,h}^{L/K,\fl}.$
	  Then we claim that each $\Spf B^{\lambdau,\tau,\circ}$ is $R$-invariant.
	  Granting this,
	  if we write $R^{\lambdau,\tau}$ for the
	  restriction of $R$ to $\Spf B^{\lambdau,\tau,\circ}$,
	  and then define
	  $\cC_{d,\ss,h}^{L/K,\fl,\lambdau,\tau}
	  := [\Spf B^{\lambdau,\tau,\circ}/R^{\lambdau,\tau}],$
	  it follows from the discussion of closed substacks
	  in Appendix~\ref{app: formal
		 algebraic stacks} 
               that 
	  $\cC_{d,\ss,h}^{L/K,\fl,\lambdau,\tau}$
	  embeds as a closed sub-formal algebraic stack of 
	  $\cC_{d,\ss,h}^{L/K,\fl}.$
	  Furthermore, the induced morphism
	  \numequation
	  \label{eqn:dominant C}
	  \coprod_{\lambdau,\tau} 
	  \cC_{d,\ss,h}^{L/K,\fl,\lambdau,\tau}
	  \to 
	  \cC_{d,\ss,h}^{L/K,\fl}
  \end{equation}
	  induces the morphism~\eqref{eqn:dominant B} after
	  pull-back via the morphism~\eqref{eqn:B to C},
	  which is representable by algebraic spaces,
	  smooth, and surjective.  Since this latter morphism is
	  scheme-theoretically dominant, so is~\eqref{eqn:dominant C}.

	  We now verify that $\Spf B^{\lambdau,\tau,\circ}$ is 
	  $R$-invariant, i.e.\ that
	  $R_0 := \Spf B^{\lambdau,\tau,\circ} \times_{\cC_{d,\ss,h}^{L/K,\fl}} \Spf \Bcirc$
	  and
	  $R_1 := \Spf \Bcirc \times_{\cC_{d,\ss,h}^{L/K,\fl}} \Spf B^{\lambdau,\tau,\circ}$
	  coincide as closed sub-formal algebraic spaces of 
	  $R = \Spf \Bcirc \times_{\cC_{d,\ss,h}^{L/K,\fl}} \Spf \Bcirc$
	  (and thus that both coincide with
	  $R^{\lambdau,\tau}
	  := \Spf B^{\lambdau,\tau,\circ} \times_{\cC_{d,\ss,h}^{L/K,\fl}}
	  \Spf B^{\lambdau,\tau,\circ}$).

	  Since $\Spf \Bcirc \to \cC_{d,\ss,h}^{L/K,\fl}$
	  is representable by algebraic spaces and smooth,
	  it is in particular representable by algebraic spaces,
	  flat, and locally of finite type, and thus so are each 
	  of the projections
	  $R_0 \to \Spf B^{\lambdau,\tau,\circ}$
	  and 
	  $R_1 \to \Spf B^{\lambdau,\tau,\circ}.$
	  It follows from Lemmas~\ref{lem: alem closed formal algebraic scheme
    adic*} and~\ref{alem: Noetherian flat affine formal algebraic
    spaces} that any affine formal algebraic space which is \'etale
	  over either of
	  $R_0$ or $R_1$
	  is of the form $\Spf \Acirc$, where $\Acirc$ is $\varpi$-adically
	  complete, topologically
	  of finite type, and flat over~$\cO$. 
	Thus, to show that $R_0$ and $R_1$ coincide as subspaces
	of $R$, it suffices to show that if $\Acirc$ is any
	$\varpi$-adically complete, topologically of finite type,
	and flat $\cO$-algebra endowed with a morphism 
	$\Spf \Acirc \to R$, then this morphism factors through
	$R_0$ if and only if it factors through $R_1$.
	Unwinding the definitions of $R,$ $R_0$, and $R_1$
	as fibre products, this amounts to showing that
	if we are given a pair of morphism
	$\Spf \Acirc \rightrightarrows \Spf \Bcirc$
	which induce the same morphism to $\cC_{d,\ss,h}^{L/K,\fl}$,
	then one factors through $\Spf B^{\lambdau,\tau,\circ}$
	if and only if the other does.

	The pair of morphisms
	$\Spf \Acirc \rightrightarrows \Spf \Bcirc$
	correspond to a pair of morphisms
	$f_0,f_1: \Bcirc \rightrightarrows \Acirc$.
	The fact that these morphisms induce the same morphism to
	$\cC_{d,\ss,h}^{L/K,\fl}$ may be rephrased as saying 
	that $f_0^*\gMt$ and $f_1^*\gMt$ are isomorphic
	Breuil--Kisin--Fargues modules over $\Acirc$.  
        Now, applying the base change statements
	from Propositions~\ref{prop: good behaviour of inertial type}
	and~\ref{prop: good behaviour of Kisin filtration},
	we find that both
	$f_0^*\WD(\gMt_{\Bcirc})$
	and
	$f_1^*\WD(\gMt_{\Bcirc})$
	are isomorphic, as $I_K$-representations,
	to $\WD(\gMt_{\Acirc})$,
        and that both
	$f_0^*D_{\dR}(\gMt_{\Bcirc})$
	and
	$f_1^*D_{\dR}(\gMt_{\Bcirc})$
	are isomorphic, as filtered $K\otimes_{\Q_p} B$-modules,
	to $D_{\dR}(\gMt_{\Acirc})$.  Thus either, and hence both,
	of
	$f_0^*\WD(\gMt_{\Bcirc})$
	and
	$f_1^*\WD(\gMt_{\Bcirc})$
	are of inertial type $\tau$
	if and only if
	$\WD(\gMt_{\Acirc})$ is of inertial type~$\tau$,
	while either, and hence both,
	of 
	$f_0^*D_{\dR}(\gMt_{\Bcirc})$
	and
	$f_1^*D_{\dR}(\gMt_{\Bcirc})$
	are of Hodge type $\lambdau$ if and only if 
	$D_{\dR}(\gMt_{\Acirc})$ is of Hodge type~$\lambdau$.  
	Consequently $f_0$ factors through
	$\Spf B^{\lambdau, \tau,\circ}$
	if and only $f_1$ does.  This shows that
	$\Spf B^{\lambdau, \tau,\circ}$ is indeed $R$-invariant.
        
	
We now show that~ $\cC_{d,\ss,h}^{L/K,\fl,\lambdau,\tau}:= [\Spf
B^{\lambdau,\tau,\circ}/R^{\lambdau,\tau}]$ satisfies the required
property, i.e.\ that for any finite flat $\cO$-algebra~$\Acirc$,
a morphism $\Spf A^{\circ} \to \cC_{d,\ss,h}^{L/K,\fl}$
factors through $\cC_{d,\ss,h}^{L/K,\fl,\lambdau,\tau}$
if and only if 
$A^{\lambdau,\tau}  = A$.
By construction, such a factorization occurs if and only if the induced
morphism
$\Spf \Acirc \times_{\cC_{d,\ss,h}^{L/K,\fl}} \Spf B^{\circ}
	\to \Spf B^{\circ}$ 
	factors through
	$\Spf B^{\lambdau,\tau,\circ}$.

	Since the surjection $\Spf B^{\circ} \to
	\cC_{d,\ss,h}^{L/K,\fl}$ is representable by algebraic
	spaces, smooth, and surjective,
	the fibre product 
$\Spf \Acirc \times_{\cC_{d,\ss,h}^{L/K,\fl}} \Spf B^{\circ}$
is a formal algebraic space,
	and the surjection
$\Spf \Acirc \times_{\cC_{d,\ss,h}^{L/K,\fl}} \Spf B^{\circ}
	\to \Spf \Acirc$
	is representable by algebraic spaces, smooth, and surjective.
	We may thus find an open cover of 
$\Spf \Acirc \times_{\cC_{d,\ss,h}^{L/K,\fl}} \Spf B^{\circ}$
by affine formal algebraic spaces  $\Spf C^{\circ}$,
with $C^{\circ}$ being a faithfully flat $\varpi$-adically complete
$\Acirc$-algebra which is topologically of
finite type. 
	  (We are again applying Lemmas~\ref{lem: alem closed formal algebraic scheme
    adic*} and~\ref{alem: Noetherian flat affine formal algebraic
    spaces}.)
	If we write $C := C^{\circ}[1/p]$,
then the base-change property of Corollary~\ref{cor: Hodge and
  inertial decomposed together} 
        shows that
	the hypothesized factorization occurs if and only
	if $C = C^{\lambdau,\tau}$, 
	and also that $C^{\lambdau,\tau} = A^{\lambdau,\tau}\otimes_A C,$
	so that (since $C$ is faithfully flat over $A$) 
	$C^{\lambdau,\tau} = C$ if and only if $A^{\lambdau,\tau} = A$.
	In conclusion, we have shown that 
$\Spf A^{\circ} \to \cC_{d,\ss,h}^{L/K,\fl}$
factors through $\cC_{d,\ss,h}^{L/K,\fl,\lambdau,\tau}$
if and only if $A^{\lambdau,\tau} = A$, as required.

Finally, the uniqueness of~$\cC_{d,\ss,h}^{L/K,\fl,\lambdau,\tau}$ is
immediate from Proposition~\ref{prop: uniqueness of a flat stack with
  local ring points} below; the properties of being of finite
presentation,
of having affine diagonal,
and
of the morphism $\cC_{d,\ss,h}^{L/K,\fl,\lambdau,\tau}\to\cX_d$ being representable by algebraic spaces, proper, and of finite
presentation, are immediate from Proposition~\ref{prop: pst and pcrys
  stacks are padic} and Lemma~\ref{lem:finite presentation from closed immersion}.
\end{proof}

\begin{prop}\label{prop: uniqueness of a flat stack with local ring points}
	Suppose that $\cY$ and $\cZ$ are two Noetherian
	formal algebraic stacks, both lying over $\Spf \cO$ and both flat
	over $\cO$, and both embedded as closed substacks
	of a stack $\cX$ over $\Spec \cO$ 
        {\em (}in the usual sense
       that the inclusions of substacks $\cY, \cZ \hookrightarrow \cX$
       are representable by algebraic spaces, and induce closed immersions
when pulled back over any scheme-valued point of $\cX${\em )}.      
Suppose also that each of $\cY_{\red}$ and $\cZ_{\red}$ are locally
of finite type over $\F$,
	and suppose further that, for any finite flat $\cO$-algebra~$\Acirc$,
	a morphism $\Spf \Acirc \to \cX$ factors through ~$\cY$ if
	and only if it factors through~ $\cZ$.  
	Then $\cY$ and $\cZ$ coincide.
\end{prop}
\begin{proof}
	The $\cO$-flat part $\cW := (\cY\times_{\cX} \cZ)_{\fl}$
        of the $2$-fibre product~$\cY\times_{\cX}\cZ$ embeds
	as a closed substack of each of $\cY$ and~$\cZ$,
	and it suffices to show that each of those embeddings
        are isomorphisms.
	The stack $\cW$
	inherits all of the hypotheses shared by $\cY$ and~$\cZ$,
	including the fact that, for any finite flat $\cO$-algebra~$\Acirc$,
	the $\Acirc$-valued points of $\cW$ coincide with 
	the $\Acirc$-valued points of each of $\cY$ and~$\cZ$.
        Thus, replacing the pair $\cY$, $\cZ$ by the pairs
        $\cW, \cY$ and $\cW, \cZ$ in turn, we may reduce to the case
  when~$\cY$ is actually a closed substack of~$\cZ$.

  Since~$\cZ$ is Noetherian, lies over $\Spf \cO$, and is flat over $\cO$,
  we can find a smooth
  surjection $U\to\cZ$, whose source is a disjoint union of affine
  formal algebraic spaces $\Spf\Ccirc$, where~$\Ccirc$ is a flat
  $I$-adically complete Noetherian $\cO$-algebra, for some ideal
  $I$ which contains $\varpi$ and defines the topology on $\Ccirc$.

  The assumption on $\cZ_{\red}$ further implies that $\Ccirc/I$ 
  is a finite type $\F$-algebra.
  Then~$\Spf\Ccirc\times_{\cZ}\cY$ is a closed sub-formal
  algebraic space of~$\Spf\Ccirc$, and (since it is also~$\cO$-flat)
  it follows from Lemmas~\ref{lem: alem closed formal algebraic scheme
    adic*} and~\ref{alem: Noetherian flat affine formal algebraic
    spaces} that it is of the form~$\Spf\Bcirc$, where~$\Bcirc$ is an $\cO$-flat
  quotient of~$\Ccirc$, 
  again endowed with the $I$-adic topology.
  Since the property of being an isomorphism can be checked {\em fppf}
  locally on the target, we may replace the closed 
  embedding $\cY \hookrightarrow \cZ$ 
  by $\Spf \Bcirc \hookrightarrow \Spf \Ccirc.$
  In summary, we have a surjection $\Ccirc \to \Bcirc$ of flat,
  $I$-adically complete 
  $\cO$-algebras, with $\Ccirc/I$ (and hence also $\Bcirc/I$) being 
  of finite type over~$\F$,
  and with the further property that any morphism $\Ccirc \to \Acirc$
  of $\cO$-algebras in which $\Acirc$ is finite flat over $\cO$ factors
  through $\Bcirc$; we then have to prove that this surjection is an
  isomorphism.

  Write~$B=\Bcirc[1/p]$, $C=\Ccirc[1/p]$;
  since~$\Bcirc$ and~$\Ccirc$ are flat over~$\cO$, it is enough to
  check that the surjection $C\to B$ is also injective. Now, we can
  embed~$C$ into the product of its localizations at all maximal
  ideals, and since~$C$ is Noetherian, it in fact embeds into the
  product of the completions of these localizations. It therefore
  embeds into the product of all of its local Artinian quotients~$A$. We
  claim that any such quotient is in fact obtained from a morphism
  ~$\Ccirc\to\Acirc$,
  where~$\Acirc$ is a finite flat~$\cO$-algebra; since any such
  morphism~$\Ccirc\to\Acirc$ factors through~$\Bcirc$ by assumption,
  it follows that any such~$A$ is a quotient of~$B$, so that 
  the surjection $C \to B$ is indeed injective.

It remains to prove the claim. Let~$\Acirc$ be the image of the
composite $\Ccirc~\to~C\to~A$; we need to show that~$\Acirc$ is a
finite $\cO$-algebra.  This follows from Lemma~\ref{lem:order in the court}
below.\end{proof}

The following lemma is no doubt standard, but we recall a proof
for the sake of completeness.

\begin{lemma}
\label{lem:order in the court}
If $R \to A$ is a morphism of $\cO$-algebras, with $R$ being $p$-adically complete
and Noetherian {\em (}e.g.\ a complete local $\cO$-algebra
with finite residue field{\em )}, and $A$ being 
a finite-dimensional $E$-algebra, then the image of $R$ in $A$
is finite over~$\cO$.
\end{lemma}
\begin{proof}
Since $R$ is Noetherian 
and $p$-adically complete, the Artin--Rees lemma shows
that the same is true of its image in~$A$ (see e.g.\ \cite[Prop.\ 10.13]{MR0242802}).  Thus,
by replacing $R$ by its image in $A$ and $A$ itself by the $E$-span
of this image, we are reduced to checking the following statement:
if $M$ is a $p$-adically complete and separated torsion free
$\cO$-module, and if $M[1/p]$ is finite-dimensional,
then $M$ is finite over $\cO$.  There are many ways to see this,
of course; here is one. 

Since $V:= M[1/p]$ is finite dimensional, we may find a finite spanning
set for this vector space contained in $M;$ if $L$ denotes its
$\cO$-span, then $L$ is a lattice in~$V$.  Since $L$ is compact
with respect to its $p$-adic topology,
the $p$-adic topology on $M$ induces the $p$-adic topology on $L$,
and in particular $L$ is closed in
the $p$-adic topology on~$M$. Thus $M/L$ is a $p$-adically separated submodule
of $V/L$.   Note that $V/L$ consists entirely of $p$-power torsion elements,
and thus that the same is true of $M/L$.

This implies that $\bigcap_n p^n\bigl((M/L)[p^{n+1}]\bigr) \subseteq \bigcap_n p^n(M/L) = 0.$
On the other hand,
$p^n\bigl((M/L)[p^{n+1}]\bigr) \subseteq  (M/L)[p] \subseteq (V/L)[p] \iso L/pL,$
and so $\{p^n\bigl((M/L)[p^{n+1}]\bigr)$ is a decreasing sequence of finite sets,
which thus eventually stabilizes.  We conclude that $p^n\bigl((M/L)[p^{n+1}]\bigr) = 0$
for some sufficiently large value of~$n$;
equivalently, $(M/L)[p^n] = (M/L)[p^{n+1}]$, and thus in fact
$(M/L)[p^n] = (M/L)[p^{\infty}] =  M/L.$
Consequently $L \subseteq M \subseteq p^{-n}L,$
showing that $M$ itself is finite over~$\cO$, as claimed.
\end{proof}

\begin{defn}
  \label{defn: our actual pcris pst stacks}
  Let~$\tau$ be an inertial type,
  and let  $\lambdau$ be a Hodge type.
  To begin with, assume in addition 
  that~$\lambdau$ is effective, and
  choose~$h\ge 0$ such
  that~$\lambdau$ is bounded by~$h$, as well as a finite Galois
  extension ~$L/K$ 
  for which the kernel of~$\tau$ contains~$I_{L}$. 
  We then define~$\cX_{K,d}^{\crys,\lambda,\tau}$ \index{$\cX_{K,d}^{\crys,\lambda,\tau}$}
  to be the scheme-theoretic image of
  the morphism $\cC_{d,\crys,h}^{L/K,\fl,\lambdau,\tau}\to\cX_{K,d}$,
  and~$\cX_{K,d}^{\semis,\lambdau,\tau}$ to be the scheme-theoretic
  image of \index{$\cX_{K,d}^{\semis,\lambda,\tau}$}
  the morphism $\cC_{d,\semis,h}^{L/K,\fl,\lambdau,\tau}\to\cX_{K,d}$.
  (The characterization of these closed substacks of $\cX_{K,d}$
  given in Theorem~\ref{thm: existence of ss stack}
  will show that they are well-defined substacks of $\cX_{K,d}$,
  independently of the auxiliary choices of $h$ and $L/K$ that were used
  in their definition.)

  In general (i.e.\ if $\lambdau$ is not effective), then we choose an
  integer~$h'$ so
  that~$\lambda'_{\sigma,i}:=\lambda_{\sigma,i}+h'\ge 0$ for
  all~$\sigma$ and all~$i$, and define
  $\cX_{K,d}^{\semis,\lambdau,\tau}$ to be the unique substack
  of~$\cX_{K,d}$ with the property that if~$B$ is a
  $\cO/\varpi^a$-algebra for some~$a\ge 1$, then a morphism $\rho:\Spec B\to \cX_{K,d}$
  factors through~ $\cX_{K,d}^{\semis,\lambdau,\tau}$ if and only if the
  morphism $\rho\otimes\epsilon^{h'}:\Spec B\to \cX_{K,d}$ factors
  through ~ $\cX_{K,d}^{\semis,\lambdau',\tau}$. (Here we are using the
  notation of Section~\ref{subsec: tensor product
  of phi gamma and duality}; more precisely, we are twisting by the
\'etale $(\varphi,\Gamma)$-module corresponding to~$\epsilon^{h'}$.)
  We define
  $\cX_{K,d}^{\crys,\lambdau,\tau}$ in the same way. The point of this
  definition is that if~$A^\circ$ is a finite flat $\cO$-algebra, then a
  representation $\rho:G_K\to\GL_d(\Acirc[1/p])$ is potentially
  semistable (resp.\ crystalline) of Hodge type~$\lambdau$ and
  inertial type~$\tau$ if and only if~$\rho\otimes\epsilon^{h'}$ is potentially
  semistable (resp.\ crystalline) of Hodge type~$\lambdau'$ and
  inertial type~$\tau$; again, Theorem~\ref{thm: existence of ss stack} below
  shows that these stacks are defined independently of the choices
  made in this definition, and in particular independently of the
  choice of~$h'$.
\end{defn}
\begin{rem}
  \label{rem: could use non effective BKF modules but no}Rather than
  introducing the twist by~$\epsilon^{h'}$, it would perhaps be
  more natural to work throughout with Breuil--Kisin and
  Breuil--Kisin--Fargues modules which are not necessarily
  $\varphi$-stable, as in~\cite{2016arXiv160203148B}. However, it
  would still frequently be useful to make the corresponding Tate
  twists (see for example the proof of ~\cite[Lem.\
  4.26]{2016arXiv160203148B}), and in particular we would need to make
  such twists in order to use the results
  of~\cite{MR2562795,EGstacktheoreticimages}, so we do not see any
  advantage in doing so.
\end{rem}

The rest of this section is devoted to proving some fundamental
properties of the stacks ~$\cX_{K,d}^{\crys,\lambda,\tau}$ and
~$\cX_{K,d}^{\semis,\lambdau,\tau}$.
%
We begin with an analysis of versal
rings. To this end, fix a point $\Spec \F'\to\cX_{K,d}$  
for some finite extension $\F'$ of $\F$, giving rise to a continuous representation $\rhobar: G_K \to \GL_d(\F')$.
Let $\cO'=W(\F')\otimes_{W(\F)}\cO$, the ring of integers in a finite extension $E'$ of $E$
having residue field $\F'$. Then the versal morphism $\Spf
R_{\rhobar}^{\square,\cO'}\to\cX_{K,d}$ of Proposition~\ref{prop:versal rings}
induces  
morphisms $\Spf R^{\crys,\underline{\lambda},\tau,\cO'}_{\rhobar}\to\cX_{K,d}$
and $\Spf R^{\semis,\underline{\lambda},\tau,\cO'}_{\rhobar}\to\cX_{K,d}$
(where these deformation rings are as in Section~\ref{subsec: notation and conventions}).

\begin{prop}
  \label{prop: versal rings for pst stacks}The morphisms  $\Spf R^{\crys,\underline{\lambda},\tau,\cO'}_{\rhobar}\to\cX_{K,d}$
and $\Spf R^{\semis,\underline{\lambda},\tau,\cO'}_{\rhobar}\to\cX_{K,d}$ factor
through versal morphisms $\Spf
R^{\crys,\underline{\lambda},\tau,\cO'}_{\rhobar}\to\cX_{K,d}^{\crys,\underline{\lambda},\tau}$
and  $\Spf
R^{\semis,\underline{\lambda},\tau,\cO'}_{\rhobar}\to\cX_{K,d}^{\semis,\underline{\lambda},\tau}$ respectively.
\end{prop}
\begin{proof}By definition, we can and do assume that~$\lambdau$ is
  effective. We give the proof for
  $\Spf R^{\crys,\underline{\lambda},\tau,\cO'}_{\rhobar}$, the
  argument in the semistable case being identical. 
  We begin by introducing 
  the fibre product $\widehat{C} :=
  \cC_{d,\crys,h}^{L/K,\fl,\lambdau,\tau}\times_{\cX_{K,d}}\Spf R_{\rhobar}^{\square,\cO'}$.
  Since 
  $\cC_{d,\crys,h}^{L/K,\fl,\lambdau,\tau}\to {\cX_{K,d}}$
  is proper and representable by algebraic spaces, 
  we see that $\widehat{C}$ is a formal algebraic space,
  and that the morphism
  \numequation
  \label{eqn:C to versal R}
  \widehat{C} \to
  \Spf R_{\rhobar}^{\square,\cO'}
  \end{equation}
  is proper and representable by algebraic spaces. The latter
  condition can be reexpressed by saying that~\eqref{eqn:C to versal R}
  is an {\em adic} morphism; thus we see that
  $\widehat{C}$ is a proper
  $\mathfrak m_{R_{\rhobar}^{\square,\cO'}}$-adic
  formal algebraic space over $\Spf R_{\rhobar}^{\square,\cO'}$.

  If we let $\Spf R$ denote the scheme-theoretic image of the
  morphism~\eqref{eqn:C to versal R}, 
  then
  by  Propositions~\ref{prop:X is an Ind-stack}, ~\ref{prop: HT
    weights are a closed condition} and Lemma~\ref{alem: scheme theoretic
      image of versal is versal}, 
the versal morphism
  $\Spf R_{\rhobar}^{\square,\cO'}\to\cX_{K,d}$ induces a versal morphism
  $\Spf R\to\cX_{K,d}^{\crys,\underline{\lambda},\tau}$,
  and so the assertion of the proposition may be rephrased as the claim that
  $R$ and $R^{\crys,\lambdau,\tau,\cO'}_{\rhobar}$
  coincide as quotients of $R_{\rhobar}^{\square,\cO'}.$

  Since
  $\cC_{d,\crys,h}^{L/K,\fl}$ is flat over~$\Spf\cO$, it
  follows from Lemma~\ref{lem:versal flatness} 
  that $R$ is an~$\cO$-flat quotient of~$
  R_{\rhobar}^{\square,\cO'}$. Since
  $R^{\crys,\underline{\lambda},\tau,\cO'}_{\rhobar}$ is also an
  $\cO$-flat quotient of~$R_{\rhobar}^{\square,\cO'}$, we are reduced
  by Proposition~\ref{prop: uniqueness of a flat stack with local ring
            points}
  to showing that if~$\Acirc$ is any finite flat $\cO'$-algebra, then an
  $\cO'$-homomorphism $R_{\rhobar}^{\square,\cO'}\to\Acirc$ factors
  through~$R$ if and only if it factors
  through~$R^{\crys,\underline{\lambda},\tau,\cO'}_{\rhobar}$.

  To this end, note that if the morphism factors
  through~$R^{\crys,\underline{\lambda},\tau,\cO'}_{\rhobar}$, then by
  Theorem~\ref{thm: admits all descents if and only if semistable}
  there is an order~$(\Acirc)'\supseteq \Acirc$ in~$\Acirc[1/p]$ such that the induced
  morphism $\Spf(\Acirc)'\to \Spf R_{\rhobar}^{\square,\cO'}$ lifts to~$\widehat{C}$.
  Consequently the morphism
  $\Spf(\Acirc)'\to \Spf R_{\rhobar}^{\square,\cO'}$ factors
  through~$\Spf R$, and since $\Spf(\Acirc)'\to\Spf\Acirc$ is
  scheme-theoretically dominant, the morphism $\Spf\Acirc\to \Spf
  R_{\rhobar}^{\square,\cO'}$ also factors through~$\Spf R$.

  It remains to prove the converse, namely that if 
  $\Spf \Acirc \to \Spf R_{\rhobar}^{\square,\cO'}$ factors through
  $\Spf R$, then it factors 
  through~$R^{\crys,\underline{\lambda},\tau,\cO'}_{\rhobar}$.
  If we 
  write $S := R\otimes_{R_{\rhobar}^{\square,\cO'}} R^{\crys,
		  \underline{\lambda},\tau,\cO'}_{\rhobar}$
	  (a quotient of $R$),
	  then equivalently,
	  we wish to show that 
	  any morphism $\Spf \Acirc \to \Spf R$
	  in fact factors through $\Spf S$.  


  To do this, we will apply Lemma~\ref{lem: criterion for Artin to map to scheme theoretic
    image}, 
  not in the context of $\cO'$-algebras, but rather in the context of
  $E'$-algebras.  That is, we will invert $p$, and work with the rings
  $R[1/p]$ and $S[1/p]$ (or rather on certain completions of these
  rings).  We would also like to invert $p$ on the object
  $\widehat{C}$, but since the latter is a formal algebraic space,
  doing so directly would lead us into considerations of rigid
  analytic geometry that we prefer to avoid.  Thus, we begin by
  observing that by Lemma~\ref{lem: algebraising C} (together
  with~\cite[Thm.\ V.6.3]{MR0302647}, to pass from
  $\cC_{d,\semis,h}^{L/K}$ to its closed
  substack~$\cC_{d,\crys,h}^{L/K,\fl}$), $\widehat{C}$ can be promoted
  to a projective scheme over $\Spec R$.  That is, there is a
  projective morphism of schemes
  $C \to \Spec R_{\rhobar}^{\square,\cO'}$, with scheme-theoretic
  image equal to $\Spec R,$ whose
  $\mathfrak m_{R_{\rhobar}^{\square,\cO'}}$-adic completion is
  isomorphic to $\widehat{C}$.  Write $C[1/p] := C\otimes_{\cO'} E'.$
  Since scheme-theoretic dominance is preserved by flat base-change,
  the morphism $C[1/p] \to \Spec R[1/p]$ is proper and
  scheme-theoretically dominant.
  
  Returning to our argument,
  we want to show, for any finite flat $\cO'$-algebra~$\Acirc$,
  that a morphism $R\to \Acirc$ necessarily factors through~$S$,
  or, equivalently,
  that the induced morphism
  $R[1/p] \to A := \Acirc[1/p]$ 
  factors through~$S[1/p]$.
  Since $A$ is a product of Artinian local $E$-algebras,
  we can check this assertion factor-by-factor;
  replacing $A$ by one of these factors (and $\Acirc$ by the image
  of $\Acirc$ in this factor)
  we may thus assume that $A$ is local.  
  We then let $\widehat{R[1/p]}$ denote the completion of
  $R$ at the kernel of its map to $A$, 
  and we let $\widehat{S[1/p]}$ denote the corresponding 
  completion of $S[1/p]$;
  by Artin--Rees, we also have that $\widehat{S[1/p]} = \widehat{R[1/p]} 
  \otimes_{R[1/p]} S[1/p]$, and so the map $R[1/p] \to A$ factors through
  $S[1/p]$
  if and only if the induced morphism $\widehat{R[1/p]} \to A$
  factors through $\widehat{S[1/p]}$.

  To show the desired factorization, 
  we see from Lemma~\ref{lem: criterion for Artin to map to scheme theoretic
    image} 
  that it suffices to show that if $\widehat{R[1/p]} \to B$
  is a morphism to a finite Artinian $E$-algebra
  for which $C_B \to \Spec B$ admits a section,
  then this morphism factors through~$\widehat{S[1/p]}$.

Write~$B^\circ$ to denote the image of the composite
$R_{\rhobar}^{\square,\cO'}\to \widehat{R[1/p]} \to B$; by Lemma~\ref{lem:order in
the court}, we see that $B^{\circ}$ is
an order in $B$.  Let $Z\hookrightarrow C_{B^{\circ}}$ denote the scheme-theoretic
  closure of the section $\Spec B\to C_B$  in $C_{B^{\circ}}$. 
  This is proper over $\Spec \Bcirc$, flat over $\cO$,
  and irreducible of dimension one.
  Thus $Z$ is finite over $\Spec \Bcirc$, and hence $Z = \Spec (\Bcirc)'$
  for some order $(\Bcirc)'$ in $B$ containing~ $\Bcirc$.
  Thinking of $Z$ as a section of $C$ over $(\Bcirc)'$,
  we find that the induced morphism $\Spf (\Bcirc)' \to \cX_{K,d}$
  lifts to a morphism $\Spf (\Bcirc)' \to 
  \cC_{d,\crys,h}^{L/K,\fl,\lambdau, \tau}$. 
  Theorem~\ref{thm: admits all descents if and only if
    semistable} 
  then implies that the morphism $R_{\rhobar}^{\square, \cO'} 
  \to \Bcirc$ factors through $R_{\rhobar}^{\crys,\lambdau,\tau,\cO'}$,
  and thus the morphism $\widehat{R[1/p]} \to B$ does indeed
  factor through $\widehat{S[1/p]}$, as required.
%
%
%
\end{proof}

\begin{thm}\label{thm: existence of ss stack}Let~$\tau$ be an inertial type, and let  $\lambdau$ be a
  Hodge type. Then the closed substacks~$\cX_{K,d}^{\crys,\lambdau,\tau}$ and
  ~$\cX_{K,d}^{\semis,\lambdau,\tau}$ of~$\cX_{K,d}$ are  $p$-adic formal
  algebraic stacks 
  which
  are of finite type and flat over~$\Spf\cO$, and  are
  uniquely 
  determined as $\cO$-flat
  closed substacks of~$\cX_{K,d}$ by the following property:
  if~$\Acirc$ is a finite flat~$\cO$-algebra, then
  $\cX_{K,d}^{\semis,\lambda,\tau}(\Acirc)$ \emph{(}resp.\
  $\cX_{K,d}^{\crys,\lambda,\tau}(\Acirc)$\emph{)} is precisely the
  subgroupoid of~$\cX_{K,d}(\Acirc)$ consisting of
  $G_K$-representations which become potentially semistable
  \emph{(}resp.\ potentially crystalline\emph{)} of Hodge
  type~$\lambdau$ and inertia type~$\tau$ after inverting~$p$.
\end{thm}
\begin{proof}
We can and do assume that~$\lambdau$ is effective. By Propositions~\ref{prop:X is an Ind-stack}, ~\ref{prop: HT weights are a closed condition} and~\ref{aprop: scheme theoretic image is p adic finite
    type II}, $\cX_{K,d}^{\crys,\lambdau,\tau}$ and
  ~$\cX_{K,d}^{\semis,\lambdau,\tau}$ of~$\cX_{K,d}$ are  $p$-adic formal
  algebraic stacks  which
  are of finite type and flat over~$\Spf\cO$.  

  We now verify the claimed description of the $\Acirc$-valued points
  of $\cX_{K,d}^{\semis,\lambda,\tau}(\Acirc)$ and
  $\cX_{K,d}^{\crys,\lambda,\tau}(\Acirc)$, for finite flat $\cO$-algebras
  $\Acirc$. Since the argument is identical in either case, we give it
  in the semistable case. Any such algebra~$\Acirc$ is a product of
  finitely many finite flat local $\cO$-algebras, and so we immediately
  reduce to verifying the claim in the case when $\Acirc$ is
  furthermore local.  Since $\Acirc$ is then a complete local finite
  flat $\cO$-algebra, the claimed description follows from
  Proposition~\ref{prop: versal rings for pst stacks}. Indeed, the
  residue field~$\F'$ is a finite extension of~$\F$, and if as above
  we set~$\cO'=W(\F')\otimes_{W(\F)}\cO$, then the morphism
  $\Spf\Acirc\to\cX_{K,d}$ factors
  through~$\cX_{K,d}^{\semis,\underline{\lambda},\tau}$ if and only if it
  factors through the versal morphism
  $\Spf
  R^{\semis,\underline{\lambda},\tau,\cO'}_{\rhobar}\to\cX_{K,d}^{\semis,\underline{\lambda},\tau}$
  of Proposition~\ref{prop: versal rings for pst stacks}. In turn,
  this happens if and only if the corresponding $G_K$-representation
  is potentially semistable of Hodge type~$\lambdau$ and inertial
  type~$\tau$, by the defining property
  of~$R^{\semis,\underline{\lambda},\tau,\cO'}_{\rhobar}$. 

It remains to show the claimed uniqueness statement. We prove in Corollary~\ref{cor:Xd is formal algebraic} below
	that $\cX_{K,d}$ is a Noetherian formal algebraic stack,
	and in Theorem~\ref{thm:Xdred is algebraic} below
	that $\cX_{d,\red}$ is of finite presentation over $\F$.
	The reader can easily confirm that
	Chapter~\ref{sec: families of extensions}
	contains no citation to the present chapter,
	and thus that those results are independent of the present ones;
	in particular, it is safe to invoke them here.
        These results imply that any closed substack $\cY$ of $\cX_{K,d}$
	is Noetherian (so that it makes sense to speak of $\cY$ being
	flat over $\cO$), and that $\cY_{\red}$ is of finite type
	over $\cF$.  The claimed uniqueness then follows from
	Proposition~\ref{prop: uniqueness of a flat stack with local ring
    points}.
\end{proof}

We will make use of the following corollary in Chapter~\ref{sec:properties}.
\begin{cor}
  \label{cor: crystalline deformation rings are effective
    versal}
For any Hodge type $\underline{\lambda}$, and
  for any~$a\ge 1$,
  the corresponding morphism
  $\Spf R^{\crys,\underline{\lambda},\cO'}_{\rhobar}/\varpi^a \to
  \cX_{K,d}^a$ is effective, i.e.\ is induced by a morphism $\Spec R^{\crys,\underline{\lambda},\cO'}_{\rhobar}/\varpi^a \to
  \cX_{K,d}^a$.
\end{cor}
\begin{proof}By Proposition~\ref{prop: versal rings for pst stacks},
  the morphism $\Spf R^{\crys,\underline{\lambda},\cO'}_{\rhobar}/\varpi^a \to
  \cX_{K,d}^a$ factors through a versal morphism $\Spf R^{\crys,\underline{\lambda},\cO'}_{\rhobar}/\varpi^a \to
  \cX_{K,d}^{\crys,\lambda}\times_{\Spf\cO}\Spec\cO/\varpi^a$. Since~$\cX_{K,d}^{\crys,\lambda}$
  is a $p$-adic formal algebraic stack (by Theorem~\ref{thm: existence of ss stack})
  the base-change
  $\cX_{K,d}^{\crys,\lambda}\times_{\Spf\cO}\Spec\cO/\varpi^a$ is an algebraic
  stack, and so the result follows from~\cite[\href{https://stacks.math.columbia.edu/tag/07X8}{Tag 07X8}]{stacks-project}.  
\end{proof}

Finally, we can compute the dimensions of our potentially crystalline
and semistable stacks. It is presumably possible to develop the
dimension theory of $p$-adic formal algebraic stacks in some
generality, but we do not need to do so, as we can instead work with
the special fibres, which are algebraic stacks.
\begin{thm}
  \label{thm: dimension of ss stack}  
The algebraic stacks~$\overline{\cX}_{K,d}^{\crys,\lambdau,\tau}:=\cX_{K,d}^{\crys,\lambdau,\tau}\times_{\Spf\cO}\Spec \F$ and
  ~$\overline{\cX}_{K,d}^{\semis,\lambdau,\tau}:=\cX_{K,d}^{\semis,\lambdau,\tau}\times_{\Spf\cO}\Spec \F$ are equidimensional of dimension \[\sum_{\sigma}\#\{1\le i<j\le
  d|\lambda_{\sigma,i}>\lambda_{\sigma,j}\}.\]%
In particular, if~$\lambdau$ is regular, then the algebraic stacks~$\overline{\cX}_{K,d}^{\crys,\lambdau,\tau}$ and
  ~$\overline{\cX}_{K,d}^{\semis,\lambdau,\tau}$ are equidimensional of
  dimension~$[K:\Qp]d(d-1)/2$.   
\end{thm}
\begin{proof}Once again, we give the argument in the crystalline case,
  the semistable case being identical. Write
  $d_{\lambdau}:=\sum_{\sigma}\#\{1\le i<j\le
  d|\lambda_{\sigma,i}>\lambda_{\sigma,j}\}$.
 The
  algebraic stack
  $\cX_{K,d}^{\crys,\lambdau,\tau}\times_{\Spf\cO}\Spec \F$
  is of finite type over ~$\F$, and we need to show that it is
  equidimensional of dimension~$d_{\lambdau}$.
  Let~$x:\Spec \F'\to\cX_{K,d}^{\crys,\lambdau,\tau}$ be a finite type
  point corresponding to a Galois representation~$\rhobar$,
with corresponding versal morphism
  $\Spf
  R^{\crys,\underline{\lambda},\cO'}_{\rhobar}/\varpi\to\overline{\cX}_{K,d}^{\crys,\lambdau,\tau}
  $. 
We have the fibre product 
  $$\Spf R^{\crys,\underline{\lambda},\cO'}_{\rhobar}/\varpi
  \times_{\widehat{(\overline{\cX}_{K,d}^{\crys,\lambdau,\tau})}_x} 
  \Spf R^{\crys,\underline{\lambda},\cO'}_{\rhobar}/\varpi
\iso \widehat{\GL}_{d,R^{\crys,\underline{\lambda},\cO'}_{\rhobar}/\varpi,1},$$
where $\widehat{\GL}_{d,R^{\crys,\underline{\lambda},\cO'}_{\rhobar}/\varpi,1}$
denotes the completion of 
$({\GL}_{d})_{R^{\crys,\underline{\lambda},\cO'}_{\rhobar}/\varpi}$
along the identity element in its special fibre.
  By~\cite[Lem.\ 2.40]{EGcomponents}, 
it is therefore enough to recall that since
  $ R^{\crys,\underline{\lambda},\cO'}_{\rhobar}$ is equidimensional
  of dimension~$d_{\lambdau}+1$, $ R^{\crys,\underline{\lambda},\cO'}_{\rhobar}/\varpi$ is
  equidimensional of dimension~$d_{\lambdau}$ (see~\cite[Lem.\ 2.1]{MR3248725}).
\end{proof}

\chapter{Families of extensions}\label{sec: families of extensions}
In this chapter we extend the theory of the Herr complex 
to the context of $(\varphi,\Gamma)$-modules with coefficients,
and use it to
develop a theory of families of extensions of
$(\varphi,\Gamma)$-modules. 
We then 
inductively construct families which cover the underlying reduced
stack $\cX_{K,d,\red}$, and employ obstruction theory to deduce
that~$\cX_{K,d}$ is a Noetherian formal algebraic stack.

\section{The Herr complex}\label{subsec: Herr complex}
In this chapter we consider the Herr complex; our approach is informed
by the papers~\cite{MR1693457,MR1839766,MR1645070,MR2416996,MR3117501,MR3230818}, and in
particular we follow~\cite{MR3117501,MR3230818} in considering it as an object
of the derived category. Our main technical result is to show that it
is a perfect complex, which we will do by showing that it satisfies
the following well-known criterion.

\begin{lem}
  \label{lem: criterion for perfectness of complex} Let~$A$ be a
  Noetherian commutative ring, and let~$C^\bullet$ be an object
  of~$D(A)$. Then the following two conditions are equivalent:
  \begin{enumerate}
  \item There is a quasi-isomorphism $F^\bullet\to C^\bullet$
    where~$F^\bullet$ is a complex of flat $A$-modules, concentrated
    in a finite number of degrees;  and the cohomology groups of~$C^\bullet$ are all finite
    $A$-modules.
  \item $C^\bullet$ is perfect.
  \end{enumerate}
\end{lem}
\begin{proof}
  This follows
  from~\cite[\href{http://stacks.math.columbia.edu/tag/0658}{Tag
    0658}]{stacks-project},
  \cite[\href{http://stacks.math.columbia.edu/tag/0654}{Tag
    0654}]{stacks-project} and
  \cite[\href{http://stacks.math.columbia.edu/tag/066E}{Tag
    066E}]{stacks-project}.
\end{proof}

Let $A$ be a $p$-adically complete $\cO$-algebra, and let $M$ be an
$A$-module with commuting $A$-linear endomorphisms $\varphi,\Gamma$
(in our main applications, $M$ will be a
projective \'etale $(\varphi,\Gamma)$-module
with~$A$-coefficients, but it will be convenient in our arguments to be
able to consider more general possibilities, such as subquotients of
base changes of such modules). 
Then the~\emph{Herr complex} of~$M$ is by \index{Herr complex}
definition the complex of~$A$-modules $\cC^\bullet(M)$ in degrees
$0,1,2$ given by
\[\xymatrix{M\ar[rr]^{(\varphi-1,\gamma-1)}&&
      M\oplus M\ar[rr]^{(\gamma-1)\oplus(1-\varphi)}&&M.}\]

 We have the
following useful interpretation of the cohomology groups of the Herr complex.
\begin{lem}
  \label{lem: Yoneda for H0 H1 of Herr complex}Let $M_1$, $M_2$ be
  projective \'etale $(\varphi,\Gamma)$-modules with
  $A$-coefficients. Then for $i=0,1$ there are natural
  isomorphisms \[H^i\bigl(\cC^\bullet(M_1^\vee\otimes M_2)\bigr)\cong\Ext^i_{\A_{K,A},\varphi,\Gamma}(M_1,M_2). \]
\end{lem}
\begin{proof}
  This is straightforward; in
  particular the case $i=0$ is immediate. When $i=1$, any extension
  $M$ of~$M_1$ by~$M_2$ splits on the level of the underlying projective
  $\A_{K,A}$-modules, and such a splitting is unique up to an element of
  $\Hom_{\A_{K,A}}(M_1,M_2)=M_1^\vee\otimes M_2$. Given such a
  splitting of~$M$, we obtain two elements~$f,g$ of $M_1^\vee\otimes M_2$ by
  writing \[\varphi_M=
    \begin{pmatrix}
      \varphi_{M_2}&\varphi_{M_2}\circ f\\0&\varphi_{ M_1}
    \end{pmatrix},\ \gamma_M=
    \begin{pmatrix}
      \gamma_{M_2}&\gamma_{M_2}\circ g\\0&\gamma_{M_1}
    \end{pmatrix}.
\]The condition that~$\varphi_M$, $\gamma_M$ commute shows that
$(f,g)$ determines a class in $H^1\bigl(\cC^\bullet(M_1^\vee\otimes M_2)\bigr)$
(this class is easily seen to be well-defined, by our above remark
about the ambiguity in the choice of splitting). We leave the
verification that this gives the claimed bijection to the sufficiently
enthusiastic reader.
\end{proof}

In order to show that the Herr complex is a perfect complex, we will
need to develop a little of the theory of the~$\psi$-operator
on~$(\varphi,\Gamma)$-modules. We will only need this in the case
that~$A$ is an~$\F$-algebra, and we mostly work in this context from
now on. Note that for any~$K$, if $A$ is an~$\F$-algebra, then
$\A_{K,A}^+$ is $(\varphi,\Gamma)$-stable, and~$(\A_{K,A}')^+$ is
$(\varphi,\Gammat)$-stable; indeed we have~$\varphi(T)=T^p$
and~$\gamma(T)-T\in T^2\A^+_{K,A}$ by Lemma~\ref{lem: gamma on T}, and
the stability under the action of~$\Gammat$ follows by an identical
argument to the proof of Lemma~\ref{lem: gamma on T} (that is, it
follows from the continuity of the action of~$\Gammat$).

We claim that for any~$K$, ~$1,\varepsilon,\dots,\varepsilon^{p-1}$ is
a basis for~$\bE'_{K}$ as a $\varphi(\bE'_{K})$-vector
space. 
To see this, note that since the
extension $\bE'_K/\varphi(\bE'_K)$ is inseparable of degree~$p$, while
$\bE'_K/\bE'_{\Qp}$ is a separable extension, it is enough to show the
result for~$\bE'_{\Qp}$. In this case we have
$\bE'_{\Qp}=\Fp((\varepsilon-1))$, so $\varphi(\bE'_{\Qp})=\Fp((\varepsilon^p-1))$
and the claim is clear. Then for any~$x\in\bE'_K$ we can write
$x=\sum_{i=0}^{p-1}\varepsilon^i\varphi(x_i)$, and we define
$\psi:\bE'_K\to\bE'_K $ by 
$\psi(x)=x_0$. By definition~$\psi$ is continuous, $\Fp$-linear, commutes
with~$\Gammat$, and satisfies $\psi\circ\varphi=\id$. Since~$\psi$
commutes with~$\Delta$, we have an induced map $\psi:\bE_K\to\bE_K$,
which again is continuous, $\Fp$-linear, commutes
with~$\Gamma$, and satisfies $\psi\circ\varphi=\id$.


If~$A$ is an $\Fp$-algebra, then since~$\psi$ is continuous we can
extend scalars from~$\Fp$ to~$A$ and complete to obtain continuous
$A$-linear maps $\psi:\A'_{K,A}\to\A'_{K,A}$ and $\psi:\A_{K,A}\to\A_{K,A}$. Again, these are continuous
and commute with~$\Gammat$ and~$\Gamma$, and satisfy $\psi\circ\varphi=\id$. We
have \[\A'_{K,A}=\oplus_{i=0}^{p-1}\varepsilon^i\varphi(\A'_{K,A}),\]and
if $x\in\A'_{K,A}$ and we
write $x=\sum_{i=0}^{p-1}\varepsilon^i\varphi(x_i)$, then
$\psi(x)=x_0$.

\begin{prop}
  \label{prop: existence of psi on phi Gamma modules}Let $A$ be an
  $\F$-algebra, 
  and let $M$ be a projective
  \'etale $(\varphi,\Gamma)$-module with $A$-coefficients. Then there
  is a continuous and open $A$-linear surjection~$\psi:M\to M$ such that $\psi$
  commutes with~$\Gamma$, and we
  have \[\psi(\varphi(a)m)=a\psi(m), \] \[\psi(a\varphi(m))=\psi(a)m\]
  for any $a\in \A_{K,A}$, $m\in M$.
\end{prop}
\begin{proof}We define $\psi:M\to M$ to be the
  composite
  \[M\stackrel{\Phi_M^{-1}}{\to} \A_{K,A}\otimes_{\A_{K,A},\varphi}M
    \stackrel{\psi\otimes 1}{\to}\A_{K,A}\otimes_{\A_{K,A}}M=M.  \]
  The relations
  $\psi(\varphi(a)m)=a\psi(m)$, $\psi(a\varphi(m))=\psi(a)m$ are
  immediate from the definitions.
That ~$\psi$ is continuous follows
  from the continuity of $\Phi_M^{-1}$ and  the continuity of~ $\psi$
  on~$\A_{K,A}$, 
  while the
   surjectivity of~$\psi$ is immediate from the relation $\psi(\varphi(m))=m$.
 That~$\psi$ is open follows from Lemma~\ref{lem: bounds on psi} below.
 \end{proof}


  \begin{lem}\label{lem: bounds on psi}Let~$A$ be an
    $\F$-algebra, let~$M$ be a projective \'etale
    $(\varphi,\Gamma)$-module with $A$-coefficients, and let~$\gM$ be
    a $\varphi$-stable lattice in~$M$. 
    Then there is an integer~$h\ge 0$ such that 
    for any integer~$n$, we
    have \[\psi(T^{h+np}\gM)\subseteq T^n\gM\subseteq \psi(T^{np}\gM). \]
    
  \end{lem}
  \begin{proof}
Since
    $\psi\circ\varphi=\id$, for any~$n$ we
    have \[\psi(T^{np}\varphi(\gM))=T^n\gM.\]
    This implies that \[T^n\gM=\psi(T^{np}\varphi(\gM))\subseteq
      \psi(T^{np}\gM). \]
For the other direction,
choose~$h$ such that \[T^h\gM\subseteq\Phi_{M}(\varphi^*\gM)\subseteq\gM.\]
Then we have
 \[\psi(T^h\gM)\subseteq\psi(\Phi_{M}(\varphi^*\gM))=
    \gM.\]
  It follows that \[\psi(T^{h+np}\gM)=T^n\psi(T^h\gM)\subseteq
    T^n\gM, \]as required.
  \end{proof}

  \begin{lem}
    \label{lem: existence of phi gamma lattice}Let~$A$ be a Noetherian
    $\F$-algebra, and let~$M$ be a projective \'etale
    $(\varphi,\Gamma)$-module with $A$-coefficients. Then~$M$ contains
    a  $(\varphi,\Gamma)$-stable lattice.
  \end{lem}
  \begin{proof}By~\cite[Lem.\ 5.2.15]{EGstacktheoreticimages}, $M$ contains
    a $\varphi$-stable lattice~$\gM$. Let~$\Gamma\gM$ be the
    $\A^+_{K,A}$-submodule of~$M$ generated by the elements~$\gamma m$
    with~$\gamma\in\Gamma$, $m\in\gM$. We claim that~$\Gamma \gM$ is a
    lattice; note that since~$\Gamma$ and~$\varphi$ commute,~$\Gamma \gM$ is
    $\varphi$-stable, and it is $\Gamma$-stable because~$\A^+_{K,A}$
    is $\Gamma$-stable. To see that it is a lattice, it is
    enough (since $\A^+_{K,A}$ is Noetherian) to show that it is contained
    in a lattice (as it certainly spans~$M$). In particular, it is
    enough to show that $\Gamma\gM\subseteq T^{-n}\gM$ for some $n\ge
    0$. This follows easily from the compactness of~$\Gamma$; for
    example, if $e_1,\dots,e_m$ are generators of the finitely
    generated $\A^+_{K,A}$-module $\gM$, then we have a continuous map
    $j:\Gamma\to M^m$, $\gamma\mapsto (\gamma e_1,\dots,\gamma e_d)$,
    and  $\Gamma=\cup_nj^{-1}(T^{-n}\gM^m) $ is an open cover
    of~$\Gamma$. Since this has a finite subcover we are done.
  \end{proof}

    \begin{lem}
    \label{lem: gamma on T unit}Let~$A$ be an~$\F$-algebra. Then for any~$i\in \Z$ we have
    $\gamma(T^{i})\in T^i(\A_{K,A}^+)^\times$. 
  \end{lem}
  \begin{proof}
    This follows from Lemma~\ref{lem: gamma on T}. Indeed, we can
    write~$\gamma(T)=T+\lambda$ with $\lambda\in T^2\A_{K,A}^+$, so for
    any~$i\ge 1$, we have $\gamma(T^i)=(T+\lambda)^i$, and
    $\gamma(T^i)-T^i\in T^{i+1}\A_{K,A}^+$, whence $\gamma(T^{i})\in
    T^i(\A_{K,A}^+)^\times$. Since~$\gamma(T^i)\gamma(T^{-i})=1$, the
    result then also holds for~$i\le 0$.
  \end{proof}
\begin{cor}
    \label{cor: existence of psi gamma lattice}Let~$A$ be a Noetherian
    $\F$-algebra, and let~$M$ be a projective \'etale
    $(\varphi,\Gamma)$-module with $A$-coefficients. Then~$M$ contains
    a lattice~$\gM'$ such that for all~$m\ge 0$, $T^{-m}\gM'$ is $(\psi,\Gamma)$-stable.
  \end{cor}
  \begin{proof}
    By Lemma~\ref{lem: existence of phi gamma lattice} there is a
    $(\varphi,\Gamma)$-stable lattice~$\gM$, and by Lemma~\ref{lem:
      bounds on psi}, the lattice $T^{-n}\gM$ is
    $\psi$-stable for all $n\gg 0$. It is also $\gamma$-stable for
    all~$n\ge 0$ by Lemma~\ref{lem: gamma on T unit}, so we may take $\gM'=T^{-n}\gM$ for any sufficiently
    large value of~$n$.
  \end{proof}


  If~$M$ is an \'etale $(\varphi,\Gamma)$-module with $A$-coefficients, then we
  write \[M':=\A'_{K,A}\otimes_{\A_{K,A}}M,\] an \'etale
  $(\varphi,\Gammat)$-module with $A$-coefficients.
  \begin{lem}
    \label{lem: ker psi is a direct summand}Let~$A$ be an~$\F$-algebra, and let~$M$ be a projective \'etale
    $(\varphi,\Gamma)$-module with $A$-coefficients. Then we have a
    decomposition of
    $\varphi(\A'_{K,A})$-modules $M'=\ker(\psi)\oplus\varphi(M')$,
    and we may write $\ker(\psi)=\oplus_{i=1}^{p-1}\varepsilon^i\varphi(M')$.
  \end{lem}
  \begin{proof}
    Since $\psi\circ \phi$ is the identity on~$M'$, the endomorphism $\phi\circ \psi$
    of $M'$ is idempotent, and so $M'$ decomposes as the direct sum of its
    kernel and image.  Since $\phi$ is injective, while $\psi$ is surjective,
    we may express this decomposition as $M' = \ker(\psi) \oplus \varphi(M')$,
    as claimed.

    As~$M'$ is \'etale, we have $\A'_{K,A}\otimes_{\varphi(\A'_{K,A})}
\varphi(M') \iso M'$. Since
    $\A'_{K,A}=\oplus_{i=0}^{p-1}\varepsilon^i\varphi(\A'_{K,A})$, we
    find that
    $M'=\oplus_{i=0}^{p-1}\varepsilon^i\varphi(\A'_{K,A})\varphi(M')=\oplus_{i=0}^{p-1}\varepsilon^i\varphi(M')$. Now,
    for~$1\le i\le p-1$ we have
    $\psi(\varepsilon^i\varphi(M'))=\psi(\varepsilon^i)M'=0$, so
    $\oplus_{i=1}^{p-1}\varepsilon^i\varphi(M')\subseteq\ker(\psi)$. Since
    we have
    $\oplus_{i=1}^{p-1}\varepsilon^i\varphi(M')\oplus\varphi(M')=\ker(\psi)\oplus\varphi(M')$,
    it follows that this inclusion is an equality, as required.
  \end{proof}

  The proof of the following result follows the approach
  of~\cite{MR1645070}, although one difference in our situation is
  that because we are working in characteristic~$p$, we do not need to
  reduce to the case that~$K/\Qp$ is unramified. 
  \begin{prop}
  \label{prop: inverting 1-gamma}Let~$A$ be a Noetherian 
  $\F$-algebra. 
  If~$M$ is a projective \'etale
  $(\varphi,\Gamma)$-module with $A$-coefficients, then $(1-\gamma)$
  is bijective on~$\ker(\psi)$.
\end{prop}
\begin{proof}
We begin with some preliminaries on the action of~$\gamma$ on~$\varepsilon\in(\A'_{K,A})^+$. We have $\gamma(\varepsilon)=\varepsilon^{\chi(\gamma)}$, and we 
write $\chi(\gamma)=1+p^Nu$ where $N\ge 1$ and~$u\in\Z_p^\times$. Then
for each $n\ge 1$ we can write $\chi(\gamma^{p^{n-1}})=1+p^{n+N-1}u_n$
with $u_n\in\Z_p^\times$, and if~$0\le r\le n+N-1$ and~$i\in\Z$ then we have
\numequation\label{eqn: moving epsilon into phi}
\gamma^{p^{n-1}}(\varepsilon^i)=\varepsilon^i\varphi^r(\varepsilon^{ip^{n+N-1-r}u_n}). \end{equation}Note
that since~$u_n$ is a $p$-adic unit, $(\varepsilon^{u_n}-1)/(\varepsilon-1)$ is a
unit in~$(\A'_{K,A})^+$, and we
can
write \numequation\label{eqn: gamma-1 on epsilon minus 1 equation}(\gamma^{p^{n-1}}-1)(\varepsilon-1)=\varepsilon(\varepsilon^{u_n}-1)^{p^{n+N-1}}\in
  (\varepsilon-1)^{p^{n+N-1}}((\A'_{K,A})^+)^\times. \end{equation}It follows that for any~$r\in\Z$, we
have \numequation\label{eqn: gamma minus 1 on epsilon minus 1}(\gamma^{p^{n-1}}-1)((\varepsilon-1)^r)\in
  (\varepsilon-1)^{r+p^{n+N-1}-1}(\A'_{K,A})^+. \end{equation} Indeed
if~$r\ge 0$ then we may
write \[(\gamma^{p^{n-1}}-1)((\varepsilon-1)^r)=(\gamma^{p^{n-1}}-1)((\varepsilon-1))\cdot\sum_{j=0}^{r-1}(\gamma^{p^{n-1}}(\varepsilon-1)^{r-1-j})\cdot(\varepsilon-1)^j, \]so
\eqref{eqn: gamma minus 1 on epsilon minus 1} follows from~\eqref{eqn:
  moving epsilon into phi} and the fact that
$\gamma^{p^{n-1}}(\varepsilon-1)\in(\varepsilon-1)(\A'_{K,A})^+$,
which in turn follows from~\eqref{eqn: gamma-1 on epsilon minus 1
  equation}. The case~$r\le 0$ then follows from the result for~$-r$ by
writing \[(\gamma^{p^{n-1}}-1)((\varepsilon-1)^r)=-\frac{(\gamma^{p^{n-1}}-1)((\varepsilon-1)^{-r})}{(\gamma^{p^{n-1}}(\varepsilon-1)^{-r})(\varepsilon-1)^{-r}},\]and
noting that $(\gamma^{p^{n-1}}(\varepsilon-1))/(\varepsilon-1)=(\varepsilon^{1+p^{n+N-1}u_n}-1)/(\varepsilon-1)$ is a unit.

   Since both~$\psi$ and~$\gamma$
  commute with the action of~$\Delta$, in order to prove the
  proposition it suffices to show
  that~$(1-\gamma)$ is bijective on the kernel of $\psi:M'\to
  M'$.  By Lemma~\ref{lem: ker psi is a direct summand}, this kernel equals
$\ker(\psi)=\oplus_{i=1}^{p-1}\varepsilon^i\varphi(M')$. 
It follows from~\eqref{eqn: moving epsilon into phi} that $\gamma$
preserves~$\varepsilon^i\varphi(M')$ for any integer~$i$, so it is
enough to prove that if $(i,p)=1$  that~$(1-\gamma)$ acts
invertibly on $\varepsilon^i\varphi(M')$.

In fact, we will prove the 
stronger statement that
for each~$n\ge 1$, each $0\le t\le N-1$, and each~$(i,p)=1$, the operator $(1-\gamma^{p^{n-1}})$ acts
invertibly on~$\varepsilon^i\varphi^{n+t}(M')$ (note that
$(1-\gamma^{p^{n-1}})$ acts on this space by~\eqref{eqn: moving epsilon into phi}).

Writing \numequation\label{eqn: yet another eqn in psi invertible gamma proof}\varepsilon^i\varphi^{n+t-1}(M')=\varepsilon^i\varphi^{n+t-1}(\oplus_{j=0}^{p-1}\varepsilon^j\varphi(M'))=\oplus_{j=0}^{p-1}\varepsilon^{i+p^{n+t-1}j}\varphi^{n+t}(M'),\end{equation}we
see that if~$t>0$ then the statement for some pair~$(n,t)$ (and
all~$i$) implies the statement for~$(n,t-1)$. Writing
$(1-\gamma^{p^{n-1}})=(1-\gamma^{p^{n-2}})(1+\gamma^{p^{n-2}}+\dots+\gamma^{p^{n-2}(p-1)})$,
and again using~\eqref{eqn: yet another eqn in psi invertible gamma proof},
we see also that the statement for~$(n,t)$ implies the statement
for~$(n-1,t)$. We can therefore assume that~$n$ is arbitrarily large
and that~$t=N-1$.  

Let~$\gM$ be a
$(\psi,\Gamma)$-stable lattice in~$M$ (which exists by
Corollary~\ref{cor: existence of psi gamma lattice}), let~$\gM'=(\A'_{K,A})^+\otimes_{\A_{K,A}}\gM$, and choose~$n$
large enough that 
$p^{n+N-1}\ge 3$
and \numequation\label{eqn: squeezing by 2 on M}(\gamma^{p^{n-1}}-1)(\gM')\subseteq (\varepsilon-1)^2\gM'.\end{equation}It follows that \numequation\label{eqn:gamma -1 squeezing}(\gamma^{p^{n-1}}-1)((\varepsilon-1)^r\gM')\subseteq (\varepsilon-1)^{r+2}\gM'\end{equation}
for all $r\in \Z$, because if $m\in \gM'$
then \[(\gamma^{p^{n-1}}-1)((\varepsilon-1)^rm)=(\gamma^{p^{n-1}}-1)((\varepsilon-1)^r)\cdot\gamma^{p^{n-1}}(m)+(\varepsilon-1)^r(\gamma^{p^{n-1}}-1)(m),\]so
that ~\eqref{eqn:gamma -1 squeezing} follows from ~\eqref{eqn: gamma
  minus 1 on epsilon minus 1} and~\eqref{eqn: squeezing by 2 on M},
together with our assumption that~$p^{n+N-1}\ge 3$.


By~\eqref{eqn: moving epsilon into phi}, we
have \[(\gamma^{p^{n-1}}-1)(\varepsilon^i\varphi^{n+N-1}(x))=\varepsilon^{i}\varphi^{n+N-1}(\varepsilon^{iu_n}\gamma^{p^{n-1}}(x)-x), \]so
it is enough to check that the map $f:M'\to M'$ given
by \[f(x)=\varepsilon^{iu_n}\gamma^{p^{n-1}}(x)-x \]is
invertible.

Let $\alpha=\varepsilon^{iu_n}-1$, so that~$\alpha/(\varepsilon-1)$ is a unit
in~$(\A'_{K,A})^+$. Then for any $r\in\Z$ and~$x\in (\varepsilon-1)^r\gM'$,  
it follows from~(\ref{eqn:gamma -1 squeezing})
that
\[\left(\frac{1}{\alpha}f-1\right)(x)=\frac{\varepsilon^{iu_n}}{\alpha}(\gamma^{p^{n-1}}-1)(x)\in
  (\varepsilon-1)^{r+1}\gM'. \] In particular the
sum
\[g(x):=\sum_{j=0}^\infty\left(1-\frac{1}{\alpha}f\right)^j(x) \]
converges. Since $f,g$ are additive by definition, we have \[\left(1-\frac{1}{\alpha}f\right)\circ g(x)=g\circ
\left(1-\frac{1}{\alpha}f\right)(x)=g(x)-x, \] so that the function $g:M\to M$ is a
left and right inverse to~$\frac{1}{\alpha}f$. Thus~$f$ is
invertible, as required. 
\end{proof}

Suppose that~$A$ is an $\F$-algebra. Using the~$\psi$ operator, we can give an alternative description of
the Herr complex, which will be important in establishing that it is a
perfect complex. Let $\cC_\psi^\bullet(M)$ be the complex in degrees
$0,1,2$ given by \[\xymatrix{M\ar[rr]^{(\psi -1,\gamma-1)}&&
      M\oplus M\ar[rr]^{(\gamma-1)\oplus(1 - \psi)}&&M. }\]

  We have a morphism of complexes
  $\cC^\bullet(M)\to\cC^\bullet_\psi(M)$ given by
  \numequation
  \label{eqn:Herr quasi-iso}
  \xymatrix{M\ar[d]^{1}\ar[rr]^{(\varphi-1,\gamma-1)}&&
      M\oplus M\ar[d]^{(-\psi,1)}\ar[rr]^{(\gamma-1)\oplus(1-\varphi)}&&M\ar[d]^{-\psi} \\ M\ar[rr]^{(\psi-1,\gamma-1)}&&
      M\oplus M\ar[rr]^{(\gamma-1)\oplus(1-\psi)}&&M} 
\end{equation}
(That this is a morphism of complexes follows from the facts that
  $\psi\circ\varphi=\id$, and that~$\psi$ commutes with~$\gamma$.)

\begin{prop}
  \label{prop: Herr complex in terms of psi}
  Let~$A$ be a Noetherian 
  $\F$-algebra, 
  and let $M$ be a projective
  \'etale $(\varphi,\Gamma)$-module with $A$-coefficients. Then
  the morphism of complexes
  $\cC^\bullet(M)\to\cC^\bullet_\psi(M)$
  defined by~{\em (\ref{eqn:Herr quasi-iso})}
  is a quasi-isomorphism.
\end{prop}
\begin{proof}
The cokernel of~(\ref{eqn:Herr quasi-iso})  
is the zero complex, while its 
kernel is the complex \[\xymatrix{0\ar[r]&
     \ker(\psi)\ar[r]^{(\gamma-1)}&\ker(\psi), }\]which is acyclic by Proposition~\ref{prop: inverting
  1-gamma}. 
\end{proof}

In fact, we require a slightly stronger statement than the quasi-isomorphism
of the preceding lemma.  Namely, we need a statement that takes into
account topologies, which will make use of the following lemma.  

\begin{lemma}\label{lem:open maps}
	If $f:X \to Y$ is a continuous open map of topological spaces, 
	and if $Y'\subseteq Y$ is the inclusion of a subspace
	into $Y$ {\em (}i.e.\ $Y'$ is a subset of $Y$
	endowed with the induced topology{\em )},
	then the base-changed map $X' := f^{-1}(Y') \to Y'$ is again open,
	when $X'$ is endowed with the topology induced by that of $X$.
\end{lemma}
\begin{proof}
	Let $U'$ be an open subset of $X'$; then we may write
	$U' = U \cap X'$, for some open subset $U$ of $X$.
	Thus $$f(U') = f(U \cap X') = f\bigl(U\cap f^{-1}(Y')\bigr)
	= f(U) \cap Y'.$$
	Since $f$ is open, we conclude that $f(U')$ is an open
	subset of $Y'$, as required.
\end{proof}

We now have the following strengthening of
Proposition~\ref{prop: Herr complex in terms of psi},
which takes into account the topologies on the complexes.

\begin{prop}
	\label{prop:topological quasi-iso}
  Let~$A$ be a countable Noetherian  $\F$-algebra
  and let $M$ be a projective
  \'etale $(\varphi,\Gamma)$-module with $A$-coefficients. Then the morphism
  of complexes
  $\cC^\bullet(M)\to\cC^\bullet_\psi(M)$
  defined by~{\em (\ref{eqn:Herr quasi-iso})}
  induces topological isomorphisms on each of the associated cohomology modules.
\end{prop}
\begin{proof}
  Proposition~\ref{prop: Herr complex in terms of psi} shows
  that~(\ref{eqn:Herr quasi-iso}) is a quasi-isomorphism. By
  definition, the induced map on $H^0$ is
  the identity map, and is therefore a homeomorphism. Since~$\psi$
  is continuous and open, we see that the maps $\cC^i(M)\to\cC^i_\psi(M)$ are
  continuous and open for each~$i$; in particular, the
 isomorphism on~$H^2$ is induced from the continuous open map
 $-\psi:\cC^2(M)\to\cC^2_\psi(M)$, and is therefore a homeomorphism.
 (Here and below
 we use the standard fact that a quotient morphism of topological
 groups is necessarily open.)

 The case of~$H^1$ requires a little more work. Since the maps
 $\cC^0(M)\to\cC^0_\psi(M)$ and $H^1\bigl(\cC^{\bullet}(M)\bigr) \to 
 H^1\bigl(\cC_{\psi}^{\bullet}(M)\bigr)$ are isomorphisms,
 we see that the continuous morphism
 $\cC^1(M)\to\cC^1_\psi(M)$ induces a continuous isomorphism of the
 modules of cocycles $Z^1(M)\to Z^1_\psi(M)$, and it suffices to show
 that this map is open (since this will imply that the induced
 bijection on $H^1$ is both continuous and open, and thus is an isomorphism
 of topological groups).

 If we let $\tZ^1(M)$ denote the preimage of~$Z^1_\psi(M)$ in $\cC^1(M)$,
and if we let $K\subseteq \cC^1(M)=M\oplus M$ denote $\ker(\psi)\oplus 0$,
 then we have inclusions $Z^1(M), K \subseteq \tZ^1(M)$,
 and hence a continuous morphism
 \numequation
 \label{eqn:continuous bijection}
 Z^1(M) \oplus K \to \tZ^1(M),
 \end{equation}
 which is in fact a bijection.
 Each of $Z^1(M)$, $K$, and $\tZ^1(M)$ is closed
 subgroup of $\cC^1(M) = M\oplus M$.  Lemma~\ref{lem:polish}
 shows that this latter topological group is Polish,
 and thus so is any of its closed subgroups.  
 Corollary~\ref{cor:open polish} then shows that~(\ref{eqn:continuous bijection})
 is in fact a homeomorphism,
while Lemma~\ref{lem:open maps} shows that the morphism $\tZ^1(M) \to 
Z^1_\psi(M)$ is open.  Consequently, we find that the morphism
$Z^1(M) \to Z^1_{\psi}(M)$ is open, as required.
\end{proof}

\begin{lem}
  \label{lem: discrete quotients of projective}Let~$A$ be a Noetherian
  $\F$-algebra,
  and let~$X$ be an
  $A$-module   subquotient of a finitely generated projective
  $A[[T]]$-module, endowed with its natural subquotient
  topology. If the topology on~$X$ is discrete, then $X$ is finitely
  generated as an $A$-module.
\end{lem}
\begin{proof}Write $X=Y/Z$, where~$Y$ is an $A$-submodule of a
  finitely generated projective $A[[T]]$-module~$\gM$. Since the
  topology on~$X$ is discrete, $Z$ is open in~$Y$, and thus
  $Z\supseteq U\cap Y$, for some open neighbourhood $U$ of zero in $\gM$.
  Such a neighbourhood $U$ contains
  $T^n\gM$ for some $n\gg 0$, and so we find that~$X$ is a subquotient of the
  finitely generated $A$-module $\gM/T^n\gM$, as required.
\end{proof}

\begin{thm}
  \label{thm: Herr complex is perfect}\leavevmode Let~$A$ be a Noetherian 
    $\cO/\varpi^a$-algebra such that $A/\varpi$ is
    countable, 
    and let~$M$ be a projective \'etale $(\varphi,\Gamma)$-module with
    $A$-coefficients. 
  \begin{enumerate}
  \item  The Herr complex~$\cC^\bullet(M)$ is a
    perfect complex of $A$-modules, concentrated in
    degrees~$[0,2]$.

  \item   If either (i) $B$ is a finite $A$-algebra;
or (ii) $B$ is a finite type $A$-algebra and $A$ itself 
is of finite type over~$\cO/\varpi^a$;
    then there is a natural isomorphism in the derived
    category
    \[\cC^\bullet(M)\otimes^{\mathbb{L}}_AB\isoto\cC^\bullet(M\otimes_{\A_{K,A}}\A_{K,B}). \]
    In particular, there is a natural isomorphism
    \[H^2(\cC^\bullet(M))\otimes_AB\isoto
      H^2(\cC^\bullet(M\otimes_{\A_{K,A}}\A_{K,B})). \]
  \end{enumerate}

\end{thm}
\begin{proof}Since~$A$ is Noetherian, $\A_{K,A}$ is a flat~$A$-algebra
  (being a localisation of the power series ring~$\A^+_{K,A}$), so
  ~$\cC^\bullet(M)$ is a complex of flat $A$-modules.  By
  \cite[\href{https://stacks.math.columbia.edu/tag/07LU}{Tag
    07LU}]{stacks-project} (for part~(1), taking~$R=\cO/\varpi^a$), and \cite[Prop.\
  2.2.2]{pilloniHidacomplexes} (for part~(2), again taking~$R=\cO/\varpi^a$), it suffices to prove
  the result in the case that~$A$ is an $\F$-algebra, which we assume
  from now on. We begin with~(1). By Lemma~\ref{lem: criterion for
    perfectness of complex}, we need only check that the cohomology
  groups of~$\cC^\bullet(M)$ are finitely generated
  $A$-modules. 
In order to do this, we will make two truncation arguments. 

Firstly, by Lemma~\ref{lem: existence of phi gamma lattice}, we can choose a
 $(\varphi,\Gamma)$-stable lattice $\gM\subseteq M$. 
%
%
We claim that for every~$n\ge 1$ the Herr complex
$\cC^\bullet(T^n\gM)$ is acyclic; consequently,  the natural morphism of
complexes $\cC^\bullet(M)\to\cC^\bullet(M/T^n\gM)$ is a
quasi-isomorphism.

To see the claim, it suffices to show that $(1-\varphi):T^n\gM\to T^n\gM$
  is an isomorphism; indeed, the
  exactness of~$\cC^\bullet(T^n\gM)$ is a formal consequence of this. We begin by checking injectivity. If $m\in T^n\gM$ and
$(1-\varphi)(m)=0$, then we have
$m=\varphi(m)\in T^{pn}\gM$ 
so that in particular 
we have $m\in T^{n+1}\gM$. By induction, we see
that $m\in T^n\gM$ for all~$n$, so that~$m=0$, as required.

We now prove surjectivity. If $m\in T^n\gM$, then we have seen in
the previous paragraph that~$\varphi(m)\in T^{n+1}\gM$,
$\varphi^2(m)\in T^{n+2}\gM$, and so on. Since $\gM$ is $T$-adically
complete, we can set $x=\sum_{i\ge 0}\varphi^i(m)\in T^n\gM$; then
$(1-\varphi)(x)=m$, as required.  


We now turn to our other truncation argument. Let $\gM'$ be as in
Lemma~\ref{cor: existence of psi gamma lattice}, so that for each $n'\le 0$, 
$T^{n'}\gM'$ is a $(\psi,\Gamma)$-stable lattice in~$M$. 
In particular, for each~$n'\le 0$, $\cC_\psi^\bullet(T^{n'}\gM)$ is a
subcomplex of $\cC_\psi^\bullet(M)$. We claim that if
$n'\le -2$ then $\cC_\psi^\bullet(M/T^{n'}\gM')$ is acyclic,
so that the natural morphism of complexes
$\cC_\psi^\bullet(T^{n'}\gM')\to\cC_\psi^\bullet(M)$ is a
quasi-isomorphism. 


As above, it is enough to show that~$(1-\psi)$ is bijective on
$M/T^{n'}\gM'$. Note firstly that if $r\le n'$
then $\psi(T^r\gM')\subseteq T^{r+1}\gM'$. Indeed, 
for each integer $r\le 0$ we
have \[\psi(T^r\gM')\subseteq\psi(T^{p\lfloor r/p\rfloor}\gM')
  =T^{\lfloor r/p\rfloor}\psi(\gM')\subseteq
  T^{\lfloor r/p\rfloor}\gM',\] and the claim follows since
if~$r\le -2$ then 
$\lfloor r/p\rfloor\ge r+1$.

We now show that~$(1-\psi)$ is injective on $M/T^{-n}\gM'$. Suppose that $m\in M$ with
$(1-\psi)(m)\in T^{n'}\gM'$. Then $m\in T^r\gM'$ for some $r\ll 0$. If
$r\ge n'$ then we are done, and if not then since $\psi(m)\in
T^{r+1}\gM'$ and $(1-\psi)(m)\in T^{n'}\gM\subseteq T^{r+1}\gM$, we
have $m\in T^{r+1}\gM$. It follows by induction that $m\in
T^{n'}\gM'$, as required.

For surjectivity, take $m\in M$, and choose $r\ll 0$ such that $m\in
T^r\gM'$. If $r< n'$ then $\psi(m)\in T^{r+1}\gM'$, so by induction we
see that there is some $s\ge 0$ such that $\psi^s(m)\in
T^{n'}\gM'$. If we set $x=\sum_{i=0}^{s-1}\psi^i(m)$, then
$(1-\psi)(x)=m-\psi^{s}(m)\in m+ T^{n'}\gM$, as required.

 We now consider the quasi-isomorphisms
\[\cC^\bullet_\psi(T^{n'}\gM')\to\cC^\bullet_\psi(M)\leftarrow
  \cC^\bullet(M)\to\cC^\bullet(M/T^n\gM).  \] 

Proposition~\ref{prop:topological quasi-iso}
shows that the middle quasi-isomorphism
induces a topological isomorphism on cohomology modules,
and thus altogether we obtain morphisms
$$H^i\bigl( \cC^{\bullet}_{\psi}(T^{n'}\gM')\bigr) \to H^i\bigl(
\cC^{\bullet}(M/T^n\gM)\bigr)$$
which are isomorphisms of $A$-modules, and continuous
with respect to the natural topologies on the source
and target.
Since the target is endowed with the discrete topology,
we find that the source is also endowed with the discrete
topology. By Lemma~\ref{lem: discrete quotients of projective}, it
follows that 
these cohomology modules are indeed finitely generated
over $A$, and hence so are the cohomology modules of
$\cC^{\bullet}(M)$, as required.
%

We now turn to~(2), where we follow the proof of~\cite[Thm.\
4.4.3]{MR3230818}. Since both~(i) and~(ii) in particular require~$B$
to be a finite type $A$-algebra, we see
that~$B$ is Noetherian, and~$B$ is countable, so that in particular
$\cC^\bullet(M\otimes_{\A_{K,A}}\A_{K,B})$ is a perfect complex of
$B$-modules by part~(1) (replacing $A$ by~$B$, and~$M$ by $M\otimes_{\A_{K,A}}\A_{K,B}$). Since~$\cC^\bullet(M)$ is a perfect complex
of $A$-modules, $\cC^\bullet(M)\otimes^{\mathbb{L}}_AB$ is also a perfect
complex of $B$-modules, and the natural map $\A_{K,A}\otimes_AB\to\A_{K,B}$
induces a morphism \numequation\label{eqn: base change for Herr complex}\cC^\bullet(M)\otimes^{\mathbb{L}}_A B\to
  \cC^\bullet(M\otimes_{\A_{K,A}}\A_{K,B}).\end{equation}

In case~(i),
when~$B$ is in fact finite as an $A$-module, the natural map
$\A_{K,A}\otimes_AB\to\A_{K,B}$ is an isomorphism (indeed, we have
$A[[T]]\otimes_AB=B[[T]]$, because $A[[T]]\otimes_AB$ is finitely
generated over~$A[[T]]$ and thus $T$-adically complete and separated
by the Artin--Rees lemma), so~(\ref{eqn: base change for Herr
  complex}) is certainly an isomorphism in this case.

We now turn to  
case~(ii),
when~$B$ is only assumed to be a finite type $A$-algebra,
but~$A$, and hence also $B$, is furthermore of finite type over~$\cO/\varpi^a$.
This implies in particular that $A$ is Jacobson,
so that if~$\m$ is any
maximal ideal of~$B$, then~$B/\m$ is finite as an $A$-module.
From this, and the already-proved case~(i),
we conclude that (\ref{eqn: base change for Herr
  complex}) is an isomorphism if we replace~$B$ by~$B/\m$. 

It follows that we have a chain of
quasi-isomorphisms \[(\cC^\bullet(M)\otimes^{\mathbb{L}}_A
  B)\otimes^{\mathbb{L}}_B B/\m \isoto \cC^\bullet(M)\otimes^{\mathbb{L}}_A B/\m \isoto \cC^\bullet(M\otimes_A B/\m)\]\[\isoto
  \cC^\bullet(M\otimes_{\A_{K,A}}\A_{K,B})\otimes^{\mathbb{L}}_B B/\m, \]
where the final quasi-isomorphism comes from replacing~$A$ by~$B$, $B$
by~$B/\m$, and~$M$ by $M\otimes_{\A_{K,A}}\A_{K,B}$ in~(\ref{eqn: base change
  for Herr complex}). (That this is indeed a quasi-isomorphism follows
by case~(i), 
since~$B/\m$ is a finite $B$-module, while we have
$(M\otimes_{\A_{K,A}}\A_{K,B})\otimes_{\A_{K,B}}\A_{K,B/\m}
=M\otimes_{\A_{K,A}}\A_{K,B/\m}=M\otimes_{\A_{K,A}}(\A_{K,A}\otimes_AB/\m)=M\otimes_AB/\m$).  

Thus~ (\ref{eqn: base change for Herr complex}) becomes a
quasi-isomorphism after applying~$\otimes^{\mathbb{L}}_B B/\m$, for
any maximal ideal~$\m$ of~$B$. It follows from~\cite[Lem.\
4.1.5]{MR3230818} that (\ref{eqn: base change for Herr complex}) is a
quasi-isomorphism, as required. Finally, the compatibility of the formation of~$H^2$ with base change
is an immediate consequence of this isomorphism in the derived
category, together with the vanishing of all of the higher degree
cohomology groups.
\end{proof}

\begin{cor}
  \label{cor: Herr complex is perfect for p-adic A}Let $A$ be a
  $p$-adically complete Noetherian $\cO$-algebra such that~$A/\varpi$ is countable, and let $M$ be a projective
  \'etale $(\varphi,\Gamma)$-module. Then the Herr complex
  $\cC^\bullet(M)$ is a perfect complex concentrated in degrees $[0,2]$.
\end{cor}
\begin{proof}This follows
  from~\cite[\href{https://stacks.math.columbia.edu/tag/0CQG}{Tag 0CQG}]{stacks-project} and Theorem~\ref{thm: Herr complex is perfect}.
(We apply the base-change statement of part~(2) of the theorem
to the surjective (and hence finite) morphisms
$A/\varpi^{a+1} \to A/\varpi^a$.)
\end{proof}

We also note the following rather technical corollary
Theorem~\ref{thm: Herr complex is perfect},
which we will use in our discussion of families of extensions below.

\begin{cor}
\label{cor:Herr complex length}
If $A$  is of finite type over $\cO/\varpi^a$ for some $a  \geq 1$,
and $M$ is a projective \'etale  $(\varphi,\Gamma)$-module,
then the Herr complex $\cC^{\bullet}(M)$ can be represented
by a complex $C^0 \to C^1 \to C^2$ of finite projective $A$-modules
in degrees~$[0,2]$.
\end{cor}
\begin{proof}
Theorem~\ref{thm: Herr complex is perfect}
shows that $\cC^{\bullet}(M)$  is perfect,
and so 
by~\cite[\href{https://stacks.math.columbia.edu/tag/0658}{Tag 0658}]{stacks-project}, 
the present corollary will follow  provided that we show that
$\cC^{\bullet}(M)$ has tor-amplitude in~[0,2].
Since any $A$-module is a filtered direct limit of finite $A$-modules,
it suffices to  show that
$\cC^{\bullet}(M)\Lotimes_A M$
has cohomological amplitude lying in~$[0,2]$
for any finite $A$-module~$M$.   
If we let $B  = A  \oplus M$, thought of as an $A$-algebra by  declaring
$M$ to be  a square-zero  ideal,
then it is  equivalent to show that
$\cC^{\bullet}(M)\Lotimes_A B$
has cohomological amplitude lying in~$[0,2]$.
But Theorem~\ref{thm: Herr complex  is perfect}
shows  that this latter complex is isomorphic (in the derived category of~$B$-modules)
to
$\cC^{\bullet}(M\otimes_{\A_{K,A}}  \A_{K,B}),$
whose cohomological amplitude  lies in~[0,2] by construction.
\end{proof}

We now explain the comparison between the Herr complex and Galois
cohomology, and the relationship to Tate local duality. These results
are essentially due to Herr and are proved
in~\cite{MR1693457,MR1839766} but we follow~\cite{MR3230818} and
formulate them as statements in the derived category.


Exactly as in~\cite[Defn.\ 2.3.10]{MR3230818}, if $M_1$, $M_2$ are
projective \'etale $(\varphi,\Gamma)$-modules with $A$-coefficients,
then there is a cup
product \[\cC^\bullet(M_1)\otimes_A\cC^\bullet(M_2)\to\cC^\bullet(M_1\otimes_{\A_{K,A}}M_2).\]

More precisely, the cup product arises from the following
generalities. If we have two complexes of $A$-modules~$C^\bullet$
and~$D^\bullet$, then the tensor product $C^\bullet\otimes D^\bullet$
has differential given by \[d(x\otimes y)=dx\otimes y+(-1)^i x\otimes
  y\] if $x\in C^i$, $y\in D^j$. If $f^\bullet:C^\bullet\to
C^\bullet$, then we write $\Fib(f|C^\bullet):=\Cone(f)[-1]$, which by definition
has $\Fib(f)^i=C^i\oplus C^{i-1}$ and
$d^i((x,y))=(d^i(x),-d^{i-1}(y)-f^i(x))$.

As a special case of~\cite[Lem.\ 2.3.9]{MR3230818}, if
$f_1:C^\bullet_1\to C^\bullet_1$, $f_2:C^\bullet_2\to C^\bullet_2$,
then we have a natural
morphism \numequation\label{eqn: cup product of fibre complex}\Fib(1-f_1|C_1^\bullet)\otimes \Fib(1-f_2|C_2^\bullet)
  \to
  \Fib(1-f_1\otimes f_2|C_1^\bullet\otimes C^\bullet_2) \end{equation}

Then by definition we have
$\cC^\bullet(M)=\Fib\bigl(1-\gamma|\Fib(1-\varphi|M)\bigr)$, so that by multiple
applications of~(\ref{eqn: cup product of fibre complex}) we have morphisms
\begin{align*}
  \cC^\bullet(M_1)\otimes_A\cC^\bullet(M_2)&=\Fib\bigl(1-\gamma|\Fib(1-\varphi|M_1)\bigr)\otimes_A\Fib\bigl(1-\gamma|\Fib(1-\varphi|M_2)\bigr)\\
  &\to \Fib\Bigl(1-\gamma|\bigl(\Fib(1-\varphi|M_1)\otimes_A\Fib(1-\varphi|M_2)\bigr)\Bigr)\\& \to
  \Fib\Bigl(1-\gamma|\Fib\bigl(1-\varphi|(M_1\otimes_A M_2)\bigr)\Bigr)\\&\onto  \Fib\Bigl(1-\gamma|\Fib\bigl(1-\varphi|(M_1\otimes_{\A_{K,A}} M_2)\bigr)\Bigr)\\
&=\cC^\bullet(M_1\otimes M_2)
\end{align*}
whose composite defines the cup product.


\begin{lem}\label{lem: H2 of Tate twist}
  If~$A$ is a finite type $\cO/\varpi^a$-algebra, then there is an isomorphism
  \[H^2(\cC^\bullet(\A_{K,A}(1)))\cong A,\] compatible with base change.
\end{lem}
\begin{proof}
  Since by Theorem~\ref{thm: Herr complex is perfect}~(2) the formation of~$H^2(\cC^\bullet(M))$ is
  compatible with base change, the result follows from the
  case~$A=\cO/\varpi^a$, which is immediate from Theorem~\ref{thm: Herr complex
    and Galois reps} below (that is, from the natural isomorphism
  $H^2\bigl(G_K,(\cO/\varpi^a)(1)\bigr)~\cong~\cO/\varpi^a$). 
\end{proof}

If $M$ is a projective \'etale $(\varphi,\Gamma)$-module 
over a finite type $\cO/\varpi^a$-algebra,
then we  define the Tate duality
pairing between the Herr complexes of~$M$ and of its Cartier dual~$M^*$ as the following composite of the cup product, truncation, and the
isomorphism of Lemma~\ref{lem: H2 of Tate twist}: \[\cC^\bullet(M)\times\cC^\bullet(M^*)\to
  \cC^\bullet(\A_{K,A}(1)))\to H^2(\cC^\bullet(\A_{K,A}(1)))[-2]\cong A[-2].\]


\begin{prop}
  \label{prop: Herr complex Tate duality and Euler char}Let $A$ be a
  finite type $\cO/\varpi^a$-algebra, and let $M$ be a projective
  \'etale $(\varphi,\Gamma)$-module. 
  \begin{enumerate}
  \item The Tate duality pairing induces a quasi-isomorphism \[\cC^\bullet(M)\isoto\RHom_A(\cC^\bullet(M^*),A))[-2]. \]
  \item If $A$ is a finite extension of~$\F$, then the Euler characteristic~$\chi_A(\cC^\bullet(M))$ is equal to
    $-[K:\Qp]d$.

  \end{enumerate}
\end{prop}
\begin{proof}For part~(1), by~\cite[Lem.\ 4.1.5]{MR3230818}, it is enough to treat
  the case that~$A$ is a field, in which case~$A$ is a finite extension
  of~$\F$. Then both parts  follow from Theorem~\ref{thm: Herr complex
    and Galois reps} below (that is, from the corresponding statements for
  Galois representations).
%
\end{proof}

Finally we recall the relationship between the Herr complex and Galois
cohomology. As in~\cite[\S2]{MR3117501} it is possible to upgrade the
following theorem to an isomorphism in the derived category, but as we
do not need this we do not give the details here. If~$A$ is a complete local Noetherian
  $\cO$-algebra with finite residue field, and~$M$ is a formal
  projective \'etale $(\varphi,\Gamma)$-module with $A$-coefficients,
  then we can define the Herr complex~$\cC^\bullet(M)$ exactly as for
  \'etale $(\varphi,\Gamma)$-modules, so that by definition we have
  $\cC^\bullet(M)=\varprojlim_n\cC^\bullet(M_{A/\m_A^n})$. By \cite[\href{https://stacks.math.columbia.edu/tag/0CQG}{Tag 0CQG}]{stacks-project}, 
  $\cC^\bullet(M)$ is a perfect complex, concentrated in degrees~$[0,2]$.
\begin{thm}
  \label{thm: Herr complex and Galois reps}If~$A$ is a complete local Noetherian
  $\cO$-algebra with finite residue field, and $T$ is a finitely
  generated projective $A$-module with a continuous action
  of~$G_K$, 
  then there
  are isomorphisms of $A$-modules \[H^i(G_K,T) \isoto H^i(\cC^\bullet(\mathbb{D}_A(T)))
    \] 
  which are
  functorial in~$T$ and compatible with cup products and duality.
\end{thm}
\begin{proof}This is~\cite[Prop.\ 3.1.1]{MR1805474}, which is deduced from the results
  of~\cite{MR1839766,MR1693457} by passage to the limit. 
%
  %
\end{proof}
\begin{rem}
  \label{rem: writing Herr complex as Galois cohomology}In accordance
  with our general convention of writing $\rho_T$ for a family
  $T\to\cX_d$, if~$T$ is an affine scheme of finite type
  over~$\cO/\varpi^a$, then we write $H^2(G_K,\rho_T)$ for the
  pull-back to~$T$ of the cohomology group $H^2(\cC^\bullet(M))$
  on~$\cX_d$, where~$M$ is the \'etale $(\varphi,\Gamma)$-module
  corresponding to the morphism $T\to\cX_d$. By Theorem~\ref{thm: Herr
    complex is perfect}, $H^2(G_K,\rho_T)$ is a coherent sheaf, whose
  formation is compatible with arbitrary finite type base-change, and
  so in particular by Theorem~\ref{thm: Herr complex and Galois reps}
  its specialisations at $\Fpbar$-points coincide with the usual
  Galois cohomology groups. This compatibility with base-change also
  allows us to extend the definition of $H^2(G_K,\rho_T)$ to arbitrary
  (not necessarily affine) schemes~$T$ of finite type over~$\cO/\varpi^a$.
\end{rem}
\begin{lem}\label{lem: Tate duality for H2 to H0 in families}
Let~$T$ be a scheme of finite type over~$\cO/\varpi^a$, 
and let~$\rho_T$, $\rho'_T$ be families of Galois representations
  over~$T$. Then sections of~$H^2(G_K,\rho_T\otimes(\rho'_{T})^\vee(1))$
  are in natural bijection with homomorphisms $\rho_T\to\rho'_{T}$.
\end{lem}
\begin{rem}
  \label{rem: clarifying what a morphism of families means}Here by a
  homomorphism $\rho_T\to\rho'_T$ we mean a homomorphism of the
  corresponding $(\varphi,\Gamma)$-modules.
\end{rem}
\begin{proof}
  [Proof of Lemma~\ref{lem: Tate duality for H2 to H0 in
    families}]
This follows from   Proposition~\ref{prop: Herr complex Tate duality
  and Euler char} and Lemma~\ref{lem: Yoneda for H0 H1 of Herr complex}.
\end{proof}
\subsection{Obstruction theory}\label{subsubsec: obstruction
  theory}
We now show that the Herr complex provides~$\cX_d$ with a nice
obstruction theory in the sense of Definition~\ref{df:obstruction}.

Let $A$ be a $p$-adically complete $\cO$-algebra, and let $x:\Spec A\to
\cX_d$ be a morphism, corresponding to a projective $(\varphi,\Gamma)$-module
$M$. We wish to consider the problem of deforming~$M$ to a square zero
thickening of~$A$. Specifically, if  \[0\to I\to A'\to A\to 0\] is a
square zero extension, then we let $\Lift(x,A')$ be the set of
isomorphism classes of projective $(\varphi,\Gamma)$-modules $M'$ with
$A'$-coefficients which have the property that $M'\otimes_{A'}A\cong M$. 

For any such thickening~$A'$, we define a corresponding obstruction
class as follows. The underlying $\A_{K,A}$-module of~$M$ has a unique (up
to isomorphism) lifting to a projective $\A_{K,A'}$-module~$\tM$, and
we may lift $\varphi,\gamma$ to semi-linear endomorphisms~$\varphit,\gammat$
of~$\tM$. (Indeed, $\A_{K,A'}$ is a square zero thickening of~$\A_{K,A}$,
and a finite projective module~$P$ over any ring~$R$ deforms
uniquely through any square zero extension $R'\to R$, as its
deformations are controlled by
\[H^1(\Spec R,\ker(R'\to R)\otimes\End_R(P)),\] which vanishes (as
$\Spec R$ is affine). To see that $\varphi,\gamma$ lift, think of them
as $\A_{K,A}$-linear maps $\varphi^*M\to M$, $\gamma^*M\to M$.)

However, there is no guarantee that we can find lifts
~$\varphit,\gammat$ which commute. To measure the obstruction to the
existence of such lifts, let $\ad M=\Hom_{\A_{K,A}}(M,M)$
be the adjoint of~$M$. This naturally has the structure of
a~$(\varphi,\Gamma)$-module; indeed, we have a natural identification $\varphi^*\ad
M=\Hom_{\A_{K,A}}(\varphi^*M,\varphi^*M)$, and we define $\Phi_{\ad M}:\varphi^*\ad
M\to\ad M$ by \[(\Phi_{\ad M}(f))(x)=\Phi_M(f(\Phi_M^{-1}(x))). \]We
define the action of~$\Gamma$ in the analogous way. Then we let $o_x(A')$ be the image in
$H^2(\cC^\bullet(\ad M))\otimes_A I=H^2(\cC^\bullet(\ad M\otimes_A I))$ of
\[\varphit\gammat\varphit^{-1}\gammat^{-1}-1\in \ad M\otimes_A
I=\cC^2(\ad M\otimes_AI).\]

\begin{lem}
  \label{lem: obstruction class well defined and controls existence of
  lifts}
  The cohomology class $o_x(A')$ is well-defined independently
  of the choice of~$\varphit,\gammat$, and vanishes if and only if
  $\Lift(x,A')\ne 0$.
\end{lem}
\begin{proof}Any other liftings $\varphit',\gammat'$ are  obtained from
  our given liftings $\varphit,\gammat$ by setting
  $\varphit'=(1+X)\varphit$, $\gammat'=(1-Y)\gammat$, for some $X,Y\in\ad
  M\otimes_A I$. A simple computation shows
  that \[\varphit'\gammat'(\varphit')^{-1}(\gammat')^{-1}-\varphit\gammat\varphit^{-1}\gammat^{-1}=(1-\gamma)X+(1-\varphi)Y, \]which
  shows that the cohomology class $o_x(A')$ is well-defined, and that it vanishes if and
  only if we can choose $\varphit',\gammat'$ so that
  $\varphit'\gammat'(\varphit')^{-1}(\gammat')^{-1}=1$, which is in
  turn equivalent to $\Lift(x,A')\ne 0$, as required.  
\end{proof}

If~$F$ is an $A$-module, we let~$A[F]:=A\oplus
F$ be the $A$-algebra with multiplication given
by \[(a,m)(a',m'):=(aa',am'+a'm).\] This is a square zero thickening
of~$A$, and $\Lift(x,A[F])\ne 0$, because we have the trivial lifting
given by $M\otimes_AA[F]$.

\begin{lem}
  \label{lem: H1 of Herr complex of adjoint controls lifts}
 Suppose that~$F$ is a finitely generated $A$-module. Then
    there is a natural isomorphism of $A$-modules $\Lift(x,A[F])\isoto
    H^1(\cC^\bullet(\ad M)\otimes_A^{\mathbb{L}}F)$. 
\end{lem}
\begin{proof}We begin by constructing the isomorphism on the level of sets. Note that $\cC^\bullet(\ad M)\otimes^{\mathbb{L}}F$ is
  computed by $\cC^\bullet(\ad M\otimes_AF)$. Liftings of~$M$
  to~$A[F]$ are determined by the corresponding liftings
  $\varphit,\gammat$ of $\varphi,\gamma$, and we obtain a class in
  $H^1(\cC^\bullet(\ad M\otimes_AF))$ by taking the image of
  of \[(\varphit\varphi^{-1}-1,\gammat\gamma^{-1}-1)\in\ad M\otimes_A F\oplus\ad M\otimes_A F =\cC^1(\ad
    M\otimes_AF).\] An elementary calculation shows that
  $(\varphit\varphi^{-1}-1,\gammat\gamma^{-1}-1)$ is in the kernel of $\cC^1(\ad
    M\otimes_AF)\to\cC^2(\ad M\otimes_AF)$ if and only if
    $\varphit\gammat=\gammat\varphit$, so it only remains to check
    that the lifting given by $\varphit,\gammat$ is trivial if and
    only if the corresponding cohomology class vanishes. 
  
    To see this, note that the endomorphisms of the trivial lifting
    are of the form $1+X$, for some $X\in\ad M\otimes F$. The corresponding
    $\varphit,\gammat$ are given by $\varphit=(1+X)\varphi(1+X)^{-1}$,
    $\gammat=(1+X)\gamma(1+X)^{-1}$, which is equivalent to
    $(\varphit\varphi^{-1}-1,\gammat\gamma^{-1}-1)=((1-\varphi)
    X,(1-\gamma) X)$, which by definition is equivalent to
    $(\varphit\varphi^{-1}-1,\gammat\gamma^{-1}-1)$ being a
    coboundary, as required.

To compare the $A$-module structures, recall that by definition the
$A$-module structure on $\Lift(x,A[F])$ is defined as follows
(see~\cite[\href{http://stacks.math.columbia.edu/tag/07Y9}{Tag
  07Y9}]{stacks-project}). If $r\in A$, then we have a homomorphism
$f_r:A[F]\to A[F]$ given by $f_r(a,f)=(a,rf)$, and given a lifting
$\tM$ of~$M$, we let~$r\tM$ be the base change of~$\tM$
via~$f_r$. Explicitly, this means that we replace~$(\varphit-\varphi)$
and~$(\gammat-\gamma)$ by ~$r(\varphit-\varphi)$
and~$r(\gammat-\gamma)$.

Similarly, the addition map
$\Lift(x,A[F])\times\Lift(x,A[F])\to \Lift(x,A[F])$ comes from the
obvious identification
$\Lift(x,A[F])\times\Lift(x,A[F])=\Lift(x,A[F\times F])$ together with
base change via the homomorphism $A[F\times F]\to A[F]$ given by
$(a,f_1,f_2)\mapsto (a,f_1+f_2)$. With obvious notation, this amounts
to setting $\varphit_1\boxplus\varphit_2:=\varphit_1+\varphit_2-\varphi$, 
$\gammat_1\boxplus\gammat_2:=\gammat_1+\gammat_2-\gamma$.

On the other hand, by definition the $A$-module structure on
$H^1(\cC^\bullet(\ad M)\otimes^{\mathbb{L}}F)$ is given by the obvious
$A$-module structure on the pairs
$(\varphit\varphi^{-1}-1,\gammat\gamma^{-1}-1)$. This is obviously the
same as the $A$-module structure that we have just explicated on~$\Lift(x,A[F])$.
%
\end{proof}

\begin{prop}
  \label{prop: Xd has a good obstruction theory}$\cX_d$ admits a nice
  obstruction theory in the sense of Definition~{\em \ref{df:obstruction}}.
\end{prop}
\begin{proof}
  Since~$\cX_d$ is limit preserving by Lemma~\ref{lem:X is limit
    preserving}, 
it follows from Theorem~\ref{thm: Herr complex is perfect} and 
  Lemmas~\ref{lem: obstruction class
    well defined and controls existence of lifts} and~\ref{lem: H1 of
    Herr complex of adjoint controls lifts} that the Herr
  complex~$\cC^\bullet(\ad M)$ provides the required obstruction
  theory.
\end{proof}

\section{Residual gerbes and isotrivial families}
\label{subsec:isotrivial}
In this section we briefly discuss the notion 
of isotrivial families of $(\varphi,\Gamma)$-modules
over reduced $\F$-schemes; i.e.\ of families
which are pointwise constant.   
The language of residual gerbes (see Appendix~\ref{app:residual gerbes}) provides
a convenient framework for doing this.
%
In the following discussion we allow $\Fnew$ to denote
any algebraic extension of $\F$. 

We have seen that $\cX_d$ is a quasi-separated Ind-algebraic stack,
which by Proposition~\ref{prop:X is an Ind-stack} may be written as the inductive limit of algebraic stacks 
with closed immersions (and so, in particular, monomorphisms) as
transition morphisms.  Thus the
discussion at the end of Appendix~\ref{app:residual gerbes} applies,
and in particular, for each point $x \in |\cX_d|,$
the residual gerbe $\cZ_x$ at $x$ in  $\cX_d$ exists.
Furthermore, if $x$ is a finite type point, then the canonical
monomorphism $\cZ_x \hookrightarrow \cX_d$ is an immersion.


If we consider an $\Fnew$-valued point 
$x: \Spec \Fnew \to \cX_{d}$
(by abuse of notation we use $x$ to denote both this point,
and its image $x \in |\cX_d|$, which is a finite type point of~$\cX_d$),
then we may base-change the residual gerbe to $\cX_{d}$ at $x$
(which is a gerbe over~$\F$) 
over $\Fnew$ via the composite $\Spec \Fnew \buildrel x \over
\longrightarrow \cX_{d} \to \Spec \F$; equivalently,
we may regard $x$ as a point of the base-change $(\cX_{d})_{\Fnew}$, 
and then consider the residual gerbe of $(\cX_{d})_{\Fnew}$ at
this point.
This residual gerbe is then of the form $[\Spec \Fnew/G]$,
for a finite type affine group scheme $G$ over $\Fnew$.
(By~\cite[\href{http://stacks.math.columbia.edu/tag/06QG}{Tag
  06QG}]{stacks-project} the group $G$ may be described as the fibre
product $x \times_{\cX_{d}} x$,
and hence its claimed properties follow from the
fact that by Proposition~\ref{prop:X is an Ind-stack}, the diagonal of $\cX_{d}$ is affine and of finite
type.) 
%
In fact, we have the following more precise result regarding $G$.

\begin{lemma}
	\label{lem:automorphisms are smooth}
	If $x$
        is an $\Fnew$-valued point of $(\cX_d)_{\Fnew}$,
	then $\Aut(x)$
	is an irreducible smooth closed algebraic subgroup
	of~$\GL_{d/\Fnew}$.
\end{lemma}
\begin{proof}As already noted, it follows from Proposition~\ref{prop:X
    is an Ind-stack}  that~$\Aut(x)$ is a finite type affine group
  scheme over~$\F'$.
	Let $D$ denote the \'etale $(\varphi,\Gamma)$-module over $\Fnew$
	corresponding to $x$,
	and let $R:= \End_{\A_{K,\F'},\varphi,\Gamma}(D)$.
	Then, if $\rhobar:G_K \to \GL_d(\Fnew)$
	denotes the Galois representation corresponding to $x$,
        we see that also $R = \End_{G_K}(\rhobar)$,
	and hence that $R$ is an $\Fnew$-subalgebra
	of $M_d(\Fnew)$. Furthermore, if $A$ is any finite
        type $\Fnew$-algebra,
	and if $D_A$ denotes the base-change of $D$ over~$A$,
	then the natural morphism
	$R\otimes_{\F'} A \to \End_{\A_{K,\F'},\varphi,\Gamma}(D_A)$
	is an isomorphism by Theorem~\ref{thm: Herr complex is
          perfect}~(2) and Lemma~\ref{lem: Yoneda for H0 H1 of Herr
          complex} (note that~$A$ is automatically flat over~$\F'$).
	Thus $$\Aut(x)(A) :=
	\Aut_{\A_{K,A},\varphi,\Gamma}(D_A) = (R\otimes_{\Fnew} A)^{\times}
	= (R\otimes_{\Fnew} A) \bigcap \GL_d(A),$$
	so that $\Aut(x)$, as a scheme over $\Fnew$,
	is precisely the open subscheme $R\cap \GL_d$, where
	we think of $R$ as an affine subspace of the affine space
	$M_d$.  In particular, we see that $\Aut(x)$ is an open
	subscheme of an affine space, and thus smooth and irreducible.
\end{proof}


Suppose that $S$ is a reduced $\Fnew$-scheme of finite type,
and that $S \to \cX_d$ is a morphism such that every closed point of
$S$ maps to some fixed $\F'$-valued point $x \in |\cX_d|$.  We can think of
this morphism as classifying a family of Galois representations
$\rhobar_S$ over $S$ whose fibre at each closed point is isomorphic
to the fixed representation $\rhobar: G_K \to \GL_d(\Fnew)$
classified by $x$.  If $G := \Aut_{G_K}(\rhobar)$,
then Lemma~\ref{lem:isotrivial} shows
that the morphism $S \to \cX$ factors through the residual gerbe
$[\Spec \Fnew/G],$ and thus corresponds to an {\em \'etale}
locally trivial $G$-bundle $E$ over $S$.
({\em A priori}, the $G$-bundle $E$ over $S$ is {\em fppf} locally
trivial.  However, since $G$ is smooth, by Lemma~\ref{lem:automorphisms
	are smooth}, we see that $E$ is in fact smooth locally trivial.
By taking an \'etale slice of a smooth cover
over which $E$ trivializes, we see that $E$ is in fact \'etale
locally trivial.)
We can then describe the family $\rhobar_S$ concretely
as a twist by $E$ of the constant family 
$S \times \rhobar,$ i.e.\ $\rhobar_S = E\times_G \rhobar$.


\section{Twisting families}\label{subsec: twists of
  families}We now begin our study of the dimensions of certain families of
$(\varphi,\Gamma)$-representations, and the behaviours of these
dimensions under the operations of twisting by families of
1-dimensional representations, and forming extensions of families.

We will need in particular to be able to twist families of representations by
unramified characters. Given an $\Fp$-algebra~$A$ and an element~$a\in
A^\times$, we have a $(\varphi,\Gamma)$-module $M_a$ whose underlying
$\A_{K,A}$-module is free of rank~$1$, generated by some~$v\in M_a$ for which
$\varphi(v)=av$ and $\gamma(v)=v$. 
If~$A=\F$
then the corresponding representation of~$G_K$ is the unramified
character~$\ur_a$ taking a geometric Frobenius to~$a$. The universal instance of this
construction comes by taking $a = x \in A = \F[x,x^{-1}]$;
the corresponding $(\varphi,\Gamma)$-module $M_x$ is then classified by 
a morphism which we denote ~$\ur_x:\Gm := \Spec \F[x,x^{\pm 1}] \to \cX_1$,
which evidently factors through~$\cX_{1,\red}.$

Given any morphism $T \to \cX_{d,\red},$ with $T$ a reduced finite
type $\F$-scheme, corresponding to a family $\rhobar_{T}$ of
$G_K$-representations over $T$, we may consider the 
family ~$\rhobar_T \boxtimes\ur_x$ over $T\times\Gm$, as in Section~\ref{subsec: tensor product
  of phi gamma and duality}. We refer to this operation on
$(\varphi,\Gamma)$-modules as \index{unramified twisting} \emph{unramified
          twisting}.  

\begin{df}
	\label{def:twistable} \index{twistable} \index{essentially twistable}
We say that~$\rhobar_T$ is \emph{twistable} if whenever
$\rhobar_{t}\cong\rhobar_{t'}\otimes\ur_a$ where $t,t'\in T(\Fpbar)$
and~$a\in\Fpbartimes$, then~$a=1$. We say that it is
\emph{essentially twistable} if for each~$t\in T(\Fpbar)$,
 the set of~$a\in\Fpbartimes$ for which there exists~$t'\in
T(\Fpbar)$  with $\rhobar_{t}\cong\rhobar_{t'}\otimes\ur_a$ is finite.
\end{df}

\begin{lem}\label{lem: dimension of Gm twist}
  If the dimension of the scheme-theoretic image of~$T$
  in~$\cX_{d,\red}$ is~$e$, then the dimension of the scheme-theoretic
  image of~$T\times\Gm$ in~$\cX_{d,\red}$ is at most~$e+1$. If~$T$ 
   contains a dense open subscheme~$U$ such
  that~$\rhobar_U$ is essentially twistable, then equality holds.
\end{lem}
\begin{rem}
  \label{rem: twistable is enough}Since a twistable
  representation is essentially twistable, we see that if~$T$ 
   contains a dense open subscheme~$U$ such
  that~$\rhobar_U$ is twistable, then equality holds in Lemma~\ref{lem: dimension of Gm twist}.
\end{rem}
\begin{proof}[Proof of Lemma~\ref{lem: dimension of Gm twist}]
  We may assume that~$T$ is irreducible. Write~$f$ for the morphism
  $T\to\cX_{d,\red}$ and~$g$ for the morphism
  $T\times\Gm\to\cX_{d,\red}$.  By~\cite[\href{https://stacks.math.columbia.edu/tag/0DS4}{Tag 0DS4}]{stacks-project}, we
  may, after possibly replacing~$T$ by a nonempty open subscheme,
  assume that for each $\Fpbar$-point $t$ of~$T$ we have
  $\dim T_{f(t)}=\dim T-e$.

  Let~$v=(t,\lambda)$ be an $\Fpbar$-point of~$T\times\Gm$. Then
  $(T\times \Gm)_{g(v)}$ contains~$T_{f(t)}\times\{\lambda\}$, 
so that
  $\dim (T\times \Gm)_{g(v)}\ge \dim T_{f(t)}=\dim T-e=\dim (T\times
  \Gm)-(e+1)$, so the first claim follows from another application
  of~\cite[\href{https://stacks.math.columbia.edu/tag/0DS4}{Tag 0DS4}]{stacks-project}.

  For the last part, we may replace~$T$ by~$U$, and then the
  hypothesis that~$\rhobar_T$ is essentially twistable means that we can identify
  $(T\times \Gm)_{g(v)}$ with a finite union of fibres $T_{f(t')}$
  (indexed by the finitely many~$a$ for which there is a ~$t'$ and an isomorphism $\rhobar_{t}\cong\rhobar_{t'}\otimes\ur_a$), so that equality holds in
  the above inequality, and the result again follows
  from~\cite[\href{https://stacks.math.columbia.edu/tag/0DS4}{Tag 0DS4}]{stacks-project}.
\end{proof}

Recall that for each embedding~$\sigmabar:k\into\F$, we have a
corresponding fundamental character
$\omega_{\sigmabar}:I_K\to\F^\times$, 
which
corresponds via local class field theory (normalised to take
uniformizers to geometric Frobeneii) to the composite of~$\sigmabar$
and the natural map $\cO_K^\times\to k^\times$. It is well-known,
and easy to check,
that
$\prod_{\sigmabar}\omega_{\sigmabar}^{e(K/\Qp)}=\epsilonbar|_{I_K}$. If~$\Fnew$
is an algebraic extension of~$\F$, then we
can and do identify the embeddings $k\into\F$ and $k\into\Fnew$.

Let~$\underline{n}=(n_{\sigmabar})_{\sigmabar:k\into\F}$ be a tuple of integers $0\le
n_{\sigmabar}\le p-1$. 
The characters~$\omega_{\sigmabar}$ can all be extended to~$G_K$, and the restriction
to~$I_K$ of any character $G_K\to\Fpbartimes$ is equal
to~$\prod_{\sigmabar}\omega_{\sigmabar}^{-n_{\sigmabar}}$ for 
some~$\underline{n}$ (which is unique unless the character is unramified). We write~$\psi_{\underline{n}}:G_K\to\F^\times$ for a
fixed choice of an extension
of~$\prod_{\sigmabar}\omega_{\sigmabar}^{-n_{\sigmabar}}$
to~$G_K$. If~$n_{\sigmabar}=0$ for all~$\sigmabar$, or
$n_{\sigmabar}=p-1$ for all~$\sigmabar$ (these are exactly
the cases corresponding to unramified characters), then we fix the choice
$\psi_{\underline{n}}=1$. We can and do make  our choice so that
if~$\underline{n}$, $\underline{n}'$ are such that
$(\psi_{\underline{n}}\psi_{\underline{n}'}^{-1})|_{I_K}=\epsilonbar|_{I_K}$,
then in fact $\psi_{\underline{n}}\psi_{\underline{n}'}^{-1}=\epsilonbar$.

Somewhat abusively, we will also
write~$\psi_{\underline{n}}$ for the constant family of
$(\varphi,\Gamma)$-modules $\mathbb{D}(\psi_{\underline{n}})$ over any
$\Fpbar$-scheme~$T$.

 \begin{rem}\label{rem: irreducible representations finite up
  to twist} The irreducible representations
          $\alphabar: G_K \to \GL_d(\Fpbar)$ are easily classified;
          indeed, since the wild inertia subgroup must act trivially,
          they are tamely ramified, so that the restriction
          of~$\alphabar$ to~$I_K$ is diagonalizable. Considering the
          action of Frobenius on tame inertia, it is then easy to see
          that~$\alphabar$ is induced from a character of the
          unramified extension of~$K$ of degree~$d$. It follows that
          each irreducible~$\alphabar$ is absolutely
          irreducible.
          
          As we have just seen, there are only finitely many
          such characters up to unramified twist, so
          that 
          in particular there are only finitely many such~$\alphabar$
          up to unramified twist. It follows that there are up to unramified twist
          only finitely many irreducible $(\varphi,\Gamma)$-modules
          over~$\Fpbar$ of any fixed rank. 
        \end{rem}

It is convenient to extend the notation~$\rho_T$ to algebraic stacks.
To this end, if~$\cT$ is a reduced algebraic stack of finite type
over~$\Fpbar$, we will denote
a morphism
$\cT\to(\cX_{d,\red})_{\Fpbar}$ 
by~$\rhobar_{\cT}$, and will sometimes abusively refer
to $\rho_{\cT}$ as a family of $G_K$-representations parameterized by~$\cT$.

\begin{df}
	\label{def:maximally nonsplit}We say that a
representation~$\rhobar:G_K\to\GL_d(\Fpbar)$ is \emph{maximally
  nonsplit of niveau~$1$} \index{maximally
  nonsplit of niveau~$1$} if it has a unique filtration by~$G_K$-stable
$\Fpbar$-subspaces such that all of the graded pieces are one-dimensional
representations of~$G_K$. We say that the family~$\rhobar_{\cT}$ is
\emph{maximally nonsplit of niveau~$1$} if~$\rhobar_t$ is maximally
nonsplit of niveau~$1$ for
all~$t\in T(\Fpbar)$, and that~$\rhobar_T$ is
\emph{generically maximally nonsplit of niveau~$1$} if there is a
dense open substack~$\cU$ \index{generically maximally nonsplit of niveau~$1$} 
of~$\cT$ such that~$\rhobar_{\cU}$ is maximally nonsplit of niveau~$1$.

 In particular, we say that an algebraic substack~$\cT$ of~$(\cX_{d,\red})_\Fpbar$ is
\emph{maximally nonsplit of niveau~$1$} if~$\rhobar_t$ is maximally
nonsplit of niveau~$1$ for
all~$t\in \cT(\Fpbar)$, and that~$\cT$ is
\emph{generically maximally nonsplit of niveau~$1$} if there is a dense open substack~$\cU$
of~$\cT$ such that~$\cU$ is maximally nonsplit of niveau~$1$.
\end{df}
\begin{rem}
  \label{rem: everything is maximally nonsplit niveau 1}We will
  eventually see in Theorem~\ref{thm:reduced dimension} below that 
~$(\cX_{d,\red})_{\Fpbar}$ itself is an algebraic stack of finite type
over~$\Fpbar$, and is generically maximally nonsplit of
  niveau 1.
\end{rem}

Our next goal is to prove a structure theorem (Proposition~\ref{prop: structure of maximally
    nonsplit}) for families which are
generically maximally nonsplit of niveau~$1$. We begin with the
following technical lemma.

\begin{lemma}
\label{lem:Fitting fibres}
Let  $T$ be a reduced scheme of finite type over an algebraically closed field~$k$,
let $M$ be a coherent sheaf on $T\times_k \Gm$, and for each $t \in T(k)$,
let $M_t$ denote the restriction of $M$ to
$\Gm \iso t\times_k \Gm \hookrightarrow T \times_k \Gm$.
Suppose that for each point $t \in T(k)$, $M_t$ is a length one skyscraper sheaf supported at a  single point
of~$\Gm$.  Then there is a dense open subscheme $U$ of
$T$ such that, over  $U$, the composite $\Supp(M) \hookrightarrow   T\times_k \Gm
\to  T$  pulls back to  an isomorphism,
while the pullback of $M$ over $U\times_k \Gm$ is locally free of rank one over
its support.
\end{lemma}
\begin{proof}
Let $\Fitt(M)\subseteq \cO_{T\times \Gm}$  denote the Fitting ideal sheaf of $M$,
and write $Z := \Spec \cO_{T\times \Gm}/\Fitt(M);$ recall that  $\Supp(M)$ 
is a closed subscheme of~$Z$, supported on  the same underlying  closed subset
of~$T\times_k \Gm.$
For each $t  \in T(k)$, we find that the fibre $Z_t$  of $Z$  over $t$
is equal to
$\Spec  \cO_{\Gm}/\Fitt(M_t)$ (recall that the formation
of Fitting  ideals is compatible  with  base-change), which, by assumption,
is a single reduced point of $\Gm$.
Thus the composite $Z \to T\times_k \Gm \to T$  is a  morphism with
reduced singleton fibres,
and hence we may find a dense open subscheme $U$ of
$T$ such that this morphism restricts to an isomorphism over~$U$.
In particular, the pull-back of $Z$ over $U$ is reduced (since $U$ is),
and thus coincides with the pull-back of $\Supp (M)$  over $U$.
This shows that the morphism $\Supp (M) \to T$  pulls back to an isomorphism
over~$U$.

If we use the  inverse of the isomorphism just constructed
to identify $U$ with an open subscheme of $\Supp(M)$,
then our assumption on the nature of $M_t$  for  $t  \in T(k)$ implies
that the fibre of $M$ over each $k$-point of $U$ is one-dimensional.  Since
$U$  is reduced, we find that the restriction of $M$ to $U$ is  furthermore locally
free of rank one, as claimed.
\end{proof}

\begin{prop}\label{prop: structure of maximally
    nonsplit}If~$\cT$ is a reduced finite type algebraic stack over~$\Fpbar$, and~$\rhobar_{\cT}$ is generically maximally
  nonsplit of niveau~$1$, then there exist:
  \begin{itemize}
  \item
    a dense 
open substack~$\cU$ of~$\cT$;
  \item
tuples $\underline{n}_i$,    $1\le i\le d$;
  \item morphisms~$\lambda_i:\cU\to\Gm$, $1\le i\le d$;
  \item morphisms~$x_i:\cU\to(\cX_{i,\red})_{\Fpbar}$, $0\le i\le d$, corresponding to
    families of $i$-dimensional representations~$\rho_i$ over~$\cU$;
  \end{itemize}
such that $x_d$ is the restriction of~$\rhobar_{\cT}$ to~$\cU$, and for
each~$0\le i\le d-1$, we have short exact
sequences of families over~$\cU$ \[0\to\rho_i\to\rho_{i+1}\to
  \ur_{\lambda_{i+1}}\otimes\psi_{\underline{n}_{i+1}}\to 0.\] 
In particular, for each~$t\in \cU$, we
have \numequation\label{eqn:maximally nonsplit structure}\rhobar_t\cong \begin{pmatrix}
      \ur_{\lambda_1(t)}\otimes\psi_{\underline{n}_1} &*&\dots &*\\
      0& \ur_{\lambda_2(t)}\otimes\psi_{\underline{n}_2}&\dots &*\\
      \vdots&& \ddots &\vdots\\
      0&\dots&0& \ur_{\lambda_d(t)}\otimes\psi_{\underline{n}_d}\\
    \end{pmatrix}.  \end{equation}%
  
\end{prop}
\begin{rem}
  \label{rem: maximally nonsplit doesn't have a universal family}
  We cannot necessarily attain~(\ref{eqn:maximally nonsplit
    structure}) over all of~$\cT$. For example, there are families of two-dimensional representations of~$G_K$ which
  are generically maximally nonsplit of niveau one, but which also specialise to irreducible
  representations. (The existence of such families follows
  from Theorem~\ref{thm:reduced dimension} below.)
\end{rem}

\begin{proof}[Proof of Proposition~\ref{prop: structure of maximally
    nonsplit}]
We begin by assuming that~$\cT$ is a scheme, say~$\cT=T$.  
Without loss of generality, we can replace by $T$ by a dense open subscheme,
and thus assume that~$\rhobar_{T}$ is maximally nonsplit of
  niveau~$1$, i.e.\ that~$\rhobar_t$ is maximally nonsplit of
  niveau~$1$ for each~$t\in T(\Fpbar)$, and also assume
that $T$ is a disjoint union of finitely many integral open and closed subschemes.  Replacing
$T$ by each of these subschemes in turn, we may assume
that $T$ itself is integral. 

For each $\underline{n}$, we may form the family
  $$\rho_T \boxtimes (\ur_x^{-1} \otimes \psi_{\underline{n}}^{-1})(1)$$ 
  over 
  $T\times_{\Fpbar} \Gm$ (where as above $x$ denotes the variable on 
  $\Gm$, so that $\Gm := \Spec \Fbar_p[x,x^{-1}]$),
  and then consider (in the notation of Remark~\ref{rem: writing Herr complex as Galois cohomology})
  \[H_{\underline{n}} := H^2\bigl(G_K, \rho_T\boxtimes (\ur_x^{-1} \otimes \psi_{\underline{n}}^{-1})(1)\bigr),\]
  a coherent sheaf on $T\times_{\Fbar_p} \Gm$. 
  Since the formation of $H^2(G_K,\text{--})$ is compatible
  with  arbitrary finite type base-change, we see that  for any point
$ t \in T(\Fbar_p)$, the pull-back of
$H_{\underline{n}}$ over $t \times_{\Fbar_p} \Gm \cong \Gm$
admits the description
\[(H_{\underline{n}})_t\cong H^2\bigl(G_K, \rho_t\otimes (\ur_{x}^{-1} \otimes \psi_{\underline{n}}^{-1})(1)\bigr),\]
and that the fibre $(H_{\underline{n}})_{(t,x_0)}$ of this pull-back
at a point $x_0 \in \Gm(\Fbar_p)$,
admits the description
  $$(H_{\underline{n}})_{(t,x_0)} \cong 
  H^2\bigl(G_K, \rho_t\otimes (\ur_{x_0}^{-1} \otimes \psi_{\underline{n}}^{-1})(1)\bigr)
	  \cong \Hom_{G_K}(\rho_t, \ur_{x_0}\otimes \psi_{\underline{n}})
	  $$ 
(the second isomorphism following
from Tate local duality).
	  The assumption that $\rho_t$ is maximally non-split
	  shows that this fibre is non-zero 
	  for exactly one choice of $\underline{n}$ and one choice of $x_0$,
	  and is then precisely one-dimensional.
In particular, we see, for this distinguished choice 
of~$\underline{n}$, that $(H_{\underline{n}})_t$ is set-theoretically supported at~$x_0$.
          Now let $\xi_0:\Spec  \Fbar_p[\epsilon]/(\epsilon^2) \hookrightarrow
\Gm$ be the (unique up to scalar) non-trivial tangent vector to $x_0$, 
and let
$(H_{\underline{n}})_{(t,\xi_0)}$
denote the pull-back of $(H_{\underline{n}})_t$
over $\xi_0$.  A similar computation using Tate local duality and the maximal non-splitness
of $\rho_t$ shows that 
$(H_{\underline{n}})_{(t,\xi_0)}$
is also one-dimensional; thus $(H_{\underline{n}})_t$  is in fact a  skyscraper sheaf of
length one supported at~$x_0$.
For all other choices of~$\underline{n}$,
we see that $(H_{\underline{n}})_t$ vanishes (since all its fibres do).

          Applying Lemma~\ref{lem:Fitting fibres} to
          $M  := \bigoplus_{\underline{n}} H_{\underline{n}},$ 
          and replacing $T$ by an appropriately chosen dense open subscheme,
          we find that $\Supp(M)$ maps isomorphically to~$T$,
          and that $M$ is locally free of rank one over its support.
          Taking into account the definition of $M$ as a direct sum,
          and the irreducibility of~$T$,
          we find that in fact $M = H_{\underline{n}}$  for exactly one choice
          $\underline{n}_d$ of~$\underline{n}$,
          and that $H_{\underline{n}} = 0$ for all other
          possible choices. 
	  Let $\lambda: T \to \Gm$ denote the composite of the inverse
          isomorphism $T \iso \Supp(M)$
	  with the projection $\Supp(M) \to \Gm$.
	  Again using the fact that the formation of $H^2$ commutes
	  with base-change, we find that the pull-back $M_T$
	  of $M$ to $T \cong \Supp(M)$ may be identified with
	  $$H^2\bigl( \rho_T \boxtimes (\ur_{\lambda}^{-1}\otimes \psi_{\underline{n}_d}^{-1})(1)\bigr).$$
          We have already remarked that $M_T$ is locally free of rank one;
	  replacing $T$ by a non-empty open subscheme once more, we may
	  in fact assume that $M_T$ is free of rank one, and
	  so choose a nowhere zero section of~ $M_T$, which by Lemma~\ref{lem: Tate duality for H2 to H0 in families} we 
	  may and do regard as a surjection 
	  $$ \rho_T \onto \ur_{\lambda}\otimes
          \psi_{\underline{n}_d}.$$ The kernel of this surjection is a
          rank~$(d-1)$ family of projective \'etale
          $(\varphi,\Gamma)$-modules, so the result for~$T$ follows by
          induction on~$d$.

          We now return to the general case that~$\cT$ is a reduced
          finite type algebraic stack over~$\Fpbar$. Since~$\cX_d$ is
          Ind-algebraic by Proposition~\ref{prop:X is an Ind-stack},
          we can form the scheme-theoretic image~$\cT'$ of the
          morphism $\cT\to(\cX_{d,\red})_{\Fpbar}$, which is again a reduced finite
          type algebraic stack over~$\Fpbar$. The
          family~$\rhobar_{\cT'}$ is again generically maximally
          nonsplit of niveau~$1$ (since by the stacky version of
          Chevalley's theorem (see~\cite[App.\ D]{MR2818725}),
          the image of $\cT$ in $\cT'$ is constructible, and so contains
          a dense open subset of $\cT'$), and since the morphism $\cT\to\cT'$
          is scheme-theoretically dominant, any dense open substack
          of~$\cT'$ pulls back to a dense open substack of~$\cT$. It
          therefore suffices to prove the result
          for~$\cT'$. Replacing~$\cT$ by~$\cT'$, then, we can and do
          assume that~$\cT$ is a closed substack of~$(\cX_{d,\red})_{\Fpbar}$.

          Let $T\to\cT$ be a smooth cover by a (necessarily) reduced
          and finite type $\Fpbar$-scheme (we can ensure that $T$ 
is finite type over $\Fbar_p$, since $\cT$ is so), and let~$\rhobar_T$ be the
          induced family. We wish to deduce the result for~$\cT$ from
          the result for~$T$, which we have already proved.  We thus
replace $T$ by a dense open subscheme (and correspondingly replace $\cT$ by
a dense open substack, namely the image of this open subscheme of $T$),
          so that the family $\rhobar_T$ is maximally non-split,
and so that the morphisms $\lambda_i:T\to\Gm$
          and ~$x_i:T\to(\cX_{i,\red})_{\Fpbar}$ exist.  We will then
show that these morphisms factor through the smooth
          surjection~$T\to\cT$. This amounts to showing that we can
          factor the morphism  of groupoids
\numequation
\label{eqn:groupoid to product}
T\times_{\cT}T\to T\times T
\end{equation}
through a morphism of groupoids $T\times_{\cT} T \to T\times_{\Gm}T$ or
$T\times_{\cT} T \to T\times_{(\cX_{i,\red})_{\Fpbar}}T$ respectively.

Since $\cT$ is a substack of $\cX$, the natural morphism induces an isomorphism
$T\times_{\cT} T \iso T\times_{\cX} T$.  Thus for any test scheme $T'$,
a $T'$-valued point of $T\times_{\cT} T$ consists of a pair of morphisms
$f_0,f_1: T' \rightrightarrows T$, and an isomorphism of families
\numequation
\label{eqn:family isomorphism}
f_0^*\rho_{T} \iso f_1^*\rho_T
\end{equation}
(where $f_i^* \rho_T$ denotes the pull-back of $\rho_T$ 
to $T'$ via~$f_i$).

In the first case, in which we consider one of the morphisms $\lambda_i$,
we need to show (employing the notation just introduced)
 that the morphism $f_0 \times f_1: T'\times T' \to T\times T$ necessarily factors through
$T \times_{\G_m} T$ (this latter fibre product formed using the morphisms $\lambda_i$).
More concretely, this amounts to checking that $\lambda_i \circ f_0 = \lambda_i \circ f_1$.
Since the source of~\eqref{eqn:groupoid to product} is reduced and of finite type
over $\Fbar_p$, it suffices to check this for reduced test schemes
$T'$ of finite type over~$\Fbar_p$, and thus to check that these morphisms coincide
at the $\Fbar_p$-valued points $t'$ of $T'$.   This follows from
the assumed isomorphism~\eqref{eqn:family isomorphism}, and the uniqueness 
of the filtration at each $\Fbar_p$-valued point of the family.


          In the second case,
in which we consider one of the morphisms $x_i$,
the morphism
          $T\times_{(\cX_{i,\red})_{\Fpbar}}T\to T\times T$ is no longer a monomorphism,
and so we have to actually construct a morphism of groupoids
          $T\times_{\cT}T\to T\times_{(\cX_{i,\red})_{\Fpbar}}T$
(rather than simply verify a factorization).
%
%
%
%
So, given $f_0, f_1: T' \to T,$ and an isomorphism~\eqref{eqn:family isomorphism},
we have to construct a corresponding isomorphism
\numequation
\label{eqn:desired isomorphism}
f_0^*(\rho_T)_i \iso f_1^*(\rho_T)_i,
\end{equation}
where $(\rho_T)_i$ denotes the family of $i$-dimensional subrepresentations
of $\rho_T$ corresponding to the morphism~$x_i$.
We will define the desired isomorphism~\eqref{eqn:desired isomorphism}
simply to be the restriction of
the isomorphism~\eqref{eqn:family isomorphism}.
For this definition to make sense, we have to show that~\eqref{eqn:family isomorphism}
does in fact restrict to an isomorphism~\eqref{eqn:desired isomorphism}.
In fact, it suffices to check that~\eqref{eqn:family isomorphism}
induces an embedding 
$f_0^*(\rho_T)_i \hookrightarrow f_1^*(\rho_T)_i;$
reversing the roles of $f_0$ and $f_1$ then allows us to promote this
to an isomorphism.

Again, it suffices to check this on reduced test schemes $T'$ of finite type
over~$\Fbar_p$, and then it suffices to show that the composite
$$f_0^*(\rho_T)_i \hookrightarrow f_0^*\rho_T \buildrel
\text{\eqref{eqn:family isomorphism}} \over \longrightarrow
f_1^*\rho_T \to f_1^*\rho_T/ f_1^*(\rho_T)_i$$ 
vanishes.  This can be checked on $\Fbar_p$-valued points,
where it again follows from the uniqueness of the filtrations that define
the~$(\rho_T)_i$.
%
\end{proof}


\section{Dimensions of families of extensions}\label{subsec:
  dimensions of families of extensions}
Suppose that we have a morphism $T \to (\cX_{d,\red})_{\Fpbar},$
with $T$ a reduced finite type $\Fpbar$-scheme,
corresponding to a family $\rhobar_{T}$ of $G_K$-representations over 
$T$.  
Fix a representation $\alphabar: G_K \to \GL_a(\Fpbar)$.

\begin{lemma}\label{lem: Ext2 locally free} If $\Ext^2_{G_K}(\alphabar,\rhobar_t)$ is of constant rank $r$ 
	for all $t \in T(\Fbar_p)$, then $H^2(G_K,\rhobar_T\otimes \alphabar^{\vee})$ is
	locally free of rank $r$ as an $\cO_T$-module.
\end{lemma}
\begin{proof}
	Since the formation of $H^2$
	is compatible with arbitrary finite type base-change (by Theorem~\ref{thm: Herr complex is perfect}~(2)),
	we see that $H^2(G_K,\rhobar_T\otimes\alphabar^{\vee})$
	is a coherent sheaf on $T$ of constant
	fibre rank $r$.  Since $T$ is reduced, the lemma follows.
\end{proof}

By Corollary~\ref{cor:Herr complex length},
if $T$ is affine, 
then we can and do choose a good complex (that is, a bounded
complex of finite rank locally free $\cO_T$-modules) 
$$C^0_T \to C^1_T \to C^2_T$$ 
computing $H^{\bullet}(G_K,\rhobar_T\otimes\alphabar^{\vee})$. Suppose that we are in the context of the preceding lemma,
i.e.\ that $\Ext^2_{G_K}(\alphabar,\rhobar_T)$ is locally free of some rank $r$. It
follows that the truncated complex
$$C^0_T \to Z^1_T$$ is again good (here $Z^1_T:=\ker(C^1_T\to
C^2_T)$). As in Remark~\ref{rem: writing Herr complex as Galois
  cohomology}, we write \[H^1(G_K,\rhobar_T\otimes\alphabar^\vee)\] for
the cohomology group of this complex in degree~$1$; its formation is
compatible with arbitrary finite type base-change, and its specialisations to
finite type points of~$T$ agree with the usual Galois cohomology.

In particular, if we choose another integer $n$,
we may use the surjection $Z^1_T \to
H^1(G_K,\rhobar_T\otimes\alphabar^{\vee})$ (together with Lemma~\ref{lem: Yoneda for H0 H1 of Herr complex})
to construct a universal family of extensions 
\numequation\label{eqn: universal family of extensions}0 \to \rhobar_T \to \cE \to \alphabar^{\oplus n} \to 0\end{equation}
parameterized by the vector bundle $V$ corresponding to the finite rank
locally free sheaf $(Z_T^1)^{\oplus n}$,
giving rise to a morphism
\numequation
\label{eqn:induced morphism}
V \to (\cX_{d+ a n,\red})_{\Fpbar}.
\end{equation}(See~\cite[\S 4.2]{CEGSKisinwithdd} for a similar construction.)

Write $G_{\alphabar} := \Aut_{G_K}(\alphabar),$
thought of as an affine algebraic group over $\Fpbar$.

\begin{prop}Maintain the notation and assumptions above; so we assume in particular
  that $\Ext^2_{G_K}(\alphabar,\rhobar_t)$ is of constant rank $r$ 
	for all $t \in T(\Fbar_p)$.
	\label{prop:dimension control}\leavevmode
        \begin{enumerate}
        \item If $e$ denotes the dimension of the scheme-theoretic
          image 
          of $T$ in $(\cX_{d,\red})_{\Fpbar},$ then the dimension of the
          scheme-theoretic image,
          with respect to the morphism~{\em (\ref{eqn:induced
              morphism})}, 
          of $V$ in $(\cX_{d+an,\red})_{\Fpbar}$ is bounded above by
          \[e + n([K:\Q_p] ad +r) -n^2 \dim G_{\alphabar}.\]

        \item  Suppose that:
          \begin{enumerate}
          \item $n=1$;
          \item  $\alphabar$ is one-dimensional;
            
          \item $\rhobar_T$ is generically maximally nonsplit of
            niveau~$1$;
          \item\label{item: weird condition over Qp} if~$K=\Qp$, then
after replacing $T$ by a dense open subscheme,
so that the tuples $\underline{n}_i$ exist,
and so that the morphisms $\lambda_i$ and $x_i$ as in
the statement of Proposition~{\em \ref{prop: structure of maximally nonsplit}}
are defined on~$T$, 
at least one of the following conditions holds at each~$t~\in~T(\Fbar_p)$:
  \begin{enumerate}
  \item
    $\ur_{\lambda_{d-1}(t)}\otimes\psi_{\underline{n}_{d-1}}\ne\alphabar(1)$.
  \item
    $\ur_{\lambda_{d}(t)}\otimes\psi_{\underline{n}_{d}}=\alphabar$.
  \item $\ur_{\lambda_{d}(t)}\otimes\psi_{\underline{n}_{d}}=\alphabar(1)$.
  \end{enumerate}
        \end{enumerate}
Then equality
          holds in the inequality of~(1). Furthermore the family of
          extensions~$\cE_V$ corresponding to~$V$ is generically
          maximally nonsplit of niveau~$1$.
\item Suppose that we are in the setting of~(2), and that 
  after replacing $T$ by a dense open subscheme, the following conditions hold:
\begin{enumerate}
\item
For each
  $t~\in~T(\Fbar_p)$ the representation $\rhobar_t$ is maximally nonsplit,
and furthermore the unique quotient character of~$\rhobar_t$
  {\em(}which exists by the assumption that~$\rhobar_T$ is generically
  maximally nonsplit of niveau~$1${\em)} is equal
  to~$\alphabar(1)$.
\item
Furthermore, if we write $$0 \to \rbar_t \to
  \rhobar_t \to \alphabar(1) \to 0$$ for the corresponding filtration of $\rhobar_t$,
and write~$\gamma_t$ for the unique quotient character of~$\rbar_t$ 
then for each $t\in T(\Fpbar)$ either
  \begin{enumerate}
  \item we have~$\gamma_t\ne \alphabar(1)$, or
\item the cyclotomic character~$\epsilonbar$ is trivial, we have
  $\gamma_t=\alphabar(1)=\alphabar$, and the extension of~$\alphabar$
  by~$\gamma_t=\alphabar(1)$ induced by~$\rhobar_t$ is tr\`es ramifi\'ee.
  \end{enumerate}
\end{enumerate}
  Then, after replacing~$V$ by a dense open subscheme, 
we have that for all
  $t\in V(\Fpbar)$ the extension of $\alphabar$ by~$\alphabar(1)$
  induced by the extension of~$\alphabar$ by $\rhobar_t$ is tr\`es
  ramifi\'ee.
\end{enumerate}
\end{prop}
\begin{rem}
  \label{rem: see below for weird condition}See Remark~\ref{rem: further explanation of weird condition} below for a
  further discussion of condition~(\ref{item: weird condition over Qp}).
\end{rem}

\begin{proof}[Proof of Proposition~\ref{prop:dimension control}]
	Writing $T$ as the union of its irreducible components,
	we obtain a corresponding decomposition of $V$ into the
	union of its irreducible components, and we can prove the
	theorem one component at a time.  Thus we may assume
	that $T$ (and hence $V$) is irreducible.
        Replacing $T$ by a non-empty affine open subscheme, we may furthermore
        assume that it is affine, so that we may find a good complex $C^{\bullet}_T$
        supported in degrees $[0,2]$ that is
        isomorphic in the derived category to the Herr complex, as above.

	Let $f: T \to (\cX_{d,\red})_{\Fpbar}$ be the given morphism, and let
	$g: V \to (\cX_{d+ a n,\red})_{\Fpbar}$ denote the
	morphism~(\ref{eqn:induced morphism}).
        In order to relate the scheme-theoretic
	image of $g$ to the scheme-theoretic image of $f$,
	we will also consider (in a slightly indirect manner)
	the scheme-theoretic image of the morphism
	$V \to (\cX_{d,\red})_{\Fpbar} \times(\cX_{d+ a n,\red})_{\Fpbar}$
	(the first factor being the composite
	of $f$ and the projection from $V$ to $T$, and the second factor
	being $g$).  If $(\cX_{d,\red})_{\Fpbar}$ were a scheme or an algebraic space,
	then the dimension of this latter scheme-theoretic image
	would be greater than or equal to the dimension of the former;
	but since $(\cX_{d,\red})_{\Fpbar}$ is a stack, we have to be slightly careful
	in comparing the dimensions of these two images.


	Let $v:\Spec \Fbar_p \to V$ be an $\Fbar_p$-point, and let $t$ 
	denote the composite $\Spec \Fbar_p \to V \to T$, so that
	$t$ is an $\Fbar_p$-point of $T$.  We also write $f(t)$ for
	the composite $f\circ t$, and $g(v)$ for the composite $g\circ v$.
	We write $T_{f(t)}$ for the fibre product
	$T\times_{(\cX_{d,\red})_{\Fpbar}} \Spec \Fbar_p$ (with respect
        to the morphisms $f: T \to (\cX_{d,\red})_{\Fpbar}$
	and $f(t): \Spec \Fbar_p \to (\cX_{d,\red})_{\Fpbar}$), and write $V_{g(v)}$ 
	and $V_{(f(t),g(v))}$ for the evident analogous fibre
	products. 

	We begin with~(1). In order to prove the required bound on the dimension of the
	scheme-theoretic image of $g$, it suffices (e.g.\
	by~\cite[\href{https://stacks.math.columbia.edu/tag/0DS4}{Tag 0DS4}]{stacks-project}) to show
	that 
\[\dim V_{g(v)} \buildrel ? \over \geq
	\dim V - e - n( [K:\Q_p] ad +r) + n^2 \dim G_{\alphabar}\]
	for $v$ lying in some non-empty open subscheme of $V$.
	As already alluded to above, what we will actually be able to
	do is to estimate the dimension of the fibre $V_{(f(t),g(v))},$
	and so our first job is to compare the dimension of this fibre
	with that of $V_{g(v)}$.  

	Let $\rhobar_{f(t)}$ denote the Galois
	representation corresponding to $f(t)$, and
	let $G_t  := \Aut(\rhobar_{f(t)}).$  Then,
        by the discussion of Section~\ref{subsec:isotrivial},
        the morphism
	$t: \Spec \Fbar_p \to (\cX_{d,\red})_{\Fpbar}$ induces an immersion
	$$[\Spec \Fbar_p/G_t] \hookrightarrow (\cX_{d,\red})_{\Fpbar},$$
	which in turn induces a monomorphism
\numequation
\label{eqn:residual gerbe}
	[\Spec \Fbar_p/G_t]\times_{(\cX_{d,\red})_{\Fpbar}} V_{g(v)} \hookrightarrow V_{g(v)}.
\end{equation}
	Since $\Spec \Fbar_p$ is a $G_t$-torsor over $[\Spec \Fbar_p/G_t],$
	we see that $V_{(f(t),g(v))}$ is a $G_t$-torsor over
	$[\Spec \Fbar_p/G_t]\times_{(\cX_{d,\red})_{\Fpbar}} V_{g(v)},$
	and hence that
		\numequation\label{eqn: first dimension inequality}\dim V_{(f(t),g(v))}  \leq \dim G_t + \dim V_{g(v)}.\end{equation}
	Thus it suffices to prove the inequality
	$$\dim V_{(f(t),g(v))} \buildrel ? \over \geq
	\dim V - e - n( [K:\Q_p] ad +r) + n^2\dim G_{\alphabar} + \dim G_t.$$
	For $t$ lying in some non-empty open subscheme of $T$,
	we have that $e = \dim T - \dim T_{f(t)},$
	and so, replacing $T$ by this non-empty open subscheme 
	and $V$ by its preimage,
	we have to show that
	$$\dim V_{(f(t),g(v))} \buildrel ? \over \geq
	\dim V - \dim T + \dim T_{f(t)} - n( [K:\Q_p] ad +r) + n^2\dim G_{\alphabar} + \dim G_t$$
	for $v$ lying in some non-empty open subscheme of $V$.
	In fact we will show that
	this inequality holds for every $\Fbar_p$-point
	$v$ of $V$.

	The right hand side of our putative inequality can be rewritten
	as
	$$n \rk Z_T^1 + \dim T_{f(t)} - n( [K:\Q_p] ad + r) +
	n^2\dim G_{\alphabar} + \dim G_t,$$
	while the local Euler characteristic formula
	shows that
	$$H^0(G_K,\rhobar_{f(t)}\otimes\alphabar^{\vee}) - H^1(G_K,\rhobar_{f(t)}\otimes\alphabar^{\vee}) + r = 
	- [K:\Q_p] ad.$$
	Thus we may rewrite our desired inequality as 
	\begin{multline*}
	\dim V_{(f(t),g(v))}
	\buildrel ? \over \geq
	n\bigl(\rk Z_T^1 - \dim H^1(G_K, \rhobar_{f(t)}\otimes \alphabar^{\vee})\bigr) \\ +
	n \dim H^0(G_K, \rhobar_{f(t)}\otimes \alphabar^{\vee}) + \dim T_{f(t)} + n^2 \dim G_{\alphabar} + \dim G_t.
\end{multline*}

	There is a canonical isomorphism
	$V_{(f(t),g(v))} \cong T_{f(t)}\times_T V_{g(v)},$
	and so, writing $S := (T_{f(t)})_{\red}$,
	we find that $\dim V_{(f(t),g(v))} = \dim S\times_T V_{g(v)}.$
	If we let $\rhobar_S$ denote
	the pull-back of the family $\rhobar_T$ over $S$,
	then the Galois representations attached to all the
	closed points of $S$ are isomorphic to $\rhobar_{f(t)}$.
	The discussion of Section~\ref{subsec:isotrivial} 
	shows we may find an \'etale cover $S' \to S$ such
	that the pulled-back family $\rhobar_{S'}$ is constant.  We may replace
	$S$ by $S'$ without changing the dimension
	we are trying to estimate, and thus we may furthermore
	assume that $\rhobar_S$ is the trivial family with fibres
	$\rhobar_{f(t)}$.

	If we let $C^0_S \to Z^1_S$
	denote the pull-back of the good complex $C^0_T \to Z^1_T$
	to $S$, then this complex is also the pull-back to
	$S$ of the good complex
	$C^0_t \to Z^1_t$ corresponding to the representation $\rhobar_t$.
	If we let $W$ denote the vector space 
	$H^1(G_K,\rhobar_t\otimes\alphabar^{\vee})^{\oplus n}$, 
	thought of as an affine space over $\Fbar_p$,
	then there is a projection $S \times_T V \to 
	S\times_{\Fbar_p} W$, whose kernel is a (trivial) vector bundle
	which we denote by~$V'$.

	The affine space
	$W$ parameterizes a universal family of extensions
        $0 \to \rhobar_{f(t)} \to \rhobar_W \to \alphabar^{\oplus n} \to 0$,
	giving rise to a morphism $h: W \to (\cX_{d+ a n,\red})_{\Fpbar},$
        and the composite morphism
	$S\times_T V \to V \buildrel g \over \longrightarrow 
	(\cX_{d+ a n,\red})_{\Fpbar}$ admits an alternative 
	factorization as
       	$$S \times_T V \to S\times_{\Fbar_p} W \to W \buildrel h \over 
	\longrightarrow (\cX_{d+ a n,\red})_{\Fpbar}.$$
	Thus 
	$$S \times_T V_{g(v)} = S\times_T V \times_W W_{h(w)}$$
	(where $w$ denotes the image of $v$ in $W$ under the projection),
	and hence
	\begin{multline*} \dim S \times_T V_{g(v)} 
	= \rk V' + \dim S + \dim W_{h(w)} \\  =
	n\bigl(\rk Z_T^1 - \dim H^1(G_K,\rhobar_{f(t)}\otimes \alphabar^{\vee})\bigr)
	+ \dim S + \dim W_{h(w)},\end{multline*}
	so that our desired inequality becomes equivalent to the inequality
	\numequation\label{eqn: second dimension inequality}\dim W_{h(w)} 
	\buildrel ? \over \geq
	n \dim H^0(G_K, \rhobar_{f(t)}\otimes \alphabar^{\vee}) + n^2\dim G_{\alphabar} + \dim G_t.\end{equation}
        There is an action of $G_t \times \Aut_{G_K}(\alphabar^{\oplus n})$ on
	$$W = H^1(G_K,\rhobar_{f(t)}\otimes \alphabar^{\vee})^{\oplus n} \cong \Ext^1_{G_K}(\alphabar^{\oplus n},
	\rhobar_t),$$
	which lifts to an action on the family $\rhobar_W.$
	There is furthermore a (unipotent)
	action of $H^0(G_K,\rhobar_t\otimes \alphabar^{\vee})^{\oplus n}$
	on $\rhobar_W$ (lying over the trivial action on $W$). 
	Altogether, we find that the morphism $h$
	factors through $[W/ \bigl( H^0(G_K,\rhobar_t \otimes \alphabar^{\vee})^{\oplus n} 
	\rtimes (G_t \times \Aut_{G_K}(\alphabar^{\oplus n})) \bigr) ]$,
	implying that 
	$$\dim W_{h(w)} \geq n\dim H^0(G_K,\rhobar_t\otimes \alphabar^{\vee}) + n^2\dim G_{\alphabar} + \dim G_t,$$
	as required.

        We now turn to~(2), so that we assume that~$n=1$,
        that~$\alphabar$ is one-dimensional, and
        that~$\rhobar_T$ is generically maximally nonsplit of niveau~$1$. By
        Proposition~\ref{prop: structure of maximally nonsplit}, after
        possibly replacing~$T$ with a non-empty open subscheme, we can and do
        assume that~$\rhobar_T$ is maximally nonsplit of niveau one, that we can
        write~$\rhobar_T$ as an extension of
        families \numequation\label{eqn: filtration on rhobarT}0\to\rbar_T\to\rhobar_T\to\betabar_T\to 0\end{equation}where the
        family~$\betabar_T$ arises from twisting a
        character~$\betabar:G_K\to\Fpbartimes$ via a morphism
        $T\to\Gm$, and that if~$K=\Qp$, then~(\ref{item: weird
          condition over Qp}) holds for all finite type points~$t$
        of~$T$. (Indeed, $\betabar_T$ is the family corresponding to
        the character~$\ur_{\lambda_d}\otimes\psi_{\underline{n}_d}$
        in the notation of Proposition~\ref{prop: structure of maximally
    nonsplit}.)

We want to show that~$\cE_V$ is generically maximally nonsplit of niveau~$1$.  Assume that this is the case; 
we now show that we then have equality in the bound of part~(1). After
replacing~$V$ with a non-empty open subscheme, we can and do suppose that~
$\cE_V$ is maximally nonsplit of niveau~$1$, so 
that in particular, for every finite type point~$v$ of~$V$, the extension
\[0\to\rhobar_{f(t)}\to\cE_{g(v)}\to\alphabar\to 0\]is maximally
nonsplit. From the proof of part~(1), we see that it
is enough to show that (after shrinking
~$V$ as we have) 
equality holds in both~(\ref{eqn: first dimension inequality})
and~(\ref{eqn: second dimension inequality}).




We begin by considering~(\ref{eqn: first dimension
  inequality}). For equality to hold here, it is enough to check
that the immersion~\eqref{eqn:residual gerbe}
induces an isomorphism on underlying reduced substacks
$$([\Spec \Fbar_p/G_t]\times_{(\cX_{d,\red})_{\Fpbar}} V_{g(v)})_{\red}
\iso (V_{g(v)})_{\red}.$$
We can check this on the level of $\Fbar_p$-valued points,
for which it suffices
to show that the
representation~$\cE_{g(v)}$ has a unique $d$-dimensional
subrepresentation, namely~$\rhobar_{f(t)}$; but this is immediate from
$\cE_{g(v)}$ being maximally non-split of niveau~$1$.

Similarly, to show that equality holds in ~(\ref{eqn: second dimension
  inequality}), it is enough to show that the morphism
\[[W/ \bigl( H^0(G_K,\rhobar_ {f(t)}\otimes \alphabar^{\vee}) \rtimes \bigl(G_t
  \times \Aut_{G_K}(\alphabar)\bigr) \bigr) ]\to(\cX_{d+a,\red})_{\Fbar_p}\]is a
monomorphism. To see this, note that  it follows from the definition of ``maximally nonsplit'' that any
automorphism of $\cE_{g(v)}$ induces automorphisms of~$\rhobar_{f(t)}$
and~$\alphabar$, so that we
have \[W\times_{(\cX_{d+a,\red})_{\Fbar_p}}W = W \times_{\Fbar_p}
\Bigl( H^0(G_K,\rhobar_t \otimes \alphabar^{\vee}) \rtimes \bigl(G_t
  \times \Aut_{G_K}(\alphabar)\bigr)\Bigr).\] Hence it is enough to show that the
morphism $[W/W\times_{(\cX_{d+a,\red})_{\Fbar_p}}W]\to(\cX_{d+a,\red})_{\Fbar_p}$ is a
monomorphism; this follows from Lemma~\ref{lem:stacks as quotients}.

In order to prove~(2), it remains to prove that~$\cE_V$ is generically maximally nonsplit of niveau~$1$. We need to show that after possibly
shrinking~$V$, for each $v\in V$ the image in
$\Ext^1_{G_K}(\alphabar,\betabar_{f(t)})$ of the element
of~$\Ext^1_{G_K}(\alphabar,\rhobar_{f(t)})$ corresponding to~$\cE_{g(v)}$ is
nonzero. To this end,
note that,
after possibly shrinking~$T$ (and remembering that $T$ is reduced),
we may assume that
  $H^2(G_K,\betabar_T\otimes\alphabar^\vee)$ is locally free of some
constant  rank. We can therefore repeat the construction that we
carried out after Lemma~\ref{lem: Ext2 locally free}, and find a good
complex \[C^0_T(\betabar)\to Z^1_T(\betabar)\] 
whose cohomology in degree~$1$, which we denote by~$H^1_T(\betabar)$,  is compatible with base change and
computes $H^1(G_K,\betabar_t\otimes\alphabar^\vee)=\Ext^1_{G_K}(\alphabar,\betabar_t)$ at finite type
points~$t$ of~$T$. By~\cite[\href{http://stacks.math.columbia.edu/tag/064E}{Tag
  064E}]{stacks-project} we have a morphism of complexes
\[(C^0_T\to Z^1_T)\to(C^0_T(\betabar)\to Z^1_T(\betabar)),\] compatible
with the natural maps on the cohomology groups
induced by the corresponding morphism of Herr complexes.

The kernel of~$Z^1_T\to H^1_T(\betabar)$ is a coherent sheaf, so after
possibly shrinking the reduced scheme~$T$, we can suppose that it is a
vector subbundle of~$Z^1_T$. By definition, we see that if we delete
from~$V$ the corresponding subbundle
 then the
required condition holds. It therefore suffices to show that
we are deleting a {\em proper} subbundle of~$V$.  
If this were not
the case, then (considering the fibre of $V$ and the subbundle under
consideration over some closed point~$t$), we would have an exact
sequence \[0\to\Ext^1_{G_K}(\alphabar,\betabar_t)\to\Ext^2_{G_K}(\alphabar,\rbar_t)\to\Ext^2_{G_K}(\alphabar,\rhobar_t)\to\Ext^2_{G_K}(\alphabar,\betabar_t)\to
  0; \]equivalently, by Tate local duality we would have an exact
sequence \[0\to\Hom_{G_K}(\betabar_t,\alphabar(1))\to\Hom_{G_K}(\rhobar_t,\alphabar(1))\to\Hom_{G_K}(\rbar_t,\alphabar(1))\to\]\[\Ext^1_{G_K}(\betabar_t,\alphabar(1))\to
  0. \]Assume
for the sake of contradiction that this is the case. Since $\rhobar_t$ is maximally nonsplit, the map
$\Hom_{G_K}(\betabar_t,\alphabar(1))\to\Hom_{G_K}(\rhobar_t,\alphabar(1))$ is an
isomorphism; so it suffices to show that we cannot have an isomorphism
$\Hom_{G_K}(\rbar_t,\alphabar(1))\to\Ext^1_{G_K}(\betabar_t,\alphabar(1))$. 

Since~$\rbar_t$ is maximally nonsplit, $\Hom_{G_K}(\rbar_t,\alphabar(1))$ is
at most 1-dimensional with equality if and only if $\ur_{\lambda_{d-1}(t)}\otimes\psi_{\underline{n}_{d-1}}=\alphabar(1)$; while by the local Euler characteristic
formula, ~$\Ext^1_{G_K}(\betabar_t,\alphabar(1))$ has
dimension at least $[K:\Qp]$, with equality if and only
if~$\betabar_t\ne\alphabar$ and~$\betabar_t\ne\alphabar(1)$. 
It follows that we must have~$K=\Qp$, $\ur_{\lambda_{d-1}(t)}\otimes\psi_{\underline{n}_{d-1}}=\alphabar(1)$, ~$\betabar_t\ne\alphabar$ and~$\betabar_t\ne\alphabar(1)$; 
but this case was excluded by assumption~(\ref{item: weird condition
  over Qp}), and we have our a contradiction.

Finally, suppose that we are in the situation of~(3). We wish to argue
in the same way as part~(2), but rather than deleting the split locus,
we need to delete the peu ramifi\'ee locus.
To be more precise, note that the surjection $\rhobar_T \to \alphabar(1)_T$
(where $\alphabar(1)_T$ denotes the constant rank one family given
by spreading out $\alphabar(1)$  over~$T$)
induces a composite morphism $Z^1_T \to H^1(G_K, \rhobar_T\otimes \alphabar^{\vee})
\to H^1(G_K, \epsilonbar_T)$
(where $\epsilonbar_T$ denotes
the spreading out of the  mod $p$ cyclotomic character  $\epsilonbar$
over~$T$),  
which in turn induces 
a morphism $V \to Y$ of total spaces of vector bundles
over $T$, with $Y$ denoting the total space of the trivial bundle
over~$T$ having fibre~$H^1(G_K, \epsilonbar)$.   Now $Y$ contains
a trivial subbundle $Y_{\mathrm{p.r.}}$ of codimension~$1$, parameterizing 
the peu ramifi\'ee classes.    The preimage of $Y_{\mathrm{p.r.}}$ is
a closed subscheme $V_{\mathrm{p.r.}}$ of~$V$,
and we wish to show that $V_{\mathrm{p.r.}}$ is a proper closed subscheme
of~$V$.

For this, it is enough to show that for each~$t$ (perhaps after replacing
$T$ by a non-empty open subset),
there is some extension of~$\alphabar$ by~$\rhobar_{f(t)}$ with the
property that the induced extension of~$\alphabar$
by~$\betabar_t=\alphabar(1)$ is tr\`es ramifi\'ee. We have an exact sequence
\[\Ext^1_{G_K}(\alphabar,\rhobar_t)\to
  \Ext^1_{G_K}(\alphabar,\betabar_t)\to\Ext^2_{G_K}(\alphabar,\rbar_t),\]so
we are done if $\Ext^2_{G_K}(\alphabar,\rbar_t)=0$. Since~$\rbar_t$ is
maximally nonsplit, this means that we are done unless~$\rbar_t$
admits~$\alphabar(1)$ as its unique quotient character. 

By hypothesis, this implies that~$\epsilonbar$ is trivial, 
and we can assume that the extension
of~$\betabar_t=\alphabar(1)=\alphabar$ by~$\alphabar(1)=\alphabar$
induced by~$\rhobar_t$ is tr\`es ramifi\'ee. Write~$c_t$ for the
corresponding class in~$H^1(G_K,\epsilonbar)$. As in
Remark~\ref{rem: further explanation of weird condition} below, it is
enough to show that there is a tr\`es ramifi\'ee extension
of~$\alphabar$ by~$\alphabar(1)$, corresponding to a
class~$d\in H^1(G_K,\epsilonbar)$ such that the cup product of~$c_t$
and~$d$ vanishes. If this is not the case, then~$c_t$ must be an
unramified class in $H^1(G_K,1)=H^1(G_K,\epsilonbar)$; but the
extension of~$K$ cut out by~$c_t$ considered as a class in
$H^1(G_K,1)$ is an extension given by adjoining a $p$th root of a
uniformizer, so is in particular a ramified extension, and we are done. 
\end{proof}
\begin{rem}
  \label{rem: further explanation of weird condition}The curious
  looking condition Proposition~\ref{prop:dimension control}~(\ref{item: weird condition over Qp})  is not
  merely an artefact of our arguments. Indeed, if we fix
  characters~$\alphabar,\betabar:G_{\Qp}\to\Fpbartimes$, with
  $\betabar\ne\alphabar$, $\betabar\ne\alphabar(1)$, and let~\[\rbar=
  \begin{pmatrix}
    \alphabar(1)&*\\0&\betabar
  \end{pmatrix}\]be a non-split extension, then any extension of~$\alphabar$ by~$\rbar$ induces
the trivial extension of~$\alphabar$ by~$\betabar$. To see this, note that if we fix
an extension of~$\alphabar$ by~$\betabar$, then the condition that we
can form an extension  of~$\alphabar$ by~$\rbar$ realising this
extension  of~$\alphabar$ by~$\betabar$ is that the cup product of the corresponding classes
in~$H^1(G_{\Qp},\alphabar\betabar^{-1}(1))$
and~$H^1(G_{\Qp},\betabar\alphabar^{-1})$ vanishes. But by Tate
local duality, this cup product is a perfect pairing
of~$1$-dimensional vector spaces, so one of the two classes must
vanish, as required.\end{rem}

\section{$\cX_{\lowercase{d}}$ is a formal algebraic stack}\label{subsec: formal
  algebraic stack}We now use the theory of families of extensions to
show that~$\cX_d$ is a formal algebraic stack. The key ingredient is
Theorem~\ref{thm:Xdred is algebraic}, which we prove by induction
on~$d$. We begin by setting up some useful terminology, which is
motivated by the generalisations of the weight part of Serre's
conjecture formulated in~\cite{2015arXiv150902527G}. Recall that the
notion of a \emph{Serre weight} was defined in Section~\ref{subsec:
  notation and conventions}. Since~$\F$ is assumed to contain~$k$, we
can and do identify embeddings $k\into\Fpbar$ and $k\into\F$.


\begin{defn}\label{defn: max nonsplit weight k}
  If~$\underline{k}$ is a Serre weight, and~$\rhobar:G_K \to \GL_d(\Fbar_p)$ is
  maximally nonsplit
  of niveau~$1$, then we say that~$\rhobar$
  is of weight~$\underline{k}$ if we can \index{maximally nonsplit of
    niveau 1!of weight~$\underline{k}$}
  write 
  \[\rhobar\cong
    \begin{pmatrix}
      \chi_1 &*&\dots &*\\
      0&\chi_2&\dots &*\\
      \vdots&& \ddots &\vdots\\
      0&\dots&0&\chi_d\\
    \end{pmatrix}
  \] where
  \begin{itemize}
  \item
    $\chi_i|_{I_K}=\epsilonbar^{1-i}\prod_{\sigmabar:k\into\F}\omega_{\sigmabar}^{-k_{\sigmabar,d+1-i}}$, 
    and
  \item if $(\chi_{i+1}\chi_{i}^{-1})|_{I_K}=\epsilonbar^{-1}$,
    then 
    $k_{\sigmabar,d-i}-k_{\sigmabar,d+1-i}=p-1$ for all $\sigmabar$ if
    and only if
    $\chi_{i+1}\chi_{i}^{-1}=\epsilonbar^{-1}$ and the element
    of~$\Ext^1_{G_K}(\chi_i,\chi_{i+1})=H^1(G_K,\epsilonbar)$
    determined by~$\rhobar$ is
    tr\`es ramifi\'ee  (and otherwise
    $k_{\sigmabar,d-i}-k_{\sigmabar,d+1-i}=0$ for all~$\sigmabar$).
  \end{itemize}
\end{defn}


Note that each maximally nonsplit representation of niveau one is of
weight~$\underline{k}$ for a unique~$\underline{k}$. For each~$\underline{k}$,
there exists a~$\rhobar$ which is maximally nonsplit of niveau~$1$ and
weight~$\underline{k}$; 
the existence of such a~$\rhobar$ follows by an easy
induction, and is in any case an immediate consequence of Theorem~\ref{thm:Xdred
  is algebraic} below.

\begin{defn}
  \label{defn: shifted weight}We say that a weight~$\underline{k}'$ is
  a \emph{shift} of a weight~$\underline{k}$ \index{shift (of a Serre weight)}
  if~$k_{\underline{\sigma},d}'=k_{\underline{\sigma},d}$ for
  all~$\sigmabar$, and if for each~$1\le i\le d-1$, we either have
  $k'_{\sigmabar,i}-k'_{\sigmabar,i+1}=k_{\sigmabar,i}-k_{\sigmabar,i+1}$
  for all~$\sigmabar$, or we have
  $k'_{\sigmabar,i}-k'_{\sigmabar,i+1}=p-1$ and
  $k_{\sigmabar,i}-k_{\sigmabar,i+1}=0$ for all~$\sigmabar$.
(Note in particular that any weight is a shift of itself.)
\end{defn}

\begin{defn}
  \label{defn: abuse of terminology crystalline lattice}We say that a
  representation~$\rho:G_K\to\GL_d(\cO)$ or~$\rho:G_K\to\GL_d(\Zpbar)$
  is \emph{crystalline}, or~\emph{crystalline of weight~$\lambda$},
  if~$\rho[1/p]$ is crystalline, resp.\ crystalline of
  weight~$\lambda$. We say that~$\rho$ is a \emph{crystalline lift} of
  $\rhobar:G_K\to\GL_d(\F)$ (resp.\ of $\rhobar:G_K\to\GL_d(\Fpbar)$).
\end{defn}


The
motivation for Definitions~\ref{defn: max nonsplit weight k}
and~\ref{defn: shifted weight} is the following lemma.
\begin{lem}
  \label{lem: maxnonsplit weight k implies crystalline lifts to
    shifts}If~$\rhobar$ is maximally nonsplit of niveau~$1$ and
  weight~$\underline{k}$, if~$\underline{k}'$ is a shift
  of~$\underline{k}$, and~$\underline{\lambda}$ is a lift
  of~$\underline{k}'$, then~$\rhobar$ has a crystalline lift {\em (}defined
over $\Zbar_p${\em )} of 
  weight~$\underline{\lambda}$. Furthermore, this lift can be chosen
  to be ordinary in the sense of Definition~{\em \ref{defn: ordinary}} below.
\end{lem}
\begin{proof}Write~$\rhobar$ as in Definition~\ref{defn: max nonsplit
    weight k}, and assume (by replacing~$E$ by a finite unramified extension if
  necessary) that both~$\rhobar$ and the filtration in Definition~\ref{defn: max nonsplit
    weight k} are defined over~$\F$.  Note firstly that by~\cite[Lem.\
  5.1.6]{2015arXiv150902527G} (which uses the opposite conventions for
  the sign of the Hodge--Tate weights to those of this book),
  together with the observation made above that
  $\prod_{\sigmabar}\omega_{\sigmabar}^{e(K/\Qp)}=\epsilonbar|_{I_K}$,
  we can find crystalline characters~$\psi_i:G_K\to\cO^\times$ such
  that~$\psibar_i=\chi_i$ and~$\psi_i$ has Hodge--Tate weights given
  by~$\lambda_{\sigma,d+1-i}$.

  It is enough to prove that~$\rhobar$ has a crystalline lift which is
  a successive extension of characters~$\psi_i$ of the form \[    \begin{pmatrix}
      \psi'_1 &*&\dots &*\\
      0&\psi'_2&\dots &*\\
      \vdots&& \ddots &\vdots\\
      0&\dots&0&\psi'_d\\
    \end{pmatrix},
  \]where each~$\psi'_i$ is an unramified twist of~$\psi_i$; note that
  any such representation is by definition ordinary in the sense of
  Definition~\ref{defn: ordinary}.  We prove this by induction on~$d$,
  the case $d=1$ being trivial. In the general case, write $\rhobar$
  as an extension \[0\to\rbar\to\rhobar\to\chibar_d\to 0,\]and suppose
  inductively that~$\rbar$ has a lift of the required form. Write
  ~$r:G_K\to\GL_{d-1}(\cO)$ for this lift; then we inductively seek a
  lift of~$\rhobar$ the form \numequation\label{eqn: rho lift is an
    extension}0\to r\to \rho\to \psi'_d \to 0, \end{equation}
where~$\psi'_d$ is an unramified twist of~$\psi_d$. We insist also
that~\numequation\label{eqn: ratio of crystalline chars isn't
  cyclotomic}\psi'_{d-1}(\psi'_d)^{-1}\ne \epsilon;\end{equation}
note that this could only fail
if~$k'_{\sigmabar,d}=k'_{\sigmabar,d-1}$ for all~$\sigmabar$, and in
this case we are only ruling out a single unramified twist.  Then
since~$\rhobar$ is maximally nonsplit, we have
$\Hom_{G_K}(r,\psi_d\epsilon)=0$, and it follows as in the proof of
Lemma~\ref{lem: extensions are crystalline} that such a~$\rho$, if it exists,
is automatically
crystalline.

It is therefore enough to show that there is a lift of the
form~\eqref{eqn: rho lift is an extension} satisfying~\eqref{eqn:
  ratio of crystalline chars isn't cyclotomic}.  As in the proof
of~\cite[Thm.\ 2.1.8]{2015arXiv150601050G}, we have an exact sequence
\nummultline\label{eqn: lifting ses H1
  H2}H^1\bigl(G_K, (\psi'_d)^{-1}\otimes r\bigr) \to 
H^1(G_K, \chi_d^{-1} \otimes \rbar) 
\stackrel{\delta}{\to}
H^2\bigl(G_K, (\psi'_d)^{-1} \otimes r\bigr) 
\end{multline}
Thus,
if~$c\in H^1(G_K, \chi_d^{-1} \otimes \rbar)$ 
is the class determined
by~$\rhobar$, it is enough to show that we can choose the unramified
twist~$\psi'_d$ so that~$\delta(c)=0$. Now,
$H^2\bigl(G_K, (\psi'_d)^{-1}\otimes r\bigr)$ 
is dual to
$H^0\bigl(G_K, r^{\vee} \otimes (\psi'_d \epsilon) \otimes E/\cO\bigr),$
and
since~$\rbar$ is maximally nonsplit, we see that this latter group is equal
to $H^0\bigl(G_K,  (\psi'_{d-1})^{-1}\psi'_d\epsilon \otimes E/\cO\bigr).$ 
In other
words, to show that~$\delta(c)=0$ in~\eqref{eqn: lifting ses H1 H2}, it
is enough to show the corresponding statement where~$r$ is replaced
by~$\psi'_{d-1}$; so we need only show that we can choose ~$\psi'_{d}$
so that the extension of~$\chi_d$ by~$\chi_{d-1}$ lifts to a
crystalline extension of~$\psi'_d$ by~$\psi'_{d-1}$. This is an
immediate consequence of the weight part of Serre's conjecture
for~$\GL_2$, as proved in~\cite{MR3324938}; see in
particular~\cite[Thm.\ 6.1.18]{MR3324938} and the references therein
(and note that these results also hold for~$p=2$,
by~\cite{2017arXiv171109035W}).
\end{proof}

\begin{remark}
We emphasize that the preceding lemma applies only to $\rhobar$ that
are maximally split of niveau~$1$.  
In Theorem~\ref{thm: strong existence of crystalline lifts}
below we prove an analogous result for arbitrary~$\rhobar$.
\end{remark}

Suppose that~$\rhobar_{\cT}$ is
generically maximally nonsplit of niveau~$1$, where~$\cT$ is some
integral finite type $\Fpbar$-stack. We say that~$\rhobar_{\cT}$ is
generically of weight~$\underline{k}$ if there is a dense open
substack~$\cU$ of~$\cT$ such that every $\Fpbar$-point of~$\rhobar_{\cU}$ is
maximally nonsplit of niveau~$1$, and is of weight~$\underline{k}$.
\begin{lem}
\label{lem:weighing families}
  If~$\rhobar_{\cT}$ is generically maximally nonsplit of
niveau~$1$, then~$\rhobar_{\cT}$ is generically of weight~$\underline{k}$ for a unique
Serre weight~$\underline{k}$.
\end{lem}
\begin{proof}
  Consider a dense open substack~$\cU\subseteq \cT$,  tuples~$\underline{n}_i$ and 
  morphisms~$\lambda_i:U\to\Gm$ as in Proposition~\ref{prop: structure of maximally
    nonsplit}. We inductively determine the~$k_{\sigmabar,i}$ as
  follows. We let~$k_{\sigmabar,d}$ be the unique choice with~$p-1\ge
  k_{\sigmabar,d}\ge 0$ not all equal to~$p-1$ and
  satisfying \[\psi_{\underline{n}_1}|_{I_K}=\prod_{\sigmabar:k\into\Fpbar}\omega_{\sigmabar}^{-k_{\sigmabar,d}}.\]
  Inductively, if $k_{\sigmabar,i+1}$ has been determined, then we
  demand that~$k_{\sigmabar,i}$ satisfies~$p-1\ge
  k_{\sigmabar,i}-k_{\sigmabar,i+1}\ge 0$  and
  \[(\psi_{\underline{n}_{d-i}}\psi_{\underline{n}_{d+1-i}}^{-1})|_{I_K}=\epsilonbar\prod_{\sigmabar:k\into\Fpbar}\omega_{\sigmabar}^{k_{\sigmabar,i}-k_{\sigmabar,i+1}}.\]This
  determines~$k_{\sigmabar,i}$, except in the case
  that~$(\psi_{\underline{n}_{d-i}}\psi_{\underline{n}_{d+1-i}}^{-1})|_{I_K}=\epsilonbar$.
  In this case, we consider the
  character~$\lambda_{d+1-i}\lambda_{d-i}^{-1}:U\to\Gm$. Since~$\cU$ is
  integral, this character is either constant or has open image. After
  possibly shrinking~$\cU$, we may therefore assume that either
  ~$\lambda_{d+1-i}(t)\ne\lambda_{d-i}(t)$ for all~$t\in \cU$, in which
  case we set $k_{\sigmabar,d+1-i}-k_{\sigmabar,d-i}=0$ for
  all~$\sigmabar$, or that the character is trivial.

  If the character is trivial, we may consider the locus in~$\cU$ over which
  the extension between~$\chi_{i+1}$ and~$\chi_i$ is peu
  ramifi\'ee. The argument used in in the proof of part~(3) of
  Proposition~\ref{prop:dimension control}
  shows that this is a constructible subset of~$\cU$; since~$\cU$ is integral,
  after shrinking~$\cU$ we can thus suppose that it
  is either equal to~$\cU$, or
  that it is empty. In the former
  case we set $k_{\sigmabar,d+1-i}-k_{\sigmabar,d-i}=0$ for
  all~$\sigmabar$, and in the latter we take
  $k_{\sigmabar,d+1-i}-k_{\sigmabar,d-i}=p-1$ for
  all~$\sigmabar$. 
  It follows from the construction that every $\Fpbar$-point
  of~$\rhobar_{\cU}$ is maximally nonsplit of niveau~$1$ and
  weight~$\underline{k}$, as required.
  \end{proof}

 It will be convenient to use the following variant of the notation
 established in Proposition~\ref{prop: structure of maximally
    nonsplit}. 
  If~$\rhobar_{\cT}$ is generically maximally nonsplit of
  niveau one and weight~$\underline{k}$, then we
  set~$\nu_i:=\lambda_{d+1-i}$, and write~$\omega_{\underline{k},i}$
  for the character~$\epsilonbar^{d-i}\psi_{\underline{n}_{d+1-i}}$,
  so that \[\omega_{\underline{k},i}|_{I_K}=\prod_{\sigmabar:k\into\Fpbar}\omega_{\sigmabar}^{-k_{\sigmabar,i}},\] and~\eqref{eqn:maximally nonsplit structure} can be rewritten as
   \numequation\label{eqn:weight k version of maximally nonsplit structure}\rhobar_t\cong \begin{pmatrix}
      \ur_{\nu_d(t)}\otimes\omega_{\underline{k},d} &*&\dots &*\\
      0& \ur_{\nu_{d-1}(t)}\otimes\epsilon^{-1}\omega_{\underline{k},d-1}&\dots &*\\
      \vdots&& \ddots &\vdots\\
      0&\dots&0& \ur_{\nu_1(t)}\otimes\epsilon^{1-d}\omega_{\underline{k},1}\\
    \end{pmatrix}.  \end{equation}
We have the \emph{eigenvalue morphism} 
~$\underline{\nu}:\cU\to(\Gm)^d$ given
by~$(\nu_1,\dots,\nu_d)$. \index{eigenvalue morphism}

\begin{df}
If $\underline{k}$ is a Serre weight, we let
$(\Gm)^d_{\underline{k}}$ denote the closed subgroup scheme of $(\Gm)^d$ parameterizing 
tuples~$(x_1,\dots,x_d)$ for which~$x_i=x_{i+1}$
whenever~$k_{\sigmabar,i}-k_{\sigmabar,i+1}=p-1$ for all~$\sigmabar$. \index{$(\Gm)^d_{\underline{k}}$}
\end{df}

Note that $(\Gm)_{\underline{k}}^d$ is closed under the simultaneous multiplication
action of $\Gm$ on $(\Gm)_{\underline{k}}^d$ (i.e.\ the action $a \cdot (x_1,
\ldots,x_d) := (a x_1,\ldots, a x_d)$).

It follows from the definitions that
the eigenvalue morphism $\underline{\nu}$
is valued in~$(\Gm)^d_{\underline{k}}.$

We have seen in Proposition~\ref{prop:X is an Ind-stack}
that~$\cX_{d}$ is an Ind-algebraic stack, which is in fact the $2$-colimit 
of algebraic stacks with respect to transition morphisms that
are closed immersions.   We let~$\cX_{d,\red}$ be the
underlying reduced substack of~$\cX_d$; it is then a closed substack of $\cX_d$
with the same structure:  namely, it is the $2$-colimit of closed algebraic stacks (which
can even be assumed to be reduced) with respect to transition morphisms that are closed
immersions. 
\begin{thm} 
	\label{thm:Xdred is algebraic}\leavevmode
        \begin{enumerate}
        \item 	The Ind-algebraic stack $\cX_{d,\red}$
is an algebraic stack, of finite presentation over $\F$.
\item  We can write~$(\cX_{d,\red})_{\Fpbar}$ 
as a union of closed algebraic substacks of finite
presentation
over~$\Fpbar$ 
\[(\cX_{d,\red})_{\Fpbar}=\cX_{d,\red,\Fpbar}^{\negligible}\cup\bigcup_{\underline{k}}\cX_{d,\red,\Fpbar}^{\underline{k}},\]
where:
\begin{itemize}
\item $\cX_{d,\red,\Fpbar}^{\negligible}$ is empty if~$d=1$, and otherwise is
  non-empty of dimension
  strictly less than~$[K:\Qp]d(d-1)/2$. 
\item each $\cX_{d,\red,\Fpbar}^{\underline{k}}$ is a closed 
irreducible substack of dimension 
  ~$[K:\Qp]d(d-1)/2$, and is generically maximally nonsplit of niveau
  one and weight~$\underline{k}$.
  The corresponding eigenvalue morphism 
is dominant {\em (}i.e.\ has dense image in~$(\Gm)^d_{\underline{k}}${\em )}.
\end{itemize}

\item  If we fix an irreducible representation
$\alphabar: G_K \to \GL_a(\Fbar_p)$ {\em (}for some $a \geq 1${\em )},
then the locus of $\rhobar$ in
$\cX_{d,\red}(\Fbar_p)$ for which
$\dim \Hom_{G_K}(\rhobar,\alphabar) \geq r$ {\em (}for any $r \geq 1${\em )}
is {\em (}either empty,
or{\em )} of dimension at most
\[[K:\Qp]d(d-1)/2- \lceil {r\bigl((a^2+1)r-a\bigr)}/{2} \rceil.\] 
Furthermore, the locus
of $\rhobar$ in
$\cX_{d,\red,\Fpbar}^{\negligible}(\Fbar_p)$ for which
$\dim \Hom_{G_K}(\rhobar,\alphabar) \geq r$ is of dimension strictly
less than this.

\item  
If we fix an irreducible representation
$\alphabar: G_K \to \GL_a(\Fbar_p)$ {\em (}for some $a \geq 1${\em )},
then the locus of $\rhobar$ in
$\cX_{d,\red}(\Fbar_p)$ for which
$\dim \Ext^2_{G_K}(\alphabar,\rhobar) \geq r$ is of dimension at most
\[[K:\Qp]d(d-1)/2-r.\]
        \end{enumerate}
\end{thm}
\begin{rem}
  Note that in parts~(3) and~(4), the locus of points in question
  corresponds to a closed substack of~$(\cX_{d,\red})_\Fpbar$, by
  upper-semicontinuity of fibre dimension. 
\end{rem}

\begin{rem}\label{rem: different k give different components}
  Since the~$\cX_{d,\red,\Fpbar}^{\underline{k}}$ are irreducible, have
  dimension equal to that of~$(\cX_{d,\red})_\Fpbar$, and 
  have pairwise disjoint open substacks (corresponding to maximally
  nonsplit representations of niveau~$1$ and weight~$\underline{k}$),
  they are in fact distinct irreducible components of~$(\cX_{d,\red})_\Fpbar$.
\end{rem}

\begin{remark}
	We can, and will, be much more precise about the structure of $\cX_{d,\red}$.
	Namely, in Chapter~\ref{sec:properties} we combine
        Theorem~\ref{thm:Xdred is algebraic} with additional arguments
        to show that $\cX_{d,\red}$ is equidimensional of dimension
        $[K:\Q_p]d(d-1)/2$; 
        accordingly, the
        irreducible components of~$(\cX_{d,\red})_\Fpbar$ are
        precisely the
        ~$\cX_{d,\red,\Fpbar}^{\underline{k}}$, and in particular
        are in bijection with the Serre
        weights~$\underline{k}$. We also show that these irreducible
        components are all defined over~$\F$. 
      \end{remark}
      \begin{rem}\label{rem: open eigenvalue morphism implies ratio
          not cyclo}
        Note that since the eigenvalue morphism has dense image, it
        follows that~$\cX_{d,\red}^{\underline{k}}$ contains a dense
        open substack with the property that $\nu_i(t)\ne
        \nu_{i+1}(t)$ unless $k_{\sigmabar,i}-k_{\sigmabar,i+1}=p-1$
        for all~$\sigmabar$.
      \end{rem}
	\begin{remark}
          As will be evident from the proof, the upper bound of
          Theorem~\ref{thm:Xdred is algebraic}~(3) is quite crude when
          $[K:\Q_p] > 1$, although it is reasonably sharp in the case
          $K = \Q_p$.  However, it suffices for our purposes, and
          indeed we will only use the (even cruder) consequence Theorem~\ref{thm:Xdred is
            algebraic}~(4) in the rest of the book. 
	\end{remark}

\begin{proof}[Proof of Theorem~{\em \ref{thm:Xdred is algebraic}}] 
 Recall that a closed immersion of
  reduced algebraic stacks that are locally of finite type over $\F_p$
  which is surjective on finite-type points is necessarily an
  isomorphism. As recalled above, $\cX_{d,\red}$ is an inductive limit of such
  stacks (indeed, by Lemma~\ref{alem: underlying reduced of Ind algebraic}, we have~$\cX_{d,\red}=\varinjlim\cX^a_{d,h,s,\red}$,
where the~$\cX^a_{d,h,s}$ are as in Section~\ref{subsec: X is Ind
  algebraic}), and so 
if we produce closed algebraic substacks
  $\cX_{d,\red,\Fpbar}^{\negligible}$ and $\cX_{d,\red,\Fpbar}^{\underline{k}}$
 of $\cX_d$,
  the union of whose $\Fpbar$-points exhausts those
  of~$\cX_{d,\red}$, then $\cX_{d,\red,\Fpbar}$ will in fact
be an algebraic stack which is the union of its closed substacks
 $\cX_{d,\red,\Fpbar}^{\negligible}$ and $\cX_{d,\red,\Fpbar}^{\underline{k}}$. 
Thus~(1) is an immediate consequence of~(2) (where the ``union'' statement
in~(2) is now to be understood on the level of $\Fbar_p$-points).

Claim~(4) follows
  from~(3) (with~$\alphabar$ replaced by
  $\alphabar\otimes\epsilonbar$) 
  by Tate local duality and the easily verified inequality
  \[\lceil {r\bigl((a^2+1)r-a\bigr)}/{2} \rceil\ge r.\]  Thus it is enough to
  prove~(2) and~(3), which we do simultaneously by induction
  on~$d$. 

As recalled in Remark~\ref{rem: irreducible representations finite up
  to twist}, there are up
to twist by unramified characters only finitely many irreducible
$\Fpbar$-representations of~$G_K$ of any fixed dimension. 
Accordingly,
we let $\{\alphabar_i\}$ be a finite set of irreducible continuous
representations $\alphabar_i: G_K \to \GL_{d_i}(\Fpbar)$, such that
any irreducible continuous representation of $G_K$ over $\Fbar_p$ of
dimension at most $d$ arises as an unramified twist of exactly one of
the~$\alphabar_i$. We let the 1-dimensional representations in this set be the
characters  ~$\psi_{\underline{n}}$ defined in Section~\ref{subsec: twists of
  families}. 

Each~$\alphabar_i$ 
  corresponds to a finite type point
  of~$\cX_{d_i,\red}$, whose associated residual gerbe
  is a substack of~$\cX_{d_i,\red}$ of dimension~$-1$:
  the morphism $\Spec\Fpbar\to\cX_{d_i,\red}$ corresponding
  to $\alphabar_i$ factors through
  a monomorphism $[\Spec\Fpbar/\Gm]\to\cX_{d_i,\red}$.
  It follows from
  Lemma~\ref{lem: dimension of Gm twist} that for each~$\alphabar_i$
  there
  is an irreducible closed zero-dimensional
  algebraic substack of~$\cX_{d_i,\red}$ of finite presentation
  over~$\Fpbar$ 
 whose $\Fpbar$-points are exactly the unramified
  twists of~$\alphabar_i$.

 In particular, if $d=1$, then 
  we let $\cX_{1,\red}^{\underline{k}}$ be the
  zero-dimensional stack constructed in the previous paragraph, whose $\Fpbar$-points
  are the unramified twists of~$\psi_{\underline{k}}$; this satisfies the
  required properties by definition, so~(2) holds when~$d=1$. For~(3),
  note that if~$r>0$ then
  we must have~$a=1$, and then the locus where
  $\Hom_{G_K}(\rhobar,\alphabar)$ is non-zero (equivalently, $1$-dimensional) is exactly the closed
  substack of dimension~$-1$ corresponding to~$\alphabar$, so the
  required bound holds. 

  We now begin the inductive proof of~(2) and~(3) for~$d>1$. In fact,
  it will be helpful 
  to simultaneously prove another statement~(2'), which is as follows:
  for each~$\underline{k}$, there is a closed irreducible algebraic
  substack~$\cX_{d,\red,\Fpbar}^{\underline{k},\textrm{fixed}}$
  of~$(\cX_{d,\red})_\Fpbar$ of finite presentation over~$\Fpbar$ and
  dimension~$[K:\Qp]d(d-1)/2-1$, which is generically maximally nonsplit of niveau~$1$ and
  weight~$\underline{k}$, and furthermore has the property that the
  corresponding character~$\nu_1$ is trivial. 
  The discussion of the previous paragraph also proves~(2')
  when~$d=1$. The point of 
  hypothesis~(2') is that we can use Proposition~\ref{prop:dimension
    control}~(2) to
  construct~$\cX_{d,\red,\Fpbar}^{\underline{k},\textrm{fixed}}$ by using~(2)
  and~(2') in dimension $(d-1)$, and then
  construct~$\cX_{d,\red,\Fpbar}^{\underline{k}}$ from it using
  Lemma~\ref{lem: dimension of Gm twist}.

  We now prove the inductive step, so we assume that~(2), (2') and~(3)
  hold in dimension less than~$d$. We begin by constructing
  the closed substacks~$\cX_{d,\red}^{\underline{k},\textrm{fixed}}$
  and~$\cX_{d,\red,\Fpbar}^{\underline{k}}$.  Let~$\underline{k}_{d-1}$ be the Serre weight in
  dimension~$(d-1)$ obtained by deleting the first entry
  in~$\underline{k}$. 
  We
  set~$\alphabar:=\epsilonbar^{1-d}\omega_{\underline{k},1}$. If~$k_{\sigmabar,1}-k_{\sigmabar,2}=p-1$
  for all~$\sigmabar$ then we say that we are in the \index{tr\`es
    ramifi\'ee} 
  \emph{tr\`es ramifi\'ee case}, and we let~$\cT$ denote
  $\cX_{d-1,\red,\Fpbar}^{\underline{k}_{d-1},\textrm{fixed}}$; otherwise,
  we let~$\cT$ denote $\cX_{d-1,\red,\Fpbar}^{\underline{k}_{d-1}}$. In the
  latter case, it follows from our inductive
  hypothesis (more precisely, from the assumption that~(2) and (3) hold in dimension~$(d-1)$), together with Tate
  local duality, that after replacing~$\cT$ with an open substack we
  can and do assume that for each $\Fpbar$-point~$t$ of~$\cT$, we have
  $\Ext^2_{G_K}(\alphabar,\rbar_t)=0$, where~$\rbar_t$
  is the $(d-1)$-dimensional representation corresponding to~$t$. If
  we are in the tr\`es ramifi\'ee case then it follows similarly 
  that after replacing~$\cT$ with an open substack we
  can and do assume that for each $\Fpbar$-point~$t$ of~$\cT$,
  $\Ext^2_{G_K}(\alphabar,\rbar_t)$ is $1$-dimensional. In either case
  we can in addition assume that~$\rbar_t$ is maximally nonsplit of
  niveau one and weight~$\underline{k}_{d-1}$.

  We let~$T$ be an irreducible scheme which smoothly covers~$\cT$, and
  we let~$\cX_{d,\red,\Fpbar}^{\underline{k},\textrm{fixed}}$ be the
  irreducible closed
  substack of $(\cX_{d,\red})_{\Fbar_p}$ constructed as the scheme-theoretic
image of $V$ in the notation of Proposition~\ref{prop:dimension
    control}. 
  Part~(2) of that proposition, together with
  the inductive hypothesis, implies
  that~$\cX_{d,\red,\Fpbar}^{\underline{k},\textrm{fixed}}$ has the
  claimed 
  dimension. (Note in particular that if~$K=\Qp$, then condition~\eqref{item: weird
    condition over Qp} of
  Proposition~\ref{prop:dimension control} holds. Indeed, if we are in the tr\`es
  ramifi\'ee case then condition~ \eqref{item: weird
    condition over Qp}(iii) holds, and otherwise the inductive hypothesis that
  the image of the eigenvalue morphism is dense in $(\Gm)^{d-1}_{\underline{k}_{d-1}}$
implies that condition \eqref{item: weird
    condition over Qp}(i) holds after shrinking~$T$.) We then let $\cX_{d,\red,\Fpbar}^{\underline{k}}$ be the
  substack obtained from~$\cX_{d,\red}^{\underline{k},\textrm{fixed}}$
  by twisting by unramified characters, which has the claimed
  dimension by Lemma~\ref{lem: dimension of Gm twist} (and the fact
  that by Proposition~\ref{prop:dimension control}~(2), there is a
  dense open substack of $\cX_{d,\red}^{\underline{k},\textrm{fixed}}$
  whose $\Fpbar$-points correspond to representations which are of
  niveau~$1$ and are maximally nonsplit, so that in particular the
  corresponding family is twistable). 
  
By construction, we see that the
stacks~$\cX_{d,\red,\Fpbar}^{\underline{k},\textrm{fixed}}$
and~$\cX_{d,\red,\Fpbar}^{\underline{k}}$ satisfy the properties required
of 
them, except possibly that in the tr\`es ramifi\'ee case, we need to check
that ~$\cX_{d,\red,\Fpbar}^{\underline{k},\textrm{fixed}}$ (and consequently
~$\cX_{d,\red,\Fpbar}^{\underline{k}}$) is generically maximally nonsplit of
niveau~$1$ and weight~$\underline{k}$. More precisely, in this case we need to check
the condition that the final extension is not just nonsplit, but is
(generically) tr\`es ramifi\'ee. This follows from
Proposition~\ref{prop:dimension control}~(3). Indeed by our
inductive assumption 
that~$\cX_{d-1,\red,\Fpbar}^{\underline{k}_{d-1}}$ is generically maximally nonsplit of niveau~$1$ and
weight~$\underline{k}_{d-1}$, and the image of the eigenvalue morphism is
dense,  we see that we are in case~(3)(b)(i) unless ~$\epsilonbar$ is trivial and
furthermore we
have~$k_{\sigmabar_2}-k_{\sigmabar_3}=p-1$ for all~$\sigmabar$, in
which case we satisfy~(3)(b)(ii).

 To complete the proof of~(2), we need only
construct~$\cX_{d,\red,\Fpbar}^{\negligible}$. By
  Proposition~\ref{prop:dimension control}, Tate local duality, upper
  semi-continuity of the fibre dimension, and the inductive
  hypothesis, we see that for each~$1\le i\le d$, each~$\alphabar_i$
  (of dimension~$a_i$, say),  and each~$s\ge 0$
  there is a finitely presented closed algebraic substack $\cX_{s,\alphabar_i,\Fpbar}$
  of~$(\cX_{d,\red})_\Fpbar$, whose $\Fpbar$-points contain all the
  representations of the form
  $0\to\rhobar_{d-{a_i}}\to\rhobar\to\alphabar_i\to 0$ for which
  $\dim_{\Fpbar}\Ext^2_{G_K}(\alphabar_i,\rhobar_{d-{a_i}})=s$, and whose
  dimension is at most 
\[[K:\Qp](d-{a_i})(d-{a_i}-1)/2-\lceil
  s(({a_i}^2+1)s-{a_i})/2\rceil+[K:\Qp]{a_i}(d-{a_i})+s-1; \] furthermore, the locus
where~$\rhobar_{d-{a_i}}$ is an $\Fpbar$-point
of~$\cX_{d-{a_i},\red,\Fpbar}^{\negligible}$ is of dimension strictly less
than this. These stacks are only nonzero for finitely many values
of~$s$. For fixed~${a_i}$, we see that as a function of~$s$, this quantity
is maximised by~$s=0$, as well as by~$s=1$ when ${a_i}=1$. (To see this, we have to maximise
the quantity $s-\lceil s(({a_i}^2+1)s-{a_i})/2\rceil$. Suppose firstly that
${a_i}>1$. Then if~$s=0$ we have $0$, while if~$s>0$ we
have
$s-\lceil s(({a_i}^2+1)s-{a_i})/2\rceil\le s-s(({a_i}^2+1)s-{a_i})/2\le
s-s({a_i}^2+1-{a_i})/2\le s-3s/2<0$. Meanwhile if~${a_i}=1$, then for~$s=0$ we
have $0$, for~$s=1$ we have $1-\lceil 1/2\rceil=1$,
while for ~$s>1$ we have
$s-\lceil s(2s-1)/2\rceil\le s-s(2s-1)/2 \le s-3s/2<0$.) It follows that
  as a function of~${a_i}$ the bound is maximised at ${a_i}=1$ and $s=0$
  or~$1$, when it is equal to $[K:\Qp]d(d-1)/2-1$, and it is otherwise
  strictly smaller.

  By Lemma~\ref{lem: dimension of Gm twist}, it follows that the locus
  in~$(\cX_{d,\red})_\Fpbar$ of representations of the form
  $0\to\rhobar_{d-{a_i}}\to\rhobar\to\alphabar'\to 0$, with~$\alphabar'$
  an unramified twist of~$\alphabar_i$ for which
  $\dim_{\Fpbar}\Ext^2_{G_K}(\alphabar',\rhobar_{d-{a_i}})=s$, is of
  dimension at most~$[K:\Qp]d(d-1)/2$, with equality holding only if
  ${a_i}=1$ and $s=0$ or~$1$. Furthermore, the locus of those
  representations for which ~$\rhobar_{d-{a_i}}$ is an $\Fpbar$-point
  of~$\cX_{d-{a_i},\red,\Fpbar}^{\negligible}$ is of dimension strictly less
  than~$[K:\Qp]d(d-1)/2$. Putting this together, we see that~(2) holds
  in dimension~$d$ if we take~$\cX_{d,\red,\Fpbar}^{\negligible}$
to be the union of the twists by unramified characters of the
substacks~$\cX_{s,\alphabar_i}$ for which $\dim\alphabar_i>1$ or~$s>1$, together
with the union of the twists by unramified characters of the
substacks of the~$\cX_{s,\alphabar_i,\Fpbar}$ for which~$\dim\alphabar_i=1$, $s=0$
or~$1$, and~$\rhobar_{d-1}$ is an $\Fpbar$-point
  of~$\cX_{d-1,\red,\Fpbar}^{\negligible}$.

Finally we prove~(3) in dimension~$d$. 
	In the case $r = 0$, there is nothing to prove.
	If $\dim \Hom_{G_K}(\rhobar,\alphabar) \geq r \geq 1,$
	then we may place $\rhobar$ in a short exact sequence
	$$0 \to \thetabar \to \rhobar \to \alphabar^{\oplus r} \to 0,$$
	where $\thetabar$ is of dimension $d- ra < d$.
        We may apply  part~(2) 
		so as to find that $\cX_{d-ar,\red,\Fpbar}$ 
has dimension at most $[K:\Q_p](d-ar)(d-ar-1)/2$.
		Let $\cU_s$ be the locally closed substack 
		of $\cX_{d-ar,\red,\Fpbar}$ over which
		$\dim H^2(G_K,\thetabar\otimes \alphabar^{\vee}) =
                s$; 
		by the inductive hypothesis, this locus has dimension at most
		$[K:\Qp](d-ar)(d-ar-1)/2-s\bigl((a^2+1)s-a\bigr)/2$, 
		and over this locus
		we may construct a universal family of extensions 
		$$0 \to \thetabar \to \rhobar_{\cU_s}
		\to \alphabar^{\oplus r} \to 0.$$
		The locus of $\rhobar$ we are interested in is
		contained in the scheme-theoretic
		image of this family in $(\cX_{d,\red})_\Fpbar$,
		and Proposition~\ref{prop:dimension control}
		shows that this scheme-theoretic image 
		has dimension bounded above by
                \[  [K:\Q_p](d-ar)(d-ar-1)/2 -s((a^2+1)s-a)/2 +r([K:\Qp]a(d-ar)+s) -r^2\]
\begin{align*}
	= &[K:\Qp]d(d-1)/2- (r(a^2+1)r-a)/2-(r-s)^2/2-(ar(ar-1))/2\\
\le &[K:\Qp]d(d-1)/2- (r(a^2+1)r-a)/2.
 \end{align*}

                Since this conclusion holds for each of the finitely
                many values of~$s$ (and since the dimension is an
                integer, allowing us to take the floor of this upper
                bound), we are done.
\end{proof}


\begin{cor}\label{cor:Xd is formal algebraic}
	$\cX_d$ is a  Noetherian formal algebraic
        stack.
\end{cor}
\begin{proof}
	Theorem~\ref{thm:Xdred is algebraic} shows that
	$\cX_{d,\red}$ is an algebraic stack that 
	is of finite presentation over $\F$ (and hence
	quasi-compact and quasi-separated).  
	Together with Proposition~\ref{prop:X is an Ind-stack},
	which shows that $\cX_d$ is isomorphic to the inductive limit
	of a sequence of finitely presented algebraic stacks with respect to transition
	morphisms that are closed immersions, 
	this implies that $\cX_d$
      is indeed a (locally countably indexed) formal algebraic stack,
      by Proposition~\ref{prop: formal stack 6.6}. By
      Remark~\ref{rem:Ind-loc. f.t. criterion}, $\cX_d$ is Ind-locally
of finite type over $\cO$,
and it then follows from
      Proposition~\ref{prop: Xd has a good obstruction theory}
      and Theorem~\ref{thm:noetherian criterion} that~$\cX_d$ is
      locally Noetherian, and hence Noetherian (as we have already
      seen that it is quasi-compact and quasi-separated).
      \end{proof}

\chapter{Crystalline lifts and the finer structure of $\cX_{\lowercase{d,\red}}$}\chaptermark{Crystalline lifts and $\cX_{\lowercase{d,\red}}$}
\label{sec:properties}
Up to this point, we have shown that $\cX_d$ is
a Noetherian formal algebraic stack. 
In this chapter,
we will additionally show that $(\cX_d)_\red$ is
equidimensional of dimension
$[K:\Q_p]d(d-1)/2$, and enumerate its irreducible components. 
In the course of doing this, we will
also prove a result on the existence of crystalline lifts
which is of independent interest. We furthermore determine the closed points
of~$(\cX_d)_\red$, and describe the maximal substack of~$\cX_d$ over which the
universal $(\varphi,\Gamma)$-module gives rise to a continuous $G_K$-representation.

\section{The fibre dimension of $H^2$ on crystalline deformation rings}
Let $\rhobar: G_K \to \GL_d(\F')$ be a continuous representation,
for some finite extension $\F'$ of $\F$, and let $\cO'$ be the ring 
of integers in a finite extension $E'/E$, with residue field~ $\F'$.
Let $\underline{\lambda}$ be  a Hodge type,
and consider the lifting ring $R :=
R_{\rhobar}^{\crys,\underline{\lambda},\cO'}$.
Let $M$ denote the universal \'etale 
$(\varphi,\Gamma)$-module
over $R$ (note that there is a universal \'etale
$(\varphi,\Gamma)$-module over~$R$, and not merely a universal
\emph{formal} \'etale $(\varphi,\Gamma)$-module,  as a consequence of Corollary~\ref{cor: crystalline deformation rings are effective
    versal}). Let $\alpha^{\circ}: G_K \to \GL_a(\cO')$ be a representation
	with $\alphabar: G_K \to \GL_a(\F')$ being absolutely irreducible,
and let~$N$ be the base change to~$R$ of the \'etale $(\varphi,\Gamma)$-module
over~$\cO'$ which corresponds to~$\alpha^{\circ}$. 

Let $\cC^{\bullet}(N^\vee\otimes M)$ denote the Herr complex
associated to $N^\vee\otimes M$, as in Section~\ref{subsec: Herr complex}; this
is a perfect complex of~$R$-modules, with cohomology concentrated in degrees~$[0,2]$,
as a consequence of
~\cite[\href{https://stacks.math.columbia.edu/tag/0CQG}{Tag
  0CQG}]{stacks-project} and Theorem~\ref{thm: Herr complex is
  perfect} (applied to the quotients~$R/\m_R^a$).
Write $\Ext^2 := H^2\bigl(\cC^{\bullet}(N^\vee\otimes M)\bigr)$.
For each $r \geq 1,$
write $$X_r := \{x \in \Spec R \, | \, \dim \kappa(x)\otimes_R \Ext^2 \geq r \}.$$
The theorem on upper semi-continuity of fibre dimension for coherent sheaves
shows that $X_r$ is a closed subset of $\Spec R$.

\begin{thm}
  \label{thm: codimension H2 in regular weight crystalline deformation ring}
  If $\Ext^2$
  is $\varpi$-power torsion,
  then for each~$r\ge 1$, $X_r$ is of codimension $\geq r+1$ in $\Spec R$.
\end{thm}
\begin{proof}
	The assumption that $\Ext^2$ is $\varpi$-power torsion
	implies that $\Ext^2$ is (set-theoretically) supported
	on the closed subscheme $\Spec (R/\varpi)_{\red}$ of $\Spec R$.
        Write $\overline{M}$ and $\overline{N}$ to denote the
        pull-backs of $M$ and~$N$ respectively,
	over $(R/\varpi)_{\red}$,
	write $\overline{\Ext}^2 := \Ext^2\bigl(\cC^{\bullet}(\overline{N}^\vee\otimes\overline{M})\bigr),$
	and define
        $$\overline{X}_r := \{x \in \Spec (R/\varpi)_{\red} \, | \,
	\dim \kappa(x)\otimes_{(R/\varpi)_{\red}} \overline{\Ext}^2 \geq r \}.$$
	Since~$\Ext^2$ is compatible with base change, and is supported
        on~$\Spec (R/\varpi)_{\red}$, while~$R$ is~$\cO$-flat,  it suffices to show that
the codimension of	$\overline{X}_r$  in
	$\Spec (R/\varpi)_{\red}$ is at least~$r$.
	Since $\Spec (R/\varpi)_{\red}$
	is equidimensional of dimension $[K:\Q_p]d(d-1)/2 + d^2,$
	it suffices in turn to show that $\overline{X}_r$
	is of dimension at most $([K:\Q_p]d(d-1)/2) + d^2 - r.$

	By Proposition~\ref{prop:versal rings},
	we have a versal morphism
	$\Spf R^{\square,\cO'}_{\rhobar}
	\to \cX_d$.
	Consider the fibre product
	$\Spf R^{\square,\cO'}_{\rhobar} \times_{\cX_d} \cX_{d,\red}.$
	This is a closed formal algebraic subspace of 
	$\Spf R^{\square,\cO'}_{\rhobar}$,
	and so by Lemma~\ref{lem: alem closed formal algebraic scheme
    adic*} 
is of the form $\Spf S$ for some quotient $S$ of 
	$R^{\square,\cO'}_{\rhobar}$. 
	Its description as a fibre product
	shows that the morphism $\Spf S \to \cX_{d,\red}$ is versal.
	Since $\cX_{d,\red}$ is an algebraic stack,
	this morphism is effective, i.e.\ arises from a morphism
	\numequation
	\label{eqn:S effective}
	\Spec S \to \cX_{d,\red}
\end{equation}
       	(\cite[\href{https://stacks.math.columbia.edu/tag/0DR1}{Tag 0DR1}]{stacks-project}), 
which is furthermore 
flat, by~\cite[\href{https://stacks.math.columbia.edu/tag/0DR2}{Tag 0DR2}]{stacks-project}.
	Let $x_0 = \Spec \F'$, thought of as the closed point of $\Spec S$,
	and endowed with a morphism $x_0 \to \cX_{d,\red}$
	corresponding to the representation~$\rhobar$.
	Since $\Spf S \to \cX_{d,\red}$
	is versal,
	it follows from the Artin Approximation Theorem that
	we may find a finite type $\F$-scheme $U$, equipped with a point
	$u_0$ with residue field $\F'$, such that the complete local
	ring of $U$ at $u_0$ is isomorphic to $S$,
	and 
	such that~\eqref{eqn:S effective} may be promoted 
	to  a smooth morphism
	\numequation
	\label{eqn:U promotion}
	U\to \cX_{d,\red},
\end{equation}
in the sense that on $u_0 (= \Spec \F')$,
the morphism~\eqref{eqn:U promotion}
induces the given morphism $x_0 (= \Spec \F') \to \cX_{d,\red}$,
and on $\Spec \widehat{\cO}_{U,u_0} \cong \Spec S$, 
the morphism~\eqref{eqn:U promotion}
induces the morphism~\eqref{eqn:S effective}.
	(See~\cite[\href{https://stacks.math.columbia.edu/tag/0DR0}{Tag 0DR0}]{stacks-project}.)


	It follows from Proposition~\ref{prop:versal rings},
	by pulling back over $\cX_{d,\red}$,
	that
	$\Spf S \times_{\cX_{d,\red}} \Spf S$
	is isomorphic to
        $\widehat{\GL}_{d,S},$
	a certain completion
	of~$(\GL_{d})_S$,  
	and thus that
	$\Spf S \times_{\cX_{d,\red}} x_0$
	is isomorphic to
        $\widehat{\GL}_{d,\F'},$
	a certain completion
	of~$(\GL_{d})_{\F'}$.  
	This latter fibre product may be identified
	with the completion of
	$\Spec S \times_{\cX_{d,\red}} x_0$ along
	its closed subspace $x_0\times_{\cX_{d,\red}} x_0,$
	and so we deduce that the morphism~\eqref{eqn:S effective}
	has relative dimension $d^2$ at the point $x_0$ in its domain.
	Thus the morphism~\eqref{eqn:U promotion} 
	has relative dimension $d^2$ at the point $u_0$ in its domain.
	Since this latter morphism is also smooth, its relative dimension
	is locally constant on its
	domain~\cite[\href{https://stacks.math.columbia.edu/tag/0DRQ}{Tag
    0DRQ}]{stacks-project},
       	and thus (since $S$
	is a local ring, so that $\Spec S$ is connected),
	we see that~\eqref{eqn:S effective}
	in fact has relative dimension~$d^2$.


	Let $\widetilde{M}$ denote the universal \'etale
	$(\varphi,\Gamma)$-module over $S$, and~$\widetilde{N}$ denote
        the base change of~$N$ to~$S$.
	Write $\widetilde{\Ext}^2 := H^2\bigl(\cC^{\bullet}(\widetilde{N}^\vee\otimes\widetilde{M})\bigr),$
	and 
write $$Y_r := \{x \in \Spec S \, | \, \dim \kappa(x)\otimes_S \widetilde{\Ext}^2
\geq r \}.$$ We claim that~$Y_r$ has dimension at most
$([K:\Q_p]d(d-1)/2) + d^2 - r. $ In order to see this, we note
that~$(Y_r)_{\Fpbar}$ is the pull-back to~$\Spec S$ via the flat morphism~\eqref{eqn:S effective}
of the locus considered in Theorem~\ref{thm:Xdred
  is algebraic}~(4) 
which has dimension at most~$[K:\Q_p]d(d-1)/2 - r. $
The required bound on the dimension of $Y_r$ follows, for example, from
  \cite[\href{https://stacks.math.columbia.edu/tag/02RE}{Tag
    02RE}]{stacks-project}.


	The composite
	$\Spf R/\varpi \to 
	\Spf R^{\square,\cO'}_{\rhobar}
	\to \cX_d$
	is effective
        by Corollary~\ref{cor: crystalline deformation rings are effective
    versal},
        i.e.\ arises from a morphism
	$\Spec R/\varpi \to \cX_d.$
	Thus the composite
	$\Spf (R/\varpi)_{\red} \to
	\Spf R/\varpi \to \cX_d$ 
	arises from a morphism
	$\Spec (R/\varpi)_{\red} \to \cX_d,$
	and hence factors through
	$\cX_{d,\red}$.  
	Consequently, we find that the surjection
	$R^{\square,\cO'}_{\rhobar} \to (R/\varpi)_{\red} $
	factors through $S$.
	The closed immersion $\Spec (R/\varpi)_{\red} \hookrightarrow \Spec S$
	then induces a closed immersion $\overline{X}_r \hookrightarrow Y_r$
	for each $r \geq 0$.  The upper bound on the dimension 
	of~ $Y_r$ computed in the preceding paragraph then provides
	the desired upper bound on the dimension of $\overline{X}_r$,
	completing the proof of the theorem.
%
%
%
%
%
\end{proof}


\section{Two geometric lemmas}
In this section we establish two lemmas which will be used in the proof
of Theorem~\ref{thm:crystalline lifts} below. 
We begin with the following simple lemma. 

\begin{lemma}
	\label{lem:torsion support}
	Suppose that $0 \to \cF \to \cG \to \cH \to 0$ is 
	a short exact sequence of coherent sheaves on a 
	reduced Noetherian scheme $X$, with $\cF$ and $\cG$ being locally free,
	and $\cH$ being torsion, in the sense that its support
	does not contain any irreducible component of $X$.
	Then there is an effective Cartier divisor~$D$ which contains
	the scheme-theoretic support of $\cH$, and whose 
	set-theoretic support coincides with the set-theoretic support
	of $\cH$,
	with the property that if $f: T \to X$
	is any morphism which meets $D$ properly,
	in the sense that the pull-back $f^*\cO_X \to f^*\cO_X(D)$ 
	of the canonical morphism $\cO_X \to \cO_X(D)$ is injective,
	then $\mathbb L_i f^*\cH = 0$ if $i  > 0$.
\end{lemma}
\begin{proof}
	Replacing $X$ by each of its finitely many connected components
	in turn, we see that it is no loss of generality to assume
	that $X$ is connected, and we do so; this ensures that
	each locally free sheaf on $X$ is of constant rank.
	If  $\eta$ is any generic point of $X$, then $\cH_{\eta} = 0$
	by assumption, and so $\cF_{\eta} \to \cG_{\eta}$ is an isomorphism.
	In particular, $\cF$ and $\cG$ are of the same rank,
        say~$r$, 
	and the induced morphism
	\numequation
	\label{eqn:wedge map}
	\wedge^r \cF \to \wedge^r \cG
\end{equation}
	is an injection of invertible sheaves. (Indeed, it follows
        from~\cite[III, \S7, Prop.\ 3]{MR1727844} that a map of
        locally free sheaves of rank~$r$ is injective if and only if
        the map on~$\wedge^r$ is injective.)  
	If we let $D$ denote the support of the cokernel of~(\ref{eqn:wedge
	       map}),	
        then $D$ is an effective Cartier divisor,
	and twisting by $\wedge^r \cF^{\vee}$
	identifies~(\ref{eqn:wedge map}) with the canonical section
	$\cO_X \hookrightarrow \cO_X(D)$.   The scheme-theoretic
	support of $\cH$ is contained in $D$, and the set-theoretic
	supports of $\cH$ and of $D$ coincide (this is a special
	case of the general fact that the Fitting ideal is
	contained in the annihilator, and has the same radical as it).

	Now suppose that $f:T \to X$
	meets $D$ properly in the sense described in the statement
	of the theorem.  Then we see that
	$$\wedge^r f^*\cF \to \wedge^r f^*\cG$$
	is injective, and thus that
	$f^*\cF \to f^*\cG$ is also injective. 
	Consequently
	$$0 \to f^*\cF \to f^*\cG \to f^*\cH \to 0$$ is
	short exact, and so $\mathbb L_i f^* \cH = 0$ for $i > 0,$
	as claimed.
\end{proof}

We now introduce some notation related to the second of the lemmas
(Lemma~\ref{lem:sections} below).

\begin{hyp}\label{hyp: codimension hypothesis}
  Let $X$ be a Noetherian scheme, and suppose that we have a short
  exact sequence of coherent sheaves \numequation
  \label{eqn:exact sequence}
  0 \to \cF \to \cG \to \cH \to 0,
\end{equation}
such that
\begin{itemize}
\item $\cG$ is locally free; 
\item $\cH$ 
  has the property that, for each $r\geq 1,$ the locus
  \[X_r:= \{x \in X \, | \, \dim \kappa(x) \otimes_{\cO_X} \cH \geq
    r\}\] is of codimension $\geq r+1$ in $X$.
\end{itemize}
\end{hyp}
In the setting of Hypothesis~\ref{hyp: codimension hypothesis}, if $\pi:\tX \to X$ is a morphism of Noetherian
schemes, then
we write $\tcG := \pi^* \cG$ and $\tcH := \pi^*\cH$,
and we let $\tcF$ denote the image of the pulled-back morphism $\pi^* \cF \to \pi^* \cG = \tcG$,
so that we again have a short exact sequence
\numequation
\label{eqn:tilde ses}
0 \to \tcF \to \tcG \to \tcH \to 0.
\end{equation}

\begin{remark} If $\cE$ is a locally free coherent sheaf on a 
	Noetherian scheme~ $X$, and $s:\cO_X \to \cE$ is a morphism,
	then we can think of $\cE$ as being the sheaf of sections
	of a vector bundle over $X$, and think of $s$ as being 
	a section of this vector bundle.  We can then speak of the zero locus
	of $s$; it is the closed subscheme of $X$ which locally,
        if we choose an isomorphism $\cE \cong \cO_X^r$ and write
$s = (a_1,\ldots,a_r),$ is cut out by the ideal sheaf $(a_1,\ldots,a_r)$;
more globally, it is cut out by the ideal sheaf which is the
image of the composite
$$\cE^{\vee}=\cO_X\otimes_{\cO_X}\cE^{\vee} \buildrel s\otimes \id \over \longrightarrow 
\cE\otimes_{\cO_X}\cE^{\vee} \to \cO_X,$$
where the second arrow is the canonical pairing.

We note that $s$ is nowhere-vanishing, i.e.\ the zero locus of $s$ is empty,
if and only if  $s$ is a split injection (so that
$s$ gives rise to a morphism of vector bundles, rather than merely
of sheaves).
We also note that
if $\cE$ has rank $r$, then the Hauptidealsatz shows that
the zero locus of a section of $\cE$ has codimension at most $r$ around
each of its points. 
\end{remark}

\begin{lem}
	\label{lem:sections}
	Suppose that we are in the setting of Hypothesis~{\em \ref{hyp:
          codimension hypothesis}}, that the $\pi: \tX \to X$ is a surjective
morphism whose domain is a Noetherian scheme,
and that $\tcF$ is locally
	free.  Suppose further that $\cO_X \to \cF$ is a morphism with
	the property that the pulled-back morphism $\cO_{\tX} \to \tcF$
is a nowhere-vanishing section.
	Then the composite 
	$\cO_X \to \cF \to \cG$ is nowhere-vanishing. 
\end{lem}
\begin{proof} The statement may be checked by working
	locally at each of the points of~$X$; more precisely, we may replace
$X$ by $\Spec \cO_{X,x}$ and $\tX$ by its pull-back over $\Spec \cO_{X,x}$.  
We have the stratification
	of $X$ by the closed subsets $X_r$ (where in the
case $r = 0$, we declare $X_0:= X$),
and if $x \in X$, then we
	set $r(x) := \dim \kappa(x)\otimes_{\cO_X} \cH$, or equivalently,
	the maximal value of $r$ for which $x \in X_r$; we then prove the
	theorem by induction on $r(x)$.  In fact, we assume that 
        the statement is true  after localizing at points $x$ for which $r(x) < r$ for {\em all} exact sequences~(\ref{eqn:exact
		sequence}) which satisfy the hypotheses of the lemma, and prove it for our given
	exact sequence after localizing at an $x$ for which $r(x) =
        r$. 

In the case $r = 0$, the sheaf $\cH$ is zero, 
so that $\cF \iso \cG$. 
Thus the induced morphism $\cO_{\tX} \to \tcG$
gives a nowhere-vanishing section.  Since $\pi$ is surjective, this ensures
that the induced morphism $\cO_X \to \cG$ is also nowhere-vanishing, as required.

We now consider the case $r \geq 1$.
	Since we have replaced $X$ by $\Spec \cO_{X,x}$,
we may assume that $\cG$ is free.
By Nakayama's lemma and our assumption
        that~$r(x)=r$, we may choose a surjection
%
%
%
%
%
	\numequation
	\label{eqn:surjection}
        \cO_X^r \to \cH,
\end{equation}
	which we may then lift to a morphism $\cO_X^r \to \cG$.
        If we let $\cK$ denote the kernel of~(\ref{eqn:surjection}), 
then we obtain a morphism of short exact sequences
\numequation
\label{eqn:K to F map}
\xymatrix{
0 \ar[r] & \cK \ar[r] \ar[d] & \cO_X^r \ar[r]\ar[d] & \cH \ar[r]\ar@{=}[d] & 0
\\ 
0 \ar[r] & \cF  \ar[r] & \cG \ar[r] & \cH \ar[r] & 0
}
\end{equation}
Pulling back this diagram back via $\pi$, and letting~$\tcK$ denote
the image of~$\pi^*\cK$ in~$\cO_{\tX}^r$,
we  obtain a corresponding morphism of short exact sequences
\numequation
\label{eqn:tK to tF map}
\xymatrix{
0 \ar[r] & \tcK \ar[r] \ar[d] & \cO_{\tX}^r \ar[r]\ar[d] & \tcH \ar[r]\ar@{=}[d] & 0
\\ 
0 \ar[r] & \tcF  \ar[r] & \tcG \ar[r] & \tcH \ar[r] & 0
}
\end{equation}
The morphism of short exact sequences~\eqref{eqn:K to F map}
induces a short exact sequence
\numequation
\label{eqn:K to F plus O to G}
	0 \to \cK \to \cF \oplus \cO_X^r \to \cG \to 0,
\end{equation}
	which is furthermore split (since $\cG$ is free and we are over an
	affine scheme), so that we may write
\numequation
\label{eqn:direct sum iso}
	 \cF \oplus \cO_X^r \cong \cK \oplus \cG.
\end{equation}
Correspondingly, the short exact sequence~\eqref{eqn:tK to tF map}
induces a short exact sequence
	$$0 \to \tcK \to \tcF \oplus \cO_{\tX}^r \to \tcG \to 0,$$
which can be thought of as being obtained from~\eqref{eqn:K to F plus O to G}
by pulling back via $\pi$ and then taking the quotient by the (naturally identified)
copies of $\mathbb L_1 \pi^*\cH$ sitting inside $\pi^* \cK$ and $\pi^* \cF$.   
Since~\eqref{eqn:K to F plus O to G} is split, so is this latter
short exact sequence, and so we also obtain an isomorphism
	\numequation
	\label{eqn:pulled-back iso}
	\tcF \oplus \cO_{\tX}^r \cong \tcK \oplus \tcG
\end{equation}
which can be thought of as being obtained via by pulling back the
isomorphism~\eqref{eqn:direct sum iso} along~$\pi,$
and then taking the quotient of each side by the appropriate
copy of $\mathbb L_1 \pi^* \cH$.
Since all the other summands appearing in the isomorphism~\eqref{eqn:pulled-back
iso} are locally free, we see that the same is true of~$\tcK$.

        
	Suppose, by way of obtaining a contradiction, that 
	the composite $\cO_X \to \cF \to \cG$
	is not nowhere-vanishing;
	then we see that its image lies in $\mathfrak m_x \cG.$
        Thus the pulled back morphism $\cO_{\tX} \to \tcG$ 
        vanishes at each point $\tx$ lying over $x$ (and 
there is at least one such point, since $\pi$ is surjective by assumption).
	Since the pulled-back morphism $\cO_{\tX} \to \tcF$
	{\em is} nowhere-vanishing,
	so is the induced morphism
        $\cO_{\tX}\to \tcF \oplus \cO_{\tX}^r.$
	A consideration of~(\ref{eqn:pulled-back iso}) then shows
	that the induced morphism $\cO_{\tX} \to \tcK$,
	which is pulled-back from the induced morphism $\cO_X \to \cK,$
	must be nowhere-vanishing.

	Thus, if we consider the exact sequence
	$$0 \to \cK \to \cO_X^r \to \cH \to 0,$$
	it satisfies Hypothesis~\ref{hyp: codimension hypothesis}, 
	and we are given a morphism $\cO_X\to \cK$ whose pull-back to $\tX$
	is a nowhere-vanishing section of $\tcK$.
        Applying the inductive hypothesis to this situation, we find 
	that the composite $\cO_{X} \to \cO_{X}^r$ is a section
	of a rank $r$ locally free sheaf
	whose zero locus is contained
        in~$X_r$ (because if we localise at a point not in~$X_r$, then
        the section is nowhere vanishing by the inductive hypothesis), and  is therefore of codimension
	at least $r+1$ by hypothesis.  This zero locus contains $x$ (since this
	section factors through $\cK$, and $\cH \cong \cO_X^r/\cK$
	has fibre dimension exactly $r$ at $x$, so the map
        $\cO_X^r\otimes_{\cO_X}\kappa(x)\to\cH\otimes_{\cO_X}\kappa(x)$
        is an isomorphism), and
	hence is non-empty. 
	On the other hand, as we noted above,
	the zero locus of a section of a rank $r$ locally free sheaf,
	if it is non-empty, has codimension at most~$r$. 
	This contradiction completes the proof of the lemma.
\end{proof}

\section{Crystalline lifts} 

		
Given a~$d$-tuple of labeled Hodge--Tate weights~$\underline{\lambda}$, and
a~$d'$-tuple of labeled Hodge--Tate weights~$\underline{\lambda}'$, we say
that~$\underline{\lambda}'$ is \emph{slightly greater than~$\underline{\lambda}$} (and that $\underline{\lambda}$ is
\emph{slightly less than~$\underline{\lambda}'$}) if for each
$\sigma:K\into\Qpbar$ we have $\lambda'_{\sigma,d}\ge \lambda_{\sigma,1}+1$, and
the inequality is strict for at least one~$\sigma$. 

This is not standard terminology, but it will be convenient for us; it
is motivated by the following well-known result. Here, as throughout
the book, we
slightly abuse terminology and refer to lattices in crystalline
representations as themselves being crystalline. 
\begin{lem}
  \label{lem: extensions are crystalline}Let
  $0\to\rho_1^\circ\to\rho^\circ\to\rho_2^\circ\to 0$ be an extension of
  $\Zpbar$-valued representations of~$G_K$, with~$\rho_1^\circ$
  and~$\rho_2^\circ$ being crystalline of Hodge--Tate
  weights~$\underline{\lambda}_1$, $\underline{\lambda}_2$
  respectively. If~$\underline{\lambda}_1$ is slightly less
  than~$\underline{\lambda}_2$, then~$\rho^{\circ}$ is 
  crystalline.
\end{lem}
\begin{proof}This follows easily from the formulae in~\cite[Prop.\
  1.24]{MR1263527}. Indeed, if we let~$\rho_1$, $\rho_2$ be
  the~$\Qpbar$-representations corresponding to~$\rho_1^\circ$,
  $\rho_2^\circ$, and we set $V:=\rho_1\otimes\rho_2^\vee$, then we
  need to show that $h^1_f(V)=h^1(V)$. Since the Hodge--Tate weights
  of~$V$ are all negative, and for at least one
  embedding~$\sigma:K\into\Qpbar$ the $\sigma$-labeled Hodge--Tate
  weights are all less than~$-1$, we have
  $h^1(V)-h^1_f(V)=h^2(V)=h^0(V^\vee(1))=0$, as required.
\end{proof}

Arguing inductively,
the following theorem will allow us to construct 
crystalline lifts of any given $\rhobar$. 
In particular it implies Theorem~\ref{thm:intro crystalline lifts} from
the introduction, in the more refined form of Theorem~\ref{thm: strong
  existence of crystalline lifts} below.
\begin{theorem}
	\label{thm:crystalline lifts}
	Suppose given a representation $\rhobar_d: G_K \to \GL_d(\Fbar_p)$
	that admits a lift
	$\rho_d^{\circ} : G_K \to \GL_d(\Zbar_p)$ which 
	is crystalline with labeled Hodge--Tate weights
        $\underline{\lambda}$. 
	Let $0 \to \rhobar_d \to \rhobar_{d+a} \to \alphabar \to 0$
	be any extension of $G_K$-representations over $\Fbar_p$,
	with $\alphabar: G_K \to \GL_a(\Fbar_p)$ irreducible,
	and let $\alpha^{\circ}: G_K \to \GL_a(\Zbar_p)$ be any
        crystalline lifting
	of $\alphabar$ with labeled Hodge--Tate
        weights~$\underline{\lambda}'$, which we assume to be slightly greater than~$\underline{\lambda}$.

	Then we may find a lifting of the given extension to an extension
	$$0 \to \theta^{\circ}_d \to \theta_{d+a}^{\circ} \to \alpha^{\circ} \to 0$$
	of $G_K$-representations over $\Zbar_p$, where
        $\theta^{\circ}_d: G_K \to \GL_d(\Zbar_p)$ again has the
        property that the associated $p$-adic representation
        $\theta_d: G_K \to \GL_d(\Qbar_p)$ is crystalline with
        labeled Hodge--Tate weights $\underline{\lambda}$.
        Furthermore, $\theta_{d+a}^{\circ}$ is crystalline, and we may
        choose $\theta_d^{\circ}$ to lie on the same irreducible
        component of
        $\Spec R_{\rhobar_d}^{\crys,\underline{\lambda}}$ that
        $\rho^{\circ}_d$ does.
\end{theorem}

\begin{proof}[Proof of Theorem~{\em \ref{thm:crystalline lifts}}]
We may and do choose our field of coefficients $E$ to be large enough
that the various Galois representations are defined over $\cO$ or $\F$
as the case may be.
	We write $R := R_{\rhobar}^{\crys,\underline{\lambda}}$, 
	and we let $X$ denote the component of the crystalline lifting scheme 
	$\Spec R_{\rhobar}^{\crys,\underline{\lambda}}$ on which $\rho^{\circ}_d$
	lies. As recalled in Section~\ref{subsec: notation and
          conventions}, $X$ is reduced. 
	Over $X$ we have a good complex $C^{\bullet}$ supported in degrees~$[0,2]$ computing
	$\Ext^{\bullet}_{G_K}(\alpha^{\circ},\rho^{\circ})$ for the
        universal deformation~$\rho^{\circ}$; indeed, by
a straightforward variant of Corollary~\ref{cor:Herr complex length},
we can choose a good complex
        supported in degrees~$[0,2]$ and quasi-isomorphic to the Herr
        complex for~$\rho^\circ\otimes(\alpha^{\circ})^\vee$.  Our assumption
	on the Hodge--Tate weights 
        implies that $\Ext^2$ is
	supported (set-theoretically) on the special fibre $\overline{X}$
	of $X$. Indeed, the formation of~$\Ext^2$ is compatible with
        base change, and at any closed point~$x$ of the generic fibre
        of~$R$, we have
        $\Ext^{2}_{G_K}(\alpha^{\circ},\rho^{\circ}_x)=\Hom_{G_K}(\rho^{\circ}_x,\alpha^{\circ}(1))$
        by Tate local duality, and this space vanishes by the
        assumption that~$\underline{\lambda}'$ is slightly greater than~$\underline{\lambda}$.

	As usual, we let $B^2$ denote the image of $C^1$ in $C^2$; 
	it is a subsheaf of the locally free sheaf $C^2$, and so is 
	torsion free.\footnote{Since $X$ is reduced, {\em torsion free}
		is a reasonable notion for a finitely
		generated module, equivalent to the associated
		primes lying among the minimal primes.}
	By~ \cite[\href{http://stacks.math.columbia.edu/tag/0815}{Tag
          0815}]{stacks-project}, we may find a blow up $\pi:\tX \to X$, whose centre lies
	(set-theoretically) in the special fibre of $X$, 
	such that the torsion-free quotient of $\pi^*B^2$ becomes
	locally free (recall that~$\Ext^2$ is supported on the special
        fibre of~$X$); in other words, if we let $\tC^{\bullet}$
	denote the pull-back by $\pi$ of $C^{\bullet}$, 
	the corresponding $\cO_{\tX}$-module of $2$-coboundaries
	$\tB^2$ (that is, the image of~$\tC^1$ in~$\tC^2$)
	is locally free. Thus, if we let $\tZ^1$ denote the $\cO_{\tX}$-module of
	$1$-cocycles for $\tC^{\bullet}$, then $\tZ^1$ is also locally
	free; so in particular, the complex $\tC^0 \to \tZ^1$ is a
        good complex. 

	The given extension $\rhobar_{d+a}$ is classified by a
	class in $\Ext^1_{G_K}(\alphabar,\rhobar_d)$, which is to say, a class
	in $H^1$ of the complex $\kappa\otimes C^{\bullet}$.
	(Here we write $\kappa$ to denote $\F$ thought of as the
	residue field of $R$.)
        Lifting this class to an element of $\kappa \otimes C^1,$
	and then to an element of $C^1$, we find ourselves in the 
	following situation:
	we have a morphism $c:R \to C^1$, whose image under the coboundary
	lies in $\mathfrak m_R C^2$ (reflecting the fact that we have a cochain
	which becomes a cocycle at the closed point). 

	Consider the composite $b: R \buildrel c \over \longrightarrow
	C^1 \to B^2,$ 
	which pulls back to a section $\tb: \cO_{\tX} \to \tB^2.$
	It follows from Lemma~\ref{lem:sections} above (applied with
        $\cF=B^2$, $\cG=C^2$, and $\cH=\Ext^2$),
	together with Theorem~\ref{thm: codimension H2 in regular
          weight crystalline deformation ring}, 
	that the lifted section~$\tb$ must have non-empty zero locus,
	which (being closed) must contain a point $\tx$ lying
	over the closed point ~$x\in X$.   The section $c$ itself
	pulls back to a section $\tc: \cO_{\tX} \to \tC^1$, whose value at the 
	point $\tx$ lies in the fibre of $\tZ^1$.   In other
	words, the fibre of $\tc$ at $\tx$ is a one-cocycle
	in the complex $\kappa(\tx)\otimes [\tC^0 \to \tZ^1]$, giving rise to
	a class
	$\overline{e} \in H^1\bigl(\kappa(\tx)\otimes [\tC^0 \to \tZ^1]\bigr)$
	lifting the original class in $\Ext^1_{G_K}(\alphabar,\rhobar_d)$ that classifies~
	$\rhobar_{d+a}$.

	Since $\tX$ is flat over $\Z_p$ (because~$X$ is), 
we may find a morphism
	$\tif:\Spec \Zbar_p \to \tX$
	passing through the point $\tx \in \tX(\Fbar_p)$ (in the sense
        that the closed point of $\Spec\Zbar_p$ maps to~$\tx$). 
	The morphism $\tif$ determines (and, by the valuative
	criterion of properness, is determined by) a morphism 
	$f: \Spec \Zbar_p \to X$,
	lifting the closed point $x \in X$, 
	which in turn corresponds to a representation
	$\theta_d^{\circ}: G_K \to \GL_d(\Zbar_p),$ as in the statement
	of the theorem.
	Now $H^2(\tC^{\bullet})$ is the cokernel of an inclusion of
	locally free sheaves (the inclusion $\tB^2 \hookrightarrow
        \tC^2$), and it is torsion (because it is pulled back
        from~$H^2(C^{\bullet})$, which is set-theoretically supported in the special fibre)
	and so Lemma~\ref{lem:torsion support} above shows that
	there is an effective Cartier divisor $D$ contained (set-theoretically)
	in the
        special 
	fibre of $\tX$ with the property that for any morphism
	to $\tX$ that meets $D$ properly, the higher derived
	pull-backs of $H^2(\tC^{\bullet})$ under this morphism
	vanish.  Since the domain of $\tif$ is $p$-torsion free,
	it meets the special fibre of $\tX$ properly,
	and in particular meets $D$ properly; 
	thus we infer that
	$$\mathbb L_i\tif^* H^2(\tC^{\bullet}) = 0$$
	for $i > 0.$ From this vanishing we deduce 
	that 
	\begin{multline*}
		\Ext^1_{G_K}(\alpha^{\circ},\theta_d^{\circ}) = H^1(f^*C^{\bullet}) =
	H^1(\tif^*\tC^{\bullet}) = \tif^*H^1(\tC^{\bullet})
	\\
	= \tif^*H^1\bigl( [ \tC^0 \to \tZ^1]) = H^1\bigl(\tif^*[\tC^0 \to
	\tZ^1]\bigr)
	.\end{multline*}
	(The first of these identifications is the general base-change
	property of the perfect complex $C^{\bullet}$, 
      the second
	follows from the fact that $f = \pi\circ \tif$,  the third
	follows from the vanishing of $\mathbb L_i\tif^* H^2(\tC^{\bullet}) $
        for $i>0$ and the base-change spectral sequence
\[E_2^{p,q} := L_{-p}\tif^*H^q(\tC^\bullet)\implies H^{p+q}(\tif^*\tC^\bullet) \] 
(see e.g.~\cite[\href{https://stacks.math.columbia.edu/tag/0662}{Tag 0662}]{stacks-project}),
	the fourth holds by definition of $H^1$, and the fifth follows
	from right-exactness of~$\tif^*$.)

	Since $\tif$ identifies the closed point of $\Spec \Zbar_p$
	with the point $\tx,$ we may find a class
	$e \in H^1\bigl(\tif^*[\tC^0 \to \tZ^1]\bigr)$ lifting the
	class $\overline{e}$, which is then identified with an element
	of $\Ext^1_{G_K}(\alpha^{\circ},\theta_d^{\circ}).$   This element classifies
	an extension $0 \to \theta_d^{\circ} \to \theta_{d+a}^{\circ} 
	\to \alpha^{\circ} \to 0$ which by construction lifts the
        given extension $0 \to \rhobar^{\circ}_d \to \rhobar_{d+a}^{\circ} \to \alphabar^{\circ} \to 0$, and so by Lemma~\ref{lem: extensions are crystalline} satisfies the conclusions of the theorem.
\end{proof}

\section{Potentially diagonalizable crystalline
  lifts}\label{subsec: precise crystalline lifts}We now use
Theorem~\ref{thm:crystalline lifts} to prove Theorem~\ref{thm:intro
  crystalline lifts}.  There are at least two differences between the
proof of Theorem~\ref{thm:intro crystalline lifts} and previous work
on the problem (in particular the results
of~\cite{MullerThesis,2015arXiv150601050G}.) One is that, rather than
working only with lifts which are extensions of inductions of characters,
we utilise extensions of more general potentially diagonalizable
representations. The other difference
is that our argument exploits our control
of the support of~$H^2$. In previous work, the only tool used to deal
with classes in~$H^2$ was twisting the various irreducible
representations by unramified characters; at points where~$H^2$ has
dimension greater than~$1$, this seems to be insufficient.

In fact, we prove various refinements of this result, allowing us to
produce potentially diagonalizable lifts; such lifts are important for
automorphy lifting theorems. We can also control the Hodge--Tate
weights of these potentially diagonalizable lifts; for example, we can
insist that the gaps between the weights are arbitrarily large, which
is useful in applications to automorphy lifting. Alternatively, we can
produce lifts whose Hodge--Tate weights lift some Serre weight,
proving a conjecture which is important for the formulation of general
Serre weight conjectures (see the discussion after~\cite[Rem.\
5.1.8]{2015arXiv150902527G}).

We begin by recalling some definitions and some basic lemmas about
extensions of crystalline representations.  Let~$\underline{\lambda}$ be a
regular $d$-tuple of Hodge--Tate weights.




\begin{defn}\label{defn: ordinary} \index{ordinary}
  A crystalline representation of weight~$\underline{\lambda}$
  is~\emph{ordinary} if it has a $G_K$-invariant decreasing 
  filtration
  whose associated graded pieces are all one-dimensional, such that
  the $\sigma$-labeled Hodge--Tate weight of the~$i$th graded piece
  is~$\lambda_{\sigma,i}$. In other words, we can write the
  representation in the form \[\begin{pmatrix}
      \chi_1 &*&\dots &*\\
      0& \chi_2&\dots &*\\
      \vdots&& \ddots &\vdots\\
      0&\dots&0& \chi_d\\
    \end{pmatrix} \]where~$\chi_i$ is a crystalline
  character whose $\sigma$-labeled Hodge--Tate weight
  is~$\lambda_{\sigma,d+1-i}$. 
\end{defn}

If~$\rho_1^\circ,\rho_2^\circ:G_K\to\GL_d(\Zpbar)$ are two crystalline
representations, then we write $\rho_1^\circ\sim\rho_2^\circ$, and say
that $\rho_1^\circ$ \emph{connects to} $\rho_2^\circ$, if and
only if the following conditions hold: $\rho_1^\circ$ and~$\rho_2^\circ$ have the same labeled
Hodge--Tate weights, have isomorphic reductions modulo~$\m_{\Zpbar}$, and
determine points on the same irreducible component of the
corresponding crystalline lifting rings. The following definition
was originally made in~\cite[\S 1.4]{BLGGT}.
\begin{defn}\label{defn: pot diag}\index{potentially diagonalizable}
  A representation~$\rho^\circ:G_K\to\GL_d(\Zpbar)$ is
  \emph{potentially diagonalizable} if there is a finite
  extension~$K'/K$ and crystalline
  characters~$\chi_1,\dots,\chi_d:G_{K'}\to\overline{\Z}_p^\times$
  such that~$\rho^\circ|_{G_{K'}}\sim\chi_1\oplus\cdots\oplus\chi_d$.
\end{defn}

\begin{lem}
  \label{lem: extensions of PD are PD}Let
  $0\to\rho_1^\circ\to\rho^\circ\to\rho_2^\circ\to 0$ be an extension of
  $\Zpbar$-valued representations of~$G_K$. If $\rho^\circ$ is
  crystalline, and~$\rho_1^\circ$ and~$\rho_2^\circ$ are potentially
  diagonalizable, then~$\rho^\circ$ is also potentially
  diagonalizable.
\end{lem}
\begin{proof}We may choose $K'/K$ such that $\rho_1^\circ|_{G_{K'}}$
  and $\rho_2^\circ|_{G_{K'}}$ both connect to direct sums of crystalline
  characters, and such that $\rhobar|_{G_{K'}}$ is trivial. It then
  follows from points~(5) and~(7) of the list before~\cite[Lem.\
  1.4.1]{BLGGT} that~$\rho^\circ|_{G_{K'}}$ connects to the direct sum
  of the union of the sets of crystalline characters for $\rho_1^\circ|_{G_{K'}}$
  and $\rho_2^\circ|_{G_{K'}}$, as required.
\end{proof}
Note that it follows in particular from Lemma~\ref{lem: extensions of
  PD are PD} that ordinary representations are potentially diagonalizable.


\begin{thm}
  \label{thm: strong existence of crystalline lifts}Let $K/\Qp$ be a
  finite extension, and let~$\rhobar:G_K\to\GL_d(\Fpbar)$ be a
  continuous representation. Then~$\rhobar$ admits a lift to a
  crystalline representation~$\rho^\circ:G_K\to\GL_d(\Zpbar)$ of some
  regular labeled Hodge--Tate weights~$\underline{\lambda}$. Furthermore:
  \begin{enumerate}
  \item $\rho^\circ$ can be taken to be potentially diagonalizable.
  \item If every Jordan--H\"older factor of~$\rhobar$ is
    one-dimensional, then~$\rho^\circ$ can be taken to be
    ordinary.
  \item $\rho^\circ$ can be taken to be potentially diagonalizable,
    and~$\underline{\lambda}$ can be taken to be a lift of a Serre weight.
  \item $\rho^\circ$ can be taken to be potentially diagonalizable,
    and~$\underline{\lambda}$ can be taken to have arbitrarily spread-out
    Hodge--Tate weights: that is, for any $C>0$, we can
    choose~$\rho^\circ$ such that for each~$\sigma:K\into\Qpbar$, we
    have $\lambda_{\sigma,i}-\lambda_{\sigma,i+1}\ge C$ for each $1\le i\le n-1$.
  \end{enumerate}

\end{thm}
\begin{proof}We prove all of these results by induction on~$d$, using
  Theorem~\ref{thm:crystalline lifts}. We begin by proving~(1), and
  then explain how to refine the proof to give each of (2)-(4).

Suppose firstly that~$\rhobar$ is irreducible (this is the base case
of the induction). Then~$\rhobar$ is of
the form~$\Ind_{G_{K'}}^{G_K}\psibar$ for some
character~$\psibar:G_{K'}\to\Fpbartimes$, where $K'/K$ is unramified
of degree~$d$. We can choose (for example by~\cite[Lem.\
7.1.1]{2015arXiv150902527G}) a crystalline
lift~$\psi:G_{K'}\to\overline{\Z}_p^\times$ of~$\psibar$ such
that~$\rho^\circ:=\Ind_{G_{K'}}^{G_K}\psi$ has regular Hodge--Tate
weights. Then~$\rho^\circ$ is crystalline (because $K'/K$ is
unramified), and it is potentially diagonalizable,
because~$\rho^\circ|_{G_K'}$ is a direct sum of crystalline
characters.

For the inductive step, we may therefore suppose that we can
write~$\rhobar$ as an extension
$0\to\rhobar_1\to\rhobar\to\rhobar_2\to 0$, with~$\rhobar_2$
irreducible, and we may assume that we have regular potentially
diagonalizable lifts~$\rho_1^\circ$, $\rho_2^\circ$ of~$\rhobar_1$,
$\rhobar_2$ respectively. Twisting~$\rho_1^\circ$ by a crystalline
character with sufficiently large Hodge--Tate weights and trivial
reduction modulo~$p$ (which exists by~\cite[Lem.\
7.1.1]{2015arXiv150902527G}), we may suppose that the Hodge--Tate
weights of~$\rho_1^\circ$ are slightly less than those
of~$\rho_2^\circ$. By Theorem~\ref{thm:crystalline lifts}, there is a
lift~$\theta_1^{\circ}$ of~$\rhobar_1$
with~$\theta_1^\circ\sim\rho_1^\circ$, and a lift~$\rho^\circ$
of~$\rhobar$ which is an
extension \[0\to\theta_1^\circ\to\rho^\circ\to\rho_2^\circ\to 0.\]
Since~$\rho_1^\circ$ is potentially diagonalizable, so
is~$\theta_1^\circ$, and it follows from Lemmas~\ref{lem: extensions
  are crystalline} and~\ref{lem: extensions of PD are PD}
that~$\rho^\circ$ is crystalline and potentially diagonalizable. This
completes the proof of~(1).

We claim that if every Jordan--H\"older factor of~$\rhobar$ is
    one-dimensional, then the lift that we produced in proving~(1) is
    actually automatically ordinary. Examining the proof, we see that
    it is enough to show that if~$\rho_1^\circ$ is ordinary
    and~$\theta_1^\circ\sim\rho_1^\circ$ then $\theta_1^\circ$ is also
    ordinary; this is immediate from~\cite[Lem.\ 3.3.3(2)]{ger}.

    To show~(3) and~(4) we have to check that we can choose the
    Hodge--Tate weights of the lifts of the irreducible pieces
    of~$\rhobar$ appropriately. This requires substantially more
    effort in the case of~(3), but that effort has already been made
    in the proof of~\cite[Thm.\ B.1.1]{2015arXiv150902527G} (which
    proves the case that~$\rhobar$ is semisimple). Indeed, examining
    that proof, we see that it produces lifts satisfying all of the
    conditions we need; we need only note that the ``slightly less''
    condition is automatic, except in the case that for each~$\sigma$
    we have (in the notation of the proof of ~\cite[Thm.\
    B.1.1]{2015arXiv150902527G}) $h_\sigma+x_\sigma=H_\sigma+1$, in
    which case we can take $h_\sigma+x_\sigma=H_\sigma+p$ for
    each~$\sigma$ instead.

Finally, to prove~(4) it is enough to check the case that~$\rhobar$ is
irreducible (because we can choose the Hodge--Tate weights of the
character that we twisted~$\rho_1^\circ$ by to be arbitrarily large). For this, note that in our application of~\cite[Lem.\
7.1.1]{2015arXiv150902527G}, we can change any labeled Hodge--Tate by
any multiple of~$(p^d-1)$, so we can certainly arrange that the gaps
between labeled Hodge--Tate weights are as large as we please.
\end{proof}

\begin{rem}
  \label{rem: bound on weights}To complete the proof of
  Theorem~\ref{thm:intro crystalline lifts}, it is enough to note that
  by definition the lifts produced in Theorem~\ref{thm: strong
    existence of crystalline lifts}~(3) have all their Hodge--Tate
  weights in the
  interval~$[0,dp-1]$. 
\end{rem}
\begin{rem}
  \label{rem: get all lifts of a Serre weight}The proof of
  Theorem~\ref{thm: strong existence of crystalline lifts}~(3)
  actually proves the slightly stronger statement that there is a
  Serre weight such that if~$\underline{\lambda}$ is any lift of that
  Serre weight, then~$\rhobar$ admits a potentially diagonalizable
  crystalline lift of weight~$\underline{\lambda}$. (Note that in the
  case that~$K/\Qp$ is ramified, there may be many such choices of a
  lift of a given Serre weight, because by definition,
  the choice of a lift of $\underline{\lambda}$ depends on a choice of
  a lift of each embedding~$k\into\Fpbar$ to an embedding
  $K\into\Qpbar$. However, it is noted in the
  first paragraph of the proof of~\cite[Thm.\
  B.1.1]{2015arXiv150902527G} that in the case that~$\rhobar$ is
  semisimple, there is a crystalline lift of
  weight~$\underline{\lambda}$ for any choice of
  lift~$\underline{\lambda}$, so the same is true in our construction.)  
\end{rem}

As a consequence of Theorem~\ref{thm: strong existence of crystalline
  lifts}, we can remove a hypothesis made in~\cite[App.\
A]{emertongeerefinedBM}, proving in particular the  following result
on the existence of globalizations of local Galois representations.
\begin{cor}
  \label{cor: existence of global lift}Suppose that~$p\nmid 2d$. Let $K/\Qp$ be a finite extension, and let
  $\rhobar:G_K\to\GL_n(\Fpbar)$ be a continuous representation. Then
  there is an imaginary CM field $F$ and a continuous irreducible
  representation $\rbar:G_{F}\to\GL_n(\Fpbar)$ such that
  \begin{itemize}
  \item each place $v|p$ of~$F$ satisfies $F_v\cong K$,
  \item for each place $v|p$ of $F$, either $\rbar|_{G_{F_v}}\cong
    \rhobar$ or $\rbar|_{G_{F_{v^c}}}\cong   \rhobar$, and
   \item $\rbar$ is automorphic, in the sense that it may be lifted
     to a representation $r:G_{F}\to\GL_n(\Qpbar)$ coming from a
     regular algebraic conjugate self dual cuspidal automorphic
     representation of~$\GL_n/F$.
  \end{itemize}
\end{cor}
\begin{proof}
  This is immediate from~\cite[Cor.\ A.7]{emertongeerefinedBM}
  (since~\cite[Conj.\ A.3]{emertongeerefinedBM} is a special case of
  Theorem~\ref{thm: strong existence of crystalline lifts}).
\end{proof}
\section{The irreducible components of $\cX_{\lowercase{d,\red}}$}\label{subsec:
  irred components}
We now complete the analysis of the irreducible
components of~$\cX_{d,\red}$ that we began in Section~\ref{subsec: formal
  algebraic stack}. Recall that Theorem~\ref{thm:Xdred is algebraic}
shows that~$\cX_{d,\red}$ is an algebraic stack of finite presentation
over~$\F$, and has dimension~$[K:\Qp]d(d-1)/2$. Furthermore, for each
Serre weight~$\underline{k}$, there is a corresponding irreducible
component~$\cX_{d,\red,\Fpbar}^{\underline{k}}$
of~$(\cX_{d,\red})_{\Fpbar}$, and the components for
different weights~$\underline{k}$ are distinct (see Remark~\ref{rem: different k give different components}).

We are now finally in a position to prove the following
result.
\begin{thm}
	\label{thm:reduced dimension}  $\cX_{d,\red}$ is
        equidimensional of dimension
	$[K:\Q_p] d(d-1)/2$, and the irreducible components of~$(\cX_{d,\red})_\Fpbar$ are
        precisely the various closed substacks~$\cX_{d,\red,\Fpbar}^{\underline{k}}$
	of~Theorem~{\em \ref{thm:Xdred is algebraic}}; in particular,
        $(\cX_{d,\red})_\Fpbar$ is maximally nonsplit of niveau~$1$.  Furthermore
        each~$\cX_{d,\red,\Fpbar}^{\underline{k}}$ can be defined
        over~$\F$, i.e.\ is the base change
        of an irreducible component~$\cX_{d,\red}^{\underline{k}}$ of~$\cX_{d,\red}$.
\end{thm}
\begin{proof}
     By Theorem~\ref{thm:Xdred is algebraic}, it is enough to prove
     that each irreducible component of 
     $(\cX_{d,\red})_\Fpbar$ is of dimension of at least
     $[K:\Q_p] d(d-1)/2$. Indeed, this shows that the irreducible
     components of~$(\cX_{d,\red})_\Fpbar$ are
        precisely the~$\cX_{d,\red,\Fpbar}^{\underline{k}}$; to see
        that these may all be defined over~$\F$, we need to show that
        the action of~$\Gal(\Fpbar/\F)$ on the irreducible components
        of~$(\cX_{d,\red})_\Fpbar$ is trivial. This follows
        immediately by considering its action on the maximally
        nonsplit representations of niveau~$1$ (since the action
        of~$\Gal(\Fpbar/\F)$ preserves the property of being maximally
        nonsplit of niveau~$1$ and weight~$\underline{k}$).

In order to see that each  irreducible component of 
     $(\cX_{d,\red})_\Fpbar$ is of dimension of at least
     $[K:\Q_p] d(d-1)/2$, it suffices to show that every finite type
     point~$x$     of $(\cX_{d,\red})_\Fpbar$ is contained in an irreducible
     substack of $(\cX_{d,\red})_\Fpbar$ of dimension at least
     $[K:\Q_p] d(d-1)/2$; so it suffices in turn to show that each~$x$
     is contained in an equidimensional substack of dimension
     $[K:\Q_p] d(d-1)/2$. But by Theorem~\ref{thm: strong existence of
       crystalline lifts}, there is a regular Hodge
     type~$\underline{\lambda}$ such that $x$ is contained in
     $\cX_d^{\crys,\underline{\lambda}}\times_{\Spf \cO}\Spec\F$
     (for~$\cO$ sufficiently large), and this stack is equidimensional
     of dimension $[K:\Q_p] d(d-1)/2$ by Theorem~\ref{thm: dimension
       of ss stack}, as required.\end{proof}
We end this section with the following result, showing that in
contrast to its substacks~$\cX_{d}^{\crys,\lambdau,\tau}$
and~$\cX_d^{\semis,\lambdau,\tau}$, the formal algebraic stack~$\cX_d$ is
not a $p$-adic formal algebraic stack.

\begin{prop}\label{prop: X is not padic formal algebraic}
   $\cX_d$ is not a $p$-adic formal algebraic stack.
\end{prop}
\begin{proof}
  Assume that~$\cX_d$ is a $p$-adic formal algebraic stack, so that
  its special fibre~$\overline{\cX}_d:=\cX_d\times_\cO\F$ is an
  algebraic stack, which is furthermore of finite type over~$\F$ (since~$\cX_d$ is a
  Noetherian formal algebraic stack, by Corollary~\ref{cor:Xd is
    formal algebraic}, and~$\cX_{d,\red}$ is of finite type
  over~$\F$, by Theorem~\ref{thm:Xdred is algebraic}). Since the underlying reduced substack
  of~$\overline{\cX}_d$ is~$\cX_{d,\red}$, which is equidimensional of
  dimension~$[K:\Qp]d(d-1)/2$, we see that~$\overline{\cX}_d$ also has
  dimension~$[K:\Qp]d(d-1)/2$.

  Consider a finite type point~$x:\Spec\F'\to\cX_{d,\red}$,
  corresponding to a representation $\rhobar:G_K\to\GL_d(\F')$.
  By Proposition~\ref{prop:versal rings} there
  is a corresponding versal morphism
  $\Spf R_{\rhobar}^\square/\varpi\to\overline{\cX}_d$, and
  applying~\cite[Lem.\ 2.40]{EGcomponents} as in the proof of Theorem~\ref{thm:
    dimension of ss stack},
we conclude
  that~$R_{\rhobar}^\square/\varpi$ must have
  dimension~$d^2+[K:\Qp]d(d-1)/2$. However, it is known that there are
  representations~$\rhobar$ for which~$R_{\rhobar}^\square/\varpi$ is
  formally smooth of dimension $d^2+[K:\Qp]d^2$ (see for
  example~\cite[Lem.\ 3.3.1]{Allen2019}; indeed, a Galois cohomology
  calculation shows that the dimension is always at
  least~$d^2+[K:\Qp]d^2$, which is all that we need). Thus we must
  have~$d^2=d(d-1)/2$, a contradiction. 
\end{proof}

\begin{remark}
\label{rem:non-effective}
A similar argument shows that the versal morphisms $\Spf R_{\rhobar}^{\square,\cO'}
\to \cX_d$ and $\Spf R_{\rhobar}^{\square,\cO'}/\varpi \to \cX_d \otimes_{\cO} \F$
of Proposition~\ref{prop:versal rings} are {\em not} effective.
\end{remark}




\section{Closed points}
\label{subsec:closed points}
Recall that if $\cY$ is  an algebraic stack, then $|\cY|$ denotes its underlying topological
space.  The points of $|\cY|$, which we also refer to as the points of~$\cY$,
 are the equivalence classes of morphisms $\Spec k \to \cY$,
with $k$ a field; two morphisms $\Spec k \to \cY$ and $\Spec l \to \cY$ are deemed to be
equivalent if we may find a field $\Omega$ and
embeddings $k, l \hookrightarrow \Omega$ such that
the pull-backs to $\Spec \Omega$ of the two given morphisms coincide.
A point of $\cY$ is called {\em finite type} if it is representable by
\index{finite type point}
a morphism $\Spec k \to \cY$ which is locally of finite type.   If $\cY$ is locally of
finite type over a base-scheme~$S$, then a point of $\cY$ is of finite type
if and only if its image in $S$ is a finite type point of $S$ in the usual
sense~\cite[\href{https://stacks.math.columbia.edu/tag/01T9}{Tag 01T9}]{stacks-project}.

We say that a point of $\cY$ is {\em closed} if the corresponding point of $|\cY|$ is closed
in the topology on~$|\cY|$.  Closed points are necessarily of finite type.
However, if $\cX$ is not quasi-DM, then finite type 
points of $\cX$ need not be closed, even if $\cX$ is Jacobson (in
contrast to the situation for
schemes~\cite[\href{https://stacks.math.columbia.edu/tag/01TB}{Tag 01TB}]{stacks-project}).

If $\cY$ is locally of finite type over an algebraically closed field~$k$,
then the map $\cY(k) \to |\cY|$ (taking a $k$-valued point of $\cY$ to the
point it represents) identifies $\cY(k)$ with the set of finite type points 
of~$\cY$.  (This is a standard consequence of the Nullstellensatz, applied
to a scheme chart of~$\cY$.)

Throughout this section,
we fix a value of $d \geq 1$, and 
write $\cX$ rather than
$\cX_d$. In order to describe the finite type and closed points
of~$(\cX_{\red})_{\Fbar_p}$, we introduce the following terminology.

\begin{defn}\label{defn: partial ss}
  Let $\rhobar, \thetabar:G_K\to\GL_d(\Fpbar)$ be continuous
  \index{partial semi-simplification}
  representations. Then we say that~$\thetabar$ is a \emph{partial
    semi-simplification} of~$\rhobar$ if  there are short exact
  sequences \[0\to \rbar_i\to\rhobar_i\to\rbar'_i\to 0\]
  for~$i=1,\dots$ such that $\rhobar_1 = \rhobar,$
  $\rhobar_{i+1}=\rbar_i\oplus\rbar_i'$,
  and $\rhobar_n=\thetabar$ for~$n$ sufficiently large.
\end{defn}
Note in particular that the semi-simplification of~$\rhobar$ is a
partial semi-simplification of ~$\rhobar$. We also consider the
following related notion.

\begin{defn}\label{defn: virtual partial ss}Let $\rhobar, \thetabar:G_K\to\GL_d(\Fpbar)$ be continuous
  representations. Then we say that~$\thetabar$ is a \emph{virtual partial
    semi-simplification} \index{virtual partial
    semi-simplification} of~$\rhobar$ if there is a short exact
  sequence of the form \[0\to\rbar\to\rbar\oplus
    \rhobar\to\thetabar\to 0.\]  
\end{defn}
It follows from an evident induction on the integer~$n$ in Definition~\ref{defn: partial ss} that if ~$\thetabar$ is a partial
semi-simplification of~$\rhobar$, then it is in particular a virtual
partial semisimplification; however, the converse does not hold in
general. (See the introduction to~\cite{MR1757882}.)

Our promised description of the finite type and closed points
of~$(\cX_{\red})_{\Fbar_p}$ is given by the following theorem.
\begin{theorem}
\label{thm:closed points}
\leavevmode
\begin{enumerate}
\item
The morphism $\cX(\Fbar_p) = \cX_{\red}(\Fbar_p) \to |(\cX_{\red})_{\Fbar_p}|$ that sends each
morphism $\Spec \Fbar_p \to \cX$ to the point
of $(\cX_{\red})_{\Fbar_p}$  that it represents is injective, and its image
consists precisely of the finite type points.  Thus the finite type points
of $(\cX_{\red})_{\Fbar_p}$ are in natural bijection with the isomorphism
classes of continuous representations $\rhobar: G_K \to \GL_d(\Fbar_p)$.
\item A finite type point
of $|(\cX_{\red})_{\Fbar_p}|$
is closed if and only if the associated Galois representation $\rhobar$ is semi-simple.
\item If $x \in |(\cX_{\red})_{\Fbar_p}|$ is a finite type point, corresponding
to the Galois representation~$\rhobar$, then the closure $\overline{\{x\}}$
contains a unique closed point, whose corresponding Galois
representation is the semi-simplification $\rhobar^{\ss}$ of~$\rhobar$.
More generally, if $x$ and $y$ are two finite type points 
of~$|(\cX_{\red})_{\Fbar_p}|$, corresponding to Galois representations
$\rhobar$ and $\thetabar$ respectively,
then $y$ lies in $\overline{\{x\}}$ if and only if $\thetabar$ is a
virtual partial
semi-simplification of~$\rhobar$. 
\end{enumerate}
\end{theorem}
\begin{proof}
Since $\Fbar_p$-valued points of $\cX$ coincide with $\Fbar_p$-valued points of $\cX_{\red}$,
and so also with $\Fbar_p$-valued points of $(\cX_{\red})_{\Fbar_p}$, the claim of~(1)
is a particular case of the more general consequence of the Nullstellensatz that we
recalled above.

Suppose now that we have a short exact sequence  \numequation\label{eqn:
  virtual ses}0\to\rbar\to\rbar\oplus
    \rhobar\to\thetabar\to 0.\end{equation} We now follow the proof
  of~\cite[Prop.\ 3.4]{MR868301}
 to show that the point~$y$ corresponding
  to~$\thetabar$ lies in the closure of the point~$x$ corresponding
  to~$\rhobar$. Denote the morphism $\rbar\to\rbar\oplus\rhobar$ in~\eqref{eqn:
  virtual ses} by~$(f,g)$. Then for~$t$ in a sufficiently small open
neighbourhood~$U$ of $0\in\A^1$, the morphism
$(f+t\mathbf{1}_{\rbar},g):\rbar\to \rbar\oplus\rhobar$ is an
injection with projective cokernel (to see this one can for example
consider a splitting of\eqref{eqn:
  virtual ses} on the level of vector spaces); so we have a  morphism
(where we as usual abuse notation by writing families of
$(\varphi,\Gamma)$-modules as families of Galois representations)
$U\to\cX$ given by $t\mapsto \rhobar_t:=
(\rbar\oplus\rhobar)/\im(f+t\mathbf{1}_{\rbar},g)$. After possibly
shrinking~$U$ further we may assume that~$f+t\mathbf{1}_{\rbar}$ is an
automorphism of~$\rbar$ for all $\Fpbar$-points $t\ne 0$ of~$U$, so
that $\rhobar_t\cong\rhobar$; while for~$t=0$ we
have~$\rhobar_t\cong\thetabar$. Thus~$\thetabar$ is indeed
in the closure of~$\rhobar$, as claimed. This proves the ``if'' direction of~(3), 
and the ``only if'' direction of~(2).

The ``if'' direction of~(2) follows from the ``only if'' direction of~(3),
and so it remains to prove this latter statement.   To this end,
we fix  $\rhobar,\thetabar: G_K \to \GL_d(\Fbar_p)$, and assume that the point $y$
corresponding to $\thetabar$ lies in the closure of the point $x$ corresponding to~$\rhobar$.
In order to avoid discussing deformation theory over $\Fbar_p$, we extend 
our coefficients $\cO$ if necessary
so that $\rhobar$ and $\thetabar$ are both defined over $\F$.
We may and do also assume that $x$ and $y$ are distinct; i.e.\ that 
$\rhobar$ and $\thetabar$ are not isomorphic (as representations defined over $\Fbar_p$,
or equivalently, as representations defined over~$\F$).
Let $x_0$ and $y_0$ denote the images of $x$ and $y$ respectively in~$\cX_{\red}$;
then $x_0$ and $y_0$ are again distinct.
Let $\cZ$ denote the closure of $\{x_0\}$, thought of as a reduced closed
substack of $\cX_{\red}$; then $y_0$ is a point of~$\cZ$. 

Let $R$ be the framed deformation ring of $\thetabar$ over~$\cO$; then $R$ is a versal
ring to $\cX$ at~$y$. 
Let $\Spf S = \Spf R \times_{\cX}\cZ$; then $S$
is a versal ring to $\cZ$ at $y_0$.   Since $\cZ$ is algebraic,
the versal ring $S$ is effective~\cite[\href{https://stacks.math.columbia.edu/tag/07X8}{Tag 07X8}]{stacks-project},  
and the corresponding morphism $\Spec S \to \cZ$
is furthermore flat~\cite[\href{https://stacks.math.columbia.edu/tag/0DR2}{Tag
  0DR2}]{stacks-project}. 
The morphism $\Spec \F \to \cZ$ corresponding to $x_0$
is quasi-compact and scheme-theoretically dominant;
thus the base-changed morphism
\numequation
\label{eqn:dominant}
\Spec S \times_{\cZ, x_0} \Spec \F \to \Spec S
\end{equation}
 is again scheme-theoretically dominant.

Let $\eta$ be a point in the image of~\eqref{eqn:dominant}, and let
$\Spec T$ be the Zariski closure of $\eta$ in $\Spec S$.  Then $T$ is a quotient
of $R$ which is an integral domain.  If $\cK$ denotes the fraction field of~$T$,
then the morphism $\Spec \cK \to \cX_{\red}$ factors through the morphism
$\Spec \F \to \cX_{\red}$ corresponding to $x_0$, and so the Galois representation
$G_K \to \GL_d(\cK)$ classified by the point $\eta$ of $\Spec T$ is isomorphic
to $\rhobar$ (see Section~\ref{subsubsec: Galois repns into certain
  fields} for the definition of this Galois representation).  
Thus, if $L$ denotes the splitting field of~$\rhobar$,
then we see that the deformations of $\thetabar$ classified by $T$ all factor
through the finite group $\Gal(L/K)$. It follows immediately from~\cite[Thm.\ 1]{MR1757882}  
that $\thetabar$ is a virtual partial semi-simplification
of~$\rhobar$, as required.
\end{proof}

\section{The substack of $G_K$-representations}\label{sec:
  comparison to CWE}
We now briefly explain the relationship between our
stacks and the stacks of $G_K$-representations constructed by Wang-Erickson
in~\cite{MR3831282}. 
Write~$\cX_d^{\Gal}$ for the formal algebraic
stack characterised by the following property: if~$A$ is an
$\cO$-algebra in which $p$ is nilpotent, then $\cX_{d}^{\Gal}(A)$ is the
groupoid of continuous morphisms \[\rho:G_{K}\to\GL_d(A)\](where~$A$
has the discrete topology, and~$G_{K}$ its natural profinite
topology). That this \emph{is} a formal algebraic stack follows
from~\cite[Thm.\ 3.8, Rem.\ 3.9]{MR3831282}. Equivalently, we can
think of $\cX_{d}^{\Gal}(A)$ as the groupoid
of rank~$d$ projective $A$-modules~$T_A$ with a continuous action of~$G_K$
(where each~$T_A$ has the discrete topology).

For any finite type $\cO/\varpi^a$-algebra~$A$, and
any~$T_A\in\cX_d^{\Gal}(A)$, we
set \[\mathbb{D}_A(T_A):=W(\C^\flat)_A\otimes_AT_A,\] which naturally has
the structure of a rank~ $d$ projective $(\varphi,G_K)$-module with
$A$-coefficients (with~$\varphi$ acting on the first factor in the
tensor product, and ~$G_K$ acting diagonally). (In the case that~$A$
is actually a finite $\cO/\varpi^a$-algebra, this agrees with the
construction given in Section~\ref{subsubsec:
Galois reps for GK phi}.)

By Proposition~\ref{prop: equivalences of categories to Ainf}, for
each finite type $\cO/\varpi^a$-algebra~$A$ the assignment
$T_A\mapsto \mathbb{D}_A(T_A)$ gives a
functor~$\cX_d^{\Gal}(A)\to\cX_d(A)$.  Since both~$\cX_d^{\Gal}$
and~$\cX_d$ are limit preserving (for $\cX_d^{\Gal}$ this follows
easily from the definition, since any~$\rho$ as above factors through
a finite quotient of~$G_K$), this defines a
morphism of stacks
\numequation
\label{eqn:Galois stacks to ours}
\cX_d^{\Gal}\to\cX_d.
\end{equation}
We now show that this
morphism is in fact a monomorphism.

\begin{thm}
  \label{thm: map from Galois stack to ours}
  The morphism~{\em \eqref{eqn:Galois stacks to ours}} is a monomorphism,
  and is furthermore versal at finite type points in the sense of Definition~{\em \ref{defn:formally smooth points}} 
below.
\end{thm}
\begin{proof}
  The versality statement follows from the fact that for any
  finite Artinian $\Z_p$-algebra~$A$, the $A$-valued points
  of $\cX_d^{\Gal}$ and $\cX_d$ coincide.  (Over such a ring~$A$,
  \'etale $(\varphi,\Gamma)$-modules do arise from Galois representations.)

  To prove the monomorphism claim, we need to show that for  any finite
  type~$\cO/\varpi^a$-algebra $A$, the functor
  $\cX_d^{\Gal}(A)\to\cX_d(A)$ is fully faithful; that is, we need to
  show that given $T_1,T_2\in \cX_d^{\Gal}(A)$, we
  have \[\Hom_{G_K,A}(T_1,T_2)\stackrel{?}{=}\Hom_{\varphi,G_K}(\mathbb{D}_A(T_1),\mathbb{D}_A(T_2)).\]
  Writing $T:=T_1^\vee\otimes_AT_2$, it suffices to show
  that \[T^{G_K}\stackrel{?}{=}\mathbb{D}_A(T)^{\varphi=1,G_K},\]or even
  that  \[T\stackrel{?}{=}\mathbb{D}_A(T)^{\varphi=1}.\] Since~$T$ is
  a finite projective $A$-module, we
  have  \[\mathbb{D}_A(T)^{\varphi=1}:=(W(\C^\flat)_A\otimes_AT)^{\varphi=1}=(W(\C^\flat)_A)^{\varphi=1}\otimes_AT,\]and
  the result follows from Lemma~\ref{lem:phi invariants}. 
\end{proof}

\begin{rem}
  \label{rem: comparison to what CWE did}In~\cite[\S4]{MR3831282},
  Wang-Erickson explains how to associate \'etale $\varphi$-modules to
  $G_{K_\infty}$-representations with open kernel. A similar argument
  would allow us to associate projective \'etale
  $(\varphi,\Gamma)$-modules to objects of~$\cX_d^{\Gal}(A)$. However,
  if we directly followed the strategy of~\cite{MR3831282},
  the resulting \'etale $(\varphi,\Gamma)$-modules
  would have as coefficients the rings $\A_{K}\otimes_{\Zp}A$,
  rather than the rings~$\A_{K,A}$ that we use throughout this book.
  Thus this approach would not obviously yield Theorem~\ref{thm: map from
    Galois stack to ours}, which is why we've preferred to make a more direct
  argument in terms of \'etale~$(\varphi,G_K)$-modules.
\end{rem}

\begin{rem}
  \label{rem: CWE stack is exactly the substack where there is a
    Galois representation}It follows immediately from
  Theorem~\ref{thm: map from Galois stack to ours} that~$\cX_d^{\Gal}$
  is largest substack of~$\cX_d$ over which the
  universal $(\varphi,\Gamma)$-module can be realised as a
  $G_K$-representation.
  
In fact, it should be possible
to strengthen Theorem~\ref{thm: map from Galois stack to ours}.
For example, \cite[Thm.\ 3.8]{MR3831282}
  shows that we can write \numequation\label{eqn: Galois breaks up
    over
    pseudorepns}\cX_d^{\Gal}=\coprod_{D}\cX_{d,D}^{\Gal},\end{equation}where~$D$
runs over the isomorphism classes of $d$-dimensional semisimple
$\Fpbar$-representations of~$G_K$, and~$\cX_{d,D}^{\Gal}(A)$ is the
groupoid of those~$\rho$ the semisimplification of whose reductions
modulo~$\varpi$ is~$D$. In view of Theorem~\ref{thm:closed points},
we expect that the monomorphism 
$(\cX_{d,D}^{\Gal})_{\red} \to \cX_{d,\red}$ (induced by the monomorphism of
Theorem~\ref{thm: map from Galois stack to ours})
will actually be a closed immersion, and hence (given the versality statement of
Theorem~\ref{thm: map from Galois stack to ours}) that
$\cX_{d,D}^{\Gal}$ will be identified with the formal completion
of $\cX_d$ along the image of this closed immersion.
In fact, this image will contain a unique closed point (corresponding
to the semisimple Galois representation~$D$),
and we imagine that $\cX_{d,D}^{\Gal}$ could also be identified
with the {\em coherent completion}
(in the sense ~\cite[Defn.\ 2.1]{MR4088350})
of $\cX_d$ a this closed point.
However, we don't pursue these ideas further here.
\end{rem}
\begin{rem}
  \label{rem: could also consider WD}As well as considering
  $G_K$-representations, it is also possible to consider a larger
  substack of~$\cX_d$ over which our $(\varphi,\Gamma)$-modules can be
  realised as ``Weil--Deligne''-representations. (Here the notion of
  Weil--Deligne representations is not quite the usual one, but rather
  is given by representations of certain discretizations of~$G_K$, as described
  in~\cite[\S 1.2]{emerton2020moduli} and the references therein.)

  In the case~$K=\Qp$ and~$d=2$, this locus is discussed in~\cite[\S
  1]{emerton2020moduli}; roughly speaking, the difference between it
  and $\cX_d^{\Gal}$ is that it contains families of semisimple
  representations given by direct sums of unramified twists of fixed
  irreducible representations (the key point being that the universal
  unramified character does not correspond to a $G_K$-representation,
  but does correspond to a representation of the Weil group~$W_K$,
  see Remark~\ref{rem: rank 1 versus rank d} below).
\end{rem}

\chapter{The rank one case}\label{sec: the rank one case}
In this chapter we describe some of our key constructions explicitly
in the case where $d = 1$.
More precisely, following the notation of
Chapter~\ref{sec: phi modules and phi gamma modules},
we give explicit descriptions of the stacks $\cR_d$, $\cR_d^{\Gamma_{\disc}}$,
and $\cX_d$, all in the case when $d = 1$ (and $K$ is arbitrary).

\section{Preliminaries}
While it may be possible to find an explicit description 
of the stacks we are interested in
directly from their definitions, this is not how 
we proceed.  Rather, we use the relationship between \'etale $\varphi$-modules
(resp.\ \'etale $(\varphi,\Gamma)$-modules) and Galois representations
to construct certain affine formal algebraic spaces $U$,
$V$, and $W$,
along with morphisms $U \to \cR_1$, $V \to \cR^{\Gamma_{\disc}}_1$,
and $W \to \cX_1$,
each satisfying the conditions of Lemma~\ref{lem:presentations} below.
An application of this lemma then yields a description
of each of $\cR_1,$ $\cR^{\Gamma_{\disc}}_1,$ and $\cX_1$.

%

\begin{defn}\index{versal at finite type points}
  \label{defn:formally smooth points}Let~$S$ be a locally Noetherian
  scheme, let $\cX$ and $\cY$ be stacks over~$S$ and let $f:\cX\to\cY$ be a
  morphism.   We say that $f$ is {\em versal
at finite type points} if for each morphism $x:\Spec k \to \cX$,
with $k$ a finite type $\cO_S$-field,
the morphism $f$ is versal at~$x$, in the sense of Definition~\ref{def:versal}.
\end{defn}

\begin{remark}
As the discussion of versality in Appendix~\ref{app: formal algebraic stacks}
should make clear,
this definition is related to the various properties considered in
\cite[Def.~2.4.4]{EGstacktheoreticimages}
and the surrounding discussion.  As noted in 
\cite[Rem.~2.4.6]{EGstacktheoreticimages}, it is somewhat complicated to define
the notion of versality at literal points of a stack, since the points of a stack
are (by definition) equivalence classes of morphisms from the spectrum of a field
to the stack, and it is not immediately clear in general whether
or not  the versality condition would be
independent of the choice of equivalence class representative.  
In the preceding definition we obviate this point, by directly requiring 
versality at all (finite type) representatives of all (finite type) points.
\end{remark}

\begin{lemma}
\label{lem:base-changing versality}
If $\cX \to \cY$ is a morphism of stacks over a locally Noetherian base $S$ which 
is versal at finite type points, 
and if $\cZ \to \cY$ is a morphism of stacks, then the induced morphism
$\cX \times_{\cY} \cZ \to \cZ$ is again versal at finite type points.
\end{lemma}
\begin{proof}
Let $\tx:\Spec k \to \widetilde{\cX} := \cX\times_{\cY} \cZ$ with $\Spec k$ of finite type over $S$;  from the definition of the $2$-fibre product $\widetilde{\cX}$,
we see that giving the morphism $\tx$ amounts to giving the 
pair of morphisms $x: \Spec k \to \cX$ and $z: \Spec k \to \cZ$ (the result
of composing $\tx$ with each of the projections),
as well as an isomorphism $y \iso y'$, where $y$ and $y'$ are the $k$-valued
points of $\cX$ obtained respectively by composing~ $x$ with the morphism $\cX \to \cY$
and by composing $z$ with the morphism $\cZ \to \cY$.
It is then a straightforward diagram chase,
working from the various definitions, 
to deduce the versal property of the morphism $\widehat{(\widetilde{\cX})}_{\widetilde{x}}
\to \widehat{\cZ}_z$ (where we use the notation of Definition~\ref{adefn: deformation category}) from the versal property of the morphism 
$\widehat{\cX}_x \to \widehat{\cY}_y.$
\end{proof}

The following result is standard, but we indicate the proof for the sake of
completeness.

\begin{lemma}
\label{lem:versality --- the algebraic case}
If $f: \cX \to \cY$ is a morphism of algebraic stacks, each of finite type
over a locally Noetherian scheme $S$,
then $f$ is versal at finite type points if and only if $f$ is smooth.
\end{lemma}
\begin{proof}
The ``if'' direction follows from an application of the infinitesimal lifting
property of smooth morphisms; see e.g.\ the implication ``$\text{(4)}
\implies  \text{(1)}$''
of
\cite[Lem.~2.4.7~(4)]{EGstacktheoreticimages}.
For the converse, choose a smooth
surjective morphism $U \to \cX$ with $U$ a scheme (necessarily
locally of finite type over~$S$); then by the direction
already proved, this morphism is also versal at finite type points,
and hence so is the composite morphism $U \to \cX \to \cY$.  To show that $f$
is smooth, it suffices to show the same for this composite.  
It follows from the implication ``(1) $\implies$ (4)'' of
\cite[Lem.~2.4.7~(4)]{EGstacktheoreticimages} that $f$ is smooth
in a neighbourhood of each finite type point of $U$; but since these points
are dense in $U$, the desired result follows.
(We have cited \cite{EGstacktheoreticimages} in this argument
purely for our own convenience; of course
the present lemma is just a version of Grothendieck's result relating smoothness
and formal smoothness.)
\end{proof}

We now strengthen the previous result to cover the case where $\cX$ and $\cY$ are
not necessarily assumed to be algebraic, but merely the map between them is
assumed to be representable by algebraic stacks.

\begin{lemma}
\label{lem:versality --- the representable case}
If $f: \cX \to \cY$ is a morphism of limit preserving
stacks over a locally Noetherian base $S$,
and~$f$ is representable by algebraic stacks, 
then $f$ is versal at finite type points if and only if $f$ is smooth.
\end{lemma}
\begin{proof}
If $f$ is smooth, then by Lemma~\ref{lem:characterizing properties}~(2)
it satisfies the infinitesimal lifting property,
which immediately implies that it is versal at finite type points.
Conversely, suppose now that $f$ is versal at finite type points.

It follows from~\cite[Cor.~2.1.8]{EGstacktheoreticimages} 
that  $f$
is limit preserving on objects,
and since it is also representable by algebraic stacks,
it is locally of finite presentation,
by Lemma~\ref{lem:characterizing properties}~(1).
We now have to verify that for any test morphism $Z \to \cY$ whose source is a scheme, the induced morphism 
\numequation
\label{eqn:morphism induced by f}
\cX\times_{\cY} Z \to Z ,
\end{equation}
which is again locally of finite presentation,
is in fact smooth. Of course it suffices to do this for affine schemes, and then,
since $\cY$ is limit preserving, for affine schemes that are of finite type
over $S$.  Since $f$ is representable by algebraic stacks,
we find that the source of~\eqref{eqn:morphism induced by f} 
is an algebraic stack.  Furthermore, since this morphism is locally of finite
presentation, its source is again locally of finite type over~$S$. 
Finally, Lemma~\ref{lem:base-changing 
versality} shows that~\eqref{eqn:morphism induced by f} 
is versal at finite type points.
The claimed result then follows from Lemma~\ref{lem:versality --- the algebraic case}.
\end{proof}

We may apply the previous lemma to obtain a criterion 
for representing a stack as the quotient of a sheaf by a smooth groupoid.
We first make another definition.

\begin{defn}
  \label{defn:surjective on f.t. points}
We say that a morphism $\cX\to \cY$ 
of stacks over a locally Noetherian scheme~$S$
is {\em surjective
on finite type points} \index{surjective
on finite type points} if for each morphism $x:\Spec k \to \cY$,
with $k$ a finite type $\cO_S$-field, we may find a field extension
$l$ of $k$ and a morphism $\Spec l \to \cX$ which makes the diagram
$$\xymatrix{ \Spec l \ar[r]\ar[d] & \cX \ar[d] \\
\Spec k \ar[r] & \cY}
$$
commutative.
\end{defn}

\begin{lemma}
\label{lem:presentations}
Suppose that $\cZ$ is a stack over a locally Noetherian scheme~$S$,
and that $f:U \to \cZ$ is a morphism from a sheaf to $\cZ$ which
is representable by algebraic spaces\footnote{Since the source
of the morphism is a sheaf,
this is equivalent to being representable by algebraic stacks.}
and surjective on finite type points {\em(}in the sense of Definition~{\em \ref{defn:surjective
on f.t. points})}.
Suppose furthermore that both $U$ and $\cZ$ are limit preserving,
and that~$f$ is 
versal at  finite type points.
Then the morphism $f$ is in fact smooth and surjective,
and, if we let $R := U\times_{\cZ} U$,
then the induced morphism $[U/R] \to \cZ$ is an isomorphism.
\end{lemma}
\begin{proof}
It follows from Lemma~\ref{lem:versality --- the representable case}
that the morphism $f$ is smooth.
We claim that it is furthermore surjective; by definition,
this means that if $T \to \cZ$ is any test morphism from an affine scheme,
we must show that the morphism of algebraic spaces
$g:T\times_{\cZ} U \to T$ induced by $f$ is surjective.
Since $\cZ$ is limit preserving, we may assume that $T$ is of finite type
over~$S$; also, since $g$ is a base-change of the smooth morphism $f$,
it is itself smooth, and so in particular has open image.  Our assumption that
$f$ is surjective on finite type points implies furthermore that the image
of $g$ contains all the finite type points of $T$; thus $g$ is indeed surjective.

Since smooth morphisms are in particular flat, the second assertion
of the lemma follows from Lemma~\ref{lem:stacks as quotients}.
\end{proof}

\begin{remark}
If we omit
the assumption in Lemma~\ref{lem:presentations} that the morphism $U \to \cZ$ be
representable by algebraic spaces, then the morphism $[U/R] \to \cZ$ need 
not be an isomorphism.

An example illustrating this is given by taking $\cZ$ to be a scheme~$Z$,
locally of finite type 
over a locally Noetherian base scheme~$S$, choosing a closed subscheme $Y$
of $Z$ which is {\em not} also open, and defining $U$ to be the formal
scheme obtained as the disjoint union $U := (Z\setminus Y) \coprod \widehat{Z}$,
where $\widehat{Z}$ denotes the completion of $Z$ along~$Y$.
We let $U \to Z$ be the obvious morphism; this is then a surjective monomorphism,
which is not an isomorphism, although it is versal at
every finite type point, and (since it is a monomorphism)
we have that $R := U\times_Z U$ coincides with the diagonal
copy of $U$ inside $U \times_S U$ (so that $[U/R] = U$).
\end{remark}

We will find the following lemmas useful in verifying the hypotheses of
Lemma~\ref{lem:presentations} in our application.

\begin{lem}
\label{lem:Artinian test}
Let $f:\cX  \to \cY$ and $g:\cY \to \cZ$ be morphisms of stacks over a locally
Noetherian scheme~$S$. Assume that~$f$ is representable by algebraic stacks
and locally of finite type, and that both $g$ and the morphism $\cZ \to S$
are limit preserving on objects.
Then $f$ is an isomorphism if and only if
the induced morphism
\numequation
\label{eqn:A pullback}
\Spec A \times_{\cZ} \cX \to \Spec A \times_{\cZ} \cY
\end{equation}
is an isomorphism
for every morphism $\Spec A \to \cZ$ with 
$A$ a locally of finite type Artinian local $\cO_S$-algebra. 
\end{lem}
\begin{proof}
The ``only if'' direction is evident,
and so we focus on the ``if'' direction. We begin by reducing to the
case $\cZ=S$.
We first recall a basic fact about fibre products:
if $T \to \cY$ is any morphism from a scheme to the stack~$\cY$,
then the base-change $T\times_{\cY} \cX$
may be described as an iterated fibre product
$T \times_{(T\times_{\cZ} \cY)} (T\times_{\cZ} \cX)$
(where we regard $T$ as lying over $\cZ$ via the composite
$T \to \cY \to \cZ$, and as lying over $T \times_{\cZ} \cY$
via the graph of the given morphism $T \to \cY$).

In order to show that $f$ is an isomorphism,
we have to show that for any morphism $T \to \cY$,
the base-changed morphism $T\times_{\cY} \cX \to T$
is an isomorphism.  By the fact about fibre products just recalled
(and remembering that the base-change of an isomorphism 
is an isomorphism),
it suffices to show that for any morphism $T \to \cZ$ whose source is
a scheme, the base-changed morphism
\numequation
\label{eqn:base-changed morphism}
T\times_{\cZ} \cX \to T\times_{\cZ} \cY
\end{equation}
is an isomorphism.  Since $\cZ \to S$ is limit preserving on objects, we may and do assume
that $T$ is locally of finite type over~$S$.
We will then replace both $\cZ$ and $S$ by~$T$ (mapping identically to itself),
and the morphism $f$ by~\eqref{eqn:base-changed
morphism}.
Note that the projection
$T\times_{\cZ} \cY \to T$ is again limit preserving on objects (being
a base-change of~$g$, which has this property). 
We now have to reinterpret our original hypothesis on $A$-valued
points of $\cZ$ (for locally of finite type Artinian local $\cO_S$-algebras~$A$)
in this new context.
To this end, note that if $\Spec A\to T$
is locally of finite type with $A$ being Artinian local,
then $\Spec A$ is also 
locally of finite type over~$S$. 
The fibre product of~\eqref{eqn:base-changed morphism} 
with $\Spec A$ over $T$ may then be identified with
the morphism~\eqref{eqn:A pullback},
and is thus an isomorphism (by assumption).
Putting all this together, we see that we have
reduced to the case when $\cZ~=~S$;
so we return to our original notation, but assume in addition 
that $\cZ = S$ (mapping to itself via the identity).

%
What we have to show is that if $T \to \cY$ is any morphism,
then
\numequation
\label{eqn:T pullback}
T\times_{\cY} \cX \to T
\end{equation}
is an isomorphism.  Since $\cY \to S$
is limit preserving, we may assume that $T$ is locally of finite type over~$S$.
Any finite type Artinian local $T$-algebra is then also finite type over $S$,
and so, replacing both $\cY$ and $S$ by $T$
(and again applying the fibre product fact recalled above, with $T$ replaced by $\Spec A$) 
we find that we may make another reduction, to the case when $\cY = S$.
The representability assumption on $f$ then implies
that $\cX$ is an algebraic stack.
In this case the lemma is standard, but we recall a proof for completeness.



Our hypothesis that~$f$ induces an isomorphism on $A$-valued points implies that the morphism $f$
is versal at finite type points, and so by Lemma~\ref{lem:versality --- the algebraic case},
it is smooth.  It also implies that $f$ contains all finite
type points in its image; since the image of $f$ is constructible, we see that $f$
is in fact surjective.  This same argument implies that the diagonal
morphism $\cX \to \cX\times_S \cX$ 
is surjective, and thus that $f$ is universally injective.  
Since $f$ is smooth, it is flat and locally of finite type.  Since it 
is also universally injective, it is an open immersion.  Since it is surjective,
it is an isomorphism, as claimed.
\end{proof}

\begin{lem}
  \label{lem: checking closed subsheaf inclusion on versal rings}
Let~$\cY$ be a stack over a locally Noetherian scheme~$S$, 
  let $\cZ$ be a closed substack of~$\cY$,
and let $f:\cX \to \cY$ be a morphism of stacks. 
Suppose that both~ $f$ and the morphism $\cY \to S$ are limit preserving on objects;
then~$f$ factors through~ $\cZ$ if and only 
  if and only if for any finite local Artinian $\cO_S$-algebra~$A$, and any
  morphism $\Spec A\to \cY$ over~$S$, the morphism
$\Spec A\times_{\cY} \cX \to \cY$ factors through $\cZ$. 
\end{lem}
\begin{proof}
The ``only if'' direction is trivial. For the converse,
consider the base-change via $f$ 
of the closed immersion $\cZ \hookrightarrow \cY$;
this is a closed immersion
$g: \cX\times_{\cY} \cZ \hookrightarrow \cX$ of stacks over~$\cY$,
which we must show is an isomorphism.
The hypothesis implies that for
any morphism $\Spec A \to \cY$ whose source is a finite type local Artinian 
$\cO_S$-algebra,  the pull-back of~ $g$ \[\Spec
  A\times_\cY(\cX\times_{\cY}\cZ)\to\Spec A\times_{\cY}\cX \] is an isomorphism; 
it then follows from Lemma~\ref{lem:Artinian test}
that $g$ is an isomorphism. 
\end{proof}

The next lemma lets us compute fibre products over stacks in certain situations.
We recall that if $\cY$ is a stack (over some base scheme~$S$),
then the inertia stack $\cI_{\cY}$ is the stack over $S$ whose $T$-valued points
(for any $S$-scheme $T$) consist of pairs $(y,\alpha)$, where $y$ is a $T$-valued
point of $\cY$, and $\alpha$ is an automorphism of~$y$. 
It is a group object in the category of stacks lying over $\cY$ 
via morphisms that are representable by sheaves.
We also remind the
reader that there is a canonical isomorphism
$$\cI_{\cY} \iso \cY\times_{\Delta_{\cY},\cY\times_S \cY, \Delta_{\cY}} \cY,$$
where the subscripts $\Delta_{\cY}$ indicate that both copies of $\cY$ are regarded 
as lying over the product $\cY\times_S \cY$ via the diagonal morphism $\Delta_{\cY}:
\cY \to\cY\times_S \cY$, and where, under this identification,
 the forgetful morphism $\cI_{\cY} \to \cY$ may be identified with projection onto
the first copy of~$\cY$.

\begin{lemma}
\label{lem:fibre product and inertia}
Let $X \to \cY$ be a morphism from a sheaf to a stack {\em(}both over some base-scheme~$S${\em )}, and suppose 
that the natural morphism $X\times_{\cY} X \to X\times_S X$ factors through the diagonal
copy of~$X$ lying in the target {\em (}so that the groupoid $X\times_{\cY} X$ is in fact
a group object in the category of sheaves over~$X${\em )}.  
Then there is a canonical isomorphism
$X\times_{\cY} X \iso X\times_{\cY} \cI_{\cY}$
of group objects over~$X$.
\end{lemma}
\begin{proof}
By assumption, we have a commutative diagram
$$\xymatrix{X\times_{\cY} X \ar[d] \ar[r] & X \ar^-{\Delta_X}[r]\ar[dr] & X\times_S X \ar[d] \\
\cY \ar^-{\Delta_{\cY}}[rr] & & \cY\times_S \cY}
$$
in which the outer rectangle is $2$-Cartesian.  
Since $X$ is a sheaf, the diagonal $\Delta_X: X \to X \times_S X$ is a monomorphism,
and so one immediately checks that the left-hand trapezoid is  also 
$2$-Cartesian.  Thus we obtain the required isomorphism
\[X\times_{\cY} X \iso X \times_{\cY\times_S \cY} \cY \iso X \times_{\cY}
(\cY\times_{\cY\times_S\cY} \cY) \iso X\times_{\cY} \cI_{\cY}. \qedhere\]
\end{proof}

\section{Moduli stacks in the rank one case}
We now give concrete descriptions of our various stacks in the rank one case.
Before doing this, we prove a result that describes the stack structure
on families of rank one objects.  (It encodes the fact that the
automorphisms of a rank one object are simply the scalars.)

\begin{lemma}
\label{lem:rank one automorphisms}
If $\cZ$ denotes either of the stacks $\cR_1$ or $\cX_1$,
then there is a canonical isomorphism
\numequation
\label{eqn:inertia stack of R1}
\cZ \times_{\cO} \widehat{\mathbb G}_m \iso \cI_{\cZ}
\end{equation}
of group objects over~$\cZ$  (where $\widehat{\G}_m$ is the $\varpi$-adic
completion of $\Gm$ over~$\cO$).
\end{lemma}
\begin{proof}
For definiteness we give the proof in the case of $\cR_1$; the proof in the $\cX_1$
case is identical.
We begin by defining the morphism~\eqref{eqn:inertia stack of R1}:
if $A$ is any $\varpi$-adically complete $\cO$-algebra,
and $M$ is an \'etale $\varphi$-module over $\A_{K,A},$
then~\eqref{eqn:inertia stack of R1} is defined on $A$-valued points
of the source lying over~$M$ by mapping an element $a \in A^{\times}$
to the automorphism of $M$ given by multiplication by~$a$.

Since $\cR_1$ has affine diagonal which is of finite presentation,
the morphism~\eqref{eqn:inertia stack of R1}
is a morphism between finite type group objects over $\cR_1$,
each of whose structure morphisms is representable by algebraic spaces,
indeed affine,
and of finite presentation.    Furthermore, $\cR_1$ itself is limit
preserving over~$\cO$. Thus 
to show that~\eqref{eqn:inertia stack of R1} is an isomorphism,
it suffices (by Lemma~\ref{lem:Artinian test})
to show that it induces an isomorphism on $A$-valued points
for $A$ an Artinian local $\cO$-algebra of finite type.
Returning to the notation of the preceding paragraph (but now assuming
that $A$ is Artinian local),
we have to show that any automorphism of $M$ is 
given by multiplication by an element of~$A^{\times}$.
Since $M$ is of rank~$1$,
we find that
$$\Hom_{\A_{K,A},\varphi}(M,M) \cong (M^{\vee}\otimes_{\A_{K,A}}M)^{\varphi=1} \cong (\A_{K,A})^{\varphi=1}
\cong A$$
(the final isomorphism following from~Lemma~\ref{lem:phi invariants}, since $\A_{K,A} \subseteq
W(\C^{\flat})_A$ by Proposition~\ref{prop: maps of coefficient rings
  are faithfully flat injections}). 
The lemma follows.
\end{proof}

\subsection{Local Galois theory}
We briefly recall the local Galois theory that is relevant to our
computation. If~$L/K$ is an algebraic extension, we write~$I_L^\ab$
for the image of~$I_L$ in~$G_L^{\ab}$ (note that this is not the
abelianization of~$I_L$).
We have short exact sequences
\numequation
\label{eqn:K ramification s.e.s.}
1 \to I_K \to G_K \to \Frob_K^{\Zhat} \to 1 
\end{equation}
and
\numequation
\label{eqn:K-cyc ramification s.e.s}
1 \to I_{K_{\cyc}} \to G_{K_{\cyc}} \to \Frob_{K}^{f \Zhat} \to 1,
\end{equation}
where $f := [k_{\infty}: k].$

These induce corresponding short exact sequences
\numequation
\label{eqn:abelian K ramification s.e.s.}
1 \to I_K^{\ab} \to G_K^{\ab} \to \Frob_K^{\Zhat} \to 1 
\end{equation}
and
\numequation
\label{eqn:abelian K-cyc ramification s.e.s}
1 \to I_{\Kcyc}^{\ab} \to G_{K_{\cyc}}^{\ab} \to \Frob_K^{f \Zhat} \to 1.
\end{equation}
%

If $G$ is any of the various profinite Galois groups appearing 
in the preceding discussion, we let $\cO[[G]]$ denote the corresponding
completed group ring over~$\cO$.  This is a pro-Artinian ring, which
gives rise to the affine formal algebraic space $\Spf \cO[[G]]$.  We endow
this space with the trivial action of $\widehat{\G}_m$ (the $\varpi$-adic
completion of $\Gm$ over~$\cO$).

\subsection{Galois lifting rings for characters}If $\F'/\F$ is a
finite extension, and $\rhobar:G_{\Kcyc}\to(\F')^\times$ is a
continuous character, then we can write
~$\cO'=\cO\otimes_{W(\F)}W(\F')$, and consider the universal lifting
$\cO$-algebra $R_{\rhobar}^{\square,\cO'}$. If $x:\Spf\F'\to\cR_1$ is
the corresponding finite type point, then we have a versal morphism
$\Spf R_{\rhobar}^{\square,\cO'}\to \cR_1$; indeed it follows exactly
as in the proof of Proposition~\ref{prop:versal rings} that we have an
isomorphism \numequation\label{eqn: deformation ring character
  torsor}\Spf R_{\rhobar}^{\square,\cO'} \times_{\cR_1} \Spf
R_{\rhobar}^{\square,\cO'} \iso
(\Gmhat)_{R_{\rhobar}^{\square,\cO'}},\end{equation} where
$(\Gmhat)_{R_{\rhobar}^{\square,\cO'}}$ denotes the completion of
$(\Gm)_{R_{\rhobar}^{\square,\cO'}}$ along $(\Gm)_{\F'}$.

\subsection{Descriptions of the rank one stacks}
We can now establish our explicit description of~$\cR_1$.

%

\begin{prop}
	\label{prop:rank 1 R}
       	There is an isomorphism
	$$
\Bigl[ \Bigl ( \Spf \cO[[I^{\ab}_{\Kcyc}]] \times \widehat{\mathbb G}_m
	\Bigr) / \widehat{\mathbb G}_m \Bigr]
\iso
	\cR_{1}$$
{\em (}where, in the formation of the quotient stack, the $\widehat{\mathbb
G}_m$-action is taken to be trivial{\em )}.
\end{prop}
\begin{proof}
We begin by constructing a morphism
	\numequation\label{eqn: covering R1}\Spf \cO[[I^{\ab}_{\Kcyc}]] \times \widehat{\mathbb G}_m
	\to \cR_1.\end{equation} For this, we choose a lift of
        $\sigma_{\cyc} \in G_{K_{\cyc}}$ of $\Frob_K^f$, splitting the
        short exact sequences~\eqref{eqn:K-cyc ramification s.e.s}
        and~\eqref{eqn:abelian K-cyc ramification s.e.s}.  If $A$ is
        any discrete Artinian quotient of $\cO[[I^{\ab}_{\Kcyc}]],$ then we
        extend the continuous morphism $\cO[[I^{\ab}_{\Kcyc}]] \to A$ to a
        continuous morphism $\cO[[G_{K_{\cyc}}^{\ab}]]\to A$ by
        mapping $\sigma_{\cyc}$ to $1 \in A$.  We may view this latter
        morphism as a rank~$1$ representation of $G_{K_{\cyc}}$ with
        coefficients in~$A$; it thus gives rise to a rank one
        projective \'etale $\varphi$-module $M_A$ over $\A_{K,A}$, and
        therefore to a morphism $\Spf \cO[[I^{\ab}_{\Kcyc}]] \to \cR_1.$ We
        extend this to the morphism~\eqref{eqn: covering R1} by
        unramified twisting; more precisely, we have an induced
        morphism $\Spf \cO[[I^{\ab}_{\Kcyc}]]\times\Gmhat \to
        \cR_1\times\Gmhat$, and we compose with the morphism
        $\cR_1\times\Gmhat\to\cR_1$ given by taking the tensor product
        with the universal unramified rank one \'etale $\varphi$-module.


By construction (and the usual explicit description of the universal
Galois deformation ring of a character), 
~\eqref{eqn: covering R1} is versal at finite type
points, as well as surjective on finite type points.   
We claim that the canonical morphism
\nummultline
\label{eqn:morphism of diagonals}
\bigl( \Spf \cO[[I^{\ab}_{\Kcyc}]] \times \widehat{\mathbb G}_m\bigr)
\times_{\cR_1} 
\bigl( \Spf \cO[[I^{\ab}_{\Kcyc}]] \times \widehat{\mathbb G}_m\bigr)
\\
\to
\bigl( \Spf \cO[[I^{\ab}_{\Kcyc}]] \times \widehat{\mathbb G}_m\bigr)
\times_{\cO} 
\bigl( \Spf \cO[[I^{\ab}_{\Kcyc}]] \times \widehat{\mathbb G}_m\bigr)
\end{multline}
factors through the diagonal copy of
$ \Spf \cO[[I^{\ab}_{\Kcyc}]] \times \widehat{\mathbb G}_m$
in the target; given this, it follows from Lemmas~\ref{lem:fibre product and inertia}
and~\ref{lem:rank one automorphisms} that there is an isomorphism of groupoids
\[\bigl ( \Spf \cO[[I^{\ab}_{\Kcyc}]] 
	\times \widehat{\mathbb G}_m
	\bigr) \times \Gmhat\isoto\bigl ( \Spf \cO[[I^{\ab}_{\Kcyc}]] 
	\times \widehat{\mathbb G}_m
	\bigr)\times_{\cR_1}\bigl ( \Spf \cO[[I^{\ab}_{\Kcyc}]] 
	\times \widehat{\mathbb G}_m
	\bigr)\]
with $\widehat{\mathbb G}_m$ acting trivially on $\Spf \cO[[I^{\ab}_{\Kcyc}]]\times
\widehat{\mathbb G}_m.$
The proposition will then follow from Lemma~\ref{lem:presentations}
provided we show that ~\eqref{eqn: covering R1} is representable by algebraic
spaces (or, equivalently, by algebraic stacks).


It remains to verify the claimed factorization of~\eqref{eqn:morphism of diagonals},
as well as the claimed representability by algebraic stacks of~\eqref{eqn: covering R1}.
We establish each of these claims in turn.

To show that~\eqref{eqn:morphism of diagonals} factors through
the diagonal diagonal copy of
$ \Spf \cO[[I^{\ab}_{\Kcyc}]] \times \widehat{\mathbb G}_m$
in the target, it suffices,
by Lemma~\ref{lem: checking closed subsheaf inclusion on versal rings},
to show that if
$$
\Spec A \to 
 \bigl(\Spf \cO[[I^{\ab}_{\Kcyc}]] \times \widehat{\mathbb G}_m\bigr)
\times
 \bigl(\Spf \cO[[I^{\ab}_{\Kcyc}]] \times \widehat{\mathbb G}_m\bigr)
$$
is a morphism whose source is a finite type Artinian local $\cO$-algebra,
then the pull back of \eqref{eqn:morphism of diagonals}
to $\Spec A$ factors through the pull-back to $\Spec A$ of the diagonal morphism
$$
 \Spf \cO[[I^{\ab}_{\Kcyc}]] \times \widehat{\mathbb G}_m
\to
\bigl( \Spf \cO[[I^{\ab}_{\Kcyc}]] \times \widehat{\mathbb G}_m\bigr)
\times_{\cO} 
\bigl( \Spf \cO[[I^{\ab}_{\Kcyc}]] \times \widehat{\mathbb G}_m\bigr)
.$$  Since morphisms $\Spec A \to
 \Spf \cO[[I^{\ab}_{\Kcyc}]] \times \widehat{\mathbb G}_m$
correspond to characters $\chi: G_{\Kcyc}^{\ab} \to A^{\times}$,
this amounts to the evident fact that if $\chi,\chi'$ are two such characters,
then the locus in $\Spec A$ over which $\chi$ and $\chi'$ become 
isomorphic coincides with the locus over which they coincide.

We now turn to proving that~\eqref{eqn: covering R1} is representable by algebraic
stacks.
To see this, it suffices to study the corresponding question modulo
some power $\varpi^a$ of the uniformizer in~$\cO$, so we work with the
stack~$\cR_{K,1}^a$ from now on.

We next recall some of the constructions we made in
Section~\ref{subsec: defn of Xd}. We have the subfield
$\Kbasic\subseteq K$ of Definition~\ref{defn: Kbasic}, and the natural
morphism $\cR_{K,1}^a\to\cR_{\Kbasic,[K:\Kbasic]}^a$ given by the
forgetful map which regards a rank one $\A_{K,A}$-module as an
$\A_{\Kbasic,A}$-module of rank~$[K:\Kbasic]$. By Lemma~\ref{lem: explicit Ind
  description for R pulled back from K basic}, we can write
$\cR_{K,1}^a=\varinjlim_h\cR_{K,1,\Kbasic,h}^a$, where the algebraic
stack $\cR_{K,1,\Kbasic,h}^a$ is the scheme-theoretic image of the
base-changed morphism
  \[\cC_{[K:\Kbasic],h}^a\times_{\cR_{\Kbasic,[K:\Kbasic]}}\cR_{K,1}\to\cR_{K,1},\]
  where $\cC_{[K:\Kbasic],h}^a$ is 
  the algebraic stack of
  $\varphi$-modules over~$\A^+_{\Kbasic,A}$ of rank~$[K:\Kbasic]$ and~$T$-height at
  most~$h$. 
  Accordingly, it suffices to show that  
 each base-change 
	$$\cYha:=\Bigl( \Spf  \cO[[I^{\ab}_{\Kcyc}]] \times \widehat{\mathbb G}_m\Bigr )
\times_{\cR_{K,1}} \cR_{K,1,\Kbasic,h}^a$$ is an algebraic stack.

By construction, $\cYha$ is a closed subsheaf of $\Spf  (\cO/\varpi^a)[[I^{\ab}_{\Kcyc}]]
\times (\Gm)_{\cO/\varpi^a}$. We claim that it is in fact a closed
subsheaf of
$\Spec  (\cO/\varpi^a)[I^{\ab}_{\Kcyc}/U] \times
(\Gm)_{\cO/\varpi^a}$,
for some open subgroup $U$ of $I^{\ab}_{\Kcyc}$.


Since~$\cY_h^a$ and~$\Spec(\cO/\varpi^a)[I^{\ab}_{\Kcyc}/U] \times
(\Gm)_{\cO/\varpi^a}$ are both
closed subsheaves of $\Spf (\cO/\varpi^a)[[I^{\ab}_{\Kcyc}]] \times
(\Gm)_{\cO/\varpi^a}$, it follows from Lemma~\ref{lem: checking closed subsheaf inclusion on versal
  rings} that 
it is enough to show that if~$A$ is a finite local Artinian
$\cO/\varpi^a$-algebra, then given any morphism $\Spec A\to \cYha$,
the composite
$\Spec A\to\cYha\to\Spf (\cO/\varpi^a)[[I^{\ab}_{\Kcyc}]] \times
(\Gm)_{\cO/\varpi^a}$ factors through
$\Spec (\cO/\varpi^a)[I^{\ab}_{\Kcyc}/U] \times
(\Gm)_{\cO/\varpi^a}$ 
for the open subgroup $U:=I^{a,h}$ 
of Proposition~\ref{prop: ramification bound in rank one case with
  coefficients and any K} below.
If~$A$
is in fact a field, then by~\cite[Lem.\
3.2.14]{EGstacktheoreticimages}, the \'etale $\varphi$-module~$M_A$
corresponding to the morphism $\Spec A\to \cR_{K,1}$ attains
$(\Kbasic,T)$-height at most~$h$
in the sense of Definition~\ref{defn: Kbasic T height at most h} below
after passing to a finite extension of~$A$,
and the required factorisation is immediate from
Proposition~\ref{prop: ramification bound in rank one case with coefficients
    and any K}. 

Having established the result in the case that $A$ is a finite type
field, in order to prove the general local Artinian case, it suffices to prove a factorisation on the level of versal rings. 
More precisely, let $\F'/\F$ be a finite extension, and let
$x:\Spec\F'\to\cYha$ be a finite type point of~$\cYha$, corresponding
to a representation $\rhobar:G_{\Kcyc}\to(\F')^\times$. As above, we
write~$R_{\rhobar}^{\square,\cO'}$ for the corresponding universal
lifting $\cO'$-algebra. Write $\Spf R^{a}_h$ for the scheme-theoretic
image of the morphism
\[
  \cC_{[K:\Kbasic],h}^a\times_{\cR_{\Kbasic,[K:\Kbasic]}}\cR_{K,1}\times_{\cR_{K,1}}\Spf
  R_{\rhobar}^{\square,\cO'}\to\Spf R_{\rhobar}^{\square,\cO'}.\] By
Lemma~\ref{alem: scheme theoretic image of versal is versal}, the
induced morphism $\Spf R^a_h\to \cR_{K,1,\Kbasic,h}^a$ is a versal
morphism at~$x$. Write~$R_{\rhobar}^{\square,\cO',U}$ for the quotient
of~$R_{\rhobar}^{\square,\cO'}$ corresponding to liftings which are
trivial on~$U$; 
then it suffices to show that $\Spf
R^a_h$ is a closed formal subscheme of ~$\Spf
R_{\rhobar}^{\square,\cO',U}$.

We now employ Lemma~\ref{lem: criterion for Artin to map to scheme
  theoretic image}. Exactly as in the proof of Lemma~\ref{lem:Artinian points of scheme theoretic images}, it
is enough (after possibly increasing~$\F'$) to show that if~$A$ is a
finite type local Artinian $R_{\rhobar}^{\square,\cO'}$-algebra with
residue field~$\F'$ for which the induced morphism
$\cC_{[K:\Kbasic],h}^a\times_{\cR_{\Kbasic,[K:\Kbasic]}}\cR_{K,1}\times_{\cR_{K,1}}\Spec
A\to \Spec A$ admits a section, then the morphism $\Spec A\to \Spf
R_{\rhobar}^{\square,\cO'}$ factors through~$\Spf
R_{\rhobar}^{\square,\cO',U}$. Since the existence of the section
implies (by definition) that the \'etale $\varphi$-module $M_A$ has
$(T,\Kbasic)$-height at most~$h$
(in the sense of Definition~\ref{defn: Kbasic T height at most h}),
we are done by Proposition~\ref{prop: ramification bound in rank one case with coefficients
    and any K}.
\end{proof}

In the previous argument, we used the following definition,
which will play an important role in the technicalities
of Section~\ref{subsec:ram bound} below.
(See in particular the statements of
Propositions~\ref{prop: ramification bound in rank one case with coefficients and any K} and~\ref{prop:
  ramification bound in rank one case with coefficients
    and any K with GK action}.)
\begin{defn}\label{defn: Kbasic T height at most h}
  Let~$A$ be a finite local Artinian $\cO/\varpi^a$-algebra for
  some~$a\ge 1$, and let~$M$ be a rank one projective \'etale $\varphi$-module
  over~$\A_{K,A}$. 
We regard~$M$ as a rank~$[K:\Kbasic]$ \'etale $\varphi$-module~$M'$
  over~$\A_{\Kbasic,A}$, and we say that~$M$ \emph{has
    $(\Kbasic,T)$-height at most~$h$} if there is a projective
  $\varphi$-module~$\gM$ over~$\A_{\Kbasic,A}^+$ of $T$-height at
  most~$h$ such that $M'=\gM[1/T]$.
\end{defn}

We next turn to describing $\cR^{\Gamma_\disc}_1$ explicitly.
We first note that
since~$I_K\cap G_{\Kcyc}=I_{\Kcyc}$, we have an exact sequence
$1 \to I_{K_{\cyc}} \to G_K \to \Frob_K^{\Zhat} \times \Gamma.$
We let $H \subseteq
\Frob_K^{\Zhat} \times \Gamma$
denote the image of $G_K$;
it is an open subgroup of 
$\Frob_K^{\Zhat} \times \Gamma$.
We write $H_{\disc} := H \cap (\Frob_K^{\Z} \times \Gamma_{\disc});$
then $H_{\disc}$ is dense in $H$, and is isomorphic to $\Z\times \Z$.
Write $\cO[H_{\disc}]$ for the group ring of $H_{\disc}$ over $\cO$,
and $\widehat{\cO[H_{\disc}]}$ for its $\varpi$-adic completion.
Then $\Spec \cO[H_{\disc}] \cong \Gm\times_\cO \Gm,$
and $\Spf \widehat{\cO[H_{\disc}]} \cong \widehat{\mathbb G}_m \times_{\cO}
\widehat{\mathbb G}_m.$  Let $I'$ denote the image of $I_{\Kcyc}$ in $G_K^{\ab}$;
equivalently, $I'$ is the quotient of $I^{\ab}_{\Kcyc}$ by
the closure of its subgroup of commutators $[ I^{\ab}_{\Kcyc},\Gamma]$.

\begin{prop}
	\label{prop:rank 1 R-Gamma-disc}
       	There is an isomorphism
	$$
\Bigl[ \Bigl ( \Spf \cO[[I']] \times \Spf \widehat{\cO[H_{\disc}]}
	\Bigr) / \widehat{\mathbb G}_m \Bigr]
\iso
	\cR_1^{\Gamma_{\disc}} $$
{\em (}where, in the formation of the quotient stack, the $\widehat{\mathbb
G}_m$-action is taken to be trivial{\em )}.
\end{prop}
\begin{proof}
We leave the construction of this isomorphism to the reader,
by combining the isomorphism of Proposition~\ref{prop:rank 1 R}
with the definition of $\cR_1^{\Gamma_{\disc}}$.
\end{proof}

\begin{remark}
	Note that the case $d=1$ is not representative of the general
	case in so far as the structure of $\cR$, and so also of
	$\cR^{\Gamma_{\disc}}$, is concerned.  Namely, when $d = 1$,
	the stacks $\cR$ and $\cR^{\Gamma_{\disc}}$ are formal
	algebraic stacks.  This will not be the case when $d > 1$.
	(They will instead be Ind-algebraic stacks whose underlying
	reduced substacks are Ind-algebraic but not algebraic.)
	The reason for this is that the $1$-dimensional mod $p$
	representations of $G_{K_{\cyc}}$ may be described
	as Frobenius twists of a finite number of characters
	when $d = 1$,
	so that in this case $\cR_{\red}$ is indeed an algebraic stack,
	while the spaces $\Ext^1_{G_{K_{\cyc}}}(
	\chi,\psi)$ are infinite dimensional, for any
	two mod $p$ characters $\chi$ and $\psi$ of $G_{K_{\cyc}}$,
	so that already in the case $d = 2$,
	the stack $\cR_{\red}$ is merely Ind-algebraic.
\end{remark}

We now give an explicit description of the stack~$\cX_1$. 
\begin{prop}
	\label{prop:rank 1 X}
       	There is an isomorphism
	$$ \Bigl[ \Bigl ( \Spf \bigl( \cO[[I_K^{\ab}]]\bigr) \times \widehat{\mathbb G}_m
	\Bigr) / \widehat{\mathbb G}_m \Bigr] \iso \cX_1$$
{\em (}where, in the formation of the quotient stack, the $\widehat{\mathbb
G}_m$-action is taken to be trivial{\em )}.
\end{prop}
\begin{proof}
We prove this in the same way as
  Proposition~\ref{prop:rank 1 R}; we leave the details to the reader,
  indicating only the key differences. We can construct a morphism
  $\Spf \bigl( \cO[[I_K^{\ab}]]\bigr) \times \widehat{\mathbb G}_m\to\cX_1$ in
  exactly the same way as in the proof of Proposition~\ref{prop:rank 1
    R} (by choosing a lift $\sigma_K \in G_{K}$ of $\Frob_K$,
  splitting the short exact sequence~\eqref{eqn:K ramification s.e.s.}),
  and we need to prove that this morphism is representable by
  algebraic spaces. This can be done by arguing exactly as in the
  proof of Proposition~\ref{prop:rank 1 R}, with the isomorphism
  $\cX_{K,1}^a=\varinjlim_{h,s}\cX_{K,1,\Kbasic,h,s}^a$ of
  Lemma~\ref{lem: explicit Ind description for X pulled back from K
    basic}  replacing the isomorphism
  $\cR_{K,1}^a=\varinjlim_h\cR_{K,1,\Kbasic,h}^a$. 

  Bearing in mind the definition of~$\cX_{K,1,\Kbasic,h,s}^a$, we find that 
  we have to show that
  there is an open subgroup~$I_K^{h,s,a}$ of~$I_K^{\ab}$ such that if~$A$ is a finite
  local Artinian $\cO/\varpi^a$-algebra, and $M$ is a rank one \'etale
  $\varphi$-module over~$\A_{K,A}$ with the property that there is a
  rank~$[K:\Kbasic]$ $\varphi$-module $\gM$ over $\A_{\Kbasic,A}^+$, such
  that:
  \begin{itemize}
  \item $\gM[1/T]=M$ (where we are regarding~$M$ as an \'etale
    $\varphi$-module over~$\A_{\Kbasic,A}$),
  \item $\gM$ is of $T$-height at most~$h$, and
  \item the action of~$\Gamma_{\Kbasic,\disc}=\Gamma_{K,\disc}$ on $\gM[1/T]$ extends the canonical action of
    $\Gamma_{\Kbasic_s,\disc}$ given by Corollary~\ref{cor: Caruso Liu Galois
      action phi Gamma discrete version},
  \end{itemize}
then the action of~$I_K^{h,s,a}$ on~$T_A(M)$ is trivial. This is 
immediate from Proposition~\ref{prop: ramification bound in rank one case with coefficients
    and any K with GK action}.
\end{proof}

\begin{rem}
  \label{rem: rank one X not closed substack}
  Note that Propositions~\ref{prop:rank 1 R-Gamma-disc} and
  \ref{prop:rank 1 X} describe $\cX$ as a certain formal completion
  of $\cR^{\Gamma_{\disc}}$ (in the case $d = 1$). 
  In particular, the monomorphism
  $\cX \hookrightarrow \cR^{\Gamma_{\disc}}$ is not a closed immersion
  (in the case $d =1$, and presumably not in the general case either).
This helps to explain why
somewhat elaborate arguments were required in
Chapter~\ref{sec: phi modules and phi gamma modules}
to deduce properties of $\cX$ from the corresponding properties of
$\cR^{\Gamma_{\disc}}$.
\end{rem}

\begin{rem}\label{rem: rank 1 versus rank d}
	It follows easily from Proposition~\ref{prop:rank 1 X} that the stack $\cX_1$ may be described
	as a moduli stack of $1$-dimensional continuous
	representations of the Weil group~$W_K$; indeed, we see that if
        $A$ is a $p$-adically complete $\cO$-algebra, then $\cX_1(\Spf
        A)$ is the groupoid of 
        continuous characters~$I_K^{\ab}\times \Z \to
        A^\times$ (if we identify $\widehat{\mathbb G}_m$ with the $p$-adically
completed group ring of $\Z$). 
Mapping the generator $1 \in \Z$ to $\sigma_K$ (a lift of Frobenius, 
as in the proof of Proposition~\ref{prop:rank 1 X})
yields an isomorphism $I_K^{\ab} \times \Z \iso W_K^{\ab}$.


        This phenomenon
	doesn't persist in the general case.  Indeed,
	as noted in Section~\ref{subsec: families of extensions}, 
	already in the case $d = 2$, there seems to be no description
        of $\cX_2$ as the moduli space of representations of 
	a group that is compatible with the description
       of its closed points in terms of representations of~$G_K$. 

One can attempt
to adapt the proof of
Proposition~\ref{prop:rank 1 X} to the case~$d>1$, by constructing a
morphism from the moduli stack of continuous $d$-dimensional Weil--Deligne representations
(suitably understood\footnote{Since we are considering representations
in rings that are $p$-power torsion,
we cannot use the usual formulation of the Weil--Deligne group.  Rather,
we choose a finitely generated dense subgroup
$\Z \ltimes \Z[1/p]$ 
of the tame Galois group of $G_K$ (by choosing a lift of Frobenius and a lift
of a topological generator of tame inertia), and form the topological
group $WD_K$ by taking its preimage in $G_K$ (and equipping 
$\Z \ltimes \Z[1/p]$ with its discrete topology).})
to the stack~$\cX_d$.
However, when $d > 1$,
this morphism of stacks will not be representable by
algebraic spaces. The precise point where the proof that we've given 
in the rank $1$ case fails to
generalise 
is that once we are not passing to abelianized Galois
groups, the upper numbered ramification
groups are not open in the inertia group.
\end{rem}

\section{A ramification bound}
\label{subsec:ram bound}
We end this chapter by proving the bounds on the ramification of a character valued in a  finite
Artinian $\cO$-algebra 
that were used in the proofs of
Propositions~\ref{prop:rank 1 R} and~\ref{prop:rank 1 X}. 
There are well-established techniques for
proving such a bound, going back to~\cite{MR807070}. We find it
convenient to follow the arguments of~\cite{MR2745530}, and indeed we
will follow some of their arguments very closely. The main differences
between our setting and that of~\cite{MR2745530} are that we are
working in the cyclotomic setting, rather than the Kummer setting;
that we are considering representations of~$G_{K_{\cyc}}$ of finite
height, rather than semistable representations of~$G_K$; and that to
define a representation of finite height, we need to pass between~$K$
and~$\Kbasic$. None of these changes make a fundamental difference
to the argument.

In fact, while~\cite{MR2745530} go to some effort to work modulo an
arbitrary power of~$p$, and to optimise the bounds that they obtain,
we are content to give the simplest proof that we can of the existence
of a bound of the kind that we need. In particular, we are able to
reduce to the case of mod~$p$ representations, 
which simplifies much of the discussion.

Recall that~$K_{\cyc}/K$ is a Galois extension with Galois
group~$\Gamma_K\cong\Z_p$. For each~$s\ge 0$, we write~$K_{\cyc,s}$
for the unique subfield of~$K_{\cyc}$ which is cyclic over~$K$ of
degree~$p^s$. We write~$e=e(K/K_0)$.

Suppose that~$K$ is 
basic, and that ~$\gM$ is a free $\varphi$-module over $\A^+_{K}/p$ of
$T$-height at most~$h$.
 Then~$M:=\A_{K}\otimes_{\A^+_{K}}\gM$ is an
\'etale $\varphi$-module over~$\A_{K}/p$, and we have a corresponding
$G_{\Kcyc}$-representation $T(M)$ as in Section~\ref{subsec:Galois
  reps}. (Note that while elsewhere in the book we have worked with
$\cO/\varpi^a$-coefficients, it is more convenient to use
$\Z/p\Z$-coefficients in most of this section.) 

In order to compare directly to the
arguments of~\cite{MR2745530}, it is more convenient to work
 with the contragredient $G_{K_{\cyc}}$-representation $T(M)^{\vee}$,
which can also be computed as
$T(M^{\vee}) = (\C^\flat\otimes_{\A_K} M^{\vee})^{\varphi=1}= \Hom_{\A_K,\varphi}\bigl(M,{\C}^\flat\bigr)= \Hom_{\A^+_K,\varphi}\bigl(M,{\C}^\flat\bigr).$
By~\cite[A.1.2.7, B.1.8.3]{MR1106901}, this can also be computed
via the functor 
\[T^{\vee}(\gM):=\Hom_{\A^+_K,\varphi}\bigl(\gM,\O_{\C}^\flat\bigr).\]
More precisely, there is a diagram
$$\xymatrixcolsep{-0.1in}
\xymatrix{
T^{\vee}(\gM):=
\Hom_{\A^+_K,\varphi}\bigl(\gM,\O_{\C}^\flat\bigr)
\ar[rd]&&
\Hom_{\A^+_K,\varphi}\bigl(M,{\C}^\flat\bigr) = T(M^{\vee})\ar[ld] \\
&\Hom_{\A^+_K,\varphi}\bigl(\gM,{\C}^\flat\bigr)&}$$
(the first arrow being induced by the inclusion of 
$\O_{\C}^\flat$ in ${\C}^\flat$, and the second by
restriction from $M$ to~$\gM$),  with both arrows in fact being
isomorphisms (by \emph{op.\ cit.}).



We will now follow the arguments of~\cite[\S 4]{MR2745530}; note that
the ring denoted~$R$ in~\cite{MR2745530} is~$\cO_\C^\flat$. For the
time being we will assume that we are in the following situation; the
rather complicated hypotheses here will allow us to interpret the
canonical actions of Corollary~\ref{cor:
      Caruso Liu Galois action phi Gamma discrete version} in terms of the constructions of~\cite[\S 4]{MR2745530}.

\begin{hyp}\label{hyp: K is basic and M has height h and s is big enough}
  Assume that~$K$ is basic, and fix some~$h\ge 1$. Let~$\gM$ be a free
  $\varphi$-module over $\A^+_{K}/p$ of $T$-height at most~$h$. Fix
  the following:
  \begin{itemize}
  \item an integer~$b\ge eh/(p-1)$. 
  \item An integer ~$N> pb/e$ such
    that $N\ge N(1,h)$, where $N(1,h)$ is as in Corollary~\ref{cor: Caruso
      Liu Galois action phi Gamma discrete version}.  
  \item An integer~$s$ such that~$p^{s-1}>b/e$, and $s\ge s(1,h,N)$,
    where $s(1,h,N)$ is as in Corollary~\ref{cor: Caruso
      Liu Galois action phi Gamma discrete version}. 
  \end{itemize}
\end{hyp}
Note in  particular that under Hypothesis~\ref{hyp: K is basic and M has height h and s is big enough} we have $p^s> ph/(p-1)>h$.

We normalize the valuation~$v$ on~$K$ so
that~$v(K^\times)=\Z$, 
and continue to write~$v$ for the unique compatible
valuation on~$\overline{K}$, and for the induced
valuation on~$\cO_\C^\flat=\varprojlim\cO_{\overline{K}}/p$. Note that
with this convention we have~$v(T)=v(p)=e$ (where we are
abusively continuing to denote the image of~$T$ in~$\cO_\C^\flat$
by~$T$; this follows
straightforwardly from the definition of~$T$ as the trace of
$T' := ([\varepsilon]-1)$, which shows (after a simple computation)
that $v(T) = (p-1)v(T')$, along with the formula
$v(\zeta_{p^s}-1)=v(p)/p^{s-1}(p-1)$, which shows that $v(T') = v(p)/(p-1)$).
We write $\mf{a}_{\cO_{\C}^\flat}^{>
  c}$ for the set of elements of~$\cO^\flat_{\C}$
of valuation
 greater than~$c$, and define $\mf{a}_{\cO_{\C}^\flat}^{\ge
  c}$ in the same way. 
For any~$c\ge 0$,
we
write 
\[T_{\cO_\C^\flat,c}^\vee(\gM):=\Hom_{\A^+_K,\varphi}\bigl(\gM,\cO_\C^\flat/\mf{a}_{\cO_\C^\flat}^{>c}\bigr),\]
so that there is a natural map $T^\vee(\gM)\to T^\vee_{\cO_\C^\flat,c}(\gM)$. For any~$c'\ge
c$ 
we write 
\[T_{\cO_\C^\flat,c',c}^\vee(\gM)=\Im(T_{\cO_\C^\flat,c'}^\vee(\gM)\to T_{\cO_\C^\flat,c}^\vee(\gM)),\]
where the    morphism is induced by the natural map $\mf{a}_{\cO_\C^\flat}^{>c'}\to
    \mf{a}_{\cO_\C^\flat}^{>c}$. 

\begin{lem}
  \label{lem: truncation for Ainf mod p representation}Assume that we are in the setting of Hypothesis~{\em\ref{hyp: K
    is basic and M has height h and s is big enough}}. Then 
the induced morphism
$T^\vee(\gM)\to T^\vee_{\cO_\C^\flat,b}(\gM)$
is injective, with image~$T^\vee_{\cO_\C^\flat,pb,b}(\gM)$.
\end{lem}
\begin{proof}


  This can be proved by a standard Frobenius amplification argument
  exactly as in the proof 
  of~\cite[Prop.\ 2.3.3]{MR2745530}. For example, to see the
  injectivity, we note that since~$\gM$ is finitely generated, any
  element of the kernel corresponds to a $\varphi$-equivariant
  morphism $\gM\to \mf{a}_{\cO_\C^\flat}^{\ge c}$ for
  some~$c>b$. Since ~$\gM$ has $T$-height at most~$h$, and
  $v(T^h)=eh$, the $\varphi$-equivariance implies that the image of
  this morphism is
  contained in $\mf{a}_{\cO_\C^\flat}^{\ge pc-eh}$. Since $c>b\ge
  eh/(p-1)$ by hypothesis, we see that the sequence
  $c,pc-eh,p(pc-eh)-eh,\dots$ goes to infinity, and the injectivity follows.

  We leave surjectivity to the reader; it is a simpler version of the
  proof of Lemma~\ref{lem: amplification from char p to char 0}
  below. (See also the proof of Lemma~\ref{lem: Frobenius
    amplification lemma over Ainf} for an almost identical
  argument.)\end{proof}

Suppose given $c \geq 0$ satisfying $p^s > c/e$.
(E.g.\ $b$ and $pb$ both satisfy this condition.)
The Frobenius $\varphi$  is an automorphism of the perfect $k$-algebra
$\cO_\C^\flat$ which multiplies valuations by~$p$,
and so $\varphi^{-s}$ 
induces the first of the following sequence of  $G_K$-equivariant 
%
%
  isomorphisms of $k$-algebras \numequation\label{eqn: isomorphism from
    truncated Ainf to truncated Kbar}
  \cO_\C^\flat/\mf{a}_{\cO_\C^\flat}^{>c}\isoto
 k\otimes_{k,\varphi^s}\cO_{\C}^\flat/\mf{a}_{\C}^{>c/p^s} 
\iso
 k\otimes_{k,\varphi^s}\cO_{\overline{K}}/\mf{a}_{\overline{K}}^{>c/p^s} ,
\end{equation}
where $\mf{a}_{\overline{K}}^{>c/p^s}$ has the 
evident meaning, and the second isomorphism
is induced by projection onto the first factor in the projective
limit
$\O_\C^\flat  = \varprojlim_{\varphi} \O_\C/p =
\varprojlim_\varphi\cO_{\overline{K}}/p$.
(The upper bound on $c$ ensures that $p \in \mf{a}_{\overline{K}}^{>c/p^s}$,
so that this projection does indeed induce the indicated morphism; it is then
straightforward to verify that it  is an isomorphism.)
This isomorphism in turn induces an isomorphism
\numequation
\label{eqn:alternative description of Tc}
T^{\vee}_{\cO_C^\flat,c}(\gM) \iso
\Hom_{\A_K^+,\varphi}(\gM, \cO_{\overline{K}}/\mf{a}_{\overline{K}}^{>c/p^s}).
\end{equation}

\begin{lem}
\label{lem:G K s action}
If $p^s > c/e \geq 0,$ then the action of $G_{K_{\cyc,s}}$ on
$\cO_\C^\flat$
induces an action of $G_{K_{\cyc,s}}$
on $T^{\vee}_{\cO_C^\flat,c}(\gM).$ 
\end{lem}
\begin{proof}
The action of $\varphi$ on
$\cO_\C^\flat$
is $G_K$-equivariant, 
and so the $G_K$-action on
$\cO_\C^\flat$
induces an action on
$\Hom_{\varphi}\bigl(\gM,\cO_\C^\flat/\mf{a}_{\cO_\C^\flat}^{>c}\bigr)$. 
The claim of the lemma is that the restriction
of this action to $G_{K_{\cyc,s}}$ preserves
the subobject
$T_{\cO_\C^\flat,c}^\vee(\gM):=\Hom_{\A^+_K,\varphi}\bigl(\gM,\cO_\C^\flat/\mf{a}_{\cO_\C^\flat}^{>c}\bigr)$
of
$\Hom_{\varphi}\bigl(\gM,\cO_\C^\flat/\mf{a}_{\cO_\C^\flat}^{>c}\bigr)$.

Rather than proving this directly,
we use the
isomorphisms~\eqref{eqn: isomorphism from truncated Ainf to truncated Kbar}
and~\eqref{eqn:alternative description of Tc}; taking these into account,
it suffices to prove that 
the  action of $G_{K_{\cyc,s}}$ on
$ \Hom_{\varphi}(\gM, \cO_{\overline{K}}/\mf{a}_{\overline{K}}^{>c/p^s})$
preserves the subobject
$\Hom_{\A_K^+,\varphi}(\gM, \cO_{\overline{K}}/\mf{a}_{\overline{K}}^{>c/p^s})$.
In other words, we have to check that the $G_{K_{\cyc,s}}$-action
preserves the property of being $\A_K^+$-equivariant.
For this, it suffices to check that $G_{K_{\cyc,s}}$ acts
trivially on the image of $T$ in $\cO_{\overline{K}}/\mf{a}_{\overline{K}}^{>c/p^s}$
under the isomorphism~\eqref{eqn: isomorphism from truncated Ainf to truncated Kbar}.
Since $T$ can be expressed as a power series in $T'$, it in fact suffices
to check this for the image of $T'$.  But this image is equal to
$\zeta_{p^{s+1}} - 1 \bmod \mf{a}_{\overline{K}}^{> c/p^s}$,
which  {\em is} fixed by $G_{K_{\cyc,s}}$ (since $\zeta_{p^{s+1}} \in K_{\cyc,s}$).
\end{proof}


\begin{cor}\label{cor: we get a canonical action of G K s}Assume that we are in the setting of Hypothesis~{\em\ref{hyp: K
    is basic and M has height h and s is big enough}}. Then the natural
  action of~$G_{K_{\cyc,s}}$
  on~$\cO_\C^\flat$ 
induces an
  extension of the action of~$G_{\Kcyc}$ on $T^\vee(\gM)$ to an action
  of~$G_{K_{\cyc,s}}$. 
\end{cor}
\begin{proof}
This follows directly from
Lemmas~\ref{lem: truncation for Ainf mod p representation}
and~\ref{lem:G K s action},
once we note that, by hypothesis, $p^s > pb/e > b/e.$
\end{proof}
By Corollary~\ref{cor: Caruso Liu Galois action phi Gamma discrete
  version} (and our hypothesized bound on~$s$), 
we can give $\gM[1/T]$ the structure of a
$(\varphi,\Gamma_{K_{\cyc,s}})$-module; this in particular induces an action
of~$G_{K_{\cyc,s}}$ on $T(M)$. 
The following corollary
expresses the compatibility of this action with the action constructed
in Corollary~\ref{cor: we get a canonical action of G K s}.

\begin{cor}
  \label{cor: the various canonical actions agree}
Assume that we are in the setting of Hypothesis~{\em\ref{hyp: K
    is basic and M has height h and s is big enough}}. Then the two actions of
  $G_{K_{\cyc,s}}$ on $T(M)^\vee$ \emph{(}from Corollary~\emph{\ref{cor: we get a
    canonical action of G K s}}, and from Corollary~\emph{\ref{cor:
      Caruso Liu Galois action phi Gamma discrete version}}\emph{)} agree.
\end{cor}
\begin{proof}By the proof of Corollary~\ref{cor: Caruso Liu Galois action phi Gamma discrete
  version}, the semilinear action of $\Gamma_{K_{\cyc,s}}$ on~$\gM$
extends to a semilinear
action of $G_{K_{\cyc,s}}$ on $\cO_{\C}^\flat\otimes_{\A^+_K}\gM$,
which is uniquely determined by the properties that it commutes
with~$\varphi$, and satisfies $(g-1)(\gM)\subset
T^N(\cO_{\C}^\flat\otimes_{\A^+_K}\gM)$ for all~$g\in
G_{K_{\cyc,s}}$. As explained in Section~\ref{subsubsec:
Galois reps for GK phi}, we have
$T(M)=(\C^\flat\otimes_{\A^+_K}\gM)^{\varphi=1}$, with the action of $G_{K_{\cyc,s}}$ on
$T(M)$ being that inherited from its action on
$\cO_{\C}^\flat\otimes_{\A^+_K}\gM$. 

As noted above, we also have
\begin{multline*}
T(M)^{\vee} = T(M^{\vee}) = T^{\vee}(\gM)
\\
=\Hom_{\A^+_K,\varphi}\bigl(\gM,\O_{\C}^\flat\bigr)
=\Hom_{{\C}^\flat,\varphi}\bigl({\C}^\flat\otimes_{\A^+_K}\gM,{\C}^\flat\bigr).
\end{multline*}
Now evaluation induces a pairing 
$$\bigl({\C}^\flat\otimes_{\A^+_K}\gM \bigr)\times
\Hom_{{\C}^\flat}\bigl({\C}^\flat\otimes_{\A^+_K}\gM,{\C}^\flat\bigr)
\to 
\C^\flat,$$
which restricts to a pairing
\numequation
\label{eqn:perfect pairing}
({\C}^\flat\otimes_{\A^+_K}\gM)^{\varphi = 1} \times
\Hom_{{\C}^\flat,\varphi}\bigl({\C}^\flat\otimes_{\A^+_K}\gM,{\C}^\flat\bigr)
\to 
(\C^\flat)^{\varphi = 1} = \Z/p\Z,
\end{equation}
and this latter pairing {\em is} the natural pairing of $G_{K_{\cyc}}$-representations
\[T(M) \times T(M)^{\vee} \to \Z/p\Z.\]
Therefore, to prove the lemma, we have to show that~\eqref{eqn:perfect pairing}
is furthermore $G_{K_{\cyc,s}}$-equivariant (with the action on~$T(M)$
being the action recalled above, coming from Corollary~\ref{cor: Caruso Liu Galois action phi Gamma discrete
  version}, and the action on~$T(M)^\vee=T^{\vee}(\gM)$ given by Corollary~\ref{cor: we get a canonical action of G K s}).

To see this, note firstly that by the same argument that we used in
the proof of Lemma~\ref{lem: truncation for Ainf mod p
  representation}, it follows from Hypothesis~\ref{hyp: K is basic and
  M has height h and s is big enough}
that \[(\C^\flat\otimes_{\A^+_K}\gM)^{\varphi=1}\subseteq (\mf{a}^{\ge b})^{-1}\otimes_{\A^+_K}\gM.\]
In addition, we
have \begin{multline*}\Hom_{{\C}^\flat,\varphi}\bigl({\C}^\flat\otimes_{\A^+_K}\gM,{\C}^\flat\bigr)=\Hom_{\cO_{\C^\flat},\varphi}\bigl(\cO_{\C^\flat}\otimes_{\A^+_K}\gM,\cO_{\C^\flat}\bigr)\\
  \subseteq
  \Hom_{\cO_{\C^\flat}}\bigl(\cO_{\C^\flat}\otimes_{\A^+_K}\gM,\cO_{\C^\flat}\bigr)
  \isofrom\Hom_{\cO_{\C^\flat}}\bigl((\mf{a}^{\ge
    b})^{-1}\otimes_{\A^+_K}\gM,(\mf{a}^{\ge
    b})^{-1}\bigr)\end{multline*} with the last
isomorphism being induced by restricting the domain of a homomorphism from $(\mf{a}^{\ge
    b})^{-1}\otimes_{\A^+_K}\gM$ to
  $\cO_{\C^\flat}\otimes_{\A^+_K}\gM$. Bearing in mind this
  isomorphism, we obtain an evaluation pairing \[\bigl((\mf{a}^{\ge
      b})^{-1}\otimes_{\A^+_K}\gM\bigr) \times \Hom_{\cO_{\C^\flat}}\bigl((\cO_{\C^\flat}\otimes_{\A^+_K}\gM,\cO_{\C^\flat}\bigr)\to (\mf{a}^{\ge
      b})^{-1} \](and taking $\varphi$-invariants again recovers the pairing~\eqref{eqn:perfect pairing}).

  Fix some $\alpha\in\C^\flat$ with $v(\alpha)=-b$, so that
  $(\alpha)=(\mf{a}^{\ge b})^{-1}$. Let~$x$ be an element of
  $({\C}^\flat\otimes_{\A^+_K}\gM)^{\varphi = 1}$, and write
  $x=\alpha\sum_{i=1}^n\lambda_i\otimes x_i$ where
  $\lambda_i\in\cO_{\C}^\flat$ and $x_i\in\gM$. Let
  $ f\in
  \Hom_{\cO_{\C}^\flat,\varphi}\bigl(\cO_{\C}^\flat\otimes_{\A^+_K}\gM,\cO_{\C}^\flat\bigr)$,
  so that $x$ and $f$ pair to
  $\alpha f(\sum_{i=1}^n\lambda_ix_i)\in (\mf{a}^{\ge b})^{-1}$; of
  course, they in fact pair to an element
  of~$\Fp\subset\cO_\C^\flat \subset (\mf{a}^{\ge b})^{-1}$.
  Let~$g\in G_{\Kcyc,s}$ be arbitrary. Then $gx$ and $gf$ pair to
  $g(\alpha) (gf)(\sum_{i=1}^ng(\lambda_i)g(x_i))$. To establish the
  claimed $G_{\Kcyc,s}$-equivariance, we need to show that
  \[g(\alpha)\cdot (gf)\bigl(\sum_{i=1}^ng(\lambda_i)g(x_i)\bigr)=g\bigl(\alpha
    f\bigl(\sum_{i=1}^n\lambda_ix_i\bigr)\bigr).\] This is an equality of
  elements of~$\Fp$, and it therefore suffices to show that it holds
  modulo~$\m_{\cO_{\C}}$ when each side is considered as an element
  of~$\cO_{\C}^\flat$. Equivalently, after multiplication
  by~$\alpha^{-1}$ it suffices to show
  that
  \[(gf)\bigl(\sum_{i=1}^ng(\lambda_i)g(x_i)\bigr)\equiv
    g\bigl(f\bigl(\sum_{i=1}^n\lambda_ix_i\bigr)\bigr) \pmod{\mf{a}^{>
        b}}.\] Recall that
  $(g-1)(\gM)\subset T^N(\cO_{\C}^\flat\otimes_{\A^+_K}\gM)$. By
  Hypothesis~\ref{hyp: K is basic and M has height h and s is big
    enough}, we have $v(T^N)=eN>pb$, so
  that in particular we have~$T^N\in \mf{a}_{\cO_{\C}^\flat}^{>b}$. Thus
  $g(x_i)\equiv x_i\pmod{\mf{a}^{> b}}$, so it suffices to show
  that
  \[(gf)\bigl(\sum_{i=1}^ng(\lambda_i)x_i\bigr)\equiv
    g\bigl(f\bigl(\sum_{i=1}^n\lambda_ix_i\bigr)\bigr) \pmod{\mf{a}^{>
        b}},\] or equivalently
  that
  \[\sum_{i=1}^ng(\lambda_i)(gf)(x_i)\equiv
    \sum_{i=1}^ng(\lambda_i)g(f(x_i)) \pmod{\mf{a}^{> b}},\] so in
  turn it is enough to show that if $x\in\gM$ then
  $(gf)(x)\equiv g(f(x))\pmod{\mf{a}^{>b}}$. But this is true by the
  very definition of~$gf$ (which was defined via the action
  of~$G_{\Kcyc,s}$ on~ $T_{\cO_\C^\flat,b}^\vee(\gM)$, which in turn
  is defined via the action of~$G_{\Kcyc,s}$ on
  $\cO_{\C}^\flat/\mf{a}^{>b}$), so we are done.
\end{proof}


Recall that we write~$k$ for the residue field
of~$K$, and since~$K$ is assumed basic we have $k_\infty=k$, so that we can
and do regard~$\gM$ as a free $\A_{K}^+/p=k[[T]]$-module. 
Write~$T_s$ for
$\tr_{K(\zeta_{p^{s+1}})/K_{\cyc,s}}(\zeta_{p^{s+1}}-1)\in\cO_{K_{\cyc,s}}$,
and note  that  we have a
homomorphism \numequation\label{eqn: A to Ks}\A_{K}^+/p=k[[T]]\to
  \cO_{K_{\cyc,s}}/p\end{equation} which sends~$T\to T_s$, and whose
restriction to~$k$ is given by $\varphi^{-s}$.

Let~$L$ be any algebraic extension of~$K_{\cyc,s}$
inside~$\overline{K}$. We 
continue to write~$v$ for the unique 
valuation on~$L$ with $v(p)=e$.  Note that we have~$v(T_s^h)=eh/p^s<e=v(p)$. For
any~$c\ge 0$ we write $\mf{a}_L^{> c}$  for the set of elements of~$L$
of valuation greater than~$c$, and  define  $\mf{a}_L^{\ge c}$ in the same way.  

Choose an (ordered) $\A_{K}^+/p$-basis $e_1,\dots,e_d$ of~$\gM$,
and 
write $\varphi(e_1,\dots,e_d)=(e_1,\dots,e_d)A$ for some~$A\in
M_d(\A_{K}^+/p)$. By our assumption that~$\gM$ has $T$-height at most~$h$,
we can write $AB=T^h\Id_d$ for some $B\in M_d(\A_{K}^+/p)$. We can
and do choose matrices $\widetilde{A},\tB\in M_d(\cO_{K_{\cyc,s}})$ which
respectively lift the images of~$A,B$ under the homomorphism
$\A_{K}^+/p\to\cO_{K_{\cyc,s}}/p$ defined above. Since the image
of~$T_s^h$ in $\cO_{K_{\cyc,s}}/p$ is nonzero, we see that after
possibly multiplying~$\tB$ by an invertible matrix which is trivial
mod~$p/T_s^h$,
we can and do assume that \numequation\label{eqn: Atilde
  times Btilde}\widetilde{A}\tB=T_s^h\Id_d.\end{equation}

We
set~\[\tT_L^*(\gM):=\{(x_1,\dots,x_d)\in\cO_L^d\mid
  (x_1^p,\dots,x_d^p)=(x_1,\dots,x_d)\widetilde{A}\}.\]
Note that since~$x\mapsto x^p$ is not a ring homomorphism on~$\cO_L$,
this is just a set with a $G_{K_{\cyc,s}}$-action, rather than
an~$\cO_L$-module.

For any~$c\in [0,ep^s)$, 
we
write 
\[T_{L,c}^*(\gM):=\Hom_{\A^+_K/p,\varphi}(\gM,(\cO_L/p)/\mf{a}_L^{>c/p^s}),\]where
the $\A^+_K$-algebra structure on $\cO_L/p$ is that induced by the
composite $\A^+_K/p \buildrel \text{\eqref{eqn: A to Ks}}\over\longrightarrow \cO_{K_s}/p\to\cO_L/p$. (The use of~$c/p^s$
rather than~$c$ in the definition of~$T_{L,c}^*$ comes from the
definition of the homomorphism~\eqref{eqn: A to Ks}; see also~\eqref{eqn: isomorphism from
    truncated Ainf to truncated Kbar}.) For any~$c'\ge
c$ 
we write 
\[T_{L,c',c}^*(\gM)=\Im(T_{L,c'}^*(\gM)\to T_{L,c}^*(\gM)),\]where the
    morphism is induced by the natural map 
    $\mf{a}_L^{>c'/p^s}\to
    \mf{a}_L^{>c/p^s}$. 

Evaluating on the members
of the basis~$e_1,\dots,e_d$ of~$\gM$ yields an
identification
\[T_{L,c}^\vee(\gM)\iso\{(x_1,\dots,x_d)\in((\cO_L/p)/\mf{a}_L^{>c/p^s})^d\mid
  (x_1^p,\dots,x_d^p)=(x_1,\dots,x_d)\overline{A}\} \]
(where~$\overline{A}$ denotes the image of~$A$
in~$M_d((\cO_L/p)/\mf{a}_L^{>c/p^s})$), so that there is a natural map
$\tT_L^\vee(\gM)\to T_{L,c}^\vee(\gM)$.

The following lemma is the analogue of~\cite[Lem.\ 4.1.4]{MR2745530}
in our setting, and the proof is essentially identical (and in fact
simpler, since we are only working modulo~$p$).
\begin{lem}\label{lem: amplification from char p to char 0}Assume that we are in the setting of Hypothesis~{\em\ref{hyp: K
    is basic and M has height h and s is big enough}}. Then the morphism $\tT_L^\vee(\gM)\to T_{L,b}^\vee(\gM)$ is injective,
  and has image $T_{L,pb,b}^\vee(\gM)$. 
\end{lem}
\begin{proof}
  We begin with injectivity. Suppose that we have two distinct
  tuples~$(x_1,\dots,x_d)$ and $(y_1,\dots,y_d)\in \tT_L^*(\gM)$ whose
  images in $T_{L,b}^\vee(\gM)$ coincide; then if we
  write~$z_i:=y_i-x_i$, we have $v(z_i)> b/p^s$ for each~$i$. 
  By~\eqref{eqn: Atilde
  times Btilde}, we
have \[((x_1+z_1)^p-x_1^p,\dots,(x_d+z_d)^p-x_d^p)\tB=(z_1,\dots,z_d)T_s^h.\]Writing
$(x_i+z_i)^p-x_i^p=z_i^p+pz_i(x_i^{p-1}+\dots)$, we see that
if~$z:=\min_iv(z_i)$, then we have
$eh/p^s+z\ge\min(pz, e+z)$. Since~$p^s>h$, we see that necessarily
$e + z  > eh/p^s + z$.   Hence it must be that
$eh/p^s+z\ge pz$, and thus that $eh/p^s(p-1)\ge z$. (This final deduction
is where we use our assumption that the tuples $(x_1,\dots,x_d)$
and $(y_1,\ldots,y_d)$ are distinct; this assumption ensures that $z < \infty$,
so that it is legitimate to subtract if from the two sides of an inequality.) 
However this contradicts our assumption that~$z>b/p^s$.

We now turn to determining the image. By definition, the morphism $\tT_L^\vee(\gM)\to T_{L,b}^\vee(\gM)$ factors through $T_{L,pb,b}^\vee(\gM)$. Choose some $(\overline{x}_1,\dots,\overline{x}_d)\in
T^\vee_{L,pb}(\gM)$. It suffices to construct $(\tx_1,\dots,\tx_d)\in \tT_L^\vee(\gM)$
with~$\tx_i$ lifting the~$\overline{x}_i$ modulo~$\mf{a}^{>b/p^s}$ by
successive approximation.

Begin
by choosing  arbitrary lifts~$x_i$ of~$\overline{x}_i$ to~$\cO_L$.

We have
$T_s^{-h}\mf{a}_L^{>b/p^{s-1}}=\mf{a}_L^{>(pb-eh)/p^s}\subseteq
\mf{a}_L^{>b/p^s}$, so by our
assumption that $(\overline{x}_1,\dots,\overline{x}_d)\in
T^\vee_{L,pb}(\gM)$, we
have \[T_s^{-h}(x_1^p,\dots,x_d^p)\tB-(x_1,\dots,x_d)\in
  (T_s^{-h}\mf{a}_L^{>b/p^{s-1}})^d\subseteq (\mf{a}_L^{>b/p^s})^d,\]
so in fact we have
\[T_s^{-h}(x_1^p,\dots,x_d^p)\tB-(x_1,\dots,x_d)\in
  (\mf{a}_L^{\ge b'/p^{s}})^d\]for some~$b'> b$. We
will construct Cauchy sequences~$(z^{(j)}_i) \in \mf{a}_L^{\ge b'/p^s}$ whose limits~$z_i$ are
such that taking~$\tx_i:=x_i+z_i$ gives the required lift.

To this end, we set~$z_i^{(0)}=0$, and define 
\[(z_1^{(j+1)},\dots,z_d^{(j+1)}):=T_s^{-h}((x_1+z_1^{(j)})^p,\dots,(x_d+z_d^{(j)})^p)\tB-(x_1,\dots,x_d).\]
To see that $z_i^{(j+1)}\in \mf{a}_L^{\ge b'/p^s}$, we  write
\[
  \begin{split}
    (z_1^{(j+1)},\dots,z_d^{(j+1)}):=T_s^{-h}((x_1+z_1^{(j)})^p-x_1^p,\dots,(x_d+z_d^{(j)})^p-x_d^p)\tB
    \\+(T_s^{-h}(x_1^p,\dots,x_d^p)\tB-(x_1,\dots,x_d)).
  \end{split}
\]  so we need
only check that each $(x_i+z_i^{(j)})^p-x_i^p\in
T_s^{h}\mf{a}_L^{\ge b'/p^s}$. Since~$(x_i+z_i^{(j)})^p-x_i^p\in(p
z_i^{(j)},(z_i^{(j)})^p)$, and~$z_i^{(j)}\in\mf{a}_L^{\ge b'/p^s}$, this holds.  

It remains to show that each~$(z^{(j)}_i)$ is
Cauchy. Set~$\epsilon=\min(e-eh/p^s,((p-1)b'-eh)/p^s)$, so that in particular $\epsilon
>0$. It is enough to show that $\min_i
v(z_i^{(j+1)}-z_i^{(j)})\ge \min_i
v(z_i^{(j+1)}-z_i^{(j)})+\epsilon$. To this end, we have
  \begin{multline*}
    (z_1^{(j+1)}-z_1^{(j)},\dots,z_d^{(j+1)}-z_d^{(j)})\\=T_s^{-h}\bigl((x_1+z_1^{(j)})^p-(x_1+z_1^{(j-1)})^p,\dots,(x_d+z_d^{(j)})^p-(x_d+z_d^{(j-1)})^p\bigr)\tB,
  \end{multline*}
so that it is enough to show that for each~$i$, we have
$(x_i+z_i^{(j)})^p-(x_i+z_i^{(j-1)})^p\in
(z_i^{(j)}-z_i^{(j-1)})T_s^{h}\mf{a}_L^{\ge\epsilon}=(z_i^{(j)}-z_i^{(j-1)})\mf{a}_L^{\ge
  eh/p^s+\epsilon}$.

For each~$i$, we
write \[(x_i+z_i^{(j)})^p-(x_i+z_i^{(j-1)})^p=\sum_{k=1}^p\binom{p}{k}x_i^{p-k}((z_i^{(j)})^k-(z_i^{(j-1)})^k),\]so
it is enough to note that for $1\le k\le p-1$ we have
$v(\binom{p}{k})=e \ge eh/p^s+\epsilon$, while for~$k=p$, we
write \[(z_i^{(j)})^p-(z_i^{(j-1)})^p=(z_i^{(j)}-z_i^{(j-1)})((z_i^{(j)})^{p-1}+\dots+(z_i^{(j-1)})^{p-1}),\]
and note that since~$z_i^{(j)},z_i^{(j-1)}\in \mf{a}_L^{\ge b'/p^s}$, the
second term here is contained in $\mf{a}_L^{\ge (p-1)b'/p^s}$, and we have
$(p-1)b'/p^s\ge eh/p^s+\epsilon$.
\end{proof}

The following corollary is the analogue of~\cite[Thm.\
4.1.1]{MR2745530} in our setting.
\begin{cor}
  \label{cor: we get everything if and only if we contain the
    splitting field}Assume that we are in the setting of Hypothesis~{\em \ref{hyp: K
    is basic and M has height h and s is big enough}}, and let~$L$ be
  an algebraic extension of~$K_{\cyc,s}$ inside~$\overline{K}$, so
  that~$G_L$ acts naturally on~$T^\vee(\gM)$ by Corollary~{\em \ref{cor: we get a
    canonical action of G K s}}. 
  Then the natural injection
  $T^\vee_{L,pb,b}(\gM)\subseteq T^\vee_{\overline{K},pb,b}(\gM)$ is an
  isomorphism if and only if the action of~$G_{L}$ on~$T^\vee(\gM)$ 
  is trivial.
\end{cor}
\begin{proof}By Lemma~\ref{lem: amplification from char p to char 0}
  (applied to each of~$L$ and~$\overline{K}$), the inclusion
  $T^\vee_{L,pb,b}(\gM)\subseteq T^\vee_{\overline{K},pb,b}(\gM)$ is an
  isomorphism if and only if the same is true of the inclusion
  $\tT_L^\vee(\gM)\subseteq \tT_{\overline{K}}^\vee(\gM)$. Directly
  from the definition of these latter sets, we see that
  $\tT_L^\vee(\gM)=(\tT_{\overline{K}}^\vee(\gM))^{G_L}$, so this containment
  holds if and only if~$G_L$ acts trivially on
  $\tT_{\overline{K}}^\vee(\gM)$, or equivalently, if and only if~$G_L$ acts
  trivially on~$T^\vee_{\overline{K},pb,b}(\gM)$.

  Now Lemmas~\ref{lem: truncation for Ainf mod p representation}
and~\ref{lem:G K s action} show that 
$T^\vee(\gM) \iso T^\vee_{\overline{K},pb,b}(\gM),$
and by its construction,
the $G_{K_{\cyc,s}}$-action on
$T^\vee(\gM)$ arises, via this isomorphism,
from the $G_{K_{\cyc,s}}$-action on~$T^\vee_{\overline{K},pb,b}(\gM)$.
The lemma follows.
%
\end{proof}

\begin{cor}
  \label{cor: the slightly weird statement about lifting
    homomorphisms} Assume that we are in the setting of Hypothesis~{\em\ref{hyp: K
    is basic and M has height h and s is big enough}}, and let $\rho$
  denote the representation $T(M)^{\vee} = T^{\vee}(\gM)$
  of~$G_{K_{\cyc,s}}$ given by Corollary~{\em\ref{cor: we get a
    canonical action of G K s}}. Let~$L$ be
  an algebraic extension of~$K_{\cyc,s}$ inside~$\overline{K}$. If
  there exists an $\cO_{K_{\cyc,s}}$-algebra homomorphism
  $\cO_{\overline{K}^{\ker \rho}}\to \cO_L/\mathfrak{a}_L^{>b/p^{s-1}}$, 
  then  $\overline{K}^{\ker \rho}\subseteq L$.
\end{cor}
\begin{proof}Suppose that we have an  $\cO_{K_{\cyc,s}}$-algebra homomorphism
  $\eta:\cO_{\overline{K}^{\ker \rho}}\to
  \cO_L/\mathfrak{a}^{>b/p^{s-1}}$. We claim that~$\eta$ induces an injective
  homomorphism \[\cO_{\overline{K}^{\ker
      \rho}}/\mathfrak{a}_{\overline{K}^{\ker \rho}}^{>b/p^{s-1}}\to
  \cO_L/\mathfrak{a}_L^{>b/p^{s-1}}.\]

To  begin with,
let $F/K_{\cyc,s}$ be the maximal unramified subextension of $\overline{K}^{\ker
\rho}/K_{\cyc,s}$, and let $k_F$ denote its residue field (which
is then also the residue field of~$\overline{K}^{\ker \rho}$).
If $k_L$ denotes the residue field of~$L$, 
then the morphism $\eta$  induces an embedding $k_F \hookrightarrow
k_L$, which in turn lifts to an embedding $W(k_F) \hookrightarrow \cO_L$.
This embedding then induces an embedding
\numequation
\label{eqn:unram embedding}
\cO_F = \cO_{K_{\cyc,s}}\otimes_{W(k)} W(k_L)\hookrightarrow \cO_L.
\end{equation}
In particular,
we find that $\cO_L$ contains
$\cO_F$, 
or, equivalently, that  $L$  contains~$F$,
although the embedding~\eqref{eqn:unram embedding}
may not be the inclusion,
but rather the composite of this inclusion with an element
of $\Gal(F/K_{\cyc,s})$.
If we regard $\cO_{\overline{K}^{\rhobar}}$ (resp.\ $\cO_L$)
as a $\cO_F$-algebra via the inclusion (resp.\ via the embedding~\eqref{eqn:unram
embedding}),
then we see (from the very construction of~\eqref{eqn:unram embedding})
that $\eta$ is a morphism of~$\cO_F$-algebras.

Now,
let~$\pi$ be a uniformizer
  of~$\cO_{\overline{K}^{\ker \rho}}$, with
  image~$\eta(\pi)\in \cO_L/\mathfrak{a}_L^{>b/p^{s-1}}$; we claim that
  $\eta(\pi)$ is non-zero, so that its valuation is well-defined,
  and that this valuation is equal to that of~$\pi$. It follows
  immediately from this that the kernel of~$\eta$ is generated by
  $\mathfrak{a}_{\overline{K}^{\ker \rho}}^{>b/p^{s-1}}$, as required.
  To verify the claim, let~$E_\pi$ be the (Eisenstein) minimal polynomial
  of~$\pi$ over~$\cO_F$. 
 Then~$\eta(\pi)$ is a root
  of~$\eta(E_\pi)$, and the non-leading coefficients of~$\eta(E_\pi)$
  all have valuation at least~$v(T_s)=e/p^s$, with equality holding
  for the constant coefficient.
  Since~$e/p^s\le eh/p^s\le (p-1)b/p^s<pb/p^s$, it follows 
  that~$v(\eta(\pi))=v(\pi)$. 

We then have a
  composite of
  injections \[ T^\vee_{\overline{K}^{\ker\rho},pb,b}(\gM)\subseteq
    T^\vee_{L,pb,b}(\gM)\subseteq T^\vee_{\overline{K},pb,b}(\gM)\](the
  second inclusion being induced by the natural map
  $\cO_L\to \cO_{\overline{K}}$). By Corollary~\ref{cor: we get everything if and only if we contain the
    splitting field} the composite is an isomorphism 
(being an injection of isomorphic finite sets),
  so all of these inclusions are isomorphisms. Applying Corollary~\ref{cor: we get everything if and only if we contain the
    splitting field} again, we conclude that $\overline{K}^{\ker
    \rho}\subseteq L$, as required. 
\end{proof}
Recall that for each finite extension~$L$ of~$\Qp$ and each real number~$v$, we have the upper numbered
ramification subgroups $G_L^v$, as defined in~\cite[Chapter IV]{SerreLF}; these are the groups denoted~$G_L^{(v+1)}$
in~\cite{MR807070} (see~\cite[Rem.\ 1.2]{MR807070}). We can now prove our
first bound on the ramification of~$T(M)$, from which our subsequent
bounds will be deduced.
\begin{cor}
  \label{cor: ramification bound in the basic mod p case}Assume that we are in the setting of Hypothesis~{\em\ref{hyp: K
    is basic and M has height h and s is big enough}}. 
If~$v> pb-1$, then~$G_{K_{\cyc,s}}^v$ acts   trivially
on~$T(M)$.
\end{cor}
\begin{proof}
Let $\rho$ denote $T(M)^{\vee}$. 
Recall that for a real number~$m$, the extension $K^{\ker \rho}/K_{\cyc,s}$
satisfies $(\mathrm{P}_m)$ if, for any subfield $L$  of $\overline{K}$
containing $K_{\cyc,s}$,
the existence of an $\cO_{K_{\cyc,s}}$-algebra homomorphism
$\cO_{K^{\ker \rho}} \to \cO_L/\mf{a}^{\ge m/p^s}$ implies that
$K^{\ker \rho} \subseteq L$.
(This property was introduced by Fontaine,
and is discussed extensively in~\cite{MR2814776},
to which we refer the reader for further explanations
and references.  Note that the denominator
$p^s$ appears in $\mf{a}^{\geq m/p^s}$ because we always compute valuations
with respect to~$K$, and $K_{\cyc,s}$ is a totally ramified
extension of $K$ of degree~$p^s$.)

Corollary~\ref{cor: the slightly
    weird statement about lifting homomorphisms}
shows that $\inf\{m \, | \, (\mathrm{P}_m) \text{ holds for  }
K^{\ker \rho}/K_{\cyc,s} \} \leq pb.$
The present corollary then follows from this,
together with~\cite[Thm.\ 1.1]{MR2814776}.
\end{proof}

We now return to the setting of a general~$K$ (i.e.\ we do not assume
that~$K$ is basic). 
Let~$A$ be a finite local Artinian $\cO/\varpi^a$-algebra for
some~$a\ge 1$, and let~$M$ be a rank one projective \'etale
$\varphi$-module over~$\A_{K,A}$, which we assume to be of
$(\Kbasic,T)$-height at most~$h$ in the sense of Definition~\ref{defn:
  Kbasic T height at most h}; that is, we regard~$M$ as a
rank~$[K:\Kbasic]$ \'etale $\varphi$-module~$M'$
over~$\A_{\Kbasic,A}$, and we assume there is a projective
$\varphi$-module~$\gM$ over~$\A_{\Kbasic,A}^+$ of $T$-height at
most~$h$ such that $M'=\gM[1/T]$.

%
%

We have the $G_{\Kbasic_{\cyc}}$-representation $T_A(M')$, which we may
compute as \[T_A(M')=(W(\C^\flat)\otimes_{\A_{\Kbasic}}M')^{\varphi=1}=(W(\C^\flat)\otimes_{\A_{\Kbasic}}M)^{\varphi=1}.\]
It follows that there is a $G_{K_{\cyc}}$-equivariant surjection
\numequation\label{eqn: equivariant surjection from T(M') to T(M)}T_A(M')=(W(\C^\flat)\otimes_{\A_{\Kbasic}}M)^{\varphi=1}\onto
(W(\C^\flat)\otimes_{\A_K}M)^{\varphi=1}=T_A(M).\end{equation}
In fact~\eqref{eqn:induction via phi mods}
shows that
$T_A(M')$ can be identified with $\Ind_{G_{K_{\cyc}}}^{G_{\Kbasic_{\cyc}}} T_A(M)$,
and the morphism~\eqref{eqn: equivariant surjection from T(M') to T(M)} then 
becomes the natural $G_{K_{\cyc}}$-equivariant surjection.


\begin{prop}
  \label{prop: ramification bound in rank one case with coefficients
    and any K}Suppose that~
  $A$ is an Artinian $\cO/\varpi^a$-algebra for some~$a\ge 1$, and $M$
  is a rank 1 projective \'etale $\varphi$-module over~$\A_{K,A}$ of
  $(\Kbasic,T)$-height at most~$h$.
  
  Then 
  there is an open subgroup~$I^{a,h}$
  of~$I_{\Kcyc}^{\ab}$ 
  depending only on
  $K,a$ and~$h$ \emph{(}and not on~$A$ or~$M$\emph{)} such that the action of ~$I^{a,h}$ on
  $T_A(M)$ is trivial.
\end{prop}
\begin{proof}
%
We choose integers $b,s,N$ satisfying the following conditions:

  \begin{itemize}
  \item $b\ge eh/(p-1)$. 
  \item $N> pb/e$ is such 
    that $N\ge N(i,h)$ for each $i \in [1,a]$, 
where the $N(i,h)$ 
are as in Corollary~\ref{cor: Caruso
      Liu Galois action phi Gamma discrete version}, with~$K$ there
    being replaced by~$\Kbasic$.
  \item $s$ is such that~$p^{s-1}>b/e$,  and such that
$s\ge s(i,h,N)$ for each $i \in [1,a]$,
    where  the $s(i,h,N)$ 
are as in Corollary~\ref{cor: Caruso
      Liu Galois action phi Gamma discrete version}, with~$K$ there
    being replaced by~$\Kbasic$.
  \end{itemize}
Note in particular that Hypothesis~\ref{hyp: K is basic and M has height h and s
  is big enough} holds, with~$K$ in the hypothesis replaced
by~$\Kbasic$. 
Note also that our assumptions imply that the hypotheses of Corollary~\ref{cor: Caruso
      Liu Galois action phi Gamma discrete version} hold (with the field
$K$ there taken to be~$\Kbasic$), so that the 
    action of~$G_{\Kbasic_{\cyc}}$ on~$T_A(M')$ has a canonical
    extension to an action of $G_{\Kbasic_{\cyc,s}}$.
This canonical action on $T_A(M')$ induces 
the canonical action on
$\varpi^i T_A(M') = T_{A}(\varpi^i M')$ 
for each $i \in [0,a],$
and hence on each of the quotients
$\varpi^{i-1} T_A(M') /\varpi^i T_A(M') = T_{A/\varpi}(\varpi^{i-1} M'/\varpi^iM'),$
for $i \in [1,a]$.

To ease notation,
write~$M_i:=\varpi^{i-1}M/\varpi^iM$, for $1\le i\le a$,
and let $M_i'$ denote $M_i$ 
regarded as an \'etale
$\varphi$-module over~$\A_{\Kbasic,A/\varpi}$. By Corollary~\ref{cor:
  ramification bound in the basic mod p case}, 
if~$v$ is sufficiently large (depending only on~$K$, $a$ and~$h$, and
our subsequent choices of~$b$, $N$ and~$s$, and in particular not
on~$A$ or~$M$), then the canonical action of~$G_{\Kbasic_{\cyc,s}}^v$
on each~$T_{A/\varpi}(M'_i)$ is trivial. It follows that there is
some~$m$ depending only on~$a$ (and the degree $[K:\Kbasic]$) such that for
any~$g\in G_{\Kbasic_{\cyc,s}}^v$, the action of~$g^{p^m}$
on~$T_A(M')$ is trivial.

  It follows from the definitions that 
  $G_{\Kbasic_{\cyc,s}}^v\cap
  G_{K_{\cyc,s}}=G_{K_{\cyc,s}}^{\psi_{K_{\cyc,s}/\Kbasic_{\cyc,s}}(v)}$,
  where~$\psi_{K_{\cyc,s}/\Kbasic_{\cyc,s}}$ is as in~\cite[Chapter IV]{SerreLF}), so if
  we take any~$w\ge \psi_{K_{\cyc,s}/\Kbasic_{\cyc,s}}(v)$,
then~$(G_{K_{\cyc,s}}^w)^{p^m}$ 
($ := \{g^{p^m} \, | \,
g\in G_{K_{\cyc,s}}^w\}$; note that this is just a subset, not a subgroup,
of~$G_{K_{\cyc,s}}$)
acts trivially on $T_A(M')$. 

The product $G_{K_{\cyc}}
(G_{K_{\cyc,s}}^w)^{p^m}$ 
is {\em a priori} a subset of $G_K$, but we claim that it is in fact
an open subgroup of~$G_K$, and thus (since it contains
$G_{K_{\cyc}}$) is equal to $G_{K_{\cyc,s'}}$  for
some $s' \geq s$ which depends only on~$s$ and $m$
(and so ultimately only on  $K$, $a$ and $h$).

To prove the claim, 
we briefly recall the relationship between
upper numbered ramification groups and the filtration of the unit
group: If~$L$ is a finite extension of~$\Qp$, then the Artin
map identifies~$\cO_L^\times$ with~$I_L^{\ab}$, and for each integer~$w\ge 1,$
identifies the subgroup~$1+(\m_L)^v$ of~$\cO_L^\times$
with the image of~$G_L^v$ in~$G_L^{\ab}$ (see~\cite[Chapter XV, \S 2,
Thm.\ 2]{SerreLF}).   It then identifies the subgroup $(1+(\m_L)^v)^{p^m}$
with the image of~$(G_L^v)^{p^m}$ in~$G_L^{\ab}$.  This former
subgroup is an open subgroup
of $\cO_L^{\times}$,
and so we find that for each~$v$, the image of $(G_L^v)^{p^m}$
in~$I_L^{\ab}$ is open.

Now the extension $K_{\cyc} / K_{\cyc,s}$ is an infinitely ramified
$\Gamma_{K_{\cyc,s}} \cong \Z_p$-extension,
and so the inertia group $I(K_{\cyc}/K)$ is finite index (equivalently,
open) in~$\Gamma_{K_{\cyc,s}}$.   Since $I(K_{\cyc}/K)$ is a quotient
of~$I_K^{\ab}$, the discussion of the preceding paragraph then shows that 
the image of
$(G_{K_{\cyc,s}}^w)^{p^m}$ 
in $\Gamma_{K_{\cyc,s}}$  is also an open subgroup. 
Since $\Gamma_{K_{\cyc,s}}$ is itself open in~$G_K$,
this implies the claim.

Since~\eqref{eqn: equivariant surjection from T(M') to T(M)}
is $G_{K_{\cyc}}$-equivariant, and since 
$(G_{K_{\cyc,s}}^w)^{p^m}$ 
acts trivially on~$T_A(M')$,
we see that the kernel 
of~\eqref{eqn: equivariant surjection from T(M') to T(M)}
is in fact~$G_{K_{\cyc,s'}}$-invariant, and thus that
the $G_{K_{\cyc}}$-action on $T_A(M)$ extends to a $G_{K_{\cyc,s'}}$-action,
which is abelian (since the $G_{K_{\cyc}}$-action
is abelian, as $T_A(M)$ has rank one over~$A$, while
$(G_{K_{\cyc,s}}^w)^{p^m}$ 
acts trivially). 

Now another application of~\cite[Chapter IV]{SerreLF}) shows that
$G_{K_{\cyc,s}}^w \cap \, G_{K_{\cyc,s'}} =
G_{K_{\cyc,s'}}^{w'}$ 
for some  $w'$ depending only on~$w$, and our preceding discussion
of abelian ramification theory shows that
$(G_{K_{\cyc,s'}}^{w'})^{p^m}$ 
has open image in $G_{K_{\cyc,s'}}^{\ab}$.
If we let $I'$  denote the preimage of this open image
(i.e.\  \[I' = \overline{[G_{K_{\cyc,s'}},G_{K_{\cyc,s'}}]} 
(G_{K_{\cyc,s'}}^{w'})^{p^m},\] where the overline
denotes closure in~$G_{K_{\cyc,s'}}$),
then $I'$ is an open subgroup  of $G_{K_{\cyc,s'}}$ which
acts trivially on~$T_A(M)$.
Thus $I'' := I' \cap G_{K_{\cyc}}$ also acts trivially on~$T_A(M)$.
Finally, we let $I^{a,h}$ denote the image of $I''$ in $I_K^{\ab}$.
(As the notation indicates, $I^{a,h}$  depends only on $a$,  $h$, and~$K$,
but not on~$A$ or~$M$, since this is true of each of the groups
$I'$ and $I''$ from which it is constructed.)
\end{proof}

The same argument allows us to prove the following variant of Proposition~\ref{prop: ramification bound in rank one case with coefficients
  and any K}.

\begin{prop}
  \label{prop: ramification bound in rank one case with coefficients
    and any K with GK action}Suppose that
  $A$ is an Artinian $\cO/\varpi^a$-algebra for some~$a\ge 1$, and $M$
  is a rank 1 projective \'etale $(\varphi,\Gamma)$-module
  over~$\A_{K,A}$. Assume furthermore that~
  $M$ has $(\Kbasic,T)$-height at most~$h$, and that ~$s$ is some
  sufficiently large integer with the property that
the action of~$\Gamma_{\Kbasic,\disc}=\Gamma_{K,\disc}$ on $M$ extends the canonical action of
    $\Gamma_{\Kbasic_s,\disc}$ given by Corollary~\emph{\ref{cor: Caruso Liu Galois
      action phi Gamma discrete version}}. {\em (}In particular, we are
  assuming that~$s$ has been chosen to be large enough that this
  canonical action is defined.{\em )}
  
  Then there is an open subgroup~$I_K^{a,h,s}$ of~$I_{K}^{\ab}$
  depending only on $K,a,h$ and~$s$ \emph{(}and not on~$A$
  or~$M$\emph{)} such that the action of ~$I_K^{a,h,s}$ on $T_A(M)$ is
  trivial.
\end{prop}
\begin{proof}
The proof is very similar to that
of Proposition~\ref{prop: ramification bound in rank one
case with coefficients and any K},
and we content ourselves with indicating the key modifications.
The main difference between the setting of that proposition and
the present one is that, in our present setting,
the morphism~\eqref{eqn: equivariant surjection from T(M') to T(M)}
is a $G_K$-equivariant morphism from a $G_{\Kbasic}$-representation
to a $G_K$-representation.  Furthermore, by assumption,
the $G_{\Kbasic_{\cyc,s}}$-action on its source is the canonical action.  

We now follow the argument 
of Proposition~\ref{prop: ramification bound in rank one
case with coefficients and any K}, and find that,
after possibly enlarging $s$ in a manner depending only on $a$, $h$,
and~$K$, some ramification group $G_{K_{\cyc,s}}^w$ acts
trivially on~$T_A(M')$ (the index $w$ also depending
only on $a$, $h$, and~$s$).    
In the present argument, there is no need to pass
to an auxiliary subgroup $G_{K_{\cyc,s'}}$, since
the morphism~\eqref{eqn: equivariant surjection from T(M') to T(M)}
is in particular already a $G_{K_{\cyc,s}}$-equivariant morphism
of $G_{K_{\cyc,s}}$-representations;
and the action of $G_{K_{\cyc,s}}$ on $T(M)$ is already abelian.
The argument then proceeds in the same manner,
to produce an open subgroup $I'$
of the inertia subgroup of $G_{K_{\cyc,s}}$ 
(and thus of $I_K$) which acts trivially on~$T_A(M)$.
(And, by construction, the group $I'$ depends only on $a$, $h$,
and~$s$, but not on~$A$ or~$M$.)
We then let
$I_K^{a,h,s}$ 
denote the image of $I'$ in~$I_K^{\ab}$.
%
%
  \end{proof}

\chapter{A geometric Breuil--M\'ezard conjecture}\label{sec: BM}
 In this chapter
we explain a (for the most part conjectural) relationship between the
geometry of our potentially semistable and crystalline moduli
stacks~$\cX_d^{\crys,\lambdau,\tau}$
and~$\cX_d^{\semis,\lambdau,\tau}$, and the representation theory
of~$\GL_n(k)$. Throughout the chapter we fix a sufficiently large
coefficient field~$E$ with ring of integers~$\cO$ and residue
field~$\F$, and we largely omit it from our notation.

The starting point for this proposed relationship is
the Breuil--M\'ezard conjecture~\cite{BreuilMezard},
which is
a conjectural formula for the Hilbert--Samuel multiplicities of the
special fibres of the lifting
rings~$R^{\crys,\underline{\lambda},\tau}_{\rhobar}$ and
~$R^{\semis,\underline{\lambda},\tau}_{\rhobar}$, and which has important
applications to proving modularity lifting theorems via the
Taylor--Wiles method~\cite{KisinFM}. The conjecture was geometrized
in~\cite{MR3248725} and~\cite{emertongeerefinedBM}, by refining the
conjectural formula for the Hilbert--Samuel multiplicity to a
conjectural formula for the underlying cycle of the special fibres of
the lifting rings ~$R^{\crys,\underline{\lambda},\tau}_{\rhobar}$ and
~$R^{\semis,\underline{\lambda},\tau}_{\rhobar}$, considered as (equidimensional)
closed subschemes of the special fibre of the universal lifting
ring~$R_{\rhobar}^{\square}$. We refer to this generalization as the
``refined Breuil--M\'ezard conjecture'', and to the original
conjecture (or rather, its generalizations to~$\GL_n$) as the
``numerical Breuil--M\'ezard conjecture''.

Our aim in this chapter, then, is to ``globalize'' the conjectures
of~\cite{emertongeerefinedBM} by formulating versions of them
for the 
stacks~$\cX_d^{\crys,\lambdau,\tau}$
and~$\cX_d^{\crys,\lambdau,\tau}$; the conjectures
of~\cite{emertongeerefinedBM} can be recovered from these conjectures
by passing to versal rings at finite type points. (We caution the
reader that there is another kind of ``globalization'' that could be
considered, namely realizing local Galois representations as the
restrictions to decomposition groups of global representations, as
used in the proofs of some of the results
of~\cite{emertongeerefinedBM} and~\cite{geekisin}.
Other than in the proof of
Theorem~\ref{thm: Gee Kisin extended to semistable} below,
we don't  consider this 
kind of globalization in the present work.) This
generalization, which we will call the ``geometric Breuil--M\'ezard
conjecture'', seems to us to be the natural setting in which to
consider the Breuil--M\'ezard conjecture, and we will use our
description of the irreducible components of~$\cX_{d,\red}$ to deduce
new results about both the refined and numerical versions 
of the Breuil--M\'ezard conjecture.

\section[The qualitative geometric Breuil--M\'ezard
  conjecture]{The qualitative geometric Breuil--M\'ezard
  conjecture\sectionmark{The qualitative geometric BM
  conjecture}}\sectionmark{The qualitative geometric BM
  conjecture}\label{subsec: qualitative BM}If~$\lambdau$ is a regular
Hodge type, and~$\tau$ is any inertial type, then by Theorems~\ref{thm: existence of ss stack} and~\ref{thm: dimension of ss stack}, the
stacks~$\cX^{\crys,\lambdau,\tau}_d$ and~$\cX^{\semis,\lambdau,\tau}_d$
are finite type $p$-adic formal algebraic stacks over~$\cO$, which are
$\cO$-flat and equidimensional of dimension~$1+[K:\Qp]d(d-1)/2$. It
follows that their special fibres $\cXbar^{\crys,\lambdau,\tau}_d$ and~$\cXbar^{\semis,\lambdau,\tau}_d$ are algebraic stacks over~$\F$ which
are equidimensional of dimension~$[K:\Qp]d(d-1)/2$. Since
~$\cX^{\crys,\lambdau,\tau}_d$ and~$\cX^{\semis,\lambdau,\tau}_d$ are
closed substacks of~$\cX_d$, $\cXbar^{\crys,\lambdau,\tau}_d$
and~$\cXbar^{\semis,\lambdau,\tau}_d$  are closed substacks of the
special fibre~$\cXbar_d$, and their irreducible components (with the
induced reduced substack structure) are therefore closed substacks of
the algebraic stack~$\cXbar_{d,\red}$ (see~\cite[\href{https://stacks.math.columbia.edu/tag/0DR4}{Tag 0DR4}]{stacks-project}
for the theory of irreducible components of algebraic stacks and their
multiplicities). 
Since~$\cXbar_{d,\red}$ is
equidimensional of dimension~$[K:\Qp]d(d-1)/2$ by
Theorem~\ref{thm:reduced dimension}, it follows that the irreducible
components of  $\cXbar^{\crys,\lambda,\tau}_d$
and~$\cXbar^{\semis,\lambda,\tau}_d$ are irreducible components
of~$\cXbar_{d,\red}$, and are therefore of the
form~$\cXbar_{d,\red}^{\underline{k}}$ for some Serre
weight~$\underline{k}$ (again by Theorem~\ref{thm:reduced
  dimension}).

For each~$\underline{k}$, we
write~$\mu_{\underline{k}}(\cXbar^{\crys,\lambdau,\tau}_d)$ and
$\mu_{\underline{k}}(\cXbar^{\semis,\lambdau,\tau}_d)$ for the
multiplicity
of~$\cXbar_{d,\red}^{\underline{k}}$ as a component
of~$\cXbar^{\crys,\lambdau,\tau}_d$
and~$\cXbar^{\semis,\lambdau,\tau}_d$. 
We
write~$Z_{\crys,\lambdau,\tau}=Z(\cXbar^{\crys,\lambdau,\tau}_d)$ and
~$Z_{\semis,\lambdau,\tau}=Z(\cXbar^{\semis,\lambdau,\tau}_d)$ for the
corresponding cycles, i.e.\ for the formal sums
\numequation\label{eqn: cris HS multiplicity stack}Z_{\crys,\lambdau,\tau}=\sum_{\underline{k}}\mu_{\underline{k}}(\cXbar^{\crys,\lambdau,\tau}_d)\cdot\cXbar_d^{\underline{k}}, \end{equation}
\numequation\label{eqn: ss HS multiplicity stack}Z_{\semis,\lambdau,\tau}=\sum_{\underline{k}}\mu_{\underline{k}}(\cXbar^{\semis,\lambdau,\tau}_d)\cdot\cXbar_d^{\underline{k}}, \end{equation}which
we regard as elements of the finitely generated free abelian group~$\Z[\cX_{d,\red}]$
whose generators are the irreducible
components~$\cXbar_d^{\underline{k}}$.

Now fix some representation~$\rhobar:G_K\to\GL_d(\F)$, corresponding to a
point $x:\Spec\F\to\cX_d$. (More generally, we could consider
representations valued in~$\GL_d(\F')$ for some finite
extension~$\F'/\F$, but in keeping with our attempt to keep the
notation in this chapter as uncluttered as possible by omitting~$\cO$,
it is convenient to suppose that~$\F$ has been chosen sufficiently
large.)  For each regular Hodge type~$\lambdau$ and inertial
type~$\tau$, we have effective versal morphisms
$\Spec
R_{\rhobar}^{\crys,\underline{\lambda},\tau}/\varpi\to\cXbar^{\crys,\lambdau,\tau}$
and
$\Spec
R_{\rhobar}^{\semis,\underline{\lambda},\tau}/\varpi\to\cXbar^{\semis,\lambdau,\tau}$
(see Corollary~\ref{cor: crystalline deformation rings are effective
  versal}),
as well as a (non-effective; see Remark~\ref{rem:non-effective}) versal morphism $\Spf R_{\rhobar}^{\square} \to \cX_d$
(see Proposition~\ref{prop:versal rings}).  


For each~$\underline{k}$,
we may consider the fibre product
$\Spf R_{\rhobar}^\square\times_{\cX_d}\cX_d^{\underline{k}}.$
This is {\em a priori} a closed formal subscheme of $\Spf R_{\rhobar}^{\square}$,
but since $R_{\rhobar}^{\square}$ is a complete local ring, it may
equally well be regarded as a closed subscheme 
of $\Spec R_{\rhobar}^{\square}$
(see Lemma~\ref{lem: alem closed formal algebraic scheme adic*}).

\begin{lemma}
\label{lem:computing cycle dimension}
The fibre product
$\Spf R_{\rhobar}^\square\times_{\cX_d}\cX_d^{\underline{k}}$,
when we
regard it as a closed subscheme of $\Spec R_{\rhobar}^{\square}$,
is equidimensional of dimension~$d^2+[K:\Qp]d(d-1)/2$.
\end{lemma}
\begin{proof}
As in the proof of Theorem~\ref{thm:reduced dimension},
we may find a regular Hodge type $\underline{\lambda}$ such that $\cX_d^{\underline{k}}$
is an irreducible component  of $(\overline{\cX}^{\crys, \underline{\lambda}})_{\red}.$
Thus
$\Spf R_{\rhobar}^\square\times_{\cX_d}\cX_d^{\underline{k}}$
is versal to a union of irreducible components of
the spectrum of the versal ring
$(R_{\rhobar}^{\crys,\underline{\lambda}}/\varpi)_{\red}$
to $(\overline{\cX}^{\crys,\underline{\lambda}})_{\red}$
(those that correspond to the various formal branches of $\cX_d^{\underline{k}}$  
passing through~$x$, in the terminology of \cite[\href{https://stacks.math.columbia.edu/tag/0DRA}{Tag
  0DRA}]{stacks-project})
and hence is of the stated dimension, since 
$R_{\rhobar}^{\crys,\underline{\lambda}}/\varpi$,
and so also
$(R_{\rhobar}^{\crys,\underline{\lambda}}/\varpi)_{\red}$,
is equidimensional of this dimension.
\end{proof}

We let $\cC_{\underline{k}}(\rhobar)$ denote the 
$d^2+[K:\Qp]d(d-1)/2$-dimensional
cycle in $\Spec R_{\rhobar}^{\square}/\varpi$ underlying the fibre product
of Lemma~\ref{lem:computing cycle dimension}.

The following theorem gives a qualitative
version of the refined Breuil--M\'ezard
conjecture~\cite[Conj.\ 4.2.1]{emertongeerefinedBM}. While its statement is purely local, we do not know
how to prove it without making use of the stack~$\cX_d$.
\begin{thm}
  \label{thm: qualitative BM}Let $\rhobar:G_K\to\GL_d(\F)$ be a
  continuous representation. Then there are finitely many cycles of
  dimension~$d^2+[K:\Qp]d(d-1)/2$ in~$\Spec R_{\rhobar}^\square/\varpi$ such that
  for any regular Hodge type~$\lambdau$ and any inertial type~$\tau$,
  each of the special fibres~$\Spec
  R_{\rhobar}^{\crys,\underline{\lambda},\tau}/\varpi$ and~$\Spec
  R_{\rhobar}^{\semis,\underline{\lambda},\tau}/\varpi$ is set-theoretically
  supported on some union of these cycles.
\end{thm}
\begin{proof}
  We have
  $\Spf R_{\rhobar}^{\crys,\underline{\lambda},\tau}/\varpi =\Spf
  R^{\square}_{\rhobar}\times_{\cX_d}\cXbar^{\crys,\lambdau,\tau}$ and
  $\Spf R_{\rhobar}^{\semis,\underline{\lambda},\tau}/\varpi =\Spf
  R^{\square}_{\rhobar}\times_{\cX_d}\cXbar^{\crys,\lambdau,\tau}$. It
  follows from~\eqref{eqn: cris HS multiplicity stack} and~\eqref{eqn:
    ss HS multiplicity stack}, together with the definition
  of~$\cC_{\underline{k}}(\rhobar)$, that we may write the underlying
  cycles as \numequation\label{eqn: cris HS multiplicity def
    ring}Z(\Spec
  R_{\rhobar}^{\crys,\underline{\lambda},\tau}/\varpi)=\sum_{\underline{k}}\mu_{\underline{k}}(\cXbar^{\crys,\lambdau,\tau}_d)\cdot\cC_{\underline{k}}(\rhobar), \end{equation}
\numequation\label{eqn: ss HS multiplicity def ring}Z(\Spec
R_{\rhobar}^{\semis,\underline{\lambda},\tau}/\varpi)=\sum_{\underline{k}}\mu_{\underline{k}}(\cXbar^{\semis,\lambdau,\tau}_d)\cdot\cC_{\underline{k}}(\rhobar). \end{equation}(Note
that by~\cite[\href{https://stacks.math.columbia.edu/tag/0DRD}{Tag
  0DRD}]{stacks-project}, the multiplicities do not change
when passing to versal rings.) The
theorem follows immediately (taking our finite set of cycles to be
the~$\cC_{\underline{k}}(\rhobar)$). 
\end{proof}
We can regard this theorem as isolating the ``refined'' part
of~\cite[Conj.\ 4.2.1]{emertongeerefinedBM}; that is, we have taken
the original numerical Breuil--M\'ezard conjecture, formulated a geometric
refinement of it, and then removed the numerical part of the
conjecture. The numerical part of the conjecture (in the optic of this
chapter) consists of relating the multiplicities
$\mu_{\underline{k}}(\cXbar^{\crys,\lambdau,\tau}_d)$ and
$\mu_{\underline{k}}(\cXbar^{\semis,\lambdau,\tau}_d)$ to the 
representation theory of~$\GL_n(k)$, as we recall in the next section.

\section{Semistable and crystalline inertial types}\label{subsec:
  type theory} We now
briefly recall the ``inertial local Langlands correspondence''
for~$\GL_d$. Let~$\rec_p$ denote the local Langlands correspondence
for~$\Qpbar$-representations of $\GL_d(K)$, normalized as
in~\cite[\S1.8]{Gpatch}; this is a bijection between the isomorphism
classes of irreducible smooth $\Qpbar$-representations of~$\GL_d(K)$
and the isomorphism classes of $d$-dimensional semisimple
Weil--Deligne $\Qpbar$-representations of the Weil group~$W_K$.  We
have the following result, which is essentially due to
Schneider--Zink \cite{MR1728541}.

\begin{thm}Let ~$\tau:I_K\to\GL_d(\Qpbar)$ be an inertial type. Then
  there are finite-dimensional smooth irreducible
  $\Qpbar$-representations $\sigma^{\crys}(\tau)$ and
  $\sigma^{\semis}(\tau)$ of~$\GL_d(\cO_K)$ with the properties that
  if ~$\pi$ is an irreducible smooth $\Qpbar$-representation
  of~$\GL_d(K)$, then the $\Qpbar$-vector space
  $\Hom_{\GL_d(\cO_K)}(\sigma^{\crys}(\tau),\pi)$ \emph{(}resp.\
the $\Qpbar$-vector space  $\Hom_{\GL_d(\cO_K)}(\sigma^{\semis}(\tau),\pi)$\emph{)} has dimension at
  most~$1$, and is nonzero precisely if~$\rec_p(\pi)|_{I_F}\cong
  \tau$, and $N=0$ on~$\rec_p(\pi)$ \emph{(}resp.\, if~$\rec_p(\pi)|_{I_F}\cong
  \tau$, and~$\pi$ is generic\emph{)}.
\end{thm}
\begin{proof}
  See~\cite[Thm.\ 3.7]{Gpatch} for~$\sigma^{\crys}(\tau)$, and~\cite[Thm.\
  3.7]{MR3769675} together with~\cite[Thm.\ 2.1, Lem.\ 2.2]{2018arXiv180302693P}  for~$\sigma^{\semis}(\tau)$.
\end{proof}Note that we do not claim that the
representations~$\sigma^{\crys}(\tau)$ and~$\sigma^{\semis}(\tau)$ are
unique; the possible non-uniqueness of these representations is
of no importance for us.

For each regular Hodge
type~$\underline{\lambda}$ we let~$L_{\lambdau}$ be the corresponding
representation of~$\GL_d(\cO_K)$, defined as follows: For
each~$\sigma:K\into\Qpbar$, we
write~$\xi_{\sigma,i}=\lambda_{\sigma,i}-(d-i)$, so that
$\xi_{\sigma,1}\ge\dots\ge\xi_{\sigma,d}$. We view each
$\xi_\sigma:=(\xi_{\sigma,1},\dots,\xi_{\sigma,d})$ as a dominant
weight of the algebraic group~$\GL_d$ (with respect to the upper
triangular Borel subgroup), and we write~$M_{\xi_\sigma}$ for the
algebraic $\cO_K$-representation of~$\GL_d(\cO_K)$ of highest
weight~$\xi_\sigma$. Then we define
$L_{\lambdau}:=\otimes_{\sigma}M_{\xi_\sigma}\otimes_{\cO_K,\sigma}\cO$.

For each~$\tau$ we let~$\sigma^{\crys,\circ}(\tau)$,
$\sigma^{\semis,\circ}(\tau)$ denote choices of $\GL_d(\cO_K)$-stable
$\cO$-lattices in~$\sigma^{\crys}(\tau)$, $\sigma^{\semis}(\tau)$
respectively (the precise choices being unimportant). Then we write
$\sigmabar^{\crys}(\lambda,\tau)$, (resp.\ $\sigmabar^{\semis}(\lambda,\tau)$)
for the semisimplification of the $\F$-representation of~$\GL_d(k)$
given by
$L_{\lambdau}\otimes_\cO\sigma^{\crys,\circ}(\tau)\otimes_\cO\F$
(resp.\
$L_{\lambdau}\otimes_\cO\sigma^{\semis,\circ}(\tau)\otimes_\cO\F$). For
each Serre weight~$\underline{k}$, we write~$F_{\underline{k}}$ for
the corresponding irreducible $\F$-representation of~$\GL_d(k)$ (see
for example the appendix to~\cite{MR2541127}). Then there are unique
integers $n_{\underline{k}}^\crys(\lambda,\tau)$ and
  $n_{\underline{k}}^\semis(\lambda,\tau)$ such
    that \[\sigmabar^{\crys}(\lambda,\tau)\cong\oplus_{\underline{k}}F_{\underline{k}}^{\oplus
        n_{\underline{k}}^\semis(\lambda,\tau)},\] \[\sigmabar^{\semis}(\lambda,\tau)\cong\oplus_{\underline{k}}F_{\underline{k}}^{\oplus n_{\underline{k}}^\semis(\lambda,\tau)}.\]
Our geometric Breuil--M\'ezard conjecture is as follows. \index{Breuil--M\'ezard conjecture}
\begin{conj}
  \label{conj: geometric BM}There are cycles~$Z_{\underline{k}}$ with
  the property that for each regular Hodge type~$\lambdau$ and each
  inertial type~$\tau$, we have
  $Z_{\crys,\lambdau,\tau}=\sum_{\underline{k}}n_{\underline{k}}^\crys(\lambda,\tau)\cdot
  Z_{\underline{k}}$,  $Z_{\semis,\lambdau,\tau}=\sum_{\underline{k}}n_{\underline{k}}^\semis(\lambda,\tau)\cdot
  Z_{\underline{k}}$.
\end{conj}

For some motivation for the conjecture (coming from the
Taylor--Wiles patching method), see for example~\cite[Thm.\
5.5.2]{emertongeerefinedBM}. Some evidence for the conjecture is given
in the following sections. 

The expressions~\eqref{eqn: cris HS multiplicity stack} and~\eqref{eqn: ss HS multiplicity
stack} describe $Z_{\crys,\lambdau,\tau}$ and $Z_{\semis,\lambdau,\tau}$
as linear combinations of (the cycles underlying) the irreducible
components~$\overline{\cX}_d^{\underline{k}'}$, and thus each cycle $Z_{\underline{k}}$
(assuming that such cycles exist) will itself be a linear combination of the
various~$\overline{\cX}_d^{\underline{k}'}.$ 
We expect that the
cycles~$Z_{\underline{k}}$ are effective, i.e.\ that each of them is a
linear
combination of the~$\cXbar_d^{\underline{k}'}$ with non-negative
(integer) coefficients. Note that the finitely many (conjectural)
cycles~$Z_{\underline{k}}$ are completely determined by the infinitely
many equations in Conjecture~\ref{conj: geometric BM} (see~\cite[Lem.\
4.1.1, Rem.\ 4.1.7(1)]{emertongeerefinedBM}).

While the original motivation for Conjecture~\ref{conj:
  geometric BM} was to understand potentially semistable deformation
rings in terms of the representation theory of~$\GL_d$, it can also be
thought of as giving a geometric interpretation of the
multiplicities~$n_{\underline{k}}^\crys(\lambda,\tau)$.

\section[{The relationship between the numerical, refined and
  geometric Breuil--M\'ezard conjectures}]{The relationship between the numerical, refined and
  geometric Breuil--M\'ezard conjectures \sectionmark{The relationship between the BM conjectures}}\sectionmark{The relationship between the BM conjectures}\label{subsec: relating
  different BM conjectures}We now explain the relationship between
Conjecture~\ref{conj: geometric BM} and the conjectures
of~\cite{emertongeerefinedBM}. Suppose firstly that
Conjecture~\ref{conj: geometric BM} holds, and fix some
$\rhobar:G_K\to\GL_d(\F)$ corresponding to a point $x:\Spec\F\to\cX_d$. For
each~$\underline{k}$, we
set \[Z_{\underline{k}}(\rhobar):=\Spf
  R_{\rhobar}^\square\times_{\cX_d}Z_{\underline{k}},\]which
we regard as a
$d^2+[K:\Qp]d(d-1)/2$-dimensional cycle
in $\Spec R_{\rhobar}^\square/\varpi$.
(As noted above, the $Z_{\underline{k}}$ will be linear combinations 
of the cycles underlying the various $\cXbar_d^{\underline{k}'}$,
and so the the $Z_{\underline{k}}(\rhobar)$ will be linear
combinations of the various cycles $\cC_{\underline{k}'}(\rhobar)$;
thus Lemma~\ref{lem:computing cycle
dimension} shows that they are indeed 
$d^2+[K:\Qp]d(d-1)/2$-dimensional cycles.)

\begin{remark}
\label{rem:cycle alternative}
Equivalently, we can define
$Z_{\underline{k}}(\rhobar)$ to be the image of~$Z_{\underline{k}}$
under the natural map from $\Z[\cX_{d,\red}]$ to the group $\Z_{d^2+[K:\Qp]d(d-1)/2}(\Spec
R_{\rhobar}^\square/\varpi)$ of $d^2+[K:\Qp]d(d-1)/2$-dimensional cycles in $\Spec
R_{\rhobar}^\square/\varpi$ which is defined as follows: We
let~$R_{\rhobar}^{\alg}$ be the quotient of~$R_{\rhobar}^{\square}/\varpi$
which is a versal ring to~$\cX_{d,\red}$ at~$x$, so that 
by~\cite[\href{https://stacks.math.columbia.edu/tag/0DRB}{Tag 0DRB},\href{https://stacks.math.columbia.edu/tag/0DRD}{Tag 0DRD}]{stacks-project} we have a 
multiplicity-preserving surjection from the set of irreducible
components of~$\Spec R_{\rhobar}^{\alg}$ to the set of irreducible
components of~$\cX_{d,\red}$ containing~$x$; we then send any
irreducible component of~$\cX_{d,\red}$ not containing~$x$ to zero,
and send each irreducible component containing~$x$ to the sum of the corresponding irreducible
components of~$\Spec R_{\rhobar}^{\alg}$  in its preimage.
\end{remark}

Exactly as in the proof of
Theorem~\ref{thm: qualitative BM}, it follows that for each regular
type~$\underline{\lambda}$ and inertial type~$\tau$, we have
\numequation\label{eqn: refined version of BM cris}Z(\Spec
  R_{\rhobar}^{\crys,\underline{\lambda},\tau}/\varpi)=\sum_{\underline{k}}n_{\underline{k}}^\crys(\lambda,\tau)\cdot
  Z_{\underline{k}}(\rhobar),\end{equation} \numequation\label{eqn: refined version of BM ss}Z(\Spec
  R_{\rhobar}^{\semis,\underline{\lambda},\tau}/\varpi)=\sum_{\underline{k}}n_{\underline{k}}^\semis(\lambda,\tau)\cdot
  Z_{\underline{k}}(\rhobar).\end{equation} The first of these statements is
\cite[Conj.\ 4.2.1]{emertongeerefinedBM} (with the cycles~$\cC_a$
there being our cycles~$Z_{\underline{k}}(\rhobar)$), and the second
is the corresponding statement for potentially semistable lifting
rings.

We now relate this conjecture to the numerical version of the
Breuil--M\'ezard conjecture. We have a homomorphism $\Z_{d^2+[K:\Qp]d(d-1)/2}(\Spec
R_{\rhobar}^\square/\varpi)\to\Z$ defined by sending each cycle to its
Hilbert--Samuel multiplicity in the sense
of~\cite[\S2.1]{emertongeerefinedBM}. 
Let $\mu_{\underline{k}}(\rhobar)$ denote the Hilbert--Samuel
multiplicity of the cycle $Z_{\underline{k}}(\rhobar)$, and write $e(\Spec
  R_{\rhobar}^{\crys,\underline{\lambda},\tau}/\varpi)$, $e(\Spec
  R_{\rhobar}^{\semis,\underline{\lambda},\tau}/\varpi)$ for the Hilbert--Samuel
  multiplicities of the indicated rings. Then it follows
  from~\eqref{eqn: refined version of BM cris} and~\eqref{eqn: refined
    version of BM ss} that we have
  \numequation\label{eqn: numerical version of BM cris}e(\Spec
  R_{\rhobar}^{\crys,\underline{\lambda},\tau}/\varpi)=\sum_{\underline{k}}n_{\underline{k}}^\crys(\lambda,\tau)\mu_{\underline{k}}(\rhobar),\end{equation} \numequation\label{eqn:
  numerical version of BM ss}e(\Spec
  R_{\rhobar}^{\semis,\underline{\lambda},\tau}/\varpi)=\sum_{\underline{k}}n_{\underline{k}}^\semis(\lambda,\tau)\mu_{\underline{k}}(\rhobar).\end{equation}
Then~\eqref{eqn: numerical version of BM cris}
is~\cite[Conj.\ 4.1.6]{emertongeerefinedBM}, and ~\eqref{eqn:
  numerical version of BM ss} is the corresponding semistable version.

Suppose now that for each $\rhobar:G_K\to\GL_d(\F)$, all regular Hodge
types~$\lambdau$ and all inertial types~$\tau$ we have~\eqref{eqn:
  numerical version of BM cris} and~\eqref{eqn: numerical version of
  BM ss} for some integers~$\mu_{\underline{k}}(\rhobar)$, but do not
assume Conjecture~\ref{conj: geometric BM} (so in particular we do not
presuppose any geometric interpretation for the
integers~$\mu_{\underline{k}}(\rhobar)$). 

For each $\underline{k}$ we choose a point
$x_{\underline{k}}:\Spec\F\to\cXbar_{d,\red}$ which is contained
in~$\cXbar^{\underline{k}}$ and not in any~$\cXbar^{\underline{k}'}$
for $\underline{k}'\ne\underline{k}$. We furthermore demand that
$x_{\underline{k}}$ is a smooth point of~$\cXbar_{d,\red}$. (Since
$\cXbar_{d,\red}$ is reduced and of finite type over~$\F$, there is a dense set
of points of~$\cXbar^{\underline{k}}$ satisfying these conditions.)
Write~$\rhobar_{\underline{k}}:G_K\to\GL_d(\F)$ for the
representation corresponding to~$X_{\underline{k}}$, and
set \numequation\label{eqn: cycle from multiplicities}Z_{\underline{k}}:=\sum_{\underline{k}'}
  \mu_{\underline{k}}(\rhobar_{\underline{k}'})\cdot\cXbar^{\underline{k}'}.\end{equation}
Then for each
regular Hodge type~$\lambdau$ and inertial type~$\tau$, it follows
from~\eqref{eqn: numerical version of BM cris}
that \begin{align*}\sum_{\underline{k}}n_{\underline{k}}^\crys(\lambda,\tau)\cdot
  Z_{\underline{k}}&=
                     \sum_{\underline{k}}e(R_{\rhobar_{\underline{k}}}^{\crys,\underline{\lambda},\tau}/\varpi)
                     \cdot\cXbar^{\underline{k}}\\
                   &=\sum_{\underline{k}}\mu_{\underline{k}}(\cXbar^{\crys,\lambdau,\tau}_d)\cdot\cXbar^{\underline{k}}\\
                   &=Z_{\crys,\lambdau,\tau}, \end{align*}where we
                 used that $x_{\underline{k}}$ is a smooth point
                 of~$\cXbar_{d,\red}$ and is only contained  in ~$\cXbar^{\underline{k}}$
to conclude that
 $\mu_{\underline{k}}(\cXbar^{\crys,\lambdau,\tau}_d)=e(R_{\rhobar_{\underline{k}}}^{\crys,\underline{\lambda},\tau}/\varpi)$. 
Similarly we have  $Z_{\semis,\lambdau,\tau}=\sum_{\underline{k}}n_{\underline{k}}^\semis(\lambda,\tau)\cdot
  Z_{\underline{k}}$, and we conclude that 
the
  geometric Breuil--M\'ezard conjecture (Conjecture~\ref{conj:
    geometric BM}) is equivalent to the numerical conjecture.
  \begin{rem}\label{rem: geom BM lets us go from sufficiently generic
      rhobar to all rhobar}
    The argument that we just made shows that the geometric conjecture
    follows from knowing the numerical conjecture for sufficiently
    generic~$\rhobar$ (indeed, it is enough to check it for a single
    sufficiently generic~$\rhobar$ for each irreducible component
    of~$\cX_{d,\red}$), while in turn the geometric conjecture implies
    the numerical conjecture for all~$\rhobar$. 
  \end{rem}

\section{The weight part of Serre's conjecture}\label{subsec: Serre
  weights}We now briefly explain the relationship between
Conjecture~\ref{conj: geometric BM} and the weight part of Serre's
conjecture. For more details, see~\cite{2015arXiv150902527G}
(particularly Section~6).

We expect that the cycles~$Z_{\underline{k}}$ will be effective, in
the sense that they are combinations of the~$\cX_d^{\underline{k}}$
with non-negative coefficients. This expectation is borne out in all
known examples (see the following sections), and in any case would be
a consequence of standard conjectures about the Taylor--Wiles
method. Indeed, the local cycles~$Z_{\underline{k}}(\rhobar)$ of Section~\ref{subsec: relating
  different BM conjectures} are conjecturally the supports of certain
``patched modules'', and in particular are effective; see~\cite[Thm.\
5.5.2]{emertongeerefinedBM}. The effectivity of
the~$Z_{\underline{k}}(\rhobar)$ would immediately imply the
effectivity of the~$Z_{\underline{k}}$.

The ``weight part of Serre's conjecture'' is perhaps more of a
conjectural conjecture than an actual conjecture: it should assign to
each~$\rhobar:G_K\to\GL_d(\F)$ a set~$W(\rhobar)$ of Serre weights, \index{$W(\rhobar)$}
with the property that if~$\rhobar$ is the restriction to a
decomposition group of a suitable
global representation (for example, an irreducible representation
coming from an automorphic form on a unitary group), then
~$\underline{k}\in W(\rhobar)$ if and only if~$\underline{k}$ is a
weight for the global representation (for example, in the sense that
the global representation corresponds to some mod~$p$ cohomology class
for a coefficient system corresponding to~$F_{\underline{k}}$).

Many conjectural definitions of the sets~$W(\rhobar)$ have been
proposed. Following~\cite{geekisin}, one definition is to assume the
Breuil--M\'ezard conjecture, for example in the form~\eqref{eqn:
  refined version of BM cris}, and define $W(\rhobar)$ to be the set
of~$\underline{k}$ for which $\mu_{\underline{k}}(\rhobar)>0$. While
this is less explicit than other definitions, it has the merit that it
would follow from standard conjectures about modularity lifting
theorems that it gives the correct set of weights; see~\cite[\S3,
4]{2015arXiv150902527G}.

Assume Conjecture~\ref{conj: geometric BM}, and assume that the
cycles~$Z_{\underline{k}}$ are effective. As explained in Section~\ref{subsec: relating
  different BM conjectures}, it follows that the numerical
Breuil--M\'ezard holds for every~$\rhobar$, with
$\mu_{\underline{k}}(\rhobar)$ being the Hilbert--Samuel multiplicity
of the cycle~$Z_{\underline{k}}(\rhobar)$, which
(since~$Z_{\underline{k}}(\rhobar)$ is effective) is positive if and
only if~$Z_{\underline{k}}(\rhobar)$ is nonzero, i.e.\ if and only
if~$Z_{\underline{k}}$ is supported at~$\rhobar$. Thus we can rephrase
the Breuil--M\'ezard version of the weight part of Serre's conjecture
as saying that $W(\rhobar)$ is the set of~$\underline{k}$ such that
$Z_{\underline{k}}$ is supported at~$\rhobar$.

Alternatively, we can rephrase this conjecture in the following way:
to each irreducible component of~$\cX_{d,\red}$, we assign the set of
weights~$\underline{k}$ with the property that~$Z_{\underline{k}}$ is
supported on this component. Then~$W(\rhobar)$ is simply the union of
the sets of weights for the irreducible components of~$\cX_{d,\red}$
which contain~$\rhobar$. As we explain in Section~\ref{subsec: CEGS},
if~$d=2$ then this description agrees with the other definitions
of~$W(\rhobar)$ in the literature, and therefore gives a
geometrization of the weight part of Serre's conjecture.
\section{The case of $\GL_2(\Qp)$}\label{subsec: GL2 Qp}
The numerical Breuil--M\'ezard conjecture for~$K=\Qp$ and~$d=2$ is completely
known, thanks in large part to Kisin's paper~\cite{KisinFM} (which
gave a proof in many cases by a mixture of local and global
techniques), and Pa{\v{s}}k{\=u}nas' paper~\cite{paskunasBM} which
reproved these results by purely local means (the $p$-adic local
Langlands correspondence), relaxed the hypotheses on~$\rhobar$, and
also proved the refined version of the correspondence. The remaining
cases not handled by these papers are proved in the
papers~\cite{HuTan,sandermultiplicities,2018arXiv180307451T},
culminating in the proof of the final cases for~$p=2$ by Tung in~\cite{2019arXiv190806174T}.

Consequently, by the discussion of Section~\ref{subsec: relating
  different BM conjectures}, Conjecture~\ref{conj: geometric BM} holds
for~$K=\Qp$ and~$d=2$. It follows easily from the explicit
description of~$\mu_{\underline{k}}(\rhobar)$ in the papers cited
above (or alternatively from the description for~$\GL_2(K)$ in
Section~\ref{subsec: CEGS} below)
that~$Z_{\underline{k}}=\cX_2^{\underline{k}}$ unless
$\underline{k}=(a+p-1,a)$ for some~$a$, in which case
$Z_{(a+p-1,a)}=\cX_2^{(a+p-1,a)}+\cX_2^{(a,a)}$.


\section{$\GL_2(K)$: potentially Barsotti--Tate types}\label{subsec: CEGS}
In this section we assume
that $p$ is odd, 
and explain some consequences of the results
of~\cite{geekisin} (which proved the numerical Breuil--M\'ezard
conjecture for $2$-dimensional potentially Barsotti--Tate
representations) and~\cite{CEGSKisinwithdd} (which studied moduli stacks of rank~$2$
Breuil--Kisin modules with tame descent data). We will take advantage
of Remark~\ref{rem: geom BM lets us go from sufficiently generic
      rhobar to all rhobar}.

We begin with the following slight extension of one of the main
results of~\cite{geekisin}. Let~$\underline{\BT}$ denote the minimal regular Hodge type,
  i.e.\ we have~$\underline{\BT}_{\sigma,1}=1$,
  $\underline{\BT}_{\sigma,2}=0$ for all $\sigma:K\into\Qpbar$.
\begin{thm}
  \label{thm: Gee Kisin extended to semistable}Let $K/\Qp$ be a finite
  extension with $p>2$, and let $\rhobar:G_K\to\GL_2(\Qpbar)$ be
  arbitrary. Then the numerical Breuil--M\'ezard conjecture holds for
  potentially crystalline and potentially semistable lifts
  of~$\rhobar$ of Hodge type~$\underline{\BT}$ and arbitrary inertial type~$\tau$.

  More precisely, there are unique non-negative
  integers~$\mu_{\underline{k}}(\rhobar)$ such that \eqref{eqn:
    numerical version of BM cris} and \eqref{eqn: numerical version of
    BM ss} both hold for~$\underline{\lambda}=\underline{\BT}$ and
  ~$\tau$ arbitrary.
\end{thm}
\begin{proof}If we remove the potentially semistable case and consider
  only potentially crystalline representations, the theorem
  is~\cite[Thm.\ A]{geekisin}, which is proved as~\cite[Cor.\
  4.5.6]{geekisin}. We now briefly explain how to modify the proofs
  in~\cite{geekisin} to prove the more general result; as writing out
  a full argument would be a lengthy exercise, and the arguments are
  completely unrelated to those of this book, we only explain the key
  points. Examining the proof of~\cite[Cor.\ 4.5.6]{geekisin}, we see
  that we just need to verify that the assertion of the first sentence
  of the proof of~\cite[Thm.\ 4.5.5]{geekisin} holds in this setting,
  i.e.\ that the equivalent conditions of~\cite[Lem.\
  4.3.9]{geekisin} hold. Exactly as in the proof of~\cite[Cor.\
  4.4.3]{geekisin}, it is enough to show that every irreducible
  component of a product of local deformation rings is witnessed by an
  automorphic representation.

  By the usual Khare--Wintenberger argument, it suffices to prove this
  after making a solvable base change, and in particular we can
  suppose that all of the residual local Galois representations are
  trivial, the mod~$p$ cyclotomic character is trivial, the trivial
  mod~$p$ representation admits a non-ordinary crystalline lift,
  and the inertial types are all trivial. Now, any non-crystalline
  semistable representation of Hodge type~$\underline{\BT}$ is
  necessarily ordinary (indeed, it follows from a direct computation
  of the possible weakly admissible modules that all such
  representations are unramified twists of an extension of the inverse
  of the cyclotomic
  character by the trivial character), and by results of Kisin and Gee
  (see ~\cite[Prop.\ 4.3.1]{MR3778977}, \cite[Cor.\
  2.5.16]{KisinModularity}, and \cite[Prop.\ 2.3]{MR2280776}), as $p>2$ and~$\rhobar$ is
  trivial, the semistable ordinary deformation ring in question is a
  domain, as are the crystalline ordinary deformation ring, and the
  crystalline non-ordinary deformation ring (all in Hodge type~$\underline{\BT}$). It
  therefore suffices to show that in the situation of the proof
  of~\cite[Cor.\ 4.4.3]{geekisin}, given any decomposition of the set~$S$
  of places of~$F_1^+$ lying over~$p$ as~$S_{\textrm{ss}}\coprod S_{\textrm{crys-ord}}
  \coprod S_{\textrm{non-ord}}$, we can arrange to have a
  congruence to an automorphic representation~$\pi''$, having the properties
  that~$\pi''$ is unramified at the places lying over a place in~$S_{\textrm{crys-ord}}
  \coprod S_{\textrm{non-ord}}$,
  is an unramified twist of the Steinberg representation at the
  places  lying over a place in~$S_{\textrm{ss}}$, and is
  furthermore ordinary (resp.\ not ordinary) at the places lying over
  a place in ~$S_{\textrm{crys-ord}}$ (resp.\ $
  S_{\textrm{non-ord}}$). 

  The existence of such a representation if~$S_{\textrm{ss}}=S$
  follows from the construction of the global
  representation~$\overline{r}$ used in the proof of numerical
  Breuil--M\'ezard conjecture in~\cite{geekisin}; more precisely, it
  follows from~\cite[Prop.\ 8.2.1]{0905.4266}, which is applied in the
  proof of~\cite[Thm.\ A.2]{geekisin} (with the type function in the
  sense of~\cite{0905.4266} being~$C$ at all places above~$p$). Thus,
  to establish the general case, it suffices to establish the
  existence of suitable ``level lowering'' congruences. This is easily
  done by switching to a group which is ramified at the places
  lying over~$S_{\textrm{ss}}$, and then making
  congruences to automorphic representations which are either
  unramified and ordinary (in the case of places lying over
  $S_{\textrm{crys-ord}}$) or have cuspidal type (in the case of
  places lying over~$S_{\textrm{non-ord}}$), and then applying the
  Khare--Wintenberger argument again at the places lying
  over~$S_{\textrm{non-ord}}$. (See e.g.\ \cite[Lem.\
  3.5.3]{KisinModularity} for a similar argument for places away
  from~$p$.)
\end{proof}

We say that a Serre weight~$\underline{k}$ for~$\GL_2$ is \index{Serre weight!Steinberg}
\emph{Steinberg} if for each~$\sigmabar$ we have
$k_{\sigmabar,1}-k_{\sigmabar,2}=p-1$. If~$\underline{k}$ is Steinberg
then we define~$\underline{\tilde{k}}$ by
$\tilde{k}_{\sigmabar,1}=\tilde{k}_{\sigmabar,2}=k_{\sigmabar,2}$. 
\begin{thm}
  \label{thm: explicit BM cycles GL2 K}
 Continue to assume that~$d=2$ and~$p>2$. Then Conjecture~{\em \ref{conj: geometric BM}} holds
  for~$\underline{\lambda}=\underline{\BT}$ and~$\tau$ arbitrary, with
  the cycles~$Z_{\underline{k}}$ being as follows: if~$\underline{k}$
  is not Steinberg, then $Z_{\underline{k}}=\cX_2^{\underline{k}}$,
  while if $\underline{k}$ is Steinberg,
  $Z_{\underline{k}}=\cX_2^{\underline{k}}+\cX_2^{\underline{\tilde{k}}}$.

\end{thm}
\begin{proof}
We use the notation of Section~\ref{subsec: relating different BM
    conjectures}. By Theorem~\ref{thm: Gee Kisin extended to semistable} and the
  discussion of Section~\ref{subsec: relating different BM
    conjectures}, we need only show that the
  cycles~$Z_{\underline{k}}$ in the statement of the theorem are those
  determined by~\eqref{eqn: cycle from multiplicities}. Suppose
  firstly that~$\underline{k}$ is not Steinberg. Then by~\cite[Thm.\
  5.2.2~(2)]{CEGSKisinwithdd}, $\cX_2^{\underline{k}}$ has a dense set of
  finite type points with the property that their only non-Steinberg
  Serre weight is~$\underline{k}$. It follows from this and~\cite[Thm.\
  5.2.2~(3)]{CEGSKisinwithdd}  that if
  neither~$\underline{k}$ nor~$\underline{k}'$ is Steinberg, then
  $\mu_{\underline{k}}(\rhobar_{\underline{k}'})=\delta_{\underline{k},\underline{k}'}$.
  By construction, if~$\underline{k}'$ is Steinberg,
  then~$\rhobar_{\underline{k}'}$ is a twist of a tr\`es ramifi\'ee
  extension of the trivial character by the mod~$p$ cyclotomic
  character. By~\cite[Lem.\ B.5]{CEGSKisinwithdd}, this implies that
  $\mu_{\underline{k}}(\rhobar_{\underline{k}'})=0$.
  The cycles
  $Z_{\underline{k}}$ are therefore as claimed if~$\underline{k}$ is
  non-Steinberg. 

  It remains to determine the values of
  $\mu_{\underline{k}}(\rhobar_{\underline{k}'})$ in the case
  that~$\underline{k}$ is Steinberg. By twisting, we can and do assume
  that~$k_{\sigmabar,2}=0$ for all~$\sigmabar$.  If we apply
  Theorem~\ref{thm: Gee Kisin extended to semistable} with $\tau$ being
  the trivial type, and recall  that $L_{\underline{\BT}}$  is the trivial representation
  and that $\sigma^{\semis}(\tau)$ is the Steinberg type, the reduction
  of which  is precisely the representation $F_{\underline{k}}$,
  we find that 
\numequation
\label{eqn:instance of BM}
Z_{\semis, \underline{\BT}, \tau} = Z_{\underline{k}}.
\end{equation}
  Thus
  $\mu_{\underline{k}}(\rhobar_{\underline{k}'})\ne 0$ if and only
  if~$\rhobar_{\underline{k}'}$ admits a semistable lift of Hodge
  type~$\underline{\BT}$. If this lift is in fact crystalline, then
  ~$\rhobar_{\underline{k}'}$ has~$\underline{\tilde{k}}$ as a Serre
  weight, so by another application of~\cite[Thm.\
  5.2.2~(2)]{CEGSKisinwithdd}, we see that either
  $\underline{k}'=\underline{\tilde{k}}$ or else that $\underline{k}'$ is
  also Steinberg.
  In the former case,
  $\rhobar_{\underline{k}'}$ is (by Remark~\ref{rem: open eigenvalue morphism implies ratio
          not cyclo}) an unramified twist of an extension of
  inverse of the mod~$p$ cyclotomic
  character by a
  non-trivial unramified character, so it does not admit a semistable non-crystalline lift of
  Hodge type~$\underline{\BT}$ (as all such lifts are unramified
  twists of an extension of the inverse of the cyclotomic  character by the trivial
  character); so we have
  $\mu_{\underline{k}}(\rhobar_{\underline{\tilde{k}}})=\mu_{\underline{\tilde{k}}}(\rhobar_{\underline{\tilde{k}}})=1$. 

  Finally, we are left with the task of
  computing $\mu_{\underline{k}}(\rhobar_{\underline{k}'})$
  when both $\underline{k}$ and $\underline{k'}$  are Steinberg.
Recall that by twisting, we are
  assuming that~$k_{\sigmabar,2}=0$ for all~$\sigmabar$.
The weight $\underline{k}'$ is then a twist of~$\underline{k}$.
We have to show that in this case we again
  have $\mu_{\underline{k}}(\rhobar_{\underline{k}'}) = \delta_{\underline{k},
\underline{k}'}.$
To see this, first note that
  since $\rhobar_{\underline{k}'}$  is  tr\`es ramifi\'ee, it does not admit a
  crystalline lift of Hodge type~$\underline{\BT}$, so all of its
  semistable lifts of Hodge type~$\underline{\BT}$ are given by
  unramified twists of  extensions of the inverse of the cyclotomic character by
  the trivial character.   Furthermore, the  reduction of any lattice in such
  an extension  is an  
  unramified twists of an extension of the inverse of the mod $\varpi$
  cyclotomic character by the trivial character, and hence  cannot equal
  $\rhobar_{\underline{k}'}$ unless $\underline{k}'$ equals $\underline{k}$
  (rather than being a non-trivial twist of it).
  If we take into account~\eqref{eqn:instance of BM}, we find
  that indeed 
  $\mu_{\underline{k}}(\rhobar_{\underline{k}'}) = 0$ when
  $\underline{k}'  \neq  \underline{k}$.
  
  Finally, any lift of $\rho_{\underline{k}}$ which is an 
  unramified twist of an
  extension of the inverse of the cyclotomic character by
  the trivial character is 
  automatically semistable of Hodge type~$\underline{\BT}$, so that
  $R^{\ss,\underline{\BT},\tau}_{\rhobar{\underline{k}}}$ is
  precisely the ordinary (framed) deformation ring parameterizing
  such lifts of $\rhobar_{\underline{k}}$.
  A standard Galois cohomology
  calculation (that we leave to the reader)
  shows that this ring is formally smooth,
  and thus that
  $R^{\ss,\underline{\BT},\tau}_{\rhobar{\underline{k}}}/\varpi$ is
  also formally smooth.
  It follows that $\mu_{\underline{k}}(\rhobar_{\underline{k}})  = 1$,
  as claimed.
%
%
%
\end{proof}
The following lemma makes precise which Galois representations occur
on Steinberg components. Note that the twist in the statement of the
lemma is determined by the values of~$k_{\underline{\sigma},2}$.

\begin{lem}
  \label{lem: points of Steinberg components are 1 by
    cyclo}If~$\underline{k}$ is Steinberg, then the $\Fpbar$-points of
  $\cX_2^{\underline{k}}$ are twists of an extension of the inverse of
  the mod~$p$ cyclotomic character by the trivial character.
\end{lem}
\begin{proof}Twisting, we may assume that~$k_{\underline{\sigma},2}=0$
  for all~$\sigmabar$. Let $\rhobar:G_K\to\GL_2(\F)$ correspond to a
  closed point of~$\cX_2^{\underline{k}}$. We claim that~$\rhobar$ is
  an unramified twist of an extension of~$\epsilonbar$ by~$1$. To
  see this, let~$R^{\underline{\BT},\operatorname{St}}_{\rhobar}$
  denote the $\Zp$-flat quotient
  of~$R^{\underline{\BT},\semis}_{\rhobar}$ determined by the
  irreducible components of the generic fibre which are not components
  of~$R^{\underline{\BT},\crys}_{\rhobar}$.

  It follows from Theorem~\ref{thm: explicit BM cycles GL2 K} that the
  cycle of
  $\Spec R^{\underline{\BT},\operatorname{St}}_{\rhobar}/\varpi$ is
  nonzero, so that in particular~$\rhobar$ must admit a semistable
  non-crystalline lift of Hodge type~$\BT$. As in the proof of
  Theorem~\ref{thm: Gee Kisin extended to semistable}, any such
  representation is an unramified twist of an extension
  of~$\epsilon^{-1}$ by the trivial character, so we are done. (We
  could presumably also phrase this argument on the level of the
  crystalline and semistable moduli stacks, but have chosen to present
  it in the more familiar setting of deformation rings.)
  \end{proof}

\begin{rem}
  \label{rem: we could compare to CEGS but we won't}The
  paper~\cite{CEGSKisinwithdd} constructs and studies moduli stacks of
  two-dimensional representations of~$G_{K_{\piflat,\infty}}$ (for a
  fixed choice of~$\piflat$) which are tamely potentially of height at
  most~$1$.  Since restriction from Barsotti--Tate representations
  of~$G_K$ to representations of~$G_{K_{\piflat,\infty}}$ of height at
  most~$1$ is an equivalence (by the results
  of~\cite{KisinModularity}), these stacks can be interpreted as
  stacks of~$G_K$-representations (and indeed are interpreted as such
  in~\cite{CEGSKisinwithdd}). It is presumably straightforward (using
  that all the stacks under consideration are of finite type, and are
  $\Zp$-flat and reduced, in order to reduce to a comparison of points
  over finite extensions of~$\Zp$) to identify them with the stacks
  considered in this book (for~$\underline{\lambda}=\underline{\BT}$
  and~$\tau$ a tame inertial type). (We were able to apply the results
  of~\cite{CEGSKisinwithdd} in the proof of Theorem~\ref{thm: explicit
    BM cycles GL2 K} without doing this because we only needed to use them on the level
  of versal rings, which are given by universal Galois lifting rings
  for both our stacks and those of~\cite{CEGSKisinwithdd}.) It may
  well be the case that the stacks of Kisin modules considered
  in~\cite{CEGSKisinwithdd} can be identified with certain of the
  moduli stacks of Breuil--Kisin--Fargues modules 
  that we defined in Section~\ref{sec: semistable stack}. 
  Since we do not need to know
  this, we leave it as an exercise for the interested reader.
\end{rem}

\subsection{A lower bound}The patching arguments of~\cite{geekisin} easily
imply a general inequality, as we now record.

\begin{prop}
  \label{prop: GL2 geometric BM inequality}Assume that~$p>2$. Let~$d=2$, and let the~$Z_{\underline{k}}$ be
  as in the statement of Theorem~{\em \ref{thm: explicit BM cycles GL2
    K}}. Then for each regular Hodge type~$\underline{\lambda}$ and
  each inertial type~$\tau$, we have 
  $Z_{\crys,\lambdau,\tau}\ge\sum_{\underline{k}}n_{\underline{k}}^\crys(\lambda,\tau)\cdot
  Z_{\underline{k}}$,  and $Z_{\semis,\lambdau,\tau}\ge\sum_{\underline{k}}n_{\underline{k}}^\semis(\lambda,\tau)\cdot
  Z_{\underline{k}}$.
\end{prop}
\begin{rem}
  The meaning of the inequalities in the statement of
  Proposition~\ref{prop: GL2 geometric BM inequality} is the obvious
  one: each side is a linear combination of irreducible components,
  and the assertion is that for each irreducible component, the
  multiplicity on the left hand side is at least the multiplicity on
  the right hand side.
\end{rem}
\begin{proof}[Proof of Proposition~{\em \ref{prop: GL2 geometric BM
    inequality}}]As in Section~\ref{subsec: relating
  different BM conjectures}, it is enough to show that for
each~$\rhobar$, we have \[e(\Spec
  R_{\rhobar}^{\crys,\underline{\lambda},\tau}/\varpi)\ge\sum_{\underline{k}}n_{\underline{k}}^\crys(\lambda,\tau)\mu_{\underline{k}}(\rhobar),\] \[e(\Spec
  R_{\rhobar}^{\semis,\underline{\lambda},\tau}/\varpi)\ge\sum_{\underline{k}}n_{\underline{k}}^\semis(\lambda,\tau)\mu_{\underline{k}}(\rhobar),\]where
the~$\mu_{\underline{k}}(\rhobar)$ are the uniquely determined
integers from~\cite{geekisin}. (In fact, it is enough to show these
inequalities when $\rhobar=\rhobar_{\underline{k}}$ for
some~$\underline{k}$, but assuming this does not simplify our
arguments.) Since~\cite{geekisin} considers only potentially
crystalline representations, we only give the proof
for~$R_{\rhobar}^{\crys,\underline{\lambda},\tau}$; the proof in the potentially
semistable case is essentially identical, and we leave it to the reader.

We now examine the proof of~\cite[Thm.\ 4.5.5]{geekisin},
setting~$\rbar$ there to be our~$\rhobar$. As we've done throughout this book, we
write~$\underline{k}$ rather than~$\sigma$ for Serre weights. If we
ignore the assumption that every potentially crystalline lift of Hodge
type~$\lambdau$ and inertial type~$\tau$ is potentially
diagonalizable, then we do not know that the equivalent conditions
of~\cite[Lem.\ 4.3.9]{geekisin} hold, but examining the proof
of~\cite[Lem.\ 4.3.9]{geekisin}, we do know that the statement of
part~(4) of \emph{loc.\ cit.} can be replaced with an inequality (with
equality holding if and only if~$M_\infty$ is a faithful
$R_\infty$-module).

Returning to the proof of~\cite[Thm.\ 4.5.5]{geekisin}, we note that
the definition of~$\mu_{\underline{k}}(\rhobar)$ is such that the
quantity $\mu'_{\sigma_{\mathrm{gl}}}(\rbar)$ is simply the product of
the corresponding~$\mu_{\underline{k}}(\rhobar)$. In particular, if we
take $\lambda_v=\lambdau$ and~$\tau_v=\tau$ for each~$v$, then
(bearing in mind the discussion of the previous paragraph), the main
displayed equality becomes an inequality
\[e(\Spec
  R_{\rhobar}^{\crys,\underline{\lambda},\tau}/\varpi)^N\ge\left(\sum_{\underline{k}}n_{\underline{k}}^\crys(\lambda,\tau)\mu_{\underline{k}}(\rhobar)\right)^N, \]where
there are~$N$ places of~$F^+$ lying over~$p$. Since each side is
non-negative, the result follows.
\end{proof}
\begin{rem}\label{rem: BM inequality well known, doesn't obviously
    work for GLd though.}
  The argument used to prove Proposition~\ref{prop: GL2 geometric BM
    inequality} is well known to the experts, and goes back
  to~\cite[Lem.\ 2.2.11]{KisinFM}. Despite the relatively formal
  nature of the argument, it does not seem to be easy to prove an
  analogous statement for~$\GL_d$ with~$d>2$; the difficulty is in the
  step where we used~\cite[Thm.\ 4.5.5]{geekisin} to replace a global
  multiplicity with a product of local multiplicities. Without this
  argument (which crucially relies on the modularity lifting theorems
  of~\cite{KisinModularity}) one only obtains inequalities involving
  cycles in products of copies of~$\cX_d$.
\end{rem}
\section{Brief remarks on $\GL_{\lowercase{d}}$, $\lowercase{d}>2$}\label{subsec: GLn}We
expect the situation for~$\GL_d$, $d>2$, to be considerably more
complicated than that for~$\GL_2$. Experience to date suggests that
the weight part of Serre's conjecture in high dimension is
consistently more complicated than is anticipated, 
and so it seems unwise to
engage in much speculation. In particular, the results of~\cite{le2020local} show that even for generic
weights~$\underline{k}$, we should not expect the
cycles~$Z_{\underline{k}}$ to have as simple a form as those
for~$\GL_2$.

In the light of Lemma~\ref{lem: maxnonsplit weight k implies crystalline lifts to
    shifts}, 
  it seems reasonable to expect~$\cX_d^{\underline{k}}$ to contribute
  to~$Z_{\underline{k}}$, and it also seems reasonable to expect
  contributions from the ``shifted'' weights, as in
  Definition~\ref{defn: shifted weight} (see also~\cite[\S
  7.4]{2015arXiv150902527G}). The comparative weakness of the
  automorphy lifting theorems available to us in dimension greater
  than~$2$ prevents us from saying much more than this, although we
  refer the reader to~\cite[\S\S 3,4,6]{2015arXiv150902527G} for a
  discussion of the conjectural general relationship between the
  numerical Breuil--M\'ezard conjecture, the weight part of Serre's
  conjecture, automorphy lifting theorems, and the stacks~$\cX_d$.

\renewcommand{\theequation}{\Alph{chapter}.\arabic{section}}
\appendix

\chapter{Formal algebraic stacks}\label{app: formal algebraic
  stacks}The theory of formal
algebraic stacks is developed in~\cite{Emertonformalstacks}. In this
appendix we briefly summarise the parts of this theory that are used in
the body of the book, and also introduce some additional terminology
and establish some additional results which we will require.

\subsection*{Algebraic stacks}We follow the terminology of~\cite{stacks-project}; in
particular, we write ``algebraic stack'' rather than ``Artin
stack''. More \index{algebraic stack}
precisely, an algebraic stack is a stack in groupoids in the \emph{fppf} topology,
whose diagonal is representable by algebraic spaces, which admits a smooth
surjection from a
scheme. See~\cite[\href{http://stacks.math.columbia.edu/tag/026N}{Tag
  026N}]{stacks-project} for a discussion of how this definition relates to
others in the literature. If~$S$ is a scheme, then by ``a stack
over~$S$'' we mean  a stack fibred in groupoids over the big
\emph{fppf} site of~$S$.

	We say that a morphism $\cX \to \cY$ of stacks over $S$
	is {\em representable by algebraic stacks} \index{representable by algebraic stacks}
	if for any morphism of stacks $\cZ \to \cY$ whose source is an algebraic
	stack, the fibre product $\cX \times_{\cY} \cZ$
	is again an algebraic stack.
(In the Stacks Project, the terminology {\em algebraic} is
used instead~\cite[\href{http://stacks.math.columbia.edu/tag/06CF}{Tag
  06CF}]{stacks-project}.) 
Note that a morphism from a sheaf to a
stack is representable by algebraic stacks if and only if it is representable
by algebraic spaces (this is easily verified,
or see e.g.\ \cite[Lem.~3.5]{Emertonformalstacks}). 

Following~\cite[\href{http://stacks.math.columbia.edu/tag/03YK}{Tag 03YK},
\href{http://stacks.math.columbia.edu/tag/04XB}{Tag 04XB}]{stacks-project},
we can define properties of morphisms representable by algebraic stacks
in the following way.

\begin{adefn}
	\label{df:properties of morphisms representable by alg stacks}
	If $P$ is a property of morphisms of algebraic stacks which is 
	{\em fppf} local on the target, and preserved by arbitrary
	base-change, then we say that a morphism $f:\cX \to \cY$ of
	stacks which is representable by algebraic stacks {\em has property
		$P$} if and only if for every algebraic stack $\cZ$
	and morphism $\cZ \to \cY$, the base-changed morphism of algebraic
	stacks $\cZ \times_{\cY} \cX \to \cZ$ has property $P$.
\end{adefn}

	When applying
	Definition~\ref{df:properties of morphisms representable by alg stacks},
	it suffices to consider the case when $\cZ$ is actually a scheme,
	or even an affine scheme.

Some properties $P$ to which we can apply
	Definition~\ref{df:properties of morphisms representable by alg stacks}
are being {\em locally of finite type}, {\em locally of finite presentation},
and {\em smooth}.   The following lemma provides alternative descriptions of
the latter two properties.

\begin{alemma}
\label{lem:characterizing properties}
Let $f:\cX \to \cY$ be a morphism of stacks which is representable by
algebraic stacks.
\begin{enumerate}
\item The morphism $f$ is locally of finite presentation if and only 
if it is {\em limit preserving on objects}, in the sense
of~\cite[\href{http://stacks.math.columbia.edu/tag/06CT}{Tag 06CT}]{stacks-project}.
\item The morphism $f$ is smooth if and only if it is locally of finite 
presentation and {\em formally smooth},
in the sense that it satisfies the usual infinitesimal lifting property:
  for every affine $\cY$-scheme $T$, and every closed
  subscheme $T_0\into T$ defined by a nilpotent ideal sheaf, the functor
  $\Hom_\cY(T,\cX)\to\Hom_\cY(T_0,\cX)$ is essentially
  surjective.
\end{enumerate}
\end{alemma}
\begin{proof}
We first recall that a morphism of algebraic stacks is locally of
finite presentation if and only if it is limit preserving on 
objects~\cite[Lem.~2.3.16]{EGstacktheoreticimages}.

Suppose now that $f$ is locally of finite presentation,
that $T$ is an affine scheme, written as a projective limit of affine schemes
$T \iso \varprojlim T_i$, and suppose given a commutative diagram
$$\xymatrix{T  \ar[r] \ar[d] & T_i \ar[d] \\ \cX \ar[r] & \cY}$$
for some value of~$i$.
We must show that we can factor the left-hand vertical arrow through~$T_{i'}$,
for some $i'\geq i$, so that the evident resulting diagram again commutes.
The original diagram induces a morphism $T \to T_i\times_{\cY} \cX,$
and obtaining a factorization of the desired type amounts to obtaining
a factorization of this latter morphism through some~$T_{i'}$.  
In other words, we must show that the morphism
$T_i\times_{\cY} \cX \to T_i$
is limit preserving on objects.  But this is a morphism of algebraic
stacks which is locally of finite presentation,
and so it is indeed limit preserving on objects.

Conversely, suppose that $f$ is limit preserving on objects.
We must show that for any morphism $T \to \cY$, the base-changed morphism
$T\times_{\cY} \cX\to T$ is locally of finite presentation.
However, this base-change is again limit preserving on
objects~\cite[\href{https://stacks.math.columbia.edu/tag/06CV}{Tag 06CV}]{stacks-project},
and since it is a morphism of algebraic stacks, it is indeed locally
of finite presentation.

We now turn to~(2). By Definition~\ref{df:properties of morphisms
  representable by alg stacks} and~\cite[Cor.\ B.9]{MR2818725} (that
is, by the result at hand in the case that~$\cY$ is itself an
algebraic stack), it is enough to show that~$f$ satisfies the
infinitesimal lifting property if and only if the base changed
morphism $Z\times_\cY\cX\to Z$ satisfies the infinitesimal lifting
property for all morphisms $Z\to\cY$ whose source is a scheme.

Assume first that ~$f$ satisfies the
infinitesimal lifting property, and let $T_0\to T$ be a nilpotent closed
immersion whose target is an affine $Z$-scheme. Given a morphism
$T_0\to Z\times_\cY\cX$, the infinitesimal lifting property lets us lift the composite $T_0\to
Z\times_\cY\cX\to \cX$ to a morphism~$T\to\cX$, and since we also have
a morphism $T\to Z$, we obtain the required morphism $T\to
Z\times_\cY\cX$.

Conversely, assumed that the base changed
morphism $Z\times_\cY\cX\to Z$ satisfies the infinitesimal lifting
property for all~$Z\to\cY$ with $Z$ a scheme, and let~$T_0\to T$ be as in the statement
of the proposition. Taking~$Z=T$, we can lift the induced morphism
$T_0\to T\times_{\cY}\cX$ to a morphism $T\to T\times_{\cY}\cX$, and
the composite of this morphism with the second projection gives us the
required lifting.\end{proof}

\subsection*{Formal algebraic spaces}
\index{formal algebraic space}
Following~\cite[\href{http://stacks.math.columbia.edu/tag/0AHW}{Tag
    0AHW}]{stacks-project}, an
affine formal algebraic space over a base scheme~ $S$
is a sheaf $X$ on the \emph{fppf} site of $S$
which admits a description as
an Ind-scheme
$X \iso \varinjlim_i X_i,$ where the $X_i$ are affine schemes
and the transition morphisms 
are thickenings 
(in the sense 
of~\cite[\href{http://stacks.math.columbia.edu/tag/04EX}{Tag
  04EX}]{stacks-project}). Here we allow the indexing set in the
inductive limit to be arbitrary; if it can be chosen to be the natural
numbers, then we say that~$X$ is \emph{countably
  indexed}.
A countably indexed affine formal algebraic space can be written in
the form $X \cong \varinjlim_n \Spec A/I_n$, where $A$ is a complete
\index{countably indexed}
topological ring equipped with a decreasing sequence $\{I_n\}_n$ of open 
ideals which are weak ideals of definition, i.e.\ consist of topologically
nilpotent elements 
(\cite[\href{https://stacks.math.columbia.edu/tag/0AMV}{Tag 0AMV}]{stacks-project}),
and which
form a fundamental basis of $0$ in~$A$ (see the discussion
of~\cite[\href{https://stacks.math.columbia.edu/tag/0AIH}{Tag 0AIH}]{stacks-project}).
We then write $X := \Spf A$
(following~\cite[\href{https://stacks.math.columbia.edu/tag/0AIF}{Tag 0AIF}]{stacks-project}).

A particular example of a countably indexed affine formal algebraic
space is an adic affine formal algebraic space. By definition, this
is of the form~$\Spf A$, where~$A$ is adic: that is, $A$ is a topological ring which is
complete (and separated, by our convention throughout this book), and
which admits an ideal of definition (that is, an ideal~$I$ whose
powers form a basis of open neighbourhoods of zero). We say that $\Spf
A$ is Noetherian if $A$ is adic and Noetherian. 
We say that~$A$
(or $\Spf A$) is adic* if~$I$ can be taken to be finitely
generated. \index{adic*} 

More generally, we say that a complete topological ring~$A$ is weakly
\index{weakly admissible}
admissible if~$A$ contains an open ideal~$I$ consisting of
topologically nilpotent elements, \cite[\href{http://stacks.math.columbia.edu/tag/0AMV}{Tag
  0AMV}]{stacks-project}). 

A formal algebraic space over $S$ is a sheaf $X$ on the \emph{fppf}
site of $S$ which receives a morphism $\coprod U_i \to X$ which is
representable by schemes, \'etale, and surjective, and whose source is
a disjoint union of affine formal algebraic spaces~$U_i$. We say
that~$X$ is \emph{locally countably indexed} if the~$U_i$ can be
\index{locally countably indexed}
chosen to be countably indexed. We say that~$X$ is \emph{locally
  Noetherian} if the~$U_i$ can be taken to be Noetherian. 

We will find the
following lemmas useful.

\begin{alem}
  \label{lem: alem closed formal algebraic scheme adic*}Let~$X$
  be an affine formal algebraic space over~$S$,
which is either countably indexed {\em (}e.g.\ an adic or adic* affine
formal algebraic space{\em )},
and hence of the form $\Spf A$ for some weakly admissible topological ring~$A$,
by~\cite[\href{https://stacks.math.columbia.edu/tag/0AIK}{Tag 0AIK}]{stacks-project},
or else is of the form $\Spf A$, for a pro-Artinian ring $A$.
  Then if~$Y\to X$ is a
  closed immersion of formal algebraic spaces over~$S$, we have that~$Y$ is of
  the form~$\Spf B$ for some quotient~$B$ of~$A$ by a closed ideal, 
  endowed with its quotient topology.  
\end{alem}
\begin{proof}
In the countably indexed case,
  this is immediate from ~\cite[\href{https://stacks.math.columbia.edu/tag/0AIK}{Tag 0AIK},\href{http://stacks.math.columbia.edu/tag/0ANQ}{Tag
  0ANQ},\href{http://stacks.math.columbia.edu/tag/0APT}{Tag
  0APT}]{stacks-project}.
If $A$ is pro-Artinian, say $A = \varprojlim A_i$ with the $A_i$ Artinian,
then $\Spec A_i \times_{\Spf A} Y$ is a closed subscheme of $\Spec A_i$ for 
each index~$i$, and hence is of the form $\Spec B_i$ for some quotient 
$B_i$ of $B$.
Thus $$Y \iso \varinjlim_i \Spec A_i\times_{\Spf A} Y \iso \varinjlim \Spec B_i
\iso \Spf B,$$
where $B = \varprojlim_i B_i$.  The general theory of pro-Artinian rings
shows that $B$ is a quotient of $A$, since each $B_i$ is a quotient
of the corresponding~$A_i$~\cite[Rem.~2.2.7]{EGstacktheoreticimages}.
\end{proof}

\begin{alem}\emph{(\cite[Lem.\ 8.18]{Emertonformalstacks})}\label{alem: Noetherian flat affine formal algebraic spaces}	A morphism 
  of Noetherian affine formal algebraic spaces $\Spf B \to \Spf A$
  which is representable by algebraic spaces is {\em (}faithfully{\em
    )} flat if and only if the corresponding morphism $A \to B$ is
  {\em (}faithfully{\em )} flat.
  
\end{alem}

\subsection*{Ind-algebraic and formal algebraic
  stacks}\index{Ind-algebraic stack}
\begin{adefn}\label{adefn: Ind algebraic stack}(\cite[Defn.\ 4.2]{Emertonformalstacks})
  An \emph{Ind-algebraic stack}
  over a scheme~$S$ is a stack~$\cX$ over~$S$ which can be written as
  $\cX\cong\varinjlim_{i\in I}\cX_i$, where we are taking the
  2-colimit in the 2-category of stacks of a 2-directed system
  $\{\cX_i\}_{i\in I}$ of algebraic stacks over~$S$.
\end{adefn}
By \cite[Rem.\
  4.9]{Emertonformalstacks}, if $\cX\to\cY$ is representable by
  algebraic spaces, and~$\cY$ is an Ind-algebraic stack, then~$\cX$ is
  also an Ind-algebraic stack.
\begin{adefn}\label{adefn: formal algebraic stack}(\cite[Defn.\ 5.3]{Emertonformalstacks})
  A \emph{formal algebraic stack} over~$S$ is a stack~$\cX$ which
  admits a \index{$p$-adic formal algebraic stack}
  morphism $U\to\cX$ which is representable by algebraic spaces,
  smooth, and surjective, and whose source is a formal algebraic
  space. If the source can be chosen to be locally countably indexed,
  then we say that~$\cX$ is locally countably indexed.
\end{adefn}

\begin{adefn}\label{adefn: p adic formal algebraic stack}(\cite[Defn.\ 7.6]{Emertonformalstacks})
  A \emph{$p$-adic formal algebraic stack} is a formal algebraic stack
  $\cX$ over~$\Spec\Zp$ which admits a morphism $\cX\to\Spf\Zp$ which
  is representable by algebraic stacks; so in particular, we may
  write~$\cX$ as an Ind-algebraic stack by writing
  $\cX\cong\varinjlim_n\cX\times_{\Spf\Zp}\Spec\Z/p^n\Z$.
\end{adefn}
The relationship between formal algebraic stacks and Ind-algebraic
stacks is discussed in detail
in~\cite[\S6]{Emertonformalstacks}. Before explaining some of this
material, we need to recall some  preliminaries on finiteness
properties and underlying reduced substacks.
\subsection*{Underlying reduced substacks}If~$\cX/S$ is any stack,
then we let $(\cX_{\red})'$ be the full subcategory of~$\cX$ whose set
of objects consists of those $T\to\cX$ for which~$T$ is a reduced
$S$-scheme.
\begin{adefn}\label{adefn: underlying reduced substack}(\cite[Lem.\ 3.27]{Emertonformalstacks})
  The \emph{underlying reduced substack} $\cX_{\red}$ of~$\cX$ is the
  intersection of all of the substacks of~$\cX$ which
  \index{underlying reduced substack}
  contain~$(\cX_{\red})'$. 
\end{adefn}

\begin{alem}\label{alem: underlying reduced of Ind algebraic}\emph{(}\cite[Lem.\ 4.16]{Emertonformalstacks}\emph{)}
  If~$\cX$ is an Ind-algebraic stack, and we write
  $\varinjlim_i \cX_i \iso \cX$, for some $2$-directed system
  $\{\cX_i\}_{i \in \cI}$ of algebraic stacks, then the induced
  morphism $\varinjlim_i (\cX_i)_{\red} \to \cX_{\red}$ is an
  isomorphism, and so in particular $\cX_{\red}$ is again an
  Ind-algebraic stack.
\end{alem}
As the following lemma records, in the case that~$\cX$ is a formal
algebraic stack, then~$\cX_{\red}$ is algebraic, and (just as in the
case of formal schemes, or formal algebraic spaces), $\cX$ is a
thickening of~$\cX_{\red}$.
\begin{alemma}\emph{(}\cite[Lem.\ 5.26]{Emertonformalstacks}\emph{)}
	\label{alem:underlying reduced}
	If $\cX$ is is a formal algebraic stack over $S$,
	then $\cX_{\red}$ is a closed and reduced algebraic substack of $\cX$, 
	and the inclusion $\cX_{\red} \hookrightarrow \cX$
	is a thickening {\em (}in the sense that its base-change 
	over any algebraic space induces a thickening of algebraic
	spaces{\em )}.   Furthermore,
	any morphism $\cY \to \cX$ with $\cY$ a reduced
	algebraic stack factors through~$\cX_{\red}.$
\end{alemma}

\subsection*{Finiteness properties}Let~$\cX$ be a formal algebraic
stack~$\cX$. We say that~$\cX$ is quasi-compact if the algebraic
stack~$\cX_{\red}$ is quasi-compact. We say that~$\cX$ is
quasi-separated if the diagonal morphism of~$\cX$ (which is
automatically representable by algebraic spaces) is quasi-compact and
quasi-separated.

We say that ~$\cX$ is locally Noetherian if and only if it admits a
morphism $U\to\cX$ which is representable by algebraic spaces, smooth,
and surjective, whose source is a locally Noetherian formal algebraic
space. Finally, we say that~$\cX$ is Noetherian if it is locally
Noetherian, and it is quasi-compact and quasi-separated.

\subsection*{The relationship between Ind-algebraic and formal algebraic stacks}
In
one direction, we have the following lemma.

\begin{alem}\label{alem: formal implies Ind}\emph{(}\cite[Lem.\ 6.2]{Emertonformalstacks}\emph{)}
  	If $\cX$ is a quasi-compact and quasi-separated
	formal algebraic stack, then $\cX \cong \varinjlim_i \cX_i$ for
	a $2$-directed system $\{\cX_i\}_{i \in I}$ of quasi-compact and quasi-separated
algebraic
	stacks in which the transition
	morphisms are thickenings. 
      \end{alem}

      A partial converse to Lemma~\ref{alem: formal implies Ind} is
      proved as~\cite[Lem.\ 6.3]{Emertonformalstacks}. In particular,
      we have the following useful criteria for being a 
      formal algebraic stack.

\begin{aprop}
  \label{prop: formal stack 6.6}\emph{(}\cite[Cor.\ 6.6]{Emertonformalstacks}\emph{)} Suppose that $\cX$ is an
  Ind-algebraic stack that can be written as the $2$-colimit
  $\cX \iso \varinjlim \cX_n$ of a directed sequence
  $(\cX_n)_{n \geq 1}$ in which the $\cX_n$ are algebraic stacks, and
  the transition morphisms are closed immersions.  If $\cX_{\red}$ is
  a quasi-compact algebraic stack, then $\cX$ is a locally countably
  indexed formal algebraic stack.
\end{aprop}

\begin{aprop}
  \label{prop: criterion for p-adic formal algebraic stack}Suppose
  that~$\{\cX_n\}_{n\ge 1}$ is an inductive system of algebraic stacks
  over~$\Spec\Zp$,
such that each~$\cX_n$ in fact lies over~$\Spec\Z/p^n\Z$, and each
of the induced morphisms
$\cX_n\to\cX_{n+1}\times_{\Spec\Z/p^{n+1}\Z}\Z/p^n\Z$
is an isomorphism. 
Then the Ind-algebraic stack $\cX:=\varinjlim_n\cX_n$
  is a $p$-adic formal algebraic stack.
\end{aprop}
\begin{proof}
  This is a special case of~\cite[Lem.\ 6.3, Ex.\ 7.8]{Emertonformalstacks}.
\end{proof}


\subsection*{Scheme-theoretic images}The notion of the
scheme-theoretic image of a (quasi-compact) morphism of algebraic stacks
is developed
in~\cite[\href{https://stacks.math.columbia.edu/tag/0CMH}{Tag
  0CMH}]{stacks-project}; equivalent alternative presentations are
given in~\cite[\S 3.1]{EGstacktheoreticimages} and in~\cite[Ex.\
9.9]{Emertonformalstacks}. In~\cite[\S 6]{Emertonformalstacks} this is
extended to a definition of the scheme-theoretic image
between formal algebraic stacks that are quasi-compact and quasi-separated.
In fact, the definition there uses 
        Lemma~\ref{alem: formal implies Ind} to write the formal stacks
in question as Ind-algebraic stacks with transition morphisms 
being closed immersions, and this is the level of generality
that is appropriate to us here.

A robust theory of scheme-theoretic images requires
a quasi-compactness assumption on the morphism whose scheme-theoretic image
is being formed, and we begin by introducing the notion of quasi-compactness
which seems appropriate to our context.

\begin{adefn}
\label{def:ind-q.c.}
Let $\cX$ be an Ind-algebraic stack
which can be written as a $2$-colimit of quasi-compact algebraic stacks
for which the transition morphisms are closed immersions, say
        $\cX \cong \varinjlim \cX_{\lambda}$.
We say that a morphism of stacks $\cX \to \cY$
is {\em Ind-representable by algebraic stacks} if each of the induced
morphisms
\anumequation
\label{eqn:restricted morphism}
\cX_{\lambda} \to \cY
\end{equation}
is representable by algebraic stacks.
We say that such a morphism is furthermore {\em Ind-quasi-compact}
if each of the morphisms $\cX_{\lambda} \to \cY$ (representable by
algebraic stacks by assumption) is quasi-compact.
\end{adefn}

Using the fact that any two descriptions of $\cX$ as an Ind-algebraic
stack as above (i.e.\ as the $2$-colimit of quasi-compact algebraic
stacks with respect to transition morphisms that are closed immersions)
are mutually cofinal, one sees that the property of a morphism
being Ind-representable by algebraic stacks (resp.\ Ind-representable
by algebraic stacks and Ind-quasi-compact) is independent of
the choice of such a description of~$\cX$.

Suppose now, in the context of the Definition~\ref{def:ind-q.c.},
that $\cY$ is also an Ind-algebraic stack,
which admits a description as the $2$-colimit of algebraic stacks
with respect to transition morphisms that are closed immersions,
say $\cY \cong \varinjlim \cY_{\mu}$.
Then the induced morphism~\eqref{eqn:restricted morphism}
factors through
$\cY_{\mu}$ for some $\mu$
(since the $\cX_{\lambda}$ are quasi-compact
	and the transition morphisms between the $\cY_{\mu}$ are monomorphisms).
Since any morphism of algebraic stacks is representable by algebraic stacks,
and since closed immersions are representable by algebraic stacks,
we find that~\eqref{eqn:restricted morphism} is necessarily representable
by algebraic stacks, and 
thus the morphism $\cX \to \cY$ is necessarily Ind-representable by algebraic
stacks.  
Furthermore, 
the morphism~\eqref{eqn:restricted morphism} is quasi-compact
if and only if the induced morphism $\cX_{\lambda} \to \cY_{\mu}$
is quasi-compact for some (or equivalently any) allowable choice
of~$\mu$.
For example, if the $\cY_{\mu}$ are all quasi-separated,
then these induced morphisms are necessarily quasi-compact
(via the usual graph argument, since each $X_{\lambda}$ 
is quasi-compact), and so in this case any morphism
$\cX\to \cY$ is necessarily Ind-quasi-compact.

Suppose now that the morphism $\cX \to \cY$ {\em is} Ind-quasi-compact,
or equivalently (as we have just explained), that the various induced morphisms 
$\cX_{\lambda} \to \cY_{\mu}$ are quasi-compact.
Then each of these induced morphisms
has a scheme-theoretic image~$\cZ_\lambda$.  One easily checks that
	 $\cZ_{\lambda}$, thought of as a closed substack of $\cY$,
	 is independent of the particular choice of the index $\mu$
	 used in its definition.
	Evidently $\cZ_{\lambda}$ is a closed substack
	of $\cZ_{\lambda'}$ if $\lambda \leq \lambda'$.
        In particular, we may form the $2$-colimit
	$\varinjlim \cZ_{\lambda}$, which is an Ind-algebraic stack.
	There is a natural morphism $\varinjlim \cZ_{\lambda} \to
        \cY$.

\begin{adefn}
\label{def:scheme-theoretic image}\index{scheme-theoretic image}
In the preceding context,
we define the scheme-theoretic image of the Ind-quasi-compact morphism
        $\cX\to\cY$ to be~$\cZ~:=~\varinjlim \cZ_{\lambda}.$
\end{adefn}

It is easily verified that the scheme-theoretic image, so defined,
is independent of the chosen descriptions of $\cX$ and~$\cY$.
        There is a canonical monomorphism $\cZ \hookrightarrow \cY$.
        (We don't claim that this monomorphism is necessarily a closed immersion
        in this level of generality.)

        A very special case of this definition is the case that~$\cX \to \cY$
        is representable by algebraic spaces and quasi-compact,
        and~$\cY=\Spf A$ for a pro-Artinian ring~$A$.
In this case the definition of the
        scheme-theoretic image coincides with~\cite[Defn.\ 3.2.15]{EGstacktheoreticimages} (i.e.\ the scheme-theoretic image can be computed via the scheme-theoretic images of the
        pull-backs of~$\cX$ to the discrete Artinian quotients of~$A$).
        In this particular case the scheme-theoretic image {\em is} a
        closed formal subspace of $\cY$, i.e.\ of the form $\Spf B$
        for some topological quotient $B$ of~$A$.

        \begin{adefn}
          \label{adefn: scheme theoretically dominant}Let
          $\cX$ and $\cY$ be 
 Ind-algebraic stacks which satisfy
the hypotheses introduced above, i.e.\ which may each be written as the $2$-colimit with
respect to closed immersions of \index{scheme-theoretically dominant}
algebraic stacks, which are furthermore quasi-compact in the case 
of~$\cX$.
Then we say that an Ind-quasi-compact morphism
$\cX\to\cY$ is \emph{scheme-theoretically dominant}
          if the induced map $\cZ\to\cY$ is an isomorphism, where
          $\cZ$ is the scheme-theoretic image of $\cX\to\cY$.
        \end{adefn}
\
\begin{aremark}
If $\cX$ is a quasi-compact and quasi-separated formal algebraic stack,
then Lemma~\ref{alem: formal implies Ind} shows
that $\cX \cong \varinjlim \cX_{\lambda}$ with the $\cX_{\lambda}$ 
being quasi-compact and quasi-separated algebraic stacks, and the morphisms
being thickenings, and so in particular closed immersions.
The preceding discussion thus applies to morphisms $\cX \to \cY$
of quasi-compact and quasi-separated formal algebraic stacks, and shows that such
morphisms are necessarily Ind-representable by algebraic stacks and Ind-quasi-compact.
In particular,
Definition~\ref{adefn: scheme theoretically dominant} applies to such
morphisms, and in this case recovers the definition of scheme-theoretic dominance 
given in~\cite[Def.~6.13]{Emertonformalstacks}.

A closely related context in which the preceding definition applies if the following:
if $A$ is an adic* topological ring, with finitely generated ideal of definition~$I$,
and $\cX,\cY \to \Spf A$ are morphisms
of quasi-compact formal algebraic stacks that are representable by algebraic stacks,
then we may write $\cX \cong \varinjlim_n \cX_n$
and $\cY \cong \varinjlim_n \cY_n$,
where $\cX_n := \cX \times_{\Spf A} \Spec A/I^n$
and $\cY_n := \cY\times_{\Spf A} \Spec A/I^n$
are quasi-compact algebraic stacks.  A morphism $\cX \to \cY$ of stacks over $\Spf A$
is then necessarily representable by algebraic stacks~\cite[Lem.~7.10]{Emertonformalstacks},
and so is Ind-quasi-compact if and only if it quasi-compact in the usual sense.
We may then define its scheme-theoretic image, following Definition~\ref{def:scheme-theoretic
image},
or speak of such a morphism being scheme-theoretically dominant.
\end{aremark}

We now show (in Proposition~\ref{aprop: scheme theoretic image is p adic finite
    type II}) that under certain hypotheses the
scheme-theoretic image of a $p$-adic formal algebraic stack is also a
$p$-adic formal algebraic stack. The deduction of this result from
those of~\cite{Emertonformalstacks} will involve the following definition.  
\begin{adefn}
  \label{adefn: locally Ind finite type}(\cite[Defn.\ 8.26, Rem.\
  8.39]{Emertonformalstacks}) We say that a formal algebraic
  stack~$\cX$ over a scheme $S$ is \index{Ind-locally of finite type} \emph{Ind-locally of finite type
    over~$S$} 
  if there exists an isomorphism
$\cX \iso \varinjlim \cX_i,$
where each $\cX_i$ is an algebraic stack locally of finite type over~$S$,
and the transition morphisms are thickenings.
\end{adefn}

\begin{aremark}
\label{rem:Ind-loc. f.t. criterion}
If $\cX$ is a quasi-compact and quasi-separated formal algebraic stack,
then in order to verify that $\cX$ is Ind-locally of finite type,
it suffices to exhibit an isomorphism
$\cX \iso \varinjlim \cX_i,$
where each $\cX_i$ is an algebraic stack locally of finite type over~$S$,
and the transition morphisms are closed immersions (not necessarily
thickenings) \cite[Lem.~8.29]{Emertonformalstacks}.
\end{aremark}

\begin{aprop}\label{aprop: scheme theoretic image is p adic finite
    type II} Suppose that we have 
  a commutative diagram $$\xymatrix{\cX \ar[r] \ar[rd] & \cY \ar[d] \\
    & \Spf \Zp}$$ in which the diagonal arrow makes $\cX$ into a
  $p$-adic formal algebraic stack of finite presentation, and where~$\cY$ is 
  an Ind-algebraic stack which can be written as  the inductive limit of
        algebraic stacks, each of finite presentation
        over~$\Spec\Z/p^a$ for some~$a\geq 1$,
with the transition maps being closed immersions.

Suppose also that
  the horizontal morphism $\cX \to \cY$ {\em (}which, by the usual graph
argument, is seen to be representable by algebraic stacks{\em )} is proper. Then the
  scheme-theoretic image~$\cZ$ of this morphism 
  is a $p$-adic
  formal algebraic stack of finite type. If~$\cX$ is flat
  over~$\Spf\Zp$, then so is~$\cZ$.
\end{aprop}
\begin{proof}
  We deduce this from~\cite[Prop.\ 10.5]{Emertonformalstacks}. 
  Note firstly that~$\cX$ is quasi-compact and
  quasi-separated, since it is of finite presentation over~$\Spf\Zp$. Examining
  the hypotheses of~\cite[Prop.\ 10.5]{Emertonformalstacks}, we need
  to show that~$\cZ$ is formal algebraic,
that it is quasi-compact and quasi-separated, that it 
  is Ind-locally of finite type over~$\Spec\Zp$, 
  and that the
  induced morphism~$\cX\to\cZ$ (which will be representable by
  algebraic stacks, by the usual graph argument) is proper.

  By hypothesis, we can write~$\cY\cong\varinjlim_\mu\cY_\mu$, where
  each~$\cY_{\mu}$ is an algebraic stack
of finite presentation over some~$\Z/p^a$ 
(with~$a$ depending on~$\mu$),
and the
  transition maps are closed immersions.
If we write~$\cX\cong\varinjlim_a\cX_a$, where $\cX_a := \Spec \Z/p^a\otimes_{
\Spf \Z_p} \cX$, then by assumption each $\cX_a$ is an algebraic stack.
By definition, then, we have that
\anumequation
\label{eqn:Z description}
  \cZ:=\varinjlim_a\cZ_a,
\end{equation}
where~$\cZ_a$ is the
  scheme-theoretic image of~$\cX_a\to\cY_\mu$, for~$\mu$
  sufficiently large. 
  It follows from Lemma~\ref{lem:thick images} below
  that each of the closed immersions $\cZ_a \hookrightarrow \cZ_{a+1}$
  is a finite order thickening.  We conclude
  from~\cite[Lem.\ 6.3]{Emertonformalstacks} that $\cZ$ is in
  fact a formal algebraic stack.

  Since~$\cY_\mu$ is of finite presentation over~$\Spec\Z_p$, each
  closed substack~$\cZ_a$ is of finite presentation over $\Z/p^a$, and in
  particular quasi-compact and quasi-separated; it follows that~$\cZ$
  is quasi-compact and quasi-separated.
  The description~\eqref{eqn:Z description} of $\cZ$ furthermore exhibits 
  $\cZ$ as being Ind-locally of finite type over $\Spec \Z_p$.

It remains to show that~$\cX\to\cZ$ is 
proper.
If $T \to \cZ$ is any morphism whose source is a scheme, then since $\cZ \to \cY$
is a monomorphism, we find that there is an isomorphism
$T\times_{\cZ} \cX \to T\times_{\cY} \cX$.   Since by 
assumption  the target of this isomorphism
is an algebraic stack, proper over $T$, 
the same is true of the source.

The claim regarding flatness follows from Lemma~\ref{lem:p-adic flatness} below.
\end{proof}
We used the following lemmas in the proof of the Proposition~\ref{aprop: scheme theoretic image is p adic finite
    type II}.
\begin{alemma}
\label{lem:thick images}
Consider a commutative diagram
of morphisms of algebraic stacks
 $$\xymatrix{ \cX \ar[r]\ar[d] & \cZ \ar[d] \\
\cX' \ar[r] & \cZ'}$$
with $\cZ'$ being quasi-compact,
in which the left-hand vertical arrow is an $n$th order thickening
for some $n \geq 1$,
the right-hand vertical arrow is a closed immersion, and
the lower horizontal arrow is quasi-compact, surjective,
and scheme-theoretically dominant.
Then the right-hand vertical arrow is then also an $n$th order thickening. 
\end{alemma}
\begin{proof}
Since $\cX \to \cX'$ and $\cX' \to \cZ'$ are surjective,
by assumption, the same is true of their composite $\cX \to \cZ'$,
and hence of the closed immersion $\cZ \to \cZ'$; thus this closed
immersion is a thickening.  To see that it is of order~$n$, we first
note that since $\cZ'$ is quasi-compact, we find a smooth surjection
$U' \to \cZ'$ whose source is an affine scheme; since the property
of a thickening being of finite order may be checked {\em fppf}
locally, and since the property of being quasi-compact
and scheme-theoretically dominant
is also preserved by flat base-change, after pulling back
our diagram over~$U'$, we may assume that $\cZ' = U'$ is an
affine scheme, and that $\cZ = U$ is a closed  subscheme.
Let $\cI$ be the ideal sheaf on $U'$ cutting out $U$,
and let $a$ be a section of $\cI$.  If $a'$ denotes the pull-back of
$a$ to $\cX'$, then $a'_{| \cX} = 0,$ and so~$(a')^n~=~0$,  
by assumption.
Since $\cX' \to U'$ is scheme-theoretically dominant, we find that~$a^n~=~0$.
Thus $U\hookrightarrow U'$ is indeed an $n$th order thickening.
\end{proof}

\begin{alemma}
\label{lem:p-adic flatness}
If $\cX \to \cY$ is a quasi-compact scheme-theoretically dominant morphism 
of $p$-adic formal algebraic stacks which are locally of finite type, 
and if $\cX$ is flat over $\Z_p$, then the same is true of $\cY$.
\end{alemma}
\begin{proof}
Let $V \to \cY$ be a morphism which is representable by algebraic spaces and smooth,
whose source is an affine formal algebraic space; so $V = \Spf B$ 
for some $p$-adically complete $\Z_p$-algebra $B$ that is topologically
of finite type.  It suffices to show that $B$ is flat over $\Z_p$.
Since $V \to \cY$ is in particular flat (being smooth),
the base-change morphism $\cX\times_{\cY} V \to V$ is again scheme-theoretically dominant.
Since it is also quasi-compact, and since $V$ is formally affine (and
so  quasi-compact),
we may find a formal algebraic space $U = \Spf A$ endowed with a 
morphism $U \to \cX\times_{\cY}V $ which is representable by algebraic spaces,
smooth, and surjective; and thus also scheme-theoretically dominant.
Thus we may replace our original situation with the composite morphism
$\Spf A \to \Spf B$.  But in this context, scheme-theoretic dominance
amounts to the morphism $B \to A$ being injective; thus $B$ is $\Z_p$-flat
if $A$ is.
\end{proof}

\subsection*{Versality and versal rings}
\label{subsubsec:versal}
We will sometimes find it useful to study scheme-theoretic images in
terms of versal rings, and so we will recall some notation and
results from~\cite[\S 2.2]{EGstacktheoreticimages} related to this topic. 

If~$\Lambda$ is a Noetherian ring,
equipped with a finite ring map $\Lambda\to k$ whose target is a
field, then we let~$\cC_{\Lambda}$ denote the category whose objects
are Artinian local $\Lambda$-algebras~$A$ equipped with an
isomorphism of~$\Lambda$-algebras $A/\m_A\isoto k$. We
let~$\pro\cC_{\Lambda}$ be the corresponding category of formal
pro-objects, which (via passage to projective limits) we identify with
the category of topological pro-(discrete Artinian) local
$\Lambda$-algebras~$A$ equipped with a $\Lambda$-algebra isomorphism
$A/\m_A\isoto k$. By~\cite[Rem.\ 2.27]{EGstacktheoreticimages}, any
morphism $A\to B$ in ~$\pro\cC_{\Lambda}$ has closed image, and
induces a topological quotient map from its source to its image, so
that in particular $A\to B$ is surjective if and only if it is induced
by a compatible system of surjective morphisms in~$\cC_\Lambda$.

Fix a locally Noetherian base scheme~$S$,
and let $k$ be a finite type $\cO_S$-field,
i.e.\ $k$ is a field equipped with a morphism $\Spec k\to S$ of
finite type.   
Choose, as we may,
an affine open subscheme $\Spec\Lambda\subseteq S$ for which $\Spec k\to S$
factors as $\Spec k\to\Spec\Lambda\to S$, with $\Lambda\to k$ being finite.
(In what follows we fix such a choice of $\Lambda$, although the notions
that we define in terms of it are independent of this choice.) We now
define a category fibred in groupoids $\widehat{\cF}_x$ in the
following way; this category is an example of a
deformation category in the sense
of~\cite[\href{https://stacks.math.columbia.edu/tag/06J9}{Tag
  06J9}]{stacks-project}. In the
notation of~\cite[\href{https://stacks.math.columbia.edu/tag/07T2}{Tag
  07T2}]{stacks-project}, our category~$\widehat{\cF}_x$ is
denoted~$\cF_{\cF,k,x}$ (with a slightly unfortunate clash of notation
in the two instances of~$\cF$).

\begin{adefn}\label{adefn: deformation category}
  If $\cF$ is a category fibred in groupoids over~$S$, and if
  $x:\Spec k \to \cF$ is a morphism,
  then we define a category~$\widehat{\cF}_{x}$, cofibred in groupoids
  over~$\cC_{\Lambda}$, as
  follows: 
  For any object $A$ of $\cC_{\Lambda}$, the objects of
  $\widehat{\cF}_x(A)$ consist of morphisms $y:\Spec A \to \cF$,
  together with an isomorphism $\alpha: x \iso \ybar$ compatible with
  the given identification of $A/\mathfrak m$ with $k$; here $\ybar$
  denotes the induced morphism $\Spec A/\mathfrak m \to \cF$.  The set
  of morphisms between two objects $(y,\alpha)$ and $(y',\alpha')$ of
  $\widehat{\cF}_x(A)$ consists of the subset of morphisms
  $\beta: y\to y'$ in $\cF(A)$ for which
  $\alpha'\circ \overline{\beta} = \alpha$; here $\overline{\beta}$
  denotes the morphism $\ybar\to \ybar'$ induced by~$\beta$.  If
  $A \to B$ is a morphism in $\cC_{\Lambda}$, then the corresponding
  pushforward $\widehat{\cF}_x(A) \to \widehat{\cF}_x(B)$ is defined
  by pulling back morphisms to $\cF$ along the corresponding morphism
  of schemes $\Spec B \to \Spec A$.
\end{adefn}
%


If $\cF$ is a category fibred in groupoids over the locally Noetherian scheme $S$,
and if $x: \Spec k \to \cF$ is a $k$-valued point of~$\cF$,
for some finite type $\cO_S$-field,
then the notion of a versal ring to $\cF$ at $x$ 
is defined in~\cite[Def.~2.2.9]{EGstacktheoreticimages}.
Rather than recalling that definition here,
we will give a definition of versality in a greater level of generality that is convenient
for us.  We will then explain how the notion of versal ring is obtained as a particular
case.

\begin{adefn}
\label{def:versal}\index{versal}
Let $f:\cF \to \cG$ be a morphism of categories fibred in groupoids over the
locally Noetherian scheme~$S$, let $k$ be a finite type $\cO_S$-field,
let $x: \Spec k \to \cF$ be a $k$-valued point of~$\cF$,
and let $y: \Spec k \to \cG$
be a $k$-valued point of~$\cG$ equipped with an isomorphism of $k$-valued points
$\alpha: f\circ x \iso y$.
The morphism $f$ induces in an evident way a morphism 
$\widehat{f}_x: \widehat{\cF}_x \to \widehat{\cG}_y$
of categories cofibred in groupoids.
We say that $f$ is {\em versal} at $x$ 
if the the morphism $\widehat{f}_x$ is {\em smooth}
in the sense 
of~\cite[\href{https://stacks.math.columbia.edu/tag/06HG}{Tag 06HG}]{stacks-project};
that is, given a commutative diagram
$$\xymatrix{\Spec B \ar[r] \ar[d] & \Spec A \ar[d] \ar@{-->}[ld]\\
\widehat{\cF}_x \ar^-{\widehat{f}_x}[r] & \widehat{\cG}_{x}}$$
in which the upper arrow is the closed immersion corresponding
to a surjection $A \to B$ in $\cC_{\Lambda}$, 
 we can fill in the
dotted arrow 
so that the diagram remains commutative.
(Clearly this notion is independent of the particular choice of $y$ and $\alpha$;
more precisely
it holds for any such choice if it holds for one such choice, such as $y = f\circ x$
and $\alpha = \id$.)
\end{adefn}

\begin{aexample}
Let $S$ be a locally Noetherian scheme,
let $k$ be a finite type $\cO_S$-field~$k$,
and let $A$ be an object of $\pro\cC_{\Lambda}$;
recall in particular then that $A$ comes equipped with a chosen isomorphism
$\Spec A/\mathfrak m \iso k.$
We let $x': \Spec k \iso \Spec A/\mathfrak m \hookrightarrow \Spf A$
denote the induced $k$-valued point of $\Spf A$.

Now let $\cF$ be a category fibred in groupoids over~$S$,
let $x:\Spec k \to \cF$ be a $k$-valued point of $\cF$,
and suppose that $f:\Spf A \to \cF$ is a morphism, for which we can
find an isomorphism $f\circ x' \iso x$ of $k$-valued points of~$\cF$.
Then the morphism $f$ is versal at the point $x'$ if and only
if $A$ is a versal ring to the $k$-valued point $x$ of $\cF$
in the sense of~\cite[Def.~2.2.9]{EGstacktheoreticimages}
\end{aexample}

\begin{aexample}
Let $S$ be a locally Noetherian scheme, let $U$ be a locally finite type $S$-scheme,
and let $f:U \to \cF$ be a  morphism to a category fibred in groupoids over~$S$.
If $u\in U$ is a finite type point, giving rise to a morphism
$\Spec \kappa(u) \to U$, then $f$ is versal at this $\kappa(u)$-valued
point of $U$ if and only if $f$ is versal at the point $u$
in the sense of~\cite[Def.~2.4.4~(1)]{EGstacktheoreticimages};
cf.~\cite[Rem.~2.4.5]{EGstacktheoreticimages}.
\end{aexample}



\begin{arem}\label{arem: versal rings for formal Ind}
  If~$\cX$ is an algebraic stack which is locally of finite
  presentation over a locally Noetherian scheme~$S$, then it admits
  (effective, Noetherian) versal rings at all finite type points~\cite[\href{https://stacks.math.columbia.edu/tag/0DR1}{Tag 0DR1}]{stacks-project}. If
  $X$ is an Ind-locally finite type algebraic space over~$S$, then it
  admits (canonical) versal rings at all finite type points,
  by~\cite[Lem.\ 4.2.14]{EGstacktheoreticimages}.
%
  We don't prove a general statement about the existence of versal
  rings for Ind-algebraic stacks, since in all the cases we consider
  in the body of the book we are able to construct them ``explicitly''
   (for example, in terms of Galois lifting rings, or lifting rings for
  \'etale $\varphi$-modules). 
\end{arem}


The following lemma and its proof are essentially \cite[Lem.\
3.2.16]{EGstacktheoreticimages}, but since the setup there is
different, we give the details here.  Before stating the lemma, 
we introduce the setup. Let~$S$ be a locally Noetherian scheme.  We suppose given a morphism
  $\cX \to \cY$ of stacks over~$S$ which is representable by algebraic
stacks and proper, 
and that $\cY$ is an Ind-algebraic stack which may
be written as the $2$-colimit of 
algebraic stacks which are of finite presentation over $S$ (and so in particular
quasi-compact and quasi-separated) with respect 
to transition morphisms that are closed immersions,  
say $\cY \cong \varinjlim \cY_{\lambda}.$
Then $\cX_{\lambda} := \cX \times_{\cY} \cY_{\lambda}$ is an algebraic stack,
and the projection $\cX_{\lambda} \to \cY_{\lambda}$ is proper, so that 
$\cX_{\lambda}$ is finite type over $S$ (and in particular quasi-compact,
and also quasi-separated, although we won't
use this latter fact).   Furthermore, we have an induced isomorphism
$\cX \cong \varinjlim_{\lambda} \cX_{\lambda}$, and so $\cX \to \cY$ 
is a morphism of Ind-algebraic stacks whose scheme-theoretic image~ $\cZ$ may be
defined.  Of course, in this context, since the Ind-structures on $\cX$
and $\cY$ are compatible, if we let $\cZ_{\lambda}$ denote the scheme-theoretic
image of $\cX_{\lambda}$ in $\cY_{\lambda}$, then this coincides with the
scheme-theoretic image of $\cX_{\lambda}$ in $\cY_{\lambda'}$ for any $\lambda' \geq
\lambda$, and we may write $\cZ \cong \varinjlim_{\lambda}
\cZ_{\lambda}$. Note that~$\cZ_\lambda$ is also of finite
presentation over~$S$, by Lemma~\ref{lem:finite presentation from closed immersion}.

\begin{alemma}
  \label{alem: scheme theoretic image of versal is versal}
Suppose that we are in the preceding situation, so that   $\cX \to
\cY$ is a morphism of stacks over a locally Noetherian base~$S$ which is representable by algebraic
stacks and proper, where $\cY$ is an Ind-algebraic stack which may
be written as the $2$-colimit of 
algebraic stacks which are of finite presentation over~$S$, with respect 
to transition morphisms that are closed immersions, and we write~$\cZ$
for the scheme-theoretic image of $\cX\to\cY$.

Suppose that~$x:\Spec k\to
  \cZ$ is a finite type point, and that $\Spf A_x\to\cY$ is a versal
  morphism for the composite $x:\Spec k\to
  \cZ\into\cY$.
Let~$\Spf B_x$ be the scheme-theoretic image of $\cX_{\Spf A_x}\to\Spf
A_x$. Then the morphism $\Spf B_x\to\cY$ factors through a versal
morphism $\Spf B_x\to \cZ$.
\end{alemma}
\begin{proof}
  By definition, we may write~ $A_x=\varprojlim A_i$,
  $B_x=\varprojlim B_i$, where the~$A_i$ are objects of~$\cC_\Lambda$,
  and $\Spec B_i$ is the scheme-theoretic image of
  $\cX_{A_i}\to \Spec A_i$. As explained immediately above, we write
  $\cX \cong \varinjlim \cX_{\lambda}$, $\cY\cong\varinjlim\cY_{\lambda}$, and
  $\cZ\cong\varinjlim \cZ_{\lambda}$, where the transition morphisms are
  closed immersions
  of algebraic stacks, and~$\cZ_\lambda$ is the
  scheme-theoretic image of the morphism
  $\cX_{\lambda}\to\cY_{\lambda}$; in particular, the morphism $\cX_\lambda\to\cZ_\lambda$ is proper and scheme-theoretically
  dominant. 

  If follows that for each~$i$, the composite $\Spec B_i\to\Spec A_i\to\cY$
  factors through~$\cZ$. Indeed, for~$\lambda$ sufficiently large the
  morphism $(\cX_\lambda)_{B_i}\to\Spec B_i$ is scheme-theoretically
  surjective, and  the morphism $(\cX_\lambda)_{B_i}\to\cY$ factors
  through $\cZ_\lambda$, so the morphism $\Spec B_i\to\cY$ also
  factors through~$\cZ_\lambda$.

We now show that a morphism $\Spec A\to\Spf A_x$, with $A$ an object of
$\cC_{\Lambda}$, 
factors through $\Spf B_x$ if and only if the composite $\Spec
A\to\Spf A_x\to\cY$ factors through~$\cZ$. In
one direction, if
$\Spec A\to\cY$ factors through~$\Spf B_x$, then it factors
through $\Spec B_i$ for some $i$, and hence 
through $\cZ$, as we saw above.

Conversely, if the composite $\Spec
A\to\Spf A_x\to\cY$ factors through~$\cZ$, then we claim 
that there exists factorisation $\Spec A\to\Spec B\to\cY$, where $B$ is
an object of~$\cC_{\Lambda}$, the morphism $\Spec
A\to\Spec B$ is a closed immersion, and $\cX_B\to \Spec B$ is
scheme-theoretically dominant. To see this, 
note firstly that the composite $\Spec A\to\cY$ factors
through~$\cZ_\lambda$ for some~$\lambda$. Since~$\cZ_\lambda$ is
an algebraic stack which is locally of finite presentation over a
locally Noetherian base, 
it admits effective Noetherian versal rings by Remark~\ref{arem: versal rings for formal Ind}, so
that $\Spec A\to\cZ_\lambda$  factors through a versal morphism $\Spec C_x\to\cZ_\lambda$ at the
finite type point of~$\cZ_\lambda$ induced by~$x$.  By ~\cite[Lem.\
1.6.3]{EGstacktheoreticimages}, we can find a factorisation  $\Spec
A\to\Spec R_x\to\Spec C_x$, where~$R_x$ is complete local Noetherian,
$\Spec R_x\to\Spec C_x$ is faithfully flat, and $\Spec A\to\Spec R_x$ is a closed
immersion. Since $\Spec R_x\to\Spec C_x$ is faithfully flat, and $\Spec
C_x\to\cZ_\lambda$ is flat
by~\cite[\href{https://stacks.math.columbia.edu/tag/0DR2}{Tag
  0DR2}]{stacks-project}, we see that the base changed morphism
$(\cX_\lambda)_{R_x}\to\Spec R_x$ is scheme-theoretically
dominant.  By~\cite[Lem.\ 3.2.4]{EGstacktheoreticimages} (applied to the proper
morphism of algebraic stacks~$\cX_\lambda\to\cZ_\lambda$), $R_x$
admits a cofinal collection of Artinian quotients~$R_i$ for which
$(\cX_\lambda)_{R_i}\to\Spec R_i$ is scheme-theoretically
dominant. The closed immersion $\Spec A\to\Spec R_x$ factors through
~$\Spec R_i$ for some~$R_i$, so we may take~$B=R_i$.

Now, by the versality of $\Spf A_x\to\cY$, we may lift the morphism
$\Spec B\to\cY$ to a morphism $\Spec B\to\Spf A_x$, which furthermore
we may factor as $\Spec B \to \Spec A_i \to \Spf A_x$, for some value
of $i$.  Since $\cX_B \to \Spec B$ is scheme-theoretically dominant,
the morphism $\Spec B \to \Spec A_i$ then factors through $\Spec B_i$,
and thus through $\Spf B_x$, as required.

It follows in particular that the composite
$\Spf B_x\to\Spf A_x\to\cY$ factors through a morphism
$\Spf B_x\to\cZ$. It remains to check that this morphism is
versal. This is formal. Suppose we are given a commutative
diagram \[\xymatrix{\Spec A_0\ar[r]\ar[d]&\Spf B_x\ar[r]\ar[d]&\Spf A_x\ar[d]\\
    \Spec A\ar[r]&\cZ\ar[r]&\cY }\]where the left hand vertical arrow
is a closed immersion, and $A_0, A$ are objects
of~
$\cC_\Lambda$. By the versality of $\Spf A_x\to\cY$, we may lift the
composite $\Spec A\to\cY$ to a morphism $\Spec A\to\Spf A_x$. Since
the composite $\Spec A\to\Spf A_x\to\cY$ factors through
$\cZ$, the morphism $\Spec A\to \Spf A_x$ then factors through $\Spf
B_x$, as required.
\end{proof}

The following lemma records a useful property of versal rings in the Noetherian
context.

\begin{alem}
\label{lem:versal flatness}
Let $\cX$ is a locally Noetherian formal algebraic stack, let $R$ be a complete
Noetherian local ring with residue field~$k$,
and let
$f:\Spf R \to \cX$ be a morphism for which the induced morphism
$x:\Spec k \to \cX_{\red}$ is a finite type point, and which
is versal to $\cX$ at $x$.   Then the morphism $f$ is flat,
in the sense of~{\em \cite[Def.~8.42]{Emertonformalstacks}}.  
\end{alem}
\begin{proof}
Since $\cX$ is locally Noetherian, by assumption,
we may find a morphism $\Spf A \to \cX$ whose source
is a Noetherian affine formal algebraic space, which is representable 
by algebraic spaces and smooth, and whose image contains the point~$x$.
By definition, we have to show that the base-changed morphism
$\Spf A \times_{\cX} \Spf R \to \Spf A$ is flat.  

The fibre product $\Spf A \times_{\cX} \Spf R$ is a formal algebraic space
(which is locally Noetherian, since $\Spf R$ is
so~\cite[\href{https://stacks.math.columbia.edu/tag/0AQ7}{Tag 0AQ7}]{stacks-project}),
and so admits a morphism $\coprod \Spf B_i \to \Spf A \times_{\cX} \Spf R$
whose source is the disjoint union of Noetherian affine formal algebraic spaces,
and which is representable by algebraic spaces, \'etale, and surjective.
Again, it suffices to show that each of the morphisms $\Spf B_i \to \Spf A$
is flat, which by definition is equivalent to each of the induced morphisms $A \to B_i$
being flat. 
If we let $I$ denote an ideal of definition of the topology on $A$,
then it suffices to show that each morphism $A/I^n \to B_i/I^n$ is
flat (see e.g.\ the proof of~\cite[Lem.\ 8.18]{Emertonformalstacks}).

For this, it suffices in turn to show, for each maximal ideal $\mathfrak n$
of $B_i$,
that the induced 
morphism $A/I^n \to
(B_i/I^n)\, \hat{}$
is flat (where $(\text{--})\, \hat{}$ denotes
$\mathfrak n$-adic completion; note that if $\mathfrak m$ denotes the maximal ideal of~$R$,
then the topology on $B_i$ is the $\mathfrak m$-adic topology,
so that any maximal ideal $\mathfrak n$ of $B_i$ contains~$\mathfrak m$,
and also contains~$IB$). 
Now the induced morphism $\Spf \widehat{B}_i \to \Spf A$ is versal 
to the induced morphism $\Spec B_i/\mathfrak n \to \Spf A$ 
(as one sees by chasing through
the constructions, beginning from the versality of $\Spf R \to \cX$,
and taking into account the infinitesimal lifting property for \'etale
morphisms),
and thus the induced morphism
$\Spf (B_i/I^n)\, \hat{} \to \Spec A/I^n$ is also versal
(to the induced morphism
$\Spec B_i/\mathfrak n \to \Spec A/I^n$).
Thus this morphism {\em is} flat, e.g.\
by~\cite[\href{https://stacks.math.columbia.edu/tag/0DR2}{Tag 0DR2}]{stacks-project},
and the lemma is proved.
\end{proof}

We will frequently find the following lemma useful.
\begin{alem}\label{lem: criterion for Artin to map to scheme theoretic
    image}
	Let $R \to S$ be a continuous surjection of objects in
	$\pro\cC_{\Lambda}$,
	and let $X\to \Spf R$ be a finite type morphism
	of formal algebraic spaces. 
        
	Make the following assumption: 
	if $A$ is any finite-type Artinian local $R$-algebra
	for which the canonical morphism $R\to A$ factors through 
	a discrete quotient of $R$,
	and for which the canonical morphism $X_A \to \Spec A$
	admits a section,
	then the canonical morphism $R \to A$ furthermore factors
        through $S$.

        Then the scheme-theoretic image of~$X\to\Spf R$ is a closed
        formal subscheme of~$\Spf S$.
\end{alem}
\begin{proof}Writing~$\Spf T$ for the scheme-theoretic image of $X\to
  \Spf R$, we need to show that the surjection $R\to T$ factors
  through~$S$. By definition, we can
write~$T$ as an inverse limit of discrete Artinian quotients~$B$
for which the canonical morphism $X_B \to \Spec B$ is
scheme-theoretically dominant. It follows from \cite[Lem.\
5.4.15]{EGstacktheoreticimages} that for any such~$B$, the surjection
~$R\to B$ 
factors through~$S$; so the surjection $R\to T$ factors
through~$S$, as required.
\end{proof}




 \subsection*{Immersed substacks} 

The following lemma is often useful for studying morphisms from a
stack in terms of morphisms from a cover of the stack.
	\begin{alemma}\cite[Lem.\ 3.18, 3.19]{Emertonformalstacks}
		\label{lem:stacks as quotients}
		Suppose that $\cX$ is a stack over a base scheme~$S$,
		and that $U\to \cX$ is a morphism over~$S$ which is
                representable by algebraic spaces and whose source is
                a sheaf. Write $R := U\times_{\cX} U$.

                \begin{enumerate}
                \item Suppose that the projections
                  $R := U\times_{\cX} U \rightrightarrows U$ {\em
                    (}which are again representable by algebraic
                  spaces, being the base-change of morphisms which are
                  so representable{\em )} are flat and locally of
                  finite presentation.  Then the morphism $U\to \cX$
                  induces a monomorphism $[U/R]\to \cX$.

                \item  Suppose that $U\to \cX$ is flat, surjective, and
                  locally of finite presentation.  Then the morphism
                  $U\to \cX$ induces an isomorphism $[U/R]\iso \cX$.
                \end{enumerate}

	\end{alemma}

In particular, if $\cX$ is a formal algebraic stack (over some base scheme $S$),
 then by definition there exists a morphism $U\to \cX$ 
 which is representable by algebraic spaces, smooth, and surjective,
 and whose source~$U$ is a formal algebraic space.
 In this case $R := U\times_{\cX}U$ is a formal algebraic
 space which is endowed in a natural way with the structure 
 of a groupoid in formal algebraic spaces over $U$, and (since smooth
 morphisms are in particular flat), by Lemma~\ref{lem:stacks as quotients} there is an isomorphism of stacks $[U/R] \iso \cX$. 
 
 Suppose now that $\cZ \hookrightarrow \cX$ is an immersion (in
the sense that it is representable by algebraic spaces,
and pulls back to an immersion over any test morphism $T \to \cX$
whose source is a scheme). 
 The
 induced morphism
 $W := U\times_{\cX} \cZ \hookrightarrow U$ is then an immersion of formal
 algebraic spaces, and $W$ is $R$-invariant, in (an evident generalization
 of) the sense 
 of~\cite[\href{http://stacks.math.columbia.edu/tag/044F}{Tag
  044F}]{stacks-project}. 
 Conversely, if $W$ is an $R$-invariant locally closed formal 
 algebraic subspace of $U$, and if we write $R_W := R\times_U W =
 W \times_U R$ for the restriction of $R$ to $W$, then
 the induced morphism $[W/R_W] \to [U/R] \iso \cX$ is an immersion.
 (In the context of algebraic stacks, this 
 is~\cite[\href{http://stacks.math.columbia.edu/tag/04YN}{Tag
  04YN}]{stacks-project}.)

\subsection*{Flat parts}Let~$\cO$ be the ring of integers in a
finite extension of~$\Qp$, and let~$\cX\to\Spf\cO$ be a $p$-adic formal
algebraic stack which is locally of finite type over~$\Spf\cO$. Then
by~\cite[Ex.\ 9.11]{Emertonformalstacks}, there is a closed substack~
$\cX_{\fl}$ of~$\cX$, the $\emph{flat part}$ of~$\cX$, which is the maximal substack of~$\cX$ which
flat over~$\Spf\cO$.  \index{$p$-adic formal algebraic stack!flat part}

More precisely, $\cX_{\fl}\to\Spf\cO$ is flat, 
and if
$\cY\to\cX$ is a morphism of locally Noetherian formal algebraic
stacks for which the composite $\cY\to\Spf\cO$ is flat, then
$\cY\to\cX$ factors through~$\cX_{\fl}$.

\subsection*{Obstruction theory}Let~$\cX$ be a limit preserving Ind-algebraic stack
over a locally Noetherian scheme $S$. If $x:\Spec A \to \cX$ is a
morphism for which the composite morphism $\Spec A \to S$ factors
through an affine open subscheme of~$S$, then
in~\cite[\href{http://stacks.math.columbia.edu/tag/07Y9}{Tag
  07Y9}]{stacks-project} there is defined a functor $T_x$ from the
category of $A$-modules to itself, whose formation is also functorial
in the pair
$(x,A)$~\cite[\href{http://stacks.math.columbia.edu/tag/07YA}{Tag
  07YA}]{stacks-project}, and such that for any finitely generated $A$-module $M$, there
is a natural identification of $T_x(M)$ with the set of lifts of $x$
to morphisms $x': \Spec A[M] \to \cX$ (where $A[M]$ denotes the square
zero extension of $A$ by $M$). (These definitions apply to~$\cX$
by~\cite[Lem.\ 4.22]{Emertonformalstacks}.)

We make the following definition, which is a special case
of~\cite[Defn.\ 11.6]{Emertonformalstacks}; that definition
incorporates an auxiliary module in its definition for technical
reasons, but in our applications this module is zero, so we have
suppressed it here. We have also incorporated~\cite[Rem.\
11.9]{Emertonformalstacks}, which allows us to restrict to the case of
finitely generated $A$-modules~$M$.

\begin{adefn}
	\label{df:obstruction}
We say that $\cX$
admits a {\em nice obstruction theory} if, \index{nice obstruction theory}
for each finite type $S$-algebra $A$
which lies over an affine open subscheme of $S$,
equipped with a morphism $x: \Spec A \to \cX$, there 
exists a complex of $A$-modules $K^{\bullet}_{(x,A)}$
such that the following conditions are satisfied:

\begin{enumerate}
	\item The complex $K^{\bullet}_{(x,A)}$ is bounded above and has 
		finitely generated cohomology modules (or, equivalently,
		is isomorphic in the derived category to a bounded above complex
		of finitely generated $A$-modules).
	\item  The formation of $K^{\bullet}_{(x,A)}$ is compatible (in
		the derived category) with pull-back. More precisely,
if $f:\Spec B \to \Spec A$, inducing the morphism $y: \Spec B
\to \Spec A \to \cX$,
then there is a natural isomorphism in the derived category
$f^*K^{\bullet}_{(x,A)} \iso K^{\bullet}_{(y,B)}$.
	\item
For any pair $(x,A)$ as above,
and for any finitely generated $A$-module $M$,
we have an isomorphism $T_x(M) \isoto H^1(K^{\bullet}_{(x,A)} \Lotimes M)$
whose formation is functorial in $M$ and in the pair $(x,A)$.

\item
	For any pair $(x,A)$ as above,
	and for any finitely generated $A$-module $M$,
	the cohomology module
$H^2(K^{\bullet}_{(x,A)}\Lotimes M)$ serves as an obstruction
module. In particular, for any square zero extension
	$$0 \to I \to A' \to A \to 0$$
	for which $A'$ is of finite type over $S$ (or, equivalently,
	for which $I$ is a finite $A$-module), 
we have a functorial obstruction element $o_x(A')
\in H^2(K^{\bullet}_{(x,A)} \Lotimes I)$, which vanishes if and only
if~$x$ can be lifted to~$\Spec A'$.
\end{enumerate}
\end{adefn}
The following result generalises the familiar fact that a
pro-representing object for a formal deformation ring whose tangent
space is finite dimensional is necessarily a complete Noetherian local
ring. 
\begin{athm}
	\label{thm:noetherian criterion}\emph{(}\cite[Thm.\ 11.13]{Emertonformalstacks}\emph{)}
	If $\cX$ is a 
locally countably indexed and Ind-locally of finite
        type formal algebraic stack over a locally Noetherian scheme
        $S$, and if~$\cX$ admits a nice obstruction theory, then $\cX$ is
        locally Noetherian.
\end{athm}

(When comparing this result with the statement of 
	\cite[Thm.\ 11.13]{Emertonformalstacks},
the reader should
bear in mind that an Ind-locally of finite type formal algebraic stack 
over~$S$ is in particular the $2$-colimit of algebraic stacks that are locally of 
finite type (or equivalently, locally of finite presentation) over $S$,
and hence is limit preserving.  Also, 
	\cite[Lem.\ 8.38]{Emertonformalstacks} ensures that an Ind-locally of
finite type formal algebraic stack over $S$ is also locally of Ind-locally finite type ---
which is the condition that appears in~\cite[Thm.\ 11.13]{Emertonformalstacks}.)

\chapter{Graded modules and rigid analysis}
\label{sec: rigid analytic perspective}
In this appendix
we establish some (mostly straightforward)
results that combine graded ring techniques 
with various results related to completions that are of a rigid analytic flavour.
\section{Associated graded algebras and modules}\label{subsec:
  associated graded}
Let $R$ be an Artinian local ring, with maximal ideal $I$ and residue field~$k$.

\begin{adefn} If $M$ is an $R$-module, then we let $\Gr^{\bullet} M$
	denote the graded $k$-vector space associated to the $I$-adic filtration
	on $M$ (so $\Gr^i M := I^iM/I^{i+1} M$).  
\end{adefn}

The formation of $\Gr^{\bullet} M$ is functorial in $M$.
If $M$ and $N$ are two $R$-modules, then there is a natural surjection
(of graded $k$-vector spaces)
\anumequation
\label{eqn:graded tensor}
\Gr^{\bullet} M \otimes_k \Gr^{\bullet} N \to \Gr^{\bullet}(M\otimes_R N).
\end{equation}
In particular, it follows that if $A$ is an $R$-algebra, then $\Gr^{\bullet} A$
is naturally a graded $k$-algebra, and if $M$ is an $A$-module,
then $\Gr^{\bullet} M$ is naturally a $\Gr^{\bullet} A$-module.

We recall the following basic lemmas.

\begin{alemma}
	\label{lem:graded Noetherian}
	If $A$ is an $R$-algebra, then $A$ is Noetherian
	if and only if $\Gr^{\bullet} A$ is Noetherian.
\end{alemma}
\begin{proof}
	This is easy and standard; see e.g.~\cite[Prop.~1.1]{MR1990669} 
	for a statement of the ``if'' direction of this result
	(which is the less obvious of the two directions)
	in a significantly more general setting.
\end{proof}

\begin{alemma}
	\label{lem:graded surjective}
	If $A$ is an $R$-algebra and if $M$ is an $A$-module,
	then the natural morphism $\Gr^{\bullet} A \otimes_{\Gr^0 A}
	\Gr^0 M \to \Gr^{\bullet} M$ is surjective.
	Furthermore, $M$ is non-zero if and only if $\Gr^0 M$ is non-zero,
	if and only if $\Gr^{\bullet} M$ is non-zero.
\end{alemma}
\begin{proof}
	The first claim is immediate.
	This implies in turn that if $\Gr^0 M$ is zero then 
	$\Gr^{\bullet} M$ is zero (the converse being evident), 
	and thus that $M$ is zero (since $\Gr^{\bullet} M$ is
	the associated graded of $M$ with respect to a finite length
	filtration beginning at $M$ and ending at $0$).
	(Essentially equivalently, one can observe directly that since $I$
	is nilpotent, if $M/IM = 0$ then $M = 0$.)
\end{proof}

\begin{alemma}
	\label{lem:graded fg}
	If $A$ is an $R$-algebra and $M$ is an $A$-module,
	then $M$ is finitely generated over $A$ if and only 
	if $\Gr^0 M$ is finitely generated over $\Gr^0 A$,
	if and only if $\Gr^{\bullet} M$ is finitely generated
	over $\Gr^{\bullet} A$.
\end{alemma}
\begin{proof} 
Clearly if $M$ is finitely generated over $A$, then
$\Gr^0 M := M/IM$ is finitely generated over $A/I$; the converse
assertion follows from Nakayama's lemma with respect to
the nilpotent ideal $I$.  If $\Gr^0 M$ is finitely generated 
over $\Gr^0 A$, then Lemma~\ref{lem:graded surjective}
shows that $\Gr^{\bullet} M$ is finitely generated over $\Gr^{\bullet} A$;
the converse assertion follows from the fact that (again by
Lemma~\ref{lem:graded surjective}) $\Gr^0 M$ may be obtained as the
quotient of $\Gr^{\bullet} M$ by the ideal $\Gr^{\bullet > 0} A$
of $\Gr^{\bullet} A$.
\end{proof}
\begin{alemma}
  \label{lem: graded iso is iso}Let $f:M\to N$ be a morphism of $A$-modules. Then~$f$ is an isomorphism if and only if the
  induced morphism $\Gr^{\bullet}f:\Gr^{\bullet}M\to\Gr^{\bullet}N$ is an
  isomorphism.
\end{alemma}
\begin{proof}
  If~$f$ is an isomorphism then certainly~$\Gr^{\bullet}f$ is an
  isomorphism. Conversely, if~$\Gr^{\bullet}f$ is an isomorphism, then
  in particular $M/IM\to N/IN$ is surjective, so~$f$ is surjective by
  Nakayama's lemma with respect to
the nilpotent ideal $I$. In addition, for each~$i$ the
  morphism~$I^iM/I^{i+1}M\to I^iN/I^{i+1}N$ is injective, so by an
  easy induction on~$i$ we see that any element of~$\ker f$ is
  contained in~$I^iM$ for all~$i$, and is thus zero, as required.
\end{proof}
\begin{alemma}
	\label{lem:graded flat}
If $A$ is an $R$-algebra, then an $A$-module $M$
is {\em (}faithfully{\em )} flat over $A$
if and only if $\Gr^{\bullet} M$ 
is {\em (}faithfully{\em )} flat over $\Gr^{\bullet} A.$
Furthermore, if any of these conditions holds,
then the natural morphism
$\Gr^{\bullet} A \otimes_{\Gr^0 A} \Gr^0 M \to
\Gr^{\bullet} M$
is an isomorphism.
\end{alemma}
\begin{proof}
This is~\cite[Lem.\ 5.5.37]{EGstacktheoreticimages}.  
\end{proof}

\begin{alemma}
	\label{lem:graded tensor}
	If $A$ is an $R$-algebra, and if $M$ and $N$ are $A$-modules,
	then the natural surjection~{\em (\ref{eqn:graded tensor})}
	induces a natural surjection 
	$$\Gr^{\bullet} M \otimes_{\Gr^{\bullet} A} \Gr^{\bullet} N \to
	\Gr^{\bullet} (M\otimes_A N).$$
\end{alemma}
\begin{proof}
	The natural surjection $M\otimes_R N \to M\otimes_A N$ induces
	a surjection $\Gr^{\bullet}(M\otimes_R N) \to \Gr^{\bullet} (M\otimes_A N).$  It is immediate that its composite with the surjection~(\ref{eqn:graded
		tensor}) factors through $\Gr^{\bullet} M \otimes_{\Gr^{\bullet}
		A} \Gr^{\bullet} N.$
\end{proof}

\begin{alemma}
	\label{lem:graded flat bis}
	If $A$ is an $R$-algebra, and if $M$ and $N$ are $A$-modules 
	with either $M$ or $N$ flat over $A$,
	then the surjection of Lemma~{\em \ref{lem:graded tensor}}
       	is a natural isomorphism
	$$\Gr^{\bullet} M \otimes_{\Gr^{\bullet} A} \Gr^{\bullet} N \iso
	\Gr^{\bullet} (M\otimes_A N).$$
\end{alemma}
\begin{proof}
	Suppose (without loss of generality) that $N$ is $A$-flat.
	Taking into account Lemma~\ref{lem:graded flat},
	we see that we have to show that natural morphism
	$\Gr^{\bullet} M \otimes_{\Gr^0 A} \Gr^0 N \to \Gr^{\bullet}(M\otimes_A N)$ is an isomorphism, i.e.\ that for each~$i \geq 0,$ the natural morphism
	$$I^iM /I^{i+1} M \otimes_{A/I} N/I \to I^i(M\otimes N)/I^{i+1}(M\otimes N)$$
	is an isomorphism.  This follows easily by tensoring $N$ over $A$
	with the various short exact sequences
$$0 \to I^n M \to I^m M \to I^mM/I^nM\to 0,$$
taking into account the flatness of $N$ over $A$.
\end{proof}

\section{A general setting}
\label{subsec:setting}
We fix an Artinian local ring $R$ with maximal ideal $I$
and residue field $k$, 
as well an $R$-algebra $C^+$, and an element $u \in C^+$, satisfying the
following properties:
\begin{enumerate}[label=(\Alph*)]
	\item $u$ is a regular element (i.e.\ a non-zero divisor) of $C^+$.
	\item $C^+/u$ is a flat $R$-algebra.
	\item $C^+/I$ is a rank one complete valuation ring, and
		the image of $u$ in $C^+/I$ (which
		is necessarily non-zero, by (A)) is of
		positive valuation (i.e.\ lies in the maximal
		ideal of $C^+/I$).
\end{enumerate}

\begin{aremark}
	\label{rem:independence of u}
	In the preceding context, if $v \in C^+$ has non-zero
	image in $C^+/I$, then
	$v^n$ divides $u^m$ in $C^+$, for some $m, n > 0$.
	(Indeed, since $C^+/I$ is a rank one valuation ring
	in which the image of $u$ has positive valuation,
	we may write $u^a = v w + x,$ for some $a>0$, some $w \in C^+$,
	and some $x \in I C^+$.  Raising both
	sides of this equation to a sufficiently large power,
	remembering that $I$ is nilpotent,
	gives the claim.)
	
	In particular, if $v \in C^+$ is another element
	satisfying conditions~(A), (B), and (C) above, then
	the powers of $u$ and of $v$ are mutually cofinal with respect
	to divisibility (i.e.\ the $u$-adic and $v$-adic topologies
        on $C^+$ coincide), and consequently $C^+[1/u] = C^+[1/v]$.
\end{aremark}

\begin{aremark}
	\label{rem:flat implies free}
	Recall  that since $R$ is an Artinian local ring,
	an $R$-module is flat if and only if it is free~\cite[\href{https://stacks.math.columbia.edu/tag/051G}{Tag 051G}]{stacks-project}.
	(We will use this fact in the proof of Lemma~\ref{lem:vanishing
		locus} below.  For a proof, note first that if $M$ 
	is any $R$-module, then any basis of $M/IM$ lifts to
	a generating set of $M$.  If $M$ is furthermore flat,
	then lifting a basis of $M/IM$, we obtain a surjection
	$F\to M$ whose source is free, with the additional
        property that, if $N$ denotes its kernel, then $N/IN = 0$.
	Thus $N = 0$, and so $F\iso M$, as required.)

\end{aremark}

We begin by establishing some simple consequences of our assumptions
on $C^+$ and~$u$.

\begin{alemma}
	\label{lem:basic properties}\leavevmode
        \begin{enumerate}[label=\normalfont(\arabic*)]
        \item  For each $n \geq 1,$ the quotient $C^+/u^n$ is flat
	over $R$.
      \item  $C^+$ itself is flat over $R$.
      \item  The ring $C^+$ is $u$-adically complete.
        \end{enumerate}


\end{alemma}
\begin{proof}
Claim~(1) follows by an evident induction from assumptions~(A) and~(B)
together with a consideration of the short exact sequences
$$
0 \to C^+/u^n \buildrel u\cdot \over \longrightarrow C^+/u^{n+1}
\to C^+/u \to 0.
$$

Now choose an increasing
filtration $(I_m)$ of $R$ by ideals so that $I_0 = 0$, and such that
each quotient $I_{m+1}/I_{m}$ is of length one (i.e.\ is isomorphic 
to $k$).
Taking into account~(1), the short exact sequences
$$0 \to I_{m} \to I_{m+1} \to I_{m+1}/I_m (\cong k = R/I) \to 0$$
give rise to short exact sequences
$$0 \to I_m\otimes_R (C^+/u^n) \to I_{m+1}\otimes_R (C^+/u^n) \to C^+/(I,u^n) \to 0.$$
Suppose that for some~$m$ we know that the natural morphism
$$I_m \otimes_R C^+ \to \varprojlim_{n} I_m \otimes_R (C^+/u^n)$$
is an isomorphism.
Then passing to the inverse limit over $n$, and taking into account 
both assumption~(C) and the fact that
the transition maps $I_m\otimes_R (C^+/u^{n+1}) \to I_m\otimes_R (C^+/u^n)$
are surjective, so that the relevant $\varprojlim^1$ vanishes,
we obtain a short exact sequence
$$0 \to I_m\otimes_R C^+ \to \varprojlim_n I_{m+1}\otimes_R (C^+/u^n)
\to C^+/I \to 0.$$
The five lemma shows that the natural morphism to this sequence
from the exact sequence
$$I_m\otimes_R C^+ \to I_{m+1}\otimes_R C^+ \to C^+/I \to 0$$
is then an isomorphism.
Proceeding by induction (the case~$m=0$ being trivial), we find (once we reach the
top of our filtration, so that $I_m = R$) that $C^+$ is $u$-adically
complete, and we also find that each of the morphisms
$$I_m\otimes_R C^+ \to C^+$$ 
is injective.  Since any ideal $J$ of $R$ can be placed in such a
filtration $(I_m)$, we find that $C^+$ is flat over $R$.  Thus~(2)
and~(3) are proved.
\end{proof}

We write $C := C^+[1/u]$. Since $C$ is a localisation of $C^+$,
which is $R$-flat by Lemma~\ref{lem:basic properties}~(2),
it is flat over $R$.
Note also that $\Gr^0C \cong (C^+/I)[1/u]$ is a field that is
complete with respect to a non-archimedean absolute value.
In particular, we can do rigid geometry over $\Gr^0 C$;
this is a key point, which we will exploit below.

\begin{adefn}
If $M$ is an $R$-module,
then we write $M\cotimes_R C^+$ to denote the $u$-adic completion
of the tensor product $M\otimes_R C^+$.  We also write
$M\cotimes_R C := (M\cotimes_R C^+)[1/u].$
\end{adefn}

\begin{aremark}
	Remark~\ref{rem:independence of u} shows
	that the formation of $M\cotimes_R C^+$ and $M\cotimes_R C$
	is independent of the choice of $u$ satisfying~(A), (B), and~(C)
	above.
\end{aremark}

\begin{aremark}
	In applications, we will apply this construction primarily 
	in the case when $M = A$ is an $R$-algebra,
	in which case $A\cotimes_R C^+$ and $A\cotimes_R C$ are again
	$R$-algebras.
\end{aremark}

\begin{aexample}
	If we take $C^+ := R[[u]]$ (for an indeterminate $u$),
	then $A\cotimes_R C^+ = A[[u]],$ and $A\cotimes_R C = A((u))$.
\end{aexample}

We next state and prove some additional
properties of $u$-adic completions that we will need.

\begin{alemma}
	\label{lem:non-zero divisor}
	If $M$ is an $R$-module, then multiplication by $u$ 
	is injective on $M\cotimes_A C^+$.
	Consequently {\em (}indeed, equivalently{\em )}, the
	natural map $M\cotimes_R C^+ \to M\cotimes_R C$ is
	injective.

	Furthermore, for each $m \geq 0$, the natural map
	$M\cotimes_R C^+ \to M\otimes_R (C^+/u^m)$
	is surjective,
	and induces an isomorphism
	$$(M\cotimes_R C^+)/u^m (M\cotimes_R C^+) \iso M\otimes_R (C^+/u^m).$$
\end{alemma}
\begin{proof}
	Note that the claim of the second paragraph
	is a general property of completion with respect to finitely 
	generated ideals
(see
e.g.~the statement and proof of
\cite[\href{http://stacks.math.columbia.edu/tag/05GG}{Tag
  05GG}]{stacks-project}). 
We will also deduce it as a byproduct of the proof of the claims
in the first paragraph, to which we now turn.

	Given the definition of $M\cotimes_R C$ as the localization
	$M\cotimes_R C^+[1/u],$ the second claim of the first
	paragraph is clearly
	equivalent to the first.
	As for this first claim,
	we note that for each $m, n \geq 1,$
	Lemma~\ref{lem:basic properties}~(1)
	shows that the terms in the short exact sequence
$$
0 \to C^+/u^n \buildrel u^m \cdot \over \longrightarrow C^+/u^{n+m}
\to C^+/u^m \to 0
$$
	are flat $R$-modules. 
	Thus, tensoring this short exact sequence
	with $M$ over $R$ yields a short exact sequence
	$$0 \to M\otimes_R (C^+/u^n) \to M\otimes_R (C^+/u^{n+1}) \to
	M\otimes_R (C^+/u)\to 0.$$
	Passing to the inverse limit over $n$,
	we obtain the exact sequence
	$$0 \to M\cotimes_R C^+ \buildrel u^m \cdot \over \longrightarrow
	M\cotimes_R C^+ \to M\otimes_R C^+/u^m C^+ $$
	from which the claims of the first paragraph follow.
	If we use the surjectivity of the transition maps in the
	inverse limit to infer
	that the relevant $\varinjlim^1$ vanishes, then we see that
	this sequence is even exact on the right, so that we also 
	obtain a confirmation of the claim of the second paragraph.
\end{proof}


\begin{alemma}
	\label{lem:completion properties}\leavevmode
        \begin{enumerate}[label=\normalfont(\arabic*)]
        \item The functors $M \mapsto M\cotimes_R C^+$ and
          $M\mapsto M\cotimes_R C$ are exact. 
        \item  If $A$ is an $R$-algebra, if $J$ is a finitely
          generated ideal in $A$, and if $M$ is an $A$-module, then
          the natural map
	$$J (M\cotimes_R C^+) \to JM \cotimes_R C^+$$
        is an isomorphism, and consequently there is a short exact
        sequence
	$$0 \to J(M\cotimes_R C^+) \to M\cotimes_R C^+
	\to (M/JM)\cotimes_R C^+ \to 0.$$
      \item If $J_1$ and $J_2$ are ideals in an $R$-algebra $A$,
        with $J_1 \supseteq J_2$, and if $M$ is an $A$-module, then
        there is a natural isomorphism
$$J_1 (M\cotimes_R C^+)/J_2 (M\cotimes_R C^+) \iso (J_1 M/J_2 M) 
\cotimes_R C^+.$$
\end{enumerate}
\end{alemma}
\begin{proof}
	If $0 \to M_1 \to M_2 \to M_3 \to 0$ is a short exact sequence
	of $R$-modules, then Lemma~\ref{lem:basic properties}~(1) 
	implies that 
	$$0 \to M_1\otimes_R (C^+/u^n) \to M_2\otimes_R (C^+/u^n) \to 
	M_3\otimes_R (C^+/u^n) \to  0$$
	is exact for each $n$.  If we pass to the inverse limit over~$n$,
	and note that the transition morphisms are evidently surjective,
	so that the relevant $\varprojlim^1$ vanishes,
	we obtain a short exact sequence
	$$0 \to M_1\cotimes_R C^+ \to M_2\otimes_R C^+ \to 
	M_3\otimes_R C^+ \to  0,$$
	proving the first exactness claim of~(1).
	The second exactness claim follows from the first,
	together with the fact that the localization
	$C := C^+[1/u]$ is flat over $C$, and that there is an
	isomorphism
	$M \cotimes_R C \iso (M\cotimes_R C^+)\otimes_{C^+}C.$



	In order to prove~(2), 
	we apply~(1) to the short exact sequence $0 \to JM \to M
	\to M/JM \to 0$, obtaining a short exact sequence
	\anumequation
	\label{eqn:J ses}
	0 \to (J M)\cotimes_R C^+ \to M\cotimes_R C^+ 
	\to (M/JM)\cotimes_R C^+ \to 0.
\end{equation}
        Thus we may (and do) 
	regard $(JM)\cotimes_R C^+$ as a submodule of $M\cotimes_R C^+.$
        There is a consequent inclusion
	\anumequation
	\label{eqn:J inclusion}
	J (M\cotimes_R C^+)
	\subseteq (JM)\cotimes_R C^+,
        \end{equation}
	whose image is $u$-adically
	dense in the target.
	Since $J$ is a finitely generated ideal,
	we see that $J(M\cotimes_R C^+)$ is furthermore $u$-adically complete,
	and thus that~(\ref{eqn:J inclusion}) is actually an equality,
	establishing the first claim of~(2).

%
	Replacing $(JM)\cotimes_R C^+$ by $J(M\cotimes_R C^+)$ in
	the short exact sequence~(\ref{eqn:J ses}) yields the short
	exact sequence whose existence is asserted in the second claim
	of~(2).

	Claim~(3) follows directly from~(2), and the exactness result
	of~(1).
\end{proof}

The following lemma will allow us to use grading techniques
to study $u$-adic completions.

\begin{alemma}
	\label{lem:passing to graded context}
	If $M$ is an $R$-module,
	then there are natural isomorphisms

	{\em (1)}
	$\Gr^{\bullet} M \otimes_{k} \Gr^0 (C^+/u^n) \iso
	\Gr^{\bullet} \bigl(M\otimes_R
	(C^+/u^n)\bigr)$
{\em (}for any $n  \geq 1${\em )},

{\em (2)}
$\Gr^{\bullet} (M\cotimes_R C^+) \iso
\Gr^{\bullet} M \cotimes_{k} \Gr^0 C^+ ,$
	and

	{\em (3)}
	$\Gr^{\bullet} (M\cotimes_R C)
	\iso
	\Gr^{\bullet} M \cotimes_{k} \Gr^0 C.$ 
\end{alemma}
\begin{proof}
	Lemma~\ref{lem:basic properties}~(1)
	shows that $C^+/u^n$ is flat over $R$,
	and~(1) follows from this,
	together with Lemmas~\ref{lem:graded flat} and~\ref{lem:graded flat bis}.

	Lemma~\ref{lem:completion
		properties}~(3) (taking the inclusion $J_2 \subseteq
	J_1$ to be the various inclusions $I^{i+1} \subseteq I^i$
	in turn) shows that each
        $\Gr^i(M\cotimes_R C^+) \iso \Gr^i M \cotimes_R  C^+.$
        Since $I$ annihilates $\Gr^i(M)$, the target of this
        isomorphism can be rewritten as $\Gr^i M \cotimes_k \Gr^0 \C^+$.
        This proves~(2).


	The $R$-algebra $C$ may be described as the direct limit
	$ C \iso \varinjlim_n u^{-n}C^+$ (the transition
	maps being the obvious inclusions).
	Since the formation of direct limits is compatible with
	tensor products, we find that
	$$\Gr^0 C \iso \varinjlim_n \Gr^0 u^{-n} C^+.$$
	The isomorphism of~(3) is thus obtained from
	the isomorphisms
	$$\Gr^{\bullet} ( M\cotimes_R u^{-n} C^+) \iso
	\Gr^{\bullet} M \cotimes_{k} \Gr^0 u^{-n} C^+$$
	of~(2) by passing to the
	direct limit over $n$.
\end{proof}

\section{Loci of vanishing and of equivariance}
We maintain the notation of Section~\ref{subsec:setting},
and prove some slightly technical results about the ``locus
of vanishing'' and ``locus of equivariance'' of morphisms 
between finitely generated projective $A\cotimes_R C$-modules.

\begin{alemma}
	\label{lem:vanishing locus}
	Let $A$ be an $R$-algebra,
	and let $M$ be a finitely generated projective module 
	over $A\cotimes_R C$.  If $m \in M$, then the ``locus
	of vanishing'' of $m$ is a Zariski closed subset of $\Spec A$.
	More precisely, there is an ideal $J$ of $A$
	such that for any morphism $\phi:A \to B$ of $R$-algebras,
	the image of $m$ in $M_B := (B\cotimes_R C) \otimes_{A\cotimes_R C}
	M$ vanishes if and only if $J$ is contained in the kernel
	of~$\phi$.
\end{alemma}
\begin{proof}
	Since $M$ is finitely generated projective over $A\cotimes_R C$,
	we may find another finitely generated projective $A\cotimes_R C$-module $N$
	such that $M \oplus N$ is free of finite rank over $A\cotimes_R C$.
	Replacing $m \in M$ by $m \oplus 0 \in M \oplus N$, we reduce
	to the case when $M$ is free of finite rank, 
	say $M = A\cotimes_R C^{\oplus r}.$  Writing 
	$m = (x_1,\ldots,x_r)$ with $x_i \in A\cotimes_R C$,
	we see that it suffices to construct a corresponding ideal
	$J_i$ for each $x_i$; the ideal $J := \sum_i J_i$ 
	will then satisfy the claim of the lemma.
	Thus we reduce to the proving the lemma in the case of an element
	$x \in A\cotimes_R C = A\cotimes_R C^+[1/u].$
	Recall from Lemma~\ref{lem:non-zero divisor} that
	$A\cotimes_R C^+$ is a subring of $A\cotimes_R C$.
	Thus, since $u$ is a unit in $A\otimes_R C$,
	we may replace $x$ by $u^n x$ for some sufficiently large
	value of $n$, and assume that $x \in A\cotimes_R C^+$.

	Note that tensoring the inclusion 
	$A\cotimes_R C^+ \hookrightarrow A\cotimes_R C$
	with $B\cotimes_R C^+$ over $A\cotimes_R C^+$
	induces the inclusion
	$$B\cotimes_R C^+ \hookrightarrow B\cotimes_R C 
	= B\cotimes_R C^+[1/u] = (B\cotimes_R C^+)\otimes_{A\cotimes_R C^+}
	A\cotimes_R C.$$
	Thus it suffices to construct an ideal $J$ in $A$
	such that the image of $x$ in $B\cotimes_R C^+$  vanishes
	if and only if the morphism $\phi:A \to B$ factors through $A/J$.
	Let $x_n$ denote the image of $x$ in $A\otimes_R (C^+/u^n)$,
	so that
	$$x = (x_n) \in A\cotimes_R C^+ = \varprojlim_n 
	A\otimes_R (C^+/u^n).$$
	If we let $y_n$ denote the image of $x_n$
	in $B\otimes_R (C^+/u^n)$,
	then the image $y$ of $x$ in $B\cotimes_R C$
	is equal to
       	$$(y_n) \in B\cotimes_R C^+ = \varprojlim_n 
	B\otimes_R (C^+/u^n).$$
	Thus it suffices to construct an ideal $J_n$ in $A$
	such that $y_n$ vanishes if and only if the morphism
	$A\to B$ factors through $A/J_n$; indeed, we may then take
	$J := \sum_n J_n$. 
	
	Lemma~\ref{lem:basic properties} shows that 
	$C^+/u^n$ is flat over $R$, and thus free over $R$.
	(See Remark~\ref{rem:flat implies free}.) 
	If we choose an isomorphism $C^+/u^n \cong R^{\oplus S},$
	then we may write 
	$x_n = (x_{n,s}) \in A^{\oplus S}.$
	If we let $J_{n,s}$ denote the ideal in $A$ generated
	by $x_{n,s}$, then $J_n := \sum_{s \in S} J_{n,s}$ 
	is the required ideal.
\end{proof}

\begin{acor}
	\label{cor:vanishing locus}
	Let $A$ be an $R$-algebra,
	and let $M$ and $N$ be finitely generated projective modules
	over $A\cotimes_R C$.  If $f: M \to N$ is an
        $A\cotimes_R C$-module homomorphism,
	then the ``locus of vanishing''
	of $f$ is a Zariski closed subset of $\Spec A$.
	More precisely, there is an ideal $J$ of $A$
	such that for any morphism $\phi:A \to B$ of $R$-algebras,
	the base-change
	$$f_B :
	M_B := (B\cotimes_R C) \otimes_{A\cotimes_R C} M\to
	N_B := (B\cotimes_R C) \otimes_{A\cotimes_R C} N
	$$
       	vanishes if and only if $J$ is contained in the kernel
	of~$\phi$.
\end{acor}
\begin{proof}
	We may regard $f \in \Hom_{A\cotimes_R C}(M,N)$
	as an element of the module
	$$M^{\vee}\otimes_{A\cotimes_R C} N.$$
	The corollary then follows from Lemma~\ref{lem:vanishing
		locus} applied to $f$ so regarded.
\end{proof}

\begin{acor}
	\label{cor:equivariance locus}
	Let $\sigma$ be a ring automorphism of $C$,
	let $A$ be an $R$-algebra,
	and let $M$ and $N$ be finitely generated projective modules
	over $A\cotimes_R C$ each endowed with a $\sigma$-semi-linear
	automorphism, denoted $\sigma_M$ and $\sigma_N$ respectively. 
	If $f: M \to N$ is an
        $A\cotimes_R C$-module homomorphism,
	then the ``locus of $\sigma$-equivariance''
	of $f$ is a Zariski closed subset of $\Spec A$.
	More precisely, there is an ideal $J$ of $A$
	such that for any morphism $\phi:A \to B$ of $R$-algebras,
	the base-change
	$$f_B :
	M_B := (B\cotimes_R C) \otimes_{A\cotimes_R C} M\to
	N_B := (B\cotimes_R C) \otimes_{A\cotimes_R C} N
	$$
       	satisfies $f_B\circ \sigma_M = \sigma_N \circ f_B$
	if and only if $J$ is contained in the kernel
	of~$\phi$.
\end{acor}
\begin{proof}
	This follows from Corollary~\ref{cor:vanishing locus},
	applied to the morphism $f  - \sigma_N \circ f \circ \sigma_M^{-1}.$
\end{proof}
We also have the following variants on these results for $A\cotimes_RC^+$-modules.
\begin{alem}
  \label{lem: inclusion lattices closed condition}Let~$A$ be an
  $R$-algebra, and let~$M$ and~$N$ be finitely generated projective
  modules over~$A\cotimes_RC^+$. Let~$f:M\to N[1/u]$ be an
  $A\cotimes_RC^+$-module homomorphism. Then the locus where
  $f(M)\subseteq N$ is a Zariski closed subset of~$\Spec A$. More
  precisely, there is an ideal  $J$ of $A$
	such that for any morphism $\phi:A \to B$ of $R$-algebras,
	the base-change
	$$f_B :
	M_B := (B\cotimes_R C) \otimes_{A\cotimes_R C} M\to
	N_B[1/u] := (B\cotimes_R C) \otimes_{A\cotimes_R C} N[1/u]
	$$
       	satisfies $f_B(M_B) \subseteq N_B$
	if and only if $J$ is contained in the kernel
	of~$\phi$.
\end{alem}
\begin{proof}We argue as in the proof of Lemma~\ref{lem:vanishing
    locus}. Since~$M,N$ are finitely generated projective
  over~$A\cotimes_RC^+$, we may find finitely generated projective
  modules~$P,Q$ over~$A\cotimes_RC^+$ such that~$M\oplus P$
  and~$N\oplus Q$ are free. Replacing~$M$ by~$M\oplus P$, $N$
  by~$N\oplus Q$, and~$f$ by~$(f,0)$, we reduce to the case that~$M$,
  $N$ are both finite free~$A\cotimes_RC^+$-modules. After choosing
  bases we may suppose that~$M=(A\cotimes_RC^+)^{\oplus r}$,
  $N=(A\cotimes_RC^+)^{\oplus s}$, and considering the matrix
  representing~$f$, we reduce to the showing that for any
  element~$x\in A\cotimes_RC$, there is an ideal~$J$ of~$A$ such that
  $(\phi\otimes 1)(x)\in B\cotimes_RC^+$ if and only if~$J$ is
  contained in the kernel of~$\phi$. This is a special case of
  Lemma~\ref{lem: element of lattice closed condition} below.
\end{proof}

\begin{alem}
  \label{lem: element of lattice closed condition}Let~$A$ be an
  $R$-algebra, and let~$M$ be a finitely generated projective
  $A\cotimes_RC^+$-module. Let~$x$ be an element of~$M[1/u]$. Then the
  locus where $x\in M$ is a Zariski closed subset of~$\Spec A$. More
  precisely, there is an ideal $J$ of $A$ such that for any morphism
  $\phi:A \to B$ of $R$-algebras, the image of~$x$ in~$M_B[1/u]$ lies
  in~$M_B$ if and only if $J$ is contained in the kernel of~$\phi$.
\end{alem}
\begin{proof}We begin by arguing  as in the proof of Lemma~\ref{lem: inclusion
    lattices closed condition}. Since~$M$ is finitely generated
  projective over~$A\cotimes_RC^+$, we may find a finitely generated
  projective module~$P,Q$ over~$A\cotimes_RC^+$ such that~$M\oplus P$
  is free. Replacing~$M$ by~$M\oplus P$ and ~$x$ by~$(x,0)$, we reduce to the case that~$M$
  is a finite free~$A\cotimes_RC^+$-module. After choosing
  a basis we may suppose that~$M=(A\cotimes_RC^+)^{\oplus r}$, so that
  we can reduce to the case that~$M=A\cotimes_RC^+$.
  
  Choose~$n$ sufficiently large so that $u^nx\in A\cotimes_RC^+$. Then
  $(\phi\otimes 1)(x)\in B\cotimes_RC^+$ if and only if the image
  of~$(\phi\otimes 1)(x)$ in
  \[u^{-n}(B\cotimes_RC^+)/(B\cotimes_RC^+)=B\otimes_R(u^{-n}C^+)/C^+)\]
  vanishes. As in the proof of Lemma~\ref{lem:vanishing locus},
  $u^{-n}C^+/C^+$ is a free~$R$-module, so we can and do choose an
  isomorphism~$u^{-n}C^+/C^+\cong R^{\oplus S}$. Then if we
  write~$x=(x_s)\in A^{\oplus S}\cong A\otimes_R (u^{-n}C^+/C^+)$, we
  can take~$J$ to be the ideal generated by the~$x_s$.
\end{proof}

\section{Extending actions}
We continue to remain in the context of Section~\ref{subsec:setting}.
We will show that various sorts of actions on $C^+$ or $C$
can be extended to the completed tensor products $A\cotimes_R C^+$
and $A\cotimes_R C$.

Throughout our discussion,
we endow $C^+$ with the $u$-adic topology, while we endow $C$
with the unique topology which makes it a topological group,
and in which $C^+$ is an open subgroup.  
More generally, if $A$ is any $R$-algebra,
we again endow
$A\cotimes_R C^+$ with the $u$-adic topology,
while we endow
$A\cotimes_R C$ 
with the unique topology which makes it a topological group,
and in which $A\cotimes_R C^+$ (endowed with its $u$-adic topology)
is an open subgroup.  

\begin{alemma}
	\label{lem:extending endomorphisms}
	Suppose that $C^+$ {\em (}resp.\ $C${\em )}
	is equipped with a continuous $R$-linear endomorphism~$\varphi$.
	Then there is a unique extension of $\varphi$
	to a continuous $A$-linear endomorphism of $A\cotimes_R C^+$
	{\em (}resp.\ $A\cotimes_R C${\em )}.
\end{alemma}
\begin{proof}
	Lemma~\ref{lem:non-zero divisor} shows
	that $u$-adic topology on $A\cotimes_R C^+$ coincides
	with its inverse limit topology,
	if we recall that
	$A\cotimes_R C^+ := \varprojlim_n A\otimes_R (C^+/u^n),$
	and we endow each of the objects appearing in the 
	inverse limit with its discrete topology.
	Given this, it is immediate that a $u$-adically continuous
	$R$-linear endomorphism
	$\varphi$ of $C^+$ extends uniquely to a $u$-adically
	continuous $A$-linear endomorphism of $A\cotimes_R C^+$.
	
	The case when $\varphi$ is a continuous $R$-linear endomorphism
	of $C$ is only slightly more involved.  Namely, to say
	that $\varphi: C \to C$ is continuous is to say that,
	for some $m \geq 0$, we have an inclusion $\varphi( u^m C^+)
	\subseteq C^+$, and that the induced morphism
	$\varphi: u^m C^+ \to C^+$ is continuous, when each of the
	source and the target are endowed with their $u$-adic topologies.
	Again taking into account
	Lemma~\ref{lem:non-zero divisor},
       we then obtain a unique $A$-linear and continuous morphism
       \begin{multline*}
A\cotimes_R C =
\bigl( \varprojlim_n A\otimes_R (u^m C^+/u^{m+n}) \bigr) [1/u] 
\\
\longrightarrow
\bigl( \varprojlim_n A\otimes_R (C^+/u^{n}) \bigr) [1/u]  
= A\cotimes_R C,
\end{multline*}
as required.
\end{proof}

If $G$ is a topological group acting continuously on~$C^+$
(i.e.\ acting in such a way that the action map $G\times C^+ \to C^+$
is continuous), then for any $g \in G$ and $m\geq 0$, we may find $n \geq 0$,
and an open neighbourhood $U$ of~$g$,
such that $U \times u^n C^+ \subseteq u^m C^+$.  (This just expresses
the continuity of the action map at the point $(g,0) \in G\times C^+$.)
If $G$ is furthermore compact, then $G$ is covered by finitely many
such open sets $U$, and thus
for any given choice of~$m$, we may in fact find $n$
such that $G \cdot u^n C^+ \subseteq u^m C^+$.
An identical remark applies in the case of a compact topological group
acting on~$C$.

We will need to introduce an additional condition on $G$-actions on~$C$,
which is related to, but doesn't seem to follow from,
the previous considerations. 

\begin{adefn}\label{adefn: bounded group action}
	If $G$ is a group,
	we say that a $G$-action on $C$
	is {\em bounded} if, for any $M\geq 0,$ there exists
	$N \geq 0$ such that $G\cdot u^{-M} C^+ \subseteq u^{-N} C^+.$
\end{adefn} 

\begin{arem}
  \label{arem: bounded doesn't depend on u}
  It follows from Remark~\ref{rem:independence of u}
  that the notion of a group action being bounded is independent
  of the choice of $u\in C^+$ satisfying conditions~(A), (B), and~(C)
  of Section~\ref{subsec:setting}.
\end{arem}

We then have the following lemma.


\begin{alemma}
	\label{lem:extending group actions}
	Suppose that $G$ is a compact topological group,
	and that $C^+$ {\em (}resp.~$C${\em )}
	is equipped with a continuous {\em (}resp.\ continuous
	and bounded{\em )} $R$-linear $G$-action. 
	Then there is a unique extension of this action 
	to a continuous $A$-linear action of $G$
	on $A\cotimes_R C^+$
	{\em (}resp.\ $A\cotimes_R C${\em )}.
\end{alemma}
\begin{proof}
The proof of this lemma is very similar to that of Lemma~\ref{lem:extending
endomorphisms}.  Indeed, applying that lemma (and taking into account
its uniqueness statement) we obtain the desired action of $G$.
It remains to confirm that this action is again continuous.

Consider first the case of a $G$-action on $C^+$.
For each $m \geq 0,$ we saw above that $G \cdot u^n C^+ \subseteq
u^m C^+$ if $n$ is sufficiently large, 
so that we obtain a well-defined map
$G\times A\cotimes_R C^+ = G\times \varprojlim_n A\otimes_R C^+/u^n 
\to A\otimes_R C^+/u^m$ which is continuous (since it is just obtained
from the continuous $G$-action on $C^+$ by tensoring with~$A$).
Passing to the projective limit in~$m$ then yields the desired continuity.

In the case of a $G$-action on $C$, we argue similarly:
namely, the assumptions of continuity and boundedness
of the $G$-action yield continuous morphisms
$$ G\times u^{-M} C^+/u^n C^+ \to  u^{-N} C^+/ u^m C^+;$$
Passing to the projective limit in $n$, then in $m$, and then to the
inductive limit in $N$, and then in $M$, we deduce the desired 
continuity of the $G$-action on $A\cotimes_R C$.
\end{proof}

\section{A rigid analytic perspective}
\label{subsec:rigid analytic}
We keep ourselves in the setting
of Section~\ref{subsec:setting}. 
We are interested in studying finitely generated projective
modules over $A\cotimes_R C$, especially their 
descent properties.
In the present general context, we don't
know whether an analogue of Drinfeld's descent result for Tate modules
(see~\cite[Thm.\ 3.3, 3.11]{MR2181808} and~\cite[Thm.\ 5.1.18]{EGstacktheoreticimages}) holds
for arbitrary $R$-algebras $A$.  However, we will see that such
a descent result is valid in the more restricted setting of finite
type $R$-algebras.  In the case of finite type $k$-algebras,
we will prove this using results from rigid analysis.  For finite type
$R$-algebras that are not necessarily annihilated by $I$,
we will use grading techniques to reduce to the case of $k$-algebras.

If $A$ is a finite type $k$-algebra,
then $\Spf (A\cotimes_k \Gr^0 C^+)$ (where the $\Spf$ 
is taken with respect to the $u$-adic topology on $A \cotimes_k
\Gr^0 C^+$) 
is a formal scheme of finite type over $\Spf \Gr^0 C^+$,
whose rigid analytic generic fibre
$\Max\Spec (A\cotimes_k \Gr^0 C)$
is an affinoid rigid analytic space over the complete
non-archimedean field $\Gr^0 C$.


As an application of this rigid analytic point-of-view (together
with graded techniques), we first establish the following
proposition.
 
\begin{aprop}\label{prop: rigid analysis FGK}\leavevmode
  
  \begin{enumerate}[label=\normalfont(\arabic*)]
  \item  If $A$ is a finite type $R$-algebra, then $A\cotimes_R C$ is
    Noetherian.
  \item  If $A \to B$ is a {\em (}faithfully{\em )} flat morphism
    of finite type $R$-algebras, then the induced morphisms
    $A\cotimes_R C^+ \to B\cotimes_R C^+$ and
    $A\cotimes_R C \to B\cotimes_R C$ are again {\em (}faithfully{\em
      )} flat.
  \end{enumerate}
\end{aprop}
\begin{proof}
	In the case when $A$ is in fact a $k$-algebra,
	it follows from the preceding discussion that
        $A\cotimes_R C$ is an affinoid algebra over~$\Gr^0 C$, and is therefore Noetherian by~\cite[Thm.\ 1, \S5.2.6]{MR746961}.
	In the general case, 
	Lemma~\ref{lem:passing to graded context}~(3)
	gives an isomorphism
	$\Gr^{\bullet} (A\cotimes_R C)
	\iso
	\Gr^{\bullet} A \cotimes_{k} \Gr^0 C.$ 
	Since $\Gr^{\bullet} A$ is a finite type $k$-algebra,
	it follows from the case already proved that
	the target of this isomorphism is Noetherian, and thus
	so is its source.  Thus $A\cotimes_R C$ is itself Noetherian,
	by Lemma~\ref{lem:graded Noetherian}.

        We explain how our rigid analytic perspective
	proves a part of~(2).
	If we suppose first that $A\to B$ is a flat
        morphism of $k$-algebras,
	then so is the morphism $A\otimes_R C \to B\otimes_R C.$
	Now the preceding discussion shows that the morphisms
	$\Max\Spec (A\cotimes_R C) \to \Max\Spec (A\otimes_R C)$
	and 
	$\Max\Spec (B\cotimes_R C) \to \Max\Spec (B\otimes_R C)$
	are open immersions of rigid analytic spaces over $\Gr^0 C$,
	so that the induced morphism 
	$\Max\Spec (A\cotimes_R C) \to \Max\Spec (B\cotimes_R C)$
	is also flat.  This proves (the $k$-algebra case of) the
	second flatness claim of~(2).  
	If $A\to B$ is a flat morphism
	of $R$-algebras that are not necessarily $k$-algebras,
	then we deduce from Lemma~\ref{lem:graded flat}
	that $\Gr^{\bullet} A \to \Gr^{\bullet} B$ is flat,
	and thus, from what we've already shown, that 
	$$\Gr^{\bullet} A \cotimes_{k} \Gr^0 C
	\to \Gr^{\bullet} B \cotimes_{k} \Gr^0 C$$
	is flat.  
	Lemma~\ref{lem:passing to graded context}~(3) then
	shows that the morphism
	\[\Gr^{\bullet} (A\cotimes_R C ) \to \Gr^{\bullet} (B\cotimes_R C)\]
        is flat,
	and one more application of Lemma~\ref{lem:graded flat}
	completes the proof of the second flatness claim of~(2).

	It is not clear to us whether one can deduce the second
	faithful flatness claim of~(2) by this style of argument.
	In order to prove this result, 
	as well as to obtain the first set of claims of~(2),
	we appeal to the results
	of~\cite{MR2774689}.  Indeed, it is clear that the
	second set of claims in~(2) follows immediately from the first,
	and it the first set of claims that we will now prove.

        We first note that if the morphism
	$A \to B$ is flat (resp.\ faithfully
	flat), then so are each of the morphisms
	$A\otimes_R C^+/u^n A \to B\otimes_R ^+C/u^n.$
	In the terminology
	of~\cite[\S 5.2]{MR2774689},
	the morphism
	$A\cotimes_R C^+ \to B\cotimes_R C^+$
	is adically flat (resp.\ adically faithfully flat); here
	we regard the source and target as being endowed with their
	$u$-adic topologies.  The discussion at the beginning
	of~\cite[\S 5.2]{MR2774689}
(see also our Proposition~\ref{prop: FGK flat})
then shows that
	$A\cotimes_R C^+ \to B\cotimes_R C^+$
	is flat.  In the
	adically faithfully flat case, we deduce in addition 
	from~\cite[Prop.~5.2.1~(2)]{MR2774689}
(again, see also Proposition~\ref{prop: FGK flat}) 
that $A\cotimes_R C^+ \to B\cotimes_R C^+$
	is faithfully flat;  
	note that
	$A\cotimes_R C^+$ and $B\cotimes_R C^+$ are Noetherian outside 
	$u$ by~(1) of the present proposition.
\end{proof}

\section{Descent for $A\cotimes_R C$-modules}
We now establish the descent result alluded to above; 
it is analogous to 
Drinfeld's \cite[Thm.\ 5.1.18~(1)]{EGstacktheoreticimages}
but is restricted to the context of finite type $R$-algebras.
We use grading techniques to reduce to the case of finite type
$k$-algebras, where we can then apply known results from rigid analysis.

\begin{atheorem}
	Let $A \to B$ be a faithfully flat morphism of 
	finite type $R$-algebras.  Then the functor
	$M \mapsto
	M\otimes_{A\cotimes_R C} (B\cotimes_R C)$
	induces an equivalence of categories between
	the category of finite type $A\cotimes_R C$-modules
	and the category of finite type $B\cotimes_R C$-modules
	equipped with descent data.
\end{atheorem}

\begin{aremark}\label{rem: what is descent data}Recall (for example
  from~\cite[\S 6.1]{MR1045822}) that if~$S\to S'$ is a ring
  homomorphism, then a descent datum for an~$S'$-module~$M'$ is an
  isomorphism of~$S'\otimes_S S'$-modules $S'\otimes_SM'\isoto
  M'\otimes_SS'$, satisfying a certain cocycle condition. Given a
  descent datum, we have a pair of morphisms of~$S$-modules $M'\to M'\otimes_SS'$, namely the
  obvious morphism and the composite of the obvious morphism $M'\to
  S'\otimes_SM'$ with the descent datum isomorphism. We let~$K$ denote
  the $S$-module given by the kernel of the difference of these two
  maps.  We
  refer to the functor $M'\mapsto K$,  from the category of
  $S'$-modules with descent data to the category of $S$-modules, as
  the \emph{kernel functor}.  (If~$S\to S'$ is faithfully flat, then the theory of faithfully
  flat descent shows that~$M'$ is the extension of scalars from~$S$ to~$S'$
  of~$K$.)
\end{aremark}
\begin{proof}
	We first treat the case when $A$ and $B$ are $k$-algebras.
	In this case, the category of
	finite type $A\cotimes_R C$-modules (resp.\ finite type $B\cotimes_R C$-modules)
	is equivalent to the category of coherent sheaves on
       the affinoid rigid analytic space $\Max\Spec (A\cotimes_R C)$ (resp.\
       $\Max\Spec (B\cotimes_R C)$) over the field $\Gr^0 C$, 
       and the statement follows from faithfully flat descent
       for rigid analytic coherent sheaves \cite[Thm.~3.1]{MR1603849}.

       We reduce the general case of $R$-algebras to the case of 
       $k$-algebras via the usual graded arguments.
       To ease notation, we let $F:\cC \to \cD$ denote the base-change functor 
	$M \mapsto
	M\otimes_{A\cotimes_R C} (B\cotimes_R C)$
	from
	the category $\cC$ of finite type $A\cotimes_R C$-modules
	to
	the category $\cD$ of finite type $B\cotimes_R C$-modules
	equipped with descent data,
	and let $G:\cD \to \cC$ denote the kernel functor.
	There are evident natural transformations
	\anumequation
	\label{eqn:adjunctions}
	G\circ F \to \id_{\cC} \quad \text{ and } \quad F \circ G \to \id_{\cD},
\end{equation}
	which we claim are natural isomorphisms (so that
	$G$ provides a quasi-inverse to $F$, proving in particular
	that $F$ induces an equivalence of categories).

	We let $\overline{\cC}$ denote the category of finite type
	$\Gr^{0} (A\cotimes_R C)$-modules, 
	and let $\overline{\cD}$ denote the category of finite type
	$\Gr^{0} (B\cotimes_R C)$-modules
	equipped with descent data to $\Gr^{0} (A\cotimes_R C)$.
	We let $\overline{F}:\overline{\cC} \to \overline{\cD}$
	denote the base-change functor
	$\text{--}\otimes_{\Gr^{0} (A\cotimes_R C)}
	\Gr^{0} (B\cotimes_R C)$,
	and let $\overline{G}:\overline{D} \to \overline{C}$
	denote the kernel functor. 
        The case of $k$-algebras that we have already treated shows that
	there are natural isomorphisms
	\anumequation
	\label{eqn:overline adjunctions}
	\overline{G}\circ \overline{F} \iso \id_{\overline{\cC}}
	\quad \text{ and } \quad
	\overline{F}\circ \overline{G} \iso \id_{\overline{\cD}}.
\end{equation}

	We think of passage to the associated graded $\Gr^{\bullet}$
	as inducing functors $\cC \to \overline{\cC}$ and also 
	$\cD\to \overline{\cD}$.
	It follows from 
	Lemmas~\ref{lem:graded flat} and~\ref{lem:graded flat bis},
        together with Proposition~\ref{prop: rigid analysis FGK},
	that there is a natural isomorphism
	$$\Gr^{\bullet} \circ F \iso \overline{F} \circ \Gr^{\bullet};$$
	there is also an evident natural isomorphism
	$$\Gr^{\bullet} \circ G \iso \overline{G} \circ \Gr^{\bullet};$$
	furthermore, these natural isomorphisms are compatible with the natural transformations~(\ref{eqn:adjunctions}) and~(\ref{eqn:overline adjunctions}).
	Since the natural transformations~(\ref{eqn:overline adjunctions})
	are isomorphisms, the same is true of the natural
        transformations~(\ref{eqn:adjunctions}), by Lemma~\ref{lem: graded iso is iso}. 
	This completes the proof that $F$ is an equivalence.
       \end{proof}



\chapter{Topological groups and modules}
\label{app:topological groups}
In this appendix we recall some more-or-less well-known facts regarding topological groups
and modules,
for which we haven't located a convenient reference.
We found the note \cite{Kaye-Polish} to be a useful
reference for the basic facts regarding Polish topological groups.

Recall the following definition.

\begin{adefn} \index{Polish group}
	A topological space is called {\em Polish} if it is separable (i.e.\
	contains a countable dense subset) and completely metrizable.
	A topological group is called {\em Polish} if its underlying 
	topological space is Polish.  Similarly,
	a topological ring is called {\em Polish} if its underlying 
	topological space is Polish.  
\end{adefn}

Our interest in Polish groups is due to the following lemma and
corollary; see
also~\cite[\href{https://stacks.math.columbia.edu/tag/0CQW}{Tag
  0CQW}]{stacks-project} for a more algebraic proof of a closely
related result.
\begin{alemma}
	\label{lem:open mapping}
	If $\phi: G \to H$ is a continuous homomorphism
        between two Polish topological groups,
	then the following are equivalent:
	\begin{enumerate}
		\item $\phi$ is open.
		\item $\phi(G)$ is not meagre in $H$.
	\end{enumerate}
\end{alemma}
\begin{proof}
	If~(1) holds, then $\phi(G)$ is a non-empty open subset
        of the Polish, and hence Baire, space $H$,
	and so is not meagre.  Thus~(1) implies~(2).
	For the converse, see e.g.~\cite[Thm.~18]{Kaye-Polish}.
\end{proof}

\begin{acor}
	\label{cor:open polish}
	A continuous surjective homomorphism
       	of Polish topological groups is necessarily
	open.
\end{acor}
\begin{proof}
	A completely metrizable space is Baire, and hence is not a
	meagre subset of itself.  The corollary thus follows from
	the implication~``(2) $\implies$ (1)'' of Lemma~\ref{lem:open mapping}.
\end{proof}

Any metrizable space is first countable (i.e.\ each point contains a countable
neighbourhood basis).  
Recall that, conversely, if $G$ is a topological group, then $G$
is first countable if and only if the identity element $1$ admits
a countable neighbourhood basis,
and in this case, if
(the underlying topological space of) $G$ is Hausdorff then it is in fact
metrizable (this is the Birkhoff--Kakutani theorem),
and even admits a (left or right) translation invariant metric.
In particular, (the underlying topological space of)
a Hausdorff topological group is separable and metrizable if and only 
if it is {\em second countable} (i.e.\ admits a countable
basis for its topology).

Recall also that a topological group $G$ admits a canonical uniform structure,
so that it makes sense to speak of $G$ being complete.
For the sake of completeness, we note that the underlying
topological space of $G$ being Polish 
forces $G$ to be complete.

\begin{alemma}
	\label{lem:complete}
	If $G$ is a Polish topological group,
	then $G$ is complete.
\end{alemma}
\begin{proof}
	Since $G$ is metrizable, it is Hausdorff, and so we may
	consider the canonical embedding of topological groups
	$G \hookrightarrow \widehat{G}$ of $G$ into its completion.
	Since $G$ is separable and dense in~$\widehat{G}$, we see
	that $\widehat{G}$ is separable.  Also, $\widehat{G}$ is a complete
	metric space (it may be identified with the metric space
	completion of $G$ with respect to any invariant
	metric inducing the topology on $G$).  Thus $\widehat{G}$ is Polish.
	Any Polish subset of a Polish space is~$G_{\delta}$ (i.e.\ a
        countable intersection of open sets),
	and so in particular $G$ is $G_{\delta}$ inside~$\widehat{G}$.
	Since $\widehat{G}$ is Polish, and thus Baire,
	we find that $G$ is not meagre in~$\widehat{G}$.
	It follows from Lemma~\ref{lem:open mapping}
	that the inclusion $G\into\widehat{G}$ is open, and thus that $G$ is an open subgroup
	of its completion~$\widehat{G}$.  Since an open subgroup
	of a topological group is also closed, we find that
	$G$ closed, as well as dense, in~$\widehat{G}$,
	and thus that $G = \widehat{G}$, which is to say, $G$
	is complete, as claimed.
\end{proof}


\begin{aremark}
	\label{rem:Polish}
It follows from Lemma~\ref{lem:complete}
and the discussion preceding it that
a topological group is Polish if and only if it is complete
and second countable.
\end{aremark}

We next
recall a result about Noetherian and Polish topological rings, 
which is inspired by a result of Grauert and Remmert
in the theory of classical Banach algebras
\cite[App.\ to \S I.5]{MR0316742}.
(See~\cite[Prop.\ 3, \S3.7.3]{MR746961}  
for the analogous result in the context of non-archimedean Banach algebras.)

\begin{aprop}
	\label{prop:module topologies}
	Let $A$ be a Polish {\em (}or, equivalently,
	a \emph{(}Hausdorff\emph{)} complete and second countable{\em )} topological ring 
	which is Noetherian {\em (}as an abstract ring{\em )},
	and which contains an open additive subgroup that
	is closed under multiplication, and consists of topologically
	nilpotent elements.
	Then any finitely generated $A$-module $M$ has a unique
	completely metrizable topology with respect to which it becomes
	a topological $A$-module.  Furthermore, $M$ is complete
	with respect to this topology, any submodule of $M$
	is closed, and any morphism of finitely generated
	$A$-modules is automatically continuous, has closed image,
	and induces an open mapping from its domain onto its image.
\end{aprop}
\begin{proof}
  To begin with, suppose that $M$ is a finitely generated $A$-module
  which is furthermore endowed with a completely metrizable topology
  which makes it a topological $A$-module.  We claim first that $M$ is
  in fact complete, and that any surjection of $A$-modules $A^n \to M$
  (for some $n \geq 1$) is continuous and open. It then follows that
  $M$ is endowed with the quotient topology via this map, and thus the
  completely metrizable topology on $M$ is uniquely determined (if it
  exists).

	To see the claim, note that since $M$ is a finitely generated $A$-module,
	we may choose a surjective homomorphism of $A$-modules
	$A^n \to M$ for some $n~\geq~1$.  Furthermore,
	since $M$ is a topological $A$-module, any such surjection
	is continuous.  Since $A$ is separable as a topological
	space by assumption, we see that $M$ is as well.  (The image
	of any countable dense subset of $A^n$ is a countable dense 
	subset of $M$.)  Thus $M$ is Polish, as is $A^n$,
	and it follows from Lemma~\ref{lem:complete} that
	$M$ is complete (as a topological module), while
        it follows from Corollary~\ref{cor:open polish}	that
	the given surjection $A^n \to M$ is open, as claimed.

	We next claim that any $A$-submodule of $M$ is closed.
	To see this, let $N$ be a submodule of $M$,
	and let $\overline{N}$ denote its closure.
	Since $A$ is Noetherian by assumption, so is its finitely 
	generated module $M$, and thus $\overline{N}$ is finitely 
	generated, say by the elements $x_1,\ldots,x_n$.
	Applying the results proved above for $M$ to
	its closed (and hence completely metrizable topological) submodule
	$\overline{N}$, we find that
	the surjection $A^n \to \overline{N}$ given by
	$(a_1,\ldots,a_n) \mapsto \sum_{i=1}^n a_i x_i$ 
	is an open mapping.
        Consequently,
	if we let $I$ be an open additive subgroup of $A$ which
	satisfies $I^2 \subseteq I,$ and which consists of topologically
	nilpotent elements (such an $I$ exists by assumption),
	then $I x_1 + \cdots + I x_n $ is an open subset of $\overline{N}$.
	Since $N$ is dense in $\overline{N}$,
	we conclude that
	$N + I x_1 + \cdots + I x_n  = \overline{N},$
        and consequently we may write
        $$x_i = \sum_{i=1}^n a_{ij} x_j + y_i$$ for each $1 \leq i \leq n$,
	for some $y_i \in N$ and $a_{ij} \in I$. 
	Rearranging, we find that
	$$ \begin{pmatrix} 1 - a_{11} & - a_{1 2} & \cdots & - a_{1 n} \\
	       -a_{21} & 1 - a_{22} & \cdots & - a_{2 n} \\
	       \vdots & \vdots & \ddots & \vdots \\
-a_{n1} & - a_{2n} & \cdots & 1 - a_{nn} \end{pmatrix}
\begin{pmatrix} x_1 \\ x_2 \\ \vdots \\ x_n\end{pmatrix}
= 
\begin{pmatrix} y_1 \\ y_2 \\ \vdots \\ y_n\end{pmatrix}
$$	       
The determinant of the matrix 
	$ \begin{pmatrix} 1 - a_{11} & - a_{1 2} & \cdots & - a_{1 n} \\
	       -a_{21} & 1 - a_{22} & \cdots & - a_{2 n} \\
	       \vdots & \vdots & \ddots & \vdots \\
-a_{n1} & - a_{2n} & \cdots & 1 - a_{nn} \end{pmatrix}$
lies in $1 + I$, and thus is a unit (since $A$ is complete
and $I$ consists of topologically nilpotent elements).
Thus we find that $x_1,\ldots,x_n$ lies in the $A$-span
of $y_1,\ldots,y_n$, so that in fact $\overline{N} \subseteq N$.
Thus $N$ is indeed closed, as claimed.

Applying the preceding result to $A^n$, we find that
any $A$-submodule of $A^n$ is closed, and hence that any
finitely generated $A$-module $M$ is isomorphic to a quotient
$A^n/N$, where $N$ is a closed submodule of $A^n$.  The quotient
topology on $A^n/N$ makes it into a complete and metrizable
topological $A$-module, and thus $M$ does indeed admit a complete
and metrizable topological $A$-module structure.   We have already
seen that this topology is unique, and that any $A$-submodule of $M$
is closed.    Finally, we see that
homomorphisms between finitely generated $A$-modules are necessarily
continuous and that their images
(being $A$-submodules of their targets)
are necessarily closed, and a final application
of Corollary~\ref{cor:open polish} shows that they induce open
mappings onto their images.
\end{proof}

\chapter{Tate modules and continuity}\label{app: Tate modules and continuity}
In this appendix we study continuity
conditions for group actions on modules over Laurent series rings. We
begin with some results on the topology of such modules, before
introducing group actions and related notions.

\section{Topologies and lattices}
We fix a finite extension
$E/\Qp$ with ring of integers~$\cO$, uniformizer~$\varpi$, and residue
field~$\F$, and we also fix a finite extension~$k/\F_p$. If $A$ is a $p$-adically complete~$\cO$-algebra, we write
$\Aplus_A:=(W(k)\otimes_{\Zp}A)[[T]]$, and we let $\AAA_A$ be the $p$-adic completion of~
$\Aplus_A[1/T]$.

\begin{arem}
  \label{rem:topologies on M}
  Any
  finitely generated projective $\AAA_A$-module $M$ has a natural topology.
  Indeed, we may write $M$ as a direct summand of $\AAA_A^n$ for some
  $n \geq~1$.  We then endow this latter module with its product topology,
  and endow $M$ with the subspace topology. 
  More intrinsically,
  since $A$ is $p$-adically complete (by assumption),
  we may write $M = \varprojlim_a M/p^a M$.   Each of the quotients
  $M/p^a M$ is then a projective $(A/p^a A)((T))$-module,
  and so has natural topology, making it a Tate $A/p^a A$-module in
  the sense \index{Tate module}
  of~\cite{MR2181808}.
  The topology on $M$ is then the projective limit of these Tate module
  topologies.
  We note that multiplication by~$T$ is topologically nilpotent
  on~$M$.

  Similarly, any Zariski locally finite free $\Aplus_A$-module $\gM$ has a natural topology.
  We may describe this topology in an analogous manner to the
  case considered in the preceding paragraph.
  Namely, such a module
  is a finitely generated and projective $\Aplus_A$-module,
  and thus is a direct summand of $(\Aplus_A)^n$ for some $n \geq 0$.
  If we endow this latter module with its product topology
  then the natural topology on $\gM$ is its corresponding subspace topology. 
  In this case, though, this topology admits a more intrinsic and succinct
  description: it is the $(p,T)$-adic topology on $\gM$.
\end{arem}

We remind the reader that a topological group $G$
is said to be {\em Polish} if its underlying
topological space is Polish, i.e.\ is separable and completely metrizable.
As explained in
Remark~\ref{rem:Polish},
this is equivalent
to $G$ being complete (as a topological group) and second countable
(as a topological space).

\begin{alemma}
	\label{lem:polish}
	If $A$ is a $p$-adically complete $\cO$-algebra
	for which $A/p$ is countable, then $\AAA_A$ is
	Polish,
	and consequently any finitely generated projective
	$\AAA_A$-module is Polish, when endowed with its canonical topology.
\end{alemma}
\begin{proof}
	We will use the fact that a countable 
	product of Polish spaces is Polish,
	as is a closed subspace of a Polish space.
	It follows from these facts,
	and from the description of the
	natural topology on a finitely generated projective
	$\AAA_A$-module given in Remark~\ref{rem:topologies on M},
	that the canonical topology on any such module is Polish,
	provided that that $\AAA_A$ itself is Polish.

	We may write $\AAA_A = \varprojlim_a \AAA_{A/p^aA}$,
	so that $\AAA_A$ is a closed subset of the (countable)
	product of the various spaces $\AAA_{A/p^aA}$.
	It suffices, then, to show that each of these 
	spaces is Polish; in other words,
	we reduce to the case when $A$ is a countable $\Z/p^a$-algebra
	for some $a \geq 1$.

	Since $A$ is countable, so is each of the quotients
	$\Aplus_A/T^n$.  Choose a subset $X_n \subseteq
	\Aplus_A$ which maps bijectively onto $\Aplus_A/T^n$,
	and write $X := \bigcup_{m,n} T^{-m} X_n.$  
	Then $X$ is a countable dense subset of $\AAA_A$,
	and thus $\AAA_A$ is separable.  
	
	The topology on $\Aplus_A$ is the $T$-adic topology,
	and thus is metrizable (as is any $I$-adic topology
	on a ring).  Since $\Aplus_A$ is $T$-adically complete,
	it is in fact completely metrizable.  The same is then
	evidently true for
	\[\AAA_A  := \Aplus_A[1/T] = \bigcup_m T^{-m} \Aplus_A.\qedhere\]
\end{proof}

\begin{arem}
	\label{rem:Polish comparision}
	Clearly $T\Aplus_A$ is an open subgroup of $\AAA_A$
that is closed under multiplication
	and consists of topologically nilpotent elements.
	Thus if $A/p$ is countable, then $\AAA_A$ satisfies
	the conditions of Proposition~\ref{prop:module topologies},
and thus the canonical topology on finitely generated projective $\AAA_A$-modules
constructed in Remark~\ref{rem:topologies on M} is a particular
case of the canonical topology constructed on any finitely generated
$\AAA_A$-module in Proposition~\ref{prop:module topologies}.
\end{arem}

We now present some additional facts related to the preceding concepts
which will be needed in the sequel.

\begin{alem}
  \label{lem: can check projectivity of Kisin and etale phi modules modulo
    p^n}A finitely generated $\Aplus_A$-module is
  projective of rank~$d$ if and only if for each~$a\ge~1$,
  the quotient
  $\gM/p^a\gM$ is projective of rank~$d$ as an
  $\Aplus_{A/p^aA}$-module. Similarly, a finitely generated 
  $\AAA_A$-module~$M$ is projective of rank~$d$ if and only if for
  each~$a\ge~1$,
  the quotient
  $M/p^aM$ is projective of rank~$d$ as an $\AAA_{A/p^aA}$-module.
\end{alem}
\begin{proof}
  This is immediate from~\cite[Prop.\
0.7.2.10(ii)]{MR3075000}, applied to the $p$-adically complete
  rings~$\Aplus_A$ and $\AAA_A$.
\end{proof}

   \begin{adefn}
    \label{defn: lattice}If~$M$ is a finitely generated $\AAA_A$-module, then a
    \emph{lattice} in~$M$ is a finitely generated $\Aplus_A$-submodule
    $\gM\subseteq M$  whose $\AAA_A$-span \index{Tate module}
    is~$M$.
  \end{adefn}

  Note that any finitely generated $\AAA_A$-module contains a lattice.
  (If $f:\AAA_A^r \to M$ is a surjection from a finitely generated 
  free $\AAA_A$-module onto the finitely generated $\AAA_A$-module $M$,
  then the image of the restriction of $f$ to $(\Aplus_A)^r$ 
  is a lattice in~$M$.)

  We record some additional lemmas and a remark which apply in the case
  when $A$ is an $\cO/\varpi^a$-algebra for some $a \geq 1$.

  \begin{alemma}
	  \label{lem:lattice properties}
	  Let $A$ be an $\cO/\varpi^a$-module for some $a \geq 1$,
	  and let $M$ be a finitely generated $\AAA_A$-module.
	  
	  \begin{enumerate}
		  \item
	  If $\gM$ and $\gN$ are two lattices
	  contained in a finitely generated $\AAA_A$-module~$M$,
	  then there exists $n \geq 0$
	  such that $T^n\gM \subseteq \gN \subseteq T^{-n} \gM$.
  \item   If in addition either $A$ is Noetherian or $M$ is projective,
	  then any lattice in $M$ is $T$-adically complete.
  \item 
	  If $A$ is Noetherian, if $\gM$ is
	  a lattice in~$M$, and if $\gN$ is an $\Aplus_A$-submodule
	  of $M$ 
	  such that $T^n\gM \subseteq \gN \subseteq T^{-n} \gM$
	  for some $n \geq 0$, 
	  then $\gN$ is a lattice in~$M$.
  \end{enumerate}
  \end{alemma}
  \begin{proof}
	  To prove~(1), it suffices to prove one of the inclusions;
	  the reverse inclusion may then be obtained (possibly after
	  increasing $n$) by switching the roles of $\gM$ and $\gN$.
          Since $\AAA_A = \Aplus_A[1/T]$, we find 
	  that $\gN \subseteq \gM[1/T] = \bigcup_{n = 0}^{\infty} T^{-n}\gM.$
	  Since $\gN$ is finitely generated as an $\A_{A}^+$-module, 
	  we obtain that $\gN \subseteq T^{-n}\gM$ for some sufficiently
	  large value of $n$, as required. 

	  To prove~(2), we note that if $A$ is Noetherian, 
	  than $\Aplus_A$ is also Noetherian, and thus any finitely generated
	  $\Aplus_A$-module is $T$-adically complete, since $\A^+_A$ itself is.
	  If $M$ is projective, then we write $M$ as a direct summand 
	  of a finitely generated free $\AAA_A$-module,
	  so that $\gM$ may then be embedded
	  into a finitely generated free $\Aplus_A$-module.   Thus $\gM$
	  is $T$-adically separated, and hence the kernel of any
	  surjection $(\Aplus_A)^r \to \gM$ is $T$-adically closed.  Since
	  $(\Aplus_A)^r$ is $T$-adically complete, we conclude that
	  the same is true of $\gM$.

	  Suppose now that $A$ is Noetherian, and that we are in
	  the situation of~(3).  The inclusion $T^n \gM \subseteq \gN$
	  shows that $M = \gM[1/T] \subseteq \gN[1/T],$
	  so that $\gN$ generates $M$ as an $\AAA_A$-module.
	  We must show that $\gN$ is furthermore finitely generated
	  over $\Aplus_A$. For this, we first note that, since 
	  $\gM$ is $T$-adically complete (by (2)),
	  the inclusion $T^n \gM \subseteq \gN \subseteq T^{-n}\gM$
	  shows that $\gN$ is also $T$-adically complete.
	  It also shows that $\gN/T\gN$ is a subquotient
	  of the finitely generated $A$-module $T^{-n}\gM/T^{n+1}\gM$,
	  and thus is a finitely generated $A$-module (as $A$ is Noetherian).
	  Since $\gN$ is $T$-adically complete, an application
	  of the topological Nakayama lemma shows that~$\gN$ is finitely
	  generated over $\Aplus_A$, as required.
          %
  \end{proof}

    \begin{arem}
    \label{rem: lattice and Drinfeld}By~\cite[Thm.\
    5.1.14]{EGstacktheoreticimages} (a theorem of Drinfeld),
    if $A$ is an $\cO/\varpi^a$-algebra for some $a \geq 1$,
    then we may
    think of a finitely generated projective~$\AAA_A$-module as a Tate $A$-module $M$
    together with a topologically nilpotent automorphism~$T$. In this
    optic a lattice in our sense is precisely
    a lattice in the Tate module~$M$
    which is also an $\Aplus_A$-submodule (by definition,
    a lattice~$L\subseteq M$
    is an open submodule with the property that for every open
    submodule~$U\subseteq L$, the $A$-module $L/U$ is finitely
    generated), 
    as the following lemma shows.
  \end{arem}

  \begin{alemma}\label{lem:lattice equivalences}
    If $A$ is an $\cO/\varpi^a$-algebra for some $a \geq 1$,
    and if $\gM$ is an $\Aplus_A$-submodule of a finitely generated
    projective $\AAA_A$-module~$M$,
	  then the following conditions on $\gM$ are equivalent:
	  \begin{enumerate}
		  \item $\gM$ is a lattice in $M$	
	  {\em (}in the sense of Definition~{\em \ref{defn: lattice})}
	  \item $\gM$ is open in~$M$,
		  and for every $\Aplus_A$-submodule $U$ of $\gM$
		  which is open in~$M$, the quotient $\gM/U$
		  is finitely generated over~$A$.
	  \item $\gM$ is open in~$M$,
		  and for every $A$-submodule $U$ of $\gM$
		  which is open in~$M$, the quotient $\gM/U$
		  is finitely generated over $A$
		  {\em (}i.e.\ $\gM$ is a lattice in~$M$,
		  when $M$ is thought of as a Tate module
		  over $A$ in the sense of Remark~{\em \ref{rem: lattice
				  and Drinfeld}).}
  \end{enumerate}
  \end{alemma}
  \begin{proof}
	  Suppose that $\gM$ is a lattice in $M$,
	  in the sense of Definition~\ref{defn: lattice}.
	  If we choose a surjection $(\Aplus_A)^r \to \gM$ for some $n \geq 0,$
	  then the induced surjection $\AAA_A^r \to M$ is a continuous
	  and  open map,
	  and so $\gM$ is open in~$M$, since $(\Aplus_A)^r$ is open in~$\AAA_A^r$.
	  Furthermore, since the submodules $T^n (\Aplus_A)^r$ ($n \geq 0$)
	  form a neighbourhood basis of $0$ in~$(\Aplus_A)^r$,
	  we see that their images $T^n\gM$ form a neighbourhood
	  basis of $0$ in~$\gM$.  Thus if $U$ is any open
	  $A$-submodule of~$\gM$, then $\gM/U$ is a quotient
	  of $\gM/T^n\gM$ for some~$n~\geq~0$,
	  and thus is finitely generated over~$A$.
	  Thus~(1) implies~(3).

	  Clearly~(3) implies~(2), and so we suppose that~(2) holds.
	  Let $\gN$ be an $\Aplus_A$-submodule of $M$ that is a lattice
	  in the sense of Definition~\ref{defn: lattice},
		  so that $M = \gN[1/T] = \bigcup_{n \geq 0} T^{-n} \gN$.
		  Then $\gM =
		  \bigcup_{n \geq 0} \gM \cap T^{-n}\gN,$
		  and thus
		 $$ \gM/(\gM \cap \gN) = 
		  \bigcup_{n \geq 0} (\gM \cap T^{-n}\gN) /
		  (\gM\cap  \gN).$$
		  But $\gM/(\gM \cap \gN)$ is finitely generated over $A$
		  by assumption,
		  and thus $$\gM/\gM \cap \gN = (\gM\cap T^{-n}\gN) /
		  (\gM\cap \gN),$$
		  for some sufficiently large value of~$n$,
		  implying that $\gM\subseteq T^{-n} \gN.$
		  In particular $\gM$ is $T$-adically complete,
		  being open (and hence closed) in
		  the lattice~$T^{-n}\gN$,
		  which is $T$-adically complete
		  by Lemma~\ref{lem:lattice properties}~(2).
	  Since $T$ is an automorphism of $M$, 
	  we find that $T\gM$ is an open submodule of $\gM$,
	  so that $\gM/T\gM$ is finitely generated over~$A$.
	  As $\gM$ is $T$-adically complete, 
	  the topological Nakayama lemma implies that $\gM$
	  is finitely generated over $\Aplus_A$.

	  Since $T$ is topologically nilpotent and $\gN$ is finitely
	  generated over $\Aplus_A$, we also find that $T^n\gN \subseteq 
	  \gM$ for some sufficiently large value of $n$.  Thus
	  $\gM[1/T] =\gN[1/T] = M.$  This completes the proof
	  that $\gM$ is a lattice in $M$
	  in the sense of Definition~\ref{defn: lattice},
	  showing that~(2) implies~(1).
  \end{proof}
  
  \begin{arem}
	  \label{rem:lattice cofinality}
	  It follows from Lemmas~\ref{lem:lattice properties}
	  and~\ref{lem:lattice equivalences}
	  that if $A$ is a~$\cO/\varpi^a$-algebra,
	  and $M$ is a finitely generated projective~$\AAA_A$-module,
	  then the lattices in $M$ form a neighbourhood basis
	  of the origin in~$M$.
  \end{arem}

  We also note the following technical lemma.

  \begin{alemma}
	  \label{lem:intersecting with lattice}
	  If $A\subseteq B$ is an inclusion of $\cO/\varpi^a$-algebras
	  for some $a \geq 1$,
	  with $A$ Noetherian,
	  if $M$ is a projective $\AAA_A$-module~$M$
	  with extension of scalars $M_B:=\A_B\otimes_{\A_A}M$ to $B$,
	  and if $\gM_B$ is a lattice in $M_B$,
	  then $\gM:=M \cap \gM_B$
	  {\em (}the intersection takes place in $M_B${\em )}
	  is a lattice in $M$,
	  having the additional property that $M \cap T^n \gM_B = T^n\gM$
	  for any integer~$n$.
  \end{alemma}
  \begin{proof}
	  Choose a projective complement to $M$,
	  i.e.\ a finitely generated projective $\AAA_A$-module
	  $N$ such that $M\oplus N$ is free over $\AAA_A$,
	  and let $\gN_B$ be a lattice in~$N_B$.
	  Then it suffices to prove the lemma with
	  $M$ replaced by $M\oplus N$ and with $\gM_B$
	  replaced by $\gM_B \oplus \gN_B$;
	  thus we may and do
	  suppose that $M$ is free over $\AAA_A$.
	  Choose a lattice $\gM'$ in $M$ which
	  is free over $\Aplus_A$,
	  and note that $\gM'_B := \Aplus_B\gM'$ is a lattice in $M_B$,
	  and that $M \cap T^n\gM'_B  = T^n\gM'$ for any integer $n$.

	  Lemma~\ref{lem:lattice properties}~(1) shows that
	  $T^n \gM'_B \subseteq \gM_B \subseteq T^{-n} \gM'_B$
	  for some sufficiently large value of $n$.
	  Intersecting with $M$, we find that
	  $T^n \gM' \subseteq \gM \subseteq T^{-n} \gM',$
	  so that $\gM$ is a lattice in $M$,
	  by Lemma~\ref{lem:lattice properties}~(3).
	  Finally,
	  by the definition of $\gM$,
	  we find that $M\cap T^n\gM_B = T^n\gM$ for any integer $n$.
  \end{proof}



\section{Group actions}
We now prove some lemmas which allow us to check whether the action of
a topological group on a Tate module is continuous, and to extend such
an action from a group to its completion.
\begin{alemma}\label{lem: continuity and lattices}  Let $G$ be a
  topological group acting on a Tate module $M$, and assume that~$G$
  admits a neighbourhood basis of the identity consisting of open
  subgroups. Then the following are equivalent:
	\begin{enumerate}
		\item The action $G\times M \to M$ is continuous.
		\item
                  \begin{enumerate}
                  \item For each~$m\in M$, the map $G\to M$, $g\mapsto
                    gm$ is continuous at the identity of~$G$,
                  \item for each~$g\in G$, the map $M\to M$, $m\mapsto
                    gm$ is continuous, and
                  \item for any lattice $L$ in $M$,
	 there is an open subgroup $H$ of $G$ that preserves~$L$.
 \end{enumerate}
                  \end{enumerate}
Furthermore, if these equivalent conditions hold, then the following stronger 
form of condition~{\em (2)(a)} holds:
\begin{enumerate}
\item[(2)(a')]
For each~$m \in M$, the map $G\to M$, $g \mapsto gm$ is equicontinuous
at the identity of~$G$.
\end{enumerate}
\end{alemma}
\begin{arem}
  \label{rem: neighbourhood basis profinite}The assumption in Lemma~\ref{lem: continuity and lattices} that~$G$
  admits a neighbourhood basis of open subgroups holds in particular
  if~$G$ is a profinite group, or if~$G$ is a subgroup of a profinite
  group with the subspace topology. In particular, it holds for the
  groups~$\Zp$ and for $\Z\subset\Z_p$.
\end{arem}
\begin{proof}[Proof of Lemma~\ref{lem: continuity and lattices}]
	Suppose first that~(1) holds; then clearly conditions~(2)(a)
	and~(2)(b) hold.  
	Let $L$ be a lattice in $M$.
	Since the action morphism $G \times M \to M$ is continuous,
	we may find an open subgroup $H$ of $G$ and an open submodule
	$U$ of~$M$, such that $H  U \subseteq L.$  Replacing $G$ by~$H$,
	we may thus suppose that $G U \subseteq L.$
	Now consider the morphism
	$G \times L/U \to M/L$  induced by the action morphism.
	Since $L$ is a lattice, the $A$-module $L/U$ is finitely generated.
	Let $\{m_i\}$ be a finite generating set.
	Since $M/L$ is discrete, for each generator~$m_i$,
	there is an open subgroup $H_i$ of $G$ such that $H_i m_i$ is
	constant, and thus equal to zero, modulo $L$.
        If we write $H := \bigcap_i H_i,$ then $H$ is an open subgroup
	of $G$ such $H m_i \subseteq L$ for every~$i$.  Consequently,
	$H L \subseteq L,$ verifying that~(2)(c) holds as well.

	Suppose conversely that~(2) holds. 
	We first note that we may strengthen~(2)(a) as follows:
	if $m \in M$, then the orbit map $g \mapsto gm$ is a continuous
	map $G \to M$.  Indeed, if $g_0 \in G$, then we may
	write this map as the composite of the continuous automorphism
	$g \mapsto g g_0^{-1}$ of $G$ and the orbit map $g \mapsto 
	g g_0 m.$  The latter orbit map is continuous at the identity
	of~$G$, by~(2)(a), and so the orbit map of $m$ is continuous
	at~$g_0$.

	We now wish to prove
	that the action map $G\times M \to M$ is continuous.
	Let $(g,m) \in G\times M.$   
	Any neighbourhood of the image $g m \in M$ contains
	a neighbourhood of the form $g m + L,$ where $L$ is a
	lattice in $M$.
	Given a lattice $L$, then, 
	we must find a neighbourhood of $(g,m)$ whose image
	lies in $g m + L.$

	Since $g$ is a continuous automorphism of $M$, by~(2)(a),
	we may find a lattice $L'$ such that $g L' \subseteq L.$
	By~(2)(c), we may find an open subgroup $H'\subseteq G$ that
	preserves~$L'.$  Then $H'$ acts on~$M/L'$, and since this quotient
	is discrete, and since the orbit maps for the action are continuous
	(by what we proved above),
        it follows that the action of $H'$ on $M/L'$ is smooth, 
	and so we may find an open subgroup $H$ of $H'$ that
	fixes the image of $m$ in $M/L'$.   Then we find that
	$(gH)(m+L') \subseteq gm + L,$ and since~$gH\times (m+L')$ is an
        open neighbourhood of $(g,m)\in G\times M$, 
        we have proved the required continuity.

Finally, suppose that conditions~(1) and~(2) hold; we must show that
the continuity condition of~(2)(a) 
can be upgraded  to the equicontinuity condition of~(2)(a'). 
To do this, we need to show that for any two
 lattices $L,L'\subseteq M$, there is an open subgroup~$H$ of~$G$ such
 that for each $m\in L$, we have $Hm\subseteq m+L'$. 

 Replacing~$L'$ by~$L\cap L'$, we can assume that~$L'\subseteq
 L$. By~(2)(c), we can choose~$H$ such that $H L'\subseteq L'$ and
 $H L\subseteq L$, so that we have an induced continuous map
 $H\times (L/L')\to L/L'$. Since~$L/L'$ is a finitely generated
 $A$-module and is discrete, after replacing~$H$ by an open subgroup we can assume
 that $H$ acts trivially on~$L/L'$, as required.
\end{proof}

\begin{alem}\label{lem: extending
  Z to Zp}
	Let $G$ be a Hausdorff topological group, 
  which admits a neighbourhood basis of the identity consisting of open
  subgroups, and suppose furthermore that for any open subgroup $H$ of~$G$,
		  the quotient $G/H$ is finite; equivalently,
		  suppose that the completion $\widehat{G}$ of $G$
		  is profinite.

		  If $M$ is a Tate module,
		  endowed with a continuous action $G \times M \to M$,
		  then 
                  the action $G\times M \to M$ extends
		       to a continuous action $\widehat{G} \times M
		       \to M$.
 \end{alem}
 \begin{proof}
We have to show that  the continuous action
 $G\times M \to M$  extends
 (necessarily uniquely, since $G$ is dense in $\widehat{G}$)
 to a continuous map $\widehat{G}\times M \to M$.  (The fact
 that this map will induce a $\widehat{G}$-action on $M$ follows
 from the corresponding fact for $G$, and the density of $G$ in $\widehat{G}$.)
 Since $M$ is the union of its lattices, it suffices to show
 that the induced map $G\times L \to M$ extends to a continuous
 map $\widehat{G} \times L \to M$, for each lattice $L~\subseteq~M$.
 For this, it suffices to show that the map $G\times L \to M$
 is uniformly continuous, for each lattice~$L$.
 That is, for any lattice $L' \subseteq M$, we have to find
 an open subgroup $H \subseteq G$ and a sublattice $L'' \subseteq L$
 such that $gH m + H L'' \subseteq gm + L'$ for all $g \in G$ and all
 $m \in L$.  

By Lemma~\ref{lem: continuity and lattices}, the conditions (2)(a)-(c)
of that lemma hold. We begin by showing that we may find a sublattice $L''$ of $L$
 such $G L'' \subseteq L'$.  For this, we first note that, by~(2)(c),
 we may find an open subgroup $H \subseteq G$ such
 that $H L' = L'$.  If let $\{g_i\}$ denote a (finite!) set
 of coset representatives for $H\backslash G,$ then
 since each $g_i$ induces a continuous automorphism of $M$,
 we may find a lattice $L_i$ such that $g_i L_i \subseteq L'$.  
 Taking into account Remark~\ref{rem:lattice cofinality},
 we may then find a lattice $L'' \subseteq L' \cap \bigcap_i L_i$,
 and by construction $G L'' \subseteq L'$.

 Condition~(2)(a') allows us to choose an open subgroup $H \subseteq G$
 such that $H m \subseteq m + L''$ for all $m \in L$.
 We then find that
 $$gH m + H L'' \subseteq gm + G L'' \subseteq gm + L'$$
 for all $g \in G$, as required.  \end{proof}





\section{$T$-quasi-linear endomorphisms}
In order to apply the preceding results to the case of interest to us
(the semilinear action of~$\Gamma$ on $(\varphi,\Gamma)$-modules), we
now introduce and study the notion of a $T$-quasi-linear endomorphism
of a finite projective $\A_A$-module. The relevance of this notion to
$(\varphi,\Gamma)$-modules is explained in
Lemma~\ref{lem:gamma minus 1 T quasi linear} below.


\begin{adefn}
If $M$ is a finite projective $\A_A$-module, 
then a {\em $T$-quasi-linear} endomorphism of $M$
is a morphism $f: M \to M$ 
which is $W(k)\otimes_{\Z_p} A$-linear, 
and which furthermore satisfies the following ($T$-quasi-linearity)
condition:
there exist power series $a(T) \in (\A_A^+)^{\times}$ and $b(T) \in (p,T)\A_A^+$
such that 
$$f(T m) = a(T) T f(m) + b(T)T m$$
for every $m \in M$.
\end{adefn}


\begin{alem}
  \label{lem: powers of T quasilinear}If $f$ is $T$-quasi-linear, then
  for all~$n\in\Z$, we may
  write \[f(T^nm)=a(T)^nT^nf(m)+b_n(T)T^{n}m\] for all~$m\in M$, where $a(T) \in (\A_A^+)^{\times}$ and $b_n(T) \in (p,T)\A^+_A$.
\end{alem}
\begin{proof}Note that the case~$n=0$ is trivial (taking $b_0(T)=0$),
  while if the claim holds for some~$n\ge 1$, then we may write
  \[f(T^{-n}m)=a(T)^{-n}T^{-n}f(m)-a(T)^{-n}b_n(T)T^{-n}m. \]
It therefore suffices to prove the result for ~$n\ge 1$.  This may be proved by induction on~$n$, the case $n=1$ being the
  definition of $T$-quasi-linearity. Indeed, if the claim holds
  for~$n$, then we have
  \begin{align*}
    f(T^{n+1}m)=f(T^n(Tm))&=a(T)^nT^nf(Tm)+b_n(T)T^{n}(Tm)\\
                          &=a(T)^nT^n(a(T) T f(m) + b(T)T
                            m)+b_n(T)T^{n+1}m\\  &=a(T)^{n+1}T^{n+1} f(m) +(b_n(T)+b(T)a(T)^n)T^{n+1}m,
  \end{align*}as required.\end{proof}

\begin{alem}\label{lem: how f behaves on lattices} 
Let~$A$ be an $\cO/\varpi^a$-algebra for some~$a\ge 1$, and let $M$ be a finite projective $\A_A$-module.  Let $f$ be a $T$-quasi-linear endomorphism
  of~$M$, and let~$\gM$ be a lattice in~$M$. Then  there is an
  integer~$m\ge 0$ such that for each~$s\in\Z$ and~$n\ge 0$ we have
  $f^n(T^s\gM)\subseteq T^{s-mn}\gM$.
\end{alem}
\begin{proof}
  Choose a finite
  set~$\{m_i\}$ of generators for~$\gM$ as an
  $\A^+_A$-module. Choose~$m$ sufficiently large that we
  have~$T^mf(m_i)\in \gM$ for each~$i$. Then by Lemma~\ref{lem: powers
    of T quasilinear} and the $W(k)\otimes_{\Zp}A$-linearity of~$f$, we
  see that for each~$s\in \Z$ we have $f(T^{s+m}\gM)\subseteq T^s\gM$,
  which gives the result in the case~$n=1$. The general case follows
  by induction on~$n$.
\end{proof}

As noted in Remark~\ref{rem:topologies on M},
any finite projective 
$\A_A$-module $M$ has a natural topology,
so it makes sense to speak of an endomorphism of $M$ being
continuous, or topologically nilpotent. 
\begin{alemma}\label{lem: quasilinear implies continuous}
	If~$A$ is an $\cO/\varpi^a$-algebra for some~$a\ge 1$, and $M$ is a finite projective $\A_A$-module,
	 then any $T$-quasi-linear endomorphism
	of $M$ is necessarily continuous.
\end{alemma}
\begin{proof}This is immediate from Lemma~\ref{lem: how f behaves on
    lattices}.
\end{proof}

\begin{alemma}
	\label{lem:topological nilpotence criteria}
	If~$A$ is an $\cO/\varpi^a$-algebra for some~$a\ge 1$, and $M$ is a finite projective $\A_A$-module, and if $f$ is a
	$T$-quasi-linear endomorphism of $M$,
	then the following are equivalent:
	\begin{enumerate}
		\item $f$ is topologically nilpotent.
		\item There exists a lattice $\gM$ in $M$
			and some $n \geq 1$
			such that $f^n(\gM) \subseteq T \gM.$
                        	\item There exists a lattice $\gM$ in $M$
			and some $n \geq 1$
			such that $f^n(\gM) \subseteq (p,T) \gM.$
                                              \item There exists a lattice $\gM$ in $M$ such
                        that for any $m \geq 1$, there exists~$n_0$
			such that  for any~$s\in\Z$ and any $n\ge
                        n_0$, we have $f^n(T^{s}\gM) \subseteq T^{s+m}
                        \gM$.
  \item For any lattice $\gM$ in $M$ and any $m \geq 1$, there exists~$n_0$
			such that  for any~$s\in\Z$ and any $n\ge
                        n_0$, we have $f^n(T^{s}\gM) \subseteq T^{s+m} \gM$.
	\end{enumerate}
\end{alemma}
\begin{proof}By Lemma~\ref{lem:lattice properties}~(1), we see
  that~$(4)\implies (5)\implies (1)\implies (2)\implies (3)$, so we
  only need to that~$(3)\implies (4)$. To this end, note firstly that
  it follows from Lemma~\ref{lem: powers of T quasilinear} and the
  $\Zp$-linearity of~$f$ that for
  each $i,j\ge 0$ and $s\in\Z$ we have
  \[f((p,T)^jT^sf^i(\gM))\subseteq (p,T)^{j}T^sf^{i+1}(\gM)+(p,T)^{j+1}T^sf^i(\gM). \] It
  follows (by induction on~$m$)
  that for each~$m\ge 0$ and $s\in\Z$ we have
  \[ f^m(T^s\gM)\subseteq \sum_{i=0}^m(p,T)^{m-i}T^sf^i(\gM).\] In
  particular, if we take $m=n$ and recall that
  $f^n(\gM)\subseteq(p,T)\gM$ by hypothesis, we find that
  \anumequation\label{eqn: base case of quasilinear induction}f^n(T^s\gM)\subseteq (p,T)T^s\gM+\sum_{i=1}^{n-1}(p,T)^{n-i}T^sf^i(\gM). \end{equation}

The same argument shows that if~$\gN$ is an $\A_A^+$-submodule of~$M$ with
  the property that
  \[\gN\subseteq \sum_{i=0}^{n-1}(p,T)^{a_i}T^sf^i(\gM)\]
      for non-negative integers $a_0,\dots,a_{n-1}$, then 
  \[f(\gN)\subseteq  \sum_{i=0}^{n-1}(p,T)^{b_i}T^sf^i(\gM)\]   
  where~$b_0=\min(a_0+1,a_{n-1}+1)$, and~$b_i=\min(a_i+1,a_{i-1})$
  if~$i>0$.  
  It follows by an easy induction on~$N$ (with the base case being
  given by~\eqref{eqn: base case of quasilinear induction}) that for all~$N\ge n$, if we
  write~$N+1=(q+1)n+r$ with $0\le r<n$, then we
  have
  \[f^N(T^{s}\gM)\subseteq \sum_{i=0}^{n-1}(p,T)^{c_i}T^sf^i(\gM)\]where
  \[(c_{n-1},\dots,c_0)=(q,q,\dots,q)+(r,r+1,\dots,n-1,1,2,\dots,r). \]
In particular we see that for all $N\ge n$ we
have   \[f^N(T^{s}\gM)\subseteq (p,T)^{\lfloor (N+1)/n\rfloor-1}T^s\sum_{i=0}^{n-1}f^i(\gM)\]

  Now, by Lemma 
  ~\ref{lem: how f
    behaves on lattices}, for any sufficiently large~$t$ we have
 ~$f^i(\gM)\subseteq T^{-t}\gM$ for $0\le i\le n-1$, so it follows
  that for $N\ge n$ we
have   \[f^N(T^{s}\gM)\subseteq (p,T)^{\lfloor
    (N+1)/n\rfloor-1}T^{s-t}\gM.\] Since~$p^a=0$ in~$A$, we also have
$(p,T)^{n+a-1}\subseteq (T^n)$ for all~$n\ge 0$, so that if $N\ge
an-1$ 
then we have \[f^N(T^{s}\gM)\subseteq T^{s+\lfloor
    (N+1)/n\rfloor-a-t}\gM.\] Since $\lfloor
    (N+1)/n\rfloor-a-t\to\infty$  as $N\to\infty$, we have~(4), as required.
\end{proof}


\begin{acor}
\label{cor: variants on topological nilpotence mod p}
Suppose that $A$ is an $\cO/\varpi^a$-algebra for
        some~$a\ge 1$.
	If $M$ is a finite projective $\A_A$-module, and if $f$ is a
	$T$-quasi-linear endomorphism of $M$,
	then the following are equivalent:
	\begin{enumerate}
		\item $f$ is topologically nilpotent.
		\item The action of~$f$ on ~$M\otimes_{\cO/\varpi^a}\F$
          is topologically nilpotent.
	\end{enumerate}
      \end{acor}
      \begin{proof}Obviously $(1)\implies (2)$. Conversely, if (2)
        holds, then by the equivalence of conditions (1) and (2) of
        Lemma~\ref{lem:topological nilpotence criteria} for the action
        of~$f$ on ~$M\otimes_{\cO/\varpi^a}\F$, we see that condition
        (3) of 
        Lemma~\ref{lem:topological nilpotence criteria} holds (for the
        action of~$f$ on~$M$), and therefore condition
        (1) of 
        Lemma~\ref{lem:topological nilpotence criteria} holds, as required.        
      \end{proof}

The following lemmas provide the key examples of $T$-quasi-linear
endomorphisms, and explain our interest in the concept.

We suppose that~$\A_A$ is endowed with a continuous action
of~$\Zp$ by $A$-algebra automorphisms which preserve~$\A^+_A$,
and suppose further that, for some 
topological generator $\gamma$  of~$\Zp$, that
\anumequation\label{eqn: condition on gamma semi-linear}\gamma(T)-T\in(p,T)T\A^+_A.\end{equation} 

\begin{alemma}
  \label{lem: all powers of gamma preserve the ideal}If~\eqref{eqn:
    condition on gamma semi-linear} holds, then $\gamma$ preserves the
  ideals~$(T)$ and~$(p,T)$ of~$\A^+_A$. Furthermore for each
  integer~$n\ge 1$ we have
  \anumequation\label{eqn: condition on powers of gamma semi-linear}\gamma^n(T)-T\in(p,T)T\A^+_A.\end{equation}
\end{alemma}
\begin{proof}The first claim follows immediately from~\eqref{eqn:
    condition on gamma semi-linear}, which shows that~$\gamma(T)$ is a
  unit multiple of~$T$, together with the fact that  $\gamma$  
preserves $\A^+_A$. Then~\eqref{eqn: condition on powers of gamma semi-linear}
    follows by induction on~$n$ (the case~$n=1$
  being~\eqref{eqn: condition on gamma semi-linear}).
\end{proof}

\begin{alemma}
	\label{lem:gamma minus 1 T quasi linear}
If
        $M$ is a finite projective $\A_A$-module which
	is endowed with an action of the subgroup  $\langle
\gamma \rangle$  of $\Zp$ which  is semi-linear with  respect
to the given action of this group on~$\A_A$ {\em(}obtained
by restricting the $\Z_p$-action{\em )}, 
	then for any integer~ $n\ge 1$, $f := \gamma^n -1$ is a $T$-quasi-linear endomorphism of~$M$.
\end{alemma}
\begin{proof}
  We have $f(Tm)=\gamma^n(T)f(m)+(\gamma^n(T)-T)m$, so it follows from~\ref{eqn: condition on powers of gamma semi-linear} that~$f$ is $T$-quasi-linear.
\end{proof}

\begin{alem}
	\label{lem:testing continuity on M mod T}
        Suppose that $A$ is an $\cO/\varpi^a$-algebra for
        some~$a\ge 1$, and that~$\A_A$ is endowed with an action of~$\Zp$
        satisfying~\eqref{eqn: condition on gamma semi-linear}. Let~$M$
        be a finite projective $\A_A$-module, equipped with a
        semi-linear action of~$\langle \gamma\rangle\subset \Zp$. Then
        the following are equivalent:
        \begin{enumerate}
        \item\label{item: gamma cts} The action of~$\langle \gamma\rangle$ extends to a continuous
          action of~$\Zp$.
        \item\label{item: gamma cts mod p} The action of~$\langle \gamma\rangle$ on~$M\otimes_{\cO/\varpi^a}\F$ extends to a continuous
          action of~$\Zp$.
        \item\label{item: gamma minus 1 squeezing M to p T n} For any lattice~$\gM\subseteq M$, and any~$n\ge 1$, there exists~$s\ge 0$
          such that $(\gamma^{p^s}-1)^i(\gM)\subseteq (p,T)^n\gM$ for
          all~$i\ge 1$.
          \item\label{item: gamma minus 1 squeezing M to T} For any lattice~$\gM\subseteq M$, there exists~$s\ge 0$
          such that $(\gamma^{p^s}-1)(\gM)\subseteq T\gM$.
           \item\label{item: gamma minus 1 squeezing some M to p T} For some lattice~$\gM\subseteq M$ and some~$s\ge 0$, we have
          $(\gamma^{p^s}-1)(\gM)\subseteq (p,T)\gM$.
        \item\label{item: gamma minus 1 top nilpt mod p}  The action of~$\gamma-1$ on~$M\otimes_{\cO/\varpi^a}\F$
          is topologically nilpotent.
          \item\label{item: gamma minus 1 top nilpt}  The action of~$\gamma-1$ on~$M$
          is topologically nilpotent.
        \end{enumerate}
  \end{alem}
  \begin{proof}

        Noting that if~$n\ge a$ then $(p,T)^n\subseteq (T)$, we see that
        \eqref{item: gamma minus 1 squeezing M to p T
          n}$\implies$\eqref{item: gamma minus 1 squeezing M to T},
        and by Lemma~\ref{lem:topological nilpotence criteria}
        and Corollary~\ref{cor: variants on topological nilpotence mod p} we
        have \eqref{item: gamma minus 1 squeezing M to
          T}$\implies$\eqref{item: gamma minus 1 squeezing some M to p
          T}$\implies$\eqref{item: gamma minus 1 top nilpt mod
          p}$\implies$\eqref{item: gamma minus 1 top nilpt}. (For \eqref{item: gamma minus 1 squeezing some M to p
          T}$\implies$\eqref{item: gamma minus 1 top nilpt mod
          p}, we also use that
        $(\gamma-1)^{p^s}\equiv(\gamma^{p^s}-1)\pmod{\varpi}$.)
        
        We next show that \eqref{item: gamma minus 1 top
          nilpt}$\implies$\eqref{item: gamma minus 1 squeezing M to p
          T n}. Suppose that~\eqref{item: gamma minus 1 top
          nilpt} holds. Recalling again that for each~$s\ge 0$ we have
        $(\gamma-1)^{p^s}\equiv(\gamma^{p^s}-1)\pmod{\varpi}$, it
        follows from Corollary~\ref{cor: variants on topological
          nilpotence mod p} and Lemma~\ref{lem:gamma minus 1 T quasi linear} that the action of $(\gamma^{p^s}-1)$
        on~$M$ is topologically nilpotent for each~$s\ge 0$. 
        We now
        argue by induction on~$a$, noting that for $a=1$, the
        implication \eqref{item: gamma minus 1 top
          nilpt}$\implies$\eqref{item: gamma minus 1 squeezing M to p
          T n} is immediate from Lemma~\ref{lem:topological nilpotence
          criteria}. We may therefore assume
        that \[(\gamma^{p^s}-1)^i(\gM)\subseteq (p,T)^n\gM
          +\varpi^{a-1}M \]for all~$i\ge 1$. It follows in particular
        that $p(\gamma^{p^s}-1)^i(\gM)\subseteq (p,T)^n\gM$ for
        all~$i\ge 1$, so that for any~$t\ge 0$ and any~$i\ge 1$, we
        have \[((\gamma^{p^{s+t}}-1)^i-(\gamma^{p^s}-1)^{ip^t})(\gM)\subseteq
          (p,T)^n\gM.\](To see this, write
        $(\gamma^{p^{s+t}}-1)=((\gamma^{p^s}-1)+1)^{p^t}-1$ and use the
        binomial theorem.) It therefore suffices to show that there is
        some~$t\ge 1$ for which~$(\gamma^{p^s}-1)^{ip^t}(\gM)\subseteq
        (p,T)^n\gM$ for all~$i\ge 1$; but we have already seen
        that~$(\gamma^{p^s}-1)$ acts topologically nilpotently on~$M$,
        so by Lemma~\ref{lem:topological nilpotence criteria} we can
        even arrange that $(\gamma^{p^s}-1)^{ip^t}(\gM)\subseteq
        T^n\gM$ for all~$i\ge 1$.


        We have shown the equivalence of conditions \eqref{item:
          gamma minus 1 squeezing M to p T n}--\eqref{item: gamma
          minus 1 top nilpt}. Suppose now that~\eqref{item: gamma cts} holds. Then Lemma~\ref{lem: continuity and
          lattices} shows that, for each lattice $\gM\subseteq M$,
we have~$\gamma^{p^s}(\gM)\subseteq\gM$
for all sufficiently large~$s\ge 0$,
and that furthermore, for each
        $m\in\gM$, there is some~ $t(m)$ such that if~$t\ge t(m)$, then
        $\gamma^{p^{t}}(m)\in m+(p,T)\gM$. Letting~$m_1,\dots,m_n$
        be generators for~$\gM$ as an $\A_{A}^+$-module, we see that
        if~$s$ is sufficiently large, then $(\gamma^{p^s}-1)(m_i)\in
        (p,T)\gM$ for each~$i$. Then if~$\lambda_i\in \A_{A}^+$, we
        have \[(\gamma^{p^s}-1)(\sum_i
          \lambda_im_i)=\sum_i\gamma^{p^s}(\lambda_i)(\gamma^{p^s}-1)(m_i)+\sum_i(\gamma^{p^s}-1)(\lambda_i)m_i.\]It
        follows from~\eqref{eqn: condition on powers of gamma semi-linear}
        that~$(\gamma^{p^s}-1)(\lambda_i)\in T\A_{A}^+$; hence $(\gamma^{p^s}-1)(\gM)\subseteq
        (p,T)\gM$, and so~\eqref{item: gamma minus 1 squeezing some M to p T} holds.

Suppose now that the equivalent conditions \eqref{item:
          gamma minus 1 squeezing M to p T n}--\eqref{item: gamma
          minus 1 top nilpt} hold. We will show 
that~\eqref{item: gamma cts} holds.  By Lemma~\ref{lem: extending Z to
          Zp} (taking the group $G$ there to be $\langle \gamma \rangle \cong \Z$
endowed with its $p$-adic topology, so that $\widehat{G}~\cong~\Z_p$),
it is enough to show that  the conditions of
Lemma~\ref{lem: continuity and lattices}~(2) hold. 

We begin with Lemma~\ref{lem: continuity and lattices}~(2)(b), the
condition that any~$g\in \langle \gamma\rangle$ acts continuously
on~$M$. It is enough to show that
if~$\gM\subseteq M$ is a lattice, then there
is a lattice~$\gN$ with $g(\gN)\subseteq \gM$.  This is in fact a general
property of semilinear automorphisms of~$\A_A$ which
preserve~$\A^+_A$. 
Indeed, let~$m_1,\dots,m_n$ be
generators of~$\gM$ as an~$\A_{A}^+$-module, and let~$\gN$ be the
$\A_{A}^+$-module generated by $g^{-1}(m_1),\dots,g^{-1}(m_n)$. This
is a lattice, because~$g$ is an automorphism of~$M$, and it
follows easily from the semi-linearity of the action of~$g$ on~$M$ that
$g(\gN)\subseteq\gM$,  as required.

We now check Lemma~\ref{lem: continuity and lattices}~(2)(c). Let~$\gM\subseteq M$ be some lattice, fix a choice of $n\ge 1$,
and then choose~$s$ as in~\eqref{item: gamma minus 1 squeezing M to p
  T n}. It suffices to show that the subgroup~$H=\langle \gamma^{p^s}\rangle$ 
of~$\langle \gamma\rangle$ preserves~$\gM$. Since $(\gamma^{p^s}-1)(\gM)\subseteq
(p,T)^n\gM\subseteq\gM$, we certainly
have~$\gamma^{p^s}(\gM)\subseteq\gM$, so it suffices to show that
$\gamma^{-p^s}(\gM)\subseteq\gM$. For this, note that
since (as recalled above) 
it follows from the congruence 
$(\gamma-1)^{p^s}\equiv(\gamma^{p^s}-1)\pmod{\varpi}$ that
$\gamma^{p^s}-1$ acts topologically nilpotently on~$M$, and
preserves~$\gM$, 
we see that for
any $m\in \gM$ we have
\[\gamma^{-p^s}(m)=(1-(1-\gamma^{p^s}))^{-1}(m)=m+(1-\gamma^{p^s})(m)+(1-\gamma^{p^s})^2(m)+\dots\in\gM, \]as
required.  

To complete the verification of the conditions of Lemma~\ref{lem:
  continuity and lattices}~(2), we need to show that for any
lattice $\gM~\subseteq~M$, 
the orbit maps $\langle \gamma \rangle \to M$, for the various $m \in \gM$,
are (equi)continuous at the identity of~$\langle \gamma \rangle$.
It is enough to show that for
each~$n\ge 1$, we can find~$s$ sufficiently large such that
$H=\langle \gamma^{p^s}\rangle$ satisfies~$Hm\subseteq m+(p,T)^n\gM$,
for each $m \in \gM$,
or equivalently,
that $(h-1)(\gM)\subseteq (p,T)^n\gM$, for all~$h\in H$. 
We accomplish this by choosing~$s$ as in~\eqref{item: gamma minus 1
  squeezing M to p T n}.

It then suffices to prove
the stronger claim that~$(\gamma^{rp^s}-1)(\gM)\subset(p,T)^n\gM$ for
all~$r\in\Z$. We have already seen that $\gamma^{rp^s}(\gM)\subseteq\gM$. 
If~$r\ge
1$ we may write
$(\gamma^{rp^s}-1)=(\gamma^{p^s}-1)(1+\gamma^{p^s}+\dots+\gamma^{(r-1)p^s})$,
and since $(1+\gamma^{p^s}+\dots+\gamma^{(r-1)p^s})(\gM)\subset\gM$,
we have~$(\gamma^{rp^s}-1)(\gM)\subseteq (\gamma^{p^s}-1)(\gM)\subset
(p,T)^n\gM$, as
required. If~$r\le 0$, then the result follows by writing~$(\gamma^{-rp^s}-1)=-(\gamma^{rp^s}-1)(\gamma^{-p^s})^r$.

        
        Finally, applying the equivalence of~\eqref{item: gamma cts} and~\eqref{item: gamma minus 1 top nilpt mod p} with~$M$
        replaced by $M\otimes_{\cO/\varpi^a}\F$, we see that~\eqref{item: gamma cts mod p}
        and~\eqref{item: gamma minus 1 top nilpt mod p} are equivalent, as required.\end{proof}




The following lemmas will allow us to reduce the problem of investigating
the topological nilpotency of a $T$-quasi-linear endomorphism
from the projective case to the free case.

\begin{alemma}
	\label{lem:quasi-linear direct sum}
	If $M_1$ and $M_2$ are Tate modules over a ring~$A$,
	endowed with continuous endomorphisms $f_1$ and $f_2$
	respectively, then $f_1$ and $f_2$ are both topologically
	nilpotent if and only if the direct sum $f := f_1 \oplus f_2$
	is a topologically nilpotent endomorphism.
\end{alemma}
\begin{proof}
This is immediate from the definitions.
\end{proof}

\begin{alemma}
	\label{lem:quasi-linear complement}
	If $M$ is a finite projective $\A_A$-module
	endowed with a $T$-quasi-linear endomorphism $f$,
	then we may find a finite projective $\A_A$-module
	$N$, endowed with a $T$-quasi-linear endomorphism
	$g$ which is furthermore topologically nilpotent,
	such that $M\oplus N$ is a free $\A_A$-module.
\end{alemma}
\begin{proof}
  By 
definition, there is a $\A_A$-module   $N$ such that $M \oplus N$ is
  free.  We may then take $g$ to be the zero endomorphism of~$N$; this
  is evidently both $T$-quasi-linear and topologically nilpotent.
\end{proof}
If $A$ is an $\cO/\varpi^a$-algebra for some~$a\ge 1$, and $M$ is a
finite projective $\A_A$-module, equipped with a $T$-quasi-linear endomorphism $f$,
and if $A \to B$ is a morphism of $\cO/\varpi^a$-algebras, then we have a
base-changed $W(k)\otimes_\Zp B$-linear endomorphism of~$B\otimes_AM$.
This base-changed endomorphism is 
$T$-adically continuous (by Lemma~\ref{lem: how f behaves on lattices}) 
and it therefore induces an
endomorphism $f_B$ of the base-changed projective $\A_B$-module $M_B$ (which by definition is  the completion of~$B\otimes_AM$ for
the $T$-adic topology). The endomorphism~$f_B$ is evidently also $T$-quasi-linear.
(Note that since~$f$ is not necessarily
$A[[T]]$-linear, we do not define~$f_B$ by viewing~$M_B$ as
$B[[T]]\otimes_{A[[T]]} M$.) 
If $f$ is furthermore topologically
nilpotent, then so is $f_B$ (for example, by Lemma~\ref{lem:topological nilpotence criteria}).

\begin{alemma}
	\label{lem:quasi-linear limit preserving}
	Let $\{A_i\}_{i \in I}$ be a directed system of $\cO/\varpi^a$-algebras,
	with $A = \varinjlim_{i \in I} A_i$, and suppose that $I$ admits
	a least element $i_0$.  Let $M$ be a finite projective $\A_{A_{i_0}}$-module,
	 let $f$ be a $T$-quasi-linear endomorphism of $M$,
	and suppose that the base-changed endomorphism $f_A$ is topologically
	nilpotent.  Then for some $i \in I,$ the base-changed
	endomorphism $f_{A_i}$ is topologically nilpotent.
\end{alemma}
\begin{proof}
	By Lemma~\ref{lem:quasi-linear complement},
	we may choose a finite projective  $\A_{A_{i_0}}$-module $N$,
	endowed with a topologically nilpotent $T$-quasi-linear
	endomorphism $g$, such that $M \oplus N$ is free.
	By Lemma~\ref{lem:quasi-linear direct sum}, the base-changed
        endomorphism $f_A \oplus g_A = (f\oplus g)_A$ of $M_A\oplus N_A$
        is topologically nilpotent, and it suffices to show 
	that $f_{A_i} \oplus g_{A_i} = (f\oplus g)_{A_i}$
	is topologically nilpotent
	for some $i \in I$.   Thus we may reduce to the case
	when $M$ is free, in which case we may also choose
	a free lattice $\gM$ contained in $M$.   Taking into account
	Lemma~\ref{lem:topological nilpotence criteria}, 
for some sufficiently large~	$n$ we have $f_A^n(\gM_A) \subseteq T\gM_A$, and it suffices to show that
	\anumequation
	\label{eqn:needed inclusion}
	f_{A_i}^n(\gM_{A_i}) \subseteq T \gM_{A_i}
\end{equation}
       	for some $i \in I$.
	Lemma~\ref{lem: how f behaves on lattices} shows that
	$f^n(T^r \gM) \subseteq T\gM$ for some $r\geq 0,$
       	and that $f^n(\gM) \subseteq T^{-s} \gM$
	for some $s \geq 0$.  Thus $f^n$ induces a morphism
	$\gM/T^r \gM \to T^{-s}\gM/T{\gM}$ of finite rank free
	$A_{i_0}$-modules,
	which, by assumption, vanishes after base-change to $A$.  
	Thus this morphism in fact vanishes after base-change
	to some $A_i,$ and consequently~(\ref{eqn:needed inclusion})
	does indeed hold for this choice of $A_i$.
\end{proof}

\chapter{Points, residual gerbes, and isotrivial
  families}\chaptermark{Points and residual gerbes}
\label{app:residual gerbes}
Recall that if $\cX$ is an algebraic stack, then the underlying set 
$|\cX|$ of points of $\cX$ is defined as the set of equivalence classes
of morphisms $\Spec K \to \cX$, with $K$ being a field; two such morphisms
are regarded as equivalent if they may be dominated by a common morphism
$\Spec L \to \cX$.  

If $X$ is a scheme (regarded as an algebraic stack), then the
set of underlying points $|X|$ is naturally identified with the
underlying 
set of points of $X$ in the usual sense,
and the equivalence class of 
morphisms representing a given point $x \in |X|$
has a canonical representative, namely the morphism $\Spec k(x) \to X$,
where $k(x)$ is the residue field of $x$.  More abstractly,
this morphism is a monomorphism, and this 
property characterizes it uniquely,
up to unique isomorphism,
among the morphisms in the equivalence class corresponding to~$x$.

If $X$ is an algebraic space, then it is not the case that every
equivalence class in $|X|$ admits a representative which is a monomorphism
(see
e.g.~\cite[\href{http://stacks.math.columbia.edu/tag/02Z7}{Tag
  02Z7}]{stacks-project}). 
However, under mild assumptions on $X$,
the points of $|X|$ do admit such representatives (which
are then unique up to unique isomorphism); this in particular
is the case if $X$ is quasi-separated.   See e.g.\
the discussion at the beginning
of~\cite[\href{http://stacks.math.columbia.edu/tag/03I7}{Tag 03I7}]{stacks-project}.  

If $\cX$ is genuinely an algebraic stack,
then it is not reasonable to expect the equivalence classes in $|\cX|$ 
to admit monomorphism representatives in general (even 
if $\cX$ is quasi-separated), since points of $\cX$ 
typically admit non-trivial stabilizers.  The aim of the theory
of residual gerbes is to provide a replacement for such representatives.

We recall the relevant
definition~\cite[\href{http://stacks.math.columbia.edu/tag/06MU}{Tag
  06MU}]{stacks-project}): 

\begin{adefn}\index{residual gerbe}
	Let $\cX$ be an algebraic stack.
	We say that the {\em residual gerbe} at a point $x \in |\cX|$ 
	exists if we may find a monomorphism
	$\cZ_x \hookrightarrow \cX$ such that $\cZ_x$ is reduced and
	locally Noetherian, $|\cZ_x|$ is a singleton,
	and the image of $|\cZ_x|$ in $|\cX|$ is equal to $x$.
	If such a monomorphism exists, 
        then $\cZ_x$ is unique up to unique isomorphism,
        and we refer to it as the {\em residual gerbe} at $x$.
\end{adefn}

In the case that $X$ is an algebraic space,
the residual gerbe $\cZ_x$ exists
at every point $x \in
|X|$~\cite[\href{http://stacks.math.columbia.edu/tag/06QZ}{Tag
  06QZ}]{stacks-project}, 
and is itself an algebraic space, called the {\em residual space} of
\index{residual space}
$X$ at~$x$~\cite[\href{http://stacks.math.columbia.edu/tag/06R0}{Tag
  06R0}]{stacks-project}. 
However, as already noted, in this case, if $X$ is furthermore
quasi-separated,
then the residual space at a point is simply the spectrum of a field.

If $\cX$ is an algebraic stack with quasi-compact diagonal
(e.g.\ if $\cX$ is quasi-separated),
then the residual gerbe exists at every point of
$|\cX|$~\cite[\href{http://stacks.math.columbia.edu/tag/06RD}{Tag
  06RD}]{stacks-project}). 

If $X$ is a scheme, then a {\em finite type point} of $X$
is a point $x \in X$ that is locally closed.  Equivalently,
these are the points for which the morphism $\Spec k(x) \to X$
is a finite type morphism, and may be characterized more abstractly
as those equivalence classes of morphisms $\Spec K \to X$ which
admit a representative which is locally of finite type.
This latter notion makes sense for an arbitrary algebraic stack,
and allows us to define the notion of a {\em finite type point} 
of an algebraic stack (see e.g.\ \cite[1.5.3]{EGstacktheoreticimages}).

We then have the following results, which provide analogues
of the topological characterization of the finite type points of a scheme
as being those points that are locally closed.

\begin{alemma}
	\label{lem:residual space immersion}
	Suppose that $X$ is an algebraic space,
	and that $x \in |X|$ is a finite type point with the property 
	that for any \'etale morphism $U \to X$ whose
	source is an affine scheme,  the fibre over $x$ 
	is finite.
	Then, if $Z_x$ denotes the residual space at $x$ in $X$,
	the canonical monomorphism $Z_x \hookrightarrow X$
	is an immersion.
\end{alemma}
\begin{proof}
	Consider the construction of $Z_x$ given
	in~\cite[\href{http://stacks.math.columbia.edu/tag/06QZ}{Tag
  06QZ}]{stacks-project}: we choose a surjective \'etale morphism
$U \to X$, form the union $$U' = \coprod_{u \in U \text{ lying
		over } x} \Spec k(u),$$
and then realize $Z_x \hookrightarrow X$ as a descent of the monomorphism
$U' \hookrightarrow U$.   In making this construction, we may replace
$U$ by any open subscheme containing $x$ in its image, 
and thus we may assume that $U$ is affine.  By assumption,
there are then only
finitely many points $u$ in $U$ lying over $x$,
these points are all finite type points of $U$,
and furthermore, none of these points are specializations of 
any of the
others~\cite[\href{http://stacks.math.columbia.edu/tag/03IM}{Tag
  03IM}]{stacks-project}.
Thus $U'\hookrightarrow U$ is in fact an immersion, and thus the same
is true of $Z_x \hookrightarrow X.$
\end{proof}

\begin{alemma}
	\label{lem:residual gerbe immersion}
	If $\cX$ is an algebraic stack whose diagonal is quasi-compact,
	if $x\in |\cX|$ is a finite type point,
	and if $\cZ_x$ is the residual gerbe of $\cX$ at $x$,
	then the canonical monomorphism
	$\cZ_x \hookrightarrow \cX$ is in fact an immersion.
\end{alemma}
\begin{proof}
	We prove this by examining the construction of
	$\cZ_x$ carried out in the proof 
        of~\cite[\href{http://stacks.math.columbia.edu/tag/06RD}{Tag
        06RD}]{stacks-project}. 
        The first step of the proof is to replace $\cX$
	by the closure of $x$ in $|\cX|$, regarded as a closed
	substack of $\cX$ with its induced reduced structure.
	The assumption that $\cX$ has quasi-compact diagonal
	implies that the morphism $\cI_{\cX} \to \cX$
	is quasi-compact, and thus that we may find a dense open substack
	$\cU$ of $\cX$ over which this morphism is flat and locally of
	finite presentation.
	Since $x$ is dense in $|\cX|$,
	we find that $x \in |\cU|$, and so replacing $\cX$
	by $\cU$, we may assume that $\cI_{\cX} \to \cX$ 
	is flat and locally of finite presentation.
        By~\cite[\href{http://stacks.math.columbia.edu/tag/06QJ}{Tag
        06QJ}]{stacks-project}, 
        this implies that $\cX$ is a gerbe over some algebraic space~$X$.
	If we let $Z_x$ denote the residual space at (the image in $|X|$ of) $x$
	in $|X|$, then $\cZ_x$ is obtained as the base-change of $Z_x$
	over the morphism $\cX \to X$.

	Note that, in the various reduction steps undertaken in the
	preceding argument, we replaced $\cX$ by an open substack
	of a closed substack; it thus suffices to verify
	that $\cZ_x \to \cX$ is an immersion after making these
	reductions.  Note that since immersions are monomorphisms,
	the diagonal of the stack obtained after these reductions
	are made is the base-change of the diagonal of the original
	stack $\cX$, and thus continues to be quasi-compact.
	Consequently, we may assume that we are in the case
	where $\cX$ is a gerbe over $X$.  Since $\cZ_x \hookrightarrow \cX$
	is then obtained as the base-change of $Z_x \hookrightarrow X,$
	it suffices to show that this latter morphism is an immersion.
	For this, if we take into account Lemma~\ref{lem:residual space
		immersion}, it suffices to show that $X$ can be chosen
	to be quasi-separated (as every \'etale morphism from an 
	affine scheme to a quasi-separated algebraic 
	space has finite fibres over every point of~$|\cX|$;
see the discussion at the beginning
of~\cite[\href{http://stacks.math.columbia.edu/tag/03I7}{Tag
  03I7}]{stacks-project}). 

If we examine the proof
        of~\cite[\href{http://stacks.math.columbia.edu/tag/06QJ}{Tag
        06QJ}]{stacks-project}, 
we see that $X$ is constructed as follows: We choose a smooth
surjective morphism $U \to \cX$ whose source is a scheme,
and write $R = U\times_{\cX} U,$ so that 
$\cX = [U/R]$.
We then factor the morphism $R \to U\times_{\Spec \Z} U$ 
through a morphism $R' \to U\times_{\Spec \Z} U$,
where $R' \to U\times_{\Spec \Z} U$
is a flat and locally of finite presentation equivalence relation,
and set $X = U/R'.$  
We note one additional aspect of the situation,
namely that the morphism $R\to R'$ is surjective.
The assumption that $\cX$ has quasi-compact diagonal then implies
that $R \to U\times_{\Spec \Z} U$ is quasi-compact, and since $R$ surjects
onto $R'$, we find that $R'\to U\times_{\Spec \Z} U$ is again quasi-compact.
Thus $X = U/R'$ is quasi-separated, as required.
\end{proof}

\begin{aexample}
	We note that quasi-separatedness (or some such hypothesis) is
	necessary for the truth of Lemma~\ref{lem:residual gerbe
	immersion}.  To illustrate this, let $G$ be an algebraic group
of positive dimension over an algebraically closed field $k$,
and let $X := G/G(k).$  Then $X$ contains a unique finite type
point $x$ --- namely, the equivalence
class of the monomorphism $\Spec k = G(k)/G(k) \hookrightarrow G/G(k)$
--- and this monomorphism realizes $\Spec k$ as the residual space
$Z_x$.  This monomorphism is {\em not} an immersion, 
since its pull-back to $G$ induces the monomorphism $G(k) \hookrightarrow
G,$ which is not an immersion.
\end{aexample}

Traditionally, in the theory of moduli problems, a family of 
some objects parameterized by a base scheme $T$ is called {\em
  isotrivial} \index{isotrivial}
if the isomorphism class of the members of the family is constant
over $T$.   From the view-point of morphisms to a moduli stack $\cX$,
this corresponds to the image of $T$ in $|\cX|$ being a singleton,
say $x$.  
We would like to conclude that the morphism $T \to \cX$
factors through the residual gerbe $\cZ_x$.
In practice, we often verify the ``constancy'' of the morphism $T \to \cX$
only at finite type points.
The following result gives sufficient conditions, under such a constancy
hypothesis
on the finite type points, for a morphism to factor through the residual
gerbe.

\begin{alemma}
	\label{lem:isotrivial}
	Let $\cX$ be an algebraic stack,
	and let $x \in |\cX|$ be a point for which
	the residual gerbe $\cZ_x$ exists,
	and for which the canonical monomorphism
	$\cZ_x \hookrightarrow \cX$ is an immersion.
	If $f:T \to \cX$ is a morphism
	whose domain is a reduced scheme,
	and for which all the finite type points of $T$
	map to the given point $x$,
	then $f$ factors through $\cZ_x$.
\end{alemma}
\begin{proof}
	Consider the fibre product $T \times_{\cX} \cZ_x$. 
	By assumption,
	the projection from this fibre product to $T$
	is an immersion 
	whose image contains every finite type point of $T$.
	A locally closed subset of $T$ that contains every finite type
	point is necessarily equal to $T$, and thus, since $T$ is
	reduced, this immersion is in fact an isomorphism.  
	Consequently the morphism $f$ factors through $\cZ_x$,
	as claimed.
\end{proof}

We end this discussion by explaining how the preceding discussion
generalizes to certain Ind-algebraic stacks.  Let $\{\cX_i\}_{i\in I}$
be a $2$-directed system of algebraic stacks, and assume that the 
transition morphisms are monomorphisms.  Let $\cX := \varinjlim_i \cX_i$
be the Ind-algebraic stack obtained as the $2$-direct limit of the $\cX_i$.
We may define the underlying set of points $|\cX|$ as equivalence classes
of morphisms from spectra of fields in the usual way, and since any
such morphism factors through some $\cX_i$, we find that 
$|\cX| = \varinjlim_i |\cX_i|$.  

Now suppose that $\cX$ has quasi-compact diagonal, or, equivalently (since
the transition morphisms are monomorphisms) that each $\cX_i$ has 
quasi-compact diagonal.  If $x \in |\cX|$, then $x \in |\cX_i|$ for 
some $i$, and the residual gerbe $\cZ_x$ at $x$ in $\cX_i$ exists.
The composite monomorphism
$\cZ_x \hookrightarrow \cX_i \hookrightarrow \cX_{i'}$
realizes $\cZ_x$ as the residual gerbe at $x$ in $\cX_{i'}$,
for any $i' \geq i$, and so we may regard $\cZ_x$ as being the 
residual gerbe at $x$ in $\cX$.
Lemmas~\ref{lem:residual gerbe immersion} and~\ref{lem:isotrivial}
immediately extend to the context of such Ind-algebraic stacks $\cX$.

\chapter{Breuil--Kisin--Fargues modules and potentially semistable
  representations (by Toby Gee and Tong Liu)}\chaptermark{BKF modules
  and $G_K$-representations}\label{app: BKF pst}In this appendix we briefly
discuss the relationship between Breuil--Kisin--Fargues modules with
semilinear Galois actions, and potentially semistable Galois
representations. As explained in Remark~\ref{rem: probably don't need all choices of pi} below, we do not expect our results to be
optimal, but they suffice for our applications in the body of the
book. The results of this appendix were originally inspired
by~\cite{MR3127808}; the recent paper~\cite{2019arXiv190508555G}
corrects a mistake in~\cite{MR3127808} and independently proves related (and in some
cases stronger) versions of some of our results, by different
methods. We do not make any use of the arguments of
either~\cite{MR3127808,2019arXiv190508555G}, but instead
combine~\cite{1302.1888} with a result of Fargues (\cite[Thm.\
4.28]{2016arXiv160203148B}). 

We use the notation introduced in the body of the book, in particular
in Section~\ref{subsec:rings}. 
The kernel
of the usual ring homomorphism~$\theta:\Ainf\to\cO_\C$ is a principal ideal~$(\xi)$; one possible choice of~$\xi$
is $\mu/\varphi^{-1}(\mu)$, where ~$\mu=[\varepsilon]-1$ (for some
compatible choice of roots of unity
$\varepsilon=(1,\zeta_p,\zeta_{p^{2}},\dots)\in\cO_\C^\flat$). Recall that $\BdR^+ $ is the $\ker(\theta)$-adic completion of $\Ainf [\frac 1 p]$ and that $t:=\log[\varepsilon]\in \BdR^+$ is a generator of $\ker(\theta)$ in $\BdR^+$. 

As usual, we let~$K$ be a finite extension of~$\Qp$ with residue
field~$k$. Recall that for each choice of uniformizer~$\pi$ of~$K$,
and each choice~$\pi^\flat\in\cO_\C^\flat$ of $p$-power roots of~$\pi$, we
write~$\gS_{\piflat}$ for~$\gS=W(k)[[u]]$, regarded as a subring
of~$\Ainf$ via $u\mapsto[\piflat]$. Write~$E_\pi(u)$ for the
Eisenstein polynomial for~$\pi$, and~$E_{\piflat}$ for its image in~$\Ainf$;
then $E_{\piflat}\in(\xi)$, because $\theta(E_{\piflat})=E_\pi(\pi)=0$.
Indeed $(E_{\pi^\flat}) = (\xi)$ (this follows for example from the
criterion given in~\cite[Rem.\ 3.11]{2016arXiv160203148B} and the
definition of an Eisenstein polynomial). 

In contrast to the body of the book, we do not use coefficients in
most of this appendix (we briefly consider $\cO$-coefficients at the end). Accordingly, we have the following definitions. 
\begin{adefn}\label{adefn: GK phi module}
  An \'etale $(\varphi,G_K)$-module is a finite free $W(\C^\flat)$-module $M$
  equipped with a $\varphi$-semilinear map $\varphi:M\to M$ which
  induces an isomorphism $\varphi^*M\isoto M$, together with a continuous
  semilinear action of~$G_K$ which commutes with~$\varphi$.
\end{adefn}
There is an equivalence of categories between the category of \'etale
~$(\varphi,G_K)$-modules $M$ of rank~$d$ and the category of
free~$\Zp$-modules~$T$ of rank~$d$ which are equipped with a
continuous action of~$G_K$. The Galois representation corresponding
to~$M$ is given by $T(M)=M^{\varphi=1}$.
(See Section~\ref{subsubsec: Galois reps for GK phi}.)

  
\begin{adefn}\label{adefn: BK module}
  Fix a choice of~$\pi^\flat$.  We define a Breuil--Kisin module of
  height at most~$h$ 
  to be a finite free
  $\gS_{\pi^\flat}$-module~$\gM$ equipped with a~$\varphi$-semi-linear
  morphism $\varphi:\gM\to\gM$, with the property that the
  corresponding morphism 
  $\Phi_{\gM}:\varphi^*\gM\to\gM$ is injective, with cokernel killed
  by~$E_{\piflat}^h$.
\end{adefn}

\begin{adefn}
  \label{adefn: BKF module}A Breuil--Kisin--Fargues module of height
  at most~$h$ is a finite free $\Ainf$-module~$\gMt$ equipped with
  a~$\varphi$-semi-linear morphism $\varphi:\gMt\to\gMt$, with the
  property that the corresponding morphism
  $\Phi_{\gMt}:\varphi^*\gMt\to\gMt$ is injective, with cokernel
  killed by~$\xi^h$.
\end{adefn}

\begin{arem}  
If~$\gM$ is a Breuil--Kisin module we write~$\gMt$ for the
Breuil--Kisin--Fargues module $\Ainf\otimes_{\gS_{\pi^\flat}}\gM$
(this is indeed a Breuil--Kisin--Fargues module, because
$\gS_{\pi^\flat}\to\Ainf$ is faithfully flat, and
$(E_{\piflat})=(\xi)$). As noted in Remark~\ref{rem: we don't twist our
    embeddings by phi}, we are not twisting the embedding
  $\gS\to\Ainf$ by~$\varphi$, and accordingly,
  in Definition~\ref{adefn: BKF module},
  we demand that the
  cokernel of~$\varphi$ is killed by a power of~$\xi$, rather than a
  power of~$\varphi(\xi)$ as in~\cite{2016arXiv160203148B}.

  Leaving this difference aside,
  our definition of a
  Breuil--Kisin--Fargues module is less general than that of
  ~\cite{2016arXiv160203148B}, in that we require~$\varphi$ to
  take~$\gM$ to itself; this corresponds to only considering Galois
  representations with non-negative Hodge--Tate weights. This
  definition is convenient for us, as it allows us to make direct
  reference to the literature on Breuil--Kisin modules. The
  restriction to non-negative Hodge--Tate weights is harmless in our
  main results, as we can reduce to this case by twisting by a large
  enough power of the cyclotomic character (the interpretation of
  which on Breuil--Kisin--Fargues modules is explained in~\cite[Ex.\
  4.24]{2016arXiv160203148B}). (We are also only considering free
  Breuil--Kisin--Fargues modules, rather than the more general
  possibilities considered in~\cite{2016arXiv160203148B}.)
\end{arem}

\begin{adefn}\label{adefn: BKF GK module}
  A Breuil--Kisin--Fargues $G_K$-module of height at most~$h$ is a Breuil--Kisin--Fargues
  module of height at most~$h$ which is equipped with a semilinear $G_K$-action which commutes
  with~$\varphi$.
\end{adefn}

\begin{arem}
Note that if~$\gMt$ is a Breuil--Kisin--Fargues $G_K$-module, then
$W(\C^\flat)\otimes_{\Ainf}\gMt$ is naturally an \'etale
$(\varphi,G_K)$-module
in the sense of Definition~\ref{adefn: GK phi module}.
\end{arem}



Recall that for each choice of~$\piflat$ and each~$s\ge 0$ we
write~$K_{\pi^\flat,s}$ for $K(\pi^{1/p^s})$,
and~$K_{\pi^\flat,\infty}$ for~$\cup_sK_{\pi^\flat,s}$. 
\begin{adefn}
\label{adefn: descending BKF to BK}Let $\gMt$ be a
  Breuil--Kisin--Fargues $G_K$-module of height at most~$h$ with a semilinear $G_K$-action. Then
  we say that~$\gMt$ \emph{admits all descents}  if the following
  conditions hold.
  \begin{enumerate}
  \item\label{item: existence of descent Zp version} For every choice of~$\pi$ and~$\piflat$, there is a
    Breuil--Kisin module~$\gM_{\pi^\flat}$ of height at most~$h$ with
    $\gM_{\pi^\flat}\subset(\gMt)^{G_{K_{\pi^\flat,\infty}}}$
    for which the induced morphism
    $\Ainf\otimes_{\gS_{\pi^\flat}}\gM_{\pi^\flat} \to \gMt$
    is an isomorphism.
  \item\label{appendix item: M mod u descends} The $W(k)$-submodule $\gM_{\piflat}/[\piflat]\gM_{\piflat}$ of
  $W(\overline{k})\otimes_{\Ainf}\gMt$ is independent of the
    choice of~$\pi$ and~$\piflat$.
  \item\label{appendix item: M mod E descends} The $\cO_K$-submodule $\varphi ^*\gM_{\piflat}/E_{\piflat}\varphi^*\gM_{\piflat}$
    of  $\cO_\C\otimes_{\theta, \Ainf}\varphi^*\gMt$  is independent of the
    choice of~$\pi$ and~$\piflat$.
  \end{enumerate}
\end{adefn}

\begin{arem}
  \label{arem: second condition is redundant}In fact
  condition~\eqref{appendix item: M mod u descends} in Definition~\ref{adefn:
    descending BKF to BK} is redundant; see
  Remark~\ref{arem: explanation that we didn't need a condition in BKF
    but that we want it for coefficients} below. However, we include
  the condition as it is useful when considering versions of the
  theory with coefficients and descent data.
\end{arem}

\begin{adefn}
  \label{adefn: crystalline descending BKF module}Let~$\gMt$ be a
  Breuil--Kisin--Fargues $G_K$-module which
  admits all descents. We say that~$\gMt$ is furthermore
  \emph{crystalline} if for each choice of~$\pi$ and~$\piflat$, and each~$g\in G_K$, we
  have \anumequation\label{appendix eqn: condition on BK GK for
    crystalline}(g-1)(\gM_{\piflat})\subset
  \varphi ^{-1}(\mu)[\piflat]\gMt.  \end{equation}
\end{adefn}

Definitions~\ref{adefn: descending BKF to BK} and~\ref{adefn:
  crystalline descending BKF module} are motivated by the
following result, whose proof occupies most of the rest of this
appendix.

\begin{athm}
  \label{athm: admits all descents if and only if semistable} Let~$M$
  be an \'etale $(\varphi,G_K)$-module. Then~$V(M)$ is semistable with
  Hodge--Tate weights in~$[0,h]$ if and only if there is a
  \emph{(}necessarily unique\emph{)}
  Breuil--Kisin--Fargues $G_K$-module $\gMt$ which is of height at
  most~$h$, which admits all descents, and which satisfies
  $M=W(\C^\flat)\otimes_{\Ainf}\gMt$.

  Furthermore, $V(M)$ is
  crystalline if and only if~$\gMt$ is crystalline.
\end{athm}
\begin{arem}
  \label{rem: probably don't need all choices of pi}
  As already noted in Remark~\ref{arem: second condition is redundant},
  condition~\eqref{appendix item: M mod u descends}
  of Definition~\ref{adefn: descending BKF to BK} is
  redundant, and it is plausible that
  condition~\eqref{appendix item: M mod E descends} is redundant as well
(i.e.\ that both condition
  condition~\eqref{appendix item: M mod u descends}
  and condition~\eqref{appendix item: M mod E descends} 
 in Definition~\ref{adefn: descending BKF to BK} are
  consequences of condition~\eqref{item: existence of descent Zp version}); but
  this does not seem to be obvious. Note though that it is not
  sufficient to demand the existence of a descent for a single choice
  of~$\pi$, as there are representations of finite height which are
  potentially semistable but not semistable (see for
  example~\cite[Ex.\ 4.2.1]{MR2558890}).
\end{arem}

We begin with some preliminary results. The following proposition and
its proof are due to Heng Du, and we thank him for allowing us to
include them here.

\begin{aprop}
  \label{aprop: BKF is de Rham if and only if GK basis mod E}Let~$\gMt$
  be a Breuil--Kisin--Fargues $G_K$-module with the property that
  $\C\otimes_{\theta , \Ainf}\varphi^*\gMt$ has a basis consisting of $G_K$-fixed 
  vectors. Let~$M=W(\C^\flat)\otimes_{\Ainf}\gMt$; then the
  $G_K$-representation ~$V(M)$ is de Rham.
\end{aprop}
\begin{arem}
  \label{arem: de Rham implies existence of lattice}Using the
  equivalence of categories of~\cite[Thm.\ 4.28]{2016arXiv160203148B} (a theorem of Fargues),
  one can easily check that Proposition~\ref{aprop: BKF is de Rham if
    and only if GK basis mod E} admits a converse: namely that if~$M$
  is an \'etale $(\varphi,G_K)$-module with the property that $T(M)$ is a
  $\Zp$-lattice in a de Rham representation of~$G_K$ with non-negative
  Hodge--Tate weights, then there is a Breuil--Kisin--Fargues
  module~$\gMt$ with the properties in the statement of
  Proposition~\ref{aprop: BKF is de Rham if and only if GK basis mod
    E}.
\end{arem}
\begin{proof}[Proof of Proposition~\ref{aprop: BKF is de Rham if and
    only if GK basis mod E}]
By~\cite[Thm.\ 4.28]{2016arXiv160203148B} (and our assumption that $\varphi(\gMt)\subseteq\gMt$), we have  injections 
  \[\BdR^+\otimes_{\Ainf}\varphi^*\gMt\hookrightarrow \BdR^+\otimes_{\Ainf}\gMt
  \hookrightarrow \BdR\otimes_{\Zp}T(M).\] To show
  that~$V(M)$ is de Rham, we need to show that the
   $\BdR$-vector space $\BdR\otimes_{\Zp}T(M)$ has a 
  basis consisting of $G_K$-fixed vectors,
  so it suffices to show that the $\BdR^+$-lattice
  $\BdR^+\otimes_{\Ainf}\varphi^*\gMt$ has a 
  basis consisting of $G_K$-fixed vectors. 
  Since~$\BdR^+$ is a complete discrete valuation ring with maximal
  ideal~$(\xi)$, it is enough to show that there are compatible
bases of
  $\BdR^+/(\xi^n)\otimes_{\Ainf}\varphi^*\gMt$ for all~$n\ge 1$,  consisting of   $G_K$-fixed vectors.

  In the case~$n=1$, since $\xi$ generates the kernel of~$\theta$,
we have such a
  basis by hypothesis. Suppose that~$\{e_i^{(n)}\}_{i=1,\dots,d}$ is a
  basis of $G_K$-fixed vectors for $\BdR^+/(\xi^n)\otimes_{\Ainf}\varphi^*\gMt$, and
  let~$\{\te_i^{(n+1)}\}_{i=1,\dots,d}$ be any basis of $\BdR^+/(\xi^{n+1})\otimes_{\Ainf}\varphi^*\gMt$
  lifting~$\{e_i^{(n)}\}$. Recall that~$t=\log[\varepsilon]\in \BdR^+$ is a
  generator of~$(\xi)$, and for each $g\in G_K$, write \[g\cdot
    (\te_1^{(n+1)},\dots,\te_d^{(n+1)})=(\te_1^{(n+1)},\dots,\te_d^{(n+1)})(1_d+t^nB_g^{(n+1)})\]where
  we can view~$B_g^{(n+1)}$ as an element of~$M_d(B^+_{\dR}/(\xi))=M_d(\C)$.

  A simple calculation shows that $g\mapsto B_g^{(n+1)}$ is a
  continuous 1-cocycle valued in $M_d(\C(n))$. Since~$n\ge 1$, the
  corresponding cohomology group $H^1(G_K,M_d(\C(n)))$ vanishes,
  because $H^1(G_K,\C(n))=0$
  by a theorem of Tate--Sen. There is therefore some $A\in M_d(\C)$ such
  that for all~$g$ we have \[B_g^{(n+1)}=\chi^n(g)g(A)-A,\]
where $\chi$ is the $p$-adic cyclotomic character.   Then \[
    (e_1^{(n+1)},\dots,e_d^{(n+1)})=(\te_1^{(n+1)},\dots,\te_d^{(n+1)})(1_d-t^nA)\]is
  the required  basis of
  $\BdR^+/(\xi^{n+1})\otimes_{\Ainf}\varphi^*\gMt$ consisting of
  $G_K$-fixed vectors
 lifting~$\{e_i^{(n)}\}$.  
\end{proof}

\begin{alem}
  \label{alemma: G K infty generate G K}There is no proper closed
  subgroup of~$G_K$ containing all of the subgroups~$G_{K_{\piflat,\infty}}$.
\end{alem}
\begin{proof}Let~$s$ be the greatest integer with the property
  that~$K$ contains a primitive $p^s$th root of unity; equivalently,
  it is the greatest integer with the property that~$K_{\piflat,s}/K$
  is Galois over~$K$ (for one, or equivalently every, choice
  of $\pi$ and~$\piflat$). Then~$K_{\piflat,s'}$ depends only on~$\pi$
if $s' \leq s$, and so we write~$K_{\pi,s'}=K_{\piflat,s'}$ for such~$s'$.

  The subgroups $G_{K_{\piflat,\infty}}$
  for any fixed choice of~$\pi$ topologically generate~$G_{K_{\pi,s}}$. It therefore
  suffices to show that as~$\pi$ varies, the various normal
  subgroups~$G_{K_{\pi,s}}$ of $G_K$
  collectively generate~$G_K$. Let $H$ be the normal subgroup that they generate;
  then for every uniformizer~$\pi$, $G_K/H$ is a quotient of $G_K/G_{K_{\pi,s}}$,  a cyclic group  of order~$p^s$.
  If $H$ were a proper subgroup of $G_K$, then necessarily $s \geq 1$,
  and the subgroups $G_{K_{\pi,1}}$ would coincide for every~$\pi$ (since they
  would coincide with  the unique index~$p$ subgroup of~$H$).
  Thus the
  extensions~$K_{\pi,1}/K$ would have to all coincide. 
By Kummer theory,
  this would imply that the ratio of any two uniformizers of~$K$ is a
  $p$th power in~$K$, which is nonsense (for example, consider the
  uniformizers~$\pi$ and~$\pi+\pi^2$).
  \end{proof}

\begin{proof}[Proof of Theorem~\ref{athm: admits all descents if and
    only if semistable}]We begin with the semistable case. If~$V(M)$ is semistable then the existence
  of~$\gMt$ is a straightforward consequence of the results
  of~\cite{1302.1888}. In particular, the existence of a unique~$\gMt$
  satisfying condition~\eqref{item: existence of descent Zp version} of
  Definition~\ref{adefn: descending BKF to BK} follows from ~\cite[Thm.\
  2.2.1]{1302.1888}. 
    
  The proofs that conditions~\eqref{appendix item: M mod u descends} and
  \eqref{appendix item: M mod E descends} hold are implicit in the proof of 
  \cite[Prop. 4.2.1]{1302.1888}, as we now explain. 
  Write
  $\overline{\gM}_{\pi ^\flat}= \gM_{\pi ^\flat}/ [\pi ^\flat]
  \gM_{\pi^\flat} $. Since~$W(\overline{k})\otimes_{W(k)}\overline{\gM}_{\pi
    ^\flat}=W(\overline{k})\otimes_{\Ainf}\gMt$ is independent
  of~$\piflat$, it suffices to show that the $K_0$-vector space
  $D_{\pi^\flat} :=\overline{\gM}_{\pi ^\flat}[\frac 1 p] $ is
  independent of $\pi^\flat$.
  Let $S_{\pi^\flat} $ be the $p$-adic
  completion of
  $\gS_{\pi ^\flat} [\frac{E_{\pi^\flat}^i}{i!}, i \geq 1]$. By
  \cite[Prop.\ 6.2.1.1]{BreuilGriffith} and  \cite[\S 2]{liulattice2},
  $D_{\pi^\flat}$ admits a section
  $s_{\pi ^\flat}: D_{\pi ^\flat} \to S_{\pi ^\flat}[\frac 1 p]
  \otimes_{\varphi, \gS_{\pi ^\flat}}\gM_{\pi ^\flat}$ so that
\anumequation
\label{eqn:strongly divisible iso}
   S_{\piflat}[\frac 1 p] \otimes _{K_0}s_{\piflat}(D_{\piflat}) =
  S_{\pi ^\flat}[\frac 1 p] \otimes_{\varphi, \gS_{\pi
      ^\flat}}\gM_{\piflat}.  
\end{equation}
In
  particular, we may regard $s_{\pi^\flat} (D_{\pi ^\flat})$ as a
  submodule of
  $\Bcris^+ \otimes_{\varphi, \gS_{\piflat}} \gM_{\piflat} \subseteq
   \Bcris^+ \otimes_{\Z_p}
   T(M)$.

 Write $\gu = \log([\pi^\flat]) \in \Bst^+= \Bcris^+ [\gu]$. Note that
 the \emph{set} $\Bst^+$ does not depend on the choice of $\piflat$,
 because, if ${\varpi^\flat}$  is another choice and we write $\gu' =
 \log ([\varpi^\flat])$, then $\Bcris^+ [\gu'] = \Bcris^+
 [\gu]$. Indeed, we have  $\log([{\varpi}^\flat]) = \log([\piflat]) +
 \lambda$ with $\lambda= \log([\frac{\varpi^\flat }{\piflat}]) \in
 \Bcris^+$. In particular, $D_\st (T(M)) : = ( \Bst^+ \otimes_{\Z_p} T(M) )^{G_K} $ is independent of~ $\piflat$.
   Furthermore,
  \cite[Prop. 2.6]{liulattice2} shows that we have a commutative diagram
   (where $i_{\piflat}$ is a $K_0$-linear isomorphism)
    $$\xymatrix{D_\st (T(M)) \ar[d]^-{i_{\piflat}}_\wr \ar@{^{(}->}[r]
      &   \Bst^+\otimes_{K_0}D_\st (T(M))\ar[d]^{\mod \gu} \ar@{^{(}->}[r] &
      \Bst^+  \otimes _{\Z_p} T(M) \ar[d]^{\mod \mathfrak u}
      \\s_{\piflat}(D_{\piflat}) \ar@{^{(}->}[r] & \Bcris^+
      \otimes_{\varphi, \gS_{\piflat}} \gM_{\piflat} \ar@{^{(}->}[r]&  \Bcris^+
      \otimes _{\Z_p} T(M)}$$    
 To show \eqref{appendix item: M mod u descends}, it therefore suffices to show that the
 image of the composite $$ \resizebox{\textwidth}{!}{\xymatrix {D_\st (T(M)) \to
   \Bst^+\otimes_{K_0} D_\st (T(M)) \ar[r]^-{\mod \gu\ \ } &
   \Bcris^+ \otimes_{\varphi, \gS_{\piflat}} \gM_{\piflat} \to W(\bar
   k)[1/p] \otimes_{\varphi, \Ainf} \gMt}} $$ is independent of
 $\piflat$. Here the last map is induced by $\nu : \Bcris^+ \to
 W(\bark)[\frac 1 p]$ which extends the natural projection $\Ainf \to
 W(\bark)$. It suffices in turn to show that the composite
 $\xymatrix{\Bst^+ \ar[r]^-{\mod \gu} & \Bcris^+ \ar[r]^-\nu &
   W(\bark ) [\frac 1 p]}$ is independent of $\piflat$. This follows
 from the fact that $ \lambda= \log([\frac{\varpi^\flat }{\piflat}])$
 is in $\ker (\nu)$ (see the proof of~\cite[Lem.\ 2.10]{liulattice2}). 
 
 To prove \eqref{appendix item: M mod E descends}, it again suffices to prove
 the statement after inverting~$p$. For any subring $A \subset
 \BdR^+$, set $F^1 A = A \cap \xi \BdR^+$. It is easy to see that $F^1 \gS_{\piflat}= E_{\piflat}\gS_{\piflat}$ and also that $S_{\piflat}/ F^1 S_{\piflat} = \cO_K$.
(Note that
the inclusion $E_{\piflat} S_{\piflat} \subseteq F^1 S_{\piflat}$ is strict.) Now we again use the isomorphism~\eqref{eqn:strongly divisible iso}.
By reducing modulo $F^1 S_{\piflat}$ on the both sides of this identification, 
we conclude that $$K \otimes _{K_0} D_{\piflat} =  K \otimes_{K_0}(s_{\piflat}(D_{\piflat}) \mod [\piflat]) =  \varphi^* \gM_{\piflat} / E_{\piflat} \varphi^* \gM_{\piflat} [\frac 1 p].  $$
On the other hand, by tensoring $\Ainf$ via $\gS_{\piflat}$ to~\eqref{eqn:strongly
divisible iso},
we obtain the isomorphism  $\Bcris^+ \otimes _{K_0} s_{\piflat}(D_{\piflat})  =
      \Bcris^+ \otimes_{\Ainf }\varphi^*\gMt$.  Modulo $F^1 \Bcris^+$
      on both sides, a similar argument to the above shows that 
      $ \C \otimes_{K_0} D_{\piflat} = (\varphi^*\gMt / \xi \varphi^* \gMt) [\frac 1 p] =\C \otimes_{W{(\bark)}}\overline{\varphi^*\gMt}
       $, where we write $ \overline{\varphi^*\gMt}
 =W(\bark)  \otimes_{\Ainf}\varphi^* \gMt
 $. In summary we see that 
 \[\resizebox{\textwidth}{!}{$(\varphi^*\gM_{\piflat}/ E_{\piflat} \varphi^*\gM_{\piflat})[\frac 1 p ] = K \otimes_{K_0}D_{\piflat} \subseteq  \C\otimes _{W(\bar k)}\overline{\varphi^*\gMt}  = (\varphi^*\gMt/ \xi \varphi^*\gMt) [\frac 1 p].$}\]  
 Since we have shown that $D_{\piflat} \subseteq \overline{\varphi^*\gMt}$ is independent of
 $\piflat$, it follows that $\varphi^*\gM_{\piflat}/ E_{\piflat}
 \varphi^*\gM_{\piflat}$ is independent of $\piflat$, as required.

  Suppose conversely that we are given~$\gMt$ as in Definition~\ref{adefn: descending BKF to BK}. Write~$\overline{\gM}$
  for the $W(k)$-module of Definition~\ref{adefn: descending BKF to
    BK}~\eqref{appendix item: M mod u descends}, and $\overline{\gM}'$
  for the $\cO_K$-module of Definition~\ref{adefn: descending BKF to
    BK}~\eqref{appendix item: M mod E descends}.   Since for
  each~$\piflat$ the action of~$G_{K_{\piflat,\infty}}$
  on~$\gM_{\piflat}$ is trivial by assumption, it follows from
  Lemma~\ref{alemma: G K infty generate G K} that   $G_K$ acts
  trivially on $\overline{\gM}$ and~$\overline{\gM}'$.

  Since $\C\otimes_{\Ainf,\theta}\varphi^*\gMt=\C\otimes_{\cO_K}\varphi^*\overline{\gM}'$,
 we see
  that the hypotheses of Proposition~\ref{aprop: BKF is de Rham if and
    only if GK basis mod E} are satisfied, so that
  $V(M)$ is de Rham, and consequently potentially
  semistable. To show that $V(M)$ is semistable, we
  may replace~$K$ with an (infinite) unramified extension and assume that~$K$ is
  a complete discretely valued field with residue field $k=\bark$. In
  particular, we now have $\overline{\gM} = W(\bark) \otimes_{\Ainf}
  \gMt $, and we claim that $\overline{\gM}[\frac 1 p] $ with its $G_K$-action is isomorphic to $D_{\pst}(V(M))$.  
  If the claim holds,  then since $G_K$ acts trivially on~$\overline{\gM}$, $G_K$ acts
  trivially on~$D_{\pst}(V(M))$, so
  that~$V(M)$ is semistable, as required.
  
  We now prove the claim. Let $L /K$ be a finite Galois extension so that
  $V(M)|_{G_L}$ is semi-stable, and let $\gM_L$ be the Breuil--Kisin
  module attached to $T(M)|_{G_L}$ for some choice of~$\piflat_L$. By~ \cite[(2.12)]{liulattice2}, 
  $ \Ainf\otimes_{\gS_{\piflat_L}}\varphi^*\gM_L  $ injects into $\Ainf
  \otimes_{\Z_p}T(M)$. Furthermore, by~\cite[Lem.\ 2.9]{liulattice2}, $\varphi^*{\gM}^{\inf}_L :=
  \Ainf \otimes_{\gS_{\piflat_L}}\varphi^*\gM_L   $ is stable under the
  $G_K$-action on $ \Ainf\otimes_{\Z_p}T(M)$, and by~\cite[Cor.\ 2.12]{liulattice2}, $ W(k) \otimes_{\gS_{\piflat_L}}\varphi^*\gM_L = 
  W(k) \otimes_{\Ainf}\varphi^*{\gM}^{\inf}_L $ (recall that $k = \bark$) together with its $G_K$-action
  is isomorphic to a $W(k)$-lattice inside $D_{\st, L} (T(M))= (\Bst^+
  \otimes_{\Z_p} T(M)) ^{G_L} $. In summary, $W(k)\otimes_{\Ainf}\varphi^*
  {\gM}^{\inf}_L [1/p] $
is isomorphic to $D_{\st, L} (T(M))$ as $G_K$-modules.

  It remains to show that $\varphi
  ^*{\gM}^{\inf}_L = \varphi ^*\gMt$. By
  \cite[Thm.\ 4.28]{2016arXiv160203148B}, it suffices to show that
  $ \BdR^+ \otimes_{\Ainf}\varphi^*{\gM}^{\inf}_L= \BdR^+\otimes_{\Ainf}  \varphi^* \gMt
  $.  By the proof of
    Proposition~\ref{aprop: BKF is de Rham if and only if GK basis mod
      E}, we see that
    \anumequation\label{eqn: first check on L Minf}
      \BdR^+\otimes_{\Ainf}\varphi^* \gMt = \BdR^+ \otimes_K
      D_\dR(T(M)).
    \end{equation}
Since $T(M)$ is semi-stable over
   $L$, it is well-known that \anumequation\label{eqn: second check on L Minf}
\Bst^+   \otimes_{\gS_{\piflat_L}} \varphi^*\gM_L = 
\Bst^+   \otimes_{W(k _L)} D_\st (T(M)|_{G_L}).\end{equation} (This can be easily seen from the proof
   of \cite[Cor.1.3.15]{KisinCrys}, or see \cite[\S 2]{liulattice2}
   for a more detailed discussion; the key point,  in terms of  the
   diagram  above, is that \cite[\S 2]{liulattice2} shows that $K_0 [\gu] \otimes_{K_0}D_\st (T(M)) =
   K_0 [\gu] \otimes_{K_0}s_{\piflat} (D_{\piflat})$.) 

From~\eqref{eqn: second check on L Minf} we obtain    
 \anumequation\label{eqn: third check on L Minf} \BdR^+\otimes_{ \gS_{\piflat_L} }\varphi^*\gM_{L} = \BdR^+ \otimes_{L}D_\dR (T(M)|_{G_L}).  \end{equation}
Comparing~\eqref{eqn: first check on L Minf} and~\eqref{eqn: third
  check on L Minf} we have  $ \BdR^+ \otimes_{\Ainf}\varphi^*{\gM}^{\inf}_L= \BdR^+\otimes_{\Ainf}  \varphi^* \gMt
  $, as required.

Finally, we turn to the crystalline case. Given the above, the result
  is a consequence of~\cite[Thm.\ 3.8]{Ozeki2018}, as we
 now explain.  In our case, $f(u)$ in~\cite{Ozeki2018} is equal to
 $u^p$, so  the assumptions of~\cite[Thm.\ 3.8]{Ozeki2018} are
 automatically satisfied (as noted at the beginning of~\cite[\S
 3.2]{Ozeki2018}).  More precisely, fix some~$\piflat$, and let $\widehat K$ be the Galois
 closure of~$K _{\piflat}$. There is a subring
 $\widehat{\mathcal R} \subseteq \Ainf$ (constructed in~\cite{Ozeki2018}) such that $\gS_{\piflat} \subset
 \widehat{\mathcal R}$ and $G_K$ acts on $\widehat{\mathcal R}$ through
 $\widehat G:=\Gal(\cup_{n\ge 1}K(\zeta_{p^n},\pi^{1/p^n})/K)$.
When $T(M)$ is crystalline, ~\cite[Thm.\
 3.8]{Ozeki2018} shows that $T(M)$ admits a $(\varphi, \widehat
 G)$-module, which by definition consists of the following data: 
  \begin{itemize}
  \item The Breuil-Kisin module $\gM_{\piflat}$ attached to $T(M)|_{G_{{\piflat},\infty}}$.
  \item A $\widehat G$-action on $\widehat{\gM}: = \widehat{\mathcal
      R} \otimes_{\gS_{\piflat}} \varphi ^* \gM_{\piflat}$ which
    commutes with the action of~$\varphi$ and satisfies~\eqref{appendix eqn: condition on BK GK for
    crystalline}.
  \item $(W(\C^\flat) \otimes_{\widehat{\mathcal R}} \widehat{\gM})^{\varphi=1} = T(M)$ as $G_K$-modules.  
  \end{itemize}(In fact \cite{Ozeki2018} uses contravariant functors,
  which can be easily translated to the covariant functors used here.)
 This proves that if $T(M)$ is crystalline then $\gMt$
 satisfies~\eqref{appendix eqn: condition on BK GK for
    crystalline}. Conversely, if $\gMt$ satisfies~\eqref{appendix eqn: condition on BK GK for
    crystalline} for \emph{one} fixed $\piflat$, then
 \cite[Lem.\ 3.15]{Ozeki2018} shows that $s_{\pi}(D_{\piflat}) \subset
 (\Bcris^+ \otimes_{\gS_{\piflat}} \varphi^*
 \gM_{\piflat})^{G_K}$. Hence $T(M)$ is crystalline, as required. \end{proof}
\begin{arem}\label{arem: explanation that we didn't need a condition
    in BKF but that we want it for coefficients}
  Note that condition~\eqref{appendix item: M mod u descends} in Definition \ref{adefn: descending BKF to
    BK} is redundant for proving Theorem \ref{athm: admits all
    descents if and only if semistable}. 
  In fact, if we
  remove~\eqref{appendix item: M mod u descends} from Definition~\ref{adefn: descending BKF to BK} then
  the necessity part of Theorem~\ref{athm: admits all descents if and
    only if semistable} of course still holds. For the sufficiency,
  note that in the above proof, we extended $K$ so that the residue
  field is $\bar k$; so we only use that $G_K$ (the inertia subgroup
  in this situation) acts trivially on
  $W(\bar k) \otimes _{\Ainf} \gMt$, and we do not need that
  $\gM_{\piflat}/ [\piflat] \gM_{\piflat} \subseteq W(\bar k) \otimes _{\Ainf} \gMt $ is independent
  of $\piflat$.

However, it is not immediately clear that the analogous condition is
redundant in the version of the theory with coefficients that we
consider in the body of the book (see Definition~\ref{defn:
  descending BKF to BK coefficients}), and this condition is used in defining our moduli stacks of potentially semistable
representations of given inertial type, so we include it here.
\end{arem}

\section{Potentially semistable representations}
It will be convenient for us to have a slight refinement of these
results, allowing us to discuss potentially semistable (and
potentially crystalline) representations. To this end, fix a finite Galois
extension~$L/K$. 

\begin{adefn}
  \label{adefn: potentially admits descents}Let $\gMt$ be a
  Breuil--Kisin--Fargues $G_K$-module of height at most~$h$. Then we say that~$\gMt$ \emph{admits all
    descents over~$L$} if the corresponding Breuil--Kisin--Fargues
  $G_L$-module (obtained by restricting the $G_K$-action on~$\gMt$ to~$G_L$) admits all descents (in the sense of
  Definition~\ref{adefn: descending BKF to BK}).
\end{adefn}

\begin{acor}
  \label{acor: admits all descents if and only if potentially semistable} Let~$M$
  be an \'etale $(\varphi,G_K)$-module. Then~$V(M)|_{G_L}$ is semistable with
  Hodge--Tate weights in~$[0,h]$ if and only if there is a \emph{(}necessarily unique\emph{)}
  Breuil--Kisin--Fargues $G_K$-module $\gMt$ which is of height at
  most~$h$, which admits all descents over~$L$, and which satisfies
  $M=W(\C^\flat)\otimes_{\Ainf}\gM$. Furthermore~$T(M)|_{G_L}$ is
  crystalline if and only if~$\gMt$ is crystalline as a
  Breuil--Kisin--Fargues $G_L$-module.
\end{acor}
\begin{proof}
  The sufficiency of the condition is immediate from
  Theorem~\ref{athm: admits all descents if and only if
    semistable}. For the necessity, by Theorem~\ref{athm: admits all
    descents if and only if semistable} there is a
  Breuil--Kisin--Fargues $G_L$-module $\gMt_L$ which admits all
  descents and satisfies $M=W(\C^\flat)\otimes_{\Ainf}\gM_L$, so we
  need only show that~$\gMt_L$ is $G_K$-stable. This follows from
  \cite[Lem.\ 2.9]{liulattice2}, or one can argue as
  follows: 
  note that the proof of Theorem~\ref{athm: admits all descents if and
    only if semistable} shows that the Breuil--Kisin--Fargues module
  $\gMt_L : =\Ainf\otimes_{\gS_{\piflat_L}}\gM_L$ corresponds via
  \cite[Thm. 4.28]{2016arXiv160203148B} to the pair
  $(T(M), \BdR^+ \otimes_{L}D_\dR(T(M)|_{G_L}))$. This pair has a
  $G_K$-action, because $V(M)$ is de Rham, hence $\gM_L^{\inf}$ is
  $G_K$-stable by \cite[Thm. 4.28]{2016arXiv160203148B}.\end{proof}


\section{Hodge and inertial types}\label{subsec: Hodge and inertial
types}
We finally recall how to interpret Hodge--Tate weights and
inertial types in terms of Breuil--Kisin modules. Fix a finite
extension $E/\Qp$ with ring of integers~$\cO$, which is sufficiently
large that~$E$ contains the image of every embedding
$\sigma:K\into\Qpbar$. The definitions above admit obvious extensions
to the case of Breuil--Kisin--Fargues modules with $\cO$-coefficients
(see Definition~\ref{defn: descending BKF to BK coefficients}),
and since $\cO$ is a finite free~$\Zp$-module, the proofs of
Theorem~\ref{athm: admits all descents if and only if semistable} and
Corollary~\ref{acor: admits all descents if and only if potentially
  semistable} go
over unchanged in this setting.
  

Suppose that~$\gMt$ is as in the statement of Corollary~\ref{acor:
  admits all descents if and only if potentially semistable}; so it is
a Breuil--Kisin--Fargues $G_K$-module of height at most~$h$, which
admits all descents over~$L$. Write~$l$ for the residue field of~$L$,
write~$\overline{\gM}$ for the $W(l)$-module
$\gM_{\piflat}/[\piflat]\gM_{\piflat}$ of Definition~\ref{adefn:
  descending BKF to BK}~\eqref{appendix item: M mod u descends}, and
$\overline{\gM}'$ for the $\cO_L$-module
$\varphi^*\gM_{\piflat}/E_{\piflat}\varphi^*\gM_{\piflat}$ of
Definition~\ref{adefn: descending BKF to BK}~\eqref{item: M mod E
  descends}. These modules are independent of the choice of~$\pi$
(which now denotes a uniformizer of~$L$) and of~$\piflat$, by definition.

The semilinear $G_K$-action on $\gMt$ induces semilinear actions
on~$\overline{\gM}$ and on~$\overline{\gM}'$, and as noted in the
proof of Theorem~\ref{athm: admits all descents if and only if
  semistable}, the action of $G_L$ on both modules is
trivial. Furthermore, by~\cite[Cor.\ 2.12]{liulattice2}, the inertial
  type $D_{\pst}(V(M))|_{I_K}$ is given by~$\overline{\gM}[1/p]$ with
  its action of~$I_{L/K}$. More precisely,  \cite[\S 2.3]{liulattice2} shows that 
  $W(\bark)[1/p] \otimes_{\Ainf}\varphi ^*\gMt$ is isomorphic to $D_{\pst}
  (T(M))\simeq W(\bark) \otimes_{W(l)}D_\st(T(M)|_{G_L})$ as $W(\bark)[1/p][G_K]$-modules. Furthermore, this isomorphism is compatible with the isomorphism (from the proof of Theorem \ref{athm: admits all descents if and only if semistable})
  \[\iota _{\piflat}: \xymatrix{\overline \gM [\frac 1 p] = D_{\piflat}
    \simeq s_{\piflat} (D_{\piflat}) & \ar[l]_-{i_{\piflat}}^-{\sim}
    D_\st(T(M)|_{G_L})}.\] Hence $\overline \gM [\frac 1 p]  \subset
  W(\bark)[1/p] \otimes_{\Ainf}\varphi ^*\gMt$ is endowed with an action of $I_{L/K}$ and is isomorphic to $D_{\pst} (V(M))$.


We now turn to the Hodge--Tate weights. We have a filtration on the $L\otimes_{\Qp}E$-module
$D_L=\overline{\gM'}[1/p]$, which is defined as follows. 
For each~$\piflat$ we write
$\Phi_{\gM_{\piflat}}:\varphi^*\gM_{\piflat}\to\gM_{\piflat}$ and $f_\pi: \varphi^*\gM_{\piflat}[\frac 1 p] \twoheadrightarrow D_L $.  
For each~$ i\geq 0$ we define
$\Fil^i\varphi^*\gM_{\piflat}=\Phi_{\gM_{\piflat}}^{-1}(E_{\piflat}^i\gM_{\piflat})$ 
and $\Fil^i D_L = f_\pi (\Fil ^i \varphi ^* \gM_{\piflat})$. 
Considering the isomorphism $\iota_{\piflat}: D_\st (V(M)|_{G_L}) \simeq D_{\piflat}$,  
we obtain an isomorphism $D_\dR(V(M)|_{G_L}) \simeq L
\otimes_{W(l)[\frac 1 p]} D_\st (V(M)|_{G_L}) \simeq D_L$. By
\cite[Cor.\ 3.2.3, Thm.\ 3.4.1]{MR2388556}, this isomorphism respects
the  filtrations on each side. In particular, $\Fil^i D_L$ is independent of the choice of~$\piflat$. 

By Hilbert 90, $D_L$ and
its filtration descend to a filtration on a
$K\otimes_{\Qp}E$-module~$D_K\simeq D_\dR(V(M))$. 
Then for each $\sigma:K\into E$, and
each~$i\in [0,h]$, the multiplicity of~$i$ in the multiset
$\HT_{\sigma}(V(M))$ is the dimension of the $i$-th graded piece of
$e_{\sigma}D_K$, where $e_\sigma\in (K\otimes_{\Qp}E)$ is the
idempotent corresponding to~$\sigma$.

  We summarise the preceding discussion in the following corollary.
\begin{acor}\label{acor: reading off inertial type and HT weights from BKF}
  In the setting of Corollary~\ref{acor: admits all descents if and
    only if potentially semistable}, the inertial type and Hodge type
  of~$V(M)$ are determined as follows: the inertial
  type~$D_{\pst}(V(M))|_{I_K}$ is given by~$\overline{\gM}[1/p]$ with
  its action of~$I_{L/K}$, while the Hodge type of~$V(M)$ is given by the
  jumps in the filtration on~$D_K$ described
  above.
\end{acor}

\emergencystretch=3em
\bibliographystyle{amsalpha}
\bibliography{universalBM}

\printindex
\end{document}